\documentclass[11pt, reqno,amsmath,amsthm,amssymb,amscd]{amsart}
\usepackage{amsmath,amssymb}
\usepackage{graphics}
\usepackage[usenames]{color}
\baselineskip 3.1pt \lineskip 3.1pt \numberwithin{equation}{section}

\def\FF{{\,\hat{\textbf{\textit{$\!$F$\ssc\,$}}}}}
\def\CA{{\textbf{CA }}}\def\CA{converges absolutely }

\def\aa{{m_0}}\def\rDS{\unrhd}

\def\lin{^{\rm lin}}

\def\ign{_{\rm igno}}
\def\igo{{\rm igo}}
\def\ign{\igo}\def\igo{{\rm ign}}

\def\inv{{\rm neg}}

\def\th{\theta}
\def\DS{\unlhd}

\def\Coeff{{\rm C_{oeff}}}

\def\ptl{\partial}\def\cC{{\textbf{\textit{c}}}}

\def\deg{{\rm deg}}
\def\max{{\rm max}}

\def \Z{\hbox{$Z\hskip -5.2pt Z$}}

\def \Q{\hbox{$Q\hskip -5pt \vrule height 6pt depth 0pt\hskip 6pt$}}
\def \R{\mathbb{R}}
\def \C{\hbox{$C\hskip -5pt \vrule height 6pt depth 0pt \hskip 6pt$}}

\def\qed{\ \ \ifhmode\unskip\nobreak\fi\ifmmode\ifinner
         \else\hskip5pt\fi\fi
 \hbox{\hskip5pt\vrule width4pt height6pt depth1.5pt\hskip 1 pt}}
\def\a{\alpha}
\def\b{\beta}
\def\d{\delta}

\def\g{\mbox{${\ssc\!}$\Large$\gamma$}}

\def\l{\lambda}

\def\si{\sigma}
\def\sc{\scriptstyle}
\def\ssc{\scriptscriptstyle}
\def\dis{\displaystyle}

\def\ol{\overline}

\def\bs{\backslash}

\def\VS#1{}

\def\PP{{\cal P}}\def\PP{\C[x]((y^{-1}))}

\def\C{\mathbb{C}}

\def\Z{\mathbb{Z}}
\def\bB{{\textbf{\textit{b}}}}
\def\zZ{{\textbf{\textit{z}}}}
\def\aA{{\textbf{\textit{a}}}}
\def\lL{{\textbf{\textit{l}}}}
\def\hH{{\textbf{\textit{h}}}}
\def\dD{{\textbf{\textit{d}}}}
\def\eE{{\textbf{\textit{e}}}}

\def\Q{\mathbb{Q}}
\def\SS{{\textbf{\textit{s}}}}
\def\bS{{\textbf{\textit{F}}}}
\def\bP{{\textbf{\textit{P}}}}\def\kk{\uu}\def\uu{{\textbf{\textit{k}}}}
\def\nn{{\textbf{\textit{n}}}}

\def\dD{{\textbf{\textit{d}}}}
\def\xX{{\textbf{\textit{x}}}}
\def\zZ{{\textbf{\textit{z}}}}

\def\yY{{\textbf{\textit{y}}}}

\def\aA{{\textbf{\textit{a}}}}
\def\bB{{\textbf{\textit{b}}}}

\def\aa{{\textbf{\textit{a}}}}
\def\ii{{\textbf{\textit{i}}}}

\newtheorem{theo}{Theorem}[section]
\newtheorem{clai}[theo]{Claim}
\newtheorem{nota}[theo]{Notation}

\newtheorem{lemm}[theo]{Lemma}
\newtheorem{fact}[theo]{Fact}

\newtheorem{prop}[theo]{Proposition}
\newtheorem{rema}[theo]{Remark}

\newtheorem{conv}[theo]{Convention}

\newtheorem{defi}[theo]{Definition}

\def\th{{\theta}}

\def\ep{{\textbf{\textit{\footnotesize$\mathbf{\mathcal{E}}$}}}}
\def\scep{{{\small\textbf{\textit{\footnotesize$\sc\mathbf{\mathcal{E}}$}}}}}

\def\th{{\theta}}
\def\equa#1#2{\begin{equation}\label{#1}\mbox{$#2$}\end{equation}}
\def\equan#1#2{\begin{equation*}\mbox{$#2$}\end{equation*}}
\def\re{_{\rm re}}\def\im{_{\rm im}}
\def\TH{\theta}\def\sTH{^2}\def\ssTH{2}
\def\TH{}\def\sTH{}\def\ssTH{}
\def\tau{}
\def\NOUSE#1{}
\def\tildev{v}
\def\tildeX{X}

\def\ZeRo{{1}}
\def\rZeRo{0}
\def\OnE{{2}}
\def\rOnE{1}
\def\dH{{h}}

\def\ln{\log}

\begin{document}
\title{Proof of two-dimensional Jacobian \vspace*{-10pt}conjecture
}
\author{Yucai Su\\
{
D\lowercase{epartment of} M\lowercase{athematics,} J\lowercase{imei} U\lowercase{niversity,} X\lowercase{iamen,} 361021, C\lowercase{hina, and}\\
S{\lowercase{chool of} M\lowercase{athematical} S\lowercase{cience,} T\lowercase{ongji} U\lowercase{niversity,}
S\lowercase{hanghai}, 200092, C\lowercase{hina}\\
\lowercase{\vspace*{-10pt}yucaisu@jmu.edu.cn}}}}

\date{\noindent  \today
\\ \indent
Supported by NSF grant  11971350 of China}

\begin{abstract}Using the local bijectivity of Keller maps, we give a proof of two-dimensional Jacobian conjecture\vspace*{-20pt}.
\end{abstract}
 \maketitle
\tableofcontents

\section{Main theorem}\label{sect1}

Let $(F,G)\in{\mathbb C}[x,y]^2$ be a pair of two polynomials. The Jacobian problem is  to give necessary and sufficient conditions that
${\mathbb C}[F,G]={\mathbb C}[x,y]$.  In geometric terms   if
\equa{Kell---}{\mbox{$\si:\C^2\to\C^2,\ p\mapsto(F(p),G(p)):=(F(a,b),G(a,b))$ \ \ for \ \ $p=(a,b)\in\C^2$,}} then the necessary and sufficient condition is that  $\si$ be invertible.\smallskip

By the chain rule an obvious necessary condition  is that the Jacobian determinant is a nonzero constant
$J(F,G):={\rm det\,}\Big[{}^{\frac{\ptl F}{\ptl x}}_{\frac{\ptl G}{\ptl x}}\,\,^{\frac{\ptl F}{\ptl y}}_{\frac{\ptl G}{\ptl y}}\Big]\in\C_{\ne0}$\vspace*{-3pt}.

If this condition is verified we call $(F,G)$ a {\em Jacobian pair} and the corresponding map $\si$ the {\em Keller map}.
The aim of this paper is to give a proof of the well known two-dimensional Jacobian conjecture (cf.~\cite{{BCW},V1}), formulated by Ott-Heinrich Keller in 1939:
\vskip5pt
\noindent
{\bf Jacobian conjecture } (cf.~\cite{{BCW},V1}): If $(F,G)$  is a {  Jacobian pair} then the {Keller map} $\si$
 is bijective, i.e.,
$F,G$ are generators of $\C[x,y]$.
 \smallskip

 It is a well known and standard fact that, in order to prove this conjecture, it is enough to prove that the map  $\si$
 is injective. This is what we do in
\begin{theo}\label{MAINT}
The {Keller map} $\si$
 is injective. Consequently, the $2$-dimensional Jacobian conjecture holds.\end{theo}
We prove this theorem by contradiction, so we can start with the converse assumption that there exists a Jacobian pair $(F,G)$ such that  the corresponding  Keller map $\si$ is not injective, i.e., there exist $p_{\ZeRo}=(x_{\ZeRo},y_{\ZeRo}),$ $p_{\OnE}=(x_{\OnE},y_{\OnE})\in{\mathbb C}^2$ with
\equa{=simag}{\mbox{$\sigma(p_{\ZeRo})=\sigma(p_{\OnE})$,  \ \ $p_{\ZeRo}\ne p_{\OnE}$.}}
For convenience, we denote
$(p_{\ZeRo},p_{\OnE})=\big((x_{\ZeRo},y_{\ZeRo}),(x_{\OnE},y_{\OnE})\big)\in\C^2\times\C^2\cong\C^4$, and
\begin{eqnarray}
\label{V=0}&\!\!\!\!\!\!\!\!\!\!\!\!\!\!\!\!\!\!&
V=\big\{(p_{\ZeRo},p_{\OnE})=\big((x_{\ZeRo},y_{\ZeRo}),(x_{\OnE},y_{\OnE})\big)\in\C^4\,\big|\,\si(p_{\ZeRo})=\si(p_{\OnE}),
\,p_{\ZeRo}\ne p_{\OnE}\big\}
.\end{eqnarray}
By assumption \eqref{=simag}, $V\ne\emptyset$ (which, to be proven in Lemma \ref{Procesi-1}, is in fact
 a smooth and closed algebraic surface
 in $\C^4$)%
.
We define the {\it height} (or the {\it norm}) of $(p_{\ZeRo},p_{\OnE})=\big((x_{\ZeRo},y_{\ZeRo}),(x_{\OnE},y_{\OnE})\big)\in V$ to be
\equa{dp0p1}{\dH_{p_{\ZeRo},p_{\OnE}}=|x_{\ZeRo}|+|y_{\ZeRo}|+|x_{\OnE}|+|y_{\OnE}|.}
To prove Theorem \ref{MAINT}, we need   two results.
Here is the first one.
\begin{theo}\label{Theo-2}
There exists some automorphism $\phi$ of $\C[x,y]$ such that the Jacobian pair\linebreak $\big(\phi(F),\phi(G)\big)$, which for convenience is still denoted by $(F,G)$, satisfies
the following: for $(p_{\ZeRo},p_{\OnE})\in V$
 when $\dH_{p_{\ZeRo},p_{\OnE}}\to\infty$,
we have $
{
|y_{\ZeRo}|+|y_{\OnE}|=o(\dH_{p_{\ZeRo},p_{\OnE}}).
}
$\end{theo}
 The reason that the result in this theorem is not symmetric on $x,y$ is because we require $(F,G)$ to satisfy
\eqref{wheraraaa}, where the variables $x$ and $y$ are not symmetric.

Once we have established Theorem \ref{Theo-2}, we fix the Jacobian pair $(F,G)$ and the variety $V$ satisfying this theorem. Then we define the projection $\pi_1:V\to\C^2$ by \equa{proj1}{\mbox{$\pi_1:(p_1,p_2)\ \ \mapsto\ \ \pi_1(p_1,p_2):=(x_1,x_2)$ \ \  for $(p_1,p_2)=\big((x_1,y_1),(x_2,y_2)\big)\in V$.}}
Then the second result can be stated as follows.

\begin{theo}\label{Theo-3}The projection $\pi_1$ is proper,   finite and surjective.
\end{theo}

We will give the proofs of the above two theorems in section  \ref{secc2}.
Then finally in section \ref{sect3}, by a rather technical study of the  properties  so far developed of the variety $V$, and with no further use of the fact that $V$ arises from the Jacobian problem, we find a contradiction to Theorem \ref{Theo-3}  thus proving  the Jacobian conjecture.

\NOUSE{%
\begin{theo}\label{Theo-3}Fix $\xi=(\xi_{\ZeRo},\xi_{\OnE})\in\C^2\backslash\pi_1(V)$. We define a continuous function $L_{p_{\ZeRo},p_{\OnE}}$ on $V$  by
\equa{Equ-Theo-3}{\dis
L_{p_{\ZeRo},p_{\OnE}}=d\big(\pi_1(p_1,p_2),\xi\big)^2=|x_{\ZeRo}-\xi_{\ZeRo}|^2+|x_{\OnE}-\xi_{\OnE}|^2\mbox{ \ for \ }(p_{\ZeRo},p_{\OnE})\in V,
}
where $d(\a,\b)$ denotes the distance between $\a$ and $\b$ for $\a,\b\in\C^2$.
Then for any $(p_{\ZeRo},p_{\OnE})\in V$, there exists $(q_{\ZeRo},q_{\OnE})=\big((\dot x_{\ZeRo},\dot y_{\ZeRo}),(\dot x_{\OnE},\dot y_{\OnE})\big)\in V$ such that
\equa{Equ-Theo-3+1}{\dis
L_{q_{\ZeRo},q_{\OnE}}<L_{p_{\ZeRo},p_{\OnE}}.
}
\end{theo}
We will give proofs of Theorem \ref{Theo-1}--\ref{Theo-3} in Section 2--4.\vskip4pt
\noindent
{\it Proof of
Theorem \ref{MAINT}.~}~The second assertion of Theorem \ref{MAINT} follows from the first by \cite{K-m1,K-m2}. To prove the first statement, assume conversely that there exists a Jacobian
 pair 
 satisfying \eqref{=simag}
.
Then we have Theorems \ref{Theo-1}--\ref{Theo-3}.
Fix any $(\bar p_{\ZeRo},\bar p_{\OnE})\in V$ and define $V_{\rOnE}=\{(p_{\ZeRo},p_{\OnE})\in V\,|\,L_{p_{\ZeRo},p_{\OnE}}\le L_{\bar p_{\ZeRo},\bar p_{\OnE}}\}$.
By \eqref{Equ-Theo-3}, using the local bijectivity of Keller maps, exactly as in the proof of Lemma \ref{S4-1}, we see that $V_{\rOnE}$ is a nonempty compact subset of $\C^4$.
Then \eqref{Equ-Theo-3+1} shows that the continuous function $\ell_{p_{\ZeRo},p_{\OnE}}$ does not have the minimal value on $V_1$, which is a contradiction.
This proves Theorem \ref{MAINT}.
\hfill$\Box$
}%

\vskip7pt\noindent{\small {\bf Acknowledgements.}  I would sincerely like to take this opportunity to thank Professor Claudio\linebreak Procesi, who
is a very kind and warmhearted professor, willing to spend more than one year of  his 
valuable time,
 selflessly, seriously, and meticulously helping me thoroughly examine every detail of the paper, discussing with me and giving a huge number of useful comments, suggestions and providing simple proofs of Theorem \ref{Theo-3}%
~and
Proposition \ref{real00-inj+1}
.
Before that, we do not personally know each other. I have great respect for him and admire his incomparably broad knowledge, mind and thoughts. Anybody who compares the present version with any version before December 2022 can see how much hard work has been done by him. Without his help, the paper cannot be the present form.
}

\section{Proofs of Theorems \ref{Theo-2} and \ref{Theo-3}}\label{secc2}

\subsection{Some preparations}

For any ring $R$, we use $R((y^{-1}))$ to denote the ring whose elements have the form $\sum_{i=-\infty}^\infty a_i y^i$ with $a_i\in R$ such that $a_i=0$ when $i\gg1$.
By $R[[y^{-1}]]$ we denote the set of usual power series $\sum_{i=-\infty}^0 a_iy^i$ with $a_i\in R$. Then $R((y^{-1}))$ is the (algebraic) localization $R((y^{-1}))=R[[y^{-1}]][y]$.

\begin{rema}\label{oss0}\rm
\begin{itemize}\item[(i)]An element $\sum_{i=0}^\infty a_iy^{-i}$ is invertible in  $R[[y^{-1}]]$ if and only if $a_0$ is invertible in $R$.
An element in $R((y^{-1}))$ can be uniquely written as $y^k\sum_{i=0}^\infty a_iy^{-i}$ with $k\in \Z$  and $a_0\neq 0$  and it is invertible in  $R((y^{-1}))$ if and only if $a_0$ is invertible in $R$.
\item[(ii)]
Given $f:=\sum_{i=1}^\infty a_iy^{-i}$ we can  define in  $R[[y^{-1}]]$ the {\em substitution homomorphism} $y^{-1}\mapsto  f$.
If  moreover $f$ is invertible  in $R((y^{-1}))$ this extends to a substitution in $R((y^{-1}))$.
\end{itemize}\end{rema}

{Recall  the following standard formula, which holds algebraically in the ring $\C[[z]]$   and converges absolutely when $z\in\C$ with $|z|<1$,
\equa{bimeformo}{(1+z)^{\b}=
 \Big(1+\sum\limits_{j=1}^\infty{\dis\binom{\b}{j}}z^j \Big),\dis\mbox{ \ where \ } \b\in\R,\ \ \binom{ \b}{j}:=\frac{ \b( \b-1)\cdots( \b-( j-1))}{j!},  }
and in general, we denote  the  {\it multi-nomial coefficient}
$
{
\binom{k}{\l_{\rOnE},\l_2,...,\l_i}=\frac{k(k-1)\cdots(k-(\l_{\rOnE}+\l_2+\cdots+\l_i-1))}{\l_{\rOnE}!\l_2!\cdots\l_i!}.}
$

\begin{rema}\label{oss}\rm
\begin{itemize}\item[(i)]
If $A$ is a commutative ring,  two formal series  $f(z),g(z)\in A[[z]]$  can be composed as if they were functions by
\begin{equation}\label{comp}
f\circ g(z)= \mbox{$\sum\limits_{j=1}^\infty$} a_j\Big(\mbox{$\sum\limits_{i=1}^\infty$} b_iz^i\Big)^j\mbox{ \ if $f(z)=\sum\limits_{i=1}^\infty a_iz^i,\ g(z)=\sum\limits_{i=1}^\infty b_iz^i.$}
\end{equation} Moreover $f(z)$ is invertible with respect to $\circ$  if and only if $a_1$ is invertible in $A$.  We  will denote by  $f^{\circ -1}$ its inverse in order to avoid confusions.
\item[(ii)]
In particular if $f$ is invertible we can find $h(z)$ so that  $g(z)=f\circ  h(z)$  [or $g(z)=h\circ  f(z)\ssc\,$].
\item[(iii)]
Finally the usual {\em chain rule}  holds.
\end{itemize}\end{rema}

We need some conventions and notations, which, for easy reference, are listed as follows.
\begin{conv}\label{conv1}\rm\begin{enumerate}\item[(1)]For $a\in\C$, we write  $a=a\re+a\im\ii$ for some
$a\re,a\im\in\R$, where $\ii=\sqrt{-1}$.
If $a^b$ appears somewhere, then we always assume $b\in\R$, and in case $a\ne0
$, we
interpret $a^b$ as the unique complex number $r^be^{b\th\ii}$ by writing
$a=re^{\th\ii}$ for some $r\in\R_{>0}$, $0\le\th<2\pi$. 
\item[\rm(2)]
Let $P=\sum_{i\in\Z_{\ge0}} p_iy^{\a-i}=y^\a\big(1+\sum_{i=1}^\infty p_iy^{-i}\big)\in\PP
$ with $\a\in\Z,\,
p_i\in\C[x]$ and $p_{\rZeRo}=1$.
\begin{itemize}\item[(i)]
Let $\b\in\Q$  with $\a\b\in\Z$. Then we can uniquely define $P^\b$ to be an element in $\PP$ as follows,
\begin{eqnarray}\label{P-as}
&\!\!\!\!\!\!\!\!\!\!\!\!\!\!\!\!\!\!\!\!\!\!&
P^{\b}=
y^{\a\b}\Big(1+\mbox{$\sum\limits_{j=1}^\infty$}{\dis\binom{\b}{j}}\Big(\mbox{$\sum\limits_{i=1}^\infty$}
p_iy^{-i}\Big)^j \Big) \in\PP,
\nonumber\\
&\!\!\!\!\!\!\!\!\!\!\!\!\!\!\!\!\!\!\!\!\!\!&
\Big(\mbox{$\sum\limits_{i=1}^\infty$}
p_iy^{-i}\Big)^j   =\mbox{$\sum\limits_{\sum  \l_i=j}$} {\dis\binom{j}{\l_{\rOnE},\l_2,... }}
  p_1^{\l_{\rOnE} }p_2^{\l_2} p_3^{\l_3}\cdots  y^{-
   \sum i\l_i} .
\end{eqnarray}
\item[(ii)]
For $Q_{\rOnE},Q_2\in\PP$, we use $P(Q_{\rOnE},Q_2)$ and $P|_{(x,y)=(Q_{\rOnE},Q_2)}$ to denote the following element [as long as it is algebraically
a well-defined element in $\PP {\sc\,}$],
\equa{P(qq)}{P(Q_{\rOnE},Q_2)=P|_{(x,y)=(Q_{\rOnE},Q_2)}=\sum\limits_{i}p_{i}(Q_{\rOnE})Q_2^{\a-i}.}
\item[(iii)]
If $Q_{\rOnE},Q_2\in\C$, 
 we also use \eqref{P(qq)} to denote a well-defined complex number as long as
the series \eqref{P(qq)} \CA$\!{\sc\!}.$
\item[(iv)]Assume $\a\ne0$.
For any $Q=\sum_{i\in
\Z_{\ge0}}q_iy^{\a_1-i}\in\PP$ with $\a_1\in\Z$, $q_i\in\C[x]$, by comparing coefficients of $y^{\a_1-i}$ in \eqref{usia} for $i\ge0$, there exists uniquely $b_i\in\C[x]$ such that
we can algebraically write
\equa{usia}{
Q=\mbox{$\sum\limits_{i=0}^\infty$}b_iP^{\frac{\a_1-i}{\a}}.}
We call $b_i$ the {\it coefficient of $P^{\frac{\a_1-i}{\a}}$ in $Q$,} and denote by $\Coeff(Q,P^{\frac{\a_1-i}{\a}})$. We also use the following notation $\Coeff(Q,x^iy^j)$:
\equa{nene34bb4b44b}{\phantom{mmmMMM}\Coeff(Q,x^iy^j)=q_{ij}\mbox{ \ \ if \ $Q$ can be written as $Q=\sum\limits_{i,j}q_{ij}x^iy^j$ with $q_{ij}\in\C$.}}
\end{itemize}
\item[(3)]
Let $\ep\to0$ be a variable. Sometimes we need to consider elements in $\PP$ which may depend on $\ep$.
For any $P\in\PP$ (or especially in $\C$)  which  depends on $\ep$, if
$P(a,b)$ \CA and $|\ep^{-i}P(a,b)|<\SS$ for some $i\in\Q_{\ge0}$ and some fixed
$\SS\in\R_{>0}$, where  $(a,b)$ is in some required region
,
then we use $O(\ep)^i$ to denote $P$:
\equa{P123O-ep}{P=O(\ep)^i.}
If $a,b$ are some variables depending on another variable $c\to c_{\rZeRo}$ (for some $c_{\rZeRo}\in\C\cup\{\infty\}$) such that
$\lim_{c\to c_{\rZeRo}}\frac{a}{b}=0$,
then we also denote \equa{a-bo}{\dis a=o(b).}
\end{enumerate}
\end{conv}

Let $P=\sum_{j}p_{j}y^{j}\in\PP,\,p_j\in\C[x] $ and $(x_{\rZeRo},y_{\rZeRo})\in\C^2$.
If  
 $
z_{\rZeRo}=\sum
_{j} |p_{j}(x_{\rZeRo})y_{\rZeRo}^{j}|$ 
converges (in particular this requires that $y_0\ne0$ if $p_j(x_0)\ne0$ for some $j<0$),
then $z_{\rZeRo}$ is 
denoted by $A_{(x_{\rZeRo},y_{\rZeRo})}(P)$ [or  by $A_{(y_{\rZeRo})}(P)$ if $P$ is independent of $x$].
\begin{defi}\label{contro}\rm
\begin{itemize}\item[(1)]
Let $P$ be as above and
$Q=\sum_iq_iy^i\in\C((y^{-1})),\,q_i\in\R_{\ge0}$, $x_{\rZeRo}\in\C$. If 
{\mbox{$|p_i(x_{\rZeRo})|\le q_i$
}}for all possible
$i$,
then we say $Q$ is a {\it controlling function} for $P$ on
$y$ at point $x_{\rZeRo}$, and denote
 \equa{MDEooooo}
 {\mbox{$P\,\DS^{x_{\rZeRo}}_y\, Q$  \ \ \ or \ \ \  $Q\,\rDS^{x_{\rZeRo}}_y\, P$
 ,}}
or  $P\,\DS_y\, Q$ 
when there is no confusion.
In particular if $P,Q$ are independent of $y$ then we write $P\DS^{x_{\rZeRo}} Q$ 
(thus $a\DS b$ for $a,b\in\C$ simply means that $|a|\le b$ with $b\ge0$).
\item[(2)]An element in $\C((y^{-1}))$ with non-negative coefficients (such as $Q$ above)
is called a {\it controlling function} on $y$.
\item[(3)]
If $Q=q_{\rZeRo}y^\a+\sum_{j>0} q_{j}y^{\a-j}\in\C((y^{-1}))$ is a controlling function on
 $y$ with 
$q_{\rZeRo}>0$,  then we always use the same symbol with
subscripts ``\,$_\igo$\,'' and
``\,$_\inv$\,'' to denote the  elements
\begin{eqnarray}
\label{P---}
\!\!\!\!\!\!\!\!\!\!\!\!\!\!\!\!\!\!\!\!\!\!\!\!\!\!&&
Q_\igo=q_{\rZeRo}^{-1}\mbox{$\sum\limits_{j>0}$} q_{j}y^{-j},
\nonumber\\
\!\!\!\!\!\!\!\!\!\!\!\!\!\!\!\!\!\!\!\!\!\!\!\!\!\!&&
Q_\inv:=q_{\rZeRo}y^\a\Big(1-q_{\rZeRo}^{-1}\mbox{$\sum\limits_{j>0} $}q_{j}y^{-j}\Big)=q_{\rZeRo}y^\a(1-Q_\igo)=2q_{\rZeRo}y^\a-Q.\!\!\!\!\!\!
\end{eqnarray}
We call $Q_\igo$ the {\it ignored  part} of $Q$, and $Q_\inv$ the {\it negative correspondence of} $Q$
[in sense of \eqref{MAMS1}
, where $a,-k$ are nonpositive].
\end{itemize}\end{defi}

 Notice that when $Q$ in \eqref{P-Q===} is a controlling function then for $k\in \R_{>0} $  we have that both $\big(1\!-\!q_{\rZeRo}^{-1}\sum_{j>0} q_{j}y^{-j}\big)^{-k}$  and  $\frac{(q_{\rZeRo}y^{\a})^k}{1-k
Q_\igo}$  are controlling functions.

\begin{lemm}\label{ds-lemm}\begin{itemize}\item[\rm(1)]
If \equa{P-Q===}{\mbox{$P=p_{\rZeRo}y^\a+\sum\limits_{j>0} p_{j}y^{\a-j}\in\PP,\ \ \
Q=q_{\rZeRo}y^\a+\sum\limits_{j>0} q_{j}y^{\a-j}\in\C((y^{-1}
)),$}} with $P\DS^{x_{\rZeRo}}_y Q$, $x_{\rZeRo}\in\C$ and $|p_{\rZeRo}(x_{\rZeRo})|=q_{\rZeRo}\in\R_{>0}$,
then for $a,b,k\in\Q$ with $a\a,b\a,k\a\in\Z$,
\begin{eqnarray}\label{MAMS1}&\!\!\!\!\!\!\!\!\!&{\rm(a)\ }\dis\frac{\ptl P}{\ptl y}\ \DS^{x_{\rZeRo}}_y\ \,  \pm{\sc\,}
\frac{d Q}{d y},\ \ \ \ \ \ \
\nonumber\\[0pt]
\!\!\!\!\!\!\!\!&\!\!\!\!\!\!\!\!\!&
{\rm(b)\ }P^a\ \DS^{x_{\rZeRo}}_y\ \, Q_\inv^a\ \DS_y\ \, (q_{\rZeRo}y^{\a})^{-b}Q_\inv^{a+b}\mbox{ for $a,b\in\Q_-$,}
\nonumber\\[0pt]
\!\!\!\!\!\!\!\!&\!\!\!\!\!\!\!\!\!&
{\rm(c)\ }Q^k\ \DS_y\ \, (q_{\rZeRo}y^{\a})^{2k}Q_\inv^{-k}\ \DS_y\ \,\left\{\begin{array}{ll}\dis
\frac{(q_{\rZeRo}y^{\a})^k}{1-k
Q_\igo}&\mbox{if }k\in\Z_{\ge1},\\[6pt]\dis
(q_{\rZeRo}y^{\a})^k\Big(1+\dis \frac{kQ_\igo}{1-Q_\igo}\Big)&\mbox{if }k\in\Q_{\ge0}\mbox{ with }k<1.
\end{array}\right.
\end{eqnarray}
where  $\eqref{MAMS1}\,(a)$ holds under the condition: either both $P,$ $Q$ are polynomials of
$y$ $($in this case the sign is ``\,$+$\,''${\ssc\,})$, or else both are power series of $y^{-1}$ $($in this case the sign is ``\,$-$\,''${\ssc\,})$.
\item[\rm(2)]If $x_{\rZeRo},y_{\rZeRo}\in\C$ 
and $P_{\rOnE}\DS^{x_{\rZeRo}}_y\, Q_{\rOnE},\ P_2\DS^{x_{\rZeRo}}_y\, Q_2$, then
\equa{ABV1}{A_{(x_{\rZeRo},y_{\rZeRo})}(P_{\rOnE}P_2)\le A_{(y_{\rZeRo})}(Q_{\rOnE})A_{(y_{\rZeRo})}(Q_2)=Q_{\rOnE}(|y_{\rZeRo}|)Q_2(|y_{\rZeRo}|).}
\end{itemize}
\end{lemm}\noindent{\it Proof.} One can see that (2)
and \eqref{MAMS1}\,(a) are obvious, and \eqref{MAMS1}\,(b),\,(c) 
are obtained by
 noting that for $a,b\in\Q_-$ and $i\in\Z_{>0}$, one has
 \begin{eqnarray}
\label{amen3n33nn3}
&\!\!\!\!\!\!\!\!\!\!\!\!\!\!\!\!\!\!\!\!\!\!\!\!\!\!\!&
\mbox{$ \dis(-1)^i\binom{a}{i}=
 \Big|\binom{a}{i}\Big|{\ssc}\le{\ssc}\Big|\binom{a+b}{i}\Big|{\ssc}={\ssc}(-1)^i\binom{a+b}{i}$,}
\nonumber\\
&\!\!\!\!\!\!\!\!\!\!\!\!\!\!\!\!\!\!\!\!\!\!\!\!\!\!\!&
\mbox{$\dis \binom{k}{i}{\ssc}\le{\ssc}\Big|\binom{-k}{i}\Big|{\ssc}\le{\ssc} \left\{\begin{array}{ll}k^i&\mbox{if  }k\in\Z_{\ge1},\\[4pt]
k&\mbox{if  }0<k\in\Q_{<1}.\!\!\!\!\end{array}\right.$}
\end{eqnarray}
This proves the lemma.\hfill$\Box$\vskip7pt

\def\bB{{\textbf{\textit{b}}}}\def\zz{y}
Throughout the rest of this section, for convenience we regard $y^{-1}$ as the variable $z$ whenever necessary.

Take, where $\tilde f_i\in\C[x]$ with $\tilde f_{\rOnE}\in\C_{\ne0},$
\equa{Faa}{\dis\tilde F=\tilde f_{\rOnE}\zz^{-1} +\mbox{$\sum\limits_{i=2}^\infty$}\tilde  f_i \zz ^{-i} \in\C[x][[\zz^{-1} ]].
}%
\def\tf{f}Then $\tilde F$ is invertible under the composition ``\,$\circ$\,'' defined in
\eqref{comp} (by regarding $y^{-1}$ as $z$). We refer $\tilde F^{\circ-1}$ to as the {\it formal inverse function} of $\tilde F$. Then $\tilde F^{\circ-1}(y^{-1})\in\C[x][[y^{-1}]]$, and thus we can write
\equa{Inv-of-T-F}{\dis
\tilde F^{\circ-1}(y^{-1})=\bB_1 y^{-1}+\mbox{$\sum\limits_{i=2}^\infty$} \bB_i y^{-i}\ \ \mbox{ for some }\ \ \bB_i\in\C[x].
}
By definition, we have (noting that $y^{-1}$ is the identity element under the composition ``\,$\circ$\,''),
\NOUSE{
Regarding $\tilde F$ as a formal function of $\zz^{-1} $ (with parameter $x$ being regarded as fixed), we have the {\it formal inverse function} denoted by
$\zz _{{\tilde F}}\in\C[x,\tilde f_{\rOnE}^{-1}][[\tilde F]]\subset\C[x][[\tilde F]]$ such that, as in \eqref{usia}, we can algebraically write [note that we denote the right-hand side below as $\zz _{{\tilde F}}(\tilde F)$ in order to emphasis that the (formal) function $\zz _{{\tilde F}}$ is defined in such a way that if $\tilde F$ is regarded as a variable then $\zz _{{\tilde F}}(\tilde F)$ is defined to be the right-hand side],
}%
\begin{eqnarray}\!\!\!\!\!\!\!\!
\label{i-Faa}
\zz^{-1} &\!\!=\!\!&
\tilde F^{\circ-1}\circ\tilde F\stackrel{{}^{\sc\rm\eqref{comp}}}{=}\tilde F^{\circ-1}(\tilde F)\stackrel{{}^{\sc\rm\eqref{Inv-of-T-F}}}{=}
\bB_{\rOnE}\tilde F+\mbox{$\sum\limits_{i=2}^\infty$} \bB_i {\tilde F}^i.
\end{eqnarray}
By the notation in Convention \ref{conv1}\,(2)\,(iv), we see from \eqref{i-Faa} that
\equa{bbi==}{\mbox{$\bB_i=\Coeff(y^{-1},\tilde F^i)\in\C[x]$ \ for \ $i\ge1$,}} which can be precisely determined by using \eqref{Faa} to substitute $\tilde F$ in \eqref{i-Faa} and comparing the coefficients of $y^{-i}$ 
in both sides of \eqref{i-Faa} such that
(we do not need to use the following explicit expression of $\bB_i$, but we only want to present that $\bB_i$'s exist),
 \begin{eqnarray}
\label{i-Faa+}
\!\!\!\!\!\!\!\!\!\!\!\!\!\!\!\!\!\!\!\!\!
\bB_1&\!\!=\!\!&
\tilde f_1^{-1}\in\C_{\ne0},
\\\nonumber
\!\!\!\!\!\!\!\!\!\!\!\!\!\!\!\!\!\!\!\!\!
\bB_i&\!\!=\!\!&
-\mbox{$\sum\limits_{j=1}^{i-1}$}\bB_j\tilde f_{\rOnE}^{j-i}
\mbox{$\sum\limits_{\ell=0}^{j}$}
{\dis\binom{j}{\ell}}\mbox{$\sum\limits_{\stackrel{\sc n\in\Z_{\ge0},\,\l_{\rOnE},\l_2,...,\l_n\ge0}{{}^{\ }_{\sc\l_{\rOnE}+2\l_2+\cdots+n\l_n={\ssc\,}i-j}}}$}{\dis\binom{\ell}{\l_{\rOnE},\l_2,...,\l_n}}
\tilde f_{\rOnE}^{-\l_{\rOnE}-\l_2-\cdots-\l_n}\tilde \tf_2^{\l_2}\tilde \tf_3^{\l_3}\cdots \tilde \tf_{n}^{\l_n},\ i\ge2.\!\!\!\!\!\!\!\!\!\!\!\!\!\!\!\!\!\!\!\!
\end{eqnarray}
\begin{lemm}\label{lemm2222}Let $\hat a_i\in\R_{\ge0}$ with $\hat a_{\rOnE}>0$, and let
\equa{AJHAHA}{\FF =\hat a_{\rOnE}\zz^{-1} +\sum\limits_{i=2}^\infty \hat a_i\zz ^{-i}\in\C[[\zz^{-1} ]]
\mbox{ \ \ and \ $\FF _\inv \stackrel{{}^{\sc\rm\eqref{P---}}}{=}\hat a_{\rOnE}\zz^{-1} -\sum\limits_{i=2}^\infty \hat a_i\zz ^{-i}\in\C[[y^{-1}]]$},} be a controlling function on $\zz $ and its negative correspondence.
Write the formal inverse function $\FF _\inv^{\circ-1}$ of  $\FF _\inv$, as
\NOUSE{
Let $[$note that we denote the right-hand side below as $\zz{\ssc\,}^\inv (\FF _\inv )$ to emphasis that the $($formal$)$ function $\zz{\ssc\,}^\inv$ is defined, which is different from $\zz _{{\tilde F}}$ defined
in \eqref{i-Faa}${\ssc\,}]$,
}
\equa{ha--f-}{\FF _\inv^{\circ-1}(y^{-1}) =
\hat \bB_{\rOnE}y^{-1} +\sum\limits_{i=2}^\infty\hat \bB_iy ^{-i}\ \ \mbox{ for some }\ \ \hat \bB_i\in\C 
.
}
Then $\hat \bB_1=\hat a_1^{-1}$ and
\begin{itemize}
\item[\rm(1)] $\FF _\inv^{\circ-1}(y^{-1})$ is a controlling function on $y^{-1}$ $($this means that the formal inverse function of the
 negative correspondence of a controlling function is a controlling function$)$,
 i.e., for $i\ge2$,
\equa{bBi>0}{\mbox{$\dis \hat \bB_i\ge0$.}}
\item[\rm(2)]If $\tilde F\DS^{x_{\rZeRo}}_\zz \FF $ with $\tilde F$ as in \eqref{Faa}
and $|\tilde f_{\rOnE}
|=\hat a_{\rOnE}$,
then by regarding $\tilde F$ as a variable, we have the following $($this means that the formal inverse function of a function can be controlled by the formal inverse function of the negative correspondence of a controlling function that controls the said function$)$,
\begin{eqnarray}
\!\!\!\!\!\!\!\!\!\!\!\!\!\!\!\!\!\!\!\!&&
\label{HAHAHAHHJ}
\zz^{-1}=\tilde F^{\circ-1}(\tilde F)\DS_{{\tilde F}}^{x_{\rZeRo}}\, \FF _\inv^{\circ-1}(\tilde F),
\mbox{ \ i.e.,}
\nonumber\\
\!\!\!\!\!\!\!\!\!\!\!\!\!\!\!\!\!\!\!\!&&
\bB_i\stackrel{{}^{\sc\rm\eqref{bbi==}}}{=}\Coeff(y^{-1},\tilde F^i)\  \DS^{x_{\rZeRo}}\ \hat \bB_i=\Coeff(y^{-1},\FF^i_\inv)
\mbox{ \ for \ }i\ge2,
\end{eqnarray}
where 
$\bB_i\DS^{x_{\rZeRo}}\,\hat \bB_i$ means that
$|\bB_i(x_{\rZeRo})|\le\hat \bB_i$.
In particular
\equa{y-ds-}{\zz^{-1} \DS_\zz \, \FF _\inv^{\circ-1} (\FF ),}
where the right side of ``\,$\DS_\zz $\,'' is regarded as a function of $\zz $ by using \eqref{AJHAHA} to substitute $\FF $.
\end{itemize}
\end{lemm}\noindent{\it{Proof.~}}~Note that (1) is a special case of (2) by taking, in \eqref{Faa}, $\tilde f_1=\hat a_1$ and $\tilde f_i=0$ for $i\ge2$, i.e., $\tilde F=\hat a_{\rOnE}\zz^{-1} $ [then by definition $\tilde F^{\circ-1}(y^{-1})$ is
$\hat a_1^{-1}y^{-1}$, i.e., $\bB_1=\hat a_1^{-1}$ and $\bB_i=0$ if $i\ge2$, and so \eqref{bBi>0} follows from \eqref{HAHAHAHHJ}$\ssc\,$].

Thus we prove (2).
Note from \eqref{Inv-of-T-F} that $\bB_i$ is the coefficient of $y^{-i}$ in the formal inverse function of $\tilde F$, or equivalently, the coefficient of $\tilde F^i$ in $y^{-1}$ when we regard $y^{-1}$ as a function of $\tilde F$ by \eqref{i-Faa}.
Thus $i!\bB_i$ is the constant term of the $i$-th partial  derivative $ \frac{\partial^i \zz^{-1} }
{{^{\ssc\,}}\partial {\tilde F}^i }$.
Therefore to prove \eqref{HAHAHAHHJ}, first we want to prove, for $i\ge1$,
\equa{to0000}{\dis\frac{\partial^i \zz^{-1} }
{{^{\ssc\,}}\partial {\tilde F}^i }\DS^{x_{\rZeRo}}_\zz \, \frac{d^i \zz^{-1} }{d \FF{} _\inv ^i }{\ssc\,},}
where the left-hand side is understood as that we first use \eqref{i-Faa} to
regard $\zz^{-1} $ as a function of $\tilde F$ (with parameter $x$) and
apply $\frac{\partial^i }{{^{\ssc\,}}\partial {\tilde F}^i }$ to it, then regard the result as a function
 of $\zz^{-1} $ by using \eqref{Faa} to substitute $\tilde F$ (and the like for the right-hand side, which does not contain the parameter $x$).
By \eqref{MAMS1}\,(a), we have
$\frac{\ptl\tilde  F}{\ptl \zz^{-1} }\DS^{x_{\rZeRo}}_\zz \,\frac{d \FF }{d \zz^{-1} }$ (here we regard $y^{-1}$ as a variable), and thus by
\eqref{MAMS1}\,(b) and definition \eqref{P---},
\equa{FiQQQ}{\dis \Big(\frac{\ptl \tilde F}{\ptl \zz^{-1} }\Big)^{-1} \stackrel{{}^{\sc\rm{\eqref{MAMS1}\,(b)}}}{\DS^{x_{\rZeRo}}_\zz}
\,\Big(\Big(\frac{d \FF }{d \zz^{-1} }\Big)_\inv\Big) ^{-1}
 \stackrel{{}^{\sc\rm{\eqref{P---}}}}{=} \Big(\frac{d\FF_\inv }{d \zz^{-1} }\Big)^{-1},} i.e., $\frac{\partial \zz^{-1} }{\partial\tilde  F}
\DS^{x_{\rZeRo}}_\zz \,
\frac{d \zz^{-1} }{d \FF_\inv }$ and  \eqref{to0000} holds for $i=1$. Inductively, by Lemma \ref{ds-lemm},
\begin{eqnarray}
&\!\!\!\!\!\!\!\!\!\!\!\!\!\!\!\!\!\!&\frac{\partial^i \zz^{-1} }{\partial \tilde F^i }
=\frac{\ptl}{\ptl \tilde F}\Big(\frac{\partial^{i-1} \zz^{-1} }
{\partial\tilde  F^{i-1} }\Big)=\frac{\ptl}{\ptl \zz^{-1} }
\Big(\frac{\partial^{i-1} \zz^{-1} }{\partial \tilde F^{i-1} }\Big)
\Big({\frac{\ptl\tilde  F}{\ptl \zz^{-1} }}\Big)^{-1}\nonumber\\[-0pt]
&\!\!\!\!\!\!\!\!\!\!\!\!\!\!\!\!\!\!&
\phantom{\frac{\partial^i \zz^{-1} }{\partial \tilde F^i }=}\!\!\!\!\!\!\!\!\!\!\!\!\!\!\!\!\!\!\!\!\!\!\!\!\!\!\!\!\!\!\!\!\!\!\!\!
\stackrel{{}^{\sc\rm inductive\ assumption\ and\ \eqref{FiQQQ}}}{\DS^{x_{\rZeRo}}_\zz}
\frac{d}{d \zz^{-1} }\Big(\frac{d^{i-1} \zz^{-1} }{d \FF{} _\inv ^{i-1} }\Big)
\Big({\frac{d \FF_\inv }{d \zz^{-1} }}\Big)^{-1}=\frac{d^i \zz^{-1} }{d \FF{}_\inv ^i }\,.
\end{eqnarray}
This proves \eqref{to0000}. Using \eqref{to0000} and noting from \eqref{i-Faa},\,\eqref{ha--f-},
we have (noting that $\frac{d^i \zz^{-1} }{d \FF{} _\inv ^i }\Big|_{\zz^{-1} =0}$ should be understood as
the constant term in $\frac{d^i \zz^{-1} }{d \FF{} _\inv ^i }$ when it is written as a series of $y^{-1}$),
\equa{MAmama}{\mbox{$\!\!\!\!\!\!\!\!\!\!\!$
$\begin{array}{llll}\dis\bB_{i
}\!\!\!&\dis
 \stackrel{{}^{\sc\rm\eqref{i-Faa}}}{=} \frac{1}{i!}\frac{\partial^i \zz^{-1} }{\partial \tilde F^i }\Big|_{\tilde F=0}
=\frac{1}{i!}\frac{\partial^i \zz^{-1} }{\partial\tilde  F^i }\Big|_{\zz^{-1} =0}
 \stackrel{{}^{\sc\rm\eqref{to0000}}}{\DS^{x_{\rZeRo}}} \,
\frac{1}{i!}\frac{d^i \zz^{-1} }{d \FF{} _\inv ^i }\Big|_{\zz^{-1} =0}=
\frac{1}{i!}\frac{d^i \zz^{-1} }{d \FF{} _\inv ^i }\Big|_{\FF _\inv =0} \stackrel{{}^{\sc\rm\eqref{ha--f-}}}{=} \hat \bB_{i
}.\end{array}\!\!\!\!$
}\mbox{$\!\!\!\!\!\!\!$}}
This proves \eqref{HAHAHAHHJ}.

Note that for any controlling functions $Q,P_2$, and any $P_1\in\C[x]((y^{-1}))$ with $P_1\DS^{x_0}_y P_2$, we always have $Q(P_1)\DS^{x_0}_y Q(P_2)$. This
with the facts that  $\tilde F\DS^{x_{\rZeRo}}_\zz \,\FF $ and that $ \FF _\inv^{\circ-1} $ is a controlling function implies
\equa{ff-thattt-that}{\mbox{ $ \FF _\inv^{\circ-1} (\tilde F)\DS^{x_{\rZeRo}}_\zz \, \FF _\inv^{\circ-1}  (\FF )$.}}
This together with \eqref{HAHAHAHHJ}
proves
\equa{y-inver-small}{
y^{-1}\stackrel{{}^{\sc\rm\eqref{HAHAHAHHJ}}}{\DS_{{\tilde F}}^{x_{\rZeRo}}}\, \FF _\inv^{\circ-1}(\tilde F)
\stackrel{{}^{\sc\rm\eqref{ff-thattt-that}}}{\DS^{x_{\rZeRo}}_\zz} \, \FF _\inv^{\circ-1}  (\FF )
,
}
i.e., we have
 \eqref{y-ds-}.\hfill$\Box$

\subsection{Proof of Theorem \ref{Theo-2}}\label{section2.2}

First we need to reformulate $(F,G$) (by applying some automorphism of $\C[x,y]$).
Fix 
a sufficiently large $\ell\in\Z_{>0}$. 
Applying the following
variable change,
\equa{newvara}{(x,y)\mapsto(y,\,y^\ell+x),}
and rescaling $F,G$, we can assume that $F,G$ have the following forms,
for some $m,n\in\Z_{>0}$, $f_{j,k},g_{j,k}\in\C$,
\begin{eqnarray}
&\!\!\!\!\!\!\!\!\!\!\!\!\!\!\!\!\!\!\!\!\!&
\label{wheraraaa}
{\rm(i)\ }F{\ssc}={\ssc}y^m{\ssc}+{\ssc}F_{\rOnE},\ \ \ F_{\rOnE}{\ssc}={\ssc}\mbox{$\sum\limits_{j=0}^{m-1}\Big(\sum\limits_{k=0}^{m-1-j}f_{j,k}x^k\Big)y^j,$}
\nonumber\\
&\!\!\!\!\!\!\!\!\!\!\!\!\!\!\!\!\!\!\!\!\!&
{\rm(ii)\ }G{\ssc}={\ssc}y^n{\ssc}+{\ssc}G_{\rOnE},\ \ \  G_{\rOnE}{\ssc}={\ssc}\mbox{$\sum\limits_{j=0}^{n-1}\Big(\sum\limits_{k=0}^{n-1-j}g_{j,k} x^k\Big)y^j.$}\!\!\!\!\!\!\!\!
\end{eqnarray}
Note that $\deg\,F_1\le m-1,\,\deg\,G_1\le n-1$ (where $\deg\,F_1$ denotes the total degree of $F_1$).
For simplicity, we can also assume $2\le n$ and $n{\ssc\,}|{\ssc\,}m$ (i.e., $\frac{m}{n}\in\Z_{>0}$)
 by replacing $(F,G)$ by \linebreak $\big(F+(G+F^k)^k,G+F^k\big)$ for some $k\in\Z_{>0}$
[the reason we assume $n|m$ will be clear in \eqref{denote-t}, \eqref{S-1-D} and \eqref{TDSP}, after all it is a reasonable choice].

Thus the new pair $(F,G)$ is in fact obtained from the original one by applying some automorphism of $\C[x,y]$ and a change of generators in $\C[F,G]$.

In the following we  consider $F,G$ as elements   in the ring $\C[x]((y^{-1}))$.

By \eqref{wheraraaa}, we can rewrite $F,G$ as, where $f_i(x)=\sum_{k=0}^{i-1}f_{m-i,k}x^k$, $g_i(x)=\sum_{k=0}^{i-1}g_{n-i,k}x^k\in\C[x]$,
and where in general for any $A\in\C[x]((y^{-1}))$ we use $\deg_x A$, called the {\it $x$-degree} of $A$, to denote  the degree of $A$ with respect to variable $x$,
\begin{eqnarray}
\label{wheraraaa+1}
&\!\!\!\!\!\!\!\!\!\!\!\!\!\!\!\!\!\!\!\!\!\!&
{\rm(i)\ }F=y^m(1+y^{-m}F_1) =y^m\Big(1+\mbox{$\sum\limits_{j=0}^{m-1}$}\Big(\mbox{$\sum\limits_{k=0}^{m-1-j}$}f_{j,k}x^k\Big)y^{j-m}\Big)
=y^m\Big(1+\mbox{$\sum\limits_{i=1}^{m }$}f_i(x) y^{-i}\Big), \!\!\!\!\!\!
\nonumber\\[-4pt]
&\!\!\!\!\!\!\!\!\!\!\!\!\!\!\!\!\!\!\!\!\!\!&
{\rm(ii)\ }G=y^n\Big(1+\mbox{$\sum\limits_{j=0}^{n-1}$}\Big(\mbox{$\sum\limits_{k=0}^{n-1-j}$}g_{j,k} x^k\Big)y^{j-n}\Big)=y^n\Big(1+\mbox{$\sum\limits_{i=1}^{n }$}g_i(x) y^{-i}\Big),
\nonumber\\\nonumber
&\!\!\!\!\!\!\!\!\!\!\!\!\!\!\!\!\!\!\!\!\!\!&
{\rm(iii)\ }\deg_xf_i\le i-1\mbox{ if }i\le m\mbox{ and }f_i=0\mbox{ if }i>m, \ \
\\[2pt]
&\!\!\!\!\!\!\!\!\!\!\!\!\!\!\!\!\!\!\!\!\!\!&
{\rm(iv)\ }\deg_xg_i\le i-1\mbox{ if }i\le n\mbox{ and }g_i=0\mbox{ if }i>n.
\end{eqnarray}
%

At this point we want to  define a choice  of an $m$-th
root for $F$ denoted $F^{\frac1m}$ by using formulas \eqref{bimeformo} and
\eqref{P-as},
\equa{root}{
F^{\frac1m}\ \stackrel{{}^{\sc\rm\eqref{P-as}}}{:=}\ y\Big(1+\sum\limits_{i=1}^{m }f_i(x) y^{-i}\Big)^{\frac1m}
\stackrel{{}^{\sc\rm\eqref{bimeformo}}}{=}y\Big(1+\sum\limits_{j=1}^\infty{\dis\binom{\frac1m}{j}}\Big(\sum\limits_{i=1}^mf_i(x)y^{-i}\Big)^j\Big).
}  This is not just an algebraic formula  but setting
\equa{Noum}{||F_{\rOnE}||:=\sum\limits_{j=0}^{m-1}\Big(\sum\limits_{k=0}^{m-1-j}|f_{j,k}x^k|\Big)|y|^j,}  it converges absolutely when
$||F_{\rOnE}|| <|y|^m.$

Set $f_0=g_0=1$, and
denote the set $A=\frac{n-\Z_{\ge0}}{m}$.
By \eqref{wheraraaa+1},  as in \eqref{usia} and \eqref{i-Faa},
we can algebraically write,
\begin{eqnarray}\label{ome-2}
&\!\!\!\!\!\!\!\!\!\!\!\!\!\!\!\!&
G=\mbox{$\sum\limits_{\a\in A}c_\a F^\a=\sum\limits_{i=0}^\infty c_{\frac{n-i}{m}}F^{\frac{n-i}{m}}
\mbox{ \  for some $c_\a\in\C[x]$},$}
\end{eqnarray}
where as in \eqref{i-Faa+}, by comparing the coefficients of $y^{n-i}$, we can inductively determine $c_{\frac{n-i}{m}}\in\C[x]$ for $i\ge0$ as follows (again we do not need the explicit expression of $c_{\frac{n-i}{m}}$):
\begin{eqnarray*}
\!\!\!\!
c_{\frac{n-i}{m}}&\!\!=\!\!&
g_i-\mbox{$\sum\limits_{j=1}^{i-1}$}c_{\frac{n-j}{m}} f_{\rOnE}^{j-i}
\mbox{$\sum\limits_{\ell=0}^{j}$}
{\dis\binom{j}{\ell}}\mbox{$\sum\limits_{\stackrel{\sc n\in\Z_{\ge0},\,\l_{\rOnE},\l_2,...,\l_n\ge0}{{}^{\ }_{\sc\l_{\rOnE}+2\l_2+\cdots+n\l_n={\ssc\,}i-j}}}$}{\dis\binom{\ell}{\l_{\rOnE},\l_2,...,\l_n}}
 f_{\rOnE}^{-\l_{\rOnE}-\l_2-\cdots-\l_n} \tf_2^{\l_2} \tf_3^{\l_3}\cdots  \tf_{n}^{\l_n}.
\end{eqnarray*}

\begin{rema}\label{oss1}\rm\begin{itemize}\item[(i)]
If $U=1+\mbox{$\sum_{i=1}^\infty$} u_iy^{-i}\in\C[x][[y^{-1}]]$  with $u_i\in\C[x],\, \deg_xu_i\le i$ and $\b\in\R$  then $U^{\b}=1+\mbox{$\sum_{i=1}^\infty$} v_iy^{-i}$ for some $v_i\in\C[x]$
with  $ \deg_xv_i\le i$.
\item[(ii)]
Furthermore, if  for all $i$ we have $ \deg_xu_i< i$  we have $ \deg_xv_i< i$.\end{itemize}\end{rema}

\noindent{\it Proof.}  This follows from a symbolic computation. Think first the $u_i$ as variables of weight $i$ then,
using formula    \eqref{bimeformo}   it is enough to prove that, for all $j\in\Z_{>0}$  we have that in  $\big(\mbox{$\sum_{i=1}^\infty$} u_iy^{-i}\big)^j$ it is a polynomial
in $y^{-1}$  where the coefficient  of $y^{-a}$  is a polynomial in the $u_i$ of weight $a$. This is immediate by writing $u_iy^{-i}=(\lambda^i  u_i)( \lambda y)^{-i}$ with $\lambda$ an auxiliary parameter.

Now by hypothesis the degree of $u_i$ is smaller  than or equal to
its weight and this property is preserved, the same if it is always smaller  than its weight. \hfill$\Box$
\begin{lemm}\label{Sect3-1}
We have, for some $a_1\in\C$,
\begin{eqnarray}\label{ome-5}
&\!\!\!\!\!\!\!\!\!\!\!\!\!\!\!\!\!\!\!\!&
{\rm(i)\ }c_\a\in\C\mbox{ if }\a>\frac{-m+1}{m},\ \ \ \ \
{\rm(ii)\ }c_{-1}=0,\ \ \ \ \ {\rm(iii)\ }\deg_xc_{\frac{-m+1-j}{m}}\le j+1\mbox{ if }j\in\Z_{\ge0},
\nonumber\\&\!\!\!\!\!\!\!\!\!\!\!\!\!\!\!\!\!\!\!\!&
{\rm(iv)\ } c_{\frac{-m+1}{m}}= -\frac{ J_{\rZeRo}}{m} (x+a_1)
.\end{eqnarray}
%
\end{lemm}

\noindent{\it Proof.}~First by \eqref{ome-2}, we have
\equa{pat-G-y}{\mbox{ $\dis \frac{\ptl G}{\ptl y}
{=}\mbox{$\sum\limits_{\a\in A}$}\a  c_\a F^{\a-1 }\frac{\ptl F}{\ptl y}.$}}
Then by \eqref{ome-2}, we obtain
\begin{eqnarray}\label{G-toF-}
&\!\!\!\!\!\!\!\!\!\!\!\!\!\!\!\!\!\!\!\!&
 \frac{\ptl G}{\ptl x}\stackrel{{}^{\sc\rm\eqref{ome-2}}}{=}
 \mbox{$\sum\limits_{\a\in A}$} \frac{dc_\a}{dx}F^\a+ \mbox{$\sum\limits_{\a\in A}$}\a  c_\a F^{\a-1}\frac{\ptl F}{\ptl x}
 \stackrel{{}^{\sc\rm\eqref{pat-G-y}}}{=} \mbox{$\sum\limits_{\a\in A}$} \frac{dc_\a}{dx}F^\a+  \frac{\ptl G}{\ptl y} \Big( \frac{\ptl F}{\ptl y}\Big)^{-1} \frac{\ptl F}{\ptl x}
.
  \end{eqnarray}
Equivalently, $ \frac{\ptl G}{\ptl x} \frac{\ptl F}{\ptl y}  = \mbox{$\sum_{\a\in A}$} \frac{dc_\a}{dx}F^\a \frac{\ptl F}{\ptl y}+  \frac{\ptl G}{\ptl y} \frac{\ptl F}{\ptl x}$, i.e.,
  \begin{eqnarray}\label{G-toF}
&\!\!\!\!\!\!\!\!\!\!\!\!\!\!\!\!\!\!\!\!&
-J_0=\mbox{$\sum\limits_{\a\in A}$} \frac{dc_\a}{dx}F^\a \frac{\ptl F}{\ptl y}.
  \end{eqnarray}
%
%
This gives (i) below, then  multiplying \eqref{ome-2} with $\frac{\ptl F}{\ptl y}$ gives (ii) below, i.e.,
\begin{eqnarray}\label{ome-3}
&\!\!\!\!\!\!\!\!\!\!\!\!\!\!\!\!\!\!\!\!\!\!\!\!\!&
{\rm(i)\ }J_{\rZeRo}\Big(\frac{\ptl F}{\ptl y}\Big)^{-1}\stackrel{{}^{\sc\rm\eqref{G-toF}}}{=}
-\mbox{$\sum\limits_{\a\in A}$} \frac{dc_\a}{dx}F^\a,\ \ \
{\rm(ii)\ }G\frac{\ptl F}{\ptl y}\stackrel{{}^{\sc\rm\eqref{ome-2}}}{=}\frac{\ptl}{\ptl y}\Big(\mbox{$\sum\limits_{-1\ne\a\in A}$}\frac{ c_\a F^{\a+1}}{\a+1}\Big)+c_{-1}F^{-1}\frac{\ptl F}{\ptl y}.
\end{eqnarray}

The $y$-degree of the left-hand side of \eqref{ome-3}\,(i) is $-m+1$, while that of the right-hand side is $m\a_0$, where $\a_0\in A$ is the maximal number with $\frac{dc_{\a_0}}{dx}\ne0$.
We obtain \eqref{ome-5}\,(i).

Since, for some $\b_i,\gamma_i\in\C[x]$,
\equa{inv}{ {\dis \frac{\ptl F}{\ptl y}}\stackrel{{}^{\sc\rm\eqref{wheraraaa+1}\,(i)}}{=}m y^{m-1}\Big(1+\sum\limits_{i=1}^\infty \b_i y^{-i}\Big)\implies
\Big( {\dis\frac{\ptl F}{\ptl y}} \Big)^{-1}={\dis\frac1m} y^{-m+1}\Big(1+\sum\limits_{i=1}^\infty \gamma_i y^{-i}\Big),} with
$(1+\sum_{i=1}^\infty \gamma_i y^{-i})=(1+\sum_{i=1}^\infty \b_i y^{-i})^{-1}$,
comparing coefficients of $y^{-m+1}$ in \eqref{ome-3}\,(i) gives that $\frac{dc_{\a_0}}{dx}=-m^{-1}J_0$ with $\a_0=\frac{-m+1}{m}$, which implies \eqref{ome-5}(iv)
for some $a_1\in\C$.

By comparing coefficients of $y^{-1}$ in \eqref{ome-3}\,(ii), we obtain \eqref{ome-5}\,(ii).

By \eqref{wheraraaa+1}\,(iii),\,(iv), we in particular have
\begin{eqnarray}\label{F-G-degree}
&\!\!\!\!\!\!\!\!\!\!\!\!\!\!\!\!\!\!\!\!\!\!&
{\rm(i)\ }\deg_x\Coeff(F,y^{m-j})\le j,\ \  \ \ \ \ \ \deg_x\Coeff(G,y^{n-j})\le j,
\nonumber\\
&\!\!\!\!\!\!\!\!\!\!\!\!\!\!\!\!\!\!\!\!\!\!&
{\rm(ii)\ }\deg_x\Coeff\mbox{\large$\Big($} \frac{\ptl F}{\ptl y} , y^{m-1-j} \mbox{\large$\Big)$} =\deg_x\Coeff(F,y^{m-j})\le j\mbox{ \ \ for \ }j\in\Z_{\ge0}.
\end{eqnarray}
Then by   Remark \ref{oss1} 
we have:
\equa{ome-4}{\dis
\deg_x\Coeff\mbox{\large$\Big($} \Big(\frac{\ptl F}{\ptl y}\Big)^{-1},y^{-m+1-j}\mbox{\large$\Big)$}\le j,\ \ \ \
 \ \ \ \deg_x\Coeff( F^{\frac{n-j}{m}},y^{n-j-k})\le k.}
%
From this and \eqref{ome-3}\,(i), we claim that \eqref{ome-5}\,(iii) holds by induction on $j$ by comparing the coefficients  of  $y^{ -m+1-j}$.

To prove the claim, from    Remark  \ref{oss1} and formula \eqref{root} one has for any $a\in\Z$,
\equa{root2}{
F^{\frac{a}{m}}=y^a\Big(1+\sum\limits_{i=1}^{\infty }\,\ell_{i,a}(x) y^{-i}\Big)   \quad  \mbox{ for some }\ell_{i,a}(x)\in\C[x]\mbox{ with }\deg_x \ell_{i,a}(x)\leq i.}       By \eqref{inv}, formula \eqref{ome-3}\,(i) becomes
 \begin{eqnarray}\label{New-===}
 \!\!\!\!\!\!\!\!\!\!\!\!\!\!\!\!\!\!\!\!\!&&
 \frac{J_0}m y^{-m+1}\Big(1+\mbox{$\sum\limits_{i=1}^\infty$} \gamma_i y^{-i}\Big)\stackrel{{}^{\rm\eqref{ome-3}\,(i),\,\eqref{root2}}}{=}
 -\mbox{$\sum\limits_{\a\in A}$} \frac{dc_\a}{dx}y^{m\a} \Big(1+\mbox{$\sum\limits_{i=1}^{\infty }$}\ell_{i,m\a}(x) y^{-i}\Big)
 \nonumber\\ \!\!\!\!\!\!\!\!\!\!\!\!\!\!\!\!\!\!\!\!\!&&
\phantom{\frac{J_0}m y^{-m+1}\Big(1+\mbox{$\sum\limits_{i=1}^\infty$} b_i y^{-i}\Big)}\ \ \ \ \ \ =\ \ \
 -\mbox{$\sum\limits_{\a\in A}$} \frac{dc_\a}{dx}  \Big(y^{m\a}+\mbox{$\sum\limits_{i=1}^{\infty }$}\ell_{i,m\a}(x) y^{m\a-i}\Big).
 \end{eqnarray}  The coefficient    of  $y^{ -m+1-j}$  of the left-hand side of the previous formula is $\frac{J_0}m \gamma_j$ if $j\geq 1$ of degree $\leq j$, it is  $\frac{J_0}m$ for $j=0$.  Instead, that of the right-hand side of the previous formula is $- \frac d {dx}  c_{ \frac{-m-(j-1)}m}$ plus  the sum for $\a,i\geq 1$  with  $m\a-i=-m+1-j$ [or $m\a =-m -(j-i-1)\ssc\,$] of  $-\frac{dc_\a}{dx}\ell_{i,m\a}(x) $.  By what we have proved only  $\a\leq \frac{-m+1}m$ contribute that is $-(j-i-1)\leq 1$ or $i\leq j$.
By induction  the degree of  $-\frac{dc_\a}{dx}$ is $\leq  j-i $ so these terms contribute to a degree $\leq j =(j-i)+i$ and the claim follows.
 The basis for the induction is for $j=0$  for which  we have \eqref{ome-5} (iv)  degree 1. 
\hfill$\Box$
\vskip5pt
We would like to mention that the following lemma [which in particular implies \eqref{nen1111}$\ssc\,$] is important for us to obtain that the $m$-th  root $\omega'$ of unity appearing in the proof of Lemma \ref{PSerConv} is equal to $1$ [cf.~\eqref{tcjne0}\,(ii) and the arguments after \eqref{T==0}$\ssc\,$].
\begin{lemm}\label{Sect3-2}
By some reformulation of $(F,G)$, we may \vspace*{-5pt}assume \equa{two-cnot=0}{\dis c_{\frac{4-m}{m}}c_{\frac{3-m}{m}}\ne0.}
\end{lemm}\noindent{\it Proof.~}~Fix any  $a_{\rZeRo}\in\C$ with $\a_0\ne3a_1$. We define $\bar F,\bar G$ below (then $\bar F,\bar G$ are still polynomials giving a counterexample to the Jacobian conjecture)
  such that $\bar F$ have the form (i) below, for  some $\bar f_{j,k}\in\C$\vspace*{-6pt},
\begin{eqnarray}
\label{ome-6}
\!\!\!\!\!\!\!\!\!\!\!\!\!\!\!\!\!\!&&
{\rm(i)\ }\bar F:=F(y,y^3+a_{\rZeRo} y^2-x)\stackrel{{}^{\sc\rm\eqref{wheraraaa}\,(i)}}{=}
(y^3+a_{\rZeRo} y^2-x)^m{\ssc}+{\ssc}\mbox{$\sum\limits_{j=0}^{m-1}\Big(\sum\limits_{k=0}^{m-1-j}f_{j,k}y^k\Big)(y^3+a_{\rZeRo} y^2-x)^j$}
\nonumber\\[-8pt]
\!\!\!\!\!\!\!\!\!\!\!\!\!\!\!\!\!\!&&
\phantom{{\rm(i)\ }\bar F}=
y^{\bar m}\Big(1+\mbox{$\sum\limits_{j=1}^{\bar m}\sum\limits_{k=0}^{j-1}$}\bar f_{j,k}y^{-j}x^k\Big),
\nonumber\\[4pt]
\!\!\!\!\!\!\!\!\!\!\!\!\!\!\!\!\!\!&&
 {\rm(ii)\ }\bar G=G(y,y^3+a_{\rZeRo} y^2-x),\!\!\!
 \end{eqnarray}
where we have used the same symbol with a bar to denote an associated element corresponding to the Jacobian pair $(\bar F,\bar G$); in particular, $\bar n=3 n,\,\bar m=3m$
.

When  $(x,y)$ is set to $(y,y^3+a_{\rZeRo}y^2-x)=\big(y,y^3(1+a_{\rZeRo}y^{-1}-xy^{-3})\big)$, we have that $y^{-1}$ is set to $y^{-3}(1+a_{\rZeRo}y^{-1}-xy^{-3})^{-1}=y^{-3}(1-a_{\rZeRo}y^{-1}+\sum_{j=2}^\infty t_j (x)y^{-j})$ for some $t_j (x)\in\C[x]$. By Remark \ref{oss1},  $\deg_xt_j<j$.

So  for $H\in  \C[x]((y^{-1}))$ the \vspace*{-5pt}substitution
\equa{phhhi}{\mbox{$H\mapsto  H|_{(x,y)=(y,y^3+a_{\rZeRo}y^2-x)}:=\phi(H)$,}} is a well-defined homomorphism which we denote for simplicity by $\phi$.
Note from \eqref{wheraraaa+1}\,(i) or \eqref{ome-6}\,(i) that $\bar F=\phi( F)$ has the form,
\begin{eqnarray}
\!\!\!\!\!\!\!\!\!\!\!\!\!\!\!\!\!\!\!\!&&
\label{form-barF}
\bar F=\phi( F)\stackrel{{}^{\sc\rm\eqref{wheraraaa+1}\,(i),\,\eqref{ome-6}\,(i)}}{=}(y^3+a_{\rZeRo} y^2-x)^m+\mbox{(terms with $y$-degree $\le 3m-2)$}
\nonumber\\\!\!\!\!\!\!\!\!\!\!\!\!\!\!\!\!\!\!\!\!&&
\phantom{\bar F}=
y^{\bar m}+m a_0y^{\bar m-1}+\mbox{(terms with $y$-degree $\le 3m-2)$}.\end{eqnarray}In particular, since $ \phi( F^{\frac{-m-i}{m}} )  $  is an $m$-th root of   $ \phi( F)^{ -m-i }    $, we
have
\equa{Mthroot}{\dis
 \phi( F^{\frac{-m-i}{m}} )  =\phi( F)^{\frac{-m-i}{m}}   =    \phi( F)^{\frac{-\bar m-3i}{\bar m}}   .}
Thus, we have the following, where the omitted terms are terms with lower $y$-degrees, or with lower powers of $\phi(F)$, 
\begin{eqnarray}\label{that-y}\label{gliF}
&\!\!\!\!\!\!\!\!\!\!\!\!\!\!\!\!\!\!\!\!&
{\rm(i)\ }\phi( F)^{\frac{1}{\bar m}}
\stackrel{{}^{\sc\rm\eqref{P-as},\,\eqref{form-barF}}}{
=}y+\frac{a_0}{3}+\cdots\implies \phi( F)^{\frac{h}{\bar m}}
=y^h+ \frac{a_0}{3}h y^{h-1}+\cdots,\ \mbox{ \ for $h\in\Z$, \ \ thus, }
\nonumber\\
&\!\!\!\!\!\!\!\!\!\!\!\!\!\!\!\!\!\!\!\!&
\nonumber
{\rm(ii)\ } y^h
=\phi( F)^{\frac{h}{\bar m}}-\frac{a_0}{3}h\phi( F)^{\frac{h-1}{\bar m}}+\cdots\ \
\in\mbox{$\sum\limits_{k=0}^\infty$}\C[x]\bar F^{\frac{h-k}{\bar m}},
\\
&\!\!\!\!\!\!\!\!\!\!\!\!\!\!\!\!\!\!\!\!&
{\rm(iii)\ }\phi\Big(c_{\frac{-m-i}{m}}F^{\frac{-m-i}{m}}\Big) =
\mbox{$\sum\limits_{k=0}^\infty$}p_{k,i}(x) \phi( F)^{\frac{2-2i-k-\bar m}{\bar m}}\quad \text{for some $p_{k,i}\in\C[x]$ and }\  i\geq -1,
\end{eqnarray}
where,  by the fact in \eqref{ome-5}\,(iii) that $c_{\frac{-m-i}{m}}$ for $i\ge-1$ has $x$-degree $\le i+2$, we see that $\phi(c_{\frac{-m-i}{m}})$ has $y$-degree $\le i+2$, and so (iii) has $y$-degree $\le i+2-3(m+i)=2-2i-\bar m$, and thus we obtain (iii) from (ii).

 Note that  one can exchange the substitution with the series in the definition of
\begin{eqnarray}
\label{MEMEMR--}
&\!\!\!\!\!\!\!\!\!\!\!\!\!\!\!\!\!\!\!\!\!\!\!\!\!\!\!\!\!\!\!\!&
 \bar G:=\phi( G)\stackrel{\rm\eqref{ome-2}}{=}\phi\Big( \mbox{$\sum\limits_{i=0}^\infty$}
  c_{\frac{n-i}{m}} F^{\frac{n-i}{m}}\Big)= \mbox{$\sum\limits_{i=0}^\infty$}
   \phi\Big(c_{\frac{n-i}{m}} F^{\frac{n-i}{m}}\Big)=
   \mbox{$\sum\limits_{i=0}^\infty$} \phi\Big(c_{\frac{-\bar m+(\bar n+\bar m -3i)}{\bar m}} F^{\frac{-\bar m+(\bar n+\bar m -3i)}{\bar m}}\Big) \!\!\!\!\!\!\!\!\!\!\!\!\!\!\!\!\!\!\!\!\!\!\!\!\!\!\!\!
\nonumber\\
&\!\!\!\!\!\!\!\!\!\!\!\!\!\!\!\!\!\!\!\!\!\!\!\!\!\!\!\!\!\!\!\!&
\phantom{\bar G} = \mbox{$\sum\limits_{-(\bar n+\bar m)\le i<\infty,\,i\in3\Z}$}
 \phi\Big(c_{\frac{-\bar m+i}{\bar m}} F^{\frac{-\bar m-i}{\bar m}}\Big).
\end{eqnarray}
Write the right-hand side of \eqref{MEMEMR--} as
\equa{the-right-hand-of}{
\mbox{ r.h.s.~of \eqref{MEMEMR--}}\ \ =\ \ \sum\limits_{i=0}^\infty \bar c_{\frac{\bar n-i}{\bar m}}\bar F^{\frac{\bar n-i}{\bar m}}\mbox{ \ for some \ $\bar c_{\frac{\bar n-i}{\bar m}}\in\C[x]$.}} Then of course we have the bar version of Lemma \ref{Sect3-1}; in particular, $\bar c_{\frac{-\bar m+i}{\bar m}}\in\C$ for $i\ge2$ by \eqref{ome-5}\,(i).
\begin{clai}\label{claim2.1}
Only the term with $i=-3$ in the right-hand side of \eqref{MEMEMR--} can contribute to
 $\bar c_{\frac{-\bar m+3}{\bar m}},\bar c_{\frac{-\bar m+4}{\bar m}}$.
\end{clai}

This can be proven as follows. If $i\le-6$ then $\frac{-\bar m-i}{\bar m}>\frac{1-m}{m}$ and $c_{\frac{-\bar m-i}{\bar m}}$ is a constant by \eqref{ome-5}\,(i) and thus the correspondent term is equal to $c_{\frac{-\bar m-i}{\bar m}}\bar F^{\frac{-\bar m-i}{\bar m}}$, which cannot contribute;
if $i=0$ the term is zero by \eqref{ome-5}\,(ii); if $i\ge3$ then  by \eqref{gliF}\,(iii) the term has $y$-degree $\le 2-\frac{2i}{3}-\bar m<-\bar m+3$, which cannot contribute. Thus the claim follows.

\NOUSE{
 Then
 the  terms from \eqref{gliF}  contribute  to $ c_{\frac{-\bar m+3}{\bar m}}, c_{\frac{-\bar m+4}{\bar m}}$ only for $i=-1,k=0$ (for $i\geq 0$ the  degree in $y$  equal to $2-2i-k-\bar m$ is   strictly less than $\bar m+3$  do not).
}%
\bigskip
Now consider the term  with $i=-3$ in the right-hand side of \eqref{MEMEMR--}. We prove
\begin{eqnarray}
\label{Next-Claim}
&\!\!\!\!\!\!\!\!\!\!\!\!\!\!\!\!\!\!\!\!\!\!\!\!\!\!\!\!\!\!\!\!\!\!\!\!&
\phi\Big(c_{\frac{-m+1}{m}}F^{\frac{-m+1}{m}}\Big)
\stackrel{{}^{^{\sc\rm\eqref{ome-5}\,(iv),\,\eqref{phhhi}}}}{=}
-\frac{ J_{\rZeRo}}{m} (y+a_1)\phi ( F)^{\frac{-\bar m+3}{\bar m}}
\nonumber\\
&\!\!\!\!\!\!\!\!\!\!\!\!\!\!\!\!\!\!\!\!\!\!\!\!\!\!\!\!\!\!\!\!\!\!\!\!&
\phantom{\phi\Big(c_{\frac{-m+1}{m}}F^{\frac{-m+1}{m}}\Big)
\ \,\ \ \ \ \ \,}
\in\ \ \,\
 -\frac{ J_{\rZeRo}}m\Big(\phi ( F)^{\frac{4-\bar m}{\bar m}}
 {\ssc\!}+{\ssc\!}
 \big(a_{\rOnE}{\ssc\!}-{\ssc\!}
 \frac{a_{\rZeRo}}{3}\big)\phi ( F)^{\frac{3-\bar m}{\bar m}}
 {\ssc\!}+{\ssc\!}\mbox{$\sum\limits_{j=0}^\infty$}\C[x]\phi ( F)^{\frac{2-\bar m-j}{\bar m}}\Big).
\!\!\!\!\!\!\!\! \end{eqnarray}
In fact,
\equa{2amsmwe}{\!\!\!\!
y\,  \phi ( F)^{\frac{-\bar m+3}{\bar m}}\stackrel{{}^{\sc\rm\eqref{that-y}\,(i)}}{=}y^{-\bar m+4 }  +(-\bar m+3)\frac{a_0}{3}y^{-\bar m+3 }    +T_1 \ \mbox{ for some }\
T_1\in\mbox{$\sum\limits_{k=-2}^\infty$}\C[x] \phi( F)^{\frac{-\bar m-k}{\bar m}}.\!\!\!\!}
From formula \eqref{that-y}  we deduce
\begin{eqnarray}
\label{that-y1}
&\!\!\!\!\!\!\!\!\!\!\!\!\!\!\!\!\!\!\!\!&y^{-\bar m+4 }
\stackrel{{}^{\sc\rm\eqref{that-y}\,(ii)}}{=}\phi( F)^{\frac{-\bar m+4}{\bar m}}-(-\bar m+4)\frac{a_0}{3}\phi( F)^{\frac{-\bar m+3}{\bar m}}+T_2\ \ \mbox{ for some }\ \
T_2\in\mbox{$\sum\limits_{k=-2}^\infty$}\C[x] \phi( F)^{\frac{-\bar m-k}{\bar m}},\!\!\!\!\!\!\!\!\!
\nonumber\\
&\!\!\!\!\!\!\!\!\!\!\!\!\!\!\!\!\!\!\!\!&y^{-\bar m+3 }
\stackrel{{}^{\sc\rm\eqref{that-y}\,(ii)}}{=}\phi(F)^{\frac{-\bar m+3}{\bar m} } +T_3\ \ \mbox{ for some }\ \
T_3\in\mbox{$\sum\limits_{k=-2}^\infty$}\C[x] \phi( F)^{\frac{-\bar m-k}{\bar m}}.
\end{eqnarray}
Thus the part ``\,$\in$\,'' in \eqref{Next-Claim} follows from \eqref{2amsmwe},\,\eqref{that-y1}, i.e., we have \eqref{Next-Claim}.

By Claim \ref{claim2.1} and \eqref{Next-Claim}, we obtain
\equa{barTwo-C-i}{\dis
\bar c_{\frac{-\bar m+4}{\bar m}}=-\frac{J_0}{m}\ne0,\ \ \ \ \ \ \bar c_{\frac{-\bar m+3}{\bar m}}=-\frac{J_0}{m}(a_1-\frac{a_0}{3})\ne0.}
 Thus by replacing $(F,G)$ by $(\bar F,\bar G)$ [by \eqref{ome-6} we still have \eqref{wheraraaa} after the replacement], we have the lemma.\hfill$\Box$
\bigskip

From now on, we fix, once and for all, the Jacobian pair $(F, G)$ satisfying \eqref{wheraraaa} and Lemma \ref{Sect3-2}
[and obviously, $(F,G)$ is obtained from the original Jacobian pair by applying some automorphism of $\C[x,y]\ssc\,$].
\bigskip

In the rest of this section we regard all elements as in the ring \equa{Ring-R}{\dis
{\mathbf R}=\C[x^{\frac1m}]((y^{-1})),}
where $x^{\frac1m}$ is regarded as a parameter such that its $m$-th  power is $x$.

Now, to be more precise, we slightly generalize notions and notations in Definition \ref{contro}.

 \begin{defi}\label{Abs-con}\rm
 Let $R=\sum_{i,j}r_{ij}x^iy^j,\,Q=\sum_{i,j}q_{ij}x^iy^j\in{\mathbf R}$ with $r_{ij}\in\C,\,q_{ij}\in\R_{\ge0}$.
 \begin{itemize}\item[(i)]If  $|r_{ij}|\le q_{ij}$ for all possible $i,j$,
then we say $R$ is controlled by $Q$ with respect to $x,y$, and
denote
$R\DS_{x,y}Q$.
\item[(ii)]
For $x_{\rZeRo},y_{\rZeRo}\in\C$, if $\sum_{i,j}|r_{ij}x_{\rZeRo}^iy_{\rZeRo}^j|$ converges, then we say the series $R$
with respect to $x,y$
converges 
{\it strongly} when $(x,y)$ is set to $(x_{\rZeRo},y_{\rZeRo})$, and denote $R|_{(x,y)=(x_{\rZeRo},y_{\rZeRo})}=\sum_{i,j}r_{ij}x_{\rZeRo}^iy_{\rZeRo}^j$.
\end{itemize}
\end{defi}

Note that
Lemmas \ref{ds-lemm} and \ref{lemm2222} can be parallelly generalized.

Denote \equa{denote-t}{\dis t=(1+x^{\frac1m})^{m-1}\in{\mathbf R}.}
%

Throughout the rest of this section,
 we always use $\SS_j\in\R_{>0}$ to denote some fixed number for all possible $j$; for instance, we can take $\SS_{\rZeRo}=\sum_{j,k}|f_{j,k}|$ in \eqref{S-1-D}.

By \eqref{wheraraaa+1}\,(iii), we can observe below that every term $x^ay^b$ appearing in $F$ (with non-negative integral $a,b$) also appears in $\bS_1$ by using \eqref{denote-t} to expand $t^j$ in $\bS_1$,
thus we  obtain,
\begin{eqnarray}
\!\!\!\!\!\!\!\!\!\!\!\!\!\!\!\!\!\!\!\!&&
\label{S-1-D}
{\rm(i)\ }
F\stackrel{{}^{\sc\rm\eqref{denote-t}}}{\DS_{x,y}}\bS_1:= y^m\Big(1+\SS_{\rZeRo}\mbox{$\sum\limits_{j=1}^{m}$}\big(ty^{-1}\big)^j\Big)\DS_{x,y}
\bS:= y^m\Big(1+\SS_{\rZeRo}\mbox{$\sum\limits_{j=1}^{\infty}$}\big(ty^{-1}\big)^j\Big)
,\mbox{ \ where \ }
\nonumber\\
\!\!\!\!\!\!\!\!\!\!\!\!\!\!\!\!\!\!\!\!&&
\dis
{\rm(ii)\ }\bS=y^m\Big(1+\frac{\SS_{\rZeRo}ty^{-1}}{1-ty^{-1}}\Big).
\end{eqnarray}
\begin{rema}\rm\label{Con-F-rema}
The importance of \eqref{S-1-D} is that though the structure of the polynomial $F$ may be complicated (and in particular we do not have any information about the coefficient of $x^iy^j$ in $F$ for general $i,j\in\Z_{>0}$), we can always use the controlling function $\bS$, which has the very simple form in \eqref{S-1-D}\,(ii), to control $F$; consequently, we are able to choose the simple controlling function $\bP$ in \eqref{P=Sm1} to control $P$, which allows us to obtain the simple form of the formal inverse function of $\bP_\inv$ in \eqref{P=Sm1+1}; then we can conveniently use Lemma \ref{lemm2222} to obtain Lemma \ref{Sect-Lemm3}, which is the key to obtain Lemma \ref{PSerConv}.
\end{rema}
Now by \eqref{S-1-D}\,(ii) and definition \eqref{P---} we have \equa{F-ign---}{\dis
\bS_\ign\stackrel{{}^{\sc\rm \eqref{P---}}}{=}\frac{\SS_{\rZeRo}ty^{-1}}{1-ty^{-1}}.}
Thus by Lemma \ref{ds-lemm}, for $j\in\Z_{\ge0}$,
\equa{S-2-D}{\dis F^{\pm\frac{j}{m}}\stackrel{{}^{\sc\rm\eqref{MAMS1}\,(b),\,(c)}}{\DS_{x,y}} y^{\pm j}\Big(1-\bS_\ign\Big)^{-\frac{j}{m}}\stackrel{{}^{\sc\rm\eqref{F-ign---}}}{=}y^{\pm j}\Big(\frac{1-ty^{-1}}{1-(1+\SS_{\rZeRo})ty^{-1}}\Big)^{\frac{j}{m}}.}

We use \eqref{P-as} to define \equa{P==defineP}{\dis
P=F^{-\frac1m},} which can be written as the form in (i) below for  some $p_j\in\C[x]$, then we
use \eqref{S-2-D} and Lemma \ref{ds-lemm} to obtain,
\begin{eqnarray}
\label{P=Sm1}
&\!\!\!\!\!\!\!\!\!\!\!\!\!\!\!\!\!\!\!\!\!\!\!\!\!\!&
\dis{\rm(i)\ } P{\ssc}\stackrel{{}^{\sc\rm\eqref{P==defineP}}}{:=}{\ssc}F^{-\frac1m}=y^{-1}\Big(1+\mbox{$\sum\limits_{j=1}^\infty$} p_j y^{-j}\Big)
\stackrel{{}^{\sc\rm\eqref{S-2-D}}}{\DS_{x,y}} \frac{y^{-1}}{(1-\bS_\ign)^{\frac1m}}{\sc}\stackrel{{}^{\sc\rm\eqref{MAMS1}\,(b)}}{\DS_{x,y}}{\sc}\bP,\mbox{ \ where,}
\nonumber\\[-0pt]&\!\!\!\!\!\!\!\!\!\!\!\!\!\!\!\!\!\!\!\!\!\!\!\!\!\!&
{\rm(ii)\ }\bP{\ssc}:={\ssc}\frac{y^{-1}}{1{\ssc}-{\ssc}\bS_\ign}{\ssc}
\stackrel{{}^{\sc\rm\eqref{F-ign---}}}
{=}{\ssc}
\frac{y^{-1}(1{\ssc}-{\ssc}ty^{-1})}{1{\ssc}-{\ssc}(1{\ssc}+{\ssc}\SS_{\rZeRo})
ty^{-1}}
{\ssc}={\ssc}y^{-1}\Big(1{\ssc}+\frac{\SS_{\rZeRo}ty^{-1}}{1{\ssc}-{\ssc}(1{\ssc}+{\ssc}
\SS_{\rZeRo})ty^{-1}}\Big)\stackrel{{}^{\sc\rm\eqref{P---}}}{=}{\ssc}y^{-1}(1{\ssc}+{\ssc}\bP_\ign)\mbox{ with}\!\!\!\!\!\!\!\!\!\!\!\!\!
\nonumber\\[-0pt]&\!\!\!\!\!\!\!\!\!\!\!\!\!\!\!\!\!\!\!\!\!\!\!\!\!\!&
{\rm(iii)\ }\bP_\ign{\ssc}={\ssc}\frac{\SS_{\rZeRo}ty^{-1}}{1{\ssc}-{\ssc}(1{\ssc}+{\ssc}
\SS_{\rZeRo})ty^{-1}}.
\end{eqnarray}
Thus by definition \eqref{P---},
\equa{P-neg=1}{\dis \bP_\inv=y^{-1}(1{\ssc}-{\ssc}\bP_\ign)
\stackrel{{}^{\sc\rm\eqref{P=Sm1}\,(iii)}}{
=}y^{-1}\Big(1-\frac{\SS_{\rZeRo}ty^{-1}}{1{\ssc}-{\ssc}(1{\ssc}+{\ssc}
\SS_{\rZeRo})ty^{-1}}\Big)=\frac{y^{-1}\big(1-(1+2\SS_{\rZeRo})ty^{-1}\big)}{1-(1+\SS_{\rZeRo})ty^{-1}}.
}

For convenience, we denote
\equa{Mememdnbrnrn}{\dis\tilde c_j:=c_{-\frac jm}=\tilde c_{{\rZeRo},j}+\tilde c_{{\rOnE},j}\mbox{ with }\tilde c_{{\rZeRo},j}\in\C,\ \tilde c_{{\rOnE},j}\in x\C[x]\mbox{ for all possible }j.} Then by \eqref{ome-5} and Lemma \ref{Sect3-2},
\equa{tcjne0}{\dis
{\rm(i)\ }\tilde c_{{\rOnE},j}=0\mbox{ \ for \ }j\le m-2,\ \ \ \ \ {\rm(ii)\ }\tilde c_{m-4}\tilde c_{m-3}\ne0.}
\begin{nota}\rm\label{notation-cij}
We denote $\hat c_{{\rOnE},j}\in x\R_{{\ssc\ge0}}[x]$ to be the polynomial of $x$ which is obtained from $\tilde c_{{\rOnE},j}\in x\C[x]$ by replacing the coefficient
$\Coeff(\tilde c_{{\rOnE},j},x^k)$ by its absolute value for all possible $j,k$.
\end{nota}
\begin{lemm}\label{Sect-Lemm3}We have, for some $\SS_{\rOnE},\SS_2\in\R_{>0}$,
\begin{eqnarray}
\label{P=Sm1+2+1}
&\!\!\!\!\!\!\!\!\!\!\!\!\!\!\!\!\!\!\!\!\!\!\!\!\!\!&
 {\rm(i)\ }G=G_{\rZeRo}+G_{\rOnE}\mbox{ \ with}
\nonumber\\[4pt]
&\!\!\!\!\!\!\!\!\!\!\!\!\!\!\!\!\!\!\!\!\!\!\!\!\!\!&
 {\rm(ii)\ }G_{\rZeRo}:=G|_{(x,P)=(0,P)}\stackrel{{}^{\sc\rm\eqref{ome-2},\,\eqref{Mememdnbrnrn}}}{=}\mbox{$\sum\limits_{j=-n}^\infty$}\tilde c_{{\rZeRo},j}P^j
\nonumber\\[-4pt]
&\!\!\!\!\!\!\!\!\!\!\!\!\!\!\!\!\!\!\!\!\!\!\!\!\!\!&
\phantom{{\rm(ii)\ }G_{\rZeRo}}
  \DS_{x,P}
P^{-n}\Big(1+\frac{\SS_{\rOnE}P}{1-\SS_2P}\Big)\DS_{x,y}
\bP_\inv^{-n}\Big(1+\frac{\SS_{\rOnE}\bP}{1-\SS_2\bP}\Big),\!\!\!\!\!\!\!\!\!\!\!\!
\nonumber\\[4pt] \nonumber
&\!\!\!\!\!\!\!\!\!\!\!\!\!\!\!\!\!\!\!\!\!\!\!\!\!\!&
{\rm(iii)\ }G_{\rOnE}:=G-G_0\stackrel{{}^{\sc\rm\eqref{ome-2},\,\eqref{Mememdnbrnrn}}}{=}\mbox{$\sum\limits_{j=-n}^\infty$}\tilde c_{{\rOnE},j}P^j
\nonumber\\[-4pt]
&\!\!\!\!\!\!\!\!\!\!\!\!\!\!\!\!\!\!\!\!\!\!\!\!\!\!&
\phantom{{\rm(iii)\ }G_{\rOnE}}
\DS_{x,P}\mbox{$\sum\limits_{j=m-1}^\infty$}x\frac{d \hat c_{{\rOnE},j}}{dx}P^j
\DS_{x,P}J_{\rZeRo}m^{-1}xP^{m-1}\Big(1-\frac{\SS_{\rOnE}tP}{1-\SS_2tP}\Big)^{-1}.
\end{eqnarray}
\end{lemm}
\noindent{\it Proof.~}Solving $y^{-1}$ from \eqref{P-neg=1} we obtain that the formal inverse function of $\bP_\inv$, which by notations in \eqref{i-Faa},\,\eqref{ha--f-}, is $\bP_\inv^{\circ-1}$, is the following
 [noting that if \eqref{P-neg=1} is regarded as a quadratic equation on variable $y^{-1}$ then we have two solutions and obviously the other solution, when expanded as a series of $y^{-1}$, contains a nonzero constant term and thus is not the one we require], 
\begin{eqnarray}
\label{P=Sm1+1}
&\!\!\!\!\!\!\!\!\!\!\!\!\!\!\!\!\!\!\!\!\!\!\!\!\!\!&
\dis y^{-1}=\bP_\inv^{\circ-1}(\bP_\inv):=\frac{1+(1+\SS_{\rZeRo})t\bP_\inv
-B}
{2(1+2\SS_{\rZeRo})t},\ \ B=(1-\b_+t\bP_\inv)^{\frac12}(1-\b_-t\bP_\inv)^{\frac12},
\end{eqnarray}where we always regard an element $(1+a)^{\frac12}$ as the unique element defined by the formula \eqref{bimeformo}, and where $\b_\pm=1+3\SS_{\rZeRo}\pm2(\SS_0+2\SS_0^2)^{\frac12}=\big(\sqrt{1+2\SS_0}\pm\sqrt{\SS_0}\big)^2\in\R_{>0}$.

Regarding $B$ defined in \eqref{P=Sm1+1} as a function of $x,\bP_\inv$, by Definition \ref{Abs-con}\,(i) and using \eqref{MAMS1} (with $k=\frac12,\,q_0=1,\a=0$, and $Q=1-\b_\pm t\bP_\inv$, and thus $Q_\inv\stackrel{{}^{\sc\rm\eqref{P---}}}{=}Q$), we obtain, where the second ``\,$\DS_{x,\bP_\inv}$\,'' is obtained by the fact that $\b_-<\b_+$, \equa{ThatAADS}{\dis
B\stackrel{{}^{\sc\rm\eqref{MAMS1}}}{\DS_{x,\bP_\inv}}(1-\b_+t\bP_\inv)^{-\frac12}
(1-\b_-t\bP_\inv)^{-\frac12}\DS_{x,\bP_\inv}(1-\b_+t\bP_\inv)^{-1}.
} Using this, we can deduce from
\eqref{P=Sm1+1} the following,
\equa{W-negg}{\dis \bP_\inv^{\circ-1}(\bP_\inv)
\stackrel{{}^{\sc\rm \eqref{P=Sm1+1},\,\eqref{ThatAADS}}}
{\DS_{x,\bP_\inv}}\bP_\inv\Big(1+\frac{\SS_3t\bP_\inv}{1-\SS_4t\bP_\inv}\Big)\mbox{ for some }\SS_3,\SS_4\in\R_{>0}.}
Thus by \eqref{MAMS1},\,\eqref{HAHAHAHHJ}, we have, by choosing sufficiently large  $\SS_{\rOnE},\SS_2\in\R_{>0}$,
\begin{eqnarray}
\label{P=Sm1+2}
&\!\!\!\!\!\!\!\!\!\!\!\!\!\!\!\!\!\!\!\!\!\!\!\!\!\!&
\dis
{\rm(i)\ } y^{-1}
\ \ \ \ \stackrel{{}^{\sc\rm\eqref{HAHAHAHHJ}}}{\DS_{x,P}}\ \
\bP_\inv^{\circ-1}(P)\stackrel{{}^{\sc\rm\eqref{W-negg}}}{\DS_{x,P}}P\Big(1+\frac{\SS_3tP}{1-\SS_4tP}\Big)\DS_{x,P}P\Big(1+\frac{\SS_{\rOnE}tP}{1-\SS_2tP}\Big),
\nonumber\\[-4pt]
&\!\!\!\!\!\!\!\!\!\!\!\!\!\!\!\!\!\!\!\!\!\!\!\!\!\!&
{\rm(ii)\ }y
\ \,\stackrel{{}^{\sc\rm\eqref{P=Sm1+2}\,(i),\,\eqref{MAMS1}}}{\DS_{x,P}} \ \,P^{-1}\Big(1-\frac{\SS_3tP}{1-\SS_4tP}\Big)^{-1}\DS_{x,P} P^{-1}\Big(1+\frac{\SS_{\rOnE}tP}{1-\SS_2tP}\Big),
\end{eqnarray}
where the second ``\,$\DS_{x,P}$\,'' in (i) is obtained from \eqref{W-negg}, and (ii) is obtained from the first two ``\,$\DS_{x,P}$\,'' of (i).

By \eqref{P=Sm1+2}, \eqref{wheraraaa}\,(ii) (with the fact that $n|m$) and  \eqref{MAMS1},\,\eqref{P=Sm1}, we can obtain, for some $\SS_5,\SS_6\in\R_{>0}$,
\equa{TDSP}{\dis G\stackrel{{}^{\sc\rm\eqref{P=Sm1+2}, \eqref{wheraraaa}}}{\DS_{x,P}}P^{-n}\Big(1+\frac{\SS_6tP}{1-\SS_5tP}\Big)\stackrel{{}^{\sc\rm\eqref{MAMS1},\,\eqref{P=Sm1}}}{\DS_{x,y}}\bP_\inv^{-n}\Big(1+\frac{\SS_6t\bP}{1-\SS_5t\bP}\Big).}
From this,  \eqref{ome-2} and definition of $G_{\rZeRo}$ in \eqref{P=Sm1+2+1}\,(ii), we obtain (note that $t|_{x=0}=1$),
\equa{TDSP-0}{\dis G_{\rZeRo}=G|_{(x,P)=(0,P)}\stackrel{{}^{\sc\rm\eqref{TDSP}}}{\DS_{x,P}}P^{-n}\Big(1+\frac{\SS_6(t|_{x=0})P}{1-\SS_5(t|_{x=0})P}\Big)
\stackrel{{}^{\sc\rm\eqref{TDSP}}}{\DS_{x,y}}\bP_\inv^{-n}\Big(1+\frac{\SS_6
(t|_{x=0})\bP}{1-\SS_5(t|_{x=0})\bP}\Big),}
i.e., we have \eqref{P=Sm1+2+1}\,(ii) (by enlarging $\SS_1,\SS_2$).

By \eqref{S-1-D} and \eqref{MAMS1}\,(a),  we obtain
\equa{df-dy==}{\dis
m^{-1}(\frac{\ptl F}{\ptl y})\DS_{x,y}y^{m-1}\Big(1+\SS_{\rZeRo}\mbox{$\sum\limits_{j=1}^{m-1}$}\big(ty^{-1}\big)^j\Big).}
Thus we have,
%
%
%
\begin{eqnarray}
\label{S-s---1}
&\!\!\!\!\!\!\!\!\!\!\!\!\!\!\!\!\!\!\!\!\!\!\!\!\!\!\!&
\mbox{$\sum\limits_{j=-n}^\infty$}\frac{d \tilde c_{{\rOnE}j}}{dx}P^j
\stackrel{{}^{\sc\rm\eqref{P==defineP},\,\eqref{Mememdnbrnrn}}}{=}
\mbox{$\sum\limits_{\a\in A}$} \frac{dc_\a}{dx}F^\a
\stackrel{{}^{\sc\rm\eqref{ome-3}\,(i)}}{=}-J_{\rZeRo}\Big(\frac{\ptl F}{\ptl y}\Big)^{-1}
\nonumber\\[-4pt]&\!\!\!\!\!\!\!\!\!\!\!\!\!\!\!\!\!\!\!\!\!\!\!\!\!\!\!&
\phantom{\mbox{$\sum\limits_{j=-n}^\infty$}
\frac{d \tilde c_{{\rOnE}j}}{dx}P^j}
\ \ \stackrel{{}^{\sc\rm\eqref{MAMS1},\,\eqref{df-dy==}}}{
\DS_{x,y}}J_{\rZeRo}m^{-1}y^{-(m-1)}\Big(1-\SS_{\rZeRo}\mbox{$\sum\limits_{j=1}^{m-1}$}(ty^{-1})^j\Big)^{-1}
\!\!\!\!\!\!\!\!\!\!\!\!\!\!\!
\nonumber\\[-4pt]&\!\!\!\!\!\!\!\!\!\!\!\!\!\!\!\!\!\!\!\!\!\!\!\!\!\!\!&
\phantom{\mbox{$\sum\limits_{j=-n}^\infty$}\frac{d \tilde c_{{\rOnE}j}}{dx}P^j}
\ \ \
\ \ \stackrel{{}^{\sc\rm\eqref{P=Sm1+2}}}{\DS_{x,P}}\ \ \
J_{\rZeRo}m^{-1}y^{-(m-1)}\Big(1-\SS_{\rZeRo}\mbox{$\sum\limits_{j=1}^{m-1}$}(ty^{-1})^j\Big)^{-1} \Big|_{y^{-1}=
P(1+\frac{\SS_3tP}{1-\SS_4tP})}
\nonumber\\[-4pt]&\!\!\!\!\!\!\!\!\!\!\!\!\!\!\!\!\!\!\!\!\!\!\!\!\!\!\!&
\phantom{\mbox{$\sum\limits_{j=-n}^\infty$}\frac{d \tilde c_{{\rOnE}j}}{dx}P^j}
\ \ \ \ \ \DS_{x,P}\ \ \ J_{\rZeRo}m^{-1}
P^{m-1}\Big(1-\frac{\SS_{\rOnE}tP}{1-\SS_2tP}\Big)^{-1},
\end{eqnarray}
where the first ``\,$\DS_{x,P}$\,'' is obtained from \eqref{P=Sm1+2}\,(i) and the second is obtained by choosing sufficiently large $\SS_{\rOnE},\SS_2$.
Now the first ``\,$\DS_{x,P}$\,'' of \eqref{P=Sm1+2+1}\,(iii) is obvious and the second ``\,$\DS_{x,P}$\,'' is obtained from \eqref{S-s---1}.
This proves the lemma.\hfill$\Box$\vskip5pt
The following proposition is a stronger version of Theorem \ref{Theo-2}.
\begin{prop}\label{Sect3-Lemm5}
There exists 
$\SS_7\in\R_{>0}$
such that for any $(p_{\ZeRo},p_{\OnE})=\big((x_{\ZeRo},y_{\ZeRo}),(x_{\OnE},y_{\OnE})\big)\in V$ with
$
{\dH_{p_{\ZeRo},p_{\OnE}}\ge\SS_7,}$
we must have
\equa{mqp1234-2}{|y_{\ZeRo}|<\tau \dH_{_{\sc p_{\ZeRo},p_{\OnE}}}^{^{\sc\frac{\ssTH m}{\ssTH m+1}}},\ \ \ \ \ \ \ \ |y_{\OnE}|<\tau  \dH_{_{\sc p_{\ZeRo},p_{\OnE}}}^{^{\sc\frac{\ssTH m}{\ssTH m+1}}}.}
\end{prop} We prove this proposition by contradiction and in several steps.  Assume thus that the statement leading to formula \eqref{mqp1234-2}  does not hold.
\bigskip

{\bf Step 1.} In this step we prove that, under the previous hypotheses, either for this Jacobian pair $(F,G)$, there exists $(p_{{\ZeRo},i},p_{{\OnE},i})$ $=\big((x_{{\ZeRo},i},y_{{\ZeRo},i}),(x_{{\OnE},i},y_{{\OnE},i})\big)\in V$
for any $i\in\Z_{\ge1}$ satisfying,
\equa{mqp1234-4---}
{\!\!\!\!\!\!\!\!{\rm(i)\ }
{\dH_{p_{{\ZeRo},i},p_{{\OnE},i}}\ge \frac 12         i,}\ \ \
 x_{{\ZeRo},i}\neq  x_{{\OnE},i}, 
 \mbox{ \ \ \ } {\rm (ii)\  }|y_{{\ZeRo},i}|\ge\frac12 \tau \dH^{^{\sc \frac{\ssTH m}{\ssTH m+1}}}_{_{\sc p_{{\ZeRo},i},p_{{\OnE},i}}}, 
\!\!\!\!\!\!\!\!}
or
for an equivalent Jacobian pair $(\tilde F,\tilde G)$
obtained by a simple linear change of coordinates and still satisfying  Lemma \ref{Sect3-2},  the tilde version of \eqref{mqp1234-4---} holds.

\noindent{\it Proof of Step 1.~}Assume \eqref{mqp1234-2} does not hold, i.e., there exists $(p_{{\ZeRo},i},p_{{\OnE},i})$ $=\big((x_{{\ZeRo},i},y_{{\ZeRo},i}),(x_{{\OnE},i},y_{{\OnE},i})\big)\in V$
for any $i\in\Z_{>0}$ satisfying
$
{\dH_{p_{{\ZeRo},i},p_{{\OnE},i}}\ge i,}$
such that
\equa{mqp1234-4}{{\rm (i)\ \ }|y_{{\ZeRo},i}|\ge \tau \dH^{^{\sc \frac{\ssTH m}{\ssTH m+1}}}_{_{\sc p_{{\ZeRo},i},p_{{\OnE},i}}}, \ \ \ \mbox{or} \ \ \ {\rm (ii)\ \ }|y_{{\OnE},i}|\ge\tau \dH^{^{\sc \frac{\ssTH m}{\ssTH m+1}}}_{_{\sc p_{{\ZeRo},i},p_{{\OnE},i}}}
.}
Since at least one of the conditions in \eqref{mqp1234-4} must hold for infinitely many $i$'s, if  necessary by switching $p_{{\ZeRo},i}$ and $p_{{\OnE},i}$, we can assume \eqref{mqp1234-4}\,(i) holds for infinitely many $i$.
If necessary
by replacing the sequence by a subsequence
[if the sequence $(p_{{\ZeRo},i},p_{{\OnE},i})$ is replaced by the subsequence  $(p_{{\ZeRo},{i_j}},p_{{\OnE},{i_j}})$, then we always have $i_j\ge j$; thus we still have $\dH_{p_{{\ZeRo},i},p_{{\OnE},i}}\ge i$ after the replacement],
 we may assume  \eqref{mqp1234-4}\,(i)  holds for all $i$.  Finally  if for infinitely many indices $i$ we have ${ x_{{\ZeRo},i}\neq  x_{{\OnE},i}}$ still passing to a subsequence we may assume that   ${ x_{{\ZeRo},i}\neq  x_{{\OnE},i}}$ holds for all $i$ and the given  pair  $(F,G)$  satisfies the conditions of Step 1.

Thus  we may assume that ${ x_{{\ZeRo},i}=  x_{{\OnE},i}}$ holds for all $i$.
%
%
%
%
Take the equivalent Jacobian pair\equa{Bar-FandG}{\dis (\tilde F,\tilde G)=\Big(F\big(x+\frac{y}{2},y\big),G\big(x+\frac{y}{2},y\big)\Big).}
which corresponds to the change of coordinates
\equa{Change-varr}{\dis(\tilde x,\tilde y)=\big(x -\frac{y }{2},y \big).  }
Note that our $(F,G)$ in \eqref{Bar-FandG} was obtained,  in Lemma  \ref{Sect3-2} by the initial choice of $(F,G)$ defined in \eqref{ome-6}, and the coordinate change $(x,y)\mapsto  (y,y^3+a_0y^2-x)$  therefore, when we change $x$ to $x+\frac{y}{2}$ in \eqref{Bar-FandG},
it amounts to that the tilde version of $(F,G)$ is   obtained from the initial  $(F,G)$ by  mapping $(x,y)\mapsto (y,y^3+a_0y^2-\frac{y}{2}-x)$.
Then one can observe that for the tilde version of $(F,G)$, we
 still have the second equality of \eqref{form-barF} and thus we still have
\eqref{that-y}, Claim \ref{claim2.1} and \eqref{barTwo-C-i}.
Hence, the proof of Lemma \ref{Sect3-2} works again for $(\tilde F,\tilde G)$.
This proves that all results before Proposition \ref{Sect3-Lemm5} hold for the tilde version Jacobian pair $(\tilde F,\tilde G)$.

Denote the points  $(p_{{\ZeRo},i},p_{{\OnE},i})$  in the new coordinates as   $(\tilde p_{1,i},\tilde p_{2,i})$ below for all $i\in\Z_{>0}$, then \eqref{Bar-FandG} together with the fact that $(p_{1,i},p_{2,i})\in V$ shows that it is in $\tilde V$
(which is the tilde version of $V$), i.e.,
\equa{bar-p-ij}{\dis(\tilde p_{1,i},\tilde p_{2,i})=\big((\tilde x_{1,i},\tilde y_{1,i}),(\tilde x_{2,i},\tilde y_{2,i})\big)
\stackrel{{}^{\sc\rm\eqref{Change-varr}}}{=}\Big(\big(x_{1,i}-\frac{y_{1,i}}{2},y_{1,i}\big),\big(x_{2,i}-\frac{y_{2,i}}{2},y_{2,i}\big)\Big)\in\tilde V.}
Since  by construction  $ p_{{\ZeRo},i}\neq p_{{\OnE},i} $ and the assumption $x_{1,i}=x_{2,i}$  we have thus $\tilde x_{1,i}\ne \tilde x_{2,i}$ for all $i $.

Note that every inequality below [except the second inequality in (ii), which follows from (i)$\ssc\,$] is obtained from one of the following three facts: (1) formula
\eqref{mqp1234-4}\,(i)  hold and $h_{p_{1,i},p_{2,i}}\ge i$ by noting from arguments after  \eqref{mqp1234-4};
 (2) $h_{p_{1,i},p_{2,i}}=|x_{1,i}|+|y_{1,i}|+|x_{2,i}|+|y_{2,i}|$ by definition \eqref{dp0p1}; (3) $\pm(|a|-|b|)\le|a-b|\le|a|+|b|$ for all $a,b\in\C$, therefore we have
\begin{eqnarray}\label{EveryIn}
&\!\!\!\!\!\!\!\!\!\!\!\!\!\!\!\!\!\!\!\!&
{\rm (i)\ }\frac{i}{2}\le\frac{h_{p_{1,i},p_{2,i}}}{2}\le|x_{1,i}|+\frac{|y_{1,i}|}{2}+|x_{2,i}|+\frac{|y_{2,i}|}{2}\le\Big|x_{1,i}-\frac{y_{1,i}}{2}\Big|+|y_{1,i}|
+\Big|x_{2,i}-\frac{y_{2,i}}{2}\Big|+|y_{2,i}|
\nonumber\\
&\!\!\!\!\!\!\!\!\!\!\!\!\!\!\!\!\!\!\!\!&
\phantom{{\rm (i)\ }\frac{i}{2}}=\tilde h_{\tilde p_{1,i},\tilde p_{2,i}}\le|x_{1,i}|+\frac{3|y_{1,i}|}{2}+|x_{2,i}|+\frac{3|y_{2,i}|}{2}\le 2h_{p_{1,i},p_{2,i}},
\nonumber\\
&\!\!\!\!\!\!\!\!\!\!\!\!\!\!\!\!\!\!\!\!&
{\rm(ii)\ }|\tilde y_{1,i}|=|y_{1,i}|\ge \dH^{^{\sc \frac{\ssTH m}{\ssTH m+1}}}_{_{\sc p_{{\ZeRo},i},p_{{\OnE},i}}}
\stackrel{{}^{\sc\rm\eqref{EveryIn}\,(i)}}{\ge} \big(\frac12\big)^{\frac{m}{m+1}} \tilde h_{_{\sc\tilde p_{1,i},\tilde p_{2,i}}}^{^{\sc\frac{m}{m+1}}}\ge \frac12 \tilde h_{_{\sc\tilde p_{1,i},\tilde p_{2,i}}}^{^{\sc\frac{m}{m+1}}}.
\end{eqnarray}
This proves that the tilde version of \eqref{mqp1234-4---} holds.

For the purpose of proving Proposition \ref{Sect3-Lemm5}, to simplify notations, in case the tilde version of \eqref{mqp1234-4---} holds, we re-denote $(\tilde F,\tilde G)$ as $(F,G)$ until the end of the proof of Proposition \ref{Sect3-Lemm5} [note that we abuse notations $F,G$ just for simplicity, but we do not change our real pair $(F,G)$ after the proof of Proposition \ref{Sect3-Lemm5}, which means that our pair $(F,G)$ is still the pair that satisfies \eqref{mqp1234-2}$\ssc\,$]. In this way, the tilde in all notations is omitted, and so in particular, the conclusion in Step 1 can be uniformly written as the following, for $i\ge1$,
\equa{mqp1234-5+}{
{\rm(i)\ } h_{ p_{1,i}, p_{2,i}}\ge \frac 12i,\ \ \
{\rm(ii)\ }
| y_{{\ZeRo},i}|\ge \frac 12\tau \dH^{^{\sc \frac{\ssTH m}{\ssTH m+1}}}_{_{\sc p_{{\ZeRo},i}, p_{{\OnE},i}}},
\ \ \ {\rm(iii)\ } x_{1,i}\ne x_{2,i}.
}

\bigskip

{\bf Step 2.} In this tep we prove that, for the previous sequence of points, we have formula \eqref{sim1aqa}  that is
$\lim_{i\to\infty}\frac{y_{{\OnE},i}}{y_{{\ZeRo},i}}=\omega,\ \ \mbox{where $\omega$ is  some $m$-th  root of unity.}$


%
\begin{nota}\rm\label{nota-sim}
Before continuing our proof, we need to use the following notations: Let $a,b,c$ be variables such that $c\to\infty$ or $c\to0$ and $a,b$ are functions of $c$, we denote
\equa{MSM1111}{\mbox{${\rm(i)\ }a\sim_{c\ssc\,}  b,\ \ \ \ {\rm(ii)\ }a\prec_{c\ssc\,} b,\ \ \ \ {\rm(iii)\ }a\preceq_{c\ssc\,} b$,}} which mean respectively, for some fixed $\SS_{\rOnE},\SS_2\in\R_{>0}$,
\equa{MSM1111+}{\mbox{$\dis{\rm(i)\ }\SS_{\rOnE}\le
 \Big|\frac{a}{b} \Big|\le\SS_2,\ \ \ \ {\rm(ii)\ }\lim\limits_{c\to\infty} \frac{a}{b}=0,\ \ \ \ {\rm(iii)\ }\Big|
\frac{a}{b}\Big|\le\SS_{\rOnE}$.}}
Then obviously, properties (i),\,(ii),\,(iii) are respectively an equivalence, a strict order, an order.
\end{nota}

Now for $k=\ZeRo,\OnE$, since $|x_{k,i}|\le \dH_{p_{{\ZeRo},i},p_{{\OnE},i}},\,|y_{k,i}|\le \dH_{p_{{\ZeRo},i},p_{{\OnE},i}}$,
by \eqref{wheraraaa}\,(i), we have, where (ii) is obtained from (i) and \eqref{wheraraaa}\,(i),
\equa{S1pq}{\dis{\rm(i)\ }
F_{\rOnE}(x_{k,i},y_{k,i})\stackrel{{}^{\sc\rm\eqref{wheraraaa}\,(i)}}{\preceq_{i\ssc\,}} \dH_{p_{{\ZeRo},i},p_{{\OnE},i}}^{m-1}\prec_{i\ssc\,} \dH_{_{\sc p_{{\ZeRo},i},p_{{\OnE},i}}}^{^{\sc\frac{m^2}{m+1}}}\stackrel{{}^{\sc\rm\eqref{mqp1234-5+}\,(ii)}}{\preceq_{i\ssc\,}} |y_{{\ZeRo},i}|^m,\ \ \ \ \ \ \ {\rm(ii)\ }F(x_{{\ZeRo},i},y_{{\ZeRo},i})\sim_{i\ssc\,} y^m_{{\ZeRo},i}.
}
Thus, where the first equality follows from the fact that $\si(p_{1,i})=\si(p_{2,i})$, 
\begin{eqnarray}\label{sim1a?qa}
&\!\!\!\!\!\!\!\!\!\!\!\!\!\!&
\dis1=\frac{F(x_{{\OnE},i},y_{{\OnE},i})}{F(x_{{\ZeRo},i},y_{{\ZeRo},i})}
=\lim_{i\to\infty}\frac{F(x_{{\OnE},i},y_{{\OnE},i})}{F(x_{{1},i},y_{{1},i})}
\stackrel{{}^{\sc\rm \eqref{wheraraaa}}}{=}\lim_{i\to\infty}
\frac{\frac{y^m_{{\OnE},i}}{y^m_{{\ZeRo},i}}+\frac{F_{\rOnE}(x_{{\OnE},i},y_{{\OnE},i})}{y^m_{{\ZeRo},i}}}{1+\frac{F_{\rOnE}(x_{{1},i},y_{{1},i})}{y^m_{{\ZeRo},i}}}
\stackrel{{}^{\sc\rm\eqref{S1pq}}}{=}\lim_{i\to\infty}\Big(\frac{y_{{\OnE},i}}{y_{{\ZeRo},i}}\Big)^m
.\end{eqnarray}
Therefore, by replacing the sequence by a subsequence, we have [we would like to remark that $\omega$ defined below is not necessarily equal to $1$, but we will prove that the $m$-th  root $\omega'$ of unity defined in \eqref{omega-p} must be $1$],
\equa{sim1aqa}{\dis
\lim_{i\to\infty}\frac{y_{{\OnE},i}}{y_{{\ZeRo},i}}=\omega,\ \ \mbox{where $\omega$ is  some $m$-th  root of unity.}}
\bigskip

{\bf Step 3.}  In this step we show that the series appearing in our algeraic identities can be evaluated at the points $p_{1,i},\,p_{2,i}$  giving rise to numerical estimates.\bigskip

Denote \equa{eps}{\dis\ep=\dH_{_{\sc p_{{\ZeRo},i},p_{{\OnE},i}}}^{^{\sc -\frac{1}{m(m+1)}}}\to0\mbox{ \ (when $i\to\infty$).}}
Let $a_i\in\C$ with $|a_i|\le \dH_{p_{{\ZeRo},i},p_{{\OnE},i}}$. 
Then  $1+|a_i|^{\frac1m} \preceq_{i\ssc\,}\dH_{_{\sc p_{{\ZeRo},i},p_{{\OnE},i}}} ^{^{\sc\frac1m}}$ so,
by \eqref{mqp1234-5+},\,\eqref{sim1aqa}, we have
\equa{tw}{\dis  (1+|a_i|^{\frac1m})^{m-1}|y_{k,i}|^{-1}\preceq_{i\ssc\,} \dH_{_{\sc p_{{\ZeRo},i},p_{{\OnE},i}}}^{^{\sc \frac{m-1}{m}}}|y_{{\ZeRo},i}|^{-1}\preceq_{i\ssc\,} \dH_{_{\sc p_{{\ZeRo},i},p_{{\OnE},i}}}^{^{\sc \frac{m-1}{m}-\frac{m}{m+1}}}=\ep\mbox{ for }k=\ZeRo,\OnE .
}
Using the notations in \eqref{Mememdnbrnrn} and the results  in formula \eqref{tcjne0}  and in formula \eqref{P=Sm1+2+1} of Lemma \ref{Sect-Lemm3}, we prove the following.
\begin{lemm}\label{PSerConv}Let $k=\ZeRo,\OnE $ and $a\in\C$ with $|a|\le \dH_{p_{{\ZeRo},i},p_{{\OnE},i}}$ and $i\gg1$.
\begin{itemize}\item[\rm(i)]
In {\rm\eqref{P=Sm1}\,(i)}, the series $P=F^{-\frac1m}$ with respect to $x,y$ converges strongly when $(x,y)$ is set to $(a,y_{k,i})$, and
\equa{PbP1}{P_{a,k}:= P|_{(x,y)=(a,y_{k,i})}=y^{-1}_{k,i}\Big(1+O(\ep)^1\Big).}
\item[\rm(ii)]
In {\rm\eqref{P=Sm1+2}\,(ii)}, the series $y$ with respect to $x,P$  converges strongly when $(x,P)$ is set to $(a,P_{a,k})$, and
\equa{Ya-y}{Y_{a,k}:=y|_{(x,P)=(a,P_{a,k})}=y_{k,i}\Big(1+O(\ep)^1\Big).}

\item[\rm(iii)]
In {\rm\eqref{P=Sm1+2+1}},
the series $G$ with respect to $x,P$ converges strongly when $(x,P)$ is set to $(a,P_{a,k})$, and
\begin{eqnarray}
&\!\!\!\!\!\!\!\!\!\!\!\!\!\!\!\!\!\!&
\label{PbP2}
{\rm(a)\ }A_{a,k,\ell}:=
\mbox{$\sum\limits_{j=\ell}^\infty$}\tilde c_{{\rZeRo}j}P_{a,k}^j\preceq_{i\ssc\,} P_{a,k}^{\ell}\sim_{i\ssc\,} y^{-\ell}_{{\ZeRo},k} \mbox{ \ for \ }\ell\ge-n,
\nonumber\\[-2pt]&\!\!\!\!\!\!\!\!\!\!\!\!\!\!\!\!\!\!&
{\rm(b)\ }B_{a,k}:=G_{\rOnE}|_{(x,y)=(a,y_{k,i})}\preceq_{i\ssc\,} \dH_{p_{{\ZeRo},i},p_{{\OnE},i}}y^{-(m-1)}_{k,i}\prec_{i\ssc\,} y^{-(m-3)}_{{\ZeRo},i}.
\end{eqnarray}
\item[\rm(iv)]
The series $\big(\frac{\ptl F}{\ptl y}\big)^{-1}$ with respect to $x,P$ converges strongly when $(x,P)$ is set to $(a,P_{a,k})$, and
\begin{eqnarray}
&\!\!\!\!\!\!\!\!\!\!\!\!\!\!\!\!\!\!&
\label{PbP3}
m\Big(\frac{\ptl F}{\ptl y}\Big)^{-1}\Big|_{(x,P)=(a,P_{a,k})}=y^{-(m-1)}_{k,i}\Big(1+O(\ep)^1\Big).
\end{eqnarray}
\item[\rm(v)]
$P_{\ZeRo}:=P_{x_{{\ZeRo},i},1}=P(x_{{\ZeRo},i},y_{{\ZeRo},i})=P(x_{{\OnE},i},y_{{\OnE},i})=P_{x_{{\OnE},i},2}$.

\end{itemize}\end{lemm}
\noindent{\it Proof.~}From \eqref{P=Sm1}\,(ii), we see that $\bP_{\ign}|_{(x,y)=(|a|,|y_{k,i}|)}$ converges to $O(\ep)^1$ by \eqref{tw}, thus \eqref{PbP1} follows from
\eqref{P=Sm1}\,(i). Similarly, we obtain \eqref{Ya-y}--\eqref{PbP3} from \eqref{P=Sm1+2},\,\eqref{P=Sm1+2+1},\,\eqref{S-s---1}.
This proves (i)--(iv).

To prove (v), we have $P(x_{{\ZeRo},i},y_{{\ZeRo},i})^{-m}
\stackrel{{}^{\sc\rm\eqref{P==defineP}}}{=}F(x_{{\ZeRo},i},y_{{\ZeRo},i})=F(x_{{\OnE},i},y_{{\OnE},i})
\stackrel{{}^{\sc\rm\eqref{P==defineP}}}{=}P(x_{{\OnE},i},y_{{\OnE},i})^{-m}$. Thus
\equa{omega-p}{\mbox{$P(x_{{\OnE},i},y_{{\OnE},i})$ $=\omega' P(x_{{\ZeRo},i},y_{{\ZeRo},i})$,}}
for some $m$-th  root $\omega'$ of unity (which may depend on $i$). Assume there exists $j\le m-3$ such that
\equa{suchthatttt}{\mbox{$\tilde c_j
\stackrel{{}^{\sc\rm\eqref{tcjne0}}}{=}\tilde c_{{\rZeRo}j}\ne0$ \ \ but \ \ $P(x_{{\ZeRo},i},y_{{\ZeRo},i})^j\ne P(x_{{\OnE},i},y_{{\OnE},i})^j\stackrel{{}^{\sc\rm\eqref{omega-p}}}{=}\omega'^jP(x_{{\ZeRo},i},y_{{\ZeRo},i})^j$.}}
Let $j_{\rZeRo}\le m-3$ be the minimal such $j$. Then $\omega'^{j_{\rZeRo}}\ne1$ and $|1-\omega'^{j_{\rZeRo}}|>\d$ for some fixed $\d>0$ (since $\omega'^{j_{\rZeRo}}$, which may though depend on $i$, is an $m$-th  root of unity).
By 
\eqref{P=Sm1+2+1} and Lemma \ref{PSerConv}\,(iii), we have [note that $\tilde c_{{\rZeRo},j_{\rZeRo}}\in\C_{\ne0}$ is a number independent of $i$ and $P(x_{{\ZeRo},i},y_{{\ZeRo},i})\sim_{i\ssc\,} y^{-1}_{{\ZeRo},i}$ by \eqref{PbP1}$\ssc\,$],
\begin{eqnarray}
&\!\!\!\!\!\!\!\!\!\!\!\!\!\!\!\!\!\!&
\label{T==0}
0\ \ \ \ =\ \ \ \ \ G(x_{{\ZeRo},i},y_{{\ZeRo},i})-G(x_{{\OnE},i},y_{{\OnE},i})
\nonumber\\&\!\!\!\!\!\!\!\!\!\!\!\!\!\!\!\!\!\!&\phantom{0}\!\!\!
\stackrel{{}^{\sc\rm\eqref{P=Sm1+2+1},\,\eqref{PbP2} }}{=}\tilde c_{{\rZeRo},j_{\rZeRo}}(1-\omega'^{j_{\rZeRo}})P(x_{{\ZeRo},i},y_{{\ZeRo},i})^{j_{\rZeRo}}+A_{x_{{\ZeRo},i},1,j_{\rZeRo}+1}-A_{x_{{\OnE},i},2,j_{\rZeRo}+1}+B_{x_{{\ZeRo},i},1}-B_{x_{{\OnE},i},2}
\nonumber\\&\!\!\!\!\!\!\!\!\!\!\!\!\!\!\!\!\!\!&\phantom{0}\!\!\!
\stackrel{{}^{\sc\rm\eqref{PbP1},\,\eqref{PbP2}}}{\sim_{i\ssc\,}}
P(x_{{\ZeRo},i},y_{{\ZeRo},i})^{j_{\rZeRo}}
\stackrel{{}^{\sc\rm\eqref{PbP1}}}{\sim_{i\ssc\,}} y^{-j_{\rZeRo}}_{{\ZeRo},i},
\end{eqnarray}
which is a contradiction. This proves the following crucial fact,
\equa{crucial-fact}{\mbox{$\omega'^j=1$ for all $j\le m-3$ with $\tilde c_j\ne0$.}} In particular by \eqref{tcjne0}\,(ii),
$\omega'^{m-4}=1$,\,$\omega'^{m-3}=1$, which   implies that $\omega'=1$. This proves (v) and the lemma.\hfill$\Box$
\vskip7pt
{\bf Step 4 (conclusion).~}
Now we can continue our proof of Proposition \ref{Sect3-Lemm5}.
%
By \eqref{ome-2},\,\eqref{P=Sm1+2+1} and Lemma \ref{PSerConv}\,(iv),\,(v), we have the following
,
where the second equality follows from \eqref{ome-2},\,\eqref{Mememdnbrnrn} and Lemma \ref{PSerConv}\,(iii),\,(v), the third equality follows from the fact  that $\tilde c_j$'s are polynomials in $x$   [the integration path is taken as $x(t):=x_{1,i}(1-t)+tx_{2,i} , \ 0\le t\le1{\ssc\,}$], the fourth follows from the fact that
the series there converges absolutely and  uniformly when  $|x|\le h_{p_{1,i},p_{2,i}}$ by Lemma \ref{PSerConv}, and the last two equalities can be observed from \eqref{S-s---1} (see Remark \ref{FifthEqrema} for the reason why we have the fifth equality),
\begin{eqnarray}
\label{GasFG-set}\!\!\!\!\!\!\!\!\!\!\!\!&\!\!\!\!\!\!\!\!\!\!\!\!\!\!\!\!\!\!\!\!\!\!\!\!&
0= G(x_{1,i},y_{{\ZeRo},i})-G(x_{2,i},y_{{\OnE},i})\stackrel{{}^{\sc\rm\eqref{ome-2},\,\eqref{Mememdnbrnrn}}}{=}-\mbox{$\sum\limits_{j=-n}^\infty$} \big( \tilde c_j(x_{2,i})-\tilde c_j(x_{1,i})\big)P_{\ZeRo}^j
\nonumber\\[-4pt]
&\!\!\!\!\!\!\!\!\!\!\!\!\!\!\!\!\!\!\!\!\!\!\!\!&\phantom{0}
=
-\mbox{$\sum\limits_{j=-n}^\infty$}\int_{x_{1,i}}^{x_{2,i}} \frac{d \tilde c_{j}}{dx}P_{\ZeRo}^jdx
=-\int_{x_{1,i}}^{x_{2,i}}\mbox{$\sum\limits_{j=-n}^\infty$} \frac{d \tilde c_{j}}{dx}P_{\ZeRo}^jdx
\stackrel{{}^{\sc\rm Remark\ \ref{FifthEqrema}}}{=}\int_{x_{1,i}}^{x_{2,i}}J_{\rZeRo}\Big(\frac{\ptl  F}{\ptl y}\Big)^{-1}\,\Big|_{(x,P)=(x,P_{\ZeRo})}dx\!\!\!\!\!\!
\nonumber\\[0pt]
&\!\!\!\!\!\!\!\!\!\!\!\!\!\!\!\!\!\!\!\!\!\!\!\!&\phantom{0}\!\!\!\!\!
\stackrel{{}^{\sc\rm\eqref{S-s---1}}}{=}(x_{2,i}-x_{1,i})J_{\rZeRo}m^{-1}y^{-(m-1)}_{1,i}\Big(1+O(\ep)^1\Big)\ne0,
\end{eqnarray}
which is a contradiction. This proves
Proposition 
\ref{Sect3-Lemm5}.
\hfill$\Box$
\begin{rema}\rm\label{FifthEqrema}
Note that the fifth equality of \eqref{GasFG-set} should be understood in this way: first we regard $\big(\frac{\ptl F}{\ptl y}\big)^{-1}$ as a power series of $P$ with coefficients in $\C[x]$, i.e., $J_0\big(\frac{\ptl  F}{\ptl y}\big)^{-1}=
\sum_{j=-n}^\infty\frac{d c_{j}}{dx}P^j$ by \eqref{S-s---1}, then we set $P$ to $P_1$ (and keep $x$ unchanged with $|x|\le h_{p_{1,i},p_{2,i}}$), and then Lemma \ref{PSerConv}\,(iv) shows that the series converges absolutely and uniformly, and so we have the fifth equality.
\end{rema}

\noindent{\it Proof of Theorem \ref{Theo-2}.~}~It follows immediately from Proposition \ref{Sect3-Lemm5}.\hfill$\Box$
\begin{rema}\rm\label{Rema-procesi+}
To prove Theorem \ref{MAINT}, we always need to make the whole use of the following three facts about $(F,G)$:
 \begin{itemize}\item[(F1)$\!\!\!\!\!\!\!\!$] \ \ \ \ $F\in\C[x,y^{\pm1}],$ $G\in\C[x]((y^{-1}))$; \item[(F2)$\!\!\!\!\!\!\!\!$] \ \ \ \ $F,G$ are locally holomorphic functions of $x,y$ on $\C^2$, and so  $F,G\in\C[x,y]$; \item[(F3)$\!\!\!\!\!\!\!\!$] \ \ \ \ $J(F,G)$ is nonzero anywhere in $\C^2$, and so $J(F,G)\in\C_{\ne0}$.\end{itemize}
However,
some results may not require all of  the above facts; for instance, our approach above shows that Theorem \ref{Theo-2} holds for any pair $(F,G)$
 satisfying (F1) as long as  the following holds [of course, when $(F,G)\notin\C[x,y]^2$, the map $\si$ in \eqref{Kell---} can be only defined on some subset of $\C^2$, and in case $G\notin\C[x,y^{\pm1}]$, we need some extra condition on $G$ so that it can be controlled somehow, cf.~\eqref{P=Sm1+2+1}$\ssc\,$]:
\begin{itemize}
 \item[(i)] $J(F,G)\in\C_{\ne0}$;
\item[(ii)]
 the homogenous part of $F$ with the highest degree contains  only one term, i.e., $y^m$ with $m>1$;
\item[(iii)]let $c_\a$ be defined in \eqref{ome-2}, then,
\equa{nen1111}{\dis
{\rm gcd}\Big(m,m\a\,\Big|\,\a\in\frac1m\Z,\a\ge\frac{3-m}{m},\,c_\a\ne0\Big)=1,}
where the symbol ``\,gcd\,'' denotes the greatest common divisor [cf.~\eqref{two-cnot=0} and statements after \eqref{T==0}, and noting that $\omega'$ appeared in \eqref{T==0} is some $m$-th  root of unity].
 \end{itemize}
 In particular, Theorem \ref{Theo-2} holds for the pair $(F,G)=(y^{m},y^{n}+cy^{k}+xy^{-m+1})$ for any $m,n,k\in\Z,$ $c\in\C_{\ne0}$ with $m>1$, $n>k\ge3-m$ and  ${\rm gcd}(m,n,k)=1$, which can be also easily proven directly; in this case, $c_{\frac{n}{m}}=1,$ $c_{\frac{k}{m}}=c\ne0.$ If \eqref{nen1111} does not hold, then the $m$-th  root $\omega'$ of unity  in \eqref{T==0}
is not necessarily equal to $1$ [if we denote the left-hand side of \eqref{nen1111} as $d$ then the proof of Lemma \ref{PSerConv}\,(v) shows that $\omega'^d=1$], and we can easily construct a counter example like this: $(F,G)=(y^6,y^4+xy^{-5})$; in this case, we can choose, say, $(p_{1i},p_{2i})=\big((x_{1i},y_{1i}),(x_{2i},y_{2i})\big)=\big((1,i),(-1,-i)\big)\in V$ for $i\in\Z_{>0}$ such that $\omega'=-1$ and Theorem \ref{Theo-2} does not hold.
\end{rema}

\subsection{Proof of Theorem \ref{Theo-3} provided by 
Claudio Procesi}
We  thank Professor Claudio Procesi for providing us the following simple proof by using a little geometry (for your reference, we also present our original  proof of surjectivity of $\pi_1$  without using geometry in Appendix \ref{sect7}).

\begin{lemm}\label{Procesi-1}
The set $V$ defined in \eqref{V=0} is a smooth and closed algebraic surface in $\mathbb C^4=\mathbb C^2\times \mathbb C^2$.
\end{lemm}
\noindent{\it Proof.~}~The map  $\phi:\C^4\to \mathbb C^2,$
$$\phi:\ (p_1,p_2)\mapsto\big(F(p_1)-F(p_2),G(p_1)-G(p_2)\big),$$ has Jacobian matrix  always of maximal rank  so that $\phi^{-1}(0,0)=V\cup D$ is a smooth closed surface where $D$ is the diagonal which is closed.  By any argument (the implicit functions theorem, or smoothness) $D$ is also  open in $\phi^{-1}(0,0)$ and the lemma follows.
\hfill$\Box$\vskip4pt
We need a simple  estimate.
\begin{lemm}\label{1}
Let $a,b\in \mathbb R,\ a,b>2$.  If
 $b<2(a+b)^{\frac m{m+1}}$  then $b<2(2 a ) ^{m}$.
\end{lemm}
\noindent{\it Proof.~}~We have  $ab-a-b=(\frac a2-1)b+(\frac b2-1)a>0$, so
$$ b<2(a+b)^{\frac m{m+1}}\implies b^{m+1}<2^{ m+1}(a+b)^{ m}<2^{ m+1}(a b)^{ m} \implies b<2(2a)  ^{m}.\eqno{\Box}$$
\vskip4pt

We apply this to $a=|x_1|+|x_2|,\ b=|y_1|+|y_2|$ with $\big((x_1,y_1),(x_2 , y_2)\big)\in V$ and to the projection $\pi_1$ [defined in \eqref{proj1}$\ssc\,$] and
deduce the following.
\begin{prop}\label{pro-also}
If $A\subset  \mathbb C^2$ is bounded  then also $\pi_1^{-1}(A)\subset V$ is bounded.
\end{prop}
\noindent{\it Proof.~}~This follows from Lemma \ref{1}, and Proposition \ref{Sect3-Lemm5}.
More precisely, fix any $\SS_8>\SS_7$ (which is the number in Proposition \ref{Sect3-Lemm5}) such that for all $(x_1,x_2)\in A$,
\equa{a-smalllThan}{\mbox{ $a:=|x_1|+|x_2|<\SS_8$.}}
 Let $(p_1,p_2)=\big((x_1,y_1),(x_2,y_2)\big)\in\pi_1^{-1}(A)$. We claim
\equa{b-is-smaller}{\dis
b:=|y_1|+|y_2|<2(2\SS_8)^m.}
Assume conversely that $b\ge2(2\SS_8)^m$. Then $h_{p_1,p_2}\stackrel{{}^{\sc\rm\eqref{dp0p1}}}{\ge} b\ge 2(2\SS_8)^m>\SS_7$. Thus we can apply Proposition \ref{Sect3-Lemm5} to obtain, \equa{b---that}{\mbox{ $b=|y_1|+|y_2|<2h_{_{\sc p_1,p_2}}^{^{\sc\frac m{m+1}}}\stackrel{{}^{\sc\rm\eqref{dp0p1},\,\eqref{a-smalllThan},\,\eqref{b-is-smaller}}}{=}2(a+b)^{\frac m{m+1}}.$}}
Then by Lemma \ref{1}, $b<2(2a)  ^{m}\stackrel{{}^{\sc\rm\eqref{a-smalllThan}}}{<}2(2\SS_8)  ^{m}$, which is a contradiction with the assumption. Thus we have \eqref{b-is-smaller} and the proposition.
\hfill$\Box$\vskip4pt
As a consequence we have  that the fibres of $\pi_1$,  being compact affine algebraic varieties, are finite   hence the closure, denoted $\overline{\pi_1(V)}$, of $\pi_1(V)$ is a 2-dimensional subvariety of   $ \mathbb C^2$ so it is  $ \mathbb C^2$.
\vskip5pt

\noindent{\it Proof of Theorem \ref{Theo-3}.~}~Now we can proceed the proof of Theorem \ref{Theo-3} as follows. We already see that $\pi_1$ is proper and finite. Now take any point $q\in \mathbb C^2=\overline{\pi_1(V)}$.  By the previous  discussion there is a sequence of points $\pi_1(p_i)$ for $p_i\in V$  with $\lim_{i\to \infty} \pi_1(p_i)=q$.  But the sequence $p_i$, by the previous proposition,  is bounded and $V$ is closed so we can extract a subsequence converging to some $p\in V$ and  $\pi_1(p)=q$.\hfill$\Box$\vskip5pt


\section{Proof of Theorem \ref{MAINT}}\label{sect3}
First we recall after the proof of Lemma \ref{Sect3-2} that we have fixed, once and for all, the Jacobian pair $(F,G)$ satisfying \eqref{wheraraaa} and Lemma \ref{Sect3-2}.

The main purpose of this section is  to prove that the projection $\pi_1:V\to\C^2$ defined in \eqref{proj1} is in fact not surjective (which  contradicts Theorem \ref{Theo-3} and thus proves Theorem \ref{MAINT}). To help understanding, we first give some explanations.

Throughout the whole section, (though $V$ may have some other properties) we only need to use  the following three properties (C1)--(C3) satisfied by $V$ [especially, (C2) will be frequently used], which do not necessarily depend on the fact that $V$ arises from the Jacobian problem. As in \eqref{V=0}, an element in $V$ is usually denoted as $(p_1,p_2)=\big((x_1,y_1),(x_2,y_2)\big)$.
\begin{itemize}
\item[(C1)] $V$ is a nonempty closed subset of ${\C^4}$ (by Lemma \ref{Procesi-1});
\item[(C2)] $|y_1|+|y_2|=o(h_{p_1,p_2})$ when $h_{p_1,p_2}:=|x_1|+|y_1|+|x_2|+|y_2|\to\infty$ (by Theorem \ref{Theo-2});
\item[(C3)]
for each $(\tilde p_1,\tilde p_2)\in V$,  there exists a neighborhood ${\mathcal O}_{\tilde p_1,\tilde p_2}\subset\C^4$ of $(\tilde p_1,\tilde p_2)$ such that for all $(p_1,p_2)=((x_1,y_1),(x_2,y_2))\in V\cap {\mathcal O}_{\tilde p_1,\tilde p_2}$, $x_1,y_1$ are holomorphic functions of $x_2,y_2$ (by the local bijectivity of Keller maps).
\end{itemize}
\noindent Note that the local bijectivity of Keller maps also implies that $V$ satisfies,
\begin{itemize}
\item[(C4)]
there exists a non-constant polynomial $\theta(x_1,y_1,x_2,y_2)$ such that $V\not\subset S_\theta$, where $S_\theta$ is the set of zeros of the polynomial $\theta$,
$$\phantom{mama}S_\theta=\{((\tilde x_1,\tilde y_1),(\tilde x_2,\tilde y_2))\in\C^4\,|\,\theta(\tilde x_1,\tilde y_1,\tilde x_2,\tilde y_2)=0\},$$
 and further, for any $(\tilde p_1,\tilde p_2)\in V\bs S_\theta$, there is some neighborhood ${\mathcal O}_{\tilde p_1,\tilde  p_2}\subset\C^4$ of $(\tilde p_1,\tilde  p_2)$ such that
for all $(p_1,p_2)=((x_1,y_1),(x_2,y_2))\in V\cap {\mathcal O}_{\tilde p_1,\tilde p_2}$,  $x_1,x_2$ are holomorphic functions of  $y_1,y_2$.
\end{itemize}
\begin{rema}\label{Easy-P-rem}\rm
We would like to mention that if we change (C4) to the following slightly stronger form (C4)$'$, then the proof (attached in Appendix \ref{sect8}) becomes very easy by using some geometry [since we already have (C4), if we change (C4)$'$ to that the Jacobian matrix of $\pi_1$ at any $(\tilde p_1,\tilde p_2)\in V\cap S_\theta$ is of maximal rank, then we also get that $\pi_1$ is not surjective]:
\begin{itemize}
\item[(C4)$'$]
 for any $(\tilde p_1,\tilde p_2)\in V$, there is some neighborhood ${\mathcal O}_{\tilde p_1,\tilde  p_2}\subset\C^4$ of $(\tilde p_1,\tilde  p_2)$ such that
for all $(p_1,p_2)=((x_1,y_1),(x_2,y_2))\in V\cap {\mathcal O}_{\tilde p_1,\tilde  p_2}$,  $x_1,x_2$ are holomorphic functions of  $y_1,y_2$.
\end{itemize}
\end{rema}

Thus (C4) means that (C4)$'$ holds for all elements in $V$ except the zeros of $\theta$.
Because of Remark \ref{Easy-P-rem}, we do think that it is very reasonable to have that $\pi_1$ is not surjective. However since we only have (C4) instead of (C4)$'$, we cannot treat $V$ in a uniform way. This is why the proof in this section looks so complicated and why we have to consider case by case. As a matter of fact, in this section,  we will try to avoid choosing elements of $V$ which are in $S_\theta$ [for example the element $(\bar p_1,\bar p_2)\in A_{\kk,\kk}$ satisfying \eqref{TaKa} is an element not in $S_\theta$ when $\kk\gg1$].

\subsection{Two propositions}\label{Sect3.1}

We want to build a proposition, namely, Proposition \ref{real00-inj}, which will immediately imply another proposition, namely, Proposition \ref{real00-inj+1}, that will give us a contradiction (thus proving Theorem \ref{MAINT}).

 First we need to introduce some notations.
For any
 $k_{\ZeRo},k_{\OnE}\in\R_{\ge0}$, by Theorem \ref{Theo-3}, 
 the following is a nonempty compact subset of $V$,
\begin{eqnarray}
\label{Ak=}&\!\!\!\!\!\!\!\!\!\!\!\!\!\!&
A_{k_{\ZeRo},k_{\OnE}}=\big\{(p_{\ZeRo},p_{\OnE})=\big((x_{\ZeRo},y_{\ZeRo}),(x_{\OnE},y_{\OnE})\big)\in V\ \big|\ |x_{\ZeRo}|=k_{\ZeRo},\,|x_{\OnE}|=k_{\OnE}\big\}\ne\emptyset.
\end{eqnarray}
Thus
the following is a well-defined 
function of $k_{\ZeRo},k_{\OnE}\in\R_{\ge0}$,
\begin{eqnarray}
\label{Ak=1}&\!\!\!\!\!\!\!\!\!\!\!\!\!\!\!\!\!\!\!\!&
\g_{k_{\ZeRo},k_{\OnE}}=\max \,\big\{|x_{\OnE}+y_{\OnE}|\ \big|\ (p_{\ZeRo},p_{\OnE})=\big((x_{\ZeRo},y_{\ZeRo}),(x_{\OnE},y_{\OnE})\big)\in A_{k_{\ZeRo},k_{\OnE}}\big\}.
\end{eqnarray}

We fix some choices of positive numbers
satisfying,
\begin{eqnarray}
\label{MSmde33333}
\!\!\!\!\!\!\!\!\!\!\!\!\!\!\!\!\!\!\!\!\!\!\!\!\!\!\!\!\!\!\!\!\!&&
1
{\ssc\,}\ll{\ssc\,}\ell_{\rZeRo}:=\d_{\rZeRo}^{-1}
{\ssc\,}\ll{\ssc\,}\ell_{\rOnE}:=\d_{\rOnE}^{-1}
{\ssc\,}\ll{\ssc\,}\ell_2:=\d_2^{-1}
{\ssc\,}\ll{\ssc\,}\ell:=\d^{-1}
{\ssc\,}\ll{\ssc\,}\kk
{\ssc\,}\ll{\ssc\,}\ep^{-1}
,
\ \ \ \ \
\ell_0
\in\Z_{>0}
.\!\!\!\!\!\!\!\!\!
\end{eqnarray}
\begin{rema}\label{rema3.1}\rm \begin{itemize}\item[(i)]
The choices may depend on the situation we encounter, and further, the choice of a number listed later may depend on all choices of numbers listed earlier; for instance, we may require that $\ell_1\gg\ell_0^{\ell_0},\,\ell_2\gg\ell_1^{\ell_1}
$, etc%
. We always require $\kk$ to be sufficiently larger than $\SS_7$,
 with $\SS_7$ being the number appeared in Proposition \ref{Sect3-Lemm5}:
 \equa{kk-bbb}{\dis\kk\gg\SS_7.}
\item[(ii)] We need to use the following convention: we can regard elements in \eqref{MSmde33333} as parameters to be fixed  upon our requirement in the situation we encounter. Sometimes we need to compare some number in \eqref{MSmde33333} with other numbers, in this case, for convenience we may regard some element in  \eqref{MSmde33333} as a variable;
     for instance, if we regard $\ell$ as a variable, then we need to
regard $\ell_0,
\ell_1,
\ell_2
    $ as  fixed numbers  and $\kk,%
\ep
$ as
%
variables
%
such that we have [cf.~\eqref{MSM1111}$\ssc\,$],
\equa{MSmde33333+}{\dis
\!\!\!\!\!\!1
{\ssc\,}\sim_{\ell\ssc\,}{\ssc\,}\ell_{\rZeRo}
{\ssc\,}\sim_{\ell\ssc\,}{\ssc\,}\ell_{\rOnE}
{\ssc\,}\sim_{\ell\ssc\,}{\ssc\,}\ell_2
{\ssc\,}\prec_{\ell\ssc\,}{\ssc\,}\ell
{\ssc\,}\prec_{\ell\ssc\,}{\ssc\,}\kk
{\ssc\,}\prec_{\ell\ssc\,}{\ssc\,}\ep^{-1}
.\!\!\!\!\!}
\end{itemize}\end{rema}

Whenever the parameter $\kk\gg\ell$ is chosen, we fix 
an element $(\bar p_{\ZeRo},\bar p_{\OnE})=\big((\bar x_{\ZeRo},\bar y_{\ZeRo}),(\bar x_{\OnE},\bar y_{\OnE})\big)\in A_{\kk,\kk}$ satisfying
 \begin{eqnarray}\label{TaKa}
&\!\!\!\!\!\!\!\!\!\!\!\!\!\!&
|\bar x_{\ZeRo}|=
|\bar x_{\OnE}|=\kk
,\ \ \ \ \
|\bar x_{\OnE}+\bar y_{\OnE}|=\g_{\kk,\kk}.\end{eqnarray}

\begin{nota}\label{notaXY}\rm\begin{itemize}\item[(i)]
For any $x_{\ZeRo},x_2,y_{\ZeRo},y_{\OnE}\in \C$, we  denote,
\begin{eqnarray}
\label{SimMMSMS}
&\!\!\!\!\!\!\!\!\!\!\!\!\!\!\!\!\!\!\!\!\!\!\!\!\!\!\!\!\!&{\rm(a)\ }
 X_{\ZeRo}=\frac{x_{\ZeRo}}{\bar x_{\ZeRo}},
\ \ \ \ \ {\rm(b)\ }
X_{\OnE}=\frac{x_{\OnE}}{\bar x_{\OnE}}, \ \  \ \ \
{\rm(c)\ }
Z=\frac{x_{\OnE}+y_{\OnE}}{\bar x_{\OnE}+\bar y_{\OnE}}
.
\end{eqnarray}
\item[(ii)]
For any subset of $\C^4$, we always use the same symbol with an overbar to denote its closure in $\C^4$, that is, if $S$ is a subset of $\C^4$,
then $\ol S$ denotes its closure in $\C^4$.
\end{itemize}
\end{nota}

The reason we 
define
\eqref{SimMMSMS} 
is in order that  when $(p_1,p_2)=\big((x_1,y_1),(x_2,y_2)\big)$ is close to $(\bar p_1,\bar p_2)$ we have that $
X_1,
X_2,
Z$ are close to $1$, which will make our arguments easier. 
%

{
\begin{defi}\rm\label{def-S0-S1}
\NOUSE
{We define the following Laurent polynomials
of $x_1,x_2$ $($or $X_1,X_2{\ssc\,})$,
\equa{Poly-x1-x2}{\dis
{\rm(i)\ }A_1=\frac{X_2^{\ell_0(1+\d_0)}}{X_1^{\ell_0^2}}\Big(\d_0 + (1-\d_0)X_1^{\ell_0^2(1+\d_0)}\Big),\ \ \
{\rm(ii)\ }H_0=\ell_0^2+1-\ell_0^2X_1^{\ell_0}X_2^{\ell_0(1+\d_0)}.
}
}%
Denote by $S_1$  the  open subset of $\C^2$ %
consisting of all elements
$(x_1,x_2)$ %
whose coordinates  satisfy
\begin{eqnarray}
&\!\!\!\!\!\!\!\!\!\!\!\!\!\!\!\!\!\!\!\!\!\!\!\!\!\!\!&
\label{meme}
{\rm(i)\ }
\d_2^{4}<|X_1|<\ell_2^{4}
,\ \
\ \ \ \ \
{\rm(ii)\ }
\d_2^{4}<|X_2|<\ell_2^{4}
,\ \ \ \ \ \ \ {\rm(iii)\ }\Big|2- \frac1{X_1}\Big|>\d_2^4
.
\end{eqnarray}
Set $S_2=\pi_1^{-1}(S_1)$
.
\end{defi}
Note that 
$(X_1,X_2)=(1,1)$ satisfies \eqref{meme}, thus by  \eqref{TaKa},\,\eqref{SimMMSMS},
 $(\bar p_1,\bar p_2)\in S_2$, i.e., $S_2$ is non-empty. Moreover, 
$S_2$ is open in $V$ and is a bounded set by Proposition \ref{pro-also}.

The reason we define $S_2$ is to ensure that 
$A_1,A_2
$ 
to be defined in Definition \ref{ABC}
~are well-defined 
functions in $\ol S_2$%
.
{The importance of condition \eqref{meme}
~is
that we can obtain the following lemma, esp., \eqref{equa-Case6-lemm}\,(iv), which is very crucial in obtaining Proposition \ref{real00-inj} in case $V_0=V_2$.
\begin{lemm}\label{Case6-lemm}
When $(p_1,p_2)\in \ol S_2\subset V$
, we have
\begin{eqnarray}
\label{equa-Case6-lemm}
\!\!\!\!\!\!\!\!\!\!\!\!\!\!\!\!\!\!\!\!\!\!\!\!\!\!\!&&
{\rm(i)\ }\d_2^4\kk\le|x_1|,|x_2|\le\ell_2^4\kk,\ \ \ \ \ \ \ \
{\rm(ii)\ }2\d_2^4\kk\le h_{p_1,p_2}\le2(\ell_2^4+\d^3)\kk,
\nonumber\\ \!\!\!\!\!\!\!\!\!\!\!\!\!\!\!\!\!\!\!\!\!\!\!\!\!\!\!&&
{\rm(iii)\ }x_2+y_2=x_2\Big(1+O(\d)^3\Big),\ \ \ \ \ \ \ \
 {\rm(iv)\ }Z{\ssc}={\ssc}X_2\Big(1+O(\d)^3\Big)=X_2{\ssc}+{\ssc}O(\d)^3,
\nonumber\\ \!\!\!\!\!\!\!\!\!\!\!\!\!\!\!\!\!\!\!\!\!\!\!\!\!\!\!&&
 {\rm(v)\ }(\d_2^4{\ssc}-
{\ssc}\d^2)\kk{\ssc}\le{\ssc}|x_2{\ssc}+{\ssc}y_2|
{\ssc}\le{\ssc}(\ell_2^4{\ssc}+{\ssc}\d^2)\kk,\ \ \ \ \ \ \  {\rm(vi)\ }\d_2^4{\ssc}-{\ssc}\d^2
{\ssc}\le{\ssc}|Z|{\ssc}\le{\ssc}\ell_2^4{\ssc}+{\ssc}\d^2.
\end{eqnarray}
\end{lemm}
\noindent{\it Proof.~}Formula \eqref{equa-Case6-lemm}\,(i) follows from \eqref{TaKa}--\eqref{meme}.
We have
\equa{msmsem999999}{\dis
2\d_2^4\kk\stackrel{{}^{\sc\rm\eqref{equa-Case6-lemm}\,(i)}}{\le}|x_1|+|x_2|\stackrel{{}^{\sc\rm\eqref{dp0p1}}}{\le}h_{p_1,p_2}.
}
From this we obtain
\equa{h-p-1-p-2===}{\dis
h_{_{\sc p_1,p_2}}^{^{\sc-\frac1{m+1}}}\stackrel{{}^{\sc\rm\eqref{msmsem999999}}}{\le}(2\d_2^4\kk)^{-\frac1{m+1}}
\stackrel{{}^{\sc\rm\eqref{MSmde33333}}}{<}\d^3.
}
Thus
\equa{Thendn---}{\dis
h_{p_1,p_2}\stackrel{{}^{\sc\rm\eqref{dp0p1}}}{=}|x_1|+|x_2|+|y_1|+|y_2|
\stackrel{{}^{\sc\rm\eqref{equa-Case6-lemm}\,(i),\,\eqref{mqp1234-2}}}{\le}
2\ell_2^4\kk+2h_{_{\sc p_1,p_2}}^{^{\sc\frac{m}{m+1}}}\stackrel{{}^{\sc\rm\eqref{h-p-1-p-2===}}}{<}2\ell_2^4\kk+2\d^3 h_{p_1,p_2}.}
This with \eqref{msmsem999999} proves \eqref{equa-Case6-lemm}\,(ii). We have
\begin{eqnarray}
\label{equa-Case6-lemm+1}
\!\!\!\!\!\!\!\!\!\!\!\!\!\!\!\!\!\!\!\!\!\!\!\!\!\!\!&&
{\rm(i)\ }
\Big|\frac{y_2}{x_2}\Big|\stackrel{{}^{\sc\rm\eqref{mqp1234-2},\,\eqref{equa-Case6-lemm}\,(i)}}{ \le}
\frac{h_{_{\sc p_1,p_2}}^{^{\sc\frac m{m+1}}}}{\d_2^4\kk}\stackrel{{}^{\sc\rm\eqref{equa-Case6-lemm}\,(ii),\,\eqref{MSmde33333}}}{<}\d^3,
\nonumber\\ \!\!\!\!\!\!\!\!\!\!\!\!\!\!\!\!\!\!\!\!\!\!\!\!\!\!\!&&
{\rm(ii)\ }x_2+y_2=x_2\Big(1+\frac{y_2}{x_2}\Big)\stackrel{{}^{\sc\rm\eqref{equa-Case6-lemm+1}\,(i)}}{=}
x_2\Big(1+O(\d)^3\Big),
\nonumber\\ \!\!\!\!\!\!\!\!\!\!\!\!\!\!\!\!\!\!\!\!\!\!\!\!\!\!\!&&
{\rm(iii)\ }Z\stackrel{{}^{\sc\rm\eqref{SimMMSMS}\,(c)}}{=}\frac{x_2+y_2}{\bar x_2+\bar y_2}
\stackrel{{}^{\sc\rm\eqref{equa-Case6-lemm+1}\,(ii)}}{=}\frac{x_2\Big(1+O(\d)^3\Big)}{\bar x_2\Big(1+O(\d)^3\Big)}
\stackrel{{}^{\sc\rm\eqref{SimMMSMS}\,(b)}}{=}X_2\Big(1+O(\d)^3\Big)\stackrel{{}^{\sc\rm\eqref{meme}\,(i)}}{=}X_2+O(\d)^3.
\end{eqnarray}
From the above we obtain \eqref{equa-Case6-lemm}\,(iii),\,(iv). Then \eqref{equa-Case6-lemm}\,(v) follows from \eqref{equa-Case6-lemm}\,(i),\,(iii). Similarly by \eqref{meme}\,(i),\,\eqref{equa-Case6-lemm}\,(iv), we have \eqref{equa-Case6-lemm}\,(vi).
\hfill$\Box$
}}%
\begin{defi}\rm\label{ABC}
{
Let $
b_\kk>0
$ be the
number
to be defined in \eqref{@suc2hthat=4}.
\NOUSE{Let $\lL$  be the small integer bigger than
$a_\kk+b_\kk$, i.e.,
\equa{ll-dd}{\dis
\lL=a_\kk+b_\kk+\l_0\in\Z_{>0}\mbox{ \ for some \ $\l_0\in\R_{>0}$ \ with \ }\l_0\le1,\ \ \ \mbox{ and denote \ }\ \dD=\lL^{-1}.}
We will see in \eqref{Akk-bkk} that $\lL\gg\ell$.
}%
\NOUSE{For $i=1,2,...,5$, let $\l_i\in\R_{\ge}0$ with $\l_i<1$ such that the following are all positive rational numbers,
\equa{ration}{\dis\!\!\!\!\!\!\lL^{6-\d_1}{\sc\!}-{\sc\!}\l_1,\,
\frac{\lL^3(1-\dD)}{2}{\sc\!}-{\sc\!}\l_1\dD^{\lL}, \,
\frac{3\lL^8}{2}{\sc\!}+{\sc\!}\l_3\dD^{\lL},\,
\frac{\lL^8}{2}{\sc\!}-{\sc\!}\l_4\dD^{\lL},\,
\dD^5{\sc\!}-{\sc\!}\l_5\dD^{\lL}{\sc\!}\in{\sc\!}\Q_{>0},\mbox{   and   }
\lL^{6-\d_1}{\sc\!}-{\sc\!}\l_1{\sc\!}\in{\sc\!}\Z_{>0}.\!\!\!\!
}%
}%
\NOUSE{Let $\lL$  be the smallest integer bigger than $10b_\kk-9a_\kk$, i.e.,
\equa{In-lLL}{\dis
\lL=10b_\kk-9a_\kk+\l_0\mbox{ with $\l_0\in\R$ and $0<\l_0\le1$, and denote }\dD=\lL^{-1}.
}
We d}Denote
\begin{eqnarray}
\label{Tx-B0}\label{tX1==}
&\!\!\!\!\!\!\!\!\!\!\!\!\!\!\!\!\!\!\!\!\!\!\!&
\widetilde X_1=
\a_0X_1
,\  \ \ \
\a_0=1+\b_0\ep
,\ \ \ \ \ \
\b_0=
-1+\d^3b_\kk
\stackrel{{}^{\sc\rm
\eqref{Akk-bkk}}}{>}0
.
\end{eqnarray}
%
}\NOUSE
{
Denote, 
\equa{tX1==}{\dis \widetilde X_1=\a_0X_1,\ \ \ \ \a_0=1+\b_0\ep,\ \ \ \ \b_0=1-(a_\kk+\b_\kk)\dD
\stackrel{{}^{\sc\rm \eqref{ll-dd}}}{=}\l_1\dD>0.}
 \NOUSE{and denote%
,
\equa{tX11111}{\dis
{\rm(i)\ }
\widetilde X_1=(\a_0X_1)^{\ell_0}
,\ \ \ \ \ \
{\rm(ii)\ }\widetilde X_2=X_2^{\ell_0}
,\ \ \ \ \ \
{\rm(iii)\ }\widetilde Z=Z^{\ell_0}
.
}
D and }}%
For any $(p_1,p_2)\in \ol S_2$, we define%
~[noting that \eqref{meme}\,(iii) with \eqref{tX1==}
shows that $A_1$ is invertible
]%
,
%
%
\NOUSE{ (we will see that all $\l_i$ do not affect our arguments, thus one can regard them as zero; further
without loss of generality, one can simply regard $\b_1,\b_2,\b_3$ as $1$),
}\NOUSE{ [note that the following only contain integral powers of $X_1,X_2,Z$ by defining $\widetilde X_1^{\d_0},\widetilde
X_2^{\d_0},\widetilde Z^{\d_0}$ to be $\a_0X_1,X_2,Z
$ respectively]%
}
\begin{eqnarray}
\label{LetNSoOP----1}
&\!\!\!\!\!\!\!\!\!\!\!\!\!\!\!\!\!\!\!\!\!\!\!\!\!\!\!\!\!\!\!\!\!\!\!\!\!\!\!\!\!\!\!\!\!\!\!\!\!\!\!
&
{\rm(i)\ }A_1
=
X_2^{\ell_0} \Big(2-\frac1{\widetilde X_1}\Big),
\ \ \ \ \ \
{\rm(ii)\ }A_2
=
\frac{X_2^{\ell_0+1}}{Z}\Big (2-\frac{\widetilde  X_1}{A_1^2}\Big) 
.
\!\!\!\!\!\!\!\!\!\!\!\!\!\!\!\!\!\!\!\!\!\!\!\!
\end{eqnarray}
\end{defi}
\NOUSE
{%
We refer readers to Remark \ref{rema-def-A1} and \eqref{B1-C1-rema} to see why we define
$A_1,A_2,A_3,B_3
$ in such a way. Here we first give the following remark.
\begin{rema}\label{B1-C1}\rm
By \eqref{LetNSoOP----1}\,(iii), we obtain
\equa{MMSMSD000000000000}{\dis{\rm(i)\ }
\frac43=\frac{\widetilde X_1^6}{3}+\frac1{A_3\widetilde X_1^3}
=\frac{C_3}{3B_3}+B_3
,\ \ {\rm(ii)\ }C_3=B_3\widetilde X_1^6=\frac{\widetilde X_1^3}{A_3}.}
\NOUSE{Thus the following are well-defined in $\ol V_2$,
\begin{eqnarray}
\label{LetNSoOP----1+CB}
&\!\!\!\!\!\!\!\!\!\!\!\!\!\!\!\!\!\!\!\!\!\!\!\!\!\!\!\!\!\!\!\!\!\!\!\!\!\!\!\!\!\!\!\!\!\!\!\!\!\!\!
&
{\rm(i)\ }B_1=\frac{ A_3}{A_1\widetilde X_1^2},\ \ \ \
{\rm(ii)\ }C_1=\frac{A_3^{\ell_0^2+1}}{A_1 \widetilde X_1^3},
\ \ \ \
{\rm(iii)\ }B_2=
\frac{A_2 A_3^3}{A_1\widetilde X_1^3 X_2},\ \ \ \
{\rm(iv)\ } C_2=\frac{A_2^2 A_3^5}{A_1\widetilde X_1^4 X_2^2}.
\!\!\!\!\!\!\!\!\!\!\!\!\!\!\!
\end{eqnarray}
}%
\NOUSE{which 
can be rewritten as
\begin{eqnarray}
\label{More-de-1}
&\!\!\!\!\!\!\!\!\!\!\!\!\!\!\!\!\!\!\!\!\!\!\!\!\!\!\!
&{\rm(i)\ }\frac{16}{15}
=\frac{\widetilde X_1^{172}}{15}+\frac{X_2^{11}}{A_1 A_2 A_3\widetilde X_1^{14} Z^{11}}
=\frac{C_1}{15B_1}+B_1,\ \ \ {\rm(ii)\ }C_1=B_1\widetilde X_1^{172}=
\frac{\widetilde X_1^{158} X_2^{11}}{A_1 A_2 A_3 Z^{11}}.
\!\!\!\!
\end{eqnarray}
}Multiplying \eqref{LetNSoOP----1+CB}\,(i) 
by $3 B_1$
,
we immediately obtain%
,
\begin{eqnarray}
\label{Re-Writtt}
&\!\!\!\!\!\!\!\!\!\!\!\!\!\!\!\!\!\!\!\!\!\!\!\!\!\!\!\!\!\!&
3B_3^2-4B_3+C_3=0.
\end{eqnarray}
\NOUSE{
Similarly, \eqref{LetNSoOP----1}\,(iii) can be rewritten as
\equa{andthooo}{\dis
\frac{5B_2^2}{2}-\frac{7B_2}{2}+C_2=0 \ \mbox{with  }
B_2=\frac{A_2 A_3^{\ell_0}\widetilde X_1^{876} Z^{64}}{X_2^{438}}, \ \
C_2=\frac{A_2\widetilde X_1^{886} Z^{64}}{A_3^{\ell_0^5 +\ell_0}X_2^{438}}.
}
}\NOUSE{
Similarly, by \eqref{LetNSoOP----1}\,(iii), we have
\begin{eqnarray}
\label{More-de-1+2}
&\!\!\!\!\!\!\!\!\!\!\!\!\!\!\!\!\!\!\!\!
1&\ \, \
\stackrel{{}^{\sc\rm\eqref{LetNSoOP----1}\,(iii)}}{=}
\frac{\widetilde X_1 X_2}{A_2A_3^2}\Big (\frac15+\frac{4A_1 \widetilde X_1^2}{5A_3}\Big)
=
\frac{A_2 A_3^3}{A_1\widetilde X_1^3 X_2}
\frac{A_1\widetilde X_1^4 X_2^2}{A_2^2 A_3^5}
 \mbox{\Large$\Big($}
\frac15 +\frac45
\frac{A_2^2 A_3^5}{A_1\widetilde X_1^4 X_2^2}
\Big(
\frac{A_1\widetilde X_1^3 X_2}{A_2 A_3^3}
\Big)^2
\mbox{\Large$\Big)$}
\!\!\!\!\!\!\!\!\!\!\!\!
  \nonumber\\
 &\!\!\!\!\!\!\!\!\!\!\!\!\!\!\!\!\!\!\!\!\!\!\!\!\!\!\!\!\!\!\!\!&
\stackrel{{}^{\sc\rm\eqref{LetNSoOP----1+CB}\,(iii),\,(iv)}}{=}
\ \frac{B_2}{C_2}
\big(\frac15+\frac{4 C_2}{5B_2^2}\big)=\frac{B_2}{5C_2}+\frac{4}{5B_2},
\end{eqnarray}
and so
\begin{eqnarray}
\label{Re-Writtt+2}
&\!\!\!\!\!\!\!\!\!\!\!\!\!\!\!\!\!\!\!\!\!\!\!\!\!\!\!\!\!\!&
4B_2^{-2}-5B_2^{-1}+C_2^{-1}=0.
\end{eqnarray}
%
}\NOUSE{Regarding the above as a quadratic equation on $B_1^{-1}$ and solving it, we obtain
the following [later on we will prove that another solution does not satisfy condition \eqref{C+ToSayas+1}\,(c)$\ssc\,$].
\begin{eqnarray}
\label{1=bbb-remark}
&\!\!\!\!\!\!\!\!\!\!\!\!\!\!\!\!\!\!\!\!\!\!\!\!\!\!\!\!\!\!\!\!\!\!\!\!\!&
B_1^{-1}=\frac{5}{8}\mbox{\Large$\Big($}1+\Big(1-\frac{16}{25C_1}\Big)^{\frac12}
\mbox{\Large$\Big)$},\!\!\!\!\!
\end{eqnarray}
where \equa{Convvv}{\dis
\Big(1-\frac{16}{25C_1}\Big)^{\frac12}
=
1+
\mbox{$\sum\limits_{i=1}^\infty$}\binom{\frac12}{i}\Big(-\frac{16}{25C_1}
\Big)^i,} is the unique element defined by \eqref{bimeformo},
which converges absolutely  by the fact from \eqref{C+ToSayas+1}\,(c) that
$\frac{16}{25|C_1|}\le\frac{16(1+\eE_1^2)}{25}<1$.
We remark that by using purely algebraic methods [without any geometry but a little analysis that the series \eqref{Convvv} converges absolutely], we will in fact prove that
formula \eqref{1=bbb-remark} holds for any element in $W_2$ [which is the set to be defined in \eqref{C+ToSayas+1}$\ssc\,$].
The importance of  \eqref{1=bbb-remark} is that
we can use it to prove that condition \eqref{ToSayas}\,(c) holds for any element in $\ol V_2$
[cf.~\eqref{ImMpP}\,(2)$\ssc\,$].
}%
\NOUSE{
First observe from \eqref{LetNSoOP----1}\,(ii),\,(iii) that we have the following relation between $A_1$ and $A_2$,
\equa{Related-A1-A2}{\dis
A_2=
\frac{X_2^{ 2\ell_0}}{A_3^{\ell_0}\widetilde X_1^{3\ell_0} Z^{\ell_0}}
\Big (\frac75-\frac{2A_1}{5A_3^{50\ell_0^2}\widetilde X_1^{10\ell_0}}
\Big)
,
}
which can be rewritten as
}
Thanks are due to Professor Claudio Procesi, who observes that formulas \eqref{Re-Writtt}
~can be  used to give an elegent proof of
Lemma \ref{procesi-lemm3}, which in particular implies that condition
\eqref{ToSayas+1}\,(c) is satisfied by elements in $\ol V_2$%
.
\end{rema}
}%
\NOUSE
{%
We need some explanations about definition \eqref{LetNSoOP----1} below as we need any number that appears as a power of a complex number to be integral [however we remind that using Convention \ref{conv1}\,(1), powers that appear outside any absolute value symbols
(or appear inside the symbols but there is only one term inside the symbols) do not need to be integral (but need to be in $\R$)$\ssc\,$]. First note that in the following we always expand a term $(a^{\hat a} b^{\hat b} c^{\hat c})^h$ into the following,
\equa{abc-d}{\dis (a^{\hat a} b^{\hat b} c^{\hat c})^h=a^{\hat ah}b^{\hat bh}c^{\hat ch}\mbox{ \ for \ }a,b,c\in\C,\ \hat a,\hat b,\hat c\in\R,\ h\in\Z.}
Recall that we have required $\ell_0=
\d_0^{-1}$ to be a positive integer.
First
if we choose $\b_1,\b_2$ to be integral, then all powers that appear in  \eqref{LetNSoOP----1} and
\eqref{ToSayas+1}\,(e)--(g) are rational.
Now we can take $\b_1=\ell_0^4$, then after applying the above formula we see that $B_1$ defined below
is an element in $\C[\tildeX_1^{\pm1},X_2^{\pm1},Z^{\pm1}]$. Next we take $\b_2>\b_1$ to be a common multiple of denominators of all rational numbers (except $\b_2$ itself)
that appear as powers  in any expression of \eqref{LetNSoOP----1} and \eqref{ToSayas+1}\,(e)--(g),
then after applying formula \eqref{abc-d} we see that all powers of $\tildeX_1,X_2,Z,B_1,B_2,B_3$ become integers.
Since $\b_1,\b_2$ only depend on $\ell_0$, we can always assume $\ell_1$ in \eqref{MSmde33333} also satisfies that $\ell_1\gg\b_2$. In this way, we can give the following definition.
\begin{lemm}\rm\label{defi-B123}
The above are well-defined 
rational functions
on $\ol S_2$
$($the closure of $S_2$ in $\C^4)$.
Further, 
we have
\begin{eqnarray}
\label{LetNSoOP----2}
&\!\!\!\!\!\!\!\!\!\!\!\!\!\!\!\!\!\!\!\!\!\!\!\!\!\!\!\!\!\!\!\!\!\!\!\!\!\!
&
{\rm(i)\ }A_1=\frac{\ell_0+1- \ell_0\frac{ X_1^{\ell_0^3(1-\d_0)}}{ X_2^{\ell_0^3}}}{X_1^{\ell_0^3}}+O(\d)^2
=\frac{\ell_0+1-\ell_0\frac{X_1^{\ell_0^7+\ell_0^3-\ell_0^2}}{A_2^{\ell_0^3}}}{X_1^{\ell_0^3}}+O(\d)^2,
\!\!\!\!
\\\nonumber
&\!\!\!\!\!\!\!\!\!\!\!\!\!\!\!\!\!\!\!\!\!\!\!\!\!
&{\rm(ii)\ }A_2=\frac1{A_2^{\ell_0^3}\widetilde X_1^{\ell_0^2(1+\d_0^2)}(1+\d_0-\d_0A_1\widetilde X_1^{\ell_0^3})}
=
\frac{X_2^{\ell_0^3}}{A_2^{\ell_0^3}X_1^{\ell_0^3+1}}+O(\d)^2=\frac1{X_1^{\ell_0^7+\ell_0^3+1}}+O(\d)^2,
\!\!\!\!
\\\nonumber
&\!\!\!\!\!\!\!\!\!\!\!\!\!\!\!\!\!\!\!\!\!\!\!\!\!
&{\rm(iii)\ }A_2=X_1^{\ell_0^4}X_2+O(\d)^2
.\!\!\!\!\!\!\!\!\!\!\!\!\!
\end{eqnarray}
\end{lemm}
\noindent{\it Proof.~}Let $(p_1,p_2)\in\ol S_2$. By \eqref{meme}(i),\,(ii) and \eqref{equa-Case6-lemm},
$X_1,X_2,Z\ne0$, we see that $A_1,A_2,A_3$ are well-defined.
By \eqref{equa-Case6-lemm}\,(iv) and  applying formula \eqref{Menenen343333} [noting
that $\widetilde X_1=X_1+O(\d)^2$ as $\a_0$ is a $1+O(\ep)^1$ element], we obtain
\eqref{LetNSoOP----2}, where, we have made use of \eqref{LetNSoOP----2}\,(iii) in the last equalities of \eqref{LetNSoOP----2}\,(i),\,(ii), and  the first equality of
\eqref{LetNSoOP----2}\,(ii)
 is obtained by simply substituting $A_1$ by
\eqref{LetNSoOP----1}\,(i).
%
%
In general when conditions \eqref{meme} and \eqref{C+ToSayas+1} hold, we need frequently use formula \eqref{Menenen343333}. Denote
\begin{eqnarray}
\label{Deno-T0}
\!\!\!\!\!\!\!\!\!\!\!\!\!\!\!\!\!\!\!\!\!\!\!\!\!\!\!\!\!\!&&
T_0=\{\mbox{$A_1,A_2,A_3$ and all factors appearing in \eqref{LetNSoOP----1},\,\eqref{LetNSoOP----2}\,}
\NOUSE{X_1,\,\widetilde X_1,\,X_2,\,Z,\,A_1,\,A_2,\,A_3,\,2 - \frac{X_2^{\ell0}}{\widetilde X_1Z^{\ell_0}},\,2 - \frac{1}{X_1},\,
\frac32-\d_0^5-\Big(\frac12-\d_0^5\Big)\frac{X_2^{\ell_0}}{\widetilde X_1 Z^{\ell_0}},\,\frac34+\frac{X_2^{\ell_0}}{4\widetilde X_1Z^{\ell_0}},\!\!\!\!\!\!\!\!\!\!\!\!\!\!\!\!\!\!\!\!\!\!\!\nonumber\\
\!\!\!\!\!\!\!\!\!\!\!\!\!\!\!\!\!\!\!\!\!\!\!\!\!\!\!\!\!\!&&
\phantom{T_0=\Big\{}\!\!\!
\frac12{\sc\!}+{\sc\!}\d_0^5{\sc\!}+{\sc\!}\Big(\frac12{\sc\!}-{\sc\!}\d_0^5\Big)\frac{\widetilde X_1^3}{A_1A_3^4X_2^2},\frac54-\frac{\widetilde X_1^3}{4A_1A_3^4X_2^2},\frac12{\sc\!}+{\sc\!}\d_0^5{\sc\!}+{\sc\!}
\Big(\frac12{\sc\!}-{\sc\!}\d_0^5\Big)\frac{X_1^{3\ell_0^2}}{A_1A_3^4X_2^2},
\frac54{\sc\!}-{\sc\!}\frac{X_1^{3\ell_0^2}}{4A_1A_3^4X_2^2}
}\}.
\end{eqnarray}
By
\eqref{meme},\,\eqref{equa-Case6-lemm},\,\eqref{C+ToSayas+1}, one can  observe
[cf.~\eqref{ImMpP}$\ssc\,$] that
$|a|^{\pm1}<\ell_2^{\ell_1}$ for all $a\in T_0$. Since $0<\d\ll\d_2^{\ell_1}$ by Remark \ref{rema3.1}\,(i), we have
\equa{Menenen343333}{\dis
\d a=O(\d)^1,\ \ \ a^{0+O(\d)^1}=O(\d)^1\mbox{ \ \ for \ \ }a\in T_0\mbox{ \ or \ }a^{-1}\in T_0.
}
\NOUSE{
From this and \eqref{LetNSoOP----1}, we obtain \eqref{LetNSoOP----2}, where the first equality of
\eqref{LetNSoOP----2}\,(iii) follows from \eqref{LetNSoOP----2}\,(i), (ii) [in fact  our original
idea is to use \eqref{LetNSoOP----2}\,(iii) as the definition of $A_2$, which can be rewritten as the form in \eqref{LetNSoOP----1}\,(iii)
by using \eqref{LetNSoOP----1}\,(i)$\ssc\,$].
}
\hfill$\Box$
}%
\NOUSE
{%
We understand that the above definitions look quite strange, but we would like to give some remarks.
Firstly, using \eqref{LetNSoOP----1}\,(iii) in $A_4$, we see that the equation $A_4=0$ defines an implicit function $A_1$ on variables $\tildeX_1,X_2,Z$.
Secondly, we will see later on that the powers appearing  in \eqref{LetNSoOP----1} will play important roles in our later proofs, but at this point we do not need to know
why we define in such ways; in fact we do not need to use
the precise forms of $A_1,A_2,A_3$ in the proof of Proposition \ref{real00-inj+1}, instead we only need to keep in mind the only thing that they are locally
holomorphic functions and nonzero everywhere when $(p_1,p_2)$ is in some set we will define.
}%
\NOUSE
{%
\begin{rema}\rm\label{B123-rema}
\begin{itemize}
\item[(i)] When we prove Proposition \ref{real00-inj}, the first thing we need to do is to prove $V_0\ne\emptyset$, which is done by choosing some suitable element $(p_1,p_2)\in V_0$. We regard such an element as the ``initial stage''. We will see from our arguments in this section that the ``initial stage'' controls globally the growths of $|x_2|,\,|x_2+y_2|$ on the whole set $V_0$.
\item[(ii)]
We understand that  definition \eqref{LetNSoOP----1} looks quite strange, but
our reason to give such complicated definition is the following:
We require $A_1,A_2,A_3$ to satisfy Proposition \ref{real00-inj}. More precisely,
\begin{itemize}\item[(a)] Firstly, we will see in the proof of Lemma \ref{Prop(ii)Holds}, as will be explained in Remark \ref{Initial-Rema}, that we have to choose so complicated powers that appear in
\eqref{LetNSoOP----1}\,(i)--(iii), which will play key roles in the proof of Proposition \ref{real00-inj}\,(ii).
\item[(b)]Secondly, we need to prove
in Lemma \ref{V0NOT0} that $V_2\ne0$ by choosing the ``initial stage'' $(p_1,p_2)$, sufficiently close to $(\bar p_1,\bar p_2)$, satisfying \eqref{ToSayas+1},\,\eqref{ToSayas+2}, such that \mbox{all} $A_1,A_2,A_3,\tildeX_1, X_2,Z$ are  $1+O(\ep)^1$ elements of the following form (where $\ep>0$ is sufficiently small),
\begin{eqnarray}
&
&
{\rm(i)\ }\tildeX_1=1+\tilde s\ep+O(\ep)^2,\ \ \tilde s=2\ell_0^4,\ \ \ \ \ \ {\rm(ii)\ }X_2=1+u\ep,\ \ u=-2\ell_0^4+\ell_0+1+\d_0,\!\!\!\!\!\!\!\!\!\!\!\!\!\!\!\!\!\!\!\!\!\!\!\!
\nonumber\\&\!\!\!\!\!\!\!\!\!\!\!\!\!\!\!\!\!\!\!\!\!\!\!\!\!\!\!\!&
{\rm(iii)\ }Z=1+\tilde v\ep,\ \ \tilde v=u-\d_0,\ \ \ \ \ \ {\rm(iv)\ }A_1=1-c_1\ep+O(\ep)^2,\ \ c_1=2(\ell_0^4-2\ell_0)a_1^{-1},\!\!\!\!\!\!\!\!\!\!\!\!\!\!\!\!\!\!\!\!\!\!\!\!
\nonumber
\\
&\!\!\!\!\!\!\!\!\!\!\!\!\!\!\!\!\!\!\!\!\!\!\!\!\!\!\!\!\!\!\!\!\!\!\!\!\!\!\!\!\!\!\!\!\!\!\!\!\!\!\!\!&
{\rm(v)\ }A_2{\ssc}={\ssc}1{\ssc}-{\ssc}c_2\ep{\ssc}+{\ssc}O(\ep)^2,\ \ c_2{\ssc}={\ssc}(\ell_0^4{\ssc}-{\ssc}2\ell_0)(2{\ssc}
-{\ssc}2\d_0{\ssc}-{\ssc}\d_0^2)\b_2,\!\!\!\!\!\!\!\!\!\!\!\!\!\!\!\!\!\!\!\!\!\!\!\!
\nonumber
\\
&\!\!\!\!\!\!\!\!\!\!\!\!\!\!\!\!\!\!\!\!\!\!\!\!\!\!\!\!\!\!\!\!\!\!\!\!\!\!\!\!\!\!\!\!\!\!\!\!\!\!\!\!\!\!\!\!\!\!\!\!\!\!\!\!\!\!\!\!\!\!\!\!\!\!\!\!\!\!\!\!\!\!\!\!\!\!\!\!\!\!\!\!\!\!\!\!\!\!\!\!\!\!\!\!\!\!\!\!\!\!\!\!\!\!\!\!\!\!\!\!\!\!\!\!\!\!\!\!&
\label{B123X12Z}
{\rm(vi)\ }A_3{\ssc}={\ssc}1{\ssc}-{\ssc}\ep{\ssc}+{\ssc}O(\ep)^2,\end{eqnarray}
where \eqref{B123X12Z}\,(iv)--(vi) are obtained from \eqref{LetNSoOP----1}\,(i)--(iii) and \eqref{B123X12Z}\,(i)--(iii) (see the proof of Lemma \ref{V0NOT0} for detail).
\item[(c)]
Thirdly, in order to obtain \eqref{ToSayas+1}\,(e)--(h),\,\eqref{ToSayas+2} and to guarantee that the first inequality in \eqref{ToSayas+1}\,(c) cannot become an equality, we need to have \eqref{X0sm,m-1},\,\eqref{A2====}, which also causes our definition \eqref{LetNSoOP----1} to be complicated (cf.~Remark \ref{Twoformula-rema})\vspace*{-4pt}.
\end{itemize}
\item[(iii)]If we change \eqref{LetNSoOP----1}\,(iii) to $A_2=\big(Z(X_{\OnE}^{1+\d_{\rZeRo}^2+\d_{\rZeRo}^3}\tildeX_{\ZeRo}A_3^{a_3})^{-1}\big)^{\b_2}
 \big(1-a_4+a_4
(\tildeX_{\ZeRo}A_1^{a_1})^{\frac{(2-\d_0^2)\b_2}{2a_4}}\big)$ for any $a_4\in\Q_{\ne\frac12}$ with $0<a_4<1$ (and that $A_2$ only contains integral powers), then \eqref{B123X12Z} still holds and
all arguments can be perfectly conducted in parallel except in \eqref{A2====}, where we have to choose $a_4=\frac{2-\d_0^2}{4(1-\d_0^2)}$ in order to have \eqref{A2====}. This is why we define $A_2$ in \eqref{LetNSoOP----1}\,(iii).
\item[(iv)] We explain why we can guarantee to have \eqref{X0sm,m-1} [which implies that  the first inequality in \eqref{ToSayas+1}\,(c) cannot become  equality]: At the ``initial stage'', by \eqref{B123X12Z}\,(i),\,(iv), we have $\tildeX_1A_1^{a_1}=1+4\ell_0\ep+O(\ep)^2>1$, thus by \eqref{LetNSoOP----1}\,(iii), $\a':=A_2^{-1}\big(Z(X_{\OnE}^{1+\d_{\rZeRo}^2+\d_{\rZeRo}^3}\tildeX_{\ZeRo}A_3^{a_3})^{-1}\big)^{\b_2}$ [which is $\a_1$ defined in \eqref{ImMpP}\,(5.b), up to $O(\d)^2$$\ssc\,$] must be smaller than $1$ at the ``initial stage''. Hence by conditions in \eqref{ToSayas+1}, we can control $|\a'|$ to keep it not to grow too big such that $|\a'|\le1+O(\d)^2$. Then when $|\tildeX_1A_1^{a_1}|=1-\d$, we have \eqref{X0sm,m-1}.
    \end{itemize}\end{rema}
}
\noindent 
\vskip5pt

Now we can state our first proposition.
\begin{prop}\label{real00-inj}
There exist some  $\kappa_i\in\R_{>0}$ and some $\eta_i\in\R_{\ne0}$ such that we have  the following.

Denote by $V_1$  the  subset of $V$  consisting of  all
 elements $(p_{\ZeRo},p_{\OnE})=\big((x_{\ZeRo},y_{\ZeRo}),(x_{\OnE},y_{\OnE})\big)$ whose coordinates $x_1,x_2,y_2$ satisfy,
\begin{eqnarray}
\label{ToSayas}
\!\!\!\!\!\!\!\!\!\!\!\!&&
{\rm(a)\ }\dis1
<{\ssc}
(\kappa_1|x_2|)^{\eta_1}\le\kappa_2|x_1|^{\eta_2}\le(\kappa_1|x_2|)^{\eta_3}<\kappa_3,\ \ \ \ \ \
\nonumber\\\!\!\!\!\!\!\!\!\!\!\!\!&\!\!\!\!\!\!\!\!\!\!\!\!\!\!\!\!\!\!\!\!\!\!\!\!\!\!\!\!\!\!\!\!\!\!\!\!\!&
{\rm(b)\ }
\ell_{p_{\ZeRo},p_{\OnE}}:=\kappa_4|x_2|^{\eta_4}\cdot|x_2+y_2|+|x_{\OnE}|+|x_{\OnE}+y_{\OnE}|
\ge\kappa_5.
\end{eqnarray}
Denote by  $V_2$ the  subset of $S_2$  consisting of  all
 elements $(p_{\ZeRo},p_{\OnE})=\big((x_{\ZeRo},y_{\ZeRo}),(x_{\OnE},y_{\OnE})\big)$ whose coordinates $x_1,x_2,y_2$
satisfy%
~$[$we will prove later on that $
A_2
$
is
invertible in the set $W_2$ to be defined in \eqref{C+ToSayas+1}${\ssc\,}]$
, 
\NOUSE{
that there
exists 
some choice of an $\ell^4$-th root of $A_2$
,
 such that $A_1,\,A_1,\,A_2$ can be defined in
{\rm\eqref{LetNSoOP----1}} satisfying
}%
\begin{eqnarray}
\label{ToSayas+1}
\!\!\!\!\!\!\!\!\!\!\!\!\!\!\!\!\!\!\!\!\!\!\!\!\!\!\!\!\!\!\!\!\!\!\!\!\!\!\!\!\!\!\!\!&\!\!\!\!\!\!\!\!\!\!\!\!\!\!\!\!\!\!
&
{\rm(a)\ }\dis1
{\ssc}<{\ssc}
|A_{\rOnE}|^{-\kappa_1}{\ssc}\le{\ssc}
|A_2|^{-1}
\le{\ssc}|A_1|^{-\kappa_2}{\ssc}<{\ssc}\kappa_3,\ \ \ \ \
{\rm(b)\ }\ell_{p_{\ZeRo},p_{\OnE}}{\ssc}:={\ssc}|X_2|\cdot|A_1|^{\eta_1}
{\ssc}+{\ssc}|x_{\OnE}|{\ssc}+{\ssc}|x_{\OnE}{\ssc}+{\ssc}y_{\OnE}|
{\ssc}\ge{\ssc}\kappa_4
,\!\!\!\!\!\!\!\!\!\!\!\!\!\!\!\!\!\!\!\!\!\!
\nonumber\\
\!\!\!\!\!\!\!\!\!\!\!\!\!\!\!\!\!\!\!\!\!\!\!\!&\!\!\!\!\!\!\!\!\!\!\!\!\!\!\!\!\!\!\!\!\!\!\!\!\!\!\!\!\!\!\!\!\!\!\!\!\!&
{\rm(c)\ }
(1-\d)|X_2|\cdot|A_1|^{2+\d_0}
<|\widetilde X_1|<
(1+\d)|X_2|^{-1}\cdot|A_1|^{-\d_0}
.
\!\!\!\!\!\!\!\!\!\!\!\!\!\!\!\!\!\!\!\!\!\!\!\!\!\!\!\!\!\!\!
\end{eqnarray}
Then there exists $V_0$ which is either equal to $V_1$ or else equal to $V_2$ such that $V_0$ is a nonempty compact
subset of $\C^4$, and further when $(p_1,p_2)\in V_0$, we have
\equa{-EiathA0}{\dis x_1,x_2,x_2+y_2\ne0
.
}
\NOUSE{, and when $(p_1,p_2)\in V_0$, we have,
\begin{eqnarray}\label{-EiathA0}
&\!\!\!\!&{\rm(1)\ }
\kappa_{10}{}\le{}|x_{\ZeRo}|,|x_{\OnE}|{}\le{}\kappa_{11},
 \ \ \ {\rm(2)\ }|x_{\OnE}+y_{\OnE}|{}\ge\kappa_{12}, \ \
{\rm(3)\ }|A_3|>0\mbox{ in case }V_0=V_2.
\end{eqnarray}
\item[\rm(iii)]For any $(p_{\ZeRo},p_{\OnE})\in V_{\rZeRo}$, no equality can occur in the first or last inequality of
{\rm\eqref{ToSayas}}\,{\rm(a)} or {\rm\eqref{ToSayas+1}\,(a)},
or in  any
inequality of {\rm\eqref{ToSayas+1}}\,{\rm(c)},\,{\rm(d)};
further, two equalities cannot simultaneously occur in the second and third inequalities of {\rm\eqref{ToSayas}\,(a)} or {\rm\eqref{ToSayas+1}\,(a)}.
}
\end{prop}

\begin{defi}\rm\label{Clo-rem}
Denote by $W_1$  the  subset of $V$  consisting of  all
 elements $(p_{\ZeRo},p_{\OnE})=$\linebreak $\big((x_{\ZeRo},y_{\ZeRo}),(x_{\OnE},y_{\OnE})\big)$ whose coordinates $x_1,x_2,y_2$
 satisfy,
\begin{eqnarray}
\label{C+ToSayas}
\!\!\!\!\!\!\!\!\!\!\!\!&&
{\rm(a)\ }\dis1
\le{\ssc}
(\kappa_1|x_2|)^{\eta_1}\le\kappa_2|x_1|^{\eta_2}\le(\kappa_1|x_2|)^{\eta_3}\le\kappa_3,\ \ \ \ \ \
\nonumber\\\!\!\!\!\!\!\!\!\!\!\!\!&\!\!\!\!\!\!\!\!\!\!\!\!\!\!\!\!\!\!\!\!\!\!\!\!\!\!\!\!\!\!\!\!\!\!\!\!\!&
{\rm(b)\ }
\ell_{p_{\ZeRo},p_{\OnE}}:=\kappa_4|x_2|^{\eta_4}\cdot|x_2+y_2|+|x_{\OnE}|+|x_{\OnE}+y_{\OnE}|
\ge\kappa_5.
\end{eqnarray}
Denote by  $W_2$ the  subset of $V$  consisting of  all
 elements $(p_{\ZeRo},p_{\OnE})=\big((x_{\ZeRo},y_{\ZeRo}),(x_{\OnE},y_{\OnE})\big)$ whose coordinates $x_1,x_2,y_2$ satisfy,
\begin{eqnarray}
\label{C+ToSayas+1}
\!\!\!\!\!\!\!\!\!\!\!\!\!\!\!\!\!\!\!\!\!\!\!\!\!\!\!\!\!\!\!\!&\!\!\!\!\!\!\!\!\!\!\!\!\!
&
{\rm(a)\ }\dis1
{\ssc}\le{\ssc}
|A_{\rOnE}|^{-\kappa_1}{\ssc}\le{\ssc}
|A_2|^{-1}
\le{\ssc}|A_1|^{-\kappa_2}{\ssc}\le{\ssc}\kappa_3,\, \ \
{\rm(b)\ }\ell_{p_{\ZeRo},p_{\OnE}}{\ssc}:={\ssc}|X_2|\cdot|A_1|^{\eta_1}
{\ssc}+{\ssc}|x_{\OnE}|{\ssc}+{\ssc}|x_{\OnE}{\ssc}+{\ssc}y_{\OnE}|
{\ssc}\ge{\ssc}\kappa_4
,\!\!\!\!\!\!\!\!\!\!\!\!\!\!\!\!\!\!\!\!\!\!
\nonumber\\
\!\!\!\!\!\!\!\!\!\!\!\!\!\!\!\!\!\!\!\!\!\!\!\!&\!\!\!\!\!\!\!\!\!\!\!\!\!\!\!\!\!\!\!\!\!\!\!\!\!\!\!\!\!\!\!\!\!\!\!\!\!&
{\rm(c)\ }
(1-\d)|X_2|\cdot|A_1|^{2+\d_0}
\le|\widetilde X_1|\le
(1+\d)|X_2|^{-1}\cdot|A_1|^{-\d_0}
.
\!\!\!\!\!\!\!\!\!\!\!\!\!\!\!\!\!\!\!\!\!\!\!\!\!\!\!\!\!\!\!
\end{eqnarray}
\end{defi}
\begin{lemm}\label{lemm-condition-XZ}
We have
\equa{Dmememe}{
{\rm(i)\ }\ol V_1\subseteq W_1,\ \ \ \ \ \ \ {\rm(ii)\ }\ol V_2\subseteq W_2\subseteq S_2.
}
\end{lemm}
\noindent{\it Proof.~}We prove the more complicated case \eqref{Dmememe}\,(ii).
{Let $(p_1,p_2)\in W_2$.
{Since
negative
powers appear
in \eqref{C+ToSayas+1}\,(a), we need to interpret
\eqref{C+ToSayas+1}\,(a) as
\equa{interpret--as}{\dis\!\!\!\!\!\!\!
\kappa_3^{-1}\le|A_1|^{\kappa_2}\le|A_2|\le|A_1|^{\kappa_1}\le1
%
%
.
\!\!\!\!\!\!\!}
In particular,%
}%
~$A_1,
A_2$
are
invertible.
}%
We have $\ol V_2\subset V$ as $V$ is closed by Lemma \ref{Procesi-1}.
\NOUSE{
We have
\begin{eqnarray}
\label{invvvvv}
&&\!\!\!\!\!\!\!\!\!\!\!\!\!\!\!\!\!\!\!\!\!\!\!\!\!\!\!
{\rm(i)\ }
\Big|\frac15 + \frac{4 X_2^{\ell_0-2}}{5\widetilde X_1^2 Z^{\ell_0}}\Big|
\stackrel{{}^{\sc\rm\eqref{equa-Case6-lemm}\,(iv),\,\eqref{tX1==}}}{=}
\Big|\frac15 + \frac{4}{5X_1^2 }\Big|+O(\d)^3>\d_2^{200}+O(\d)^3>0
,
\nonumber\\
&&\!\!\!\!\!\!\!\!\!\!\!\!\!\!\!\!\!\!\!\!\!\!\!\!\!\!\!
{\rm(ii)\ }
\Big|
2 - \d_0 - \frac{(1 - \d_0) X_2^{\ell_0-2}}{
\widetilde X_1^2 Z^{\ell_0}}\Big|
\stackrel{{}^{\sc\rm\eqref{equa-Case6-lemm}\,(iv),\,\eqref{tX1==}}}{=}
\Big|
2 - \d_0 - \frac{1 - \d_0}{
X_1^2 }\Big|
+O(\d)^3>\d_2^{200}+O(\d)^3>0
.
\end{eqnarray}
Thus, $\frac15 + \frac{4 X_2^{\ell_0-2}}{5\widetilde X_1^2 Z^{\ell_0}}$ is invertible and $A_3$ is well-defined.
}%
We will see in \eqref{A1-A2-cond},\,\eqref{ImMpP} that \eqref{C+ToSayas+1} 
implies \eqref{meme}.
\NOUSE{
\eqref{ImMpP} that $|A_3|>\d_2$, thus $A_3,\,2 - \d_0 - \frac{(1 - \d_0) X_2^{\ell_0-2}}{
\widetilde X_1^2 Z^{\ell_0}}$ are
invertible, and $A_1,A_2$ are well-defined, and further $A_1,A_2$ are invertible by \eqref{C+ToSayas+1}\,(a).
\NOUSE{
We immediately see from  \eqref{C+ToSayas+1}\,(a) that
$A_1,A_2
$ are invertible
.
%
Further, by \eqref{meme}\,(iii),\,\eqref{equa-Case6-lemm}\,(iv),\,\eqref{tX1==},
$|\frac43 - \frac{X_2^{84} Z^8}{3\widetilde X_1^{88}}|=|\frac43 - \frac{X_2^{92}}{3X_1^{88}}|+O(\d)^2\ge\d_2^{200}+O(\d)^2>0$.
Thus $A_3$ is well-defined. We will also see from \eqref{ImMpP} that $|A_3|>1$, and so $A_3$ is invertible.
}%
We have $\ol V_2\subset V$ as $V$ is closed by Lemma \ref{Procesi-1}.
We see from \eqref{ToSayas+1}\,(a),\,(c) that $|A_1|,|A_2|,|\widetilde X_1|$ are $1+O(\d)^1$ elements, which with
\eqref{ToSayas+1}\,(d) shows that
\NOUSE{
in \eqref{A1-A2-cond},\,\eqref{ImMpP},\,\eqref{ImMpP+1} that we have the following,
\equa{Forllw}{\dis\!\!\!\!\!\!\!\!\!
{\rm(i)\,}\d_2{\sc\!}<{\sc\!}|A_3|,  \ {\rm(ii)\,}\d_2{\sc\!}<{\sc\!}|A_2|{\sc\!}<{\sc\!}\ell_2,
 \ {\rm(iii)\,}\d_2^{100}{\sc\!}<{\sc\!}\Big|\frac15 {\sc\!}+{\sc\!} \frac{4 X_2^{\ell_0-2}}
{5\widetilde X_1^2 Z^{\ell_0}}\Big|
,\
{\rm(iv)\,}\d_2^{100}{\sc\!\!}<{\sc\!\!}\Big|
2 {\sc\!\!}-{\sc\!\!} \d_0 {\sc\!\!}- {\sc\!\!}\frac{(1 {\sc\!\!}-{\sc\!\!} \d_0) X_2^{\ell_0-2}}{
\widetilde X_1^2 Z^{\ell_0}}\Big|
.
\!\!\!\!\!\!\!}
From this,\,\eqref{equa-Case6-lemm}\,(iv),\,\eqref{tX1==} and \eqref{C+ToSayas+1}\,(d),\,(e), we see that
}%
}%
%
~Hence, $W_2\subset S_2$.

Now let
 $(\hat p_1,\hat p_2)=\big((\hat x_1,\hat y_1),(\hat x_2,\hat y_2)\big)$ be any element in $\ol V_2$. By definition of $\ol V_2$, there \mbox{exists} a sequence
 $(p_{1,i},p_{2,i})=\big((x_{1,i},y_{1,i}),(x_{2,i},y_{2,i})\big)\in V_2$ for $i\in\Z_{\ge1}$
 converging to  $(\hat p_1,\hat p_2)$. Then \eqref{meme} 
 is satisfied by $(p_{1,i},p_{2,i})$, i.e.
,
 \equa{AboutXZ}{\dis
{\rm(i)\ } \d_2^4<|X_{1,i}|<\ell_2^4,\ \  \ \ \ \ \ {\rm(ii)\ }\d_2^4<|X_{2,i}|<\ell_2^4
 ,\ \ \ \ \ \ \ {\rm(iii)\ }\Big|2- \frac1{X_{1,i}}\Big|
  >\d_2^4
.} 
By taking the limit of the sequence, we see from \eqref{AboutXZ} that the hat versions
$\hat X_1,\hat X_2
$ of
$X_1,X_2
$, which are limits of $ X_{1,i},X_{2,i}
$, satisfy the following,
\equa{0wele,}{\dis
{\rm(i)\ } \d_2^4\le|\hat X_{1}|\le\ell_2^4,\ \ \ \ \ \ \ {\rm(ii)\ }\d_2^4\le|\hat X_{2}|\le\ell_2^4
   ,\ \ \ \ \ \ \ {\rm(iii)\ }\Big|2 - \frac1{\hat X_1}\Big|
   \ge\d_2^4
.}
Thus $\hat X_1,\hat X_2,\hat Z$ [cf.~\eqref{equa-Case6-lemm}\,(iv)$\ssc\,$]
are invertible, and
 $A_1
$ 
is 
well-defined at $(\hat p_1,\hat p_2)\in\ol V_2$%
,  and further, $A_1$ is invertible by
noting from \eqref{tX1==} that the tilde version of $\hat X_1$ is $\widetilde {\hat X}_1=\a_0\hat X_1$ and
then noting from \eqref{0wele,}\,(iii) that
$\big|2-\frac1{{\widetilde{\hat X}}{}_1}\big|$ is invertible%
. Thus
$
A_2 
$
is
well-defined
. Since \eqref{ToSayas+1} is satisfied by $(p_{1,i},p_{2,i})$,
by taking the limit of the sequence, we immediately see that \eqref{C+ToSayas+1} is satisfied by $(\hat p_1,\hat p_2)$. By definition of $W_2$, we have $(\hat p_1,\hat p_2)\in W_2$.
This completes the proof of \eqref{Dmememe}\,(ii).
\hfill$\Box$

\begin{rema}\rm\label{AboutX1}\begin{itemize}
\NOUSE{%
\item[(i)]
Note that $V_2$ is algebraically well-defined, namely, an element $(p_1,p_2)$ is in $V_2$ if and only if there
exists 
a choice of $A_2^{\d^4}$ 
(possibly there are more choices
), such that \eqref{LetNSoOP----1} and \eqref{ToSayas+1} hold.
When we say $(p_1,p_2)\in V_2$, we always take it for granted that some $A_1^{\d^4}
$ is
 already chosen.
By \eqref{LetNSoOP----1}\,(a), $|A_2|\ge\kappa_4^{-\kappa_2^{-1}}$
,
 we see that $A_2^{\d^4}
$ is
always  locally holomorphic
functions of $\tildeX_1,X_2,Z$.
}%
\item[(i)]
We emphasis that
\NOUSE{
Proposition \ref{real00-inj}\,(iii) is highly nontrivial, which plays the key role in
the proof of Proposition \ref{real00-inj+1}, and that, 
as  will be seen in the proof of Proposition \ref{real00-inj+1},
}%
the last two terms of $\ell_{p_1,p_2}$ in \eqref{ToSayas}\,(b) or \eqref{ToSayas+1}\,(b) play extremely important roles which guarantee the inequation  \eqref{GSoo+} has solution $(q_1,q_2)$ as required  [cf.~\eqref{Awww}$\ssc\,$].
\item[(ii)] When we prove Proposition \ref{real00-inj}, the first thing we need to do is to prove $V_0\ne\emptyset$, which will be done by choosing some suitable element $(p_1,p_2)\in V_0$. We regard such an element as the ``initial stage''. We will see from our arguments in this section that we observe the following important fact,
\begin{fact}\label{fact-initial}
    The ``initial stage'' controls globally the growths of $|x_2|,\,|x_2+y_2|$ on the whole set $V_0$ if we define $V_0$ by \eqref{ToSayas}.
\end{fact}
\item[(iii)]
We use the strange way in  \eqref{LetNSoOP----1} to define $A_1,A_2$ in order to satisfy Proposition \ref{real00-inj}. More precisely,
\begin{itemize}\item[(a)] Firstly, we will see in the proof of Lemma \ref{Prop(ii)Holds} 
that 
powers that appear in
\eqref{LetNSoOP----1}\,(i)
~play key roles
.
\item[(b)]Secondly, we need to prove
in Lemma \ref{V0NOT0} that $V_2\ne0$ by choosing the ``initial stage'' $(p_1,p_2)$, sufficiently close to $(\bar p_1,\bar p_2)$, satisfying \eqref{ToSayas+1}, such that \mbox{all} elements in \eqref{LetNSoOP----1}
 are  $1+O(\ep)^1$  numbers (though they are not necessarily real numbers) and
$X_2,Z$ have  are positive numbers of
the following form (where $\ep>0$ is sufficiently small),
\begin{eqnarray}
&\!\!\!\!\!
&
{\rm(i)\ }X_2=1
,\
\ \ \ \ \
{\rm(ii)\ }Z=1-\d^3\ep
.
\label{B123X12Z}
\NOUSE{
A_1{\ssc}={\ssc}1{\ssc}+{\ssc}
\d_0(1-\d_0)^2\ep+{\ssc}O(\ep)^2,
\!\!\!\!\!\!\!\!\!\!\!\!\!\!\!\!\!\!\!\!
\\
&\!\!\!\!\!\!\!\!\!\!\!\!\!\!\!\!\!\!\!\!\!\!\!\!\!\!\!\!\!\!\!\!\!\!\!\!\!\!\!\!\!\!\!\!\!\!\!\!\!\!\!\!\!\!\!\!\!\!\!\!\!\!\!\!\!\!\!\!\!\!\!\!\!\!\!\!\!\!\!\!\!\!\!\!\!\!\!\!\!\!\!\!\!\!\!\!\!\!\!\!\!\!\!\!\!\!\!\!\!\!\!\!\!\!\!\!\!\!\!\!\!\!\!\!\!\!\!\!&
\nonumber
{\rm(v)\ }A_2{\ssc}={\ssc}1{\ssc}+{\ssc}(1 + \d_0^2 + \d_0^3 + \d_0^4)\ep{\ssc}+{\ssc}O(\ep)^2,
\ \ \ {\rm(vi)\ }A_3{\ssc}={\ssc}1{\ssc}-\frac{\ell_0(1 + \d_0 - \d_0^2)}{1-\d_0}\ep+{\ssc}O(\ep)^2,
}\!\!\!\!\!\!\!\!\!\!\!\!\!\!\!\!\!\!\!\end{eqnarray}
\end{itemize}
\item[(iv)]We would like to point out that because of our choice of the ``initial stage'' in \eqref{B123X12Z} that $|Z|$ is smaller than $|X_2|$, which is the key points to obtain a contradiction when the last inequality of \eqref{C+ToSayas+1}\,(a) becomes  equality as there will exist an inconsistence in \eqref{C+ToSayas+1}%
. 
\NOUSE
{%
\item[(vi)] We explain why we can guarantee to have \eqref{X0sm,m-1} [which implies that  the first inequality in \eqref{ToSayas+1}\,(c) cannot become an equality]: At the ``initial stage'', by \eqref{B123X12Z}, we have $\frac{\tildeX_1^{20} X_2^{40}}{A_1^5 Z^{10}}=1+50\ep+O(\ep)^2>1$, thus by \eqref{LetNSoOP----1}\,(iii), $\a':=\frac{Z^5}{A_2A_3^{11}\tildeX_1^{10} X_2^{10}}$ [which is $\a_1$ defined in \eqref{ImMpP}\,(5.b), up to $O(\d)^2$$\ssc\,$] must be smaller than $1$ at the ``initial stage''. Hence by conditions in \eqref{ToSayas+1}, we can control $|\a'|$ to keep it not to grow too big such that $|\a'|\le1+O(\d)^2$. Then when $|\tildeX_1A_1^{\frac12}|=1-\d$ [noting that $\frac{\tildeX_1^{20} X_2^{40}}{A_1^5 Z^{10}}=\tildeX_1^{10}A_1^5$ by \eqref{LetNSoOP----1}\,(i)$\ssc\,$], we can have \eqref{X0sm,m-1}.
}%
\item[(v)] Note that we do not have any information about $y_1$ in \eqref{ToSayas},\,\eqref{ToSayas+1}. This is because we cannot put too many conditions in the definition of $V_0$. We can only use Theorem \ref{Theo-2} to control $|y_1|$ as long as we have some control on $|x_1|,|x_2|$
.
\end{itemize}\end{rema}

\NOUSE
{%
First we need some remarks.
\begin{rema}\rm\label{B123-unique}\begin{itemize}\item[(i)]
We define $V_2$ purely algebraically (without using any geometry). It says that an element $(p_1,p_2)\in V$ is in $V_2$
 if and only if there are some choices of complex numbers $A_1,A_2,A_3,\ell_{p_1,p_2}$ and
 some choices of rational powers of numbers that appear
 in \eqref{LetNSoOP----1}  such that \eqref{LetNSoOP----1},\,\eqref{ToSayas+1} and \eqref{ToSayas+2} hold
 [thus as far as Proposition \ref{real00-inj} is concerned we do not consider $A_1,A_2,...$, as functions, but
 instead as numbers related to $(p_1,p_2){\ssc\,}$].
Therefore $V_2$ is a well-defined subset of $V$.
\item[(ii)] We will prove that
 there exist $\tilde X_1,X_2,Z$ and the said choices, all close to $1$, such that
 \eqref{LetNSoOP----1},\,\eqref{ToSayas+1} and \eqref{ToSayas+2}  hold, i.e., $V_2$ contains an element.
\item[(iii)]
For a given $(\check p_1,\check p_2)\in V_2$, of course it is possible that there exist different choices of the
above satisfying \eqref{LetNSoOP----1},\,\eqref{ToSayas+1} and \eqref{ToSayas+2}.
However, when we say an element $(\check p_1,\check p_2)$ is in $V_2$,
we always take it for granted that all said choices are already given.
\item[(iv)]We do not require that the complex numbers $A_1,A_2,A_3,\ell_{p_1,p_2}$ (and said rational powers of numbers) are functions defined globally
[the fact is that we never need to use any global property about these numbers even in the proof of Theorem \ref{MAINT} at the end of subsection \ref{subsect3.2} except the fact that
 $L$ defined in \eqref{L===ss} is a bounded set, all our arguments are only concerned with local properties].
However they are locally multi-valued functions, and each branch of any of them is a locally holomorphic function (cf.~Lemma \ref{holo}), and since all powers appearing in
\eqref{LetNSoOP----1} are rational numbers, all said locally holomorphic multi-valued functions have only finite branches
(in fact there exists a fixed integer $N\in\Z_{>0}$ which only depend on $\ell_0,\ell$ such that all multi-valued functions have no more than $N$ branches). Further,
for any  $(\check p_1,\check p_2)\in V_2$,  there exists a small neighbourhood
such that different branches of any said locally-defined function can separate because there exists a
fixed number $c>0$ [which may depend on  $(\check p_1,\check p_2)\ssc\,$] such that associated with the element $(\check p_1,\check p_2)$,
the difference of any two choices
of any element $A_1,A_2,A_3,\ell_{p_1,p_2}$ (and said rational powers of numbers)  has an absolute value larger than $c$.
\item[(v)] For any $(\check p_1,\check p_2)\in V_2$, we will prove that there exists $(q_1,q_2)\in V_2$, close to  $(\check p_1,\check p_2)$,
such that
all numbers $A_1,A_2,A_3,\ell_{q_1,q_2}$, etc., associated with $(q_1,q_2)$, are chosen to be the same branches as that
associated with $(\check p_1,\check p_2)$ (thus they are
all holomorpic), and Proposition \ref{real00-inj+1} holds.
\end{itemize}
\end{rema}
}\NOUSE{%
Before stating Proposition \ref{real00-inj+1}, we first prove the following
\begin{lemm}
\label{V2-compact}If $V_0\ne\emptyset$ then it is compact in $V$ $($thus in $\C^4)$.
\end{lemm}
\noindent{\it Proof.~}We consider the more complicated case with $V_0=V_2$.
First by \eqref{-EiathA0}\,(1) and
Proposition \ref{pro-also}, $V_2$ is bounded [thus all $A_1,A_2,A_3,\ell_{p_1,p_2}$, etc., are bounded
by condition \eqref{ToSayas+1}$\ssc\,$].
Let $(p_{1i},p_{2i})=\big((x_{1i},y_{1i}),(x_{2i},y_{2i})\big)\in V_2$
be a sequence converging to $(\tilde p_1,\tilde p_2)=\big((\tilde x_{1},\tilde y_{1}),(\tilde x_{2},\tilde y_{2})$.
First by Lemma \ref{Procesi-1}, we have $(\tilde p_1,\tilde p_2)\in V$.
Let $A_{2i}^{\d^4}$ 
be the chosen $\ell^4$-th root of $A_2$ 
associated with $(p_{1i},p_{2i})$.
Then clearly, there exists a subsequence $(p_{1,i_j},p_{2i_j})$ such that  $A_{2,i_j}^{\d^4}
$
converges
; say,
it 
converges to $\tilde A_2^{\d^4}$, which must be an $\ell^4$-th root of  $A_2|_{(p_1,p_2)=(\tilde p_1,\tilde p_2)}$ since by definition $(A_{2,i_j}^{\d^4})^{\ell^4}=A_{2,i_j}$
.
Then
 $(\tilde p_1,\tilde p_2)$ associated with the chosen $\tilde A_2^{\d^4}
$ satisfies
\eqref{LetNSoOP----1},\,\eqref{ToSayas+1}. By definition, $(\tilde p_1,\tilde p_2)\in V_2$, i.e., $V_2$ is closed.\hfill$\Box$
}%

\begin{prop}\label{real00-inj+1}
For any $(\check p_{\ZeRo},\check p_{\OnE})\in V_{\rZeRo}$, in every neighborhood of $(\check p_{\ZeRo},\check p_{\OnE})$ there exists $(q_{\ZeRo},q_{\OnE})=\big((\dot x_{\ZeRo},\dot y_{\ZeRo}),(\dot x_{\OnE},\dot y_{\OnE})\big)\in V_{\rZeRo}$
sufficiently close to $(\check p_1,\check p_2)$ such that \equa{GSoo+}{\dis \ell_{q_{\ZeRo},q_{\OnE}}> \ell_{\check p_{\ZeRo},\check p_{\OnE}}.}
\end{prop}


We will first use  Proposition \ref{real00-inj} to give a proof of Proposition \ref{real00-inj+1}.
\NOUSE
{%
Before the proof, we would like to remark that
our proof is standard in the sense that it only uses some ``formal arguments'' which do not require the precise forms of $A_1,A_2,A_3$ but we only need to know the only fact that they are locally holomorphic and nonzero when $(p_1,p_2)\in V_0$, in order to obtain \eqref{dotB} below.
}%

\subsection{Proof of Proposition \ref{real00-inj+1} provided by 
Claudio Procesi  and proof of  Theorem \ref{MAINT}}\label{subsect3.2}
We thank Professor Claudio Procesi for providing us the material after \eqref{AOnnn} until the end of proof of Proposition \ref{real00-inj+1}  (for your reference we also  present our original arguments in Appendix \ref{our-prop33}).

\vskip5pt
\noindent{\it Proof of Proposition \ref{real00-inj+1}.~}Let $(\check p_{\ZeRo},\check p_{\OnE})=\big((\check x_{\ZeRo},\check y_{\ZeRo}),(\check x_{\OnE},\check y_{\OnE})\big)\in  V_{\rZeRo}$.
Note that \eqref{-EiathA0}
implies that  $\check x_{\ZeRo},\check x_{\OnE},\check x_{\OnE}+\check y_{\OnE}\ne0$.
We  need to find an element $(q_1,q_2)\in V_0$ which is sufficiently close to $(\check p_1,\check p_2)$. Thus we can assume $q_1=(\dot x_1,\dot y_1),\,q_2=(\dot x_2,\dot y_2)$ such that their coordinates have the following forms, for some sufficiently small $\ep_1>0$, and $u,v,s,t\in\C$ (for later convenient use, it is no harm to write the coordinate of $\dot x_2+\dot y_2$ instead of $\dot y_2$),
\equa{1++??+q0q1}{\!\!\!\!\!\!\!\!\!\!\!\!
\dot x_{\ZeRo}=\check x_{\ZeRo}(1{\ssc}+{\ssc} s\ep_1),\ \
\dot y_1=\check y_{\ZeRo}{\ssc}+{\ssc}t\ep_1,\ \
\dot x_{\OnE}=
\check x_{\OnE}(1{\ssc}
+{\ssc}u
\ep_1),\ \
\dot x_2{\ssc}+{\ssc} \dot y_{\OnE}{\ssc}={\ssc}
(\check x_2{\ssc}+{\ssc}\check y_{\OnE})(1{\ssc}+{\ssc}v\ep_1).
\!\!\!\!\!\!}

The local bijectivity of Keller maps says that for any $u,v\in\C$ in a bounded set there always exists sufficiently small $\ep_1>0$ such that,
there exist $s,t\in\C$
with $(q_{\ZeRo},q_{\OnE})\in V$, namely, $q_1\ne q_2$ [which is guaranteed by the fact that $(q_1,q_2)$ is sufficiently close to $(\check p_1,\check p_2)$ and $\check p_1\ne \check p_2$] and
\equa{detmm}{\mbox{$\big(F(q_1),G(q_1)\big)=\big(F(q_2),G(q_2)\big)$.}}
Since $F,G$ are polynomials, we can always solve $s$ (throughout the paper we do not need to use $t$ so we omit the solution of $t$) from \eqref{1++??+q0q1},\,\eqref{detmm} to obtain that $s$ is a power series of $\ep_1$ such that each coefficient of $\ep_1^k$ is a homogenous polynomial of $u,v$ with degree $k+1$ for all $k\ge0$; in particular, we can write $s$ as the following form, for some $a,b,\tilde a,\tilde b,\tilde c\in\C$,
\equa{s1t=}{\dis \!\!\!\!s{\ssc}={\ssc}s_{\ZeRo}{\ssc}+{\ssc}(\tilde a u^2{\ssc}+{\ssc}\tilde b u\tildev
{\ssc}+{\ssc}\tilde c\tildev ^2)\ep_1{\ssc}+{\ssc}O(\ep_1)^2
, \ \ \ s_{\ZeRo}{\ssc}={\ssc} au{\ssc}+{\ssc} b\tildev ,
\ \  \ (a,b)\ne(0,0).\!\!\!\!
}

First note that since $(q_1,q_2)$ is sufficiently close to $(\check p_1,\check p_2)$, we see that
\eqref{meme},\,%
\eqref{ToSayas+1}\,(c)%
~are
automatically satisfied by $(q_1,q_2)$  in case $V_0=V_2$%
.
Further when \eqref{GSoo+} holds for $(q_1,q_2)$, we always have
\equa{Alwaysmsm}{\mbox{ $\ell_{q_1,q_2}>\ell_{\check p_1,\check p_2}
\stackrel{{}{\sc\rm\eqref{ToSayas}\,(b)
}}{\ge} \kappa_5$
 \ \ (or respectively $ \stackrel{{}{\sc\rm\eqref{ToSayas+1}\,(b)}}{\ge}\kappa_4$)
,}} i.e., \eqref{ToSayas}\,(b) [or respectively \eqref{ToSayas+1}\,(b)$\ssc\,$] automatically holds for $(q_1,q_2)$.
Thus in order for $(q_1,q_2)$ to be in $V_0$ and to satisfy \eqref{GSoo+}, we only need to require  \eqref{ToSayas}\,(a) [respectively \eqref{ToSayas+1}\,(a)
$\ssc\,$] and \eqref{GSoo+} to be satisfied by $(q_1,q_2)$.

Observe from the first strict inequality of  {\rm\eqref{ToSayas}\,(a)} or {\rm\eqref{ToSayas+1}\,(a)} that
two equalities cannot simultaneously occur in the second and third inequalities of {\rm\eqref{ToSayas}\,(a)} or {\rm\eqref{ToSayas+1}\,(a)}.
Therefore%
, we only need to consider the following  two possible cases.
\vskip4pt \noindent{\it Case 1: Assume that, for $(\check p_1,\check p_2)$, all inequalities of {\rm\eqref{ToSayas}\,(a)} {\rm[}or respectively {\rm\eqref{ToSayas+1}\,(a)
$\ssc\,$]} are strict inequalities.}
\vskip4 pt

Then the same must hold for  $(q_1,q_2)$ (since $\ep_1$ is infinitesimal), i.e., $(q_1,q_2)$ is automatically in $V_0$ if \eqref{GSoo+} holds%
. Thus in this case we only need to consider one inequation, i.e.,
\eqref{GSoo+}. Therefore we see that  this case is easier than the case  we encounter below (as we always have no need to worry about $E_1$ appeared below).

\vskip4pt \noindent{\it Case 2:
Assume that, for $(\check p_1,\check p_2)$,
 equality occurs in either the second or else the third inequality of {\rm\eqref{ToSayas}\,(a)} $[$or respectively {\rm\eqref{ToSayas+1}\,(a)}${\ssc\,}]$.}
\vskip4pt

Accordingly, we need to consider two inequations for $(q_1,q_2)$
: one  is the inequation for $(q_1,q_2)$ corresponding to the part of
\eqref{ToSayas}\,(a) [or respectively \eqref{ToSayas+1}\,(a)$\ssc\,$], where  equality occurs for $(\check p_1,\check p_2)$; another is  \eqref{GSoo+}.

First assume we have case \eqref{ToSayas} and the second equality of \eqref{ToSayas}\,(a) holds for $(\check p_1,\check p_2)$ (the case for the third equality is similar), i.e., $\kappa_2\kappa_1^{-\eta_1}|\check x_1^{\eta_2}\check x_2^{-\eta_1}|=1$. Then the two inequations we need to consider for $(q_1,q_2)$ are
\equa{case1-2equ}{\mbox{$(\kappa_1|\dot x_2|)^{\eta_1}\le\kappa_2|\dot x_1|^{\eta_2}$ \ and \ $\ell_{q_1,q_2}
:=\kappa_4|\dot x_2|^{\eta_4}|\dot x_2+\dot y_2|+|\dot x_{\OnE}|+|\dot x_{\OnE}+\dot y_{\OnE}|>\ell_{\check p_1,\check p_2}$,}} which, by  \eqref{ToSayas}\,(b),\,\eqref{1++??+q0q1}
, can be rewritten as the following forms,
\begin{eqnarray}
\label{noB-form}
&\!\!\!\!\!\!\!\!\!\!\!\!\!\!\!\!&
E_1:=\kappa_2\kappa_1^{-\eta_1}|\dot x_1|^{\eta_2}\cdot|\dot x_2|^{-\eta_1}-1\stackrel{{}^{\sc\rm\eqref{1++??+q0q1}}}{=}|1+s\ep_1|^{\eta_2}\cdot|1+u\ep_1|^{-\eta_1}-1\ge0,
\nonumber
\\
&\!\!\!\!\!\!\!\!\!\!\!\!\!\!\!\!&
E_2{\ssc}:={\ssc}(\ell_{q_1,q_2}{\ssc}-{\ssc}\ell_{\check p_1,\check p_2})|\check x_2|^{-1}
\nonumber
\\
&\!\!\!\!\!\!\!\!\!\!\!\!\!\!\!\!&
\phantom{E_2}\!\!\!\!\!
\stackrel{{}^{\sc\rm\eqref{ToSayas}\,(b)}}{=}
\kappa_4|\dot x_2|^{\eta_4}\cdot|(\dot x_2{\ssc}+{\ssc}\dot y_2)\check x_2^{-1}|{\ssc}+{\ssc}
|\dot x_{\OnE}\check x_2^{-1}|{\ssc}+{\ssc}|(\dot x_{\OnE}{\ssc}+{\ssc}\dot y_{\OnE})\check x_2^{-1}|
\nonumber
\\
&\!\!\!\!\!\!\!\!\!\!\!\!\!\!\!\!&
\phantom{E_2}\ \ \ \ \ \ \
{\ssc}-{\ssc}\Big(\kappa_5|\check x_2|^{\eta_4}\cdot|(\check x_2
{\ssc}+{\ssc}\check y_2)\check x_2^{-1}|{\ssc}+{\ssc}1{\ssc}+{\ssc}|(\check x_{\OnE}{\ssc}+{\ssc}\check y_{\OnE})\check x_2^{-1}|\Big)\!\!\!\!\!\!
\nonumber
\\
&\!\!\!\!\!\!\!\!\!\!\!\!\!\!\!\!&
\phantom{E_2}\!\!{\ssc\!}
\stackrel{{}^{\sc\rm\eqref{1++??+q0q1}}}{=}
\kappa'_1|1{\ssc}+{\ssc}u\ep_1|^{\eta_4}\cdot|1{\ssc}+{\ssc}\tildev \ep_1|{\ssc}+{\ssc}|1{\ssc}
+{\ssc}u\ep_1|{\ssc}+{\ssc}\kappa'_2|1{\ssc}+{\ssc}\tildev \ep_1|{\ssc}-{\ssc}(\kappa'_1{\ssc}+{\ssc}1{\ssc}+{\ssc}\kappa'_2){\ssc}>{\ssc}0,
\end{eqnarray}
where, $\kappa'_1=\kappa_4|\check x_2|^{\eta_4}\cdot|\check x_2+\check y_2)\check x_2^{-1}|$,
$\kappa'_2=|(\check x_2+\check y_2)\check x_2^{-1}|$.

It may be worth mentioning that we have the following simple formulas,
for any $a_i\in\C,\,b_i\in\R$ which are independent of $\ep_1$,
\begin{eqnarray}
\label{MSMSMSMEN939393}
&\!\!\!\!\!\!\!\!\!\!\!\!\!\!\!\!\!\!\!\!\!\!\!\!&
{\rm(i)\ }
(1+a_1\ep_1)^{b_1}(1+a_2\ep_1)^{b_2}=1+(a_1b_1+a_2b_2)\ep_1+O(\ep_1)^2,
\nonumber\\
&\!\!\!\!\!\!\!\!\!\!\!\!\!\!\!\!\!\!\!\!\!\!\!\!&
{\rm(ii)\ }
|1+a_1\ep_1|^{b_1}=|1+a_1b_1\ep_1+O(\ep_1)^2|.\end{eqnarray}
Using \eqref{s1t=},\,\eqref{MSMSMSMEN939393}, and expanding each element inside the absolute value symbols in $E_1$ and $E_2$  in \eqref{noB-form} as a power series of $\ep_1$ (such that every coefficient of $\ep_1^k$ is a homogenous polynomial of $u,\tildev $ with degree $k$ for $k\ge1$), we see that  \eqref{noB-form}
can be always rewritten  as
the following forms,
where
 $f(u\ep_1,\tildev \ep_1),\ g(u\ep_1,\tildev \ep_1)$ are two  holomorphic functions
of $u,\tildev $ (when $\ep_1$ is regarded as fixed) in some neighbourhood of $(0,0)$ and vanishing at $(0,0)$,
 \begin{eqnarray}
\label{i.e.1@such111that=2}
\!\!\!\!\!\!\!\!\!\!\!\!\!\!\!\!\!\!\!\!\!\!\!\!\!
&
\!\!\!\!\!\!\!\!\!\!&
{\rm(i)\ }E_{\rOnE}{\ssc}:={\ssc}
\big|1+f(u\ep_1,\tildev \ep_1) 
\big|
-1 \ge0,
\nonumber\\[0pt] 
\!\!\!\!\!\!\!\!\!\!\!\!\!\!\!\!\!\!\!\!\!\!\!\!\!\!\!\!&
\!\!\!\!\!\!\!\!\!\!&
{\rm(ii)\, }E_2{\sc}:={\sc}
\kappa_{\ZeRo}'\big|1{\sc}+g(u\ep_1,\tildev \ep_1)
\big|
{\sc}+{\sc}|1{\sc}+{\sc}u\ep_1|{\sc}
+{\sc}\kappa'_{\OnE}|1{\sc}+{\sc}\tildev \ep_1|{\sc}-{\sc}
(\kappa'_{{\ZeRo}}{\sc}+{\sc}1+\kappa'_{\OnE}){\sc} >{\sc}0.\!\!\!\!\!\!\!\!\!\!\!\!\!\!\!\!\!
\end{eqnarray}

For the  case \eqref{ToSayas+1}, since $A_i\ne0$ for $i=1,2
$ by \eqref{ToSayas+1},\,\eqref{-EiathA0}, and $\ep_1>0$ is sufficiently small, we can always write
 $A_i|_{(p_1,p_2)=(q_1,q_2)}$ as, for some $\a_{1,i},\a_{2,i}\in\C$,
\equa{BB1222}{\dis A_i|_{(p_1,p_2)=(q_1,q_2)}=A_i|_{(p_1,p_2)=(\check p_1,\check p_2)}\Big(1+(\a_{1,i}u+\a_{2,i}\tildev )\ep_1+O(\ep_1)^2\Big).}
Thus by using \eqref{MSMSMSMEN939393}, we can also write the two inequations we need to consider as the forms in \eqref{i.e.1@such111that=2}.
\NOUSE{%
We prove that in case \eqref{ToSayas+1} we can also write the two inequations we need to consider as the forms in \eqref{i.e.1@such111that=2}.
Thus assume we have case \eqref{ToSayas+1} and the third equality of \eqref{ToSayas+1}\,(a) holds for $(\check p_1,\check p_2)$  (the case for the second equality is similar), i.e., $|\check A_1^{\kappa_3}\check A_2^{-\kappa_2}|=1$, where in general, we denote
\equa{bar=dot-B}{\mbox{$\check A_i=A_i(\check x_1,\check x_2,\check y_2)$, \ \ $\dot A_i=A_i(\dot x_1,\dot x_2,\dot y_2)$ \ \ for $i=1,2,3$.}} Then the two inequations we need to consider
 for $(q_1,q_2)$ are $|\dot A_2|^{\kappa_2}\le|\dot A_1|^{\kappa_3}$ and $\ell_{q_1,q_2}>\ell_{\check p_1,\check p_2}$, which can be rewritten as
 the following forms,
\begin{eqnarray}
\label{Bformmm}
&\!\!\!\!\!\!\!\!\!\!\!\!\!\!\!\!&
E_1:=|\dot A_1^{\kappa_3}\dot A_2^{-\kappa_2}|-1\ge0,
\\\nonumber
&\!\!\!\!\!\!\!\!\!\!\!\!\!\!\!\!&
E_2{\ssc\!}:={\ssc\!}(\ell_{q_1,q_2}{\ssc\!}-{\ssc\!}\ell_{\check p_1,\check p_2})|\check x_2|^{-1}{\ssc\!}={\ssc\!}\Big(|\dot A_1^{\eta_1}\dot A_3|{\ssc\!}+{\ssc\!}|\dot x_2|{\ssc\!}
+{\ssc\!}|\dot x_2{\ssc\!}+{\ssc\!}\dot y_2|{\ssc\!}-{\ssc\!}(|\check A_1^{\eta_1}\check A_3|{\ssc\!}+{\ssc\!}|\check x_2|{\ssc\!}+{\ssc\!}|\check x_2+\check y_2|)\Big)|\check x_2|^{-1}{\ssc\!}>{\ssc\!}1.
\end{eqnarray}
Since $A_i$ for $i=1,2,3$ is a local holomorphic function of variables $x_1,y_1,x_2$ when $(p_1,p_2)\in V_0$ (cf.~Lemma \ref{holo}), which does not take the zero value [by  \eqref{ToSayas+1}\,(a),\,\eqref{-EiathA0}$\ssc\,$], by \eqref{s1t=},\,\eqref{bar=dot-B}, we can expand $\dot A_i$ as a power series of $\ep_1$ (such that each coefficient of $\ep_1^k$ is a homogenous polynomial of $u,\tildev $ with degree $k$ for $k\ge1$) to obtain that $A_i$ has the following form, for some $a_{ki}\in\C$ (from the following we see that the fact that $\check A_i\ne0$ is important to us),
\equa{dotB}{\!\!\!\!\!\!\!\dot A_i{\ssc\!:}={\ssc\!}A_i(\dot x_1,\dot x_2,\dot y_2)
{\ssc\!}={\ssc\!}\check A_i\big(1{\ssc\!}+{\ssc\!}(a_{1i}u{\ssc\!}+{\ssc\!}a_{2i}\tildev )\ep_1{\ssc\!}+{\ssc\!}(a_{3i}u^2{\ssc\!}+{\ssc\!}a_{4i}u\tildev {\ssc\!}+{\ssc\!}a_{5i}\tildev ^2)\ep_1^2{\ssc\!}+{\ssc\!}O(\ep_1)^3\big).\!\!\!\!\!\!}
No matter what are $a_{ki}$'s in \eqref{dotB}, we can always use \eqref{dotB} in \eqref{Bformmm} and expand each of $E_1$ and $E_2$ as a power series of $\ep_1$, so that  \eqref{Bformmm} can be also rewritten as the forms in  \eqref{i.e.1@such111that=2}, where in this case $\kappa'_1=|\check A_1^{\eta_1}\check A_3\check x_2^{-1}|$.
Thus in
 any case, the two inequations we need to consider can be always written as the forms in \eqref{i.e.1@such111that=2}
.
For the  case \eqref{ToSayas+1}, since $A_2\ne0$ we can always write
 $A_2|_{(p_1,p_2)=(q_1,q_2)}$ as, for some $\a_1,\a_2\in\C$,
\equa{BB1222}{\dis A_2|_{(p_1,p_2)=(q_1,q_2)}=A_2|_{(p_1,p_2)=(\check p_1,\check p_2)}\Big(1+(\a_1u+\a_2\tildev )\ep_1+O(\ep_1)^2\Big).}
Thus we can also write the two inequations we need to consider as the forms in \eqref{i.e.1@such111that=2}.
}%

First observe that for any $\a=\a\re+\a\im\ii\in\C$ [recall Convention \ref{conv1}\,(1)$\ssc\,$], we have, 
\begin{eqnarray}
\label{AOnnn}
&\!\!\!\!\!\!\!\!\!\!\!\!\!\!\!\!\!\!\!\!\!\!\!\!&
|1+\a\ep_1|=\sqrt{(1+\a\re\ep_1)^2+(a\im\ep_1)^2}=1+\a\re\ep_1+\frac{(\a\im)^2}2\ep_1^2+O(\ep_1)^3.
\end{eqnarray}
Therefore the   coefficient of the $\ep_1$  term of $E_i$ in \eqref{i.e.1@such111that=2}  has the form
 $\ell_i(u,\tildev )\re$  for some linear homogeneous form $\ell_i(u,\tildev )$ in $u,\tildev $ for $i=1,2$.  If these two forms are linearly independent of course we can choose $u,\tildev $ so that the two forms take positive real values (and the proof of Proposition \ref{real00-inj+1} is then completed), otherwise there exist two complex numbers $\a,\b$  not both zero so that  setting  $u:=\a z,\, \tildev :=\b  z$ we can have that $\ell_i(\a z,\b  z)=0,\ i=1,2$.

Therefore we consider the later case and assume for instance that  $\a  $   is nonzero and set $w:=   \a z$.
Then
\NOUSE{%
$g(u\ep_1,\tildev \ep_1)$ becomes the form
\equa{gu-v===}{\mbox{$g(u\ep_1,\tildev \ep_1)=\g w^2\ep_1^2+O(\ep_1)^3$ for some $\gamma\in\C$,}} and thus
}%
by \eqref{AOnnn}
the $\ep_1^2$  coefficient  in  \eqref{i.e.1@such111that=2}\,(ii)    is of the form $A(w)$ below
for some $\gamma\in\C$ [as mentioned in Remark \ref{AboutX1}\,(i), here is the reason why we need the last two terms of $\ell_{p_1,p_2}$, which guarantee that we always have the last two non-negative terms of $A(w)$ and at least one of both is positive (noting that in general it is possible that $\a\ne0,\,\b=0$ or $\a=0,\,\b\ne0)\ssc\,$],
\equa{Awww}{\dis E_2=A(w)\ep_1^2+O(\ep_1)^3,\ \ \ A(w):=(\gamma w^2)\re+\frac{(w\im)^2}2+\kappa'_{\OnE}\frac{\big((\b\a^{-1} w)\im\big)^2}2 . }
Further
$E_1$ is either zero or   of the form
  $ (c w^k)\re\ep_1^k +O(\ep_1)^{k+1}$ for some $c\in\C_{\ne0}$, $ k\geq 2$.
Assume we have the later case as otherwise $E_1\ge0$ holds trivially. Namely,
\equa{C1====wlwl}{\dis E_1=(c w^k)\re\ep_1^k +O(\ep_1)^{k+1}.}
If $\g=0$ then we can easily choose $w\in\C$ with $w\im\ne0$ and $(c w^k)\re>0$ so that $E_1>0$ and  $E_2>0$ by \eqref{Awww} and \eqref{C1====wlwl}.

Thus assume $\g\ne0$.  We may take $w\in\C$ with $|w|=1,\, w=e^{\theta i}$ for $\th\in\R$ with $0\le\th<2\pi$.

  Notice now that for any integer $h_1\geq 1$  and any nonzero constant complex number $h_2$  the two sets  where $ (h_2 w^{h_1})\re >0$ respectively    $ (h_2 w^{h_1})\im >0$   are two open sets  of the circle $|w|=1$ of arclength  $\pi$.

  \begin{lemm}
The set where    $A(w)>0$  is   an open set  of the circle $|w|=1$  of arclength  $\ge\pi+\d$ for some $\d>0$ so there is a nonempty open set of the circle    $|w|=1$
  of arclength  $\ge\d$  where both inequalities
$A(w)>0$ and $(cw^k)\re>0$ hold.
\end{lemm}
\noindent{\it Proof.~}~The set  where   $(\g w^2)\re\le0$  is formed by two opposite arcs of  total arclength $\frac{\pi}2$  so at least two of the $4$ endpoints  of these arcs, where   $(\g w^2)\re=0$,  are different from $\pm 1$  where $w\im=0$, thus in a neighbourhood of these two points  $(\g w^2)\re+\frac{(w\im)^2}2>0$.
\hfill$\Box$
\vskip5pt
Note that we can always choose sufficiently small $\ep_1>0$ such that $0<\ep_1\ll\d$. This proves Proposition \ref{real00-inj+1}.\hfill$\Box$\vskip5pt

\noindent{\it Proof of Theorem \ref{MAINT}.~}~Now we use Proposition \ref{real00-inj+1} to prove Theorem \ref{MAINT}. \
Proposition \ref{real00-inj+1} (which says that the continuous function $\ell_{p_1,p_2}$ does not have the maximal value on the nonempty compact set $V_0$) immediately gives a contradiction, which proves
the first statement of Theorem \ref{MAINT}. The second statement follows from  \cite{K-m1,K-m2}.
\NOUSE
{%
Now we use Proposition \ref{real00-inj+1} to prove Theorem \ref{MAINT}. By
 Lemma \ref{V2-compact},
 $V_{\rZeRo}$ is a compact subset of $\C^4$. We consider the more complicated case with $V_0=V_2$.
Let [note that  for a given $(p_1,p_2)\in V_2$ if there exist different choices of $\ell_{p_1,p_2}$ satisfying Proposition \ref{real00-inj+1} then all
choices are in the set $L$),
\equa{L===ss}{\dis
L=\{\ell_{p_1,p_2}\,|\,(p_1,p_2)\in V_2\}.
}
 Then $L$ is bounded (cf.~proof of Lemma \ref{V2-compact}).
Let $\tilde \ell={\rm sup\,}L$ be the supremum of $L$. This means that there exists a sequence $(p_{1i},p_{2i})\in V_2$ (with chosen $A_{1i},A_{2i},A_{3i},
\ell^{(i)}_{p_{1i},p_{2i}}$, etc.) such that $\lim_{i\to\infty}\ell^{(i)}_{p_{1i},p_{2i}}=\tilde \ell$.
Choosing a subsequence $(p_{1,i_j},p_{2,i_j})$ such that $(p_{1,i_j},p_{2,i_j})$ and all $A_{1,i_j},A_{2,i_j}$, $A_{3,i_j},
\ell^{(i_j)}_{p_{1,i_j},p_{2,i_j}}$, etc., converge; say, they converge to  $(\tilde p_{1},\tilde p_{2})$ and $\tilde A_1,\tilde A_2,\tilde A_3$,
$\tilde \ell_{\tilde p_1,\tilde p_2}:=\tilde \ell$, etc. Then $(\tilde p_1,\tilde p_2)\in V_2$ with chosen  $\tilde A_1,\tilde A_2,\tilde A_3$,
$\tilde \ell_{\tilde p_1,\tilde p_2}:=\tilde \ell$, etc. Thus $\ell\in L$, which contradicts Proposition \ref{real00-inj+1}. This proves
the first statement of Theorem \ref{MAINT}. The second statement follows from  \cite{K-m1,K-m2}.
}%
\hfill$\Box$

\subsection{Proof of 
Proposition \ref{real00-inj}}
It remains to  prove Proposition \ref{real00-inj}.
We would like to mention that because we require $V_0$ to be closed, it seems to us that in general there is no way to define a system of inequations \eqref{ToSayas} or \eqref{ToSayas+1} satisfying  that $V_0$ is closed. We need some ``extra fact''.
By considering all possible different cases, each of which provides us some different ``extra fact'', we will be able to achieve the goal.

Before proceeding our proof of Proposition \ref{real00-inj}, let us give some further explanations. We aim to construct a nonempty subset $V_0$ satisfying, say for example, \eqref{ToSayas}
.
To do this, usually we can start with some given element
$$(\tilde p_1,\tilde p_2)=\big((\tilde x_1,\tilde y_1),(\tilde x_2,\tilde y_2)\big)\in V,$$
and of course we may choose some different $(\tilde p_1,\tilde p_2)$ for a different case.
Look at the condition \eqref{ToSayas}: if, for example, we take
$$\mbox{$\kappa_0=1,$ \ \ \ \
$\kappa_1=|\tilde x_2|^{-1}$, \ \ \ \ $\kappa_2=|\tilde x_1|^{-\eta_2},$ 
}$$ then \eqref{C+ToSayas}\,(a) can be satisfied by $(\tilde p_1,\tilde p_2)$ by choosing some suitable $\kappa_3>1.$

We always need to choose $\kappa_5$ in \eqref{ToSayas}\,(b) to be as big as possible
otherwise $V_0$ can possibly contain some extra elements which we may be unable to control somehow. In order to do this,  we should choose $(\tilde p_1,\tilde p_2)\in V$ such that $|\tilde x_2+\tilde y_2|$
is maximal (i.e., $\g_{|\tilde x_1|,|\tilde x_2|}$). In this case we can choose $\kappa_5$ to be
$$\b_0:=\kappa_4|\tilde x_2|^{\eta_4} \g_{|\tilde x_1|,|\tilde x_2|}+|\tilde x_2|+
\g_{|\tilde x_1|,|\tilde x_2|}.$$ In this way, we may have that $(\tilde p_1,\tilde p_2)\in \ol V_0$ (the closure of $V_0$). However then,  equality can  occur in the
first  inequality of \eqref{C+ToSayas}\,(a) and thus \eqref{ToSayas}\,(a) is not satisfied [i.e.,
$(\tilde p_1,\tilde p_2)\notin V_0$]. To avoid this, we have to enlarge $\kappa_5$ [this is why we choose $\a$ in \eqref{GMSMSMSMS}\,(b), respectively \eqref{GMSMSMSMS??}\,(b), to be as in \eqref{GMSMSM0999+1}, respectively
\eqref{GMSMSM0999}$\ssc\,$]. Then $(\tilde p_1,\tilde p_2)$ cannot be in $V_0$.
This means that we cannot choose the given element $(\tilde p_1,\tilde p_2)$ as the ``initial stage''.
We have to find some other element  as the ``initial stage'', which is the reason why we require some ``extra fact''.
\vskip5pt

First we need a lemma.

\begin{lemm}\label{g-1-lemm}
The  $\g_{k_{\ZeRo},k_{\OnE}}$ defined in $\eqref{Ak=1}$ is a  strictly  increasing function of $k_{\ZeRo}\in\R_{\ge0}$ when $k_{\OnE}\in\R_{\ge0}$ is fixed, i.e.,
\begin{eqnarray}
\label{wePPPP1+}
&&\!\!\!\!\!\!\!\!\!\!\!\!\!\!\!\!\!\!\!\!\!\!\!\!
\g_{k'_{\ZeRo},k_{\OnE}}>\g_{k_{\ZeRo},k_{\OnE}} \ \ \mbox{ if \ }k'_{\ZeRo}>k_{\ZeRo}\ge0,\ k_{\OnE}\ge0.
\end{eqnarray}
Consequently, $\g_{k_1,k_2}>0$ for any $k_1\in\R_{>0},\,k_2\in\R_{\ge0}$.
\end{lemm}
\noindent{\it Proof.~}%
For any $k'_{\ZeRo}>0$, let
\begin{eqnarray}\label{LetA=}\!\!\!\!\!\!\!\!\!\!\!\!\!\!\!\!
\b\!\!&:=&\!\!\max\{\g_{k_{\ZeRo},k_{\OnE}}\,|\,0\le k_{\ZeRo}\le k'_{\ZeRo}\}\nonumber\\[0pt]
\!\!\!\!\!\!\!\!\!\!\!\!\!\!\!\!\!\!\!\!\!\!\!\!&\stackrel{{}^{\sc\rm\eqref{Ak=1}}}{=}&\!\!\max\big\{\,|x_{\OnE}+y_{\OnE}|\ \big|\ (p_{\ZeRo},p_{\OnE})=
\big((x_{\ZeRo},y_{\ZeRo}),(x_{\OnE},y_{\OnE})\big)\in V,\, 0\le|x_{\ZeRo}|\le k'_{\ZeRo},\, |x_{\OnE}|= k_{\OnE}\,\big\}\ge0.\end{eqnarray}
 Assume
conversely that
there exists $k_{\ZeRo}<k'_{\ZeRo}$ with $k_{\ZeRo}\ge0$ such that $\g_{k_{\ZeRo},k_{\OnE}}=\b$.
We will use the local bijectivity of Keller maps  to obtain a contradiction.
 Let \equa{MSAMSMS}
{\mbox{$(\tilde p_{\ZeRo},\tilde p_{\OnE})=\big((\tilde x_{\ZeRo},\tilde y_{\ZeRo}),(\tilde x_{\OnE},\tilde y_{\OnE})\big)\in V$ \ \ \ with \ \ $|\tilde x_{\ZeRo}|=k_{\ZeRo},\ \ |\tilde x_{\OnE}|=k_{\OnE},\ \  |\tilde x_{\OnE}+\tilde y_{\OnE}|=\b$.}}
Using $ x_{\OnE}, y_{\OnE}$ as local coordinates in $V$ we choose $(q_1,q_2)\in V$ sufficiently close to $(\tilde p_1,\tilde p_2)$, satisfying, for some $v\in\C$ with $|v|$ being sufficiently small,
\begin{eqnarray}
\label{q0q1}
q_1=(\dot x_{ 1}, \dot y_1),\ \ q_{\OnE}:=(\dot x_2,\dot y_2)=
(\tilde x_{\OnE}
 ,\tilde y_{\OnE}+v),
\end{eqnarray}  such that
\begin{eqnarray}
\label{@suchthat=2-for}&\!\!\!\!\!\!\!\!\!\!\!\!\!\!\!\!\!\!\!\!\!\!\!\!&
 |\dot x_{\OnE}+\dot y_{\OnE}|=|\tilde x_{\OnE}+\tilde y_{\OnE}+v|>|\tilde x_{\OnE}+\tilde y_{\OnE}| .\end{eqnarray}
Since $|\tilde x_{\ZeRo}|=k_{\ZeRo}<k'_{\ZeRo}$, and locally $x_1$ is a holomorphic, hence a continuous function of  $(x_2,y_2)$, we have $|\dot x_{\ZeRo}|<k'_{\ZeRo}$ when $|v|>0$ is  sufficiently small.
This means that we can choose $(q_{\ZeRo},q_{\OnE})\in V$ with
$$\mbox{$|\dot x_{\ZeRo}|<k'_{\ZeRo},\ \ \ |\dot x_{\OnE}|=k_{\OnE}$, \ \ but \ $|\dot x_{\OnE}+\dot y_{\OnE}|>\b$,}$$ which is a contradiction with the definition of $\b$ in \eqref{LetA=}. This proves the first assertion of the lemma. Then for any $k_1>0$, $k_2\ge0$, we have $$\g_{k_1,k_2}
\stackrel{{}^{\sc\rm\eqref{wePPPP1+}}}{>}\g_{\frac{k_1}{2},k_2}\ge0,$$ proving the second.\hfill$\Box$
\vskip7pt

\noindent{\it Proof of Proposition   \ref{real00-inj}.~}~Now we proceed the proof of Proposition
\ref{real00-inj} case by case.
\begin{rema}\rm\label{P-remak}
 In the first five 
 cases, we are able to define $\ell_{p_1,p_2}$ such that we can choose the ``initial stage'',
which controls the growths of $|x_2|,\,|x_2+y_2|$ as
  mentioned in Fact \ref{fact-initial}, to guarantee that
$|x_2+y_2|$ can grow faster (or decrease slower) than $|x_1|+|x_2|$  when $|x_1|$ or $|x_2|$ goes to 
infinity [see for example 
\eqref{yyy6666}$\ssc\,$];
 thus by Theorem \ref{Theo-2}, $|x_2|$ cannot go too far from the correspondent value of the ``initial stage'',
which allow us to choose some 
$\kappa_4$ in \eqref{ToSayas}\,(a) so that for any element in $\ol V_0$,  equality cannot occur in the last inequality of
\eqref{C+ToSayas}\,(a).
\end{rema}
{\noindent\it {Case 1}:~Assume $\g_{k_1,k_2}$ is not a weakly increasing function of $k_2$, i.e., $\g_{k_{\ZeRo},k_{\OnE}'}>\g_{k_{\ZeRo},k_{\OnE}}$ for some $k_{\ZeRo}>0,\, k_{\OnE}> k_{\OnE}'>0$. }
\vskip4pt

This provides us the following ``extra fact'' \eqref{GMSMSM0999+1}:
\begin{itemize}
\item[(1)]  There exists a  $\d>0$  such that
\equa{k1-k1-isiis}{\dis(k_{\OnE}^{-1}k'_{\OnE})^\d\g_{k_{\ZeRo},k_{\OnE}'}>\g_{k_{\ZeRo},k_{\OnE}},}  since strict inequality   holds when $\d=0$.

\item[(2)]  Once such $\d$ is fixed, there is a constant $\kk_0$ so that for  $\kk>\kk_0$ we have:
\equa{GMSMSM0999+1}{\dis
\a:=(k_{\OnE}^{-1}k'_{\OnE})^\d\g_{k_{\ZeRo},k_{\OnE}'}+(k'_{\OnE}+\g_{k_{\ZeRo},k_{\OnE}'})\kk^{-1}>\g_{k_{\ZeRo},k_{\OnE}}+(k_{\OnE}+\g_{k_{\ZeRo},k_{\OnE}})\kk^{-1}.
}
\end{itemize}
We will see in \eqref{y1===0} what is the role the above ``extra fact'' plays.

We define $V_1=V_1(\kk)$, depending on the parameter $\kk>\kk_0$,  to be  the subset of $V$  consisting of    all
 elements $(p_{\ZeRo},p_{\OnE})\!=\!\big((x_{\ZeRo},y_{\ZeRo}),(x_{\OnE},y_{\OnE})\big)$  such that its coordinates $x_1,x_2,y_2$ satisfy,
\begin{eqnarray}
\label{GMSMSMSMS}\dis\!\!\!\!\!\!\!\!\!\!\!\!\!\!\!\!\!\!&&
{\rm(a)\ }1<(k_2|x_2|^{-1})^{\kk^{-1}}\le k_1^{-1}|x_1|\le(k_2|x_2|^{-1})^{\kk^{-1}+\kk^{-2}}<\kk^{\kk^{-1}+\kk^{-2}},
\nonumber\\[0pt]
\dis\!\!\!\!\!\!\!\!\!\!\!\!\!\!\!\!\!\!&&
{\rm(b)\ }\ell_{p_1,p_2}:=
(k_{\OnE}^{-1}|x_{\OnE}|)^\d|x_{\OnE}+y_{\OnE}|+(|x_{\OnE}|+|x_{\OnE}+y_{\OnE}|)\kk^{-1}\ge\a.
\end{eqnarray}

\begin{rema}\label{d-k-n}\rm\begin{itemize}\item[(i)]As mentioned in Remark \ref{rema3.1}\,(ii), we treat $ \kk  $  as a parameter  which will be fixed upon our requirement in the course of the proof.
\item[(ii)]
We put $\kk^{-1}$ in \eqref{GMSMSMSMS}\,(b) in order for the element $(|x_{\OnE}|+|x_{\OnE}+y_{\OnE}|)\kk^{-1}$ to be sufficiently smaller than $\a$ when \eqref{GMSMSMSMS} holds.
\item[(iii)]
The reason we put the factor $(k_{\OnE}^{-1}|x_{\OnE}|)^\d$ in  \eqref{GMSMSMSMS}\,(b) is to guarantee that when the last strict inequality of \eqref{GMSMSMSMS}\,(a) becomes an equality, $|x_2|$ will decrease, but  $|x_2+y_2|$ will grow.
\item[(iv)]Note that the number in the last term of \eqref{GMSMSMSMS}\,(a) cannot be too big, otherwise for an element in $\ol V_1$, when the last strict inequality of \eqref{GMSMSMSMS}\,(a) becomes  equality, $|x_1|$ will grow faster than $|x_2+y_2|$
 and  then we cannot apply Theorem \ref{Theo-2}.
\end{itemize}\end{rema}

First, multiplying \eqref{GMSMSMSMS}\,(b) by $\kk$ and re-denoting $\kk\ell_{p_1,p_2}$ as $\ell_{p_1,p_2}$, we see that
 \eqref{GMSMSMSMS} can be rewritten in the form   \eqref{ToSayas}.

\begin{rema}\label{prii}\rm
For  a point in $V_1$, by \eqref{GMSMSMSMS}\,(a), we have
 $1<(k_2|x_2|^{-1})^{\kk^{-1}}$, $(k_2|x_2|^{-1})^{\kk^{-1}+\kk^{-2}}<\kk^{\kk^{-1}+\kk^{-2}},$ which implies (ii) below,
\equa{Remak-equa}{\dis{\rm(i)\ } k_1 \le |x_1|<\kk^{\kk^{-1}+\kk^{-2}}k_1,\ \ \ \ \ {\rm(ii)\ }
k_2\kk^{-1}<|x_2|<k_2.}where (i) is obtained from the part  $1< k_1^{-1}|x_1|<\kk^{\kk^{-1}+\kk^{-2}}$ in  \eqref{GMSMSMSMS}\,(a).
 \end{rema}
Now we prove that $V_1(\kk)\ne\emptyset$ (for all $\kk>\kk_0$)  by choosing some suitable ``initial stage'' $(\check p_{\ZeRo},\check p_{\OnE})$ mentioned in Remark \ref{AboutX1}\,(ii).
By hypothesis, $ k_{\OnE}k'^{-1}_{\OnE}  >1$, so denote \equa{kPPP}{\dis
\check k_{\ZeRo}:=(k_{\OnE}k'^{-1}_{\OnE})^{\kk^{-1}}k_{\ZeRo}>k_{\ZeRo}.}
By definition, there exists,
\equa{GEXiSTS}{\mbox{ $(\check p_{\ZeRo},\check p_{\OnE})=\big((\check x_{\ZeRo},\check y_{\ZeRo}),(\check x_{\OnE},\check y_{\OnE})\big)\in V$ \ \ with \
$|\check x_{\ZeRo}|=\check k_{\ZeRo},\ \ |\check x_{\OnE}|=k'_{\OnE}$, \ $|\check x_{\OnE}+\check y_{\OnE}|=\g_{\check k_{\ZeRo},k'_{\OnE}}$.}}
When $(p_1,p_2)$ is set to $(\check p_1,\check p_2)$, we denote the middle three terms in \eqref{GMSMSMSMS}\,(a) by $T_1,T_2,T_3$ respectively.
Then we have the following three facts:
\begin{itemize}\item[(a)] $T_1:=(k_2|\check x_2|^{-1})^{\kk^{-1}}\stackrel{{}^{\sc\rm\eqref{GEXiSTS}}}{=}(k_2k_2'^{-1})^{\kk^{-1}}\stackrel{
{}^{\sc\rm\eqref{kPPP}}}{>}1$
 by \eqref{GEXiSTS} and \eqref{kPPP} respectively;
\item[(b)]$T_2:=k_1^{-1}|\check x_1|\stackrel{{}^{\sc\rm\eqref{GEXiSTS},\, \eqref{kPPP}}}{=}(k_2k_2'^{-1})^{\kk^{-1}}\stackrel{{}^{\sc\rm(a)}}{=}T_1$  by \eqref{GEXiSTS}, \eqref{kPPP} and (a) respectively;
\item[(c)]
$T_3:=(k_2|\check x_2|^{-1})^{\kk^{-1}+\kk^{-2}}\stackrel{{}^{\sc\rm(a)}}{=}T_1^{1+\kk^{-1}}>T_1\stackrel{{}^{\sc\rm(b)}}{=}T_2$ by (a) and (b) respectively.
\end{itemize}
Thus we obtain
\begin{eqnarray}
\label{Veriiify}
&\!\!\!\!\!\!\!\!\!\!\!\!\!\!\!\!\!\!\!\!\!\!\!\!\!\!\!\!\!\!&
1<T_1=T_2<T_3<\kk,
\end{eqnarray}
i.e., \eqref{GMSMSMSMS}\,(a) is satisfied by $(\check p_{\ZeRo},\check p_{\OnE})$.

As for   \eqref{GMSMSMSMS}\,(b),  recall that,  from Lemma \ref{g-1-lemm}, we have $\g_{\check k_{\ZeRo},k'_{\OnE}}>\g_{k_{\ZeRo},k'_{\OnE}}$.
Using \eqref{GEXiSTS} and the definition of $\a$ given by \eqref{GMSMSM0999+1}, we obtain: \begin{eqnarray}
\label{Veriiify+1}
&\!\!\!\!\!\!\!\!\!\!\!\!\!\!\!\!\!\!\!\!\!\!\!\!\!\!\!\!\!\!&
(k_{\OnE}^{-1}|\check x_{\OnE}|)^\d|\check x_{\OnE}{\ssc\!}+{\ssc\!}\check y_{\OnE}|
{\ssc\!}+{\ssc\!}(|\check x_{\OnE}|
{\ssc\!}+{\ssc\!}|\check x_{\OnE}
{\ssc\!}+{\ssc\!}\check y_{\OnE}|)\kk^{-1}
\stackrel{{}^{\sc\rm\eqref{GEXiSTS}}}{=}(k_{\OnE}^{-1}k'_{\OnE})^\d\g_{\check k_{\ZeRo},k'_{\OnE}}+(k'_{\OnE}+\g_{\check k_{\ZeRo},k'_{\OnE}})\kk^{-1}
\nonumber\\&\!\!\!\!\!\!\!\!\!\!\!\!\!\!\!\!\!\!\!\!\!\!\!\!\!\!\!\!\!\!&
\phantom{(k_{\OnE}^{-1}|\check x_{\OnE}|)^\d|\check x_{\OnE}+\check y_{\OnE}|+(|\check x_{\OnE}|+|\check x_{\OnE}+\check y_{\OnE}|)\kk^{-4}
}
\!\!\!\!\!\!\!\!\!\!\!\!
\stackrel{{}^{\sc\rm\eqref{wePPPP1+},\,\eqref{kPPP}}}{
>}(k_{\OnE}^{-1}k'_{\OnE})^\d\g_{k_{\ZeRo},k'_{\OnE}}{\ssc\!}+{\ssc\!}(k'_{\OnE}{\ssc\!}+{\ssc\!}\g_{k_{\ZeRo},k'_{\OnE}})\kk^{-1}
\stackrel{{}^{\sc\rm\eqref{GMSMSM0999+1}}}{=}\a,\end{eqnarray}
namely, \eqref{GMSMSMSMS}\,(b) is satisfied by $(\check p_{\ZeRo},\check p_{\OnE})$. Hence we see that the ``initial stage'' $(\check p_1,\check p_2)$ is in $V_1$. We take $V_0=V_1\ne\emptyset$.
\medskip

Observe that
$$\lim_{\kk\to\infty} \kk^{\kk^{-1}+\kk^{-2}}=1 . $$     So we can take $\kk_0$ so that if $\kk>\kk_0$ we have     $\kk^{\kk^{-1}+\kk^{-2}}<2$.      Thus by \eqref{Remak-equa}, we obtain, when $\kk>\kk_0$   and  $(p_{\ZeRo},p_{\OnE})\!=\!\big((x_{\ZeRo},y_{\ZeRo}),(x_{\OnE},y_{\OnE})\big)\in V_1$,\equa{k1k2--that}{\mbox{(i) $k_1\le|x_1|\le k_1\kk^{\kk^{-1}+\kk^{-2}}<2k_1$, \ \ \ \ \ (ii) $k_2\kk^{-1}\le|x_2|\le k_2$.}}
 In particular,
we have $x_1,x_2\ne0$.
%
Assume \eqref{GMSMSMSMS} holds with $x_{\OnE}+y_{\OnE}=0$.
Then \eqref{GMSMSMSMS}\,(b) shows
that $|x_2|\ge\a\kk>k_2$ by \eqref{GMSMSM0999+1} (when $\kk$ is sufficiently larger than $\kk$ as $\g_{k_{\ZeRo},k_{\OnE}}>0$), a contradiction with \eqref{k1k2--that}\,(2). Thus we have \eqref{-EiathA0} [we do not have $A_3$  in the present case]. Further, $V_0$ is bounded by \eqref{k1k2--that} and Proposition \ref{pro-also}.

Now we want to prove that $V_0$ is closed. Thus we let $(p_{\ZeRo},p_{\OnE})\in\ol V_{\rZeRo}$.
By Lemma \ref{lemm-condition-XZ}, first assume that the first strict inequality of
\eqref{GMSMSMSMS}\,(a) becomes an equality for $(p_1,p_2)$. Then we immediately obtain that $|x_{\OnE}|=k_{\OnE},\,|x_{\ZeRo}|=k_{\ZeRo}$.

For this point thus we must have
\equa{ell-p1p2===}{\dis
\ell_{p_1,p_2}\stackrel{{}^{\sc\rm\eqref{GMSMSMSMS}\,(ii)}}{=}
 |x_{\OnE}+y_{\OnE}|+( k_2 +|x_{\OnE}+y_{\OnE}|)\kk^{-1}\stackrel{{}^{\sc\rm\eqref{GMSMSMSMS}\,(ii)}}{\ge}\a.} By definition \eqref{Ak=1}, we have $\g_{k_{\ZeRo},k_{\OnE}}\ge|x_{\OnE}+y_{\OnE}|,$
  which gives \begin{eqnarray}
\label{y1===0}
\!\!\!\!\!\!\!\!\!\!\!\!\!\!\!\!&&
\g_{k_{\ZeRo},k_{\OnE}}+(k_{\OnE}+\g_{k_{\ZeRo},k_{\OnE}})\kk^{-1}\stackrel{{}^{\sc\rm\eqref{Ak=1},\,\eqref{ell-p1p2===}}}{\ge}
\ell_{p_1,p_2}\stackrel{{}^{\sc\rm\eqref{ell-p1p2===}}}{ \ge } \a.
\end{eqnarray}
%
The above is a contradiction with the ``extra fact'' \eqref{GMSMSM0999+1}, which proves that the first strict inequality of
\eqref{GMSMSMSMS}\,(a) is satisfied by $(p_1,p_2)$.

This also implies that the second and third inequalities cannot be both equalities.
\\ \indent
Continuing, now assume the last strict inequality of \eqref{GMSMSMSMS}\,(a) becomes an  equality for $(p_1,p_2)$, i.e., $k_2|x_2|^{-1}=\kk$,
we denote \equa{denn-beta}{\dis
\b=\frac12(k_{\OnE}^{-1}k'_{\OnE})^\d\g_{k_{\ZeRo},k_{\OnE}'}>0.}
 Then by \eqref{GMSMSM0999+1},
\equa{b-alpha}{\mbox{ $2\b<\a,\ \text{ \ or equivalently, \ }\ \a-\b>\b$.}}
First assume \equa{fir-ASS}{\mbox{ $\a-(|x_{\OnE}|+|x_{\OnE}+y_{\OnE}|)\kk^{-1}\ge\b$.}}
Then by \eqref{GMSMSMSMS}\,(b),
\equa{by-wehavEVEV}{\mbox{$|x_2+y_2|\stackrel{{}^{\sc\rm\eqref{GMSMSMSMS}\,(b),\,\eqref{fir-ASS}}}
{\ge}\b(k_2^{-1}|x_2|)^{-\d}=\b\kk ^{\d}$.}}
We want to prove that \eqref{by-wehavEVEV} holds in general. Thus assume now we do not have \eqref{fir-ASS}, i.e., \begin{eqnarray}
\label{fir-ASS+1}
\!\!\!\!\!\!\!\!\!\!\!\!\!\!\!\!\!\!\!\!\!\!\!\!&&
\a-(|x_{\OnE}|+|x_{\OnE}+y_{\OnE}|)\kk^{-1}<\b,
\ \ \ \implies\ \ \ \ |x_2|+|x_2+y_2|>(\a-\b)\kk\stackrel{{}^{\sc\rm\eqref{b-alpha}}}{>}\b\kk,\ \ \ \implies
\nonumber\\
\!\!\!\!\!\!\!\!\!\!\!\!\!\!\!\!\!\!\!\!\!\!\!\!&&
|x_2+y_2| > \b\kk-|x_2|=  \b \kk-k_2\kk^{-1}
.\end{eqnarray}
Again for $\kk$ sufficiently large  we have  $ \b \kk-k_2\kk^{-1}> \b \kk^\d$, thus \eqref{by-wehavEVEV} holds in general. 

Now then $|y_2|\ge |x_2+y_2|-|x_2|\ge \b \kk^{\d}-k_2\kk^{-1}\ge \frac{\b}{2}\kk^\d$, which,
%
%
 together with the definition of $h_{p_1,p_2}$ in \eqref{dp0p1}  implies
\equa{h-yxppp}{\dis h_{p_1,p_2}\stackrel{{}^{\sc\rm\eqref{dp0p1}}}{\ge}|y_{\OnE}|
>\frac{\b}{2}\kk^{\d}.}
Note from \eqref{MSmde33333}, Remark \ref{rema3.1}\,(ii) and Remark \ref{d-k-n} that 
we can choose $\kk$ as large as we wish. When $\kk$   tends to 
infinity also  $h_{p_1,p_2}\to\infty$ by \eqref{denn-beta},\,\eqref{h-yxppp}. Thus we can apply Theorem \ref{Theo-2} to obtain that $h_{p_1,p_2}\sim_{\kk\ssc\,}|x_1|+|x_2|$, which with \eqref{k1k2--that} shows that $h_{p_1,p_2}\sim_{\kk\ssc\,}1$, a contradiction with \eqref{h-yxppp}
[by \eqref{k1k2--that},\,\eqref{by-wehavEVEV} and by Remark \ref{Remmma} (to be given later), one can see why we can get a contradiction]. This proves that \eqref{GMSMSMSMS} is satisfied by $(p_1,p_2)$. By definition, $(p_1,p_2)\in V_0$. Thus $\ol V_0=V_0$.


The above shows that with $V_1=V_1(\kk)$ being defined by \eqref{GMSMSMSMS}, we have, for $\kk$ sufficiently large, Proposition\,{\rm\ref{real00-inj}}. Case 1 is now completed.\vskip7pt

 Thus from now on, we can assume that $\g_{k_{\ZeRo},k_{\OnE}}$ is a weakly increasing function of $k_{\OnE}\in\R_{>0}$ when $k_{\ZeRo}\in\R_{>0}$ is fixed, i.e.,
\begin{eqnarray}
\label{wePPPP1}
&&\!\!\!\!\!\!\!\!\!\!\!\!\!\!\!\!\!\!\!\!\!\!\!\!
\g_{k_{\ZeRo},k_{\OnE}}\ge\g_{k_{\ZeRo},k_{\OnE}'} \ \ \mbox{ if \ }k_{\ZeRo}>0,\ k_{\OnE}> k_{\OnE}'>0.
\end{eqnarray}
\vskip4pt

{\noindent\it Case 2:~Assume
$\g_{\bar k^{1+\d'}k_{\ZeRo},\bar kk_{\OnE}}\ge \bar k\sTH{\ssc\,}\g_{k_{\ZeRo},k_{\OnE}}$ for some $\d'\in\R_{\ge0},\,\bar k,k_{\ZeRo},k_{\OnE}\in\R_{>0}$ with $\bar k>1,\,\d'<\frac1m.$}
\vskip4pt

By choosing $\d''$ with $\d'<\d''<\frac1{\ssTH m}$ and by Lemma \ref{g-1-lemm} [or \eqref{wePPPP1+}$\ssc\,$], we have
\equa{mwnenenernr}{\mbox{$\g_{\bar k^{1+\d''}k_{\ZeRo},\bar kk_{\OnE}}
\stackrel{{}^{\sc\rm\eqref{wePPPP1+}}}{>}
\g_{\bar k^{1+\d'}k_{\ZeRo},\bar kk_{\OnE}}\ge \bar k\sTH{\ssc\,}\g_{k_{\ZeRo},k_{\OnE}}$.}}
This, as in Case 1, provides us the following ``extra fact'' with $\d,\kk,\kk$ satisfying \eqref{MSmde33333} [here we require that $\ell\gg\max\{\bar k^{1+\d''}k_{\ZeRo},\bar kk_{\OnE}\}$, cf.~Remark \ref{rema3.1}\,(i)$\ssc\,$],
\equa{GMSMSM0999}{\dis\!\!\!\!\!
\a:=(\bar kk_{\OnE})^{-(1+\d)}\g_{\bar k^{1+\d''}k_{\ZeRo},\bar k k_{\OnE}}+(\bar k k_{\OnE}+\g_{\bar k^{1+\d''}k_{\ZeRo},\bar kk_{\OnE}})\nn^{-1}
>k_{\OnE}^{-(1+\d)}\g_{k_{\ZeRo}, k_{\OnE}}+( k_{\OnE}+\g_{k_{\ZeRo},k_{\OnE}})\nn^{-1}
.\!\!\!}
We define $V_1$ to be the subset of $V$  consisting of    all
 elements $(p_{\ZeRo},p_{\OnE})\!=\!\big((x_{\ZeRo},y_{\ZeRo}),(x_{\OnE},y_{\OnE})\big)$  such that its coordinates $x_1,x_2,y_2$ satisfy (cf.~Remark \ref{d-k-n}),
\begin{eqnarray}
\label{GMSMSMSMS??}\dis\!\!\!\!\!\!\!\!\!\!\!\!\!\!\!\!\!\!&&
{\rm(a)\ }1<\big(k_{\OnE}^{-1}|x_{\OnE}|\big)^{1+\d''-\kk^{-1}}\le k_{\ZeRo}^{-1}|x_{\ZeRo}|\le
\big(k_{\OnE}^{-1}|x_{\OnE}|\big)^{1+\d''+\kk^{-1}}<\kk^{1+\d''+\kk^{-1}},\ \ \ 
\nonumber\\[-2pt]\dis\!\!\!\!\!\!\!\!\!\!\!\!\!\!\!\!\!\!&&
{\rm(b)\ }\ell_{p_1,p_2}:=\frac{|x_{\OnE}+y_{\OnE}|}{\,|x_{\OnE}|^{\ssTH
1+\d
}\,}+(|x_{\OnE}|+|x_{\OnE}+y_{\OnE}|)\nn^{-1}\ge\a
.
\end{eqnarray}
As in Case 1,
we can rewrite \eqref{GMSMSMSMS??} as the form in \eqref{ToSayas}.
\begin{rema}\rm\label{FirRemma}
The reason we put the power $1+\d$, which is bigger than $1$, in the denominator of \eqref{GMSMSMSMS??}\,(b),
is to ensure that $|x_2+y_2|$ will grows faster than $|x_2|$ when  the last strict inequality of
\eqref{GMSMSMSMS??}\,(a) becomes  equality  for an element in $\ol V_1$.
\end{rema}

Now we prove that $V_1\ne\emptyset$ by choosing some suitable ``initial stage'' $(\check p_{\ZeRo},\check p_{\OnE})$ mentioned in Remark \ref{AboutX1}\,(ii).
By definition, there exists\equa{GEXiSTS+}{\!\!\!\!\!\!\!\!\!\!\mbox{ $(\check p_{\ZeRo},\check p_{\OnE})=\big((\check x_{\ZeRo},\check y_{\ZeRo}),(\check x_{\OnE},\check y_{\OnE})\big)\in V$  with
$|\check x_{\ZeRo}|=\bar k^{1+\d''}k_{\ZeRo},\ |\check x_{\OnE}|=\bar kk_{\OnE}$,
 $|\check x_{\OnE}+\check y_{\OnE}|=\g_{\bar k^{1+\d''}k_{\ZeRo},\bar kk_{\OnE}}$.}\!\!\!\!}
When $(p_1,p_2)$ is set to $(\check p_1,\check p_2)$, we denote the middle three terms in \eqref{GMSMSMSMS??}\,(a) by $T_1,T_2,T_3$ respectively.
Then using \eqref{GEXiSTS+} and definition of $\a$ in \eqref{GMSMSM0999}, one can verify
\begin{eqnarray}
&\!\!\!\!\!\!\!\!\!\!\!\!\!\!\!\!\!\!\!\!\!\!\!\!\!&
{\rm(a)\ }1<T_1:=(k_2^{-1}|\check x_2|)^{1+\d''-\kk^{-1}}\stackrel{{}^{\sc\rm\eqref{GEXiSTS+} }}{=}\bar k^{1+\d''-\kk^{-1}}<T_2:=k_1^{-1}|\check x_1|\stackrel{{}^{\sc\rm\eqref{GEXiSTS+} }}{=}\bar k^{1+\d''}
\nonumber\\
&\!\!\!\!\!\!\!\!\!\!\!\!\!\!\!\!\!\!\!\!\!\!\!\!\!&
\phantom{{\rm(a)\ }1}<T_3:=
(k_2^{-1}|\check x_2|)^{1+\d''+\kk^{-1}}\stackrel{{}^{\sc\rm\eqref{GEXiSTS+} }}{=}\bar k^{1+\d''+\kk^{-1}}<\kk^{1+\d''+\kk^{-1}},
\nonumber\\
&\!\!\!\!\!\!\!\!\!\!\!\!\!\!\!\!\!\!\!\!\!\!\!\!\!&
{\rm(b)\ }\frac{|\check x_2+\check y_2|}{|\check x_2|^{1+\d}}+(|\check x_2|+|\check x_2+\check y_2|)\nn^{-1}
\stackrel{{}^{\sc\rm\eqref{GEXiSTS+} }}{=}\frac{\g_{\bar k^{1+\d''}k_{\ZeRo},\bar k k_{\OnE}}}{(\bar kk_{\OnE})^{1+\d}}+(\bar k k_{\OnE}+\g_{\bar k^{1+\d''}k_{\ZeRo},\bar kk_{\OnE}})\nn^{-1}\stackrel{{}^{\sc\rm\eqref{GMSMSM0999} }}{=}\a,
\end{eqnarray}
i.e.,  \eqref{GMSMSMSMS??} is satisfied by $(\check p_{\ZeRo},\check p_{\OnE})$, namely,  the ``initial stage'' $(\check p_{\ZeRo},\check p_{\OnE})$ is in $V_1$.
We take $V_0=V_1\ne\emptyset$.

Exactly similar to Case 1, we can obtain that $V_0$ is a bounded set and \eqref{-EiathA0} holds.


Let $(p_{\ZeRo},p_{\OnE})\in\ol V_{\rZeRo}$. By Lemma \ref{lemm-condition-XZ},
first assume  the first strict inequality of
\eqref{GMSMSMSMS??}\,(a) becomes  equality for $(p_1,p_2)$. Then we obtain that $|x_{\OnE}|=k_{\OnE},\,|x_{\ZeRo}|=k_{\ZeRo}$.
We have, where the first inequality follows from
the definition of $\g_{k_{\ZeRo},k_{\OnE}}$ in \eqref{Ak=1}, while the second
from \eqref{GMSMSMSMS??}\,(b),
\begin{eqnarray}
\label{y1===0+}
\!\!\!\!\!\!\!\!\!\!\!\!\!\!\!\!&&
k_{\OnE}^{-(1+\d)}\g_{k_{\ZeRo}, k_{\OnE}}+( k_{\OnE}+\g_{k_{\ZeRo},k_{\OnE}})\nn^{-1}
\stackrel{{}^{\sc\rm\eqref{Ak=1}}}{\ge}
\frac{|x_{\OnE}+y_{\OnE}|}{\,|x_{\OnE}|^{\ssTH
1+\d
}\,}+(|x_{\OnE}|+|x_{\OnE}+y_{\OnE}|)\nn^{-1}\stackrel{{}^{\sc\rm\eqref{GMSMSMSMS??}\,(b)}}{\ge}\a,\end{eqnarray}
which is a contradiction with the ``extra fact'' \eqref{GMSMSM0999}.

Now assume  
the last strict inequality of \eqref{GMSMSMSMS??}\,(a) becomes  equality for $(p_1,p_2)$. Then using \eqref{MSM1111}, \eqref{MSmde33333} [cf. Remark \ref{rema3.1}\,(ii) and Remark \ref{d-k-n}], one
obtains, where (ii) is obtained by
 the second and third inequalities of  \eqref{GMSMSMSMS??}\,(a),
\equa{MSnenenene}{\mbox{${\rm(i)\ } |x_{\OnE}|=k_2\kk\sim_{\kk\ssc\,}\kk$; \ \ (ii) $|x_{\ZeRo}|
=k_1\big(k_{\OnE}^{-1}|x_{\OnE}|\big)^{1+\d''+O(\kk^{-1})^1}\stackrel{{}^{\sc\rm\eqref{GMSMSMSMS??}\,(a)}}{\sim_{\kk\ssc\,}}\kk^{1+\d''}$.}} 
Note from Theorem \ref{Theo-2} or \eqref{mqp1234-2} that
\equa{Thatwww}{\mbox{$\dH_{p_{\ZeRo},p_{\OnE}}\stackrel{{}^{\sc\rm\eqref{dp0p1},\,\eqref{mqp1234-2} }}{\sim_{\kk\ssc\,}}|x_{\ZeRo}|+|x_{\OnE}|\stackrel{{}^{\sc\rm\eqref{MSnenenene}}}{\sim_{\kk\ssc\,}}\kk^{1+\d''}$ when $\kk\gg1$,}}
but by \eqref{GMSMSMSMS??}\,(b), as in Case 1
, we have
\equa{mxxxxxx-ssss}{\mbox{
$|x_{\OnE}+y_{\OnE}|\stackrel{{}^{\sc\rm\eqref{GMSMSMSMS??}\,(b)}}{\succeq_{\kk\ssc\,}}|x_{\OnE}|^{\ssTH
1+\d
}\stackrel{{}^{\sc\rm\eqref{MSnenenene}}}{\sim_{\kk\ssc\,}} \kk^{\ssTH
1+\d
}\stackrel{{}^{\sc\rm\eqref{MSnenenene}}}{\succ_{\kk\ssc\,}}|x_2|$,}} and thus
 (noting that here is the only place we make use of the fact that $\d''<\frac1m$),
\equa{Thus-mwnene}{\dis |y_{\OnE}|\stackrel{{}^{\sc\rm\eqref{mxxxxxx-ssss}}}{\sim_{\kk\ssc\,}}|x_{\OnE}+y_{\OnE}|
\stackrel{{}^{\sc\rm\eqref{mxxxxxx-ssss}}}{\succ_{\kk\ssc\,}}\kk\succ_{\kk\ssc\,}\kk^{\frac{(1+\d'')m}{1+m}}
\stackrel{{}^{\sc\rm\eqref{Thatwww}}}{\sim_{\kk\ssc\,}} \dH_{_{\sc\,p_{\ZeRo},p_{\OnE}}}^{^{\sc\frac{m}{m+1}}},}
 a contradiction with
\eqref{mqp1234-2} [by \eqref{MSnenenene},\,\eqref{mxxxxxx-ssss} and by Remark \ref{Remmma} (to be given later), one can see why we can get a contradiction].

This proves that \eqref{GMSMSMSMS??} is satisfied by $(p_1,p_2)$, i.e., $(p_1,p_2)\in V_0$, and so $V_0$ is closed.
 Proposition
{\rm\ref{real00-inj}} holds. Case 2 is now completed.\vskip7pt

Hence from now on, we can assume
\equa{AssG}{\dis
\g_{\bar k^{1+\d'}k_{\ZeRo},\bar kk_{\OnE}}< \bar k\sTH{\ssc\,}\g_{k_{\ZeRo},k_{\OnE}}\mbox{ for any $\d'\in\R_{\ge0},\,\bar k,k_{\ZeRo},k_{\OnE}\in\R_{>0}$ with $\bar k>1,\,\d'<\frac1m.$}}
\vskip5pt
{Then we can prove
\begin{lemm}\label{Anlemmmmmm}For any $k_1,k_2\in\R_{>0}$, we have $\g_{k_1,k_2}>
k_2\sTH.$
\end{lemm}\noindent\noindent{\it Proof.~}Assume $\g_{k_1,k_2}\le k_2$ for some $k_1,k_2\in\R_{>0}$.
Choosing $k'_1\in\R_{>0}$ with $k'_1<k_1$, by Lemma \ref{g-1-lemm}, we obtain that
$
{\mbox{$\g_{k'_1,k_2}<\g_{k_1,k_2}\le
k_2\sTH$}}$. Then $\a:=\frac{\gamma_{k'_1,k_2}}{k_2
}<
1
$. By \eqref{AssG}, we have
$\g_{\kk k'_1,\kk k_2}< \kk\sTH{\ssc\,} \g_{k'_1,k_2}=\kk\sTH  k_2
\a$ for all $\kk\gg1$. Let
\equa{meleTTT}{\mbox{$\dis (\tilde p_1,\tilde p_2)\in V$ \ \ with \ \ $|\tilde x_1|=\kk k'_1$, \ \
$|\tilde x_2|=\kk k_2$, \ \ $|\tilde x_2+\tilde y_2|=\g_{\kk k'_1,\kk k_2}<\kk\sTH
k_2
\a$.}}
By Proposition \ref{Sect3-Lemm5} or \eqref{mqp1234-2}, we have [using notation \eqref{MSM1111}$\ssc\,$]
\equa{thahahaha}{\mbox{ $h_{\tilde p_1,\tilde p_2}\stackrel{{}^{\sc\rm\eqref{dp0p1},\,\eqref{mqp1234-2}}}{\sim_{\kk\ssc\,}}|\tilde x_1|
+|\tilde x_2|\stackrel{{}^{\sc\rm\eqref{meleTTT}}}{\sim_{\kk\ssc\,}}\kk\sTH $,}} but then
\equa{ButTherne}{|\tilde y_2|\ge|\TH \tilde x_2\sTH|-|\tilde x_2+\tilde y_2|\stackrel{{}^{\sc\rm\eqref{meleTTT}}}{>}(
1
 -\a)\kk\sTH k_2
\stackrel{{}^{\sc\rm\eqref{thahahaha}}}{ >}\tau h_{_{\sc \tilde p_1,\tilde p_2}}^{^{\sc\frac{\ssTH m}{\ssTH m+1}}}\mbox{ \ (when $\kk\gg1$)},} which is a contradiction with
\eqref{mqp1234-2}. 
\hfill$\Box$
\vskip7pt
}%

Recall that we fix
some choices of positive numbers
satisfying
\eqref{MSmde33333} (cf.~Remark \ref{rema3.1}) and that the element $(\bar p_1,\bar p_2)$ satisfies
\eqref{TaKa}.
We choose $(q_1,q_2)=\big((\dot x_1,\dot y_1),(\dot x_2,\dot y_2)\big)\in V$ sufficiently close to $(\bar p_1,\bar p_2)$ such that its coordinates have the form as in \eqref{1++??+q0q1}, i.e.,
\equa{1++??+q0q1-bar}{\!\!\!\!\!\!\!\!\!\!\!\!
\dot x_{\ZeRo}=\bar x_{\ZeRo}(1{\ssc}+{\ssc} s\ep),\ \
\dot y_1=\bar y_{\ZeRo}{\ssc}+{\ssc}t\ep,\ \
\dot x_{\OnE}=
\bar x_{\OnE}(1{\ssc}
+{\ssc}u
\ep),\ \
\dot x_2{\ssc}+{\ssc} \dot y_{\OnE}{\ssc}={\ssc}
(\bar x_2{\ssc}+{\ssc}\bar y_{\OnE})(1{\ssc}+{\ssc}v\ep).
\!\!\!\!\!\!}
and as in \eqref{s1t=}, 
we can solve from \eqref{detmm} to obtain, for some $a_\kk,b_\kk\in\C$ [here for later convenience, we re-denote $a,b$ in \eqref{s1t=}
as $- a_\kk, b_\kk$],
\begin{eqnarray}
\label{@suc2hthat=4}&\!\!\!\!\!\!\!\!\!\!\!\!\!\!\!\!\!\!\!\!\!\!\!\!&
s=-a_\kk u+b_\kk \tildev +O(\ep)^1,
\ \ \ (a_\kk,b_\kk)\ne(0,0).
\end{eqnarray}
%
{%
We will see that the numbers $a_\kk,b_\kk$ play crucial roles in our definition of $V_0$. First we prove
\begin{lemm}\label{anLemmA}
We have $a_\kk\ge0,\,
b_\kk\ge0
$.\end{lemm}
\noindent\noindent{\it Proof.~}%
{Assume
$a_{\kk\,\rm im}\ne0$ or $a_{\kk\,\rm re}<0$ 
[cf.~Convention \ref{conv1}\,(1)$\ssc\,$]. Then
we can choose
$u,v\in\C$
in the following way: we always choose $u\re<0,\,v\re>0$ so that
\eqref{suchThaT}\,(ii),\,(iii) hold; if $a_{\kk{\rm\,im}}\ne0$, we can choose $u\im\ne0$ with
$a_{\kk\,\rm im}u\im$ being sufficiently large 
to guarantee
\eqref{u-v====}\,(iii) holds; if $a_{\kk\,\rm re}<0$
, we can simply take $u\im=
\tildev \im=
0$, $v\re=1$ and $u\re\ll-1$  to guarantee
\eqref{u-v====}\,(iii) [thus \eqref{suchThaT}\,(i) holds],  i.e.,
}%
\equa{u-v====}{\dis {\rm(i)\ }u\re
<0,\ \ {\rm(ii)\ }\tildev \re>0
,\  \mbox{ and such that \ }{\rm(iii)\ }s\re
\stackrel{{}^{\sc\rm\eqref{@suc2hthat=4}}}{=}(-a_\kk u+b_\kk \tildev )\re+O(\ep)^1<0.}
Then by \eqref{1++??+q0q1-bar} and \eqref{u-v====}, we have
\begin{eqnarray}
\label{suchThaT}&\!\!\!\!\!\!\!\!\!\!\!\!\!\!\!\!\!\!\!\!\!\!\!\!\!\!\!\!\!\!\!\!\!\!\!\!\!\!&
{\rm(i)\ }0<k_{\ZeRo}:=|\dot x_{\ZeRo}|=\kk|1+s\ep|
\stackrel{{}^{\sc\rm\eqref{u-v====}\,(iii)}}{<}\kk,\ \ \ \ \ \ \ {\rm(ii)\ }0<k_{\OnE}:=|\dot x_{\OnE}|=\kk|1+u\ep|<\kk,
\nonumber\\[0pt]
&\!\!\!\!\!\!\!\!\!\!\!\!\!\!\!\!\!\!\!\!\!\!\!\!\!\!\!\!&
{\rm(iii)\ }|\dot x_{\OnE}+\dot y_{\OnE}|=
\g_{\kk,\kk}|1+\tildev \ep|>\g_{\kk,\kk}.\end{eqnarray}
This means that $0<k_{\ZeRo}<\kk$ and $0<k_{\OnE}<\kk $ with the following [where the first inequality follows from definition \eqref{Ak=1}, the second from \eqref{suchThaT}\,(iii)$\ssc\,$], \equa{Assss}{\dis \g_{k_{\ZeRo},k_{\OnE}}=\g_{|\dot x_1|,|\dot x_2|}\stackrel{{}^{\sc\rm \eqref{Ak=1}}}
{\ge} |\dot x_{\OnE}+\dot y_{\OnE}|
\stackrel{{}^{\sc\rm \eqref{suchThaT}\,(iii)}}
{>}\g_{\kk,\kk },}
which is a contradiction with the fact obtained from Lemma \ref{g-1-lemm} and \eqref{wePPPP1} that $\g_{k_1,k_2}$ is an increasing function on both variables. Thus $a_\kk \ge0$
.

Similarly, if
$b_{\kk\,\rm im}\ne0$ or $b_{\kk\,\rm re}<0$, 
then we can choose
$u,v\in\C$
in the following way: again we always choose $u\re<0,\,v\re>0$ so that
\eqref{suchThaT}\,(ii),\,(iii) hold; if $b_{\kk{\rm\,im}}\ne0$, we can choose $v\im\ne0$ with
$-b_{\kk\,\rm im}v\im$ being sufficiently large 
to guarantee
\eqref{u-v====}\,(iii) holds; if $b_{\kk\,\rm re}<0$
, we can simply take $u\im=
\tildev \im=
0$, $v\re\gg1$ and $u\re=-1$  to guarantee
\eqref{u-v====}\,(iii) holds.
\hfill$\Box$
}%

{%
\begin{lemm}\label{YYYy1==}Let $\d'\in\R_{>0}$ be such that $\d'>\ln(\kk)^{-1}$ $[$where $\ln(\cdot)$ is the natural logarithmic function$]$, we have
\equa{ggggg}{\dis
\kk\sTH<\g_{\kk,\kk}<(1+\d')
\kk\sTH.}
\end{lemm}\noindent\noindent{\it Proof.~}The first inequality follows from Lemma \ref{Anlemmmmmm}. By \eqref{dp0p1},\,\eqref{mqp1234-2}, when $\kk\gg1$, we have $h_{\bar p_1,\bar p_2}\sim_{\kk\ssc\,}|\bar x_1|+|\bar x_2|\sim_{\kk\ssc\,}\kk$. If $\g_{\kk,\kk}\ge(1+\d')
 \kk\sTH$, then 
\equa{Msme22229}{\mbox{ $\dis |\bar y_2|\ge|\TH\bar x_2\sTH+\bar y_2|-|\TH\bar x_2\sTH|
\stackrel{{}^{\sc\rm\eqref{TaKa}}}{\ge}\d'\kk\sTH\ge\ln(\kk)^{-1}\kk\sTH\succ\kk^{\frac{\ssTH m}{\ssTH m+1}}\sim h_{_{\sc \bar p_1,\bar p_2}}^{^{\sc\frac{\ssTH m}{\ssTH m+1}}},$}} a contradiction with \eqref{mqp1234-2}.\hfill$\Box$
}

{%
\vskip7pt
{\noindent\it Case 3: Assume 
for any $\SS\in\R_{>0}$ there exists $\kk>\SS$ such that
$\ssTH b_\kk<1+\frac1m
+ a_\kk$
.}
\vskip4pt
%
Then we can choose sufficiently small $\d'>0$ [which can
depend on $\kk$ but is independent of $\nn$ when $\nn\gg\kk$, cf.~notation \eqref{MSmde33333} and
Remark \ref{rema3.1}\,(i)$\ssc\,$]
such that we have the following ``extra fact'' (note that the following holds when $\d'=0$ thus
holds when $\d'>0$ is sufficiently small),
\equa{NSNENNENE}{\ssTH (1+\d')b_\kk<1+\frac1m-\d'+ a_\kk.}
We define  $V_1$ to be the subset of $V$  consisting of    all
 elements $(p_{\ZeRo},p_{\OnE})\!=\!\big((x_{\ZeRo},y_{\ZeRo}),(x_{\OnE},y_{\OnE})\big)$  such that its coordinates $x_1,x_2,y_2$ satisfy  (cf.~Remark \ref{d-k-n}),
\begin{eqnarray}
\label{NSoOP}&&\!\!\!\!\!\!\!\!\!\!\!\!\!\!\!\!\!\!\!\!\!\!\!\!\!\!\!\!
\dis{\rm(i)\ }
1{\sc}<{\sc}\big(\kk^{-1}|x_{\OnE}|\big)^{1+\frac1m-\d'-\nn^{-1}}{\sc}\le{\sc}
\kk^{-1}|x_{\ZeRo}|{\sc}\le{\sc}
(\kk^{-1}|x_{\OnE}|)^{1+\frac1m{\sc}}<{\sc}\nn^{1+\frac1m},\ \ \ 
\nonumber\\[-2pt]&&\!\!\!\!\!\!\!\!\!\!\!\!\!\!\!\!\!\!\!\!\!\!\!\!\!\!\!\!
{\rm(ii)\ }\ell_{p_1,p_2}:=
\frac{\g_{\kk,\kk}^{-1}|x_{\OnE}+y_{\OnE}|}{\big(\kk^{-1}|x_{\OnE}|\big)^{\ssTH
1+\d'}}
{\sc}+(|x_{\OnE}|+|x_{\OnE}+y_{\OnE}|)\ep^3\ge{\sc}1{\sc}+{\sc}\ep^2.
\end{eqnarray}
Similar to \eqref{GMSMSMSMS??}, we put the power $1+\d'$, which is bigger than $1$, in the denominator of \eqref{NSoOP}\,(ii)
is to ensure that $|x_2+y_2|$ will grows faster than $|x_2|$ when  equality occurs in the last strict inequality of
\eqref{NSoOP}\,(ii) for an element in $\ol V_1$.

We will prove $V_1\ne\emptyset$, but firstly, as in Cases 1 and 2, we can rewrite \eqref{NSoOP} as the form in \eqref{ToSayas},
and  $V_1$ is bounded such that $x_1,x_2\ne0$ if $(p_1,p_2)\in V_1$.
Secondly, we explain the reason  
we put $\ep^2,\ep^3$ in \eqref{NSoOP}\,(ii) is that we want to guarantee the following
\equa{gAuaua}{\dis (|x_{\OnE}|+|x_{\OnE}+y_{\OnE}|)\ep^3<\ep^2\mbox{ when }(p_1,p_2)\in V_1.}
which is obtained by Theorem \ref{Theo-2} and the facts from \eqref{NSoOP}\,(a) that $|x_1|,|x_2|<\nn^2$ and $\ep^{-1}\gg\nn$.
Then \eqref{NSoOP}\,(ii) with \eqref{gAuaua} gives
\equa{x2+ya2}{\dis |x_2+y_2|
\stackrel{{}^{\sc\rm\eqref{NSoOP}\,(ii),\,\eqref{gAuaua} }}{>}\g_{\kk,\kk}\big(\kk^{-1}|x_{\OnE}|\big)^{\ssTH
1+\d'},}
which is also satisfied by any element in $\ol V_0$ by Lemma \ref{lemm-condition-XZ}.
In particular, we have
\eqref{-EiathA0}.

Now we prove that $V_1\ne\emptyset$ by choosing some suitable ``initial stage'' $(q_1,q_2)$ mentioned in Remark \ref{AboutX1}\,(ii).
 With $(q_1,q_2)$, sufficiently close to $(\bar p_1,\bar p_2)$, being defined in \eqref{1++??+q0q1-bar}, we want to choose suitable $u,v$ such that
\eqref{NSoOP} is satisfied by $(q_{\ZeRo},q_{\OnE})$, i.e.,
\begin{eqnarray}
\label{1NSoOP}
&&\!\!\!\!\!\!\!\!\!\!\!\!\!\!\!\!\!\!\!\!\!\!\!\!\!\!\!\!
\dis{\rm(i)\ }
1<|1+u\ep|^{1+\frac1m-\d'-\nn^{-1}}\le|1+s\ep|\le
|1+u\ep|^{1+\frac1m}<\nn^{1+\frac1m},
\nonumber\\[-4pt]&&\!\!\!\!\!\!\!\!\!\!\!\!\!\!\!\!\!\!\!\!\!\!\!\!\!\!\!\!
{\rm(ii)\ }
\frac{\Big|1+\tildev \ep\Big|}{|1+u\ep|^{\ssTH
1+\d'
}}+O(\ep)^3\ge1+\ep^2.
\end{eqnarray}
We take $u,v\in\R_{>0}$ such that
\equa{LetF3rst}{\mbox{$\dis u=1,\ \ \
\tildev =
\frac{a_\kk+1+\frac1m-\d'}{b_\kk},$
\ \ \ \ then \ \ $\dis s\stackrel{{}^{\sc\rm\eqref{@suc2hthat=4}}}{=}1+\frac1m-\d'+O(\ep)^1$,}}
where the last equation is obtained  from \eqref{@suc2hthat=4}.
Then the coefficients of $\ep$ in the middle three terms of \eqref{1NSoOP}\,(i) are respectively
\equa{mdmdmrmrmrm}{\mbox{$\dis 1+\frac1m-\d'-\nn^{-1},$ \ \ \ \ $\dis 1+\frac1m-\d'$, \ \ \ \ $\dis 1+\frac1m$,}} i.e.,
all inequalities  in \eqref{1NSoOP}\,(i) are strict inequalities
by recalling from  notation \eqref{MSmde33333} that we  have complete freedom in choosing $\ep$ with
$0<\ep\ll\eE_1=\nn_1^{-1}$ independently of all other choices of the parameters.
Further, the coefficient of $\ep^1$ in the left hand-side of \eqref{1NSoOP}\,(ii) is
\equa{bebebr3848}{\dis v-(1+\d')u
\stackrel{{}^{\sc\rm\eqref{LetF3rst}}}{=}\frac{a_\kk+1+\frac1m-\d'}{b_\kk}-\ssTH (1+\d'),}
which is bigger than $0$ by the ``extra fact'' \eqref{NSNENNENE}. Thus the ``initial stage'' $(q_1,q_2)$ is in $V_1$.
We take $V_0=V_1\ne\emptyset$%
.

Next we let $(p_{\ZeRo},p_{\OnE})\in \ol V_{\rZeRo}$.
By Lemma \ref{lemm-condition-XZ}, first assume the first strict inequality of \eqref{NSoOP}\,(i) becomes  equality for $(p_1,p_2)$. Then $|x_1|=|x_2|=\kk$, but by
\eqref{x2+ya2} (which also holds for elements in $\ol V_0$), we have $|x_{\OnE}+y_{\OnE}|>\g_{\kk,\kk}$, a contradiction with the definition of $\g_{\kk,\kk}$ in \eqref{Ak=1}.

Now assume 
the last strict inequality of \eqref{NSoOP}\,(i) becomes  equality for $(p_1,p_2)$.
Then we  obtain,  when $\nn\gg\kk$ [using \eqref{MSM1111},\,\eqref{MSmde33333+}, cf.~Remark \ref{rema3.1}\,(ii)$\ssc\,$],
 \equa{ansnen0000}{\mbox{$|x_{\OnE}|=\kk\nn\sim_{\nn\ssc\,}\nn$, \ \ \ \
 $|x_{\ZeRo}|\le\kk\nn^{1+\frac1m}\sim_{\nn\ssc\,}
 \nn^{1+\frac1m}$,}} but by
\eqref{x2+ya2},
\equa{x2y2-x2}{\dis
|x_{\OnE}+y_{\OnE}|\stackrel{{}^{\sc\rm\eqref{x2+ya2}}}{\succeq_{\nn\ssc\,}}|x_2|^{1+\d'}
\stackrel{{}^{\sc\rm\eqref{ansnen0000}}}{\sim_{\nn\ssc\,}}\nn^{1+\d'}
\stackrel{{}^{\sc\rm\eqref{ansnen0000}}}{\succ_{\nn\ssc\,}}|x_{\OnE}|.} Thus \equa{y2=a=a=a}{\mbox{$|y_{\OnE}|\stackrel{{}^{\sc\rm\eqref{x2y2-x2}}}{\sim_{\nn\ssc\,}}|x_{\OnE}+y_{\OnE}|
\stackrel{{}^{\sc\rm\eqref{x2y2-x2}}}{\succ_{\nn\ssc\,}}\nn^{1+\d'}.$}}
By \eqref{ansnen0000} and \eqref{mqp1234-2}, we have  \equa{hdhdhdhhhh}{\mbox{$\dH_{p_{\ZeRo},p_{\OnE}}
\stackrel{{}^{\sc\rm\eqref{dp0p1},\,\eqref{mqp1234-2}}}{\sim_{\nn\ssc\,}}
|x_{\ZeRo}|+|x_{\OnE}|
\stackrel{{}^{\sc\rm\eqref{ansnen0000}}}{\preceq_{\nn\ssc\,}}\nn^{1+\frac1m}
\stackrel{{}^{\sc\rm\eqref{y2=a=a=a}}}{\prec_{\nn\ssc\,}}|y_{\OnE}|^{\frac{m+1}{m}}$,}} a contradiction with
\eqref{mqp1234-2}  [by \eqref{MSnenenene},\,\eqref{mxxxxxx-ssss} and by Remark \ref{Remmma} (to be given later), one can see why we can get a contradiction].

This proves that \eqref{NSoOP} is satisfied by $(p_1,p_2)$, i.e., $(p_1,p_2)\in V_0$, and so $V_0$ is closed.
 Proposition
{\rm\ref{real00-inj}} holds. Case 3 is now completed.
\vskip7pt
Hence from now on, we can assume
\equa{amamsn}{\mbox{ $\dis\ssTH b_\kk\ge1+\frac1m
+ a_\kk$
\ \ when $\kk\gg\ell$.}}
}
{%
\vskip4pt
{\noindent\it Case 4: Assume there exists a fixed but sufficiently small $\d'\in\R_{>0}$ such that for any $\SS\in\R_{>0}$ there exists $\kk>\SS$
satisfying $(1-\d')b_\kk>1+a_\kk$ $($the ``extra fact''$)$.
}
\vskip4pt
First we remark that the reason we require $\d'$ to be independent of $\kk$ is to guarantee that we have \eqref{Amsnenenene}.
We define  $V_1$ to be the subset of $V$  consisting of    all
 elements $(p_{\ZeRo},p_{\OnE})\!=\!\big((x_{\ZeRo},y_{\ZeRo}),(x_{\OnE},y_{\OnE})\big)$  such that its coordinates $x_1,x_2,y_2$ satisfy  (cf.~Remark \ref{d-k-n}),
\begin{eqnarray}
\label{NSoMMMOP}
\!\!\!\!\!\!\!\!\!\!\!\!\!\!\!\!\!\!\!\!\!
&&
\dis{\rm(i)\ }
1<\big(\kk^{-1}|x_2|\big)^{-1+\kk^{-1}}\le(\kk^{-1}|x_1|)^{-1}\le
\big(\kk^{-1}|x_2|\big)^{-1-\kk^{-1}}<(1-\d')^{-1-\kk^{-1}},
\nonumber\\[-4pt]
\!\!\!\!\!\!\!\!\!\!\!\!\!\!\!\!\!\!\!\!\!\!\!\!\!\!\!\!
&&
 {\rm(ii)\ }\ell_{p_1,p_2}:=
\frac{\g_{\kk,\kk}^{-1}|x_2+y_2|}{\big(\kk^{-1}|x_2|\big)^{\ssTH 
1-\d'
}}+(|x_2|+|x_2+y_2|)\ep^3\ge1+\ep^2.\end{eqnarray}
As in the previous case, we have \eqref{ToSayas}, 
 \eqref{-EiathA0} and $V_1$ is bounded
.
\begin{rema}\rm\label{FirRemma+}\begin{itemize}\item[(i)]
In contrary to \eqref{NSoOP}, in the denominator of \eqref{NSoMMMOP}\,(ii), we put the power
$1-\d'$, which is smaller than $1$ and is a number independent of $\kk$
[that is important otherwise we cannot obtain \eqref{Amsnenenene}$\ssc\,$],
is to ensure that $|x_2+y_2|$ will decrease slower than $|x_2|$ when the last strict inequality of
\eqref{NSoMMMOP}\,(i) becomes  equality for an element in $\ol V_1$.
\item[(ii)]Note that the number in the last term of \eqref{NSoMMMOP}\,(i) cannot be too big otherwise
for an element in $\ol V_1$,
when the last strict inequality of
\eqref{NSoMMMOP}\,(i) becomes  equality,
$|x_2|$ will fall to a too small number and then we cannot use Theorem \ref{Theo-2} or  \eqref{mqp1234-2} (for example
 we cannot choose the number to be $\kk$).
\end{itemize}\end{rema}

With $(q_1,q_2)$ being defined in \eqref{1++??+q0q1-bar}, we want to  choose suitable $u,v$ such that
\eqref{NSoMMMOP}\,(i),\,(ii) hold for $(q_1,q_2)$, i.e.,
\begin{eqnarray}
\label{1NSoMMMOP}&&\!\!\!\!\!\!\!\!\!\!\!\!\!\!\!\!\!\!\!\!\!\!\!\!\!\!\!\!
\dis{\rm(i)\ }
1<|1+u\ep|^{-1+\kk^{-1}}\le|1+s\ep|^{-1}\le
|1+u\ep|^{-1-\kk^{-1}}<(1-\d')^{-1-\kk^{-1}},
\nonumber\\[-4pt]&&\!\!\!\!\!\!\!\!\!\!\!\!\!\!\!\!\!\!\!\!\!\!\!\!\!\!\!\!
{\rm(ii)\ }
\frac{\Big|1+v
\ep\Big|}{|1+u\ep|^{\ssTH 
1-\d'
}}+O(\ep)^3\ge1+\ep^2.
\end{eqnarray}
Take $u,v\in\R_{<0}$ such that [the last equation is obtained from \eqref{@suc2hthat=4}$\ssc\,$\vspace*{-2pt}],
\equa{LetF3rMMMst}{\mbox{$\dis u=-1,\ \ \ v=
-\frac{1+a_\kk}{b_\kk},$
\ \ \ \ \ and \ \ $\dis s\stackrel{{}^{\sc\rm\eqref{@suc2hthat=4}}}{=}-1+O(\ep)^1$.}}
Then the coefficients of $\ep^1$ in the middle three terms in \eqref{1NSoMMMOP}\,(i) are respectively $1-\kk^{-1},\,1,\,1+\kk^{-1}$, i.e.,
%
all inequalities  in \eqref{1NSoMMMOP}\,(i) are strict inequalities.
Further, the coefficient of $\ep^1$ in the left hand-side of \eqref{1NSoMMMOP}\,(ii) is
$\ssTH 
1-\d'
-\frac{1+a_\kk}{b_\kk}$, which is positive by the ``extra fact''. Thus $(q_1,q_2)\in V_1$
(for $\ep$  sufficiently small). We take $V_0=V_1\ne\emptyset.$

Similarly to \eqref{x2+ya2}, we obtain from \eqref{NSoMMMOP}\,(ii) the following (which also holds for any element in $\ol V_0$),
\equa{x2+ya2++++}{\dis
|x_2+y_2|>\g_{\kk,\kk}\big(\kk^{-1}|x_2|\big)^{\ssTH 
1-\d'
}.}

Let $(p_1,p_2)\in\ol V_0$. By Lemma \ref{lemm-condition-XZ}, first assume
the last strict inequality of \eqref{NSoMMMOP}\,(i) becomes  equality.
Then we   obtain that $|x_2|=(1-\d')\kk$, and from   the second inequality of  \eqref{NSoMMMOP}\,(i),
$|x_1|\le\kk$.
%
By \eqref{x2+ya2++++} and Lemma \ref{Anlemmmmmm}, we have \equa{yyy6666}{\mbox{$\dis|x_2+y_2|\stackrel{{}^{\sc\rm\eqref{x2+ya2++++}}}{>}(1-\d')^{\ssTH 
1-\d'
}\g_{\kk,\kk}\stackrel{{}^{\sc\rm Lemma\ \ref{Anlemmmmmm}}}{>}\Big(1-\d'+\d'^{2}+O(\d')^{3}\Big)
\kk\sTH$.}}
As before we have \equa{hp1ip2i}{\mbox{$h_{p_1,p_2}
\stackrel{{}^{\sc\rm\eqref{dp0p1},\,\eqref{mqp1234-2}}}{\sim_{\kk\ssc\,}}
|x_1|+|x_2|\sim_{\kk\ssc\,}\kk,$}}
and we obtain \equa{Amsnenenene}{\dis\mbox{
$|y_2|\ge|x_2+y_2|-|\TH x_2\sTH |\stackrel{{}^{\sc\rm\eqref{yyy6666}}}{>}\big(\d'^{2}+O(\d')^{3}\big)\kk\sim_{\kk\ssc\,}
\kk\sTH \succ_{\kk\ssc\,}
\kk^{\frac{\ssTH m}{\ssTH m+1}}\stackrel{{}^{\sc\rm\eqref{hp1ip2i}}}{\sim_{\kk\ssc\,}} h_{_{\sc p_1,p_2}}^{^{\sc\frac{\ssTH m}{\ssTH m+1}}}$,}} a contradiction with
Theorem \ref{Theo-2} or  \eqref{mqp1234-2}.

Before continuing the proof, at this point it may be worth stating the following remark.
\begin{rema}\label{Remmma}\rm
Let $a\in\R_{>0}$ be a parameter such that $a\gg1$. If for some $\d'\in\R_{>0}$ independent of $a$ with $\d'<\frac1m$, we have $|x_1|\preceq_{a\ssc\,}a^{1+\d'}$, and further either $|x_2|\succeq_{a\ssc\,}a$ or
$|x_2+y_2|\succeq_{a\ssc\,}a$
[thus $h_{p_1,p_2}\succeq_{a\ssc\,}|x_2|\succeq_{a\ssc\,}a$ or $h_{p_1,p_2}\succeq_{a\ssc\,}|x_2+y_2|\succeq_{a\ssc\,}a$], then \eqref{Amsnenenene} has in fact proved that $|x_2|$ cannot differ from
$|x_2+y_2|$ by a factor which is independent of $a$ [in case $|x_2|>|x_2+y_2|$, we use the inequation $|y_2|\ge|x_2|-|x_2+y_2|$ in \eqref{Amsnenenene}$\ssc\,$].
\end{rema}

Now assume the first strict inequality  of \eqref{NSoMMMOP}\,(i),
becomes  equality for $(p_1,p_2)$.
  Then $|x_1|=|x_2|=\kk$, but by
\eqref{x2+ya2++++}, $|x_2+y_2|>\g_{\kk,\kk}$, a contradiction with definition \eqref{Ak=1}.

This proves that \eqref{NSoMMMOP} is satisfied by $(p_1,p_2)$, i.e., $(p_1,p_2)\in V_0$, and so $V_0$ is closed.
 Proposition
{\rm\ref{real00-inj}} holds.   Case 4 is now completed.\vskip7pt

Hence from now on, we assume
\equa{(1-dp")}{(1-\d')b_\kk\le 1+a_\kk\mbox{ for any fixed sufficiently small $\d'>0$ and all $\kk\gg\ell$}.}
%
%
This with \eqref{amamsn}  shows that $a_\kk\ge\frac{(1-\d')(1+\frac1m)}{\d'}$. Since $\d'>0$ is arbitrarily sufficiently small number,
we see that $a_\kk>0$ (thus also $b_\kk>0$) is unbounded, i.e.,
\equa{Akk-bkk}{\dis\lim_{\kk\to\infty}a_\kk=\lim_{\kk\to\infty}b_\kk=\infty,\mbox{ \ \ \ and }\ \ \ \lim_{\kk\to\infty}\frac{a_\kk}{\ssTH b_\kk}=1.
}
}\vskip7pt

{\noindent\it Case 5: Assume there exist some fixed $\l,\,\l_2\in\R_{>0}$ $($i.e., $\l,\l_2$ are independent of $\kk)$ with $\l<1$ such that  whenever $\kk\gg1$ there exist
  $k_1,k_2,\hat k\in\R_{>0}$ $($which can depend on $\kk{\ssc\,})$
  with $k_1,k_2,\hat k<1$
and $(\hat p_1,\hat p_2)=\big((\hat x_1,\hat y_1),(\hat x_2,\hat y_2)\big)\in V$ satisfying the following,
\begin{eqnarray}\label{TheFasss}
&\!\!\!\!\!\!\!\!\!\!\!\!\!\!\!\!\!\!\!\!\!\!&
{\rm(i)\ }|\hat x_1|=k_1\kk,\ \ |\hat x_2|=k_2\kk,\ \ |\hat x_2+\hat y_2|=\hat  k\g_{\kk,\kk},\ \
\nonumber\\
&\!\!\!\!\!\!\!\!\!\!\!\!\!\!\!\!\!\!\!\!\!\!&
{\rm(ii)\ }\l\le k_2<1,\ \ \ \ \
{\rm(iii)\ }0<k_1< 1
,\ \
\ \ \
{\rm(iv)\ }0<\hat k<1,\ \ \ \ \ {\rm(v)\ }k_2< \hat  k^{
1+\l_2}
.
%
\end{eqnarray}
}

First we remark that the reason we require $\l_2$ to be independent of $\kk$ is to guarantee that we  have \eqref{Y1>+++} and the reason we require $\l$ to be independent of $\kk$ is to guarantee that $\hat x_2$ is large enough for us to apply Theorem \ref{Theo-2} [cf.~\eqref{Mandnd}$\ssc\,$].

Fix any sufficiently small $\d>0$ (independent of $\kk$) such that $\d\ll\min\{\l,
\l_2\}$.
First assume $k_2\le1-\d$. Then we have, when $\kk\gg1$,
\begin{eqnarray}
\!\!\!\!\!\!\!\!\!\!\!\!\!\!\!\!\!\!\!\!\!\!\!\!\!\!\!\!
&&
\label{Mandnd}
{\rm(i)\ }|\hat x_2|=k_2\kk\stackrel{{}^{\sc\rm\eqref{TheFasss}\,(ii)}}{>}\l\kk\sim_{\kk\ssc\,}\kk,
\ \ \  {\rm(ii)\ }|\hat x_1|=k_1\kk\stackrel{{}^{\sc\rm\eqref{TheFasss}\,(iii)}}{<}
\kk\sim_{\kk\ssc\,}\kk
,
\nonumber\\[4pt]\!\!\!\!\!\!\!\!\!\!\!\!\!\!\!\!\!\!\!\!\!\!\!\!\!\!\!\!&&
{\rm(iii)\ }
|\hat y_2|\ge|\hat x_2+\hat y_2|-|\hat x_2|\stackrel{{}^{\sc\rm\eqref{TheFasss}}}{\ge}
k_2^{\frac1{1+\l_2}}\g_{\kk,\kk}-k_2\kk\stackrel{{}^{\sc\rm Lemma\ \ref{Anlemmmmmm}}}{>}
(k_2^{\frac1{1+\l_2}}-k_2)\kk
\nonumber\\[0pt]\!\!\!\!\!\!\!\!\!\!\!\!\!\!\!\!\!\!\!\!\!\!\!\!\!\!\!\!&&
\phantom{{\rm(iii)\ }|\hat y_2|}
\sim_{\kk\ssc\,}\ \ \kk\ \, \stackrel{{}^{\sc\rm\eqref{Mandnd}(i),\,(ii),\,\eqref{dp0p1},\,\eqref{mqp1234-2}}}{\sim _{\kk\ssc\,}}\ \,h_{\hat p_1,\hat p_2},
%
\end{eqnarray}
 where
the strict inequality in (iii) follows from the fact in Lemma \ref{Anlemmmmmm} that $\g_{\kk,\kk}>\kk$,
and where the first ``\,$\sim_{\kk\ssc\,}$\,'' in (iii) follows from the fact that
 $\min\{x^{\frac1{1+\l_2}}-x\,|\,\l\le x\le1-\d\}$ is a fixed positive number (i.e., independent of $\kk$), and the
last ``\,$\sim_{\kk\ssc\,}$\,'' follows from  \eqref{Mandnd}\,(i),\,(ii) and \eqref{dp0p1},\,\eqref{mqp1234-2}.
We obtain from \eqref{Mandnd}
a contradiction with \eqref{mqp1234-2}. Thus
\equa{k2-is-big}{\mbox{$k_2>1-\d$.}}

%
%
Since $\hat k\ne1$, we can always write, for some $\l'_1,\l'_2\in\R_{>0}$%
,
\equa{WecanWSAS}{\dis {\rm(i)\ }k_1=\hat k^{\l'_1},\ \ \ \ {\rm(ii)\ }k_2=\hat k^{\l'_2},\mbox{ \ \ \ then \ \ }
{\rm(iii)\ }\l'_2\stackrel{{}^{\sc\rm\eqref{TheFasss}\,(v)}}{>}1+\l_2
,}
where the 
inequation follows from 
\eqref{TheFasss}\,(v).

Let $\ep_0\in\R_{>0}$ be such that $\ep_0^{-1}\gg\kk$. Now we define $V_0=V_1$ to be the subset of $V$ consisting of elements $(p_1,p_2)=\big((x_1,y_1),(x_2,y_2)\big)$ satisfying,
where $\ln(\cdot)$ is the  natural logarithm function,
\begin{eqnarray}
\label{LetNSoOP-add}
&\!\!\!\!\!\!\!\!\!\!\!\!\!\!
&
\dis{\rm(i)\  }
1< (\kk^{-1}|x_2|)^{-\frac{\l_1'}{\l'_2}+\scep_0}
\le(\kk^{-1}|x_1|)^{-1}
\le(\kk^{-1}|x_2|)^{-\frac{\l_1'}{\l'_2}-\scep_0}<(1-\d)^{-\frac{\l_1'}{\l'_2}-\scep_0},
\nonumber\\[0pt]&\!\!\!\!\!\!\!\!\!\!\!\!\!\!\!\!\!\!\!\!\!\!\!\!\!\!&
{\rm(ii)\ }
\frac{\g_{\kk,\kk}^{-1}|x_2+y_2|}{(\kk^{-1}|x_2|)^{\l_2'^{-1}(1+\scep_0)}}+(|x_2|+|x_2+y_2|)\ep_0^2\ge\kappa_2,
\\[0pt]\nonumber&\!\!\!\!\!\!\!\!\!\!\!\!\!\!\!\!\!\!\!\!\!\!\!\!\!\!&
{\rm(iii)\ }\kappa_2
:=\frac{\hat k}{\hat k^{1+\scep_0}}+(k_2\kk+\hat k\g_{\kk,\kk})\ep_0^2=1-\ln(\hat k)\ep_0+O(\ep_0)^2
\stackrel{{}^{\sc\rm\eqref{TheFasss}\,(iv)}}{>}1+(\kk+\g_{\kk,\kk})\ep_0^2
,\!\!\!\!
\end{eqnarray}where the inequality in \eqref{LetNSoOP-add}\,(iii) follows by noting that
$\ln(\hat k)<0$ since $0<\hat k<1$ by \eqref{TheFasss}\,(iv) (noting that we have complete freedom in choosing $\ep_0>0$ to be sufficiently smaller than  $\kk^{-1}$).
Then as before we can rewrite \eqref{LetNSoOP-add} as the form in
\eqref{ToSayas}, 
and \eqref{-EiathA0} holds and $V_0$ is bounded%
. 

Similar to \eqref{NSoMMMOP}, in the denominator of \eqref{LetNSoOP-add}\,(ii), we put the power
$\l_2'^{-1}(1+\ep_0)$, which by \eqref{WecanWSAS}\,(iii) is, up to $O(\ep_0)^1$, smaller than $\frac1{1+\l_2}$
(that is a number independent of $\kk$)
is to ensure that $|x_2+y_2|$ will decrease slower than $|x_2|$ when the last strict inequality of
\eqref{LetNSoOP-add}\,(i) becomes  equality for an element in $\ol V_1$.

When we set $(p_1,p_2)$ to be  $(\hat p_1,\hat p_2)$ in \eqref{LetNSoOP-add},
the middle three terms in  \eqref{LetNSoOP-add}\,(i) are respectively
\equa{Middle-three}{\dis
T_1:=k_{_{\sc2}}^{{\sc-\frac{\l'_1}{\l'_2}+\scep_0}},\ \ \ \ \
T_2:=k_1^{-1}\stackrel{{}^{\sc\rm\eqref{WecanWSAS}}}{=}\hat k^{-\l'_1}\stackrel{{}^{\sc\rm\eqref{WecanWSAS}}}{=}k_{_{\sc2}}^{-\frac{\l'_1}{\l'_2}},\ \ \ \ \
T_3=k_{_{\sc2}}^{-\frac{\l'_1}{\l'_2}-\scep_0}.
}
 Since $1-\d<k_2<1$ by \eqref{TheFasss}\,(ii),\,\eqref{k2-is-big}, we see from \eqref{Middle-three} that $1<T_1<T_2<T_3<(1-\d)^{-\frac{\l_1'}{\l'_2}-\scep_0}$ so that
\eqref{LetNSoOP-add}\,(i) holds for  $(\hat p_1,\hat p_2)$.
Further the left-hand side of
\eqref{LetNSoOP-add}\,(ii) is $\frac{\hat k}{\hat k^{1+\scep_0}}+(k_2\kk+\hat k\g_{\kk,\kk})\ep_0^2=\kappa_2$ by \eqref{TheFasss}\,(i), i.e.,
\eqref{LetNSoOP-add}\,(ii) holds for  $(\hat p_1,\hat p_2)$.
\NOUSE{
it simply reads as [recalling that $k_2>\d$ and noting from \eqref{WecanWSAS}
that $\l'_2<1$ thus $1-\d<k_2<k_2^{\l'_2}$],
\equa{ReAdAS}{\mbox{
$\dis{\rm(i)\ }(1-\d)^{\frac{\l_1'}{\l'_2}+\scep_0}\le (k_2^{\l'_2})^{\frac{\l_1'}{\l'_2}+\scep_0}
\le k_2^{\l'_1}\le (k_2^{\l'_2})^{\frac{\l_1'}{\l'_2}-\scep_0}\le1$, \ \ \ \ \ $\dis{\rm(ii)\ }
\frac{k_2^{\l'_2}+\ep_0^2}{k_2^{\l_2'+\scep_0}}\ge\kappa_2$.}}
}%
Thus $(\hat p_1,\hat p_2)\in V_0$.

Let  $(p_1,p_2)\in\ol V_0$. By Lemma \ref{lemm-condition-XZ}, first assume the last strict inequality of
\eqref{LetNSoOP-add}\,(i) becomes  equality for $(p_1,p_2)$.
Then
\equa{Thehehehdxb}{\mbox{ ${\rm(i)\ }|x_1|\le\kk,$
\ \ \ \ ${\rm(ii)\ }|x_2|=(1-\d)\kk$,}}
thus by \eqref{mqp1234-2},
\equa{MSMSMememem}{\mbox{$h_{p_1,p_2}
\stackrel{{}^{\sc\rm\eqref{dp0p1},\,\eqref{mqp1234-2}}}{\sim_{\kk\ssc\,}}
|x_1|+|x_2|\sim_{\kk\ssc\,}|x_2|\stackrel{{}^{\sc\rm\eqref{mqp1234-2}}}{\sim_{\kk\ssc\,}}|x_2+y_2|
\stackrel{{}^{\sc\rm\eqref{Thehehehdxb}\,(ii)}}{\sim_{\kk\ssc\,}}\kk.$}}

Now we will evaluate \eqref{Y1>+++} up to $O(\ep_0)^1$.
First we have
\begin{eqnarray*}\label{MEMEM}
&\!\!\!\!\!\!\!\!\!\!\!\!\!\!\!\!\!\!\!\!\!\!\!\!\!&
|x_2|{\sc\!}+{\sc\!}|x_2
{\sc\!}+{\sc\!}y_2|\stackrel{{}^{\sc\rm\eqref{Thehehehdxb},\eqref{Ak=1}}}{<}(1{\sc\!}-{\sc\!}\d)\kk
{\sc\!}+{\sc\!}\g_{|x_1|,|x_2|}\stackrel{{}^{\sc\rm
\eqref{Thehehehdxb},\,\eqref{wePPPP1+},\,\eqref{wePPPP1}}}{\le}
(1{\sc\!}-{\sc\!}\d)\kk{\sc\!}+{\sc\!}\g_{\kk,(1-\d)\kk}\stackrel{{}^{\sc\rm\eqref{wePPPP1}}}{\le}
\kk{\sc\!}+{\sc\!}\g_{\kk,\kk}.\!\!\!\!\!
\end{eqnarray*}
From this and \eqref{LetNSoOP-add}\,(iii), we obtain that $\kappa_2-(|x_2|+|x_2+y_2|)\ep_0^2>1$.
This with \eqref{LetNSoOP-add}\,(ii) and \eqref{Thehehehdxb}\,(ii) gives, up to $O(\ep_0)^1$
 [so $\ep_0$ in \eqref{LetNSoOP-add} is omitted], the first inequality below; the second inequality below follows from the inequation in \eqref{WecanWSAS}\,(iii); the third inequality follows from Lemma \ref{Anlemmmmmm}.
%
Thus up to $O(\ep_0)^1$,
we obtain,
\begin{eqnarray}
\label{Y1>+++}
&\!\!\!\!\!\!\!\!\!\!\!\!\!\!\!\!\!\!\!\!\!\!\!\!\!\!\!\!\!\!\!\!\!&
|x_2+y_2|\stackrel{{}^{\sc\rm\eqref{LetNSoOP-add}\,(ii),\,(iii),\,\eqref{Thehehehdxb}\,(ii)}}{>}(1-\d)^{\l_2'^{-1}}\g_{\kk,\kk}
\stackrel{{}^{\sc\rm\eqref{WecanWSAS}\,(iii)}}{>}(1-\d)^{\frac{1}{1+\l_2}}\g_{\kk,\kk}
\nonumber\\[0pt]&\!\!\!\!\!\!\!\!\!\!\!\!\!\!\!\!\!\!\!\!\!\!\!\!\!\!\!\!\!\!\!\!\!\!\!\!\!&\phantom{|x_2+y_2|}
\ \ \ \ \ \
\stackrel{{}^{\sc\rm Lemma\ \ref{Anlemmmmmm}}}{\ge}(1-\d)^{\frac{1}{1+\l_2}}\kk=\Big(1-\frac{\d}{1+\l_2}+O(\d)^2\Big)
\sTH{\ssc\,}\kk\!\!\!\!\!\!\!\!\!\!\!\!\!\!\!\!\!\!\!\!\!\!
\\[0pt]\nonumber&\!\!\!\!\!\!\!\!\!\!\!\!\!\!\!\!\!\!\!\!\!\!\!\!\!\!\!\!\!\!\!\!\!\!\!\!\!&\phantom{|x_2+y_2|}
\ \ \ \ \ \ \ \  \stackrel{{}^{\sc\rm\eqref{amdnen333}}}{=} \ \ \ \Big(\big(1+\a\d+O(\d)^2\big)(1-\d)\Big)\kk
=
\big(1+\a\d+O(\d)^2\big)| x_2|\mbox{ \ [up to $O(\ep_0)^1$]},\!\!\!\!\!\!\!\!\!\!\!\!\!\!\!\!
\end{eqnarray}
where
\equa{amdnen333}{\mbox{$\dis \a=\frac{\l_2}{1+\l_2}>0$, which is a number independent of $\kk$.}}
By \eqref{Thehehehdxb},\,\eqref{Y1>+++}, using arguments as in \eqref{yyy6666},\,\eqref{Amsnenenene},
we obtain a contradiction (cf.~Remark \ref{Remmma}). 

Now assume the first strict inequality of
\eqref{LetNSoOP-add}\,(i) becomes  equality for $(p_1,p_2)$.
Then $|x_1|=|x_2|=\kk$. Thus $|x_2+y_2|\le\g_{|x_1|,|x_2|}=\g_{\kk,\kk}$ by definition \eqref{Ak=1}. We have
\begin{eqnarray}\label{Finsmemem}
\!\!\!\!\!\!\!\!\!\!\!\!\!\!\!\!\!\!\!\!\!\!\!\!\!\!\!\!\!\!&&
1+(\kk+\g_{\kk,\kk})\ep_0^2\stackrel{{}^{\sc\rm\eqref{LetNSoOP-add}\,(iii)}}{<}
\kappa_2\stackrel{{}^{\sc\rm\eqref{LetNSoOP-add}\,(ii)}}{\le}\frac{\g_{\kk,\kk}^{-1}|x_2+y_2|}{(\kk^{-1}|x_2|)^{\l_2'^{-1}(1+\scep_0)}}+(|x_2|+|x_2+y_2|)\ep_0^2
\nonumber\\
\!\!\!\!\!\!\!\!\!\!\!\!\!\!\!\!\!\!\!\!\!\!\!\!\!\!\!\!\!\!&&
\phantom{1+(\kk+\g_{\kk,\kk})\ep_0^2}\ \ \ \ \,
=\ \ \ \ \g_{\kk,\kk}^{-1}|x_2+y_2|+(\kk+|x_2+y_2|)\ep_0^2\le1+(\kk+\g_{\kk,\kk})\ep_0^2,
\end{eqnarray}
which is a contradiction. 
\NOUSE{%
Assume the first equality of \eqref{LetNSoOP-add}\,(i) holds.
Then $|x_1|=(1-\d)\kk$, $|x_0|\le\kk$, but by \eqref{LetNSoOP-add}\,(i), we have
(note that  $\d\ll\d_0\ssc\,$)
\begin{eqnarray}
\!\!\!\!\!\!\!\!\!\!\!\!\!\!&&
\label{yYYaha}
|y_1|\ge(1-\d)^{\frac{1+\d^4}{1+\d_1}}\g_{\kk,\kk}=\Big(1-\frac{\d}{1+\d_1}+O(\d^2)\Big)\g_{\kk,\kk}\nonumber\\[4pt]
\!\!\!\!\!\!\!\!\!\!\!\!\!\!&&
\phantom{|y_1|}
>
\Big(1-\frac{\d}{1+\d_0}+O(\d^2)\Big)\kk=\Big(1+\frac{\d_0}{1+\d_0}{\ssc\,}\d+O(\d^2)\Big)|x_1|,\end{eqnarray}
 [since $|x_0|\preceq\kk\preceq|y_1|$, we must have $h_{p_1,p_1}\sim|x_1|\sim|y_1|$ by \eqref{mqp1234-2+}$\ssc\,$].
 Thus Theorem \ref{real00-inj}\,(1)\,(i) holds.
}%
%

This proves that \eqref{LetNSoOP-add} is satisfied by $(p_1,p_2)$, i.e., $(p_1,p_2)\in V_0$, and so $V_0$ is closed.
 Proposition
{\rm\ref{real00-inj}} holds.   Case 5 is now completed.\vskip7pt

Thus from now on we assume that Case 5 does not occur.
\vskip5pt
{\noindent\it Case 6
: The remaining case.}
\vskip4pt
In this case, unfortunately, we are not able to use \eqref{ToSayas} to define $V_0$. We have to define $V_0$ in a more complicated way. The reason is below:
\begin{rema}\rm\label{Antotm}
 Because $b_\kk$ is too big [by \eqref{amamsn}$\ssc\,$], we are unable to choose suitable ``initial stage'' $(q_1,q_2)$ mentioned in Remark \ref{AboutX1}\,(ii) such that $u,\tildev ,s$ satisfy \eqref{@suc2hthat=4}, and $(q_1,q_2)$, defined in
\eqref{1++??+q0q1-bar}, satisfies \eqref{NSoOP} if we define $V_0$ as in \eqref{NSoOP}. To see this, observe that in order for \eqref{NSoOP}\,(ii) to hold, we
 have to choose $\tildev $ to be bigger than  $u$ but then \eqref{@suc2hthat=4} shows that $s$ may  possibly be too big, which
 implies that $|x_1|$ may grow too fast by 
Fact \ref{fact-initial}%
. Since we do not have any information
about $y_1$ in our definition of $V_0$, when $|x_1|\gg|x_2|$ we can not apply Theorem \ref{Theo-2}
to obtain a contradiction if the last strict inequality of \eqref{GMSMSMSMS??}\,(i) becomes  equality for an element in $\ol V_0$.
Therefore in order to able to obtain a contradiction,
we have to find some other way to control the growth of $|x_1|$.
\NOUSE{ in the following our strategy is to choose $\tildev $ to be smaller
 than $u$ [cf.~\eqref{LetF3rst-final}$\ssc\,$], then we can obtain,
 \begin{itemize}\item[(a)] on one hand,  in contrary to \eqref{x2+ya2}, $|x_2+y_2|$ or $|Z|$ will grow slower or decrease faster than $|x_2|$ (or $|X_2|$)
 by Fact \ref{fact-initial}
;
 \item[(b)]
on other other hand, 
\eqref{equa-Case6-lemm}\,(iv) shows that $|Z|$ must grow or decrease at the same speed as $|X_2|$.
\end{itemize}
Therefore we obtain an inconsistence, which allows us to prove that
 the last inequality of \eqref{LetNSoOP}\,(i) cannot become an equality for any element in $\ol V_0$.
\end{itemize}
}\end{rema}

Recall notations \eqref{MSmde33333},\,\eqref{SimMMSMS} and \eqref{LetNSoOP----1} and Remark \ref{rema3.1}; in particular, we  will frequently use the fact 
that $0<\ep\ll
\kk^{-1}\ll\d=\ell^{-1}
\ll\d_2=\ell_2^{-1}
\ll\d_1=\ell_1^{-1}%
\ll\d_0=\ell_0^{-1}
\ll1$
.
\NOUSE{ We will need the following ``extra fact''
[recalling from \eqref{ll-dd} that $\dD=\lL^{-1}\le b_\kk^{-1}$],
\equa{b-0-bigger}{\dis\b_0\stackrel{{}^{\sc\rm\eqref{LetNSoOP----1}\,(iv)}}{=}
\eE_1\tilde s+(1-\dD)b_\kk-a_\kk\stackrel{{}^{\sc\rm\eqref{ll-dd}}}{\ge}
\eE_1\tilde s+b_\kk-1-a_\kk\stackrel{{}^{\sc\rm\eqref{amamsn}}}{\ge}
\frac{1}{m}+O(\eE_1)^1>0.
}
}%
\NOUSE{\begin{eqnarray}
\label{extra-alp}
&\!\!\!\!\!\!\!\!\!\!\!\!\!\!\!\!\!\!\!\!\!\!&
\a_0=b_\kk+\d>0. \ \ \
\end{eqnarray}
}\NOUSE{
With $a_1,a_3,c_1,c_2\in\Q_{>0}$ being defined in \eqref{LetNSoOP----1}\,(iv),\,(v) and \eqref{B123X12Z}\,(iv),\,(v), we also define
[noting that we define below in order to obtain \eqref{X0sm,m+1}$\ssc\,$]%
,
\begin{eqnarray}
\label{LetNSoOP----1**}
&\!\!\!\!\!\!\!\!\!\!\!\!\!\!\!\!\!\!\!\!\!\!\!\!\!\!\!\!\!\!\!\!\!\!
&
a_2=
\frac{c_1^{-1}+\b_1^{-1}}{1 + 2 \d_0- \d_0^4}=\b_1^{-1}\big(1+O(\d_0)^1\big)
.
\end{eqnarray}
}%
\begin{defi}\rm\label{V2=?}
Recall that $A_1,A_2
$ are 
the following functions defined in 
\eqref{LetNSoOP----1}
:
\begin{eqnarray}
\label{LetNSoOP----1-re-give}
&\!\!\!\!\!\!\!\!\!\!\!\!\!\!\!\!\!\!\!\!\!\!\!\!\!\!\!\!\!\!\!\!\!\!\!\!\!\!\!\!\!\!\!\!\!\!\!\!\!\!\!
&
{\rm(i)\ }A_1
=
X_2^{\ell_0} \Big(2-\frac1{\widetilde X_1}\Big),
\ \ \ \ \ \
{\rm(ii)\ }A_2
=
\frac{X_2^{\ell_0+1}}{Z}\Big (2-\frac{\widetilde  X_1}{A_1^2}\Big)
.
\!\!\!\!\!\!\!\!\!\!\!\!\!\!\!\!\!\!\!\!\!\!\!\!
\end{eqnarray}
\NOUSE{
\begin{eqnarray}
\label{LetNSoOP----1-redefine}
&\!\!\!\!\!\!\!\!\!\!\!\!\!\!\!\!\!\!\!\!\!\!\!\!\!\!\!\!\!\!\!\!\!\!\!\!\!\!\!\!\!\!\!
&
{\rm(i)\ }A_3=\frac{ X_2^{75\ell_0+52}}{\widetilde X_1^{60\ell_0+8} Z^{75\ell_0+53}},
\ \ \ \ \ \ \ \
{\rm(ii)\ }
A_1=\frac{A_3^6\widetilde X_1^{ 25\ell_0+10} Z^{25}}{X_2^{20}}(2 -\widetilde X_1^{10}),
\nonumber\\
&\!\!\!\!\!\!\!\!\!\!\!\!\!\!\!\!\!\!\!\!\!\!\!\!\!\!\!\!\!\!\!\!\!\!\!\!\!\!\!\!\!\!\!\!\!\!\!\!\!\!\!
&
{\rm(iii)\ }A_2=
\frac{Z^{20\ell_0}}{\widetilde X_1^{25\ell_0-6}X_2^{25\ell_0}}\Big(\frac85 -\frac{3\widetilde X_1^{10}}5\Big)
.
\end{eqnarray}
}
{
}\NOUSE{\equa{a1a2a3}{\dis
a_1=-\ell_0^2(1- 3 \d_0 + 2 \d_0^4 - \d_0^5),\ \
a_2=-\d_0^2 (1- \d_0 + \d_0^4),
\ \ \
a_3=\d_0(1-\d_0).}
Now}%
We take
 $V_0=V_2$ with \eqref{ToSayas+1} 
being specified as follows%
, 
%
%
\begin{eqnarray}
\label{LetNSoOP}
\!\!\!\!\!\!\!\!\!\!\!\!\!\!\!\!\!\!\!\!\!&\!\!\!\!\!\!\!\!\!\!\!\!\!\!\!\!\!\!\!\!\!\!\!\!\!\!\!\!&
{\rm(i)\ }1< |A_{\rOnE}|^{-1+\d^2
}
\le |A_2|^{-1}
\le|A_1|^{-1-\d^2
}
<
\ell_{_{\sc 1}}^{^{\sc1+\d^2}}
,
\nonumber\\[0pt]
\!\!\!\!\!\!\!\!\!\!\!\!\!\!\!\!\!\!\!\!\!\!\!\!&\!\!\!\!\!\!\!\!\!\!\!\!\!\!\!\!\!\!\!\!\!\!\!\!\!\!\!\!&
{\rm(ii)\ }\ell_{p_1,p_2}{\ssc}:={\ssc}
|X_2|\cdot|A_1|^{-\d^2
}
{\ssc}+{\ssc}(|x_{\OnE}|{\ssc}+{\ssc}|x_{\OnE}{\ssc}+{\ssc}y_{\OnE}|)\ep_2^3{\ssc}\ge
{\ssc}1{\ssc}+{\ssc}\ep_2^2
,
\nonumber\\[0pt]
\!\!\!\!\!\!\!\!\!\!\!\!\!\!\!\!\!\!\!\!\!\!\!\!&\!\!\!\!\!\!\!\!\!\!\!\!\!\!\!\!\!\!\!\!\!\!\!\!\!\!\!\!&
{\rm(iii)\ }
(1-\d)|X_2|\cdot|A_1|^{2+\d_0}
\le|\widetilde X_1|\le
(1+\d)|X_2|^{-1}\cdot|A_1|^{-\d_0}
.\!\!\!\!\!\!\!\!\!\!\!\!\!\!\!\!\!\!\!\!
\end{eqnarray}
%
\end{defi}
\begin{rema}\rm\label{Rema-V2-def}
{The reason we put $\pm\d^2$ in the powers of the second and fourth terms of \eqref{LetNSoOP}\,(i) is to control
$|A_2|$ 
such that
 we have the following,
which holds for any element in $\ol V_0$
,
\equa{cont-B2}{\dis|A_2|
\stackrel{{}^{\sc\rm\eqref{LetNSoOP}\,(i)}}{=}
|A_1|^{1+O(\d)^{2}
}
=
|A_1|
+O(\d)^2
,
}
where the last equality follows from the fact that
$\d=\ell^{-1}\ll\d_1$%
%
%
%
.
Further, by }
\eqref{LetNSoOP}\,(i)
~and Lemma \ref{lemm-condition-XZ},
for any $(p_1,p_2)\in\ol V_0$, we have
,
\begin{eqnarray}
\label{A1-A2-cond}
&&\!\!\!\!\!\!\!\!\!\!\!\!\!\!\!\!\!\!\!\!\!\!\!\!\!\!\!\!\!\!\!\!\!\!\!\!\!\!\!\!
{\rm(a)\ }
\d_1\le
|A_1|
\le1
,
\ \ \ \ \ \ \
{\rm(b)\ }
\d_1+O(\d)^2\le
|A_2|\le
1
.\!\!
\!\!\!\!\!\!\!\!\!\!\!\!
\end{eqnarray}
\NOUSE{
\item[(ii)]
Observe from \eqref{LetNSoOP----1-redefine}\,(ii),\,(iii), we have the following
important relation between $A_1$ and $A_2$,
\equa{A2-A1-rela}{\dis
A_2=
\frac{Z^{20\ell_0}}{\widetilde X_1^{25\ell_0-6}X_2^{25\ell_0}}
\Big(\frac25 + \frac{3A_1X_2^{20}}{5A_3^6\widetilde X_1^{ 25\ell_0+10} Z^{25}}\Big).
}
\end{itemize}
}\NOUSE{\item[(ii)]
The reason we use $\ell_1$ to
control $A_1,A_2
$ in \eqref{A1-A2-cond} and use $\ell_2$ to control $\tildeX_1,X_2
$  in the definition of $S_1$ in \eqref{meme}
is to obtain
 Lemmas \ref{Case6-lemm},\,\ref{Step1} and \eqref{LetNSoOP----2},
the importance of which
is that we in particular  have \eqref{equa-Case6-lemm}\,(iv), that  is extremely important to us in the following sense.
\begin{itemize}\item[(a)]Firstly, we can obtain \eqref{LetNSoOP----2} that means that $A_1,A_2,A_3$,
which are originally  functions of three variables $\tildeX_1,X_2,Z$, become  functions of two
variables $\tildeX_1,X_2$, up to $O(\d)^2$ (which, as we shall see, will not affect our arguments below). This  allows us to obtain \eqref{ImMpP}\,(2)--(4), which says that  no equality can occur in any inequality of
\eqref{C+ToSayas+1}\,(c) 
for an element in $\ol V_0$.
 \item[(b)]Secondly, as mentioned in Remark \ref{AboutX1}\,(iv) and Remark \ref{Antotm}\,(ii),
for an element in $\ol V_0$, when the last strict inequality of \eqref{LetNSoOP}\,(i) becomes an equality,
the key point we can  obtain a contradiction is the following: we choose $|Z|$ to be smaller than $|X_2|$ at the ``initial stage'', however
 \eqref{equa-Case6-lemm}\,(iv) shows that $|X_2|$ must grow or decrease at the same speed as $|Z|$, that, as mentioned in Remark \ref{Antotm}\,(ii), will give us some
 inconsistence
.
 \end{itemize}
\NOUSE{
\item[(iii]
Using notations \eqref{TaKa},\, \eqref{SimMMSMS},
we see from \eqref{LetNSoOP}\,(vi),\,(vii) that \eqref{-EiathA0}\,(1) holds.
\NOUSE{\eqref{LetNSoOP}\,(v),\,(vi)
and  the fact from \eqref{LetNSoOP}\,(i) that $1\le|A_1|\le
\ell_1$
}
We will see after Lemma \ref{Step1} that \eqref{-EiathA0}\,(2),\,(3) also hold.
}
\end{itemize}
}\end{rema}
{%
Now we divide the proof of Proposition \ref{real00-inj} into
three
lemmas.
}%

{%
\begin{lemm}\label{V0NOT0}The set $V_{\rZeRo}$ is nonempty.
\end{lemm}
\noindent{\it Proof.~}First at this point it may be worth recalling that we have fixed $(\bar p_1,\bar p_2)=\big((\bar x_1,\bar y_1),(\bar x_2,\bar y_2)\big)\in A_{\kk,\kk}\subset V$ satisfying \eqref{TaKa}, i.e.,
\equa{TaKa++}{\mbox{$|\bar x_{\ZeRo}|=|\bar x_{\OnE}|=\kk
$, \ \ \ $|\bar x_{\OnE}+\bar y_{\OnE}|=\g_{\kk,\kk}.$}}
We need to choose a suitable ``initial stage''  $(q_1,q_2)$ as mentioned in Remark \ref{AboutX1}\,(ii). To do this, as already mentioned in \eqref{1++??+q0q1-bar},
we choose  $(q_1,q_2)=\big((\dot x_1,\dot y_1),(\dot x_2,\dot y_2)\big)\in V_0$ to be sufficiently close to $(\bar p_1,\bar p_2)$ such that
\eqref{1++??+q0q1-bar} holds, i.e.,
\equa{1++??+q0q1-bar++}{\!\!\!\!\!\!\!\!\!\!\!\!
\dot x_{\ZeRo}=\bar x_{\ZeRo}(1{\ssc}+{\ssc} s\ep),\ \ \
\dot y_1=\bar y_{\ZeRo}{\ssc}+{\ssc}t\ep,\ \ \
\dot x_{\OnE}=
\bar x_{\OnE}(1{\ssc}
+{\ssc}u
\ep),\ \ \
\dot x_2{\ssc}+{\ssc} \dot y_{\OnE}{\ssc}={\ssc}
(\bar x_2{\ssc}+{\ssc}\bar y_{\OnE})(1{\ssc}+{\ssc}v\ep),
\!\!\!\!\!\!}
and further, $s$ is determined by \eqref{@suc2hthat=4}, i.e.,
\equa{@suc2hthat=4++}{
s=s_0+O(\ep)^1,\ \ \ \ \ s_0=-a_\kk u{\sc}+{\sc}b_\kk v
,}
where $a_\kk\gg1,\,b_\kk\gg1$ satisfy \eqref{Akk-bkk}.
\NOUSE{
Note that $(1+2\d_0)a_\kk-(1+\d_0)b_\kk
\gg-(1+2\d_0)$ by
\eqref{amamsn},\,\eqref{Akk-bkk}
, we then
define
\begin{eqnarray}
&\!\!\!\!\!\!\!\!\!\!\!\!\!\!\!\!\!\!\!\!\!\!\!\!\!\!\!\!\!\!\!\!\!\!\!\!\!\!\!\!\!\!\!\!\!\!\!&
\label{a0BiggerThen}
{\rm(i)\ }\b_0:=
1+2\d_0+(1+2\d_0) a_\kk
-(1+\d_0)b_\kk
\stackrel{{}^{\sc\rm\eqref{Akk-bkk}}}{>}0,\ \ \ \ \ \
{\rm(ii)\ }\a_0{\ssc}:={\ssc}1{\ssc}+{\ssc}\b_0\ep{\ssc}>{\ssc}1.\!\!\!\!\!\!\!\!\!\!\!\!\!\!\!\!\!\!\!\!
\NOUSE
{\nonumber\\
&\!\!\!\!\!\!\!\!\!\!\!\!\!\!\!\!\!\!\!\!\!\!\!\!\!\!\!&
{\rm(iii)\ }\b_1=9\ell_0-2-2\d_0^2
,\ \ \ \ \ \ {\rm(iv)\ }\a_1=1+\b_1\ep
.}
\end{eqnarray}
The reason we define \eqref{a0BiggerThen}\,(i),\,(ii) is in order to obtain \eqref{X1X2Z===}\,(i).
We also remark that the fact that
 $\a_0>1$  is the ``extra fact'' we need that will be extremely important in obtaining the first
 inequality of
\eqref{664646}\,(ii)
.
}%
We want to take suitable $u,v$ such that
\eqref{LetNSoOP} is satisfied by $(q_{\ZeRo},q_{\OnE})$, defined in
\eqref{1++??+q0q1-bar}.
{
{
Note from 
\eqref{SimMMSMS} and \eqref{1++??+q0q1-bar++} that setting $(p_{\ZeRo},p_{\OnE})$ to $(q_{\ZeRo},q_{\OnE})$
corresponds to that $X_{\ZeRo},X_{\OnE},Z$ are respectively set to,
\begin{eqnarray}
\!\!\!\!\!\!\!\!\!\!\!\!\!\!\!\!\!\!\!\!&&
\label{XsXXX}
{\rm(i)\ }X_1{\sc}\stackrel{{}^{\sc\rm\eqref{SimMMSMS},\, \eqref{1++??+q0q1-bar++}}}{=}
{\sc}
{1{\sc}+{\sc}s\ep}
,\ \ \ {\rm(ii)\ }
X_2{\sc}\stackrel{{}^{\sc\rm\eqref{SimMMSMS},\, \eqref{1++??+q0q1-bar++}}}{=}{\sc}1{\sc}+{\sc}u\ep
,
 \ \ \ 
{\rm(iii)\ }Z{\sc}
\stackrel{{}^{\sc\rm\eqref{SimMMSMS},\, \eqref{1++??+q0q1-bar++}}}{=}{\sc}1{\sc}+{\sc}v\ep
.\!\!\!\!\!\!
\end{eqnarray}
%

At this point, we also want to remind from  notation \eqref{MSmde33333} that we  have complete freedom in
choosing the parameter $\ep
$ with $0<\ep\ll
\kk^{-1}
$ independently of all other choices of
the parameters $\ell_0,\ell_1,\ell_2,\ell,\kk$%
.
Further, it may be worth presenting the following obvious facts, for any $a,b\in\C$ which are bounded by some number independent of $\ep$,
\equa{ep-pro}{(1+a\ep)^{-1}=1-a\ep+O(\ep)^2,\ \ \ \ \ \
(1+a\ep)(1+b\ep)=1+(a+b)\ep+O(\ep)^2.}

Now we take
[we remark again that, as mentioned in Remark \ref{Antotm}\,(ii), we define $\tildev $ to be
smaller than $u$, which allows us to obtain a contradiction if
  equality occurs in the last inequality of \eqref{LetNSoOP}\,(i)
for an element in $\ol V_0$]
,
\begin{eqnarray}
\label{LetF3rst-final}
&\!\!\!\!\!\!\!\!\!\!\!\!\!\!\!\!\!\!\!\!\!\!\!\!\!\!\!\!\!\!\!\!\!\!\!&
{\rm(i)\ }u{\sc}={\sc}0,  \ \ \ \ \
{\rm(ii)\ }v
{\sc}=-\d^3,\  \ \ \ \  {\rm(iii)\ }s=s_0+O(\ep)^1\mbox{  \ with \ }
s_0\stackrel{{}^{\sc\rm\eqref{@suc2hthat=4++}}}{=}-\d^3b_\kk.
\end{eqnarray}
{%
So
[the following is the  reason we define $\b_0$ in \eqref{tX1==}$\ssc\,$],
\begin{eqnarray}
\label{so-that}\dis
&&\!\!\!\!\!\!\!\!\!\!\!\!\!\!\!\!\!\!\!\!\!\!\!\!\!\!\!\!\!\!
\b_0+s_0\stackrel{{}^{\sc\rm\eqref{tX1==}}}{=}
-1+\d^3b_\kk-\d^3b_\kk=-1
.\end{eqnarray}
}%
{%
Then by expanding in $\ep$ and recalling that $\ell_0\d_0=1$%
, up to $O(\ep)^2$,
we have, for $i=1,2
$
,
%
%
{%
\begin{eqnarray}
&\!\!\!\!\!\!\!\!\!\!\!\!\!\!\!\!\!\!\!\!\!\!\!\!\!\!\!&
\label{B1-ep}\label{X1X2Z===}
{\rm(i)\ }\widetilde X_1
=
1+\tilde s\ep+O(\ep)^2,\ \ \ \
\tilde s\stackrel{{}^{\sc\rm\eqref{Tx-B0}}}{=}
\b_0+s_0
\stackrel{{}^{\sc\rm\eqref{so-that}}}{=}
-1,
\ \ \ \
{\rm(ii)\ }
A_i=1+a_i\ep+O(\ep)^2
,
\!\!\!\!\!\!\!\!\!\!\!\!\!\!\!\!\!
\nonumber\\
&\!\!\!\!\!\!\!\!\!\!\!\!\!\!\!\!\!\!\!\!\!\!\!\!\!\!\!&
 {\rm(iii)\ }a_1
\stackrel{{}^{\sc\rm\eqref{recallll}
}}{=}
-1,
\!\!\!\!\!\!\!\!\!\!\!\!\!\!\!\!
\\\nonumber
&\!\!\!\!\!\!\!\!\!\!\!\!\!\!\!\!\!\!\!\!\!\!\!\!\!\!\!&
{\rm(iv)\ }a_2\stackrel{{}^{\sc\rm\eqref{LetNSoOP----1-re-give}\,(ii)}}{=}
(\ell_0+1)u-v-(-2a_1+\tilde s)
=
0+\d^3-(2-1)=-1+\d^3
.\!\!\!\!\!\!\!\!\!\!\!\!\!\!\!\!\!\!\!
\end{eqnarray}
where \eqref{X1X2Z===}\,(iii)  is precisely computed, up to $O(\ep)^2$, by recalling \eqref{LetNSoOP----1}\,(i),
}%
\begin{eqnarray}
\label{recallll}
\!\!\!\!\!\!\!\!\!\!\!\!\!\!\!\!\!\!\!\!\!\!\!\!\!\!\!\!\!\!\!\!\!\!\!\!\!\!\!\!\!\!\!\!\!\!\!\!&&
A_1\stackrel{{}^{\sc\rm\eqref{LetNSoOP----1-re-give}\,(i)}}{=}
 X_2^{\ell_0} \Big(2{\sc}-{\sc}\frac1{\widetilde X_1}\Big)
{\sc}={\sc}
(1{\sc}+{\sc}u\ep)^{\ell_0}
 \Big(2{\sc} -{\sc} (1{\sc}+{\sc}\tilde s\ep)^{-1}\Big)
\nonumber\\\!\!\!\!\!\!\!\!\!\!\!\!\!\!\!\!\!\!\!\!\!\!\!\!\!\!\!\!\!\!\!\!\!\!\!\!\!\!&&
\phantom{A_1\ \ \ \ }
{\sc}={\sc}1{\sc}+{\sc}\Big(\ell_0u{\sc}+\tilde s\Big)\ep
{\sc}={\sc}1{\sc}-{\sc}\ep{\sc}
\stackrel{{}^{\sc\rm\eqref{X1X2Z===}\,(iii)}}{=}{\sc}1{\sc}+{\sc}a_1\ep
.
\!\!\!\!\!\!\!\!\!\!\!\!\!\!\!
\end{eqnarray}
Similarly, we can easily compute \eqref{recallll}\,(iv)
.

\NOUSE{%
Then by \eqref{XsXXX},\,\eqref{ep-pro},\,\eqref{LetF3rst-final}, $X_1,X_2,Z$ are $1+O(\ep)^1$ elements of the following forms
\equa{X1X2Z===}{\dis X_1=1+s_1\ep+O(\ep)^2,\ \ \ \ \ X_2=1+u_1\ep+O(\ep)^2,\ \ \ \ \ Z=1+v_1\ep+O(\ep)^2,}
such that the coefficients of $\ep^1$ in them can be easily computed as follows,
\begin{eqnarray}\label{suv=====}
&\!\!\!\!\!\!\!\!\!\!\!\!\!\!\!\!\!\!\!\!\!\!\!\!&
{\rm(i)\ }s_1\stackrel{{}^{\sc\rm\eqref{XsXXX}\,(i)}}{=}\ell_0^2\ell(\a_0+s_0)
\stackrel{{}^{\sc\rm\eqref{@suc2hthat=4++},\,\eqref{a0BiggerThen},\,\eqref{LetF3rst-final}}}{=}
1-\d_0^6, \ \ \ \ \ \ \ {\rm(ii)\ }u_1\stackrel{{}^{\sc\rm\eqref{XsXXX}\,(ii)}}{=}\ell_0^2u\stackrel{{}^{\sc\rm\eqref{LetF3rst-final}\,(i)}}{=}-(1-\d_0^6),\!\!\!\!\!\!\!\!\!\!\!\!\!\!\!\!
\nonumber\\[4pt]
&\!\!\!\!\!\!\!\!\!\!\!\!\!\!\!\!\!\!\!\!\!\!\!\!&
{\rm(iii)\ }v_1\stackrel{{}^{\sc\rm\eqref{XsXXX}\,(iii)}}{=}\ell_0^2v\stackrel{{}^{\sc\rm\eqref{LetF3rst-final}\,(ii)}}{=}-1.
\end{eqnarray}
}\NOUSE{
 \eqref{XsXXX} gives that $X_1,X_2,Z$ are all $1+O(\ep)^1$ elements
 such that the coefficients of $\ep^1$ in them are respectively,
  \begin{eqnarray}
\label{W-XXXX}
&\!\!\!\!\!\!\!\!\!\!\!\!\!\!\!\!\!\!\!\!&
\tilde s=-\ell_0\ell(\a_0+s_1)\stackrel{{}^{\sc\rm\eqref{@suc2hthat=4++},\,\eqref{LetF3rst-final}}}{=}-4(1+\d_0),\ \ \ \ \ \ \tilde u=0,\ \ \ \ \ \ \tilde v=\ell_0v=-24.
\end{eqnarray}
}%
}
\NOUSE{
Now in \eqref{LetNSoOP}, we take  $(p_{\ZeRo},p_{\OnE})$ to be $(q_{\ZeRo},q_{\OnE})$.
Then  by expanding in $\ep$, we have, where $e$ is the natural number [recalling that $\nn=\ep^{-1}$,
and noting that when we take $X_1,X_2,Z$ to be as in \eqref{XsXXX} with \eqref{LetF3rst-final},
all $A_1,A_2,A_3$ are positive numbers up to $O(\ep)^2$ (though $A_1,A_2,A_3$ are not necessarily real numbers)$\ssc\,$],
\begin{eqnarray}
\label{Mnnwnw}
&\!\!\!\!\!\!\!\!\!\!\!\!\!\!\!\!\!\!\!\!\!\!\!\!\!\!&
{\rm(i)\ }
A_1^{\nn}\stackrel{{}^{\sc\rm\eqref{LetNSoOP----1}\,(i)}}{=}X_2^{\nn}
=(1+\ep)^{\ep^{-1}}=e+O(\ep)^1,
\nonumber\\
&\!\!\!\!\!\!\!\!\!\!\!\!\!\!\!\!\!\!\!\!\!\!\!\!\!\!&
{\rm(ii)\ }
A_2\stackrel{{}^{\sc\rm\eqref{LetNSoOP----1}\,(ii)}}{=}
X_1^{-3\nn}\Big (\frac15 + \frac{4 X_1^{5\nn}}5\Big)=
(1-\ep)^{-3\ep^{-1}}\Big (\frac15 + \frac{4 (1-\ep)^{5\ep^{-1}}}5\Big)
\nonumber\\
&\!\!\!\!\!\!\!\!\!\!\!\!\!\!\!\!\!\!\!\!\!\!\!\!\!\!&
\phantom{{\rm(ii)\ }A_2 \ \ \ }
=\ \ \ \Big(e^3+O(\ep)^1\Big)\Big(\frac15+\frac{4}{5e^5}+O(\ep)^1\Big)
=\frac{e^3}{5}+\frac{4}{5e^2}+O(\ep)^1,
\nonumber\\
&\!\!\!\!\!\!\!\!\!\!\!\!\!\!\!\!\!\!\!\!\!\!\!\!\!\!&
{\rm(iii)\ }
A_1^{2\nn}\stackrel{{}^{\sc\rm\eqref{LetNSoOP----1}\,(i)}}{=}X_1^{\nn}
=(1+\ep)^{2\ep^{-1}}=e^2+O(\ep)^1.
\end{eqnarray}
Observe that $1<e<\frac{e^3}{5}+\frac{4}{5e^2}<e^2$
(noting that $e\thickapprox2.7,\,\frac{e^3}{5}+\frac{4}{5e^2}\thickapprox4.1,\,e^2\thickapprox7.3$),
}%
\NOUSE
{
Then by \eqref{XsXXX},\,\eqref{LetF3rst-final},\,\eqref{X1X2Z===},
 we see from Definition \ref{defi-B123} that all $
A_1,A_2,A_3$ 
 are $1+O(\ep)^1$ elements of the forms
\equa{Ai===}{ 
A_i=1+a_i\ep+O(\ep)^2,\ \ i=1,2,3}
   such that the coefficients of $\ep^1$ in them are respectively,
%
\begin{eqnarray}
\label{B1-ep1}
\nonumber
&\!\!\!\!\!\!\!\!\!\!\!\!\!\!\!\!\!\!\!\!\!\!\!\!\!\!\!\!\!\!\!\!\!\!\!\!\!\!
&
{\rm(i)\ }a_1{\sc}\stackrel{{}^{\sc\rm\eqref{LetNSoOP----1}}}{=}
-1+b_1=
\frac{2 + 4 \d_0 + 2 \d_0^2 - \d_0^3 - 4 \d_0^4 - 2 \d_0^5 - \d_0^6}{3 (1 + \d_0)}
,
\!\!\!\!\!\!\!\!\!\!\!\!\!\!\!\!\!\!\!\!\!
\\[4pt]\nonumber
&\!\!\!\!\!\!\!\!\!\!\!\!\!\!\!\!\!\!\!\!\!\!\!\!\!\!\!\!\!\!\!\!\!\!\!\!\!\!\!\!\!\!\!\!\!\!\!&
{\rm(ii)\ }a_3\stackrel{{}^{\sc\rm\eqref{LetNSoOP----1}\,(iii)}}{=}
\d_0^2\Big(\tilde s{\sc\!}+{\sc\!}2\ell_0^2(2{\sc\!}-{\sc\!}\d_0^2)(b_1
{\sc\!}-{\sc\!}a_1)\Big){\sc\!}+{\sc\!}(2{\sc\!}+{\sc\!}\d_0^3)\tilde v
{\sc\!}-{\sc\!}2\tilde s
{\sc\!}={\sc\!}-\frac{4{\sc\!} +{\sc\!}
 8 \d_0 {\sc\!}+{\sc\!} 4 \d_0^2 {\sc\!}+{\sc\!} 4 \d_0^3 {\sc\!}+{\sc\!} 4 \d_0^4
 {\sc\!}+{\sc\!} 2 \d_0^5 {\sc\!}+{\sc\!} \d_0^6}{3(1 +  \d_0)}
,\!\!\!\!\!\!\!\!\!\!\!\!\!\!\!\!\!\!\!\!\!\!\!\!\!\!\!\!\!\!\!\!\!\!\!\!\!\!\!\!\!\!\!\!\!\!\!\!\!\!\!\!\!
\\[4pt]
&\!\!\!\!\!\!\!\!\!\!\!\!\!\!\!\!\!\!\!\!\!\!\!\!\!\!\!\!\!\!\!\!\!\!\!\!\!\!\!\!\!\!\!\!\!\!\!&
{\rm(iii)\ }a_2=
\ell_0^2(2-\d_0^2)(b_1-a_1)+\ell_0^2(1+\d_0^2)(b_3-a_3)-(3+4\d_0)\tilde u
=2 + 5 \d_0 + 2 \d_0^2 + \d_0^3.
\end{eqnarray}
%
%
}%
%
{%

Now in \eqref{LetNSoOP}, we take  $(p_{\ZeRo},p_{\OnE})$ to be $(q_{\ZeRo},q_{\OnE})$.
By \eqref{recallll}\,(iii),\,(iv),
 we can easily see that 
 the middle three terms of \eqref{LetNSoOP}\,(i) are respectively the following three numbers%
~(recalling that $0<\d\ll\d_0=\ell_0^{-1}$)%
,
\begin{eqnarray}
 \label{a1a2a3}
\!\!\!\!\!\!\!\!\!\!\!\!\!\!\!\!\!\!\!\!\!\!\!\!\!\!\!\!\!\!\!\!\!\!&&
(-1+\d^2
)a_1=1-\d^{2}
,
\ \ \ \
\ \ \ \
-a_2=1+O(\d)^3,
\ \ \ \  \ \ \
(-1-\d^2
)a_1=1+\d^{2}
.\!\!\!\!\!\!\!\!\!\!\!\!\!\!\!\!\!\!\!
\end{eqnarray}
We }
see that all
inequalities in \eqref{LetNSoOP}\,(i) are strict inequalities when $(p_1,p_2)$ is set to $(q_1,q_2)$, i.e., \eqref{LetNSoOP}\,(i) is satisfied by $(q_1,q_2)$.

When $(p_1,p_2)$ is set to $(q_1,q_2)$, one can   observe from \eqref{TaKa++}--\eqref{XsXXX} that
$(|x_2|+|x_2+y_2|)\ep^3$ is an $O(\ep)^3$ element
, and
then  from \eqref{LetNSoOP}\,(ii),\,\eqref{recallll}\,(i),\,(ii), we see  that
$\ell_{q_1,q_2}$ is a $1+O(\ep)^1$ element such that
the coefficient of $\ep^1$ in $\ell_{q_1,q_2}$ is
\NOUSE{ $v>0$,
by computing the following, 
\begin{eqnarray}
\label{Bdndn}
&&\!\!\!\!\!\!\!\!\!\!\!\!\!\!\!\!\!\!\!\!\!\!\!\!\!\!
A_3A_1^{-1}
\stackrel{{}^{\sc\rm\eqref{LetNSoOP----1}\,(i),\,(iii)}}{=}
\frac{Z}{X_2^3X_1^3}\stackrel{{}^{\sc\rm\eqref{XsXXX}}}{=}(1+v\ep)(1+u\ep)^{-3}(1+s_0\ep)^{-3}+O(\ep)^2
\nonumber\\
&&\!\!\!\!\!\!\!\!\!\!\!\!\!\!\!\!\!\!\!\!\!\!\!\!\!\!
\phantom{A_3A_1^{-1}==}=\ \ \ 1+(v-3u-3s_0)\ep+O(\ep)^2=1+v\ep+O(\ep)^2.
\end{eqnarray}
}
$u-\d^2a_1
=\d^2
>0$
,
which means
that
 the
inequality  in \eqref{LetNSoOP}\,(ii) is a strict inequality.

\NOUSE{
Note from \eqref{LetNSoOP}\,(ii),\,\eqref{TaKa++},\,\eqref{1++??+q0q1-bar++} that, when $(p_1,p_2)$ is set to $(q_1,q_2)$,
$\ell_{q_1,q_2}$ is a $1+O(\ep)^1$ element, and
the coefficient of $\ep^1$ in $\ell_{q_1,q_2}$ is, by \eqref{ep-pro},\,\eqref{Ai===},
$a_3+(-\frac14-\d^2)a_1=4\d^2>0$ (recalling that $\d\ll\d_0$), which means that
 the
inequality  in \eqref{LetNSoOP}\,(ii) is a strict inequality.
}

When $(p_1,p_2)$ is set to $(q_1,q_2)$, one can easily observe from \eqref{tX1==},\,\eqref{XsXXX}--\eqref{LetF3rst-final} 
that
 $\widetilde X_1,
 X_2,A_1$, $A_2
 $
are all $1+O(\ep)^1$ elements, thus
 all strict inequalities  in
%
\eqref{ToSayas+1}\,(c) 
hold [in particular $(p_1,p_2)\in S_2$, cf.~proof of Lemma \ref{lemm-condition-XZ}]
.
\NOUSE{
As in Remark \ref{X1tilde},
by regarding \eqref{1=bbb}\,(i) as  a quadratic equation of $\tildeX_1^{10}$,
we have two possible choices of $\tildeX_1^{10}$,
one is stated as in
\eqref{ToSayas+1}\,(e)--(h) satisfying \eqref{ToSayas+2} [note that in this case $A_2=1+O(\ep)^1$ thus \eqref{ToSayas+2}\,(ii) holds$\ssc\,$],
the other, which is not a $1+O(\ep)^1$ element, does not satisfy \eqref{W-XXXX}\,(i). Thus we have \eqref{ToSayas+1}\,(e)--(h) and \eqref{ToSayas+2}.
}%

Hence the ``initial stage'' $(q_{\ZeRo},q_{\OnE})\in V_0.$\hfill$\Box$%
}%
\vskip5pt

From now on, we always assume $(p_1,p_2)\in\ol V_0$, and for convenience, we denote
(some notations will be introduced later),
\equa{denote-t2}{\dis\!\!\!\!\!\!\!
T_2=\Big\{\tilde\xX_1=|\widetilde X_1|,\,
\zZ=|Z|,\,
\xX_i{\ssc\!}={\ssc\!}|X_i|,\,\aA_i=|A_i|
,\ \forall\,i
\Big\}.
\!\!\!\!\!\!}
Then by Lemma \ref{lemm-condition-XZ} and \eqref{C+ToSayas+1}, we see that elements of $\ol V_0\subset V$ satisfy
the following%
,
\begin{eqnarray}
\label{C+LetNSoOP}
\!\!\!\!\!\!\!\!\!\!\!\!\!\!\!\!\!\!\!\!\!&\!\!\!\!\!\!\!\!\!\!\!\!\!\!\!\!\!\!\!\!\!\!\!\!\!\!\!\!&
{\rm(i)\ }1\le \aA_{\rOnE}^{-1+\d^2
}
\le \aA_2^{-1}
\le\aA_1^{-1-\d^2
}
\le
\ell_{_{\sc 1}}^{^{\sc1+\d^2}}
,
\nonumber\\[0pt]
\!\!\!\!\!\!\!\!\!\!\!\!\!\!\!\!\!\!\!\!\!\!\!\!&\!\!\!\!\!\!\!\!\!\!\!\!\!\!\!\!\!\!\!\!\!\!\!\!\!\!\!\!&
{\rm(ii)\ }\ell_{p_1,p_2}{\ssc}:={\ssc}
\xX_2\aA_1^{-\d^2
}
{\ssc}+{\ssc}(|x_{\OnE}|{\ssc}+{\ssc}|x_{\OnE}{\ssc}+{\ssc}y_{\OnE}|)\ep_2^3{\ssc}\ge
{\ssc}1{\ssc}+{\ssc}\ep_2^2
,
\nonumber\\[0pt]
\!\!\!\!\!\!\!\!\!\!\!\!\!\!\!\!\!\!\!\!\!\!\!\!&\!\!\!\!\!\!\!\!\!\!\!\!\!\!\!\!\!\!\!\!\!\!\!\!\!\!\!\!&
{\rm(iii)\ }
(1-\d)\xX_2\aA_1^{2+\d_0}
\le\tilde \xX_1\le
(1+\d)\xX_2^{-1}\aA_1^{-\d_0}
.\!\!\!\!\!\!\!\!\!\!\!\!\!\!\!\!\!\!\!\!\!\!\!\!\!\!\!\!\!\!\!\!\!\!\!\!\!\!\!\!\!\!\!\!\!\!\!\!
\end{eqnarray}

\NOUSE
{%
Before continuing, in order to better understand why we define $A_1,A_2,A_3,B_1$ in the way in \eqref{LetNSoOP----1}
, at this point it may be worth presenting the following remark.
\begin{rema}\rm\label{rema-def-A1}\begin{itemize}\item[(1)]
If we choose $A_i,\,i=1,2,3$ to be  monomials, then it is just the same as the case \eqref{ToSayas}.
In that case, as  already mentioned in
Remark \ref{Antotm}, we have Fact \ref{fact-initial}, and we cannot obtain a contradiction. Thus we have to
choose at least one of $A_i$ 
to have more terms (of course, the easiest way is  to
 take one of $A_i$ 
 to have two terms).
Here, for example our original idea is to choose  $A_1$ to satisfy the following
[noting that  \eqref{More-de-1} can be rewritten as the following form],
\equa{A1-form1}{\dis 1=E_1(\beta_1+\beta_2E_2),}
where $E_1,E_2$ are some monomials on $X_1,X_2,Z,A_1,A_2,A_3$,
and $\b_1,\b_2\in\R_{>0}$ such that $\b_1+\b_2=1$ [which is
necessary in order that at the ``initial stage'', the right-hand side of \eqref{A1-form1} is equal to $1$].
Then we define $B_1,C_1$ to be the following [cf.~\eqref{LetNSoOP----1}\,(iv),\,(v)$\ssc\,$],
\equa{B1-C1-rema}{\dis {\rm(i)\ }B_1=
{E_1E_2},\ \ \ \ \ {\rm(ii)\ }C_1=
{E_1^2E_2}.}
In this way, we can obtain from \eqref{A1-form1} the following [cf.~\eqref{Re-Writtt}$\ssc\,$],
\equa{Nosm-Rem}{\dis\b_2B_1^{2}-B_1
+\b_1C_1
=0.}
We need to choose suitable $E_1,E_2,\b_1,\b_2$ to satisfy the following.
\begin{itemize}
\item[(a)]
The elements $B_1,C_1$ defined in \eqref{B1-C1-rema} need to satisfy
%
\eqref{LetNSoOP----1-redefine++}\,(v)%
.
Further, when equality occurs in the last inequality of \eqref{C+LetNSoOP}\,(i) for an element in $\ol V_0$,
the formula
\eqref{A1-form1} should be controlled by the second term inside the bracket. This means that
$E_2\succeq_{\ell_1\ssc\,}1$,
which is  equivalent to that $B_1
\succeq_{\ell_1\ssc\,}1$ [cf.~\eqref{Formmeme}$\ssc\,$]
.
This requires us to choose suitable $E_1,E_2$ [noting that we require condition \eqref{ToSayas+1}\,(c) so
that we  can use \eqref{Nosm-Rem} to obtain  \eqref{Formmeme}%
$\ssc\,$].
\item[(b)]
At the ``initial stage'', we need to require the coefficient of $\ep^1$ in
$A_1$ 
be some number in \eqref{B1-ep}
by adjusting $\b_1,\b_2$.
\NOUSE{%
\item[(iii)]
We also need to choose suitable $\b_1,\b_2$ so that some positive power of $B_1$ in the denominator of \eqref{CCCC===} can be used in the proof of the last inequality of \eqref{ToSayas+1}\,(c). The very crucial fact is the following.
\begin{fact}\label{Fact-power-B1}
The positive power of $B_1$ in the denominator of \eqref{CCCC===} must be
$<\b_2\b_1^{-1}-1$ $($thus $\b_2$ must be $>\b_1$, and so $\frac12<\b_2<1)$, otherwise the coefficient of $\eE_1^2$ in the right-hand side of \eqref{CCCC===}  is $\ge1$ and we cannot obtain   the last inequality of \eqref{ToSayas+1}\,{\rm(c)}. In our case in the paper, we have chosen $\b_1=\frac15,\,\b_2=\frac45$ $[$cf.~\eqref{LetNSoOP----1}\,{\rm(iv)}${\ssc\,}]$, and so $\b_2\b_1^{-1}-1=3$, and the power of $B_1$
 in the denominator of \eqref{CCCC===}  is $\frac52$, which is  $<3$.
 \end{fact}
}%
\end{itemize}\item[(2)]
Note also that by \eqref{LetNSoOP----1}\,(ii),\,(iv) in \eqref{LetNSoOP----1}\,(iii), we have
\begin{eqnarray}
\label{A2-rewrite}
&\!\!\!\!\!\!\!\!\!\!\!\!\!\!\!\!\!\!\!\!\!\!&
{\rm(i)\ }
1=
\frac{A_3\widetilde X_1^{\ell_0}}{A_2X_2}
 \Big(\frac{1 + \d_0}{1 + 2 \d_0} +\frac{\d_0D_0}{1+2\d_0}\Big)
 =\frac{(1+\d_0)B_2}{1+2\d_0}+\frac{\d_0C_2}{(1+2\d_0)B_2},
\!\!\!\!\!\!\!\!\!\!\!\!\!\!\!\!\!\!\!\!\!
\nonumber\\
&\!\!\!\!\!\!\!\!\!\!\!\!\!\!\!\!\!\!\!\!\!\!&
{\rm(ii)\ }B_2=
\frac{X_2^{\ell_0^2(1+2\d_0+\d_0^2)}}{A_2A_3^{\ell_0^2}\widetilde X_1^{\ell_0^4(1+2\d_0-\d_0^2-3\d_0^3-2\d_0^4)}}
=\frac1{B_2},  \
{\rm(iii)\ }C_3=D_0B_3^2=\frac{X_2^{2\ell_0^4(1+4\d_0+5\d_0^2+2\d_0^3)}}{A_2^2A_3^{\ell_0^{10}+2\ell_0^2}
\widetilde X_1^{\ell_0^2(2-3\d_0^2-2\d_0^3)}}
,
\!\!\!\!\!\!\!\!\!\!\!\!\!\!\!\!\!\!\!
\nonumber\\
&\!\!\!\!\!\!\!\!\!\!\!\!\!\!\!\!\!\!\!\!\!\!&
\implies\ \ {\rm(iv)\ }
\ell_0(1+\d_0)B_3^2-\ell_0(1+2\d_0)B_3+C_3=0.
\end{eqnarray}
\item[(4)]Observe that by multiplying  \eqref{LetNSoOP----1}\,(ii) with \eqref{A2-rewrite--}\,(i) we
we obtain  the following, which is the
reason why we define $A_2$ as in \eqref{LetNSoOP----1}\,(iii),
\begin{eqnarray}
\label{A2-rewrite++}
&\!\!\!\!\!\!\!\!\!\!\!\!\!\!\!\!\!\!\!\!\!\!&
{\rm(i)\ }1=\frac1{B_4}\Big(\frac{(1+\d_0)^2}{1+2\d_0}-\frac{\d_0^2D_0^2}{1+2\d_0}\Big)
,\ \ \ \implies
{\rm(ii)\ }\ell_0^2(1+2\d_0)B_4^2-\ell_0^2(1+\d_0)^2B_4+C_4=0,
\!\!\!\!\!\!\!\!\!\!\!\!\!\!\!\!\!\!\!\!\!\nonumber\\
&\!\!\!\!\!\!\!\!\!\!\!\!\!\!\!\!\!\!\!\!\!\!&
{\rm(iii)\ }B_4=\frac{B_1}{B_3}=\frac{A_1A_2A_3^{\ell_0^2-1}
\widetilde X_1^{\ell_0^4(2+4\d_0-4\d_0^2-9\d_0^3+2\d_0^4)}}
{X_2^{\ell_0^2(2+4\d_0-\d_0^2)}Z^2}
,
\nonumber\\
&\!\!\!\!\!\!\!\!\!\!\!\!\!\!\!\!\!\!\!\!\!\!&
 {\rm(iv)\ }C_4=D_0^2B_4=
\frac{A_1A_2
\widetilde X_1^{\ell_0^4(2+4\d_0-2\d_0^2-5\d_0^3-2)}}
{A_3^{\ell_0^{10}-\ell_0^2+1}X_2^{\ell_0^2(2+4\d_0-1)}Z^2}.
\end{eqnarray}
\end{itemize}
\end{rema}\vskip7pt
}%

\NOUSE{
\begin{lemm}\label{ExLemm}When $(p_1,p_2)\in\ol V_0$, we have
\begin{eqnarray}
\!\!\!\!\!\!\!\!\!\!\!\!\!\!\!\!\!\!\!\!\!\!\!\!\!\!\!&&
\label{extra-equa}
{\rm(i)\, }
\d_2<|\tildeX_1 
|
,\ \ \  \ \ \ \ \ {\rm(ii)\ }\d_3<|Z|<\ell_3
,
\nonumber\\
\!\!\!\!\!\!\!\!\!\!\!\!\!\!\!\!\!\!\!\!\!\!\!\!\!\!\!&&
{\rm(iii)\ }\d_3^{\ell_2}<|X_2|<\ell_3^{\ell_2},\ \ \ \ \ \ \
{\rm(iv)\ }\widetilde Z=\widetilde X_2+O(\d)^3
.
\end{eqnarray}
\end{lemm}
\noindent{\it Proof.~}By \eqref{C+LetNSoOP}\,(i), we in particular have the following, where the first inequality in (ii) follows from the fact in \eqref{B1-ep} that $-a_2a_1^{-1}+\d^3<1$,
\equa{+C+LetNSoOP}{
{\rm(i)\ }\d_1\le |A_1|\le1,\ \ \ \ \ \ {\rm(ii)\ }\d_1\le|A_2|\le1.
}
By \eqref{LetNSoOP----1}\,(ii), we have
\begin{eqnarray}
&\!\!\!\!\!\!\!\!\!\!\!\!\!\!\!\!\!\!\!\!\!\!\!\!\!\!\!&
\label{anAAA2-aaa}\pm\Big(
\frac{(1+\d_0^4)}{2}\Big|\frac{A_3^{2-\d_0^{11}}}{X_1^2}\Big|-\frac{1-\d_0^4}{2}
\Big)\le\Big|\frac{A_2A_3^{2-2\d_0^2}}{X_1}\Big|\le\frac{1-\d_0^4}{2}+\frac{(1+\d_0^4)}{2}\Big|\frac{A_3^{2-\d_0^{11}}}{X_1^2}\Big|.
\end{eqnarray}
Recall from Remark \ref{rema3.1} that
 we  have complete freedom in choosing $\ell_2$ with $\ell_2\gg\ell_1$ independently of  choices of $\ell_0,\ell_1$.
To prove \eqref{extra-equa}\,(i),  assume
there exist some $(p_1,p_2)\in\ol V_0$ such that \equa{First-ass}{\mbox{ $|X_1|\le\d_2$.}}
We consider 2 cases. In the first case we  assume
\equa{Case1-equuu}{\dis
\Big|\frac{A_3^{2-\d_0^{11}}}{X_1^2}\Big|\ge1.
}
Then
\equa{Case1-fi1}{\dis
\!\!\!\!\!\!\!\!\!
\d_0^4
\Big|\frac{A_3^{2-\d_0^{11}}}{X_1^2}\Big|\
\stackrel{{}^{\sc\rm\eqref{Case1-equuu}}}{\le}
\frac{(1+\d_0^4)}{2}\Big|\frac{A_3^{2-\d_0^{11}}}{X_1^2}\Big|{\ssc\!}-
{\ssc\!}\frac{1{\ssc\!}-{\ssc\!}\d_0^4}{2}\
\stackrel{{}^{\sc\rm\eqref{anAAA2-aaa}}}{ \le}\Big|\frac{A_2A_3^{2-2\d_0^2}}{X_1}\Big|
\stackrel{{}^{\sc\rm\eqref{+C+LetNSoOP}\,(ii)}}{\le}\Big|\frac{A_3^{2-2\d_0^2}}{X_1}\Big|
.\!\!\!\!\!\!}
This implies
\equa{MSMSMS--=q=2=}
{\dis |A_3|^{1-2\d_0^2+\frac{\d_0^{11}}{2}}\stackrel{{}^{\sc\rm\eqref{Case1-fi1}}}{\ge}\d_0^4
\Big|\frac{A_3^{2-\d_0^{11}}}{X_1^2}\Big|^{\frac12}
\stackrel{{}^{\sc\rm\eqref{Case1-equuu}}}{\ge}\d_0^4.
}
Then we obtain the following contradiction, which shows that the first case does not occur,
\equa{MSMSMS--=q=2=1}
{\d_0^{4+\frac{4(2\d_0^2-\d_0^{11})}{1-2\d_0^2+\frac{\d_0^{11}}{2}}}\stackrel{{}^{\sc\rm\eqref{MSMSMS--=q=2=}}}{\le}
\d_0^4|A_3|^{2\d_0^2-\d_0^{11}}\stackrel{{}^{\sc\rm\eqref{Case1-fi1}}}{\le}|X_1|
\stackrel{{}^{\sc\rm\eqref{First-ass}}}{\le}\d_2\stackrel{{}^{\sc\rm\eqref{MSmde33333}}}{<}
\d_0^{4+\frac{4(2\d_0^2-\d_0^{11})}{1-2\d_0^2+\frac{\d_0^{11}}{2}}}.
}
In the second case we assume
\equa{Next-SSS}{\dis{\rm(i)\ }
|X_1|>|A_3|^{\frac{2-\d_0^{11}}{2}}\stackrel{{}^{\sc\rm\eqref{LetNSoOP----1}\,(iii)}}{=}|Z|^{-\frac{2-\d_0^{11}}{2}},
\ \ \ \ \ {\rm(ii)\ }|Z|\stackrel{{}^{\sc\rm\eqref{Next-SSS}\,(i)}}{>}|X_1|^{-\frac{2}{2-\d_0^{11}}}
\stackrel{{}^{\sc\rm\eqref{First-ass}}}{\ge}\ell_2^{\frac{2}{2-\d_0^{11}}}.}
}\NOUSE{
Then we have
\begin{eqnarray}
&\!\!\!\!\!\!\!\!\!\!\!\!\!\!\!\!\!\!\!\!\!\!\!\!&
\label{absolute-x1}
{\rm(i)\ }|x_1|\stackrel{{}^{\sc\rm\eqref{TaKa},\,\eqref{SimMMSMS}}}{=}\kk(1+\a_0\ep)^{-1}|X_1|^{\d_0^{11}\d}
\stackrel{{}^{\sc\rm\eqref{First-ass}}}{\le}\d_2^{\d_0^{11}\d}\kk\stackrel{{}^{\sc\rm\eqref{MSmde33333}}}{<}\d_1\kk,
\nonumber\\
&\!\!\!\!\!\!\!\!\!\!\!\!\!\!\!\!\!\!\!\!\!\!\!\!&
{\rm(ii)\ }\kk\stackrel{{}^{\sc\rm Lemma\,\ref{YYYy1==}}}{<}\g_{\kk,\kk}
\stackrel{{}^{\sc\rm\eqref{Next-SSS}\,(ii)}}{<}\g_{\kk,\kk}|Z|^{\d_0^{11}}\
\stackrel{{}^{\sc\rm\eqref{TaKa},\,\eqref{SimMMSMS}}}{=}|x_2+y_2|\le|x_2|+|y_2|
\stackrel{{}^{\sc\rm\eqref{dp0p1}}}{\le} h_{p_1,p_2},
\end{eqnarray}
and
\equa{msmwemememem}{\dis
h_{_{\sc p_1,p_2}}^{^{\sc-\frac1{m+1}}}\stackrel{{}^{\sc\rm\eqref{absolute-x1}\,(ii)}}{\le}
\kk^{-\frac1{m+1}}\stackrel{{}^{\sc\rm\eqref{MSmde33333} }}{<}\d_1.}
Then
\begin{eqnarray}
&\!\!\!\!\!\!\!\!\!\!\!\!\!\!\!\!\!\!\!\!\!\!\!\!&
\label{absolute-x1++}
h_{p_1,p_2}\ \ \ \ \ \stackrel{{}^{\sc\rm\eqref{dp0p1}}}{=}\ \ |x_1|+|x_2|+|y_1|+|y_2|
\stackrel{{}^{\sc\rm\eqref{absolute-x1}\,(i),\,
\eqref{mqp1234-2}}}
{<}\d_1\kk+|x_2|+2h_{_{\sc p_1,p_2}}^{^{\sc\frac m{m+1}}}
\nonumber\\
&\!\!\!\!\!\!\!\!\!\!\!\!\!\!\!\!\!\!\!\!\!\!\!\!&
\phantom{h_{p_1,p_2}}\!\!\!
\stackrel{{}^{\sc\rm\eqref{absolute-x1}\,(ii),\,\eqref{msmwemememem}}}{<}\d_1 h_{p_1,p_2}+|x_2|+2\d_1h_{p_1,p_2}
\end{eqnarray}
The above shows
\equa{x-1====hhhh}{\dis(1+3\d_1)^{-1}h_{p_1,p_2}\le|x_2|\stackrel{{}^{\sc\rm
\eqref{dp0p1}}}{\le} h_{p_1,p_2}.
}
Then
\equa{x2+y2-smaller-xx}{\dis
|x_2|-\d_1h_{p_1,p_2}
\stackrel{{}^{\sc\rm\eqref{msmwemememem}}}{\le}|x_2|-|y_2|\le|x_2+y_2|\le|x_2|+|y_2|
\stackrel{{}^{\sc\rm\eqref{msmwemememem}}}{\le}|x_2|+\d_1h_{p_1,p_2}.
}
Thus
\equa{msmsmem22233}{\dis|x_2+y_2|\stackrel{{}^{\sc\rm \eqref{x-1====hhhh},\,\eqref{x2+y2-smaller-xx}}}{=}\Big(1+O(\d_1)^1\Big)|x_2|.}
Note that since $\kk\gg\ell_1$,  we can choose $\d'$ in Lemma \ref{YYYy1==} to be $\d_1$, which shows that $\g_{\kk,\kk}=\big(1+O(\d_1)^1\big)\kk$.
Then
we have the following, by noting that $\big(1+O(\d_1)^1\big)^{\ell_0^{11}}=1+O(\d_1)^1$ by the fact in \eqref{MSmde33333} that $\d_1\ll\d_0=\ell_0^{-1}$,
\begin{eqnarray}
&\!\!\!\!\!\!\!\!\!\!\!\!\!\!\!\!\!\!\!\!\!\!\!\!\!\!\!\!\!\!\!\!&
\label{Z===X2-simm}
|Z|
\stackrel{{}^{\sc\rm\eqref{TaKa},\,\eqref{SimMMSMS}}}
{=}\Big|\g_{\kk,\kk}^{-1}(x_2+y_2)\Big|^{\ell_0^{11}}\,\ \stackrel{{}^{\sc\rm\eqref{msmsmem22233}}}{=}
\Big|\Big(1+O(\d_1)^1\Big)\kk^{-1}x_2\Big|^{\ell_0^{11}}\,\
\stackrel{{}^{\sc\rm\eqref{TaKa},\,\eqref{SimMMSMS}}}
{=}\Big(1+O(\d_1)^1\Big)|X_2|.\!\!
\end{eqnarray}
Then
\begin{eqnarray}
&\!\!\!\!\!\!\!\!\!\!\!\!\!\!\!\!\!\!\!\!&
\label{PartA1111}
{\rm(i)\ }
\a_1:=\Big|\frac{ X_2}{A_3^{\d_0^2}X_1 Z^{1+\d_0^2}}\Big|
\stackrel{{}^{\sc\rm\eqref{LetNSoOP----1}\,(iii)}}{=}\Big|\frac{X_2}{X_1Z}\Big|
\stackrel{{}^{\sc\rm\eqref{Z===X2-simm}}}
{=}\Big(1+O(\d_1)^1\Big)|X_1|^{-1}\stackrel{{}^{\sc\rm\eqref{First-ass}}}{\ge}\Big(1+O(\d_1)^1\Big)\ell_2,
\nonumber\\
&\!\!\!\!\!\!\!\!\!\!\!\!\!\!\!\!\!\!\!\!&
{\rm(ii)\ }
\a_2:=\Big| \frac{A_2^6 X_1^2 Z^{2 + 2 \d_0^2}}
{    X_2^2}\Big|\stackrel{{}^{\sc\rm\eqref{Z===X2-simm}}}{=}\Big(1+O(\d_1)^1)\Big| 
{A_2^6 X_1^2 Z^{ 2 \d_0^2}}
{    }\Big|\stackrel{{}^{\sc\rm\eqref{LetNSoOP----1}\,(iii)}}{=}\Big(1+O(\d_1)^1\Big)
\Big|\frac{A_2^6X_1^2}{A_3^{2\d_0^2}}\Big|.
\end{eqnarray}
By \eqref{LetNSoOP----1}\,(i), we have
\equa{M54-4050}{\dis
\pm\Big(\frac{1+\d_0^4}{2}-\frac{(1-\d_0^4)\a_2}{2}\Big)\stackrel{{}^{\sc\rm\eqref{LetNSoOP----1}\,(i)}}{\le}
|A_1|\a_1^{-1}\stackrel{{}^{\sc\rm\eqref{+C+LetNSoOP}\,(i),\,\eqref{PartA1111}\,(ii)}}
{\le}\Big(1+O(\d_1)^1\Big)\d_2=O(\d_2)^1.
}
Thus
\equa{aplsle32}{\dis
\a_2\stackrel{{}^{\sc\rm\eqref{M54-4050}}}{=}\frac{1+\d_0^4}{1-\d_0^4}+O(\d_2)^1=1+O(\d_0)^2.
}
Then
\begin{eqnarray}
\label{AnaThR}
\!\!\!\!\!\!\!\!\!\!\!\!\!\!\!\!\!\!\!\!&&
\Big|\frac{A_2A_3^{2-2\d_0^2}}{X_1}\Big|\ \ \ =\ \ \
\Big|\Big(\frac{A_3^{2\d_0^2}}{A_2^6X_1^2}\Big)^{\frac12}A_2^3A_3^{2-3\d_0^2}\Big|
\stackrel{{}^{\sc\rm\eqref{PartA1111}\,(ii),\,\eqref{aplsle32}}}{=}\Big(1+O(\d_1)^1\Big)\a_2^{-\frac12}|A_3|^{2-3\d_0^2}
\nonumber\\
\!\!\!\!\!\!\!\!\!\!\!\!\!\!\!\!\!\!\!\!&&\ \ \ \ \ \ \ \ \ \ \ \ \ \ \ \,
\stackrel{{}^{\sc\rm\eqref{LetNSoOP----1}\,(iii)}}{=}
\Big(1+O(\d_1)^1\Big)\a_2^{-\frac12}|Z|^{-2+3\d_0^2}\ \,
\stackrel{{}^{\sc\rm\eqref{Next-SSS}\,(ii),\,\eqref{aplsle32}}}{=}O(\d_2)^1.
\end{eqnarray}
This together with the first inequality of \eqref{anAAA2-aaa} shows
\equa{sho-that---}{\dis\Big|\frac{A_3^{2-\d_0^{11}}}{X_1^2}\Big|=
\frac{1-\d_0^4}{1+\d_0^4}+O(\d_2)^1=1+O(\d_0)^4.}
We obtain the following contradiction, which shows that the second case does not occur either,
\equa{al2l2-that}{\dis
1+O(\d_0)^2\stackrel{{}^{\sc\rm\eqref{aplsle32}}}{=}
\a_2\stackrel{{}^{\sc\rm\eqref{PartA1111}\,(ii),\,\eqref{sho-that---}}}{=}
\Big(1+O(\d_0)^4\Big)\Big|
{A_2^6 X_1^{2-\frac{4\d_0^2}{2-\d_0^{11}}}}{}\Big|
\stackrel{{}^{\sc\rm\eqref{+C+LetNSoOP}\,(ii),\,\eqref{First-ass}}}{=}O(\d_2)^1.}
}\NOUSE{
This together with \eqref{LetNSoOP----1}\,(i) gives
\equa{Ananwne}{\dis1\stackrel{{}^{\sc\rm\eqref{PartA1111},\,\eqref{LetNSoOP----1}\,(i)}}{\sim_{\ell_2\ssc\,}}
\Big|\frac{ X_1^2 Z^{2 + 2 \d_0^2}}
{    A_2^2 X_2^2}\Big|\stackrel{{}^{\sc\rm\eqref{A1-A2-sim-ell2},\,\eqref{Z===X2-simm}}}{\sim_{\ell_2\ssc\,}}
|X_1Z^{2\d_0^2}|\stackrel{{}^{\sc\rm\eqref{Z-is-bigg}}}{\sim_{\ell_2\ssc\,}}|X_1|^{1-2\d_0^2}
\stackrel{{}^{\sc\rm\eqref{Z-is-bigg}}}{\prec_{\ell_2\ssc\,}}1,
}
a contradiction.
In this sense we can regard $\ell_0,\ell_1$ as fixed numbers, and regard $\ell_2$ as a parameter tending to infinity.
 Thus we can use \eqref{MSM1111}\,(i) in Notation \ref{nota-sim} by comparing \eqref{+C+LetNSoOP}\,(i) with \eqref{MSM1111+}\,(i) to obtain
the first ``\,$\sim_{\ell_2\ssc\,}$\,'' below, similarly we have the second ``\,$\sim_{\ell_2\ssc\,}$\,'' below:
\equa{A1-A2-sim-ell2}{|A_1|\stackrel{{}^{\sc\rm\eqref{+C+LetNSoOP}\,(i)}}{\sim_{\ell_2\ssc\,}}
1\stackrel{{}^{\sc\rm\eqref{+C+LetNSoOP}\,(ii)}}{\sim_{\ell_2\ssc\,}}|A_2|.}
Then
\equa{anAAA2-}{\dis
1\stackrel{{}^{\sc\rm\eqref{A1-A2-sim-ell2}}}{\sim_{\ell_2\ssc\,}}|A_2|
\stackrel{{}^{\sc\rm\eqref{LetNSoOP----1}\,(ii)}}{=}\Big|\frac{X_1}{ A_3^{2- 2 \d_0^2}}\Big|\cdot\Big|\frac{1 - \d_0^4}{2}
  +\frac{ (1 + \d_0^4)A_3^2}{2 X_1^2}\Big|.}
We consider case by case. First assume
 \equa{X1-sim-A3-1}{{\rm(i)\ }|X_1|\sim_{\ell_2\ssc\,}|A_3|^{2-2\d_0^2},\ \ \mbox{ or equivalently, \ \ }{\rm(ii)\ }\Big|\frac{X_1}{A_3^{2-2\d_0^2}}\Big|\sim_{\ell_2\ssc\,}1.}
Then \eqref{X1-sim-A3-1}\,(ii) shows that the right-hand side of \eqref{anAAA2-} is
\equa{anAAA2}{\dis
\mbox{r.h.s of \eqref{anAAA2-} \ }
\stackrel{{}^{\sc\rm\eqref{X1-sim-A3-1}\,(ii)}}{\sim_{\ell_2\ssc\,}}\Big|\frac{1 - \d_0^4}{2}
  +\frac{ (1 + \d_0^4)A_3^2}{2 X_1^2}\Big|,}
which with \eqref{anAAA2-} implies
\equa{Amama}{\dis\Big|\frac{A_3^2}{ X_1^2}\Big|\preceq_{\ell_2\ssc\,}1.
}
We obtain the following contradiction,
\equa{A3---3-3}{1\stackrel{{}^{\sc\rm\eqref{Amama}}}{\succeq_{\ell_2\ssc\,}}\Big|\frac{A_3^2}{ X_1^2}\Big|
\stackrel{{}^{\sc\rm\eqref{X1-sim-A3-1}\,(i)}}{\sim_{\ell_2\ssc\,}}\Big|\frac{X_1^{\frac2{2-2\d_0^2}}}{X_1^2}\Big|
=\big|X_2|^{\frac{-2+4\d_0^2}{2-2\d_0^2}}\
\stackrel{{}^{\sc\rm\eqref{First-ass}}}{\succ_{\ell_2\ssc\,}}
1.}
Next assume
 \equa{+X1-sim-A3-1}{{\rm(i)\ }|X_1|\prec_{\ell_2\ssc\,}|A_3|^{2-2\d_0^2},\ \ \mbox{ or equivalently, \ \ }{\rm(ii)\ }\Big|\frac{X_1}{A_3^{2-2\d_0^2}}\Big|\prec_{\ell_2\ssc\,}1.}
Then \eqref{+X1-sim-A3-1}\,(ii) shows that the ``\,$\sim_{\ell_2\ssc\,}$\,'' in \eqref{anAAA2} becomes
``\,$\prec_{\ell_2\ssc\,}$\,'', which implies that
the ``\,$\preceq_{\ell_2\ssc\,}$\,'' in \eqref{Amama} becomes ``\,$\succ_{\ell_2\ssc\,}$\,'', i.e.,
$\big|\frac{A_3^2}{X_2^2}\big|\succ_{\ell_2\ssc\,}1$. This with  \eqref{anAAA2-} implies the second ``\,$\sim_{\ell_2\ssc\,}$\,'' below:
\equa{+A3---3-3}{1\stackrel{{}^{\sc\rm\eqref{A1-A2-sim-ell2}}}{\sim_{\ell_2\ssc\,}}|A_1|
{\sim_{\ell_2\ssc\,}}
\Big|\frac{X_1}{A_3^{2-2\d_0^2}}\Big|\cdot\Big|\frac{A_3^2}{X_2^2}\Big|
=\Big|\frac{A_3^{2\d_0^2}}{X_1}\Big|\stackrel{{}^{\sc\rm\eqref{+X1-sim-A3-1}\,(i)}}{\succ_{\ell_2\ssc\,}}
\Big|\frac{X_1^{\frac{2\d_0^2}{2-2\d_0^2}}}{X_1}\Big|=|X_1|^{\frac{-2+4\d_0^2}{2-2\d_0^2}}\
\stackrel{{}^{\sc\rm\eqref{First-ass}}}{\succ_{\ell_2\ssc\,}}1
,}
a contradiction with \eqref{First-ass}.
Now assume
 \equa{++X1-sim-A3-1}{{\rm(i)\ }|X_1|\succ_{\ell_2\ssc\,}|A_3|^{2-2\d_0^2},\ \ \mbox{ or equivalently, \ \ }{\rm(ii)\ }\Big|\frac{X_1}{A_3^{2-2\d_0^2}}\Big|\succ_{\ell_2\ssc\,}1.}
Then \eqref{++X1-sim-A3-1}\,(ii) shows that
the  ``\,$\sim_{\ell_2\ssc\,}$\,'' in \eqref{anAAA2} becomes
``\,$\succ_{\ell_2\ssc\,}$\,'', which implies that $\big|\frac{A_3^2}{X_1^2}\big|\sim_{\ell_2\ssc\,}1$, i.e., $|X_1|\sim_{\ell_2\ssc\,}|A_3|$. Thus
  \equa{Z-is-bigg}{1\stackrel{{}^{\sc\rm\eqref{First-ass}}}{\succ_{\ell_2\ssc\,}}|X_1|\sim_{\ell_2\ssc\,}|A_3|=\frac1{|Z|}.}
}%
\NOUSE{
The above shows that we cannot have \eqref{First-ass}, i.e., we have \eqref{extra-equa}\,(i).
Now assume $|Z|\le\d_3$. Then $|A_3|\ge\ell_3$ by \eqref{LetNSoOP----1}\,(iii) and $\big|\frac{A_3^{2-\d_0^{11}}}{X_1^2}\big|\ge\d_2^2\ell_3^{2-\d_0^{11}}$ by \eqref{C+LetNSoOP}\,(iv). Thus $1-\d_0^4<\d_0^4\big|\frac{A_3^{2-\d_0^{11}}}{X_1^2}\big|$ by \eqref{MSmde33333}, and the left-hand side of \eqref{anAAA2-aaa}, when we take the positive sign, is bigger than $\frac12\big|\frac{A_3^{2-\d_0^{11}}}{X_1^2}\big|$, which with the first inequality of \eqref{anAAA2-aaa} implies
that $\frac12\big|\frac{A_3^{2-\d_0^{11}}}{X_1^2}\big|\le\big|\frac{A_2A_3^{2-2\d_0^2}}{X_1}\big|$, that implies the second inequality below,
\equa{a3AAA3}{\dis
\ell_3^{2\d_0^2-\d_0^{11}}\le|A_3|^{2\d_0^2-\d_0^{11}}
{\,\le\,} 2|A_2X_1|\stackrel{{}^{\sc\rm\eqref{C+LetNSoOP}\,(iv),\,\eqref{+C+LetNSoOP}\,(ii)}}{\le}2\ell_2
\stackrel{{}^{\sc\rm\eqref{MSmde33333}}}{<}\ell_3^{2\d_0^2}.}
This is a contradiction. Assume $|Z|\ge\ell_3$.
Then $|A_3|\le\d_3$ by \eqref{LetNSoOP----1}\,(iii) and $\big|\frac{A_3^{2-\d_0^{11}}}{X_1^2}\big|\le\ell_2^2\d_3^{2-\d_0^{11}}=O(\d_3)^1$ by \eqref{C+LetNSoOP}\,(iv). Thus we obtain the following contradiction,
\equa{MRMm4m4m5}{\dis
\frac{1-\d_0^4}{4}+O(\d_3)^2\stackrel{{}^{\sc\rm\eqref{anAAA2-aaa}}}{=}\Big|\frac{A_2A_3^{2-2\d_0^2}}{X_1}\Big|
\stackrel{{}^{\sc\rm\eqref{extra-equa}\,(i),\,\eqref{+C+LetNSoOP}\,(ii)}}{\le}
\ell_2|A_3|^{2-2\d_0^2}\le\ell_2\d_3^{2-2\d_0^2}\
\stackrel{{}^{\sc\rm\eqref{MSmde33333}}}{=}O(\d_3)^1
.
}
This proves \eqref{extra-equa}\,(ii).
Now assume $|X_2|\le\d_3^{\ell_1}$. We define $\a_3$ in (i) below, then
\equa{alpha-4444}
{\dis
{\rm(i)\ }\a_3:=\frac{ X_2^2}{A_2^6 X_1^2 Z^{2 + 2 \d_0^2}},\ \ \ \ \
{\rm(ii)\ }|\a_3|\stackrel{{}^{\sc\rm\eqref{+C+LetNSoOP}\,(ii),\,\eqref{extra-equa}\,(i),\,(ii)}}
{\le}
\d_3^{2\ell_1}\ell_1^6\ell_2^2\ell_3^{2+2\d_0^2}\ \stackrel{{}^{\sc\rm\eqref{MSmde33333}}}{=}O(\d_3)^1,
}
and we obtain the following contradiction,
\begin{eqnarray}
&\!\!\!\!\!\!\!\!\!\!\!\!\!\!\!\!\!\!\!\!\!\!\!\!\!&
\label{A1=spwl2,}\dis\!\!\!\!\!\!\!\!\!
1\stackrel{{}^{\sc\rm\eqref{+C+LetNSoOP}\,(i)}}{\ge}|A_1|
\stackrel{{}^{\sc\rm\eqref{LetNSoOP----1}\,(i),\,(iii)}}{=}
\Big|\frac{X_2}{X_1Z\a_3}\Big|\cdot\Big|\frac{1+\d_0^4}{2}+\frac{(1-\d_0^4)\a_3}{2}\Big|
\\\nonumber
&\!\!\!\!\!\!\!\!\!\!\!\!\!\!\!\!\!\!\!\!\!\!\!\!\!\!\!\!\!&
\stackrel{{}^{\sc\rm\eqref{alpha-4444}}}{=}
\Big|\frac{A_2^6X_1Z^{1+2\d_0^2}}{X_2}\Big|\Big(\frac{1+\d_0^4}{2}+O(\d_3)^1\Big)
\stackrel{{}^{\sc\rm\eqref{extra-equa}\,(i),\,(ii)}}{\ge}\d_1^6\d_2\d_3^{1+2\d_0^2}\ell_3^{\ell_1}\Big(\frac{1+\d_0^4}{2}+O(\d_3)^1\Big)
\stackrel{{}^{\sc\rm\eqref{MSmde33333}}}{>}1.\!\!\!\!\!\!\!\!\!\!\!\!\!\!\!\!\!\!\end{eqnarray}
Assume $|X_2|\ge\ell_3^{\ell_1}$. Then with $\a_3$ being defined in
\eqref{alpha-4444}\,(i), we obtain as in \eqref{alpha-4444}\,(ii), $|\a_3|^{-1}=O(\d_3)^1$, and as in
\eqref{A1=spwl2,}, we have the following contradiction,
\begin{eqnarray}
&\!\!\!\!\!\!\!\!\!\!\!\!\!\!\!\!\!\!\!\!\!\!&
\label{A1=spwl2,++}\dis\!\!\!\!\!\!\!\!\!
1\stackrel{{}^{\sc\rm\eqref{+C+LetNSoOP}\,(i)}}{\ge}|A_1|
\stackrel{{}^{\sc\rm\eqref{LetNSoOP----1}\,(i),\,(iii)}}{=}
\Big|\frac{X_2}{X_1Z}\Big|\cdot\Big|\frac{1-\d_0^4}{2}+\frac{1+\d_0^4}{2\a_3}\Big|
\nonumber\\
&\!\!\!\!\!\!\!\!\!\!\!\!\!\!\!\!\!\!\!\!\!\!&
\stackrel{{}^{\sc\rm\eqref{alpha-4444}}}{=}
\Big|\frac{X_2}{X_1Z}\Big|\Big(\frac{1-\d_0^4}{2}+O(\d_3)^1\Big)
\stackrel{{}^{\sc\rm\eqref{extra-equa}\,(i),\,(ii)}}{\ge}\d_2\d_3^{1}\ell_3^{\ell_1}\Big(\frac{1-\d_0^4}{2}+O(\d_3)^1\Big)
\stackrel{{}^{\sc\rm\eqref{MSmde33333}}}{>}1.\end{eqnarray}
This proves \eqref{extra-equa}\,(iii).
To prove \eqref{extra-equa}\,(iv),  we now have
\begin{eqnarray}
&\!\!\!\!\!\!\!\!\!\!\!\!\!\!\!\!\!\!\!\!\!\!\!\!&
\label{absolute-x1+1}
{\rm(i)\ }\Big(1+O(\ep)^1\Big)\d_2^{\d_0^{11}\d}\kk\stackrel{{}^{\sc\rm\eqref{extra-equa}\,(i)}}{\le}\kk(1+\a_0\ep)^{-1}|X_1|^{\d_0^{11}\d}
\stackrel{{}^{\sc\rm\eqref{TaKa},\,\eqref{SimMMSMS}}}{=}|x_1|
\stackrel{{}^{\sc\rm\eqref{C+LetNSoOP}\,(iv)}}{\le}\ell_2^{\d_0^{11}\d}\kk
,
\nonumber\\
&\!\!\!\!\!\!\!\!\!\!\!\!\!\!\!\!\!\!\!\!\!\!\!\!&
{\rm(ii)\ }\d_3^{\d_0^{11}\ell_1}\kk\stackrel{{}^{\sc\rm\eqref{extra-equa}\,(iii)}}{\le}\kk|X_2|^{\d_0^{11}}
\stackrel{{}^{\sc\rm\eqref{TaKa},\,\eqref{SimMMSMS}}}{=}|x_2|
\stackrel{{}^{\sc\rm\eqref{extra-equa}\,(iii)}}{\le}\ell_3^{\d_0^{11}\ell_1}\kk,
\nonumber\\
&\!\!\!\!\!\!\!\!\!\!\!\!\!\!\!\!\!\!\!\!\!\!\!\!&
{\rm(iii)\ }\d_3^{\d_0^{11}}\kk\stackrel{{}^{\sc\rm Lemma\,\ref{YYYy1==}}}{<}\d_3^{\d_0^{11}}\g_{\kk,\kk}
\stackrel{{}^{\sc\rm\eqref{extra-equa}\,(ii)}}{<}\g_{\kk,\kk}|Z|^{\d_0^{11}}
\stackrel{{}^{\sc\rm\eqref{TaKa},\,\eqref{SimMMSMS}}}{=}|x_2+y_2|\le|x_2|+|y_2|
\stackrel{{}^{\sc\rm\eqref{dp0p1}}}{\le} h_{p_1,p_2},\!\!\!\!
\nonumber\\
&\!\!\!\!\!\!\!\!\!\!\!\!\!\!\!\!\!\!\!\!\!\!\!\!&
{\rm(iv)\ }h_{_{\sc p_1,p_2}}^{^{\sc-\frac1{m+1}}}\stackrel{{}^{\sc\rm\eqref{absolute-x1+1}\,(iii)}}{\le}
(\d_3^{\d_0^{11}}\kk)^{-\frac1{m+1}}\stackrel{{}^{\sc\rm\eqref{MSmde33333} }}{<}\d^4.
\nonumber\\
&\!\!\!\!\!\!\!\!\!\!\!\!\!\!\!\!\!\!\!\!\!\!\!\!&
{\rm(v)\ }h_{p_1,p_2}
\stackrel{{}^{\sc\rm \eqref{dp0p1}}}{=}\ \ |x_1|+|x_2|+|y_1|+|y_2|
\stackrel{{}^{\sc\rm\eqref{absolute-x1+1}\,(i),\,(ii),(iv),\,\eqref{mqp1234-2}}}{\le}
 (\ell_2^{\d_0^{11}\d}+\ell_3^{\d_0^{11}\ell_1})\kk+2\d^4h_{p_1,p_2},
\end{eqnarray}
Note that by choosing $\d'=\d^4$ in  Lemma \ref{YYYy1==} we obtain that $\g_{\kk,\kk}=\big(1+O(\d)^4\big)\kk$ and that
$\d^4\widetilde Z=O(\d)^4$ by \eqref{extra-equa},\,\eqref{MSmde33333} and the fact in \eqref{SimMMSMS}\,(vi) that $|\widetilde Z|=|Z|^{\d_0^{11}}$.
Assume $|\widetilde Z-\widetilde X_2|\ge\d^3$. Then we have the following,
where the last inequality follows from the fact in \eqref{absolute-x1+1}\,(iii) that $h_{p_1,p_2}<(1-2\d^4)^{-1} (\ell_2^{\d_0^{11}\d}+\ell_3^{\d_0^{11}\ell_1})\kk$,
\begin{eqnarray}
\label{y1IsSmall}
&\!\!\!\!\!\!\!\!\!\!\!\!\!\!\!\!\!\!\!\!\!\!\!\!\!\!\!&
|y_2|\ge|x_2+y_2-x_2|
\stackrel{{}^{\sc\rm\eqref{TaKa},\,\eqref{SimMMSMS}}}{=}|\g_{\kk,\kk}\widetilde Z-\kk\widetilde X_2|
\ge(\d^3+O(\d)^4)\kk>h_{_{\sc p_1,p_2}}^{^{\sc\frac{m}{m+1}}},
\end{eqnarray}
which is a contradiction with \eqref{mqp1234-2}. This proves \eqref{extra-equa}\,(iv).
\hfill$\Box$
}}
%
\NOUSE{
\begin{rema}\rm\label{rema-nextlemma}
We give some explanations on the next lemma.
The reason we define \eqref{SimMMSMS}\,(????iv) is to obtain \eqref{ImMpP}\,(1), which controls  $|\widetilde X_1|$. From this, we can mainly use Theorem \ref{Theo-2} to obtain
\eqref{ImMpP}\,(2), which is extremely important to us in the following sense.
\begin{itemize}\item[(i)]Firstly, we can obtain \eqref{ImMpP}\,(4)--(6) that means that $A_1,A_2,A_3$,
which are originally  functions of three variables $\tildeX_1,X_2,Z$, become  functions of two variables $\tildeX_1,X_2$, up to $O(\d)^3$ (which, as we shall see, will not affect our arguments below). This  allows us to obtain \eqref{ImMpP}\,(7)--(9), which says that  no equality can occur in any inequality of
\eqref{LetNSoOP}\,(iii),\,(iv) for any element in $\ol V_0$.
 \item[(ii)]Secondly, as mentioned in Remark \ref{AboutX1}\,(iv) and Remark \ref{Antotm}\,(ii),
we are able to prove that  equality cannot occur in the last inequality of \eqref{LetNSoOP}\,(i)  for any element in $\ol V_0$%
.
\end{itemize}\end{rema}
}

\NOUSE{
Denote
\begin{eqnarray}
&\!\!\!\!\!\!\!\!\!\!\!\!\!\!\!\!\!\!\!\!\!\!\!\!\!\!\!&
\label{C1-2and}
{\rm(i)\ }C_1:=\frac{A_1}{ X_1^{\ell_0}X_2^{\ell_0(1+\d_0)}},\ \ \ \ \
\stackrel{{}^{\sc\rm\eqref{C+ToSayas+1}\,(f)}}{\implies} \ \ \ \ \ {\rm(ii)\ }|C_1|\le1+\d,
\nonumber\\
&\!\!\!\!\!\!\!\!\!\!\!\!\!\!\!\!\!\!\!\!\!\!\!\!\!\!\!&
{\rm(iii)\ }D_1:=\frac{X_2^{2\ell_0 (1 + \d_0)}}{X_1^{\ell_0^2(1-\d_0)}A_1^2}
.\end{eqnarray}
}

We will need to apply frequently formula \eqref{(T0o(eE)1}. To obtain that formula, we
require the following [noting that  $\widetilde T_0$ is the set consisting of absolute values of
non-monomial factors
appearing in \eqref{LetNSoOP----1-re-give}
$\ssc\,$],
\NOUSE{We already see from \eqref{A1-A2-cond} that
\eqref{ImMpP}\,(ii) holds for $a=\aA_1,\aA_2$, thus also holds for $a=\tilde\xX_1,\xX_1$ by \eqref{tX1==},\,\eqref{C+ToSayas}\,(c). Then by \eqref{LetNSoOP----1}\,(ii), we see that $\aA_3\le1+O(\d)^1$, from this and \eqref{ImMpP}\,(i),  $\aA_3$ is also a $1+O(\d)^1$. By \eqref{Case6-lemm}\,(iv),
$A_3=\frac{Z^3}{\widetilde X_1^{22} X_2^4}=\frac{1}{\widetilde X_1^{22}X_2}+O(\d)^2$, this implies that $\xX_2,\zZ$ are also $1+O(\d)^1$ element.
This proves \eqref{ImMpP}\,(ii). From this, we obtain
}
{
\begin{eqnarray}
\label{Amsmene}
\!\!\!\!\!\!\!\!\!\!\!\!\!\!\!\!\!\!\!\!\!\!\!\!\!\!\!\!&&
\d_2^{\ell_1}{\sc\!}<{\sc\!} a{\sc\!}<{\sc\!}\ell_2^{\ell_1}\mbox{ for }a\in
\widetilde T_2,\mbox{ where  }
\widetilde T_2{\sc\!}={\sc\!}T_2\cup\widetilde T_0,\ \ \
\widetilde T_0=\mbox{\Large
$\Big\{$}
\Big|2 -\frac1{\widetilde X_1}\Big|,
\Big |2-\frac{\widetilde  X_1}{A_1^2}\Big|
%
\mbox{\Large
$\Big\}$}
.
\end{eqnarray}
First by \eqref{meme}
,\,\eqref{equa-Case6-lemm}\,(iv),\,\eqref{tX1==},\,\eqref{LetNSoOP----1}\,(iii),\,\eqref{A1-A2-cond}, we see
that \eqref{Amsmene} holds for $a\in T_2$
. 
Then by  \eqref{LetNSoOP----1}\,(ii),\,(iii),
\,\eqref{A1-A2-cond},%
~we see that it holds for all $a\in\tilde T_2$.
Therefore, by the fact that $0<\d=\ell^{-1}\ll\d_2=\ell_2^{-1}\ll\d_1=\ell_1^{-1}\ll\d_0=\ell_0^{-1}$, we have,
\equa{(T0o(eE)1}{\dis
a^{0+O(\d)^1}=1+O(\d)^1,\ \ \ a\Big(1+O(\d)^1\Big)=a+O(\d)^1
\mbox{ \ \ for \ \ }a\in\widetilde T_2\mbox{ \ or \ }a^{-1}\in\widetilde T_2.
}

\NOUSE{
\begin{rema}\label{rem-xz-not-zero}\rm In \eqref{LetNSoOP----1},
the functions $A_1,A_2,A_3,B_1,C_1$ are defined only when  $Z,\widetilde X_1,X_2$ are different from $0$ and so a priori the set $V_0$  as in \eqref{C+LetNSoOP} should be considered only in this  open set.  But in the next Lemma we will prove that $V_0$  is contained in the closed set  where  each of  $Z,X_1,X_2$ is bounded away from 0 by some explicit constant and so in this set the functions are well defined and in fact we will show that they are non zero.  Therefore we can always think of them as invertible holomorphic functions.
%
\end{rema}
}
Now we have

\begin{lemm}
\label{Step1}
For $(p_1,p_2)\in\ol V_0$, we have
,
\begin{eqnarray}\label{ImMpP}
&\!\!\!\!\!\!\!\!\!\!\!\!\!\!\!\!\!\!\!\!\!\!\!\!\!\!\!\!\!\!
\!\!\!\!\!\!\!\!\!\!
&
{\rm(i)\  }\xX_2
>\aA_1^{\d^2}\stackrel{{}^{\sc\rm\eqref{A1-A2-cond}\,(a),\,\eqref{(T0o(eE)1}}}{=}1+O(\d)^2
,\ \ \ \ \ \ {\rm(ii)\ }
\xX_2<\ell_1^2,
\ \ \ \ \ \
{\rm(iii)\ }\d_1<\tilde\xX_1<\ell_1^3,
\!\!\!\!\!\!\!\!\!\!
\!\!\!\!\!\!\!\!\!\!\!\!\!\!\!\!\!\!\!\!
\nonumber\\
&\!\!\!\!\!\!\!\!\!\!\!\!\!\!\!\!\!\!\!\!\!\!\!\!\!\!\!\!\!\!
\!\!\!\!\!\!\!\!\!\!
\!\!\!\!\!\!\!\!\!\!\!\!\!\!\!\!\!\!\!\!
&
{\rm(iv)\ }
(1-\d)\xX_2\aA_1^{2+\d_0}
<\tilde \xX_1<
(1+\d)\xX_2^{-1}\aA_1^{-\d_0}
,\ \ \ \ \ {\rm(v)\ }
\Big|2-\frac1{\widetilde X_1}\Big|>\d_2
.\!\!\!\!\!
\end{eqnarray}
 \end{lemm}%
\noindent{\it Proof
.~}First using notation \eqref{denote-t2}, we have $|x_1|=\kk\xX_1,\,|x_2|=\kk\xX_2$ by
\eqref{TaKa},\,\eqref{SimMMSMS}.%
\NOUSE{%
We obtain
\begin{eqnarray}
\label{bsbsbs}
\!\!\!\!\!\!\!\!\!\!\!\!\!\!\!\!\!\!\!\!\!\!\!\!&&
{\rm(i)\ }
\eE_1^{\ell_1}\kk\stackrel{{}^{\sc\rm\eqref{C+LetNSoOP}\,(iii)}}{\le}|x_1|
\stackrel{{}^{\sc\rm\eqref{TaKa},\,\eqref{SimMMSMS}}}{=}
\xX_1\kk
\stackrel{{}^{\sc\rm\eqref{C+LetNSoOP}\,(iii)}}{\le}
(1+\eE_2)\aA_1^{-\frac15}\stackrel{{}^{\sc\rm\eqref{A1-A2-cond}\,(a)}}{\le}1+\eE_1,
\nonumber\\
\!\!\!\!\!\!\!\!\!\!\!\!\!\!\!\!\!\!\!\!\!\!\!\!&&
{\rm(ii)\ }(1-\d_1)\kk
\stackrel{{}^{\sc\rm\eqref{C+LetNSoOP}\,(iii)}}{\le}
\xX_2\kk=|x_2|,
\ \ \ \stackrel{{}^{\sc\rm\eqref{MSmde33333+}}}{\implies}\ \ \
{\rm(iii)\ }|x_1|\preceq_{\kk\ssc\,}\kk\preceq_{\kk\ssc\,}|x_2|.
\!\!\!\!
\end{eqnarray}
By Remark \ref{rema3.1}\,(i), $\kk\gg\SS_7$ with $\SS_7$ being the number appeared in
Proposition \ref{Sect3-Lemm5}, we can apply
Proposition \ref{Sect3-Lemm5} to obtain
\begin{eqnarray}
\label{memenncnnnr}
\!\!\!\!\!\!\!\!\!\!\!\!\!\!\!\!\!\!\!\!\!\!\!\!\!\!\!\!\!\!&&
{\rm(i)\ }
|y_1|,|y_2|\stackrel{{}^{\sc\rm\eqref{mqp1234-2} }}{<}h_{_{\sc p_1,p_2}}^{^{\sc\frac{m}{m+1}}},\ \ \ \implies\ \ \
{\rm(ii)\ }\kk\stackrel{{}^{\sc\rm\eqref{bsbsbs}}}{\preceq_{\kk\ssc\,}}
h_{p_1,p_2}\stackrel{{}^{\sc\rm\eqref{dp0p1}}}{=}
|x_1|+|x_2|+|y_1|+|y_1|\stackrel{{}^{\sc\rm\eqref{bsbsbs}}}{\preceq_{\kk\ssc\,}}|x_2|,
\!\!\!\!\!\!\!\!\!\!\!\!\!\!\!\!\!\nonumber\\\!\!\!\!\!\!\!\!\!\!\!\!\!\!\!\!\!\!\!\!\!\!\!\!\!\!\!\!\!\!&&
{\rm(iii)\ }h_{p_1,p_2}\stackrel{{}^{\sc\rm\eqref{memenncnnnr}\,(ii)}}{\sim_{\kk\ssc\,}}|x_2|
\stackrel{{}^{\sc\rm\eqref{memenncnnnr}\,(i),\,(ii)}}{\sim_{\kk\ssc\,}}|x_2+y_2|
\stackrel{{}^{\sc\rm\eqref{bsbsbs}\,(ii)}}{\succeq_{\kk\ssc\,}}\kk.
\end{eqnarray}
Write
 $\TH x_{\OnE}\sTH=(x_{\OnE}+y_{\OnE})(1+\mu_{\OnE})$ for some $\mu_{\OnE}\in\C$,
then by \eqref{mqp1234-2},
\equa{MSMSee33}{\dis 
|(x_{\OnE}+y_{\OnE})\mu_{\OnE}|=|x_{\OnE}-(x_{\OnE}+y_{\OnE})|=|y_{\OnE}|\stackrel{{}^{\sc\rm\eqref{mqp1234-2} }}{<}\dH_{_{\sc p_{\ZeRo},p_{\OnE}}}^{^{\sc\frac{\ssTH m}{\ssTH m+1}}}\,
\stackrel{{}^{\sc\rm\eqref{memenncnnnr}\,(iii)}}{\sim_{\kk\ssc\,}}\,
|x_{\OnE}+y_{\OnE}|
^{\frac{\ssTH m}{\ssTH m+1}},}
i.e., $\mu_{\OnE}\preceq_{\kk\ssc\,}|x_{\OnE}+y_{\OnE}|^{-\frac1{\ssTH m+1}}
\stackrel{{}^{\sc\rm\eqref{memenncnnnr}\,(iii)}}{\preceq_{\kk\ssc\,}}
\kk^{-\frac1{\ssTH m+1}}\ll\d^2$. Thus $|\mu_{\OnE}|=O(\d)^2$.
Similarly, we can write $\TH\bar x_{\OnE}\sTH=(\bar x_{\OnE}+\bar y_{\OnE})(1+\bar \mu_{\OnE})$ with
$|\bar\mu_{\OnE}|=O(\d)^2$. Hence we have, \equa{Heneneaa}{\mbox{$\dis
X_{\OnE}\stackrel{{}^{\sc\rm\eqref{SimMMSMS}}}{=}
\frac{(\bar x_{\OnE}^{-1}x_{\OnE}) Z }{
(\bar x_{\OnE}+\bar y_{\OnE})^{-1}(x_{\OnE}+y_{\OnE})}
=\frac{(1+\mu_{\OnE})Z }{1+\bar\mu_{\OnE}}
=Z
\Big(1+O(\d)^2\Big)$,}\!\!\!\!\!\!}
which implies \eqref{ImMpP}\,(i).
By \eqref{meme},\,\eqref{memenncnnnr},%
}\NOUSE{
First by {C+LetNSoOP}\,(iv),\,\eqref{extra-equa}\,(i) and the fact from \eqref{MSmde33333} that $0<\d\ll\d_2=\ell_2^{-1}\ll\d_0\ll1$ (and in particular the choice of $\d=\ell^{-1}$ is independent of choices of $\d_0,\d_2$ by Remark \ref{Notatta-rema}), we have [see also \eqref{A0===+1}$\ssc\,$],
\begin{eqnarray}
\label{x1-small---}
\!\!\!\!\!\!\!\!\!\!\!\!\!\!\!\!\!\!\!\!\!\!\!\!&&
1+\ln(\d_2^{\d_0^{11}})\d+O(\d)^2=(\d_2^{\d_0^{11}})^{\d}=\d_2^{\d_0^{11}\d}
\stackrel{{}^{\sc\rm\eqref{extra-equa}\,(i) }}{<}
|\tildeX_1|^{\d_0^{11}\d}
\stackrel{{}^{\sc\rm\eqref{C+LetNSoOP}\,(iv) }}{\le}\ell_2^{\d_0^{11}\d}
\nonumber\\
\!\!\!\!\!\!\!\!\!\!\!\!\!\!\!\!\!\!\!\!\!\!\!\!&&
\phantom{1+\ln(\d_2^{\d_0^{11}})\d+O(\d)^2}
=(\ell_2^{\d_0^{11}})^{\d}=1+\ln(\ell_2^{\d_0^{11}})\d+O(\d)^2,
\end{eqnarray}
i.e., $|\tildeX_1|^{\d_0^{11}\d}=1+O(\d)^1$.
Then by \eqref{SimMMSMS}\,(iv),\, and  the fact in \eqref{a0BiggerThen} that $\a_0>0$,  we obtain
\eqref{ImMpP}\,(1) as follows, 
\begin{eqnarray}
\!\!\!\!\!\!\!\!\!\!\!\!\!\!\!\!\!\!&&
\label{X11111}\dis\!\!\!\!\!\!\!\!
|\widetilde X_1|\stackrel{{}^{\sc\rm\eqref{SimMMSMS}\,(iv)}}{=}(1{\sc}+{\sc}\a_0\ep)^{-1}|
\tildeX_1|^{\d_0^{11}\d}\stackrel{{}^{\sc\rm\eqref{a0BiggerThen},\,\eqref{ep-pro}}}{=}
\Big(1{\sc}-{\sc}\a_0\ep{\sc}+{\sc}O(\ep)^2\Big)|\tildeX_1|^{\d_0^{11}\d}
\nonumber\\
\!\!\!\!\!\!\!\!\!\!\!\!\!\!\!\!\!\!&&
\phantom{|\widetilde X_1|\!}{\sc}<{\sc}\
|\tildeX_1|^{\d_0^{11}\d}{\sc}={\sc}
1{\sc}+{\sc}O(\d)^1
.\!\!\!\!\!\!\end{eqnarray}
}\NOUSE{
To prove \eqref{ImMpP}\,(2),  first by Remark \ref{rema3.1}\,(ii) and \eqref{MSmde33333+}, we have $1\sim_{\kk\ssc\,}\ell_1\sim_{\kk\ssc\,}\ell_2\sim_{\kk\ssc\,}\ell_3\sim_{\kk\ssc\,}1$.
Since the choice of $\kk$ is independent of $\d_0,\d_0,\d_2,\d_3,\d$ by Remark \ref{rema3.1}\,(i),
exactly similar to the proof of \eqref{A1-A2-sim-ell2}, we obtain from
\eqref{C+LetNSoOP}\,(iv),\,\eqref{extra-equa},\,\eqref{+C+LetNSoOP} the following
%
\equa{Bb11111}{\mbox{(i) $|A_{1}|\stackrel{{}^{\rm\eqref{+C+LetNSoOP}}}{\sim_{\kk\ssc\,}}1
\stackrel{{}^{\sc\rm\eqref{+C+LetNSoOP}}}{\sim_{\kk\ssc\,}}|A_2|$, \ \ \ \ \ \ (ii) $|\tildeX_{\ZeRo}|
\stackrel{{}^{\sc\rm \eqref{extra-equa}\,(i),\,\eqref{C+LetNSoOP}\,(iv)}}{\sim_{\kk\ssc\,}}1
\stackrel{{}^{\sc\rm\eqref{extra-equa}\,(ii)}}{\sim_{\kk\ssc\,}}|Z|$.}}
From \eqref{Bb11111}\,(ii), exactly similar to the proof of \eqref{Z===X2-simm}, we have
\equa{|X_2|sim|Z|}{|X_2|\sim_{\kk\ssc\,}|Z|\sim_{\kk\ssc\,}1.}
Thus by \eqref{X11111},\,
\eqref{TaKa},\,\eqref{SimMMSMS}%
,
we have
\equa{WeHAHAHSDde}{\mbox{${\rm(i)\ }|x_{\ZeRo}|\stackrel{{}^{\sc\rm\eqref{TaKa},\,\eqref{SimMMSMS}}}{=}\kk|\widetilde X_{\ZeRo}|
\stackrel{{}^{\sc\rm\eqref{X11111},\,\eqref{Bb11111}\,(ii)}}{\preceq
_{\kk\ssc\,}}\kk,\ \ \ \ \ \ {\rm(ii)\ }|x_2|\stackrel{{}^{\sc\rm\eqref{TaKa},\,\eqref{SimMMSMS}}}{=}\kk|\widetilde X_2|\stackrel{{}^{\sc\rm\eqref{SimMMSMS}\,(v),\,
\eqref{|X_2|sim|Z|}}}
{\sim_{\kk\ssc\,}}\kk$.}}
Hence 
we must have
\begin{eqnarray}
\label{MustHave}
\!\!\!\!\!\!\!\!\!\!\!\!\!\!\!\!\!\!\!\!&&
{\rm(i)\ }
 \dH_{p_{\ZeRo},p_{\OnE}}
\stackrel{{}^{\sc\rm\eqref{dp0p1},\,\eqref{mqp1234-2} }}{\sim_{\kk\ssc\,}}|x_{\ZeRo}|+|x_{\OnE}|\stackrel{{}^{\sc\rm\eqref{WeHAHAHSDde}}}{\sim_{\kk\ssc\,}}\kk,\ \ \ \ \ \ \ \ \
{\rm(ii)\ }
|y_{\OnE}|\stackrel{{}^{\sc\rm\eqref{mqp1234-2} }}{\prec_{\kk\ssc\,}}\dH_{p_{\ZeRo},p_{\OnE}},
\nonumber\\
\!\!\!\!\!\!\!\!\!\!\!\!\!\!\!\!\!\!\!\!&&
{\rm(iii)\ }
 |x_{\OnE}+y_{\OnE}|\stackrel{{}^{\sc\rm\eqref{WeHAHAHSDde}\,(ii),\,\eqref{MustHave}\,(i),\,(ii)}}{\sim_{\kk\ssc\,}}|x_{\OnE}|
 \stackrel{{}^{\sc\rm\eqref{WeHAHAHSDde}\,(ii)}}{\sim_{\kk\ssc\,}}\kk.
\end{eqnarray}
Write
 $\TH x_{\OnE}\sTH=(x_{\OnE}+y_{\OnE})(1+\mu_{\OnE})$ for some $\mu_{\OnE}\in\C$,
then by \eqref{mqp1234-2},
\equa{MSMSee33}{\dis 
|(x_{\OnE}+y_{\OnE})\mu_{\OnE}|=|x_{\OnE}-(x_{\OnE}+y_{\OnE})|=|y_{\OnE}|\stackrel{{}^{\sc\rm\eqref{mqp1234-2} }}{<}\dH_{_{\sc p_{\ZeRo},p_{\OnE}}}^{^{\sc\frac{\ssTH m}{\ssTH m+1}}}\,
\stackrel{{}^{\sc\rm\eqref{MustHave}}}{\sim_{\kk\ssc\,}}\,
|x_{\OnE}+y_{\OnE}|
^{\frac{\ssTH m}{\ssTH m+1}},}
i.e., $\mu_{\OnE}\preceq_{\kk\ssc\,}|x_{\OnE}+y_{\OnE}|^{-\frac1{\ssTH m+1}}
\stackrel{{}^{\sc\rm\eqref{MustHave}}}{\sim_{\kk\ssc\,}}
\kk^{-\frac1{\ssTH m+1}}\ll\d^3$. Thus $|\mu_{\OnE}|=O(\d)^3$.
Similarly, we can write $\TH\bar x_{\OnE}\sTH=(\bar x_{\OnE}+\bar y_{\OnE})(1+\bar \mu_{\OnE})$ with
$|\bar\mu_{\OnE}|=O(\d)^3$. Hence we have, \equa{Heneneaa}{\mbox{$\dis
\widetilde X_{\OnE}=
\frac{(\bar x_{\OnE}^{-1}x_{\OnE})\widetilde Z }{
(\bar x_{\OnE}+\bar y_{\OnE})^{-1}(x_{\OnE}+y_{\OnE})}
=\frac{(1+\mu_{\OnE})\widetilde Z }{1+\bar\mu_{\OnE}}
=\widetilde Z
\Big(1+O(\d)^3\Big)=\widetilde Z+O(\d)^3$,}\!\!\!\!\!\!}
 where the last equality
follows from \eqref{extra-equa}\,(ii) as
$\ell_3\ll\ell=\d^{-1}$.
This proves
}\NOUSE{
Note from \eqref{extra-equa},\,\eqref{MSmde33333},\,\eqref{SimMMSMS}\,(iii)--(vi) that we have
\equa{OSOSOSWL}{\dis
\d a=O(\d)^1\mbox{ \ for any \ }a\in T_0:=\{X_1,X_2,Z,\widetilde X_1,\widetilde X_2,\widetilde Z\}.
}
From this and by \eqref{extra-equa}\,(iv), we have \eqref{ImMpP}\,(2) as follows,
\equa{Z-2-2-2}{\dis
Z=\widetilde Z^{\ell_0^{11}}=\widetilde X_2^{\ell_0^{11}}\big(1+O(\d)^3\big)^{\ell_0^{11}}=X_2(1+O(\d)^3)=X_2+O(\d)^3.}
By \eqref{MSmde33333},\,\eqref{TaKa},\,\eqref{SimMMSMS},  we see that
 $(|x_{\OnE}|+|x_{\OnE}+y_{\OnE}|)\ep^3<\ep^2$, which with \eqref{C+LetNSoOP}\,(ii) implies \eqref{ImMpP}\,(3) [the equality there will follow from \eqref{A0===+1}$\ssc\,$].
}
\NOUSE{Denote
\begin{eqnarray}
\label{A0===}
&\!\!\!\!\!\!\!\!\!\!\!\!\!\!\!\!\!\!\!\!\!\!\!\!\!\!\!\!\!\!&
T_{1}{\ssc}:=
{\ssc}\Big\{{\ssc}\tildeX_{\ZeRo},\,{\ssc}X_{\OnE},\,
{\ssc}Z ,\,X_1',\,{\ssc}A_{\rOnE},\,A_2,\,
\d_0+\frac{(1-\d_0)X_1'^{\ell_0(1 + \d_0)}}{Z^{\ell_0^4(1 + \d_0)}},\,\frac14+\frac{3A_1^{4\ell_0^4}(\a_0X_1)^2Z^4}{4X_2^2}
\Big\}.
\end{eqnarray}
By
\eqref{C+LetNSoOP}\,(i), we have [
note  from \eqref{B1-ep} that $a_2a_1^{-1}+\d^3<3\ell_0^4\ssc\,$],
\equa{A1andA2}{
\dis{\rm(i)\ }1\le|A_1|\le\ell_1,\ \ \ \ \ {\rm(ii)\ }1\le|A_2|\le \ell_1^{3\ell_0^4}.
}
From this, Lemma \ref{Case6-lemm} and \eqref{C+LetNSoOP}\,(iv),\,(v),
we can easily observe the following,
\equa{ObSSSSS}{\mbox{$|\a|^{\pm1}\le\ell_2^{\ell_2}$ \ for any  $\a\in T_1$.}}
Thus, 
%
%
\begin{eqnarray}
\label{A0===+1}
&\!\!\!\!\!\!\!\!\!\!\!\!\!\!\!\!\!\!\!\!&
\a^{0+O(\d)^1}=1+O(\d)^1,\ \ \a\Big(1+O(\d)^1\Big)=\a+O(\d)^1\mbox{  \ \ for all  }\a\in
T_{1}\mbox{  or }\a^{-1}\in T_{1}.
\end{eqnarray}
Thus $(\a_0X_1)^{-1}{\ssc}={\ssc}X_1^{-1}{\ssc}+{\ssc}O(\ep)^1{\ssc}={\ssc}X_1^{-1}{\ssc}+{\ssc}O(\d)^3$ by \eqref{a0BiggerThen} and we have \eqref{ImMpP}\,(2),\,(3),\,(8) by \eqref{LetNSoOP----1},\,\eqref{equa-Case6-lemm}\,(iv)%
. 
}%
\NOUSE{Let  $(p_1,p_2)\in\ol V_0$.
If $\xX_1\le\ep$ or $\xX_1\ge\nn$, we obtain from
 \eqref{LetNSoOP----1}\,(ii) that $|A_2|\gg\nn$, a contradiction with \eqref{A1-A2-cond}\,(b).
This proves \eqref{ImMpP}\,(i).
Then
by }%
\NOUSE{notation \eqref{TaKa},\,\eqref{SimMMSMS}, and \eqref{A1-A2-cond}\,(a),\,\eqref{ImMpP}\,(i),
we obtain that  $|x_1|=\xX_1\kk<\nn\kk,\,|x_2|=\xX_2\kk\le\nn\kk$.
We must have $|y_2|<2\nn\kk$ by Proposition \ref{Sect3-Lemm5} (cf.~Remark \ref{Remmma}).%
}%
\NOUSE
{first 
by \eqref{C+ToSayas+1}\,(d),\,(e) and  notation \eqref{denote-t2}, we have $\eE_0\le\tilde \xX_1,\xX_2\le\nn_0$ (recalling that $\eE_0=\nn_0^{-1}$). Thus
\begin{eqnarray}
\label{x1-small}
&\!\!\!\!\!\!\!\!\!\!\!\!\!\!\!\!\!\!\!\!\!\!\!\!\!\!\!\!\!\!\!\!\!&
{\rm(i)\ }
\xX_1\stackrel{{}^{\sc\rm\eqref{LetNSoOP----1}\,(iii)}}{=}
\a_0^{-\frac1{100}}\xX_2^{\frac15}
\tilde \xX_1^{\frac1{100}}
\stackrel{{}^{\sc\rm\eqref{LetNSoOP----1}\,(iii)}}{<}\xX_2^{\frac15}
\tilde \xX_1^{\frac1{100}}
\stackrel{{}^{\sc\rm\eqref{C+ToSayas+1}\,(d),\,(e)}}{\le}
\nn_0,\!\!\!\!\!\!\!\!\!\!\!\!\!\!\!\!\!\!\!\!\!\!\!
\nonumber\\
&\!\!\!\!\!\!\!\!\!\!\!\!\!\!\!\!\!\!\!\!\!\!\!\!\!\!\!\!\!\!\!\!\!&
{\rm(ii)\ }
|x_1|\stackrel{{}^{\sc\rm\eqref{TaKa},\,\eqref{SimMMSMS}}}{=}\kk|X_1|=
\kk\xX_1\stackrel{{}^{\sc\rm\eqref{x1-small}\,(i)}}{\le}
\kk\nn_0,
\nonumber\\
&\!\!\!\!\!\!\!\!\!\!\!\!\!\!\!\!\!\!\!\!\!\!\!\!\!\!\!\!\!\!\!\!\!&
{\rm(iii)\ }|x_2|\stackrel{{}^{\sc\rm\eqref{TaKa},\,\eqref{SimMMSMS}}}{=}\kk|X_2|=\kk\xX_2
\stackrel{{}^{\sc\rm\eqref{C+ToSayas+1}\,(e)}}{\le}\kk\nn_0.
\end{eqnarray}
We claim
\equa{y1-y2-isSmall}{\dis
\max\{|y_1|,|y_2|\}\le\kk\nn_0.
}
To prove the claim, say (the proof for the case $|y_1|=\max\{|y_1|,|y_2|\}$ is exactly similar),
\equa{ma-y12}{\dis
|y_2|=\max\{|y_1|,|y_2|\}>\kk\nn_0.
}
Then
\equa{yyyyyaaa2==}
{\dis
\kk\nn_0\stackrel{{}^{\sc\rm\eqref{ma-y12}}}{<}|y_2|\le|x_1|+|x_2|+|y_1|+|y_2|\stackrel{{}^{\sc\rm\eqref{dp0p1}}}{=}h_{p_1,p_2}\stackrel{{}^{\sc\rm\eqref{x1-small}\,(ii),\,(iii),\,\eqref{ma-y12}}}{\le}4|y_2|.
}
Note that we have complete freedom in choosing $\nn_0\gg\nn$ by \eqref{MSmde33333}, in particular we may assume $\kk\nn_0>\SS_7$ with $\SS_7$ being the number in Proposition \ref{Sect3-Lemm5} so that by the fact in \eqref{yyyyyaaa2==} that $h_{p_1,p_2}>\kk\nn_0$, we can apply \eqref{mqp1234-2} to obtain
\equa{y2-conttrraa}{\dis
|y_2|\stackrel{{}^{\sc\rm\eqref{mqp1234-2} }}{<}h_{_{\sc p_1,p_2}}^{^{\sc\frac m{m+1}}}\stackrel{{}^{\sc\rm\eqref{yyyyyaaa2==}}}{\le} 4^{\frac m{m+1}}|y_2|^{\frac m{m+1}},
\ \ \ \implies\ \ \ |y_2|<4^m,
}
which is a contradiction with \eqref{ma-y12}. This proves \eqref{y1-y2-isSmall}.
%
} 
$\!\!\!\!$Then by  \eqref{meme},\,\eqref{equa-Case6-lemm}\,(iii),
we obtain immediately
that
 $(|x_{\OnE}|+|x_{\OnE}+y_{\OnE}|)\ep^3<\ep^2$
, which with \eqref{C+LetNSoOP}\,(ii)
implies \eqref{ImMpP}\,(i)%
.

Recall from \eqref{equa-Case6-lemm}\,(iv) that $Z=X_2+O(\d)^3$, 
we can rewrite $A_2$ in \eqref{LetNSoOP----1-re-give}\,(ii), up to $O(\d)^2$, as the following, 
\begin{eqnarray}
\!\!\!\!\!\!\!\!\!\!\!\!\!\!\!\!\!\!\!\!\!\!\!\!\!\!\!\!\!\!\!\!\!\!\!\!
&&
A_2
\stackrel{{}^{\sc\rm\eqref{LetNSoOP----1-re-give}\,(ii)}}{=}
\frac{X_2^{\ell_0+1}}{Z}\Big (2- \frac{\widetilde  X_1}{A_1^2}\Big)
 \stackrel{{}^{\sc\rm\eqref{equa-Case6-lemm}\,(iv),\,\eqref{(T0o(eE)1}}}{=}
X_2^{\ell_0}\Big (2- \frac{\widetilde  X_1}{A_1^2}\Big)
+O(\d)^2.
\!\!\!\!\!\!\!\!
\label{We0o0o0o++1}
\end{eqnarray}
%
%
{%
We use \eqref{C+LetNSoOP} to obtain \eqref{ImMpP}\,(ii),\,(iii) as follows
(recalling that $0<\d\ll\d_0\ll1$),
\begin{eqnarray}
\label{We0o0o0o}
\!\!\!\!\!\!\!\!\!\!\!\!\!\!\!\!\!\!\!\!\!\!\!\!\!\!\!\!\!\!\!\!&&
{\rm(i)\ }
\xX_2
\stackrel{{}^{\sc\rm\eqref{C+LetNSoOP}\,(iii)}}{\le}
\Big(\frac{(1+\d)\aA_1^{-(2+2\d_0)}}{1-\d}\Big)^{\frac1{2}}
=
\aA_1^{-(1+\d_0)}+O(\d)^1
\stackrel{{}^{\sc\rm\eqref{A1-A2-cond}\,(a)}}{<}\ell_1^2,
\!\!\!\!\!\!\!\!\!\!\!\!\!
\nonumber\\
\!\!\!\!\!\!\!\!\!\!\!\!\!\!\!\!\!\!\!\!\!\!\!\!\!\!\!\!\!\!\!\!&&
{\rm(ii)\ }
\tilde\xX_1\stackrel{{}^{\sc\rm\eqref{C+LetNSoOP}\,(iii)}}{\le}
(1+\d)\xX_2^{-1}\aA_1^{-(2+\d_0)}
\stackrel{{}^{\sc\rm\eqref{A1-A2-cond}\,(a),\,\eqref{ImMpP}\,(i)}}{\le}
\ell_1^3.
\nonumber
\\
\!\!\!\!\!\!\!\!\!\!\!\!\!\!\!\!\!\!\!\!\!\!\!\!&&
{\rm(iii)\ }
\tilde\xX_1\stackrel{{}^{\sc\rm\eqref{C+LetNSoOP}\,(iii)}}{\ge}
(1-\d)\xX_2\aA_2^{\d_0}
\stackrel{{}^{\sc\rm\eqref{A1-A2-cond}\,(a),\,\eqref{ImMpP}\,(i)}}{\ge}
\d_1.
\!\!\!\!\!\!\!\!\!\!\!\!\!\!\!\!\!
\end{eqnarray}
Thus we have \eqref{ImMpP}\,(ii),\,(iii).
}%
\NOUSE{
\begin{eqnarray}
\label{We0o0o0o}
&&\!\!\!\!\!\!\!\!\!\!\!\!\!\!\!\!\!\!\!\!\!\!\!\!\!\!\!\!\!\!\!\!
\frac{A_2}{A_1}
\stackrel{{}^{\sc\rm\eqref{LetNSoOP----1}\,(i),\,(ii)}}{=}
\frac{X_2\Big(
\frac23{\ssc\!}+{\ssc\!}\frac1{3\widetilde X_1^3X_2^3}
\Big)}{Z\widetilde X_1^{11} X_2^{11}}
\stackrel{{}^{\sc\rm\eqref{equa-Case6-lemm}\,(iv)}}{=}
\frac{
\frac23{\ssc\!}+{\ssc\!}\frac1{3\widetilde X_1^3X_2^3}
}{\widetilde X_1^{11} X_2^{11}}{\ssc\!}+{\ssc\!}O(\d)^2
\stackrel{{}^{\sc\rm\eqref{LetNSoOP----1}\,(i),\,(ii)}}{=}
\frac{
\frac43{\ssc\!}-{\ssc\!}\frac{A_1\widetilde X_1^2}{3X_2^{11}}
}{\widetilde X_1^{11} X_2^{11}}{\ssc\!}+{\ssc\!}O(\d)^2
\NOUSE{,
\nonumber\\
&&\!\!\!\!\!\!\!\!\!\!\!\!\!\!\!\!\!\!\!\!\!\!\!\!
{\rm(ii)\ }
\frac{A_1}{A_2}\stackrel{{}^{\sc\rm\eqref{equa-Case6-lemm}\,(iv),\,\eqref{LetNSoOP----1}\,(i),\,(ii)}}{=}
\frac{\widetilde X_1^{11} X_2^{11}}{\Big(
\frac23+\frac1{3\widetilde X_1^3X_2^3}
\Big)}+O(\d)^2
=
\frac{\widetilde X_1^{11} X_2^{11}}{\Big(
\frac43-\frac{A_1\widetilde X_1^2}{3X_2^{11}}
\Big)}+O(\d)^2
}
.
\end{eqnarray}
}
\NOUSE{
Now we have
\equa{A1-A2222222}{\dis\!\!\!\!\!\!\!\!
{\rm(i)\, }
\aA_1\stackrel{{}^{\sc\rm\eqref{LetNSoOP----1}\,(ii)}}{=}
\xX_2^2\xX_1^{-\nn_1}\zZ
\stackrel{{}^{\sc\rm\eqref{ImMpP}\,(i)}}{=}
\xX_2\xX_1^{-\nn_1}\Big(1{\ssc\!}+{\ssc\!}O(\d)^2\Big),\, \
{\rm(ii)\, }
\aA_2\stackrel{{}^{\sc\rm\eqref{LetNSoOP----1}\,(i),\,(iii)}}{=}\frac1{\xX_1} \Big|
\lL{\ssc\!}+{\ssc\!}1{\ssc\!} -{\ssc\!}\frac{\lL}{X_2X_1}\Big|.\!\!\!\!\!\!
}
}

Now we first assume the last inequality of \eqref{ImMpP}\,(iv) does not hold. Then
 by \eqref{C+LetNSoOP}\,(iii),
 \equa{tX1-is=}{\dis
\tilde\xX_1
=
(1+\d)\xX_2^{-1}\aA_1^{-\d_0}
.}
Using the definition of $A_1$ in \eqref{LetNSoOP----1-re-give}\,(i), up to $O(\d)^2$, we have
,
\begin{eqnarray}
\label{x2-00ro4m4m4-}
\!\!\!\!\!\!\!\!\!\!\!\!\!\!\!\!\!\!\!\!\!\!\!\!\!\!\!\!\!\!\!\!\!\!\!\!\!\!\!\!\!\!\!\!\!\!&&
{\rm(i)\ }
2
\stackrel{{}^{\sc\rm\eqref{LetNSoOP----1-re-give}\,(i)}}{\le}
\frac1{\tilde \xX_1}
+\frac{\aA_1}{\xX_2^{\ell_0}}
\stackrel{{}^{\sc\rm\eqref{tX1-is=}}}{=}
\frac{\xX_2\aA_1^{\d_0}}{1+\d}
+\frac{\aA_1}{\xX_2^{\ell_0}}
\\ \nonumber
\!\!\!\!\!\!\!\!\!\!\!\!\!\!\!\!\!\!\!\!\!\!\!\!\!\!\!\!\!\!\!\!\!\!\!\!\!\!\!\!\!\!\!\!\!\!&&
\ \ \ \ \ \ \ \ \ \ 
=(1-\d)\xX_2\aA_1^{\d_0}
+(\xX_2\aA_1^{\d_0})^{-\ell_0}\aA_1^2
\stackrel{{}^{\sc\rm\eqref{A1-A2-cond}\,(a)}}{\le}
\gamma_1(\yY):=(1-\d)\yY+\yY^{-\ell_0}
\mbox{ \ [up to $O(\d)^2\ssc\,$]},
\!\!\!\!\!\!\!\!\!\!\!\!\!\!\!\!\!\!\!\!\!\!\!\!\!\!\!\!
\nonumber\\
\!\!\!\!\!\!\!\!\!\!\!\!\!\!\!\!\!\!\!\!\!\!\!\!\!\!\!\!\!\!\!\!\!\!\!\!\!\!\!\!\!\!\!\!\!\!&&
{\rm(ii)\ }
\yY := \xX_2\aA_1^{\d_0}
.
\nonumber\!\!\!\!\!\!\!\!\!\!\!\!\!\!\!\!\!\!\!\!\!\! \!\!\!\!\!\!\!\!
\end{eqnarray}
Noting from
\eqref{A1-A2-cond}\,(a),\,\eqref{ImMpP}\,(i) that the last term
in the first line of
\eqref{x2-00ro4m4m4-}\,(i) is $\le1$ [up to $O(\d)^2\ssc\,$], we must have $\yY\ge1+\d$ [up to $O(\d)^2\ssc\,$].
Regarding $\gamma_1(\yY)$ as a function on $\yY$, by solving  $\frac{d\gamma_1}{d\yY}=0$, we see that
the function $\gamma_1(\yY)$ 
is
strictly decreasing when $1\le\yY\le\yY_0$ and strictly increasing when $\yY\ge\yY_0$, where
$\yY_0=\big(\frac{\ell_0}{1-\d}\big)^{\frac{\d_0}{1+\d_0}}$.
\NOUSE{
(recalling that $\ln(\cdot)$ is the natural  logarithmic function; also noting that $\d\ll\d_0$)
\begin{eqnarray}
\label{uuuauyyY}
\!\!\!\!\!\!\!\!\!\!\!\!\!\!\!\!\!\!\!\!\!\!&&
\yY_0=\Big(\ell_0^5(1-\d_0^7)\Big)^{\d_0^7}
=1+5\ln(\ell_0)\d_0^2+O(\d_0)^3.
\end{eqnarray}
}%
Since $\gamma_1(2-\d_0^2)=2-\d_0^2+O(\d_0)^{100}<3$ [noting that $2^{-\ell_0}\ll\d_0^{100}$ and $\d\ll\d_0\ssc\,$],
in order for \eqref{x2-00ro4m4m4-}\,(i) to hold, we must have 
(i) below, and
so the second term in \eqref{x2-00ro4m4m4-}\,(i) is an $O(\d_0)^{100}$ element, which gives (ii) below. Therefore, we have, 
up to $O(\d)^1$ [noting that we can conduct the following computations up to $O(\d)^1$ as $\d\ll\d_0$
and recalling that $\ln(\cdot)$ is the natural  logarithmic function],
\begin{eqnarray}
\label{uuuauyyY}
\!\!\!\!\!\!\!\!\!\!\!\!\!\!\!\!\!\!\!\!\!\!&&
{\rm(i)\ }\yY>2-\d_0^2,\ \ \ \ \ \ \ \ 
{\rm(ii)\ }2
\stackrel{{}^{\sc\rm\eqref{LetNSoOP----1-re-give}\,(i)}}{\ge}
\frac1{\tilde \xX_1}-\frac{\aA_1}{\xX_2^{\ell_0}}
=\yY+O(\d_0)^{100},
\\\nonumber
\!\!\!\!\!\!\!\!\!\!\!\!\!\!\!\!\!\!\!\!\!\!&&
{\rm(iii)\ }\yY\stackrel{{}^{\sc\rm\eqref{uuuauyyY}\,(i),\,(ii)}}{=}2+O(\d_0)^{100},
\ \ \ \ \ \ \ \ 
{\rm(iv)\ }\yY
\stackrel{{}^{\sc\rm\eqref{x2-00ro4m4m4-}\,(ii)}}{=}
\xX_2\aA_1^{\d_0}
\stackrel{{}^{\sc\rm\eqref{We0o0o0o}\,(i)}}{\le}\aA_1^{-1},
\\\nonumber
\!\!\!\!\!\!\!\!\!\!\!\!\!\!\!\!\!\!\!\!\!\!&&
 {\rm(v)\ }\xX_2\stackrel{{}^{\sc\rm\eqref{A1-A2-cond}\,(a)}}{\ge}\xX_2\aA_1^{\d_0}
=\yY=2+O(\d_0)^{100},
\\\nonumber
\!\!\!\!\!\!\!\!\!\!\!\!\!\!\!\!\!\!\!\!\!\!&&
{\rm(vi)\ }
\frac{\tilde\xX_1 }{\aA_1^2}
\stackrel{{}^{\sc\rm\eqref{tX1-is=}}}{=}
\aA_1^{-2-\d_0}\xX_2^{-1}
\stackrel{{}^{\sc\rm\eqref{x2-00ro4m4m4-}\,(ii)}}{=}
\aA_1^{-2-2\d_0}\yY^{-1}
\stackrel{{}^{\sc\rm\eqref{uuuauyyY}\,(iv)}}{\ge}
\yY^{1+2\d_0}
\\\nonumber
\!\!\!\!\!\!\!\!\!\!\!\!\!\!\!\!\!\!\!\!\!\!&&
\phantom{{\rm(vi)\ }\ \  }
\stackrel{{}^{\sc\rm\eqref{uuuauyyY}\,(iii)}}{=}
2^{1+2\d_0}+O(\d_0)^{100}=2+4\ln(2)\d_0+O(\d_0)^2.
\end{eqnarray}
The importance of the fact in \eqref{uuuauyyY}\,(vi) that it is a number strictly bigger than $2$
is that we can use it and the
 expression of $A_2$ in \eqref{We0o0o0o++1} to obtain, up to $O(\d)^1$,
\begin{eqnarray}
\label{Futtth}
\!\!\!\!\!\!\!\!\!\!\!\!\!\!\!\!\!\!\!\!\!\!&&
\frac12+O(\d_0)^1
\stackrel{{}^{\sc\rm\eqref{uuuauyyY}\,(iii),\,(iv)}}{\ge}
\aA_1
\stackrel{{}^{\sc \rm\eqref{cont-B2}}}{=}
\aA_2
\stackrel{{}^{\sc \rm\eqref{We0o0o0o++1}}}{\ge}
\xX_2^{\ell_0}\Big (\frac{\tilde \xX_1}{\aA_1^2}-2\Big)
\nonumber\\
\!\!\!\!\!\!\!\!\!\!\!\!\!\!\!\!\!\!\!\!\!\!&&
\ \ \ \ \ \ \ \ \ \ \ \ \ \ \ \ 
\stackrel{{}^{\sc \rm\eqref{uuuauyyY}\,(v),\,(vi)}}{\ge}
\Big(2+O(\d_0)^{100}\Big)^{\ell_0}
\Big(4\ln(2)\d_0+O(\d_0)^2\Big)\gg2
,
\end{eqnarray}
which is a contradiction.
\NOUSE{%
Equ.~\eqref{uuuauyyY}\,(iii) contradiction with \eqref{A1-A2-cond}\,(a)
[noting that though we give our above evaluation up to $O(\d)^1$, we can still obtain
that $\aA_1$ is strictly bigger than $1$].
}%
This proves the first inequality of
\eqref{ImMpP}\,(iv).
\NOUSE{%
Thus
\equa{yyIs-Big}
{\dis\frac{51}{50}<
\yY
\stackrel{{}^{\sc\rm\eqref{x2-00ro4m4m4-}\,(ii)}}{=}
\xX_2^{6\ell_0}\aA_1^{-\frac34+\d_0}
\stackrel{{}^{\sc\rm\eqref{We0o0o0o}\,(i)}}{\le}\aA_1^{\frac18+\d_0}
\mbox{ \ [up to $O(\d)^2\ssc\,$]}
.
}
We claim
\equa{msmsm---}{\dis
\yY<\aA_1^{2\d_0}.
}
Assume conversely that it does not hold, i.e., we have
\equa{xX2---2-2}{\dis
\xX_2>
\aA_1^{\frac{\d_0}8+\frac{\d_0^2}{6}}.
}
Then we can use \eqref{We0o0o0o++1} to obtain the following, up to $O(\d)^1$
[note that we can conduct the following computations up to $O(\d)^1$ as $\d\ll\d_0$,
and
 that although
$\frac{A_1 \widetilde X_1^4}{X_2^{2\ell_0}}
,\frac{\widetilde X_1 }{X_2^{3\ell_0}}$ are not necessarily  real numbers, we are able to determine  in
\eqref{uuuauyyY}\,(ii),\,(v) that
they are  positive numbers up to $O(\d_0)^{100}$; note also that the last equality of
(iii) below follows from the definition of $A_1$ in \eqref{LetNSoOP----1}\,(i)$\ssc\,$],
%
Denote $\eta_i=|\Theta_i|$ for $i=1,2$.
We want to use \eqref{tX1-is=},\,\eqref{uuuauyyY}\,(ii),\,(v) to solve
$\tilde\xX_1,\tilde\xX_2:=\xX_2^{\ell_0},\aA_1$ in terms of $\theta_1,\theta_2$ as follows [up to $O(\d)^1\ssc\,$],
\begin{eqnarray}
\label{Detmmm--}
&\!\!\!\!\!\!\!\!\!\!\!\!\!\!\!\!\!\!\!\!&
{\rm(i)\ }1
\stackrel{{}^{\sc\rm\eqref{uuuauyyY}\,(ii)}}{=}
\frac{\tilde\xX_1^4}{\th_1\aA_1\tilde\xX_2^2}=
\stackrel{{}^{\sc\rm\eqref{tX1-is=}}}{=}
\frac{\tilde \xX_2^{34}}{\th_1\aA_1^{4- 4 \d_0}},
\!\!\!\!\!\!\!\!\!\!\!\!\!\!
\\\nonumber
&\!\!\!\!\!\!\!\!\!\!\!\!\!\!\!\!\!\!\!\!&
{\rm(ii)\ }
1
\stackrel{{}^{\sc\rm\eqref{uuuauyyY}\,(v)}}{=}
\frac{\tilde\xX_1}{\th_2\tilde\xX_2^3}
\stackrel{{}^{\sc\rm\eqref{tX1-is=}}}{=}
\frac{\tilde\xX_2^6}{\th_2\aA_1^{\frac34-\d_0}},
\!\!\!\!\!\!\!\!\!\!\!\!\!\!
\nonumber\\
&\!\!\!\!\!\!\!\!\!\!\!\!\!\!\!\!\!\!\!\!&
{\rm(iii)\ }
1
\stackrel{{}^{\sc\rm\eqref{Detmmm--}\,(i),\,(ii)}}{=}
\Big(\frac{\tilde \xX_2^{34}}{\th_1\aA_1^{4- 4 \d_0}}\Big)^3
\Big(\frac{\th_2\aA_1^{\frac34-\d_0}}{\tilde\xX_2^6}\Big)^{17}
=\frac{\th_2^{17}
\aA_1^{\frac34- 5 \d_0}}{\th_1^3},
\!\!\!\!\!\!\!\!\!\!\!\!\!\!
\nonumber\\
&\!\!\!\!\!\!\!\!\!\!\!\!\!\!\!\!\!\!\!\!&
{\rm(iv)\ }
\aA_1
\stackrel{{}^{\sc\rm\eqref{Detmmm--}\,(iii)}}{=}
\Big(\frac{\th_1^3}{\th_2^{17}}\Big)^{\frac1{\frac34-5\d_0}}
=
\th_1^{\frac{12}{3-20\d_0}}\th_2^{-\frac{68}{3-20\d_0}},
\!\!\!\!\!\!\!\!\!\!\!\!\!\!
\nonumber\\
&\!\!\!\!\!\!\!\!\!\!\!\!\!\!\!\!\!\!\!\!&
{\rm(v)\ }
\tilde\xX_2
\stackrel{{}^{\sc\rm\eqref{Detmmm--}\,(ii),\,(iv)}}{=}
\Big(\th_2\Big(\th_1^{\frac{12}{3-20\d_0}}\th_2^{-\frac{68}{3-20\d_0}}\Big)^{\frac34-\d_0}\Big)^{\frac16}
=
\th_1^{\frac{3 - 4 \d_0}{2(3 - 20 \d_0)}}
\th_2^{-\frac{8 (1 - \d0)}{3 - 20 \d_0}},
\end{eqnarray}
Now we can use definition of $A_1$ in
\eqref{LetNSoOP----1}\,(i) to obtain
\begin{eqnarray}
\label{AAA111}
&\!\!\!\!\!\!\!\!\!\!\!\!\!\!\!\!\!\!\!\!\!\!\!\!\!\!\!\!\!\!\!\!\!\!\!\!\!\!\!\!\!\!\!\!\!\!\!\!\!\!\!
&
1
\stackrel{{}^{\sc\rm\eqref{LetNSoOP----1}\,(i)}}{=}
\aA_1^{-1}\tilde \xX_1^3 \tilde\xX_2
 \Big|5 - \frac{4\widetilde X_1}{ X_2^{3\ell_0}}\Big|
\stackrel{{}^{\sc\rm\eqref{tX1-is=},\,\eqref{uuuauyyY}\,(v)}}{=}
\aA_1^{-\frac{13}{4}+3\d_0}\tilde\xX_2^{28}\Big|5-4
\Big(\frac{5(1+\d_0)}{4+5\d_0}+O(\d_0)^{100}\Big)\Big|
\nonumber\\
&\!\!\!\!\!\!\!\!\!\!\!\!\!\!\!\!\!\!\!\!\!\!\!\!\!\!\!\!\!\!\!\!\!\!\!\!\!\!\!\!\!\!\!\!\!\!\!\!\!\!\!
&
\stackrel{{}^{\sc\rm\eqref{Detmmm--}\,(iv),\,(v)}}{=}
\Big(\th_1^{\frac{12}{3-20\d_0}}\th_2^{-\frac{68}{3-20\d_0}}\Big)^{-\frac{13}{4}+3\d_0}
\Big(\th_1^{\frac{3 - 4 \d_0}{2(3 - 20 \d_0)}}
\th_2^{-\frac{8 (1 - \d0)}{3 - 20 \d_0}}\Big)^{28}
\Big(\frac{5\d_0}{4 + 5 \d_0}+O(\d_0)^{100}\Big)=
\th_1\th_2^{-1}
\Big(\frac{5\d_0}{4 + 5 \d_0}+O(\d_0)^{100}\Big)
\end{eqnarray}
Obviously \eqref{uuuauyyY}\,(iv) implies that $\aA_1\le1$ [up to $O(\d)^1\ssc\,$], which is a contradiction with
the fact in \eqref{uuuauyyY}\,(i) that $\aA_1>1+\frac{\d_0^4}{20}+O(\d_0)^5$.
This proves the first inequality of
\eqref{ImMpP}\,(ii).
}\NOUSE{
Then up to $O(\d)^1$ (one can observe that the following evaluation does not depend on $\d$ as $0<\d\ll\d_0$), we also have
\begin{eqnarray}
\label{ii-LetNSoOP----1}
&\!\!\!\!\!\!\!\!\!\!\!\!\!\!\!\!\!\!\!\!\!\!\!\!\!\!\!\!\!\!\!\!\!\!\!\!\!\!\!\!\!\!\!\!\!\!\!\!\!\!\!
&
{\rm(i)\ }\tilde \xX_1^{2\ell_0^2} \xX_2^{-2\ell_0^2}-3
\stackrel{{}^{\sc\rm\eqref{tX1-is=}}}{=}\xX_2^{-2\ell_0(1+\d_0)}-3
\stackrel{{}^{\sc\rm\eqref{uuuauyyY}\,(iii)}}{=}
\Big(\frac65\Big)^{\ell_0(1+\d_0)}\Big(1+O(\d_0)^{99}\Big)-3>3,
\!\!\!\!\!\!\!\!\!\!\!\!\!\!\!\!\!\!\!\!\!\!\!\!\!\!\!
\nonumber\\
&\!\!\!\!\!\!\!\!\!\!\!\!\!\!\!\!\!\!\!\!\!\!\!\!\!\!\!\!\!\!\!\!\!\!\!\!\!\!\!\!\!\!\!\!\!\!\!\!\!\!\!
&
{\rm(ii)\ }
2\ \ \,\stackrel{{}^{\sc\rm\eqref{LetNSoOP----1}\,(ii)}}{\ge}
\aA_2^{-1}\tilde\xX_1^{\ell_0^2(1+\d_0^2)}\xX_2^{\ell_0^6(1-\d_0^4)}(\tilde\xX_1^{2\ell_0^2}\xX_2^{-2\ell_0^2}-3)
\stackrel{{}^{\sc\rm\eqref{cont-B2},\eqref{tX1-is=},\,\eqref{ii-LetNSoOP----1}\,(i)}}{>}
3\aA_1^{-1}\xX_2^{\ell_0^6(1-\d_0^5+O(\d_0)^6)}
\!\!\!\!\!\!\!\!\!\!\!\!\!\!\!\!\!\!\!\!\!\!\!\!\!\!\!\!\!\!\!\!\!\!\!
\nonumber\\
&
\!\!\!\!\!\!\!\!\!\!\!\!\!\!\!\!\!\!\!\!\!\!\!\!\!\!\!\!\!\!\!\!\!\!\!\!\!\!\!\!\!\!\!\!\!\!\!\!\!\!\!
&
\phantom{{\rm(ii)\ }}
\ \ \ \stackrel{{}^{\sc\rm\eqref{ImMpP}\,(i)}}{\ge}
3\aA_1^{-\d_0^5+O(\d_0)^6}
\stackrel{{}^{\sc\rm\eqref{cont-B2}}}{\ge}
3.
\end{eqnarray}
Equ.~\eqref{ii-LetNSoOP----1}\,(ii) is a contradiction, which proves the first inequality of
\eqref{ImMpP}\,(iv).
}\NOUSE{
 As in \eqref{x2-00ro4m4m4-}\,(i), we also have $1+\d_0\ge\d_0\yY^{2\ell_0^4}-\frac{\aA_1\tilde\xX_1^{\ell_0^4}}{\xX_2^{\ell_0^6 + 4\ell_0^3}}$, which implies that
 $\yY<\ell_0^{\d_0^4}$. Thus
 \equa{yy-is-between}{\dis
 1+O(\d_0)^4=\yY_0<\yY<\ell_0^{\d_0^4}=1+O(\d_0)^4.}
Up to $O(\d)^1$
(one can observe that the following evaluation does not depend on $\d$ as $0<\d\ll\d_0$), we also have
 as $\gamma_1\big(1-\frac{\d_0}{3}\big)=3-\frac{29\d_0^3}{54}+O(\d_0)^4<3$ by \eqref{x2-00ro4m4m4-}
[recalling from \eqref{MSmde33333} that $\d\ll\d_0$]. Thus up to $O(\d)^1$, we have
\begin{eqnarray}
\label{aA1-is-smaller}
\!\!\!\!\!\!\!\!\!\!\!\!\!\!\!\!\!\!\!\!\!\!\!\!\!\!\!\!\!\!\!\!\!\!\!\!&&
{\rm(i)\ }\aA_1<1-\frac{\d_0}{3},
\nonumber\\
\!\!\!\!\!\!\!\!\!\!\!\!\!\!\!\!\!\!\!\!\!\!\!\!\!\!\!\!\!\!\!\!\!\!\!\!&&
  {\rm(ii)\ }
|\ell_0^2+1-\ell_0^2\widetilde X_1|\ge\ell_0^2\tilde\xX_1-(\ell_0^2+1)
\stackrel{{}^{\sc\rm\eqref{tX1-is=}}}{=}\ell_0^2\aA_1^{-1+\d_0}-(\ell_0^2+1)
\stackrel{{}^{\sc\rm\eqref{aA1-is-smaller}\,(i)}}{>}\frac{\ell_0}{3}+O(\d_0)^0,
\nonumber\\
\!\!\!\!\!\!\!\!\!\!\!\!\!\!\!\!\!\!\!\!\!\!\!\!\!\!\!\!\!\!\!\!\!\!\!\!&&
{\rm(iii)\ }
1\ \ \ \ \ \ \stackrel{{}^{\sc\rm\eqref{LetNSoOP----1}\,(ii)}}{=}
\ \ \ \aA_2^{-1}
\tilde \xX_1^{\ell_0^3}\xX_2^{2\ell_0^3} |\ell_0^2+1 -\ell_0^2\widetilde X_1|
\nonumber\\
\!\!\!\!\!\!\!\!\!\!\!\!\!\!\!\!\!\!\!\!\!\!\!\!\!\!\!\!\!\!\!\!\!\!\!\!&&
\ \ \ \ \ \
\stackrel{{}^{\sc\rm\eqref{cont-B2},\,\eqref{tX1-is=},\,\eqref{aA1-is-smaller}\,(ii)}}{>}
\Big(\frac{\ell_0}{3}+O(\d_0)^0\Big)
\aA_1^{-\ell_0^3(1-\d_0)+\ell_0^3(-1+\d_0)}=\frac{\ell_0}{3}+O(\d_0)^0.
\!\!\!\!\!\!
\end{eqnarray}
Equ.~\eqref{aA1-is-smaller}\,(iii)
}\NOUSE{
Then up to $O(\d)^1$, we have [noting from \eqref{aA1-is-smaller}\,(iii) that $\yY^{1-\d_0}\ge2\big(1+O(\d_0)^1\big)\ssc\,$],
\begin{eqnarray}
\label{x2-00ro4m4m4++11}
\!\!\!\!\!\!\!\!\!\!\!\!\!\!\!\!\!\!\!\!\!\!\!\!\!\!\!\!\!\!\!\!\!&&
{\rm(i)\ }
\Big |3 - \frac{2 X_2^{\ell_0}}{A_1^{4\ell_0-3}\widetilde X_1}\Big|
\ge\frac{2 \xX_2^{\ell_0}}{\aA_1^{4\ell_0-3}\tilde \xX_1}-3
\stackrel{{}^{\sc\rm\eqref{tX1-is=}}}{=}
2\aA_1^{-4\ell_0+3+\d_0}\xX_2^{\ell_0}-3
\nonumber\\
\!\!\!\!\!\!\!\!\!\!\!\!\!\!\!\!\!\!\!\!\!\!\!\!\!\!\!\!\!\!\!\!\!\!\!\!&&
\phantom{{\rm(i)\ }\Big |3 - \frac{2 X_2^{\ell_0}}{A_1^2\widetilde X_1}\Big|}
=2(\aA_1^{-3\d_0}\xX_2^{\ell_0})^{1-\d_0}\aA_1^{-4\ell_0+3+\d_0+3\d_0(1-\d_0)}-3
=2\yY^{1-\d_0}\aA_1^{-4\ell_0+O(\d_0)^0}-3
\nonumber\\
\!\!\!\!\!\!\!\!\!\!\!\!\!\!\!\!\!\!\!\!\!\!\!\!\!\!\!\!\!\!\!\!\!\!\!\!&&
\ \ \ \ \ \ \ \ \ \ \ \ \ \ \ \
\stackrel{{}^{\sc\rm\eqref{A1-A2-cond}\,(a)}}{\ge}
4\Big(1+O(\d_0)^1\Big)\aA_1^{-4\ell_0+O(\d_0)^0}-3\aA_1^{-4\ell_0+O(\d_0)^0}
=\Big(1+O(\d_0)^1\Big)\aA_1^{-4\ell_0+O(\d_0)^0},
\!\!\!\!\!\!\!\!\!\!\!\!\!\!\!\!\!\!\!\!\!\!\!\!\!\!\!\!\!\!\!\!\!\!\!\!
\nonumber\\
\!\!\!\!\!\!\!\!\!\!\!\!\!\!\!\!\!\!\!\!\!\!\!\!\!\!\!\!\!\!\!\!\!\!\!\!&&
{\rm(ii)\ }
1\stackrel{{}^{\sc\rm\eqref{LetNSoOP----1}\,(ii)}}{=}
\aA_2^{-1}
\tilde \xX_1 \xX_2^{\ell_0^3}
\Big|3 - \frac{2 X_2^{\ell_0}}{A_1^2\widetilde X_1}\Big|
\stackrel{{}^{\sc\rm\eqref{cont-B2},\,\eqref{tX1-is=}}}{=}
\aA_1^{-1-\d_0}\xX_2^{\ell_0^3+1}
\Big|3 - \frac{2 X_2^{\ell_0}}{A_1^2\widetilde X_1}\Big|
\nonumber\\
\!\!\!\!\!\!\!\!\!\!\!\!\!\!\!\!\!\!\!\!\!\!\!\!\!\!\!\!\!\!&&\phantom{{\rm(i)\ }3\,}
\stackrel{{}^{\sc\rm\eqref{x2-00ro4m4m4++11}\,(i)}}{\ge}
\Big(1+O(\d_0)^1\Big)\aA_1^{-4\ell0+O(\d_0)^0}\xX_2^{\ell_0^3+1}
\stackrel{{}^{\sc\rm\eqref{A1-A2-cond}\,(a),\,\eqref{ImMpP}\,(i),\,\eqref{x2-00ro4m4m4-}}}{>}\yY
.
\end{eqnarray}
Equ.~\eqref{x2-00ro4m4m4++11}\,(ii)
}%
}


Next assume the first inequality of \eqref{ImMpP}\,(iv) does not hold. Then
 by \eqref{C+LetNSoOP}\,(iii)
,
 \equa{111---tX1-is=111}{\dis
\tilde\xX_1=
(1-\d)\xX_2\aA_1^{2+\d_0}.
}
Using the second expression of $A_2$ in  \eqref{We0o0o0o++1} and the fact in \eqref{cont-B2}
that $\aA_2=\aA_1$ [up to $O(\d)^2\ssc\,$],
we can compute the following, up to $O(\d)^2$
,
\begin{eqnarray}
\label{x2-00ro4m4m4-111}
\!\!\!\!\!\!\!\!\!\!\!\!\!\!\!\!\!\!\!\!\!\!\!\!\!\!\!\!\!\!
\!\!\!\!\!\!\!\!\!\!\!\!\!\!\!\!\!\!\!\!\!\!\!\!\!\!\!\!\!\!&&
2 \stackrel{{}^{\sc\rm\eqref{We0o0o0o++1}}}{\le}
\frac{\tilde \xX_1}{\aA_1^2}+
 \frac{\aA_2}{\xX_2^{\ell_0}}
\stackrel{{}^{\sc\rm
\eqref{cont-B2},\,\eqref{111---tX1-is=111}}}{=}
(1-\d)\xX_2\aA_1^{\d_0}
+\frac{\aA_1}{\xX_2^{\ell_0}}\mbox{ \ [up to $O(\d)^2\ssc\,$]}.
\NOUSE{%
\!\!\!\!\!\!\!\!\!\!\!\!\!\!\!\!\!\!\!\!\!\!\!\!\!\!\!\!\!\!\!\!\!\!\!\!\!\!\!\!\!\!
\nonumber\\
\!\!\!\!\!\!\!\!\!\!\!\!\!\!\!\!\!\!\!\!\!\!\!\!\!\!\!\!\!\!
\!\!\!\!\!\!\!\!\!\!\!\!\!\!\!\!\!\!\!\!\!\!\!\!\!\!\!\!\!\!&&
\phantom{1+\d_0^5\ \ \, {\rm(i)\ }}
\stackrel{{}^{\sc\rm\eqref{111---tX1-is=111}}}{=}
\frac{\aA_1^{5-4\d_0^2}\xX_2^{5\ell_0-1}}{(1+\d)^4}
+\frac{\d_0^5\aA_1^{-\ell_0^6+O(\d_0)^{-1}}}{(1+\d)^{\ell_0^8-5\ell_0^5}\xX_2^{\ell_0^{20}+O(\d_0)^{-19}}}
\!\!\!\!\!\!\!\!\!\!\!\!\!\!\!\!\!\!\!\!\!\!\!\!\!\!\!\!\!\!\!\!\!\!\!\!\!\!\!\!\!\!
\nonumber\\
\!\!\!\!\!\!\!\!\!\!\!\!\!\!\!\!\!\!\!\!\!\!\!\!\!\!\!\!\!\!
\!\!\!\!\!\!\!\!\!\!\!\!\!\!\!\!\!\!\!\!\!\!\!\!\!\!\!\!\!\!&&
\phantom{1+\d_0^5\ \ \ \ \ {\rm(i)\ }}
<\frac1{1+\d}\mbox{\Large$\Big($}\aA_1^{5-4\d_0^2}\xX_2^{5\ell_0-1}+
\frac{\d_0^5(\aA_1^{5-4\d_0^2}\xX_2^{5\ell_0-1})^{-\frac{4\ell_0^6}{5}}
\aA_1^{-\frac{\ell_0^6}{5}+O(\d_0)^{-5}}}{\xX_2^{\ell_0^{20}+O(\d_0)^{-19}}}\mbox{\Large$\Big)$}
\!\!\!\!\!\!\!\!\!\!\!\!\!\!\!\!\!\!\!\!\!\!\!\!\!\!\!\!\!\!\!\!\!\!\!\!\!\!\!\!\!\!
\\\nonumber
\!\!\!\!\!\!\!\!\!\!\!\!\!\!\!\!\!\!\!\!\!\!\!\!\!\!\!\!\!\!
\!\!\!\!\!\!\!\!\!\!\!\!\!\!\!\!\!\!\!\!\!\!\!\!\!\!\!\!\!\!&&
\phantom{1\!+\!\d_0^5\,}
\stackrel{{}^{\sc\rm\eqref{A1-A2-cond}\,(a),\,\eqref{ImMpP}\,(i)}}{\le}
\gamma_2(\yY):=\frac1{1+\d}\Big(\yY+\d_0^5\yY^{-\frac{4\ell_0^6}{5}}\Big),\ \ \mbox{ where}
\!\!\!\!\!\!\!\!\!\!\!\!\!\!\!\!\!\!\!\!\!\!\!\!\!\!\!\!\!\!\!\!\!\!\!\!\!\!\!\!\!\!
\\\nonumber
\!\!\!\!\!\!\!\!\!\!\!\!\!\!\!\!\!\!\!\!\!\!\!\!\!\!\!\!\!\!
\!\!\!\!\!\!\!\!\!\!\!\!\!\!\!\!\!\!\!\!\!\!\!\!\!\!\!\!\!\!&&
{\rm(ii)\ }
\yY:=\aA_1^{5-4\d_0^2}\xX_2^{5\ell_0-1}
\stackrel{{}^{\sc\rm\eqref{A1-A2-cond}\,(a),\,\eqref{ImMpP}\,(i)}}{>}1
.
}
\!\!\!\!\!\!\!\!\!\!\!\!\!\!\!\!\!\!\!\!\!\!\!\!\!\!\!\!\!\!\!\!\!\!\!\!\!\!\!\!\!\!
\end{eqnarray}
Observing that
\eqref{x2-00ro4m4m4-111} is exactly the same as the first line in
\eqref{x2-00ro4m4m4-}\,(i),
we can thus use the same arguments after \eqref{x2-00ro4m4m4-} to show
that we have
\eqref{uuuauyyY}\,(i)--(v), and further, using \eqref{111---tX1-is=111},
exactly as in \eqref{uuuauyyY}\,(vi), we have $\frac1{\tilde\xX_1}=\aA_1^{-2-\d_0}\xX_2^{-1}\ge2+4\ln(2)\d_0+O(\d_0)^2$. Then
as in \eqref{Futtth}, we can use the definition of $A_1$ in 
\eqref{LetNSoOP----1-re-give}\,(i) to obtain, up to $O(\d)^1$,
\equan{msmsmsnnn}{\dis
\frac12+O(\d_0)^1\stackrel{{}^{\sc\rm\eqref{uuuauyyY}\,(iii),\,(iv)}}{\ge}
\aA_1
\stackrel{{}^{\sc\rm\eqref{LetNSoOP----1-re-give}\,(i)}}{\ge}
\xX_2^{\ell_0} \Big(\frac1{\tilde \xX_1}-2\Big)\ge\Big(2+O(\d_0)^{100}\Big)^{\ell_0}\Big(4\ln(2)\d_0+O(\d_0)^2\Big)\gg2,
}
which is again a contradiction. This proves the last inequality of \eqref{ImMpP}\,(iv).
\NOUSE{
$\gamma_2(\aA_1)$ is decreasing when $0<\aA_1<\aA_0$ and increasing when $\aA_1>\aA_0$,
where $\aA_0=\big(\frac{\d_0}{1+\d}\big)^{\frac{1}{1+\d_0}}$.
Since $\aA_1\le1$ and $\gamma_2(1)=2-\d+O(\d)^2<2$, we must have (i) below, and so we have
}\NOUSE{ and
\gamma_2
[noting that  $(1+\frac{\d_0^5}{2})^{-\frac{4\ell_0^6}{5}}\ll\d_0^{100}$ and that  $\d\ll\d_0^{100}$],
we must have
(i) below
}%
\NOUSE{%
[again we can conduct the following computations up to $O(\d)^1\ssc\,$%
],
\begin{eqnarray}
\label{uuuauyyY++1}
\!\!\!\!\!\!\!\!\!\!\!\!\!\!\!\!\!\!\!\!\!\!&&
{\rm(i)\ }
\aA_1^{\d_0}<\frac{4}{7},\ \ \ \ \ \
\!\!\!\!\!\!\!\!\!\!\!\!\!
\\\nonumber
\!\!\!\!\!\!\!\!\!\!\!\!\!\!\!\!\!\!\!\!\!\!\!\!\!\!\!\!\!\!&&
{\rm(ii)\ }
\frac{\tilde\xX_1}{\xX_2^2}
\stackrel{{}^{\sc\rm\eqref{111---tX1-is=111}}}{=}
\aA_1^{-2+\d_0}
\frac{\aA_1^{-\ell_0^8}}{\tilde \xX_1^{\ell_0^8 - 5\ell_0^5} \xX_2^{\ell_0^{20}}}
\stackrel{{}^{\sc\rm\mbox{deducted as in }\eqref{x2-00ro4m4m4-111}\,(i)}}{\le}\yY^{-\frac{4\ell_0^6}{5}}
\stackrel{{}^{\sc\rm\eqref{uuuauyyY++1}\,(i)}}{\ll}\d_0^{100},
\!\!\!\!\!\!\!\!\!\!\!\!\!
\nonumber\\
\nonumber
\!\!\!\!\!\!\!\!\!\!\!\!\!\!\!\!\!\!\!\!\!\!\!\!\!\!\!\!\!\!&&
{\rm(iii)\ }
\frac{A_1 X_2^{\ell_0-1}}{\widetilde X_1^4}
\stackrel{{}^{\sc\rm\eqref{uuuauyyY++1}\,(ii)}}{=}
1+\d_0^5+O(\d_0)^{105}=1+\d_0^5+O(\d_0)^{10}.
\!\!\!\!\!\!\!\!\!\!\!\!\!
\end{eqnarray}
Further, using \eqref{uuuauyyY++1}\,(iii) and the definition of $A_1$ in
\eqref{LetNSoOP----1}\,(i), which implies that
$\frac{A_1}{X_2\widetilde X_1^3}=5-\frac{4\widetilde X_1}{X_2^{\ell_0}}$, we obtain, up to $O(\d)^1$,
\begin{eqnarray}
\label{uuuauyyY++1KK}
\!\!\!\!\!\!\!\!\!\!\!\!\!\!\!\!\!\!\!\!\!\!\!\!\!\!\!\!\!\!&&
{\rm(i)\ }
\frac{\widetilde X_1 }{X_2^{\ell_0}}
=\frac{\widetilde X_1^4}{A_1X_2^{\ell_0-1}}\cdot\frac{A_1}{X_2\widetilde X_1^3}
 \stackrel{{}^{\sc\rm\eqref{uuuauyyY++1}\,(iii),\,\eqref{LetNSoOP----1}\,(i)}}{=}
\Big(\frac1{1+\d_0^5}+O(\d_0)^{10}\Big)\Big(5-\frac{4\widetilde X_1}{X_2^{\ell_0}}\Big),
\!\!\!\!\!\!\!\!\!\!\!\!\!
\\
\!\!\!\!\!\!\!\!\!\!\!\!\!\!\!\!\!\!\!\!\!\!\!\!\!\!\!\!\!\!&&
{\rm(ii)\ }
\frac{\widetilde X_1 }{X_2^{\ell_0}}
\Big(1+\frac4{1+\d_0^5}+O(\d_0)^{10}\Big)
\stackrel{{}^{\sc\rm\eqref{uuuauyyY++1KK}\,(i)}}{=}
\frac{5}{1+\d_0^5}+O(\d_0)^{10},
\!\!\!\!\!\!\!\!\!\!\!\!\!
\nonumber\\
\!\!\!\!\!\!\!\!\!\!\!\!\!\!\!\!\!\!\!\!\!\!\!\!\!\!\!\!\!\!&&
{\rm(iii)\ }
\frac{\widetilde X_1 }{X_2^{\ell_0}}
\stackrel{{}^{\sc\rm\eqref{uuuauyyY++1KK}\,(iii)}}{=}
\frac{5}{5+\d_0^5}+O(\d_0)^{10}
=1-\frac{\d_0^5}{5}+O(\d_0)^{10},
\!\!\!\!\!\!\!\!\!\!\!\!\!
\nonumber\\
\!\!\!\!\!\!\!\!\!\!\!\!\!\!\!\!\!\!\!\!\!\!\!\!\!\!\!\!\!\!&&
{\rm(iv)\ }
5 - \frac{4\widetilde X_1 }{X_2^{\ell_0}}
\stackrel{{}^{\sc\rm\eqref{uuuauyyY++1KK}\,(iii)}}{=}
5 - 4\Big(1-\frac{\d_0^5}{5}+O(\d_0)^{10}\Big)
=
1+\frac{4\d_0^5}{5}+O(\d_0)^{10},
\!\!\!\!\!\!\!\!\!\!\!\!\!
\nonumber\\
\!\!\!\!\!\!\!\!\!\!\!\!\!\!\!\!\!\!\!\!\!\!\!\!\!\!\!\!\!\!&&
{\rm(v)\ }
\tilde\xX_1^3=\Big(\frac{\tilde \xX_1^4}{\aA_1 \xX_2^{\ell_0-1}}\Big)^{\frac34}
\aA_1^{\frac34}\xX_2^{\frac34(\ell_0-1)}
\stackrel{{}^{\sc\rm\eqref{uuuauyyY++1}\,(iii)}}{=}\Big(1-\frac{3\d_0^5}{4}+O(\d_0)^{10}\Big)
\aA_1^{\frac34}\xX_2^{\frac34(\ell_0-1)}.
\!\!\!\!\!\!\!\!\!\!\!\!\!
\nonumber
\end{eqnarray}
Finally, using
\eqref{uuuauyyY++1KK}\,(iv),\,(v) and the definition of $A_1$ in
\eqref{LetNSoOP----1}\,(i), up to $O(\d)^1$, we have
\begin{eqnarray}
\label{uuuauyyY++1KK-000}
\!\!\!\!\!\!\!\!\!\!\!\!\!\!\!\!\!\!\!\!\!\!\!\!\!\!\!\!\!\!&&
{\rm(i)\ }
\aA_1
\stackrel{{}^{\sc\rm\eqref{LetNSoOP----1}\,(i)}}{=}
\xX_2\tilde\xX_1^3
\Big|5 - \frac{4\widetilde X_1 }{X_2^{\ell_0}}
\Big|
\\\nonumber
\!\!\!\!\!\!\!\!\!\!\!\!\!\!\!\!\!\!\!\!\!\!\!\!\!\!\!\!\!\!&&
\phantom{{\rm(i)\ }}
\stackrel{{}^{\sc\rm\eqref{uuuauyyY++1KK}\,(iv),\,(v)}}{=}
\xX_2\aA_1^{\frac34}\xX_2^{\frac34(\ell_0-1)}\Big(1-\frac{3\d_0^5}{4}+O(\d_0)^{10}\Big)
\Big(1+\frac{4\d_0^5}{5}+O(\d_0)^{10}\Big)
\!\!\!\!\!\!\!\!\!\!\!\!\!
\\
\!\!\!\!\!\!\!\!\!\!\!\!\!\!\!\!\!\!\!\!\!\!\!\!\!\!\!\!\!\!&&
\phantom{{\rm(i)\ }\ \ \ \ \ \ \ }
=\aA_1^{\frac34}\xX_2^{\frac{3\ell_0}4+\frac14}\Big(1+\frac{\d_0^5}{20}+O(\d_0)^{10}\Big),
\!\!\!\!\!\!\!\!\!\!\!\!\!
\nonumber\\\nonumber
\!\!\!\!\!\!\!\!\!\!\!\!\!\!\!\!\!\!\!\!\!\!\!\!\!\!\!\!\!\!&&
{\rm(ii)\ }
\aA_1^{\d_0^4+O(\d_0)^5}
\stackrel{{}^{\sc\rm\eqref{We0o0o0o}\,(i)}}{\ge}
\xX_2
\stackrel{{}^{\sc\rm\eqref{uuuauyyY++1KK-000}\,(i)}}{=}
\Big(\aA_1^{\frac14}\Big(1-\frac{\d_0^5}{20}+O(\d_0)^{10}\Big)\Big)^{\frac1{\frac{3\ell_0}{4}+\frac14}}
\nonumber\\\nonumber
\!\!\!\!\!\!\!\!\!\!\!\!\!\!\!\!\!\!\!\!\!\!\!\!\!\!\!\!\!\!&&
\phantom{{\rm(ii)\ }\ \ \ \ \ \ \ \ \ \ \ \ \ \ \ \ \ }
=\ \ \aA_1^{\frac{\d_0}{3}+O(\d_0)^2}\Big(1+\frac{4\d_0^6}{15}+O(\d_0)^7\Big).
\end{eqnarray}
Equ.~\eqref{uuuauyyY++1KK-000}\,(ii) is
a contradiction with the fact in
\eqref{A1-A2-cond}\,(a) that $\aA_1\ge1$.
This proves the last inequality of
\eqref{ImMpP}\,(iv).
%
%
%
\NOUSE{%
On the other hand, we have, where $\ln(\cdot)$ is the natural logarithmic function,
\begin{eqnarray}
\label{uuuauyyY++1+XXX}
\!\!\!\!\!\!\!\!\!\!\!\!\!\!\!\!\!\!\!\!\!\!&&
{\rm(i)\ }1+O(\d_0)^{101}
\stackrel{{}^{\sc\rm\eqref{uuuauyyY++1}\,(vi)}}{=}
\tilde\xX_1^{2\ell_0}\xX_2^{2\ell_0^2}
\stackrel{{}^{\sc\rm\eqref{111---tX1-is=111}}}{=}
\frac{\aA_1^{-2\ell_0(\ell_0^2-\ell_0+1)}}{\xX_2^{6\ell_0^3-2\ell_0^2}}
\nonumber\\\!\!\!\!\!\!\!\!\!\!\!\!\!\!\!\!\!\!\!\!\!\!&&
{\rm(iii)\ }
\aA_1^{\ell_0^2(1-\d_0+\d_0^2)}=\xX_2^{-(3\ell_0^2-\ell_0)}\Big(1+\d_0+O(\d_0)^{100}\Big)^{-\frac{\d_0}{2}}=
\xX_2^{-(3\ell_0^2-\ell_0)}\Big(1-\frac{\d_0^2}{2}+O(\d_0)^3\Big).
\!\!\!\!\!\!\!\!\!\!\!\!\!\!\!\!\!\!\!\!\!\!
\nonumber\\
\!\!\!\!\!\!\!\!\!\!\!\!\!\!\!\!\!\!\!\!\!\!&&
{\rm(iv)\ }
\frac12-\frac{\ln(2)}{2}\d_0+O(\d_0)^2
\stackrel{{}^{\sc\rm\eqref{uuuauyyY++1}\,(iv)}}{=}
\tilde\yY^{-(1+\d_0-\d_0^2)}
\stackrel{{}^{\sc\rm\eqref{x2-00ro4m4m4-111}\,(ii)}}{=}
\aA_1^{2\ell_0^2(1+\d_0-\d_0^2)}
\!\!\!\!\!\!\!\!\!\!\!\!\!\!\!\!\!\!\!\!\!\!
\\\nonumber
\!\!\!\!\!\!\!\!\!\!\!\!\!\!\!\!\!\!\!\!\!\!&&
\phantom{{\rm(iv)\ }\ \ \ \ \ \ \ \ \ \ \ \ \ \ \ \ \ \ \ \ \ \ \ \ \ \ \ }
\stackrel{{}^{\sc\rm\eqref{uuuauyyY++1+XXX}\,(i),\,(iii)}}{=}
\Big(\frac12+\frac{\d_0^2}{4}+O(\d_0)^3\Big)\Big(1-\frac{\d_0^2}{2}+O(\d_0)^3\Big)=\frac12+O(\d_0)^2
.\end{eqnarray}
On the other hand, we have, where $\ln(\cdot)$ is the natural logarithmic function,
\begin{eqnarray}
\label{uuuauyyY++1+XXX-YYY}
\!\!\!\!\!\!\!\!\!\!\!\!\!\!\!\!\!\!\!\!\!\!&&
{\rm(i)\ }
\aA_1^{\ell_0^2(1+\d_0-\d_0^2)}
\stackrel{{}^{\sc\rm\eqref{We0o0o0o++1}\,(i),\,\eqref{cont-B2}}}{=}
\frac{\xX_2^{\ell_0^2}}{\tilde\xX_1^{2\ell_0-1}}
 \Big|\frac{\ell_0+1}{2} - \frac{\ell_0-1}{2}\widetilde X_1^{2\ell_0} X_2^{2\ell_0^2}\Big|
\!\!\!\!\!\!\!\!\!\!
\nonumber\\\nonumber
\!\!\!\!\!\!\!\!\!\!\!\!\!\!\!\!\!\!\!\!\!\!&&
\phantom{{\rm(vii)\ }
}
\stackrel{{}^{\sc\rm\eqref{uuuauyyY++1}\,(vi)}}{=}
\frac{\xX_2^{\ell_0^2}}{\Big(\xX_2^{-2\ell_0^2}(1{\sc\!}+{\sc\!}\d_0)
\Big(1{\sc\!}+{\sc\!}O(\d_0)^{101}\Big)\Big)^{\frac{2\ell_0-1}{2\ell_0}}}
 \Big(\frac{\ell_0{\sc\!}+{\sc\!}1}{2} {\sc\!}-{\sc\!}
  \frac{\ell_0{\sc\!}-{\sc\!}1}{2}(1{\sc\!}+{\sc\!}\d_0)\Big(1{\sc\!}+{\sc\!}O(\d_0)^{101}\Big)\Big)
\!\!\!\!\!\!\!\!\!\!
\\
\!\!\!\!\!\!\!\!\!\!\!\!\!\!\!\!\!\!\!\!\!\!&&
\phantom{{\rm(vii)\ }\ \ \ \
}=\ \ \xX_2^{3\ell_0^2-\ell_0}\Big(\frac12+\frac{\d_0^2}{4} -\frac{ \d_0^3}{8} + \frac{7 \d_0^4}{48}+O(\d_0)^5\Big)
.
\nonumber\\
\!\!\!\!\!\!\!\!\!\!\!\!\!\!\!\!\!\!\!\!\!\!&&
{\rm(iii)\ }
\frac{\aA_1^{-2\ell_0(\ell_0^2-\ell_0+1)}}{\xX_2^{6\ell_0^3-2\ell_0^2}}
\stackrel{{}^{\sc\rm\eqref{111---tX1-is=111}}}{=}
\tilde\xX_1^{2\ell_0}\xX_2^{2\ell_0^2}\stackrel{{}^{\sc\rm\eqref{uuuauyyY++1}\,(vi)}}{=}
1+\d_0+O(\d_0)^{100},\ \ \implies\
\nonumber\\\!\!\!\!\!\!\!\!\!\!\!\!\!\!\!\!\!\!\!\!\!\!&&
{\rm(iii)\ }
\aA_1^{\ell_0^2(1-\d_0+\d_0^2)}=\xX_2^{-(3\ell_0^2-\ell_0)}\Big(1+\d_0+O(\d_0)^{100}\Big)^{-\frac{\d_0}{2}}=
\xX_2^{-(3\ell_0^2-\ell_0)}\Big(1-\frac{\d_0^2}{2}+O(\d_0)^3\Big).
\!\!\!\!\!\!\!\!\!\!\!\!\!\!\!\!\!\!\!\!\!\!
\nonumber\\
\!\!\!\!\!\!\!\!\!\!\!\!\!\!\!\!\!\!\!\!\!\!&&
{\rm(iv)\ }
\frac12-\frac{\ln(2)}{2}\d_0+O(\d_0)^2
\stackrel{{}^{\sc\rm\eqref{uuuauyyY++1}\,(iv)}}{=}
\tilde\yY^{-(1+\d_0-\d_0^2)}
\stackrel{{}^{\sc\rm\eqref{x2-00ro4m4m4-111}\,(ii)}}{=}
\aA_1^{2\ell_0^2(1+\d_0-\d_0^2)}
\!\!\!\!\!\!\!\!\!\!\!\!\!\!\!\!\!\!\!\!\!\!
\\\nonumber
\!\!\!\!\!\!\!\!\!\!\!\!\!\!\!\!\!\!\!\!\!\!&&
\phantom{{\rm(iv)\ }\ \ \ \ \ \ \ \ \ \ \ \ \ \ \ \ \ \ \ \ \ \ \ \ \ \ \ }
\stackrel{{}^{\sc\rm\eqref{uuuauyyY++1+XXX}\,(i),\,(iii)}}{=}
\Big(\frac12+\frac{\d_0^2}{4}+O(\d_0)^3\Big)\Big(1-\frac{\d_0^2}{2}+O(\d_0)^3\Big)=\frac12+O(\d_0)^2
.\end{eqnarray}
\begin{rema}\rm\label{Reasso}
At this point it may be worth mentioning that the reason we can obtain the above contradiction is simply because  we
have \eqref{uuuauyyY}\,(iv), which determines $\widetilde X_1^{\ell_0}$ up to $O(\d_0)^{100}$.
\end{rema}
Equ.~\eqref{uuuauyyY++1+XXX}\,(iv) is a contradiction (cf.~Remark \ref{Reasso}), which proves the last inequality of
\eqref{ImMpP}\,(iv).
}%
\NOUSE{%
Observing that $\frac{d\gamma_3}{d\aA_1}|_{\aA_1=\tilde\aA_0}=O(\d)^1$ for
$\tilde\aA_0=\big(1+\frac{\d_0^2}{2}\big)^{\frac{\d_0^2}{4+\d_0^2}}=1+\frac{\d_0^4}{8}+O(\d_0)^5$. Exactly Similar to
the arguments after
\eqref{x2-00ro4m4m4-}, by the facts that $\gamma_3(1)=2-\d+O(\d)^2<2$, $\gamma_3(\aA_1)\ge2$, we obtain that $\aA_1\ge\tilde\aA_0$.
Thus the second term of $\gamma_2$ is $\le1-\frac{\d_0^2}{4}+O(\d_0)^3$, which implies (i) below
(recalling that $\d\ll\d_0$),
\begin{eqnarray}
\label{ooo-x2-00ro4m4m4-111}\!\!\!\!\!\!\!\!\!\!\!\!\!\!\!\!\!\!\!\!\!\!\!\!\!\!\!\!\!\!
\!\!\!\!\!\!\!\!\!\!\!\!\!\!\!\!\!\!\!\!\!\!\!\!\!\!\!\!\!\!&&
{\rm(i)\ }\yY:=\frac{\aA_1^{2\ell_0^2}}{\xX_2^{3\ell_0^5}}\ge1+\frac{\d_0^2}{4}+O(\d_0)^3,
\nonumber\\
\!\!\!\!\!\!\!\!\!\!\!\!\!\!\!\!\!\!\!\!\!\!\!\!\!\!\!\!\!\!&&
{\rm(ii)\ }
\Big(\frac{\widetilde X_1^{\ell_0^2}}{A_1^{\ell_0}X_2^{\ell_0^5+\ell_0^2 -\ell_0}}\Big)^{\ell_0^3}
\stackrel{{}^{\sc\rm\eqref{111---tX1-is=111}}}{=}
\frac{(1+\d)^{2\ell_0^5}\aA_1^{\ell_0^5(2-\d_0)}}{\xX_2^{3\ell_0^8-\ell_0^4}}
\end{eqnarray}
}
\NOUSE{
Note that the function $\gamma_1(\aA_1)$ on $\aA_1$ is
strictly increasing when $0<\aA_1\le\aA_0$ and strictly decreasing when $\aA_0\le\aA_1\le1$, where
$\aA_0=\frac{2(1-\d^2)\d_0}{(1-\d^2)^{-1}(1-\d_0)}=2\d_0+O(\d_0)^2$ as $\frac{\d\gamma_1}{d\aA_1}|_{\aA_1=\aA_0}=0$.
Since $\aA_1\le1$ by \eqref{A1-A2-cond}\,(a) and $\gamma_1(1)=3-\d^2+O(\d)^3<3$, in order for
\eqref{x2-00ro4m4m4-} to hold, we must have (recalling from \eqref{MSmde33333} that $\ell_0\d_0=1$),
\equa{aA1-is-smaller}{\dis
\aA_1<\Big(\frac34\Big)^{\ell_0}\ll\d_0^{100},
}
as $\gamma_1\big((\frac34)^{\ell_0}\big)<3$ by \eqref{x2-00ro4m4m4-}.
}%
\NOUSE{%
Exactly similar to \eqref{x2-00ro4m4m4-},\,\eqref{aA1-is-smaller}, we can obtain
\equa{aA1-is-smaller+11}{\dis\!\!\!\!\!\!\!\!\!\!\!\!
{\rm(i)\, }\tilde\yY{\ssc\!}>{\ssc\!}\frac43,\  \ \implies\ \   {\rm(ii)\, }
\tilde\yY^{-\ell_0(3-\d_0)}\ll\d_0^{100},\
 \ {\rm(iii)\, }\tilde\yY\stackrel{{}^{\sc\rm\eqref{x2-00ro4m4m4-111}\,(i)}}{\ge}\frac32{\ssc\!}+
 {\ssc\!}O(\d_0)^{100}{\ssc\!}={\ssc\!}\frac32\Big(1{\ssc\!}+{\ssc\!}O(\d_0)^{100}\Big).\!\!\!\!\!\!
\!\!\!\!\!\!}
Then as in \eqref{x2-00ro4m4m4-111}\,(i), we also have
\begin{eqnarray}
\label{x2-00ro4m4m4-111+abc}
\!\!\!\!\!\!\!\!\!\!\!\!\!\!\!\!\!\!\!\!\!\!\!\!\!\!\!\!\!\!&&
{\rm(i)\ }2(1+\d)^{-1}\tilde\yY=
\frac{2\xX_2^{\ell_0}}{\aA_1^{4\ell_0-3}\tilde\xX_1}
\stackrel{{}^{\sc\rm\eqref{LetNSoOP----1}\,(ii)}}{\le}
3+\aA_2\tilde\xX_1^{-\ell_0}\aA_3^{-\ell_0^3}
=3+O(\d_0)^{100},\ \ \ \stackrel{{}^{\sc\rm\eqref{aA1-is-smaller+11}\,(iii)}}{\implies}
\!\!\!\!\!\!\!\!\!\!\!\!\!\nonumber\\
\!\!\!\!\!\!\!\!\!\!\!\!\!\!\!\!\!\!\!\!\!\!\!\!\!\!\!\!\!\!&&
{\rm(ii)\ }\tilde\yY=\frac32\Big(1+O(\d_0)^{100}\Big).
\!\!\!\!\!\!\!\!\!\!\!\!\!\!\!\!\!\!\!\!
\end{eqnarray}
Now up to $O(\d)^1$, we have
,
\begin{eqnarray}
\label{x2-00ro4m4m4-111}
\!\!\!\!\!\!\!\!\!\!\!\!\!\!\!\!\!\!\!\!\!\!\!\!\!\!\!\!\!\!\!\!\!&&
{\rm(i)\ }
\tilde \xX_1^6 \xX_2^{2\ell_0-6}
\stackrel{{}^{\sc\rm\eqref{tX1-is=111}}}{=}
\aA_1^{-6\ell_0+18}\xX_2^{2\ell_0-6}\stackrel{{}^{\sc\rm\eqref{x2-00ro4m4m4-111}\,(ii)}}{=}\tilde\yY^{2(1-3\d_0)}
\stackrel{{}^{\sc\rm\eqref{x2-00ro4m4m4-111+abc}\,(ii)}}{=}\frac94\Big(1+O(\d_0)^1\Big),
\nonumber\\
\!\!\!\!\!\!\!\!\!\!\!\!\!\!\!\!\!\!\!\!\!\!\!\!\!\!\!\!\!\!\!\!\!\!\!\!&&
{\rm(ii)\ }\Big |\frac43 - \frac13\widetilde X_1^6 X_2^{2\ell_0-6}\Big|
\ge
\frac43- \frac13\tilde \xX_1^6 \xX_2^{2\ell_0-6}
\stackrel{{}^{\sc\rm\eqref{x2-00ro4m4m4-111}\,(i)}}{=}
\frac7{12}+O(\d_0),
\nonumber\\
\!\!\!\!\!\!\!\!\!\!\!\!\!\!\!\!\!\!\!\!\!\!\!\!\!\!\!\!\!\!\!\!\!\!\!\!&&
{\rm(iii)\ }
1\stackrel{{}^{\sc\rm\eqref{We0o0o0o}\,(iv)}}{=}
\aA_1^{-1}\tilde \xX_1 \xX_2^{\ell_0^3}\Big |\frac43 - \frac{\widetilde X_1^6 X_2^{2\ell_0-6}}{3}\Big|
\stackrel{{}^{\sc\rm\eqref{tX1-is=111},\,\eqref{x2-00ro4m4m4-111}\,(ii)}}{\ge}
\aA_1^{-\ell_0+2}\xX_2^{\ell_0^3}
\nonumber\\
\!\!\!\!\!\!\!\!\!\!\!\!\!\!\!\!\!\!\!\!\!\!\!\!\!\!\!\!\!\!\!\!\!\!\!\!&&
\phantom{{\rm(i)\ }\Big |\frac43 - \frac13\widetilde X_1^6 X_2^{2\ell_0-6}\Big|}
\stackrel{{}^{\sc\rm\eqref{A1-A2-cond}\,(a)}}{\ge}.
\end{eqnarray}
Equ.~\eqref{x2-00ro4m4m4-111}\,(v) is a contradiction with \eqref{A1-A2-cond}\,(a), which proves the last inequality of
\eqref{ImMpP}\,(iv).
Equ.~\eqref{ImMpP}\,(v) simply follows from \eqref{ImMpP}\,(iv).
Now we have
\begin{eqnarray}
\label{fafdnenen}
\!\!\!\!\!\!\!\!\!\!\!\!\!\!\!\!\!\!\!\!\!\!\!\!\!\!\!\!\!\!\!\!\!&&
{\rm(i)\ }
\tilde\xX_1\aA_3^{-\ell_0}
\stackrel{{}^{\sc\rm\eqref{ImMpP}\,(iv)}}{\ge}(1-\d^2)\aA_1^{-\d_0}\aA_3^{-1}
=(1-\d^2)(\aA_1\aA_3^{\ell_0})^{-\d_0}
\stackrel{{}^{\sc\rm\eqref{ImMpP}\,(v)}}{\ge}1+O(\d)^2,
\nonumber\\
\!\!\!\!\!\!\!\!\!\!\!\!\!\!\!\!\!\!\!\!\!\!\!\!\!\!\!\!\!\!\!\!\!&&
{\rm(ii)\ }
\aA_1^2\tilde\xX_1\aA_3^{\ell_0}
\stackrel{{}^{\sc\rm\eqref{ImMpP}\,(iv)}}{\le}(1+\d^2)\aA_1^{\d_0}\aA_3
=(1+\d^2)(\aA_1\aA_3^{\ell_0})^{\d_0}
\stackrel{{}^{\sc\rm\eqref{ImMpP}\,(v)}}{\le}1+O(\d)^2.
\end{eqnarray}
This proves
\eqref{ImMpP}\,(vi),\,(vii).
Finally assume \eqref{ImMpP}\,(viii) does not hold. Then by \eqref{C+LetNSoOP}\,(iv), we have
\equa{Now-tildeX-is}{\dis
\tilde\xX_1=(1-\d)\frac{\aA_3^{\ell_0+1}}{\aA_1^{\d_0^2}}.
}
Then
\begin{eqnarray}
\label{a----mswmw}
\!\!\!\!\!\!\!\!\!\!\!\!\!\!\!\!\!\!\!\!\!\!\!\!\!\!\!\!\!\!\!\!\!&&
{\rm(i)\ }
1+O(\d)^2\stackrel{{}^{\sc\rm\eqref{ImMpP}\,(vii)}}{\le}\frac1{\aA_1^2\tilde\xX_1\aA_3^{\ell_0}}
\stackrel{{}^{\sc\rm\eqref{Now-tildeX-is}}}{=}
\frac1{(1-\d)
\aA_1^{2-\d_0^2}\aA_3^{2\ell_0+1}}
=\frac{\aA_1^{\d_0^2}}{(1-\d)\aA_3(\aA_1\aA_3^{\ell_0})^2}
\nonumber\\
\!\!\!\!\!\!\!\!\!\!\!\!\!\!\!\!\!\!\!\!\!\!\!\!\!\!\!\!\!\!\!\!\!&&
\stackrel{{}^{\sc\rm\eqref{ImMpP}\,(v)}}{\le}
\aA_1^{\d_0^2}\aA_3^{-1}+O(\d)^1
\stackrel{{}^{\sc\rm\eqref{A1-A2-cond}\,(a),\,\eqref{ImMpP}\,(i)}}{\le}
1+O(\d)^1.
%
\!\!\!\!\!\!\!\!\!\!\!\!
\end{eqnarray}
Equ.~\eqref{a----mswmw}\,(iii) with \eqref{A1-A2-cond}\,(a) shows that $\aA_1=1+O(\d)^1$. Then \eqref{ImMpP}\,(i),\,(v) gives that
$\aA_3=1+O(\d)^1$, and so $\tilde\xX_1=1+O(\d)^1$ by \eqref{ImMpP}\,(iv). This proves, where $T_2$ is defined in \eqref{denote-t2},
\equa{a-is-O(d)-ele}{\dis
a=1+O(\d)^1\mbox{ \ for \ }a\in T_2.}
Then we have
\begin{eqnarray}
\label{a----mswmw+2222}
\!\!\!\!\!\!\!\!\!\!\!\!\!\!\!\!\!\!\!\!\!\!\!\!\!\!\!\!\!\!\!\!\!&&
3
\stackrel{{}^{\sc\rm\eqref{We0o0o0o}\,(iv)}}{\le}
2\tilde\xX_1\aA_3^{-\ell_0}+\aA_1\tilde\xX_1^{-1}\aA_3^{-\ell_0^5}
\stackrel{{}^{\sc\rm\eqref{We0o0o0o}\,(iv)}}{\le}
2(1-\d)\aA_1%
\\\nonumber
\!\!\!\!\!\!\!\!\!\!\!\!\!\!\!\!\!\!\!\!\!\!\!\!\!\!\!\!\!\!\!\!\!&&
{\rm(iii)\ }
1+O(\d)^2\stackrel{{}^{\sc\rm\eqref{ImMpP}\,(vi)}}{\le}
\tilde\xX_1\aA_3^{-\ell_0}
\stackrel{{}^{\sc\rm\eqref{Now-tildeX-is}}}{=}
(1-\d)\aA_1^{-2-2\d_0}\aA_3^{-2\ell_0-1}
\stackrel{{}^{\sc\rm\eqref{a----mswmw}\,(ii)}}{\le}
\aA_1^{1+O(\d_0)^1}+O(\d)^1.
\!\!\!\!\!\!\!\!\!\!\!\!
\end{eqnarray}
}%
\NOUSE{%
Note that $\yY=\aA_1^{-\d_0}\aA_3>\aA_1^{-\d_0+\d^3}\ge1$ by \eqref{A1-A2-cond}\,(a),\,\eqref{ImMpP}\,(i).
Since the function $\gamma_1(\yY)$ on $\yY$ is  strictly decreasing when $1\le\yY\le\yY_0$ and strictly increasing
when $\yY>\yY_0$, where $\yY_0=\big(\frac{(1-\d^2)^{-1}\ell_0(1-\d_0)}{2(1-\d^2)}\big)^{\d_0}$
(as $\frac{d\gamma_1}{d\yY}|_{\yY=\yY_0}=0$), in order for \eqref{x2-00ro4m4m4}\,(i) to hold, we must have
\equa{yy-is-bigger}{\dis
\yY>\frac43,}
by noting that $\gamma_1(\frac43)<3$.
Then up to $O(\d)^3$, we have
\begin{eqnarray}
\label{ABC-LetNSoOP----1}
&\!\!\!\!\!\!\!\!\!\!\!\!\!\!\!\!\!\!\!\!\!\!\!\!\!\!\!\!\!\!\!\!\!\!\!\!\!\!\!\!\!\!\!\!\!\!\!\!\!\!\!
&
{\rm(i)\,}
\aA_1^2\tilde\xX_1\aA_3^{\ell_0}
\stackrel{{}^{\sc\rm\eqref{tX1-is=}}}{=}(1-\d^2)\aA_1^{2-\d_0}\aA_3^{2\ell_0-1}
\stackrel{{}^{\sc\rm\eqref{x2-00ro4m4m4}\,(ii)}}{=}
(1-\d^2)\yY
\nonumber\\
&\!\!\!\!\!\!\!\!\!\!\!\!\!\!\!\!\!\!\!\!\!\!\!\!\!\!\!\!\!\!\!\!\!\!\!\!\!\!\!\!\!\!\!\!\!\!\!\!\!\!\!
&
\aA_1^3\stackrel{{}^{\sc\rm\eqref{cont-B2}}}{=}
\aA_2\stackrel{{}^{\sc\rm\eqref{LetNSoOP----1}\,(ii),\,(iii)}}
{=}\aA_3^{\ell_0^5}\tilde \xX_1^{-1}\Big |3 - \frac{2}{A_1^2\widetilde  X_1 X_2^{\ell_0}}\Big|
\ge\aA_3^{\ell_0^5}\tilde \xX_1^{-1}\aA_1^{-2}\tilde  \xX_1^{-1}\aA_3^{-\ell_0}(2-3\aA_1^2\tilde\xX_1\aA_3^{\ell_0})
\!\!\!\!\!\!\!\!\!\!\!\!\!\!\!\!\!\!\!
\end{eqnarray}
which is a contradiction.
This proves  the first inequality of \eqref{ImMpP}\,(ii).
Next assume the last  inequality of \eqref{ImMpP}\,(ii) does not hold. Then $\tilde\xX_1\xX_2=(1+\d)\aA_1^{\frac{13}{14}}$ by
 \eqref{C+LetNSoOP}\,(iii). Using the second equality of \eqref{We0o0o0o},
up to $O(\d)^2$, we have%
,
\begin{eqnarray}
\label{x2-00ro4m4m4++}
\!\!\!\!\!\!\!\!\!\!\!\!\!\!\!\!\!\!\!\!\!\!\!\!\!\!\!\!\!\!&&
3
\stackrel{{}^{\sc\rm\eqref{We0o0o0o},\,\eqref{LetNSoOP----1}\,(iii)}}{\le}
\frac{\aA_1}{\aA_2(\tilde\xX_1\xX_2)^{11}}\Big(2+\frac1{(\tilde\xX_1\xX_2)^3}\Big)
=2(1+\d)^{-11}\aA_1^{1+12-11\times\frac{13}{14}}+(1+\d)^{-14}\aA_1^{1+12-14\times\frac{13}{14}}
\!\!\!\!\!\!\!\!\!\!\!\!\!\!\!\!\!\!\!\!
\nonumber\\[-3pt]
\!\!\!\!\!\!\!\!\!\!\!\!\!\!\!\!\!\!\!\!\!\!\!\!\!\!\!\!\!\!&&\ \
\stackrel{{}^{\sc\rm\phantom{\eqref{cont-B2},\,\eqref{ImMpP}\,(i)}}}{=}
2(1+\d)^{-11}\aA_1^{\frac{39}{14}}+(1+\d)^{-14}\stackrel{{}^{\sc\rm\eqref{A1-A2-cond}\,(a)}}{<}3-
36\d\mbox{ \ [up to $O(\d)^2$]},
\!\!\!\!\!\!\!\!\!\!\!\!\!\!\!\!\!\!\!\!
\end{eqnarray}
which is a contradiction. This proves \eqref{ImMpP}\,(ii).
The first inequality of \eqref{ImMpP}\,(iii) follows from \eqref{ImMpP}\,(i).
Assume the inequality of \eqref{ImMpP}\,(iii) does not hold. Then $\xX_2=\ell_2$ by \eqref{C+LetNSoOP}\,(iv).
Then
\begin{eqnarray}
\label{x2-00ro4m4m4++2}
\!\!\!\!\!\!\!\!\!\!\!\!\!\!\!\!\!\!\!\!\!\!\!\!\!\!\!\!\!\!&&
\Big|2-\frac1{\widetilde X_1^3X_2^3}\Big|
\stackrel{{}^{\sc\rm\eqref{LetNSoOP----1}\,(i)}}{=}
\aA_1(\tilde\xX_1\xX_2)^2\xX_2^{-13}
\stackrel{{}^{\sc\rm\eqref{A1-A2-cond}\,(a),\,\eqref{C+LetNSoOP}\,(iii)}}{=}\d_2^{13}
\!\!\!\!\!\!\!\!\!\!\!\!\!\!\!\!\!\!\!\!
\nonumber\\[-3pt]
\!\!\!\!\!\!\!\!\!\!\!\!\!\!\!\!\!\!\!\!\!\!\!\!\!\!\!\!\!\!&&\ \
\stackrel{{}^{\sc\rm\phantom{\eqref{cont-B2},\,\eqref{ImMpP}\,(i)}}}{=}
2(1+\d)^{-11}\aA_1^{\frac{39}{14}}+(1+\d)^{-14}\stackrel{{}^{\sc\rm\eqref{A1-A2-cond}\,(a)}}{<}3-
36\d\mbox{ \ [up to $O(\d)^2$]},
\!\!\!\!\!\!\!\!\!\!\!\!\!\!\!\!\!\!\!\!
\end{eqnarray}
does not hold.
Then $\xX_2=1-\d$ by \eqref{C+LetNSoOP}\,(iv).
We have
\begin{eqnarray}
\label{x2-00ro4m4m4+}
\!\!\!\!\!\!\!\!\!\!\!\!\!\!\!\!\!\!\!\!\!\!\!\!\!\!\!\!\!\!&&
\aA_1^{\eE_1}
\stackrel{{}^{\sc\rm\eqref{LetNSoOP----1}\,(i)}}{=}
\tilde\xX_1^2\zZ^{2\eE_1}\xX_2^{-\eE_1}
\stackrel{{}^{\sc\rm\eqref{ImMpP}\,(i),\,(iii)}}{\le}
(1+\eE_2)\Big(1+O(\d)^2\Big)^{2\eE_1}\xX_2^{\eE_1}
\\\nonumber
\!\!\!\!\!\!\!\!\!\!\!\!\!\!\!\!\!\!\!\!\!\!\!\!\!\!\!\!\!\!&&
\phantom{\aA_1^{\eE_1}}
\stackrel{{}^{\sc\rm\eqref{ImMpP}}}{=}
(1+\eE_2)\Big(1+O(\d)^2\Big)^{2\eE_1}(1-\d)^{\eE_1}
=(1+\eE_2)\Big(1-\d\Big(1+O(\d)^1\Big)\eE_1+O(\eE_1)^2\Big)<1,\!\!\!\!\!\!\!\!
\end{eqnarray}
which is a contradiction with \eqref{A1-A2-cond}.
This proves  the first inequality of \eqref{ImMpP}\,(iv).
Now assume $\xX_2\ge\nn_1^{\nn_1}$. Then
\begin{eqnarray*}
\!\!\!\!\!\!\!\!\!\!\!\!\!\!\!\!\!\!\!\!\!\!\!\!\!\!\!\!\!\!&&
\aA_1^{\eE_1}
\stackrel{{}^{\sc\rm\eqref{LetNSoOP----1}\,(i)}}{=}
\tilde\xX_1^2\zZ^{2\eE_1}\xX_2^{-\eE_1}
\stackrel{{}^{\sc\rm\eqref{ImMpP}\,(i),\,(iii)}}{\ge}
(1-\d)^2\Big(1+O(\d)^2\Big)^{2\eE_1}\xX_2^{\eE_1}
\ge(1-\d)^2\Big(1+O(\d)^2\Big)^{2\eE_1}\nn_1
,\!\!\!\!\!\!\!\!
\end{eqnarray*}
which is a contradiction with \eqref{A1-A2-cond}.
}%
}

Finally, \eqref{ImMpP}\,(v) follows from \eqref{LetNSoOP----1}\,(i) together with
\eqref{LetNSoOP}\,(i),\,\eqref{ImMpP}\,(i)--(iii).
This completes the prove of Lemma \ref{Step1}.
\hfill$\Box$\vskip7pt

\NOUSE
{
We will need to apply frequently formula \eqref{(T0o(eE)1}. To obtain that formula, we
require the following [noting that  $\widetilde T_0$ is the set consisting of absolute values of
non-monomial factors
appearing in \eqref{LetNSoOP----1}\,(ii),\,(iii)
$\ssc\,$],
\NOUSE{We already see from \eqref{A1-A2-cond} that
\eqref{ImMpP}\,(ii) holds for $a=\aA_1,\aA_2$, thus also holds for $a=\tilde\xX_1,\xX_1$ by \eqref{tX1==},\,\eqref{C+ToSayas}\,(c). Then by \eqref{LetNSoOP----1}\,(ii), we see that $\aA_3\le1+O(\d)^1$, from this and \eqref{ImMpP}\,(i),  $\aA_3$ is also a $1+O(\d)^1$. By \eqref{Case6-lemm}\,(iv),
$A_3=\frac{Z^3}{\widetilde X_1^{22} X_2^4}=\frac{1}{\widetilde X_1^{22}X_2}+O(\d)^2$, this implies that $\xX_2,\zZ$ are also $1+O(\d)^1$ element.
This proves \eqref{ImMpP}\,(ii). From this, we obtain
}
\NOUSE
{
\begin{eqnarray}
\label{Amsmene}
&&\!\!\!\!\!\!\!\!\!\!\!\!\!\!\!\!\!\!\!\!\!\!\!\!
\d_2^{\ell_1}< a<\ell_2^{\ell_1}\mbox{ \ \ for \ \ }a\in
\widetilde T_2,\mbox{ where \ }
\widetilde T_2{\sc\!}={\sc\!}T_2\cup\widetilde T_0,\ \
\\\nonumber
&&\!\!\!\!\!\!\!\!\!\!\!\!\!\!\!\!\!\!\!\!\!\!\!\!
\widetilde T_0{\sc}={\sc}\mbox{\Large$\Big\{$}
\Big|\frac{1 + \d_0}{1 + 2 \d_0} + \frac{\d_0\widetilde X_1^{\ell_0^2}}{1 + 2 \d_0}\Big|
,
|1 + \d_0 - \d_0\widetilde X_1^{\ell_0^2}|
%
\mbox{\Large$\Big\}$}
.\!\!\!\!\!\!\!\!\!\!\!\!
\end{eqnarray}
First by \eqref{meme}
,\,\eqref{equa-Case6-lemm}\,(iv),\,\eqref{tX1==},\,\eqref{A1-A2-cond}, we see
that \eqref{Amsmene} holds for $a=\tilde\xX_1,\xX_1,\xX_2,\zZ,\aA_1,\aA_2$. Then we see from
\eqref{LetNSoOP----1}\,(i) 
that it also holds for $a=\aA_3
$.
Thus it holds for elements in $T_2$ [which contains absolute values of monomial factors
\eqref{LetNSoOP----1}$\ssc\,$].
Then by  \eqref{LetNSoOP----1}\,(ii),\,(iii),
~we see that it holds for all $a\in\tilde T_2$.
Therefore, by the fact that $0<\d=\ell^{-1}\ll\d_2=\ell_2^{-1}\ll\d_1=\ell_1^{-1}$, we have,
\equa{(T0o(eE)1}{\dis
a^{0+O(\d)^1}=1+O(\d)^1,\ \ \ a\Big(1+O(\d)^1\Big)=a+O(\d)^1
\mbox{ \ \ for \ \ }a\in\widetilde T_2\mbox{ \ or \ }a^{-1}\in\widetilde T_2.
}
Throughout the rest of the paper, we use the following notations,
\equa{Notattt}{\dis
{\rm(i)\ }\equiv\mbox{ means ``equal up to $O(\d)^2$''},
\ \ \ {\rm(ii)\ }\preccurlyeq \mbox{ means ``smaller than up to $O(\d)^2$''}.
}
Then by \eqref{(T0o(eE)1},\,\eqref{equa-Case6-lemm}\,(iv), we have
\begin{eqnarray}
\label{LetNSoOP----1-redefine++}
&\!\!\!\!\!\!\!\!\!\!\!\!\!\!\!\!\!\!\!\!\!\!\!\!\!\!\!\!\!\!\!\!\!\!
&
{\rm(i)\ }1
\stackrel{{}^{\sc\rm\eqref{LetNSoOP----1}\,(ii)}}{\equiv}
\aA_1^{-1}\tilde \xX_1^2\Big|\frac{1 + \d_0}{1 + 2 \d_0} + \frac{\d_0\widetilde X_1^{\ell_0^2}}{1 + 2 \d_0}\Big|
,
\ \ \ \ \
{\rm(ii)\ }1\stackrel{{}^{\sc\rm\eqref{LetNSoOP----1}\,(ii)}}{\equiv}
\aA_2^{-1}\aA_3\tilde \xX_1^2 |1 + \d_0 - \d_0\widetilde X_1^{\ell_0^2}|
,
\!\!\!\!\!\!\!\!\!\!\!\!\!\!\!\!\!\!\!\!\!\!
\nonumber\\
&\!\!\!\!\!\!\!\!\!\!\!\!\!\!\!\!\!\!\!\!\!\!\!\!\!\!\!\!\!\!\!\!\!\!
&
{\rm(iii)\ }
1\equiv
\frac{\aA_3\tilde\xX_1^4}{\aA_1\aA_2
}\Big|\frac{(1+\d_0)^2}{1+2\d_0}-\frac{\d_0^2\widetilde\xX_1^{2\ell_0^2}}{1+2\d_0}\Big|
.
\end{eqnarray}
Now assume the first inequality of \eqref{ImMpP}\,(ii) does not holds. Then $\tilde\xX_1=1-\d$ by
\eqref{C+LetNSoOP}\,(iii). We have (recalling that $\d_0\ell_0=1$ and $\d\ll\d_0$)
\equa{EEMdcmdmdm}{\dis
1\stackrel{{}^{\rm\sc\eqref{LetNSoOP----1-redefine++}\,(i),\,\eqref{A1-A2-cond}\,(a)}}{\preccurlyeq}
(1-\d)^2\Big(\frac{1+\d_0}{1+2\d_0}+\frac{\d_0(1-\d)^{\ell_0^2}}{1+2\d_0}\Big)
\equiv1-\Big(2+\frac{\ell_0}{1+2\d_0})\d<1-\d,
}
which is a contradiction.
Next assume the last inequality of \eqref{ImMpP}\,(ii) does not holds. Then $\tilde\xX_1=(1+\d)\aA_1^{-\frac14}$ by
\eqref{C+LetNSoOP}\,(iii). We have
\begin{eqnarray}
\label{EEMdcmdmdm+1}
\!\!\!\!\!\!\!\!\!\!\!\!\!\!\!\!\!\!\!\!\!\!\!\!\!\!\!\!\!\!\!\!\!&&
\frac{(1+\d_0)^2}{1+2\d_0}\ \ \ \
\stackrel{{}^{\rm\sc\eqref{LetNSoOP----1-redefine++}\,(iii),\,\eqref{A1-A2-cond}\,(a)}}{\preccurlyeq}
\ \ \
\frac{\d_0^2\tilde\xX_1^{2\ell_0^2}}{1+2\d_0}
+\frac{\aA_1\aA_2}{\aA_3\tilde\xX_1^{4}}
\equiv
\frac{\d_0^2(1+\d)^{2\ell_0^2}\aA_1^{-\frac{\ell_0^2}{2}}}{1+2\d_0}
+\frac{\aA_1\aA_2}{(1+\d)^{4}\aA_1^{-1}\aA_3}
\nonumber\\
\!\!\!\!\!\!\!\!\!\!\!\!\!\!\!\!\!\!\!\!\!\!\!\!\!\!\!\!\!\!\!\!\!&&
\phantom{\frac{(1+\d_0)^2}{1+2\d_0}
}
\stackrel{{}^{\sc\rm\eqref{cont-B2},\,\eqref{A1-A2-cond}\,(a),\,\eqref{ImMpP}\,(i)}}{\preccurlyeq}
\frac{\d_0^2(1+\d)^{2\ell_0^2}}{1+2\d_0}
+\frac{1}{(1+\d)^{4}}<1-\d,
\end{eqnarray}
which is again a contradiction.
This proves \eqref{ImMpP}\,(ii).
By \eqref{tX1==},\,\eqref{A1-A2-cond},\,\eqref{ImMpP}\,(ii), we see \eqref{ImMpP}\,(iii) holds for
$a=\aA_1,\aA_3,\tilde\xX_1,\xX_1$. By \eqref{LetNSoOP----1}\,(i),\,\eqref{ImMpP}\,(i), we have $\xX_2=\aA_3\ge1+O(\d)^1$.
Using \eqref{LetNSoOP----1-redefine++}\,(ii), we can obtain that $\aA_3\le1+O(\d)^1$. This proves
\eqref{ImMpP}\,(iii).\hfill$\Box$
\NOUSE{
Now assume the first inequality of \eqref{ImMpP}\,(iv) does not hold. Then $\tilde\xX_1=(1-\d)\aA_1^{-1}$ by \eqref{C+LetNSoOP}\,(iii).
Then up to $O(\d)^2$, we have,
\begin{eqnarray}
&&\!\!\!\!\!\!\!\!\!\!\!\!\!\!\!\!\!\!\!\!\!\!\!\!\!\!\!
\label{We-Hksks}
2\stackrel{{}^{\sc\rm\eqref{LetNSoOP----1-redefine++}\,(iv),\,\eqref{ImMpP}\,(i)}}{\le}
(1-\d)^{10}+(1-\d)^{275\ell_0+30}\aA_1^{1-\frac{5}{14}\times(\frac{12}{67}+O(\d_0)^1)-(275\ell_0+30)}
\nonumber\\
&&\!\!\!\!\!\!\!\!\!\!\!\!\!\!\!\!\!\!\!\!\!\!\!\!\!\!\!
\ \ \ \ \ \ \ \ \ \ \ \ \,
=\ \ \ \ \ (1-\d)^{10}+(1-\d)^{275\ell_0+30}\aA_1^{-275\ell_0(1+O(\d_0)^1)}
\nonumber\\
&&\!\!\!\!\!\!\!\!\!\!\!\!\!\!\!\!\!\!\!\!\!\!\!\!\!\!\!
\ \ \ \ \ \ \ \ \stackrel{{}^{\sc\rm\eqref{A1-A2-cond}\,(a)}}{\le}\ \ \
(1-\d)^{10}+(1-\d)^{275\ell_0+30}=1-(275\ell_0+40)\d \mbox{ \ \ [up to $O(\d)^2\ssc\,$]},
\end{eqnarray}
which is a contradiction.
Next assume the last inequality of \eqref{ImMpP}\,(iv) does not hold.
Then $\tilde\xX_1=(1+\d)\aA_1^{-\d_0^2}$ by \eqref{C+LetNSoOP}\,(iii).
Then up to $O(\d)^2$, we have,
\begin{eqnarray}
&&\!\!\!\!\!\!\!\!\!\!\!\!\!\!\!\!\!\!\!\!\!\!\!\!\!\!\!
\label{We-Hksks++}
8\stackrel{{}^{\sc\rm\eqref{LetNSoOP----1-redefine++}\,(vi),\,\eqref{cont-B2},\,\eqref{ImMpP}\,(i)}}{\le}
3(1+\d)^{10}\aA_1^{-\d_0^2}+5\frac{\aA_1^
{(\frac5{67}+O(\d_0)^1)-5\ell_0(\frac{12}{67}+O(\d_0)^1)-(300\ell_0^2+15\ell_0+6)\d_0^2}}{(1+\d)^{300\ell_0^2+15\ell_0+6}}
\nonumber\\
&&\!\!\!\!\!\!\!\!\!\!\!\!\!\!\!\!\!\!\!\!\!\!\!\!\!\!\!
\ \ \ \ \ \ \ \ \ \ \ \ \ \ \ \ \,
=\ \ \ \ 3(1+\d)^{10}\aA_1^{-\d_0^2}+5(1+\d)^{-300\ell_0^2-15\ell_0-6}\aA_1^{-\frac{60\ell_0}{67}(1+O(\d_0)^1)}
\\\nonumber
&&\!\!\!\!\!\!\!\!\!\!\!\!\!\!\!\!\!\!\!\!\!\!\!\!\!\!\!
\ \ \ \ \ \ \ \ \ \ \ \ \,
\stackrel{{}^{\sc\rm\eqref{A1-A2-cond}\,(a)}}{\le}
\ \ \
3(1+\d)^{10}+5(1+\d)^{-300\ell_0^2-15\ell_0-6}=1-(1500\ell_0^2+75\ell_0)\d \mbox{ \ \ [up to $O(\d)^2\ssc\,$]},
\!\!\!\!\!\!\!\!\!\!\!\!
\end{eqnarray}
which is again a contradiction. This proves \eqref{ImMpP}\,(iv).
By \eqref{LetNSoOP----1-redefine++}\,(ii),\,\eqref{ImMpP}\,(i),\,(iii),
we see that the last inequality of
\eqref{ImMpP}\,(v) holds.
If the first inequality does not hold, i.e., $\xX_2\le\d_2$, then $\aA_3\gg\ell_1$ by
\eqref{LetNSoOP----1-redefine++}\,(i) (recalling that $\ell_2=\d_2^{-1}\gg\ell_1$), and so $\aA_1\gg\ell_1$
by \eqref{LetNSoOP----1-redefine++}\,(iii),
a contradiction with \eqref{A1-A2-cond}\,(a). This proves
\eqref{ImMpP}\,(v).}
\NOUSE{
We already see from \eqref{A1-A2-cond},\,\eqref{C+LetNSoOP}\,(iii) that
 \eqref{ImMpP}\,(iii) holds for $a=\tilde\xX_1,\aA_1,\aA_2$.
Then we have $\aA_3\ge1+O(\d)^1$ by
 \eqref{ImMpP}\,(i).
On the other hand, we have
\equa{Smeme}{\dis
1+O(\d)^1=\aA_1\stackrel{{}^{\sc\rm\eqref{LetNSoOP----1-redefine++}\,(iii)}}{\ge}
\aA_3^{\frac{5}{14}}\Big(1+O(\d)^1\Big)\Big(2-\big(1+O(\d)^1\big)\Big)=\aA_3^{\frac{5}{14}}\Big(1+O(\d)^1\Big),
}
which implies that $\aA_3\le1+O(\d)^1$. Thus
 \eqref{ImMpP}\,(ii) hold for $a=\aA_3$. Then it also hold
 for $a=\xX_2,\zZ,\tilde\xX_1$ by
\eqref{LetNSoOP----1-redefine++}\,(ii),\,\eqref{equa-Case6-lemm}\,(iv),\,\eqref{tX1==}. This proves
 \eqref{ImMpP}\,(iii).
}%
}
\NOUSE{
\equa{MAmams--000}{\dis
A_2=\frac{\widetilde X_1^{22\ell_0+11} X_2}{A_3}\Big (2 - \frac{A_1 \widetilde X_1^2}{A_3^2}\Big)^{\ell_0}
\stackrel{{}^{\sc\rm\eqref{equa-Case6-lemm}\,(iv),\,\eqref{LetNSoOP----1}\,(i)}}{=}
\widetilde X_1^{22\ell_0+11}\Big (2 - \frac{A_1 \widetilde X_1^2}{A_3^2}\Big)^{\ell_0}+O(\d)^3.
}
By  \eqref{LetNSoOP----1}\,(iv),\,(v), up to $O(\d)^3$, we have
\begin{eqnarray}
\label{Apply-LetNSoOP----1}
&\!\!\!\!\!\!\!\!\!\!\!\!\!\!\!\!\!\!\!\!\!\!\!\!\!\!\!\!\!
&
\bB_1^{-3}\cC_1^{2}=(\aA_1
\aA_3^{-2}\tilde\xX_1^{2})^{-3}
(\aA_1
\aA_3^{-(44\ell_0+26)}\tilde\xX_1^{3})^{2}
=\aA_1^{-1}\aA_3^{-(88\ell_0+46)}\stackrel{{}^{\sc\rm\eqref{A1-A2-cond}\,(a),\,\eqref{ImMpP}}}{\le}1+O(\d)^3.
\NOUSE{
{\rm(i)\ }A_3\stackrel{{}^{\sc\rm\eqref{LetNSoOP----1}\,(i)}}{=}\frac{\widetilde X_1^{81} X_2^{46}}{Z^{45}}
=\widetilde X_1^{81} X_2+O(\d)^3,
\ \ \ \ \ \
{\rm(ii)\ }X_2\stackrel{{}^{\sc\rm\eqref{Apply-LetNSoOP----1}\,(i)}}{=}A_3\widetilde X_1^{-81}+O(\d)^3,
\!\!\!\!\!\!\!\!\!
\!\!\!\!\!\!\!\!\!
\!\!\!\!\!\!\!\!\!
\!\!\!\!\!\!\!\!\!
\nonumber\\
&\!\!\!\!\!\!\!\!\!\!\!\!\!\!\!\!\!\!\!\!\!\!\!\!\!\!\!\!\!\!\!\!\!\!\!\!\!\!\!\!\!\!\!\!\!\!\!\!\!\!\!
&
{\rm(iii)\ }
A_1\stackrel{{}^{\sc\rm\eqref{LetNSoOP----1}\,(ii)}}{=}\frac{X_2^3}{A_3^2\widetilde X_1^4 Z^3(2 -\widetilde X_1^{10})}
=\frac1{A_3^2\widetilde X_1^4(2 -\widetilde X_1^{10})}+O(\d)^3,
\!\!\!\!\!\!\!\!\!
\!\!\!\!\!\!\!\!\!
\!\!\!\!\!\!\!\!\!
\!\!\!\!\!\!\!\!\!
\nonumber\\
&\!\!\!\!\!\!\!\!\!\!\!\!\!\!\!\!\!\!\!\!\!\!\!\!\!\!\!\!\!\!\!\!\!\!\!\!\!\!\!\!\!\!\!\!\!\!\!\!\!\!\!
&
{\rm(iii)\ }A_2\stackrel{{}^{\sc\rm\eqref{LetNSoOP----1}\,(iii)}}{=}\frac{Z}{A_3\Big (\frac95 - \frac{4\widetilde X_1^{10}}{5}\Big)}
\stackrel{{}^{\sc\rm\eqref{Apply-LetNSoOP----1}\,(ii)}}{=}
\frac{1}{\widetilde X_1^{81}\Big (\frac95 - \frac{4\widetilde X_1^{10}}{5}\Big)}+O(\d)^3,
%
\!\!\!\!\!\!\!\!\!
\!\!\!\!\!\!\!\!\!
\!\!\!\!\!\!\!\!\!
\!\!\!\!\!\!\!\!\!
\nonumber\\
&\!\!\!\!\!\!\!\!\!\!\!\!\!\!\!\!\!\!\!\!\!\!\!\!\!\!\!\!\!\!\!\!\!\!\!\!\!\!\!\!\!\!\!\!\!\!\!\!\!\!\!
&
{\rm(iv)\ }
1\stackrel{{}^{\sc\rm\eqref{meme01020932}}}{=}
\frac{A_2 A_3}{Z}\Big (\frac15 +
\frac{4 X_2^3}{5 A_1 A_3^2\widetilde X_1^4 Z^3}\Big)
=A_2\widetilde X_1^{81}\Big (\frac15 +
\frac{4}{5 A_1 A_3^2\widetilde X_1^4}\Big)+O(\d)^3,
%
\!\!\!\!\!\!\!\!\!
\!\!\!\!\!\!\!\!\!
\!\!\!\!\!\!\!\!\!
\!\!\!\!\!\!\!\!\!
\nonumber\\[4pt]
&\!\!\!\!\!\!\!\!\!\!\!\!\!\!\!\!\!\!\!\!\!\!\!\!\!\!\!\!\!\!\!\!\!\!\!\!\!\!\!\!\!\!\!\!\!\!\!\!\!\!\!
&
{\rm(v)\ }|B_2|^{-1}\stackrel{{}^{\sc\rm\eqref{LetNSoOP----1}\,(iv)}}{=}
\aA_1^{-1}\aA_2 \aA_3^{-1}\tilde \xX_1^{-4}\xX_2^3 \zZ^{-4}
\stackrel{{}^{\sc\rm\eqref{equa-Case6-lemm}\,(iv),\,\eqref{Apply-LetNSoOP----1}\,(ii)}}{=}
\aA_1^{-1}\aA_2\aA_3^{-2}\tilde\xX_1^{77}+O(\d)^3,
\!\!\!\!\!\!\!\!\!
\!\!\!\!\!\!\!\!\!
\!\!\!\!\!\!\!\!\!
\!\!\!\!\!\!\!\!\!
\nonumber\\[7pt]
&\!\!\!\!\!\!\!\!\!\!\!\!\!\!\!\!\!\!\!\!\!\!\!\!\!\!\!\!\!\!\!\!\!\!\!\!\!\!\!\!\!\!\!\!\!\!\!\!\!\!\!
&
{\rm(vi)\ }|C_2|^{-1}\stackrel{{}^{\sc\rm\eqref{LetNSoOP----1}\,(v)}}{=}
\aA_1^{-1}\aA_2^2\tilde \xX_1^{-4}\xX_2^3 \zZ^{-5}
\stackrel{{}^{\sc\rm\eqref{equa-Case6-lemm}\,(iv),\,\eqref{Apply-LetNSoOP----1}\,(ii)}}{=}\aA_1^{-1}\aA_2^2\aA_3^{-2}\tilde\xX_1^{158}+O(\d)^3,
\!\!\!\!\!\!\!\!\!
\!\!\!\!\!\!\!\!\!
\!\!\!\!\!\!\!\!\!
\!\!\!\!\!\!\!\!\!
\nonumber\\[7pt]
&\!\!\!\!\!\!\!\!\!\!\!\!\!\!\!\!\!\!\!\!\!\!\!\!\!\!\!\!\!\!\!\!\!\!\!\!\!\!\!\!\!\!\!\!\!\!\!\!\!\!\!
&
{\rm(vii)\ }|B_2C_2^4|^{-1}=
\aA_1^{-5}\aA_2^9\aA_3^{-10}\tilde\xX_1^{709}
\stackrel{{}^{\sc\rm\eqref{cont-B2},\,\eqref{ImMpP}\,(i),\,\eqref{C+ToSayas+1}\,(c)}}{\le}
\aA_1^{-5+9\times\frac{27}{65}+10\times\frac{34}{65}-709\times\frac{7}{1250}}+O(\d)^2
\!\!\!\!\!\!\!\!\!
\!\!\!\!\!\!\!\!\!
\!\!\!\!\!\!\!\!\!
\!\!\!\!\!\!\!\!\!
\nonumber\\[7pt]
&\!\!\!\!\!\!\!\!\!\!\!\!\!\!\!\!\!\!\!\!\!\!\!\!\!\!\!\!\!\!\!\!\!\!\!\!\!\!\!\!\!\!\!\!\!\!\!\!\!\!\!
&
\phantom{{\rm(vii)\ }|B_2C_2^4|^{-1}}
=\aA_1^{-5+\frac{243}{65}+\frac{68}{13}-\frac{4963}{1250}}+O(\d)^2=
\aA_1^{-\frac{19}{16250}}+O(\d)^2\stackrel{{}^{\sc\rm\eqref{A1-A2-cond}\,(a)}}{\le}1+O(\d)^2.
}\end{eqnarray}
}\NOUSE{
\begin{eqnarray}
\label{Asno-LetNSoOP----1}
&\!\!\!\!\!\!\!\!\!\!\!\!\!\!\!\!\!\!\!\!\!\!\!\!\!\!\!\!\!\!\!\!\!\!\!\!\!\!\!\!\!\!\!\!\!\!\!\!\!\!\!
&
{\rm(i)\ }A_3\stackrel{{}^{\sc\rm\eqref{LetNSoOP----1}\,(i)}}{=}\frac{Z^3}{\widetilde X_1^{22} X_2^4}=\frac{1}{\widetilde X_1^{22}X_2}+O(\d)^2,
\ \ \
{\rm(ii)\ }X_2=\frac{1}{A_3\widetilde X_1^{22}}+O(\d)^2,
\!\!\!\!\!\!\!\!\!
\!\!\!\!\!\!\!\!\!
\!\!\!\!\!\!\!\!\!
\!\!\!\!\!\!\!\!\!
\nonumber\\
&\!\!\!\!\!\!\!\!\!\!\!\!\!\!\!\!\!\!\!\!\!\!\!\!\!\!\!\!\!\!\!\!\!\!\!\!\!\!\!\!\!\!\!\!\!\!\!\!\!\!\!
&{\rm(iii)\ }A_2\stackrel{{}^{\sc\rm\eqref{LetNSoOP----1}\,(ii)}}{=}\frac{X_2^2\Big(\frac{5}{3} - \frac{2\widetilde X_1^{15}}{3}\Big)}{Z\Big(\frac15 + \frac{4\widetilde X_1^{15}}{5}\Big)}
=\frac{X_2\Big(\frac{5}{3} - \frac{2\widetilde X_1^{15}}{3}\Big)}{\frac15 + \frac{4\widetilde X_1^{15}}{5}}+O(\d)^2
\!\!\!\!\!\!\!\!\!
\!\!\!\!\!\!\!\!\!
\!\!\!\!\!\!\!\!\!
\!\!\!\!\!\!\!\!\!
\nonumber\\
&\!\!\!\!\!\!\!\!\!\!\!\!\!\!\!\!\!\!\!\!\!\!\!\!\!\!\!\!\!\!\!\!\!\!\!\!\!\!\!\!\!\!\!\!\!\!\!\!\!\!\!
&\phantom{{\rm(iii)\ }A_2}
\stackrel{{}^{\sc\rm\eqref{LetNSoOP----1}\,(iii)}}{=}
\frac{X_2\Big(\frac23+\frac1{3A_1A_3\widetilde\xX_{15}}\Big)}{\frac75-\frac{2}{5A_1A_3\widetilde X_1^{15}}}+O(\d)^2
\stackrel{{}^{\sc\rm\eqref{Asno-LetNSoOP----1}\,(ii)}}{=}
\frac{\frac23+\frac1{3A_1A_3\widetilde\xX_{15}}}{A_3\widetilde X_1^{22}\Big(\frac75-\frac{2}{5A_1A_3\widetilde X_1^{15}}\Big)}+O(\d)^2
\NOUSE{,
\!\!\!\!\!\!\!\!\!
\!\!\!\!\!\!\!\!\!
\!\!\!\!\!\!\!\!\!
\!\!\!\!\!\!\!\!\!
\nonumber\\
&\!\!\!\!\!\!\!\!\!\!\!\!\!\!\!\!\!\!\!\!\!\!\!\!\!\!\!\!\!\!\!\!\!\!\!\!\!\!\!\!\!\!\!\!\!\!\!\!\!\!\!
&
{\rm(iv)\ }B_3=\frac{\widetilde X_1^4Z^{2\ell_0+1}}{A_3X_2^{2\ell_0-2}},\ \ \ \ \ \ \
{\rm(v)\ }C_3=\frac{\widetilde X_1^6Z^{3\ell_0+2}}{A_3^2X_2^{3\ell_0-2}}
}%
.
\!\!\!\!\!\!\!\!\!\!\!\!\!\!\!
\end{eqnarray}
Now assume the first inequality of \eqref{ImMpP}\,(iii) does not hold. Then $\tilde\xX_1=(1-\d)\aA_1^{\frac1{50}}$ by \eqref{C+ToSayas+1}\,(c).
Up to $O(\d)^2$, by \eqref{LetNSoOP----1}\,(ii), we have
\begin{eqnarray}
\label{A1==soso}
\!\!\!\!\!\!\!\!\!\!\!\!\!\!\!\!\!\!\!\!\!\!\!\!\!\!\!&&
3\le2\tilde\xX_1^{15}+\aA_1^{-1}\aA_3^{-1}\tilde\xX_1^{-15}
\stackrel{{}^{\sc\rm\eqref{ImMpP}\,(i)}}{<}
2(1-\d)^{15}\aA_1^{\frac3{10}}+(1-\d)^{-15}\aA_1^{-1+\frac58-\frac3{10}}
\nonumber\\
\!\!\!\!\!\!\!\!\!\!\!\!\!\!\!\!\!\!\!\!\!\!\!\!\!\!\!&&
\phantom{2}=\gamma_1(\aA_1):=2(1-\d)^{15}\aA_1^{\frac3{10}}+(1-\d)^{-15}\aA_1^{-\frac{27}{40}}
\mbox{ \ \ [up to $O(\d)^2$]}.
\end{eqnarray}
Observe that $\gamma_1(\aA_1)$ is a strictly decreasing function on $\aA_1$ in a small neighborhood at $\aA_1=1$ as $\frac{d\gamma_1}{d\aA_1}|_{\aA_1=1}=\frac{3(1-\d)^{15}}{5}-\frac{27(1-\d)^{-15}}{40}
=-\frac{3}{40}+O(\d)^1<0$. Thus by the fact from \eqref{A1-A2-cond} that $\aA_1\ge1$ and
by \eqref{ImMpP}\,(ii), we obtain from \eqref{A1==soso} the following, up to $O(\d)^2$,
\equa{A1==soso+}{
3<\gamma_1(\aA_1)<\gamma_1(1)=2(1-\d)^{15}+(1-\d)^{-15}=3-15\d
\mbox{ \ \ [up to $O(\d)^2$]},
}
which is a contradiction. This proves the first inequality of \eqref{ImMpP}\,(iii).
Now assume \eqref{ImMpP}\,(ii) does not hold. Then
$\xX_1=(1-\d^2)\aA_1^{-\frac1{65}}$ by \eqref{C+ToSayas+1}\,(c).
We will need to apply formula \eqref{(T0o(eE)1} frequently.
To obtain this formula, first we have, where $T_2$ is defined in \eqref{denote-t2} [recalling from \eqref{MSmde33333} that $\d_2=\ell_2^{-1}$],
This can be easily seen from
\eqref{equa-Case6-lemm}\,(iv),\,\eqref{tX1==},\,\eqref{LetNSoOP----1}\,(i)--(iii),\,\eqref{C+ToSayas+1}\,(d),\,(e)
and \eqref{A1-A2-cond} [see also \eqref{ImMpP+1}\,(v),\,(vi)$\ssc\,$].
\begin{eqnarray}
\label{Another-LetNSoOP----1}\label{LetNSoOP----1+ThisOne}
&\!\!\!\!\!\!\!\!\!\!\!\!\!\!\!\!\!\!\!\!\!\!\!\!\!\!\!\!\!\!\!\!\!\!\!\!\!\!\!\!\!\!\!\!\!\!\!\!\!\!\!
&
{\rm(i)\ }A_1=\frac{X_2^5}{A_3^6Z^2}=\frac{X_2^3}{A_3^6}+O(\d)^3,\ \ \ \implies\ \ \ {\rm(ii)\ }
\xX_2=\aA_1^{\frac13}\aA_3^2+O(\d)^3,
\!\!\!\!\!\!\!\!\!
\!\!\!\!\!\!\!\!\!
\!\!\!\!\!\!\!\!\!
\!\!\!\!\!\!\!\!\!
\nonumber\\
&\!\!\!\!\!\!\!\!\!\!\!\!\!\!\!\!\!\!\!\!\!\!\!\!\!\!\!\!\!\!\!\!\!\!\!\!\!\!\!\!\!\!\!\!\!\!\!\!\!\!\!
&
{\rm(iii)\ }
A_3^{-1}=\frac{X_2^{\ell_0}}{\widetilde X_1^2 Z^{\ell_0+1}}\Big(\frac15 + \frac{4 X_2^{\ell_0-2}}{5\widetilde X_1^2 Z^{\ell_0}}\Big)
=\frac1{\widetilde X_1^2X_2}\Big(\frac15 + \frac{4 }{5\widetilde X_1^2 X_2^2}\Big)+O(\d)^3,
\!\!\!\!\!\!\!\!\!
\!\!\!\!\!\!\!\!\!
\!\!\!\!\!\!\!\!\!
\!\!\!\!\!\!\!\!\!
\nonumber\\
&\!\!\!\!\!\!\!\!\!\!\!\!\!\!\!\!\!\!\!\!\!\!\!\!\!\!\!\!\!\!\!\!\!\!\!\!\!\!\!\!\!\!\!\!\!\!\!\!\!\!\!
&
{\rm(iv)\ }A_2=\frac{\widetilde X_1^2Z^{\ell_0}}{A_3^{10}X_2^{\ell_0-2}\Big(2 - \d_0 - \frac{(1 - \d_0) X_2^{\ell_0-2}}{
\widetilde X_1^2 Z^{\ell_0}}\Big)}
=
\frac{\widetilde X_1^2X_2^2}{A_3^{10}\Big(2 - \d_0 - \frac{1 - \d_0}{
\widetilde X_1^2 X_2^2}\Big)}+O(\d)^3
,
\!\!\!\!\!\!\!\!\!
\!\!\!\!\!\!\!\!\!
\!\!\!\!\!\!\!\!\!
\!\!\!\!\!\!\!\!\!
\nonumber\\
&\!\!\!\!\!\!\!\!\!\!\!\!\!\!\!\!\!\!\!\!\!\!\!\!\!\!\!\!\!\!\!\!\!\!\!\!\!\!\!\!\!\!\!\!\!\!\!\!\!\!\!
&
{\rm(iv)\ }|B_3|=\aA_3^{-1}\tilde \xX_1^4\zZ^{2\ell_0+1}\xX_2^{-2\ell_0+2}
=\aA_3^{-1}\tilde\xX_1^4\xX_2^3+O(\d)^3\stackrel{{}^{\sc\rm\eqref{Another-LetNSoOP----1}\,(ii)}}{=}
\aA_1\aA_3^5\tilde\xX_1^4+O(\d)^3
,
\!\!\!\!\!\!\!\!\!
\!\!\!\!\!\!\!\!\!
\!\!\!\!\!\!\!\!\!
\!\!\!\!\!\!\!\!\!
\nonumber\\
&\!\!\!\!\!\!\!\!\!\!\!\!\!\!\!\!\!\!\!\!\!\!\!\!\!\!\!\!\!\!\!\!\!\!\!\!\!\!\!\!\!\!\!\!\!\!\!\!\!\!\!
&
{\rm(v)\ }|C_3|^{-1}=\aA_3^2\xX_2^{3\ell_0-2}\tilde \xX_1^{-6}\zZ^{-3\ell_0-2}
=\aA_3^2\tilde\xX_1^{-6}\xX_2^{-4}+O(\d)^3=\aA_1^{-\frac43}\aA_3^{-6}\tilde\xX_1^{-6}+O(\d)^3,
\!\!\!\!\!\!\!\!\!
\!\!\!\!\!\!\!\!\!
\!\!\!\!\!\!\!\!\!
\!\!\!\!\!\!\!\!\!
\nonumber\\
&\!\!\!\!\!\!\!\!\!\!\!\!\!\!\!\!\!\!\!\!\!\!\!\!\!\!\!\!\!\!\!\!\!\!\!\!\!\!\!\!\!\!\!\!\!\!\!\!\!\!\!
&
{\rm(vi)\ }
|B_3|^{\frac65}\cdot|C_3|^{-1}
=\aA_1^{-\frac2{15}}\tilde\xX_1^{-\frac65}+O(\d)^3
\stackrel{{}^{\sc\rm\eqref{C+ToSayas+1}\,(d)}}{\le}
\aA_1^{-\frac2{15}+\frac65(\frac{25}{237}-\frac{5\d_0}{1422})}+O(\d)^2
\!\!\!\!\!\!\!\!\!
\!\!\!\!\!\!\!\!\!
\!\!\!\!\!\!\!\!\!
\!\!\!\!\!\!\!\!\!
\nonumber\\
&\!\!\!\!\!\!\!\!\!\!\!\!\!\!\!\!\!\!\!\!\!\!\!\!\!\!\!\!\!\!\!\!\!\!\!\!\!\!\!\!\!\!\!\!\!\!\!\!\!\!\!
&
\phantom{{\rm(vi)\ }
|B_3|^{\frac65}\cdot|C_3|^{-1}}
=\aA_1^{-\frac8{1185}+O(\d_0)^1}+O(\d)^2\stackrel{{}^{\sc\rm\eqref{A1-A2-cond}\,(b)}}{\le}1+O(\d)^2
.
\!\!\!\!\!\!\!\!\!
\end{eqnarray}
}
\NOUSE{Now assume \eqref{ImMpP}\,(ii) does not hold. Then
$\xX_1=(1-\d^2)\aA_1^{-\frac1{65}}$ by \eqref{C+ToSayas+1}\,(c).
Up to $O(\d)^3$, we have
\begin{eqnarray}
\label{A1-X1-small}
&&\!\!\!\!\!\!\!\!\!\!\!\!\!\!\!\!\!\!\!\!\!\!\!\!\!\!\!
1\stackrel{{}^{\sc\rm\eqref{Apply-LetNSoOP----1}\,(iv),\,\eqref{cont-B2},\,\eqref{ImMpP}\,(i)}}{\le}
(1-\d^2)^{81}\aA_1^{\frac{27}{65}-\frac{81}{65}}\Big(\frac15+\frac45(1-\d^2)^{-4}\aA_1^{-1+\frac{68}{65}+\frac{4}{65}}\Big)
\nonumber\\
&&\!\!\!\!\!\!\!\!\!\!\!\!\!\!\!\!\!\!\!\!\!\!\!\!\!\!\!
\phantom{1\ \ \ \ }
=\frac15(1-\d^2)^{81}\aA_1^{-\frac{54}{65}}+\frac45(1-\d)^{77}\aA_1^{-\frac{47}{65}}
\stackrel{{}^{\sc\rm\eqref{A1-A2-cond}}}{\le}
1-\frac{398\d^2}{5}
\ \ \mbox{ \ [up to $O(\d)^3$]},
\end{eqnarray}
which is a contradiction. This proves \eqref{ImMpP}\,(ii).
}\NOUSE{\equa{X1==33000}{\dis\tilde\xX_1=(1-\d^2)
\aA_1^{-\frac{25}{237} + \frac{5 \d_0}{1422}}.}
Then up to $O(\d)^3$, we have
\begin{eqnarray}
\label{A2-to-uu}
&\!\!\!\!\!\!\!\!\!\!\!\!\!\!\!\!\!\!\!\!\!\!\!\!\!\!\!&
2-\d_0\stackrel{{}^{\sc\rm\eqref{Another-LetNSoOP----1}\,(iv)}}{\le}
\aA_2^{-1}\aA_3^{-10}\tilde\xX_1^2\xX_2^2+(1-\d_0)\tilde\xX_1^{-2}\xX_2^{-2}
\nonumber\\
&\!\!\!\!\!\!\!\!\!\!\!\!\!\!\!\!\!\!\!\!\!\!\!\!\!\!\!&
\ \ \ \,\ \
\stackrel{{}^{\sc\rm\eqref{Another-LetNSoOP----1}\,(ii),\,\eqref{X1==33000}}}{=}
(1-\d^2)^2\aA_2^{-1}\aA_1^{2(-\frac{25}{237}+\frac{5\d_0}{1422})+\frac23}\aA_3^{-10+4}+(1-\d_0)(1-\d^2)^{-2}
\aA_1^{-2(-\frac{25}{237}+\frac{5\d_0}{1422})-\frac23}\aA_3^{-4}
\!\!\!\!\!\!\!\!\!\!\!\!\!\!\!\!\!\!\!\nonumber\\
&\!\!\!\!\!\!\!\!\!\!\!\!\!\!\!\!\!\!\!\!\!\!\!\!\!\!\!&
\ \ \ \,\ \
\stackrel{{}^{\sc\rm\eqref{cont-B2},\,\eqref{ImMpP}\,(i)}}{\le}
(1{\sc\!}-{\sc\!}\d^2)^2\aA_1^{-\frac{5(162+\d_0)}{711}+(\frac{36}{79}+\frac{5\d_0}{711})+(-6)(-\frac{9}{79})}
{\sc\!}+{\sc\!}(1{\sc\!}-{\sc\!}\d_0)(1{\sc\!}-{\sc\!}\d^2)^{-2}\aA_1^{(-\frac{36}{79}-\frac{5\d_0}{711})+(-4)(-\frac{9}{79})}
\!\!\!\!\!\!\!\!\!\!\!\!\!\!\!\!\!\!\!\nonumber\\
&\!\!\!\!\!\!\!\!\!\!\!\!\!\!\!\!\!\!\!\!\!\!\!\!\!\!\!&
\ \ \ \ \ \ \ \ \ \ \ \ \
=(1-\d^2)^2+(1-\d_0)(1-\d^2)^{-2}\aA_1^{-\frac{5\d_0}{711}}\le2-\d_0-\d_0\d^2\mbox{ \ \ [up to $O(\d)^3$]},
\end{eqnarray}
which is a contradiction. This proves \eqref{ImMpP}\,(ii).
}
\NOUSE
{Now assume the inequality of \eqref{ImMpP}\,(ii) does not hold, then
$\tilde\xX_1=(1-\d)\aA_2^{-\frac53}$ by \eqref{C+ToSayas+1}\,(c).
Up to $O(\d)^2$, we have
\begin{eqnarray}
\label{A1===sl,we,3}
\!\!\!\!\!\!\!\!\!\!\!\!\!\!\!\!\!\!\!\!\!\!\!\!&&
1\ \ \,\stackrel{{}^{\sc\rm\eqref{LetNSoOP----1}\,(iii)}}{\le}
\aA_1^{-1}\tilde\xX_1^6\xX_2^{-2}\Big(\frac15+\frac45\tilde\xX_1^5\xX_2^{-10}\Big)
\stackrel{{}^{\sc\rm\eqref{LetNSoOP----1+ThisOne}\,(iii)}}{=}
\aA_1^{-1}\tilde\xX_1^6(\aA_2 \aA_3 \tilde\xX_1^{20})^{-\frac2{21}}
\Big(\frac15+\frac45\tilde\xX_1^5(\aA_2 \aA_3 \tilde\xX_1^{20})^{-\frac{10}{21}}\Big)
\!\!\!\!\!\!\!\!\!\!\!\nonumber\\
\!\!\!\!\!\!\!\!\!\!\!\!\!\!\!\!\!\!\!\!\!\!\!\!&&
\ \ \ \ \ \ \ \ \,=\ \ \ \frac15\Big(\aA_1^{-21}\aA_2^{-2}\aA_3^{-2}\tilde\xX_1^{86}\Big)^{\frac1{21}}
+\frac45\Big(\aA_1^{-7}\aA_2^{-4}\aA_3^{-4}\tilde\xX_1^{-3}\Big)^{\frac17}
\!\!\!\!\!\!\!\!\!\!\!\nonumber\\
\!\!\!\!\!\!\!\!\!\!\!\!\!\!\!\!\!\!\!\!\!\!\!\!&&
\stackrel{{}^{\sc\rm\eqref{cont-B2},\,\eqref{ImMpP}\,(i)}}{\le}
\frac15(1-\d)^{\frac{86}{21}}\Big(\aA_2^{105-2-18-\frac{430}{3}}\Big)^{\frac1{21}}+
\frac45(1-\d)^{\frac37}\Big(\aA_2^{35-4-36+5}\Big)^{\frac17}
\!\!\!\!\!\!\!\!\!\!\!\nonumber\\
\!\!\!\!\!\!\!\!\!\!\!\!\!\!\!\!\!\!\!\!\!\!\!\!&&
\ \ \ \ \ \ \ \ \,=\ \ \
\frac15(1-\d)^{\frac{86}{21}}\Big(\aA_2^{-\frac{175}{3}}\Big)^{\frac1{21}}+\frac45(1-\d)^{-\frac37}
\stackrel{{}^{\sc\rm\eqref{A1-A2-cond}\,(b)}}{\le}1-\frac{10\d}{21}
\mbox{ \ [up to $O(\d)^2$]}.
\end{eqnarray}
This is a contradiction. This proves \eqref{ImMpP}\,(ii).\hfill$\Box$\vskip7pt
For convenience, we denote [cf.~\eqref{A1-form1}$\ssc\,$],
\equa{DeD12}{\dis
{\rm(i)\ }D_1:=\frac{X_2^5}{\widetilde X_1^4Z^2}\stackrel{{}^{\sc\rm\eqref{LetNSoOP----1+ThisOne}\,(i)}}{=}\frac{A_3^2}{\widetilde X_1^4}+O(\d)^2,
\ \ \ \
{\rm(ii)\ }D_2=\frac{X_2^7}{\widetilde X_1^{10}}\stackrel{{}^{\sc\rm\eqref{LetNSoOP----1+ThisOne}\,(i)}}{=}\frac{A_3^7}{\widetilde X_1^{10}}+O(\d)^2.
}
Then
\begin{eqnarray}
\label{AnoThA1A2}
&\!\!\!\!\!\!\!\!\!\!\!\!\!\!\!\!\!\!\!\!\!\!\!\!\!\!\!\!\!\!\!\!\!\!\!\!\!\!\!\!\!\!\!\!\!\!\!\!\!\!\!
&
{\rm(i)\ }
A_1\stackrel{{}^{\sc\rm\eqref{LetNSoOP----1}\,(ii)}}{=}
\frac{X_2^5}{\widetilde X_1^4Z^2}\Big(\frac15+\frac{4X_2^7}{5\widetilde X_1^{10}}\Big)
=D_1\Big(\frac15+\frac{4D_2}{5}\Big),
\!\!\!\!\!\!\!\!\!
\nonumber\\
&\!\!\!\!\!\!\!\!\!\!\!\!\!\!\!\!\!\!\!\!\!\!\!\!\!\!\!\!\!\!\!\!\!\!\!\!\!\!\!\!\!\!\!\!\!\!\!\!\!\!\!
&
{\rm(ii)\ }A_2\stackrel{{}^{\sc\rm\eqref{LetNSoOP----1}\,(iii)}}{=}
\frac{X_2^7}{\widetilde X_1^{10}}\Big(3-\frac{2X_2^7}{\widetilde X_1^{10}}\Big)=D_2(3-2D_2),
\ \ \ \ \ \ \
{\rm(iii)\ }
A_1^{-1}\stackrel{{}^{\sc\rm\eqref{AnoThA1A2}\,(i),\,(ii)}}{=}\frac{D_1^{-1}}{\frac75-\frac{2A_2}{5D_2}}.
\!\!\!\!\!\!\!\!\!\!\!\!\!\!
\end{eqnarray}
One can also observe from \eqref{LetNSoOP----1}\,(ii),\,(iii) that the following holds [to see the first equality, simply substitute $A_2$ by
\eqref{LetNSoOP----1}\,(ii)$\ssc\,$],
\equa{MEMEME-A-inv}{\dis\!\!\!\!\!\!\!\!\!\!\!\!
A_1^{-1}=\frac{\frac1{A_2}+9A_3^{10}Z^{90}}{10A_3^{10}X_2^{20}Z^{71}\Big(\frac{49}{45}-\frac{4X_2^{90}}{45A_2A_2^{30}Z^{90}}\Big)}
\stackrel{{}^{\sc\rm\eqref{equa-Case6-lemm}\,(iv),\,\eqref{(T0o(eE)1}}}{=}\frac{\frac1{A_2}+9A_3^{10}X_2^{90}}{10A_3^{10}X_2^{91}\Big(\frac{49}{45}-\frac{4}{45A_2A_3^{30}}\Big)}+O(\d)^2.
\!\!\!\!\!\!\!\!}
Now assume $\widetilde\xX_1\le\d_2$, i.e., $\aA_3\ge\ell_2$ by \eqref{LetNSoOP----1}\,(i). Then
\eqref{LetNSoOP----1++++++}\,(ii) with \eqref{A1-A2-cond}\,(b) shows that
$10-\frac{9X_2^{90}}{A_3^{20}}=\frac1{A_2A_3^{30}}+O(\d)^2=O(\d_2)^{30}$ (recalling that $\ell_1=\d_1^{-1}\ll\ell_2=\d_2^{-1}\ll\ell=\d^{-1}$). In particular,
\equa{msdmseme444}{\dis{\rm(i)\ }
\frac{X_2^{90}}{A_3^{20}}=\frac{10}{9}+O(\d_2)^1,\ \ \ \
{\rm(ii)\ }\frac15+\frac{4X_2^{90}}{5A_3^{20}}\stackrel{{}^{\sc\rm\eqref{msdmseme444}\,(i)}}{=}\frac{49}{45}+O(\d_2)^1,
\ \ \ \
{\rm(iii)\ }\xX_2\stackrel{{}^{\sc\rm\eqref{msdmseme444}\,(i)}}{>}\aA_3^{\frac29}.
}
Using this in \eqref{LetNSoOP----1++++++}\,(i), we obtain that $\aA_1>\xX_2>\aA_3^{\frac29}\ge\ell_2^{\frac29}\gg\ell_1$,
a contradiction with \eqref{A1-A2-cond}\,(a).
This proves the first inequality of \eqref{ImMpP}\,(ii), the second inequality follows from \eqref{ImMpP}\,(i) and the fact that
$\tilde\xX_1=\aA_3^{-1}$.
Next, assume the first inequality of \eqref{ImMpP}\,(iii) does not hold. Then $\xX_2=1-\d$ by \eqref{C+ToSayas+1}\,(d).
Therefore up to $O(\d)^2$, we have
\begin{eqnarray}
\!\!\!\!\!\!\!\!\!\!\!\!\!\!\!\!\!\!\!\!\!\!\!\!\!\!\!\!\!\!&&
\label{Up-To-+1}
1\ \ \ \,\stackrel{{}^{\sc\rm\eqref{LetNSoOP----1++++++}\,(i)}}{\le}
\ \ \
\frac{\xX_2^{91}}{\aA_1\aA_3^{20}}\Big(\frac15+\frac{4\xX_2^{90}}{5\aA_3^{20}}\Big)=
(1-\d)^{91}\aA_1^{-1}\aA_3^{-20}\Big(\frac15+\frac45(1-\d)^{90}\aA_3^{-20}\Big)
\nonumber\\
\!\!\!\!\!\!\!\!\!\!\!\!\!\!\!\!\!\!\!\!\!\!\!\!\!\!\!\!\!\!&&
\stackrel{{}^{\sc\rm\eqref{A1-A2-cond}\,(a),\,\eqref{ImMpP}\,(i)}}{\le}
(1-\d)^{91}\Big(\frac15+\frac45(1-\d)^{90}\Big)
=1-163\d\mbox{ \ [up to $O(\d)^2$]},
\end{eqnarray}
which is a contradiction.
Now assume the last inequality of \eqref{ImMpP}\,(iii) does not hold. Then $\xX_2=(1+\d)\aA_1^2$ by \eqref{C+ToSayas+1}\,(d).
Therefore up to $O(\d)^2$, we have, where the second inequality follows by
noting from \eqref{cont-B2},\,\eqref{A1-A2-cond},\,\eqref{ImMpP}\,(i) that
$\frac{\aA_1}{\aA_2\aA_3^{10}\aA_1^{182}}=\aA_2^{-1}\aA_3^{-10}\aA_1^{-181}\le1$,
$\frac{\aA_1}{\aA_1^2}=\aA_1^{-1}\le1$ and $\frac1{\aA_2\aA_3^{30}}\le1$,
\begin{eqnarray}
\!\!\!\!\!\!\!\!\!\!\!\!\!\!\!\!\!\!\!\!\!\!\!\!\!\!\!\!\!\!&&
\label{Up-To-}
1\ \ \ \ \ \ \,\stackrel{{}^{\sc\rm\eqref{MEMEME-A-inv}}}{\le}
\ \ \ \ \ \ \
\frac{\frac{\aA_1}{\aA_2}+9\aA_1\aA_3^{10}\xX_2^{90}}{10\aA_3^{10}\xX_2^{91}\Big(\frac{49}{45}-\frac{4}{45\aA_2\aA_3^{30}}\Big)}
=\frac{\frac{\aA_1}{10\aA_2\aA_3^{10}(1+\d)^{91}\aA_1^{182}}+\frac{9\aA_1}{10(1+\d)\aA_1^2}}{\frac{49}{45}-\frac4{45\aA_2\aA_3^{30}}}
\nonumber\\
\!\!\!\!\!\!\!\!\!\!\!\!\!\!\!\!\!\!\!\!\!\!\!\!\!\!\!\!\!\!&&
\stackrel{{}^{\sc\rm\eqref{cont-B2},\,\eqref{A1-A2-cond},\,\eqref{ImMpP}\,(i)}}{\le}
\frac{\frac{1}{10(1+\d)^{91}}+\frac{9}{10(1+\d)}}{\frac{49}{45}-\frac4{45}}
=1-10\d\mbox{ \ [up to $O(\d)^2$]},
\end{eqnarray}
which is again  a contradiction.
Finally assume $|10-9\widetilde X_1^{20}X_2^{90}|\le\d_2^{100}$. Then by \eqref{ImMpP}\,(ii) and the fact that
$\aA_3=\tilde\xX_1^{-1}$, we see from \eqref{LetNSoOP----1++++++}\,(ii) that $\aA_2\gg\ell_1^{\frac{48}{11}}$, a
contradiction with \eqref{A1-A2-cond}\,(b).
\NOUSE{
Now assume $\big|\ell_0+1-\frac{\ell_0}{X_2^{\ell_0}}\big|\le\d_2^{90}$. Then
by \eqref{Case6-lemm},\,\eqref{LetNSoOP----1}\,(i),\,\eqref{C+ToSayas+1}\,(d), we see that $\aA_2>\ell_2^{80}$
(recalling that $\ell_2=\d_2^{-1}$),
and the second term inside the brackets of $A_1$ in
\eqref{LetNSoOP----1}\,(ii) is an $O(\d_2)^1$ element, and
$\aA_2^6 \xX_2^{-6}\gg\ell_1$, which implies that $\aA_1\gg\ell_1$, a contradiction with \eqref{A1-A2-cond}\,(a). This proves \eqref{ImMpP}\,(ii).
}
} \NOUSE{that \eqref{Amsmene} holds for $a=\tilde\xX_1,\xX_2,\zZ,\aA_1,\aA_2$
[by noting from \eqref{MSmde33333}
that $\nn_0=\eE_0^{-1}\gg\nn=\eE^{-1}$]. It also holds for $a=\aA_3=\a_1\zZ^{16}$ by \eqref{LetNSoOP----1}\,(i),\,\eqref{C+ToSayas+1}\,(f).
Further, we  see from \eqref{LetNSoOP----1}\,(iii) [or the equality in \eqref{x1-small}\,(i)$\ssc\,$] and \eqref{C+ToSayas+1}\,(d),\,(e) that  \eqref{Amsmene} holds for $a=\xX_1$.
\NOUSE{
Further, by \eqref{LetNSoOP----1},\,\eqref{ImMpP}\,(1), we can easily see that $\zZ\ge\eE_0^{10}$, and
\equa{zzzz===}{\dis
\zZ\stackrel{{}^{\sc\rm\eqref{denote-t2}}}{=}|Z|
\stackrel{{}^{\sc\rm\eqref{TaKa},\,\eqref{SimMMSMS}}}
{=}\g_{\kk,\kk}^{-1}|x_2+y_2|\stackrel{{}^{\sc\rm\eqref{ggggg},\,\eqref{x1-small}\,(iii),\,\eqref{y1-y2-isSmall}}}{\le}\nn_0^{10}.
}%
}%
This proves \eqref{Amsmene}.
}%
\NOUSE{
Now assume \eqref{ImMpP}\,(ii) does not hold. Then by \eqref{C+ToSayas+1}\,(d), we have $\tilde\xX_1=1-\d$.
Then \eqref{LetNSoOP----1}\,(ii) gives the following, up to $O(\d)^2$ [using \eqref{(T0o(eE)1}$\ssc\,$],
\begin{eqnarray}
\label{MSMn4n4}
&&\!\!\!\!\!\!\!\!\!\!\!\!\!\!\!\!\!\!\!\!\!\!\!\!\!\!
1\stackrel{{}^{\sc\rm\eqref{LetNSoOP----1}\,(ii)}}{\ge}
\tilde\xX_1^{-1}\aA_2^{-1}\aA_3^3(1+\d_0-\d_0\tilde\xX_1^5)
\nonumber\\
&&\!\!\!\!\!\!\!\!\!\!\!\!\!\!\!\!\!\!\!\!\!\!\!\!\!\!
\ \ \  \ \ \ =\ \ (1-\d)^{-1}\aA_2^{-1}\aA_3^3\Big(1+\d_0-\d_0(1-\d)^5\Big)
\stackrel{{}^{\sc\rm\eqref{cont-B2},\,\eqref{ImMpP}\,(i)}}{\ge}(1+\d)\aA_1^{(\frac1{20}-\frac {\d_0}2)-\frac{3}{40}}(1+5\d_0\d)
\!\!\!\!\!\!\!\nonumber\\
&&\!\!\!\!\!\!\!\!\!\!\!\!\!\!\!\!\!\!\!\!\!\!\!\!\!\!
\ \
\stackrel{{}^{\sc\rm\eqref{A1-A2-cond}\,(a)}}{\ge}1+(1+5\d_0)\d\mbox{ \ \ [up to $O(\d)^2$]},
\end{eqnarray}
which is a contradiction. This proves \eqref{ImMpP}\,(ii).
\hfill$\Box$\vskip7pt
}
}
\NOUSE{In the following, we will need to use the following formulas, which are obtained from \eqref{B1-ep}
[noting from \eqref{MSmde33333},\,\eqref{ll-dd},\,\eqref{Akk-bkk} that $\ell\dD=O(\dD)^1$],
\equa{a2a1a3}{\dis a_2a_1^{-1}=\ell\Big(2\dD+\frac{2130\dD^2}{101}\Big)+O(\dD)^3,\ \ \ \ \ \ \ \ \ \ \ \ a_3a_1^{-1}=
\frac{\dD}{101}
 + \frac{1931 \dD^2}{10201}+O(\dD)^3.}
%
%
Since $\aA_1$ is invertible [cf.~\eqref{denote-t2},\,\eqref{Amsmene}$\ssc\,$], we obtain from \eqref{cont-B2},\,\eqref{ImMpP} 
the following
[recalling that $0<\eE_1=\nn_1^{-1}\ll\dD=\lL^{-1}$, cf.~\eqref{MSmde33333} and \eqref{ll-dd}$\ssc\,$],
\equa{ImMpP-(1)}{\aA_3\stackrel{{}^{\sc\rm\eqref{ImMpP}\,(1)}}{>}\aA_2^{a_3a_2^{-1}-\eE_1^{3}}\stackrel{{}^{\sc\rm\eqref{cont-B2}}}{=}\aA_1^{a_3a_1^{-1}+O(\eE_1)^3}
\stackrel{{}^{\sc\rm\eqref{a2a1a3}}}{=}
\aA_1^{\frac{\dD}{101}
 + \frac{1931 \dD^2}{10201}+O(\dD)^3}.}
}
\NOUSE{In the following we apply formula \eqref{(T0o(eE)1} frequently. First up to $O(\d)^2$, we have
\equa{X2===}{\dis{\rm(i)\ }\a_1\stackrel{{}^{\sc\rm\eqref{LetNSoOP----1}\,(ii)}}{=}1\mbox{ \ [up to $O(\d)^2$]},\ \ \
{\rm(ii)\ }\xX_2\stackrel{{}^{\sc\rm\eqref{equa-Case6-lemm}\,(iv)}}{=}\zZ\stackrel{{}^{\sc\rm\eqref{LetNSoOP----1}\,(i)}}{=}
\aA_3\tilde\xX_1^{-4}\mbox{ \ [up to $O(\d)^2$]}.
}
}\NOUSE{Recalling notations in \eqref{denote-t2}, and
noting from \eqref{LetNSoOP----1}\,(i) that $\a_1=1+O(\ep)^1=1+O(\eE_1)^2$ [thus the factor $\a_1$ can be omitted in \eqref{C1B1}$\ssc\,$] and so $\zZ=\aA_3^{\frac1{16}}+O(\eE_1)^2$,
we obtain (noting that $B_1,C_1$ are invertible as $A_1$ is invertible in $\ol V$),
}%
\NOUSE
{%
Note that we have already proved in Lemma \ref{lemm-condition-XZ} that $\aA_1$ is invertible,
thus $B_1,C_1$
are invertible. We have,
 \begin{eqnarray}
\label{C1B1}
&\!\!\!\!\!\!\!\!\!\!\!\!\!\!\!\!\!\!\!\!\!\!\!\!\!\!\!\!\!\!&
{\rm(i)\ }|B_1|^{-1}\stackrel{{}^{\sc\rm\eqref{LetNSoOP----1}\,(iii),\,\eqref{denote-t2}}}{=}
\aA_1^{-1} \aA_2^{-12}\tilde\xX_1^{18\ell_0^2+11} \xX_2^{4}
%
,\!\!\!\!\!\!\!\!\!\!\!\!\!\!\!\!\!\!\!\!\!\!
\nonumber\\
&\!\!\!\!\!\!\!\!\!\!\!\!\!\!\!\!\!\!\!\!\!\!\!\!\!\!\!\!\!\!&
{\rm(ii)\ }
|C_1|^{-1}\stackrel{{}^{\sc\rm\eqref{LetNSoOP----1}\,(iv),\,\eqref{denote-t2}}}{=}
\aA_1^{-2} \aA_2^{-6}\tilde \xX_1^{18\ell_0^2+11} \xX_2^{-2}
\NOUSE{
\nonumber\\
&\!\!\!\!\!\!\!\!\!\!\!\!\!\!\!\!\!\!\!\!\!\!\!\!\!\!\!\!\!\!&
{\rm(iv)\ }
|B_2|^{-1}\stackrel{{}^{\sc\rm\eqref{LetNSoOP----1}\,(vi),\,\eqref{denote-t2}}}{=}
\aA_2^{-1}\aA_3^{-100}\tilde\xX_1^{11}\xX_2^{-82}
=\aA_2^{-1}\aA_3^{-182}\tilde\xX_1^{93}
\mbox{ \ [up to $O(\d)^2$]},
\nonumber\\
&\!\!\!\!\!\!\!\!\!\!\!\!\!\!\!\!\!\!\!\!\!\!\!\!\!\!\!\!\!\!&
{\rm(v)\ }
|C_2|^{-1}\stackrel{{}^{\sc\rm\eqref{LetNSoOP----1}\,(vii),\,\eqref{denote-t2}}}{=}
\aA_2^{-2}\aA_3^{-100}\tilde\xX_1^{22}\xX_2^{-254}
=\aA_2^{-2}\aA_3^{-354}\tilde\xX_1^{276}\mbox{ \ [up to $O(\d)^2$]}
}
.
\end{eqnarray}
%
\NOUSE
{We need to compute \eqref{CCCC===} below. First observe from \eqref{C1B1} that
 the factor $\tilde\xX_1$ is cancelled in \eqref{CCCC===}. Then we compute the powers of $\aA_i,\,i=1,2,3$ [which is done in the second equality of \eqref{CCCC===}$\ssc\,$], and use \eqref{cont-B2} and \eqref{ImMpP-(1)} (noting that in order to use this formula, it is required that the power of $\aA_3$ is non-positive) to write the result as a power of $\aA_1$ [which is done in the second inequality of \eqref{CCCC===}$\ssc\,$] and obtain finally that the power of $\aA_1$ is negative. Namely,
 Thus up to  $O(\eE_1)^2$,
 we 
 obtain,
\begin{eqnarray}
\label{CCCC===}
\!\!\!\!\!\!\!\!\!\!\!\!\!\!\!\!\!\!\!\!\!\!\!\!\!\!\!\!\!\!&&
|B_1|^{-\frac12}\cdot|C_1|^{-1}
\stackrel{{}^{\sc\rm\eqref{C1B1}\,(i),\,(ii)}}
{=}
\aA_1^{-\frac52}\aA_3^{-12}\tilde\xX_1^{27\ell_0^2+\frac{33}{2}}
\stackrel{{}^{\sc\rm\eqref{ImMpP}\,(i),\,\eqref{cont-B2}}}{\le}
\tilde\xX_1^{-\frac52(6\ell_0^2+\frac{34}{5})-12\ell_0^2+27\ell_0^2+\frac{33}{2}}
\nonumber\\
\!\!\!\!\!\!\!\!\!\!\!\!\!\!\!\!\!\!\!\!\!\!\!\!\!\!\!\!\!\!&&
\ \ \ \ \ \ \ \ \ \ \ \ \ \ \ \ \ \ \ \ \ \ \ \ \ \
=\ \ \ \ \tilde\xX_1^{-\frac{1}{2}}\stackrel{{}^{\sc\rm\eqref{A1-A2-cond}}}{\le}1
\mbox{ \ [up to $O(\d)^2$]}.
\end{eqnarray}
}}%
\NOUSE{
We have
\equa{MSnene}{\dis
|B_2C_2^2|^{-1}\stackrel{{}^{\sc\rm\eqref{LetNSoOP----1}\,(iv),\,(v)}}{=}
\aA_2^{-5}\stackrel{{}^{\sc\rm\eqref{A1-A2-cond}\,(b)}}{\le}1.
}
}%
\NOUSE
{\noindent
We sincerely thank Professor Claudio Procesi, who observes that formula \eqref{Re-Writtt}
with \eqref{LetNSoOP----1-redefine++}\,(v)  and
condition \eqref{C+LetNSoOP}\,(iii)%
~can be used to give an elegant proof of the following lemma.
}\NOUSE
{%
Now assume equality occurs in the first inequality of \eqref{C+ToSayas+1}\,(c), i.e.,
\equa{MSMSMS000000}{
\tilde\xX_1=\eta\aA_2^{\lL^5(1-\dD)\b_2^{-1}}
\mbox{ \ with \ }\eta=1-\eE_1.}
First assume $\aA_2\ge\eta^{\lL\b_2}$. Then up to $O(\eE_1)^2$, we obtain
,
\begin{eqnarray*}
&\!\!\!\!\!\!\!\!\!\!\!\!\!\!\!\!\!\!\!\!\!\!\!\!\!\!\!\!\!\!&
\xX_2\stackrel{{}^{\sc\rm\eqref{LetNSoOP----1}\,(ii)}}{=}
\aA_2^{-\beta_2^{-1}}\tilde\xX_1^{\dD^5+O(\dD)^N}
\stackrel{{}^{\sc\rm\eqref{B1-ep}\,(ii),\,(iii),\,\eqref{MSMSMS000000}}}{=}\ \
\eta^{\dD^5+O(\dD)^N}\aA_2^{(-\dD+O(\dD)^N)\b_1^{-1}}
>
\eta^{\lL^{10}(1+O(\dD)^1)},\!\!\!\!\!\!\!\!\!\!\!\!\!\!\!
\end{eqnarray*}
which is  bigger than $1$, $\xX_1$ [by \eqref{x1-small}\,(i)$\ssc\,$] and $
\a_1^{\lL^5(2 + \dD)}\tilde\xX_1^{2\lL^4(2-\dD+\dD^2)}$.
From this
and 
\eqref{LetNSoOP----1}\,(i), \eqref{ImMpP}\,(1) (and noting that $a_3a_2^{-1}>0$), we obtain that $\zZ>\xX_2$, thus also $\zZ>\xX_1$.
Then we obtain a contradiction as in \eqref{Conndndm}. This proves that the following holds for $\aA_2$
[thus also holds for $\aA_1$ by \eqref{cont-B2} and holds for $\tilde\xX_1,\,\eta$ by \eqref{MSMSMS000000}$\ssc\,$],
\equa{OSOSOSO}{
a=1+O(\eE_1)^1\mbox{ \ for \ }a\in T_0:=\{\aA_1,\aA_2,\aA_3,\tilde\xX_1,\xX_2,\eta\}.
}
To see \eqref{OSOSOSO} holds for $\xX_2$, first by
\eqref{LetNSoOP----1}\,(ii),\,\eqref{ImMpP}\,(1), we have $\xX_2\ge1+O(\eE_1)^1$. If $\xX_2>
\a_1^{\lL^5(2 + \dD)}\eta^{2\lL^4  (2 - \dD + \dD^2)} $,
then
\eqref{MSMSMS000000} with \eqref{LetNSoOP----1}\,(i),\,\eqref{ImMpP}\,(1)
shows that $\zZ>\xX_2$ and again we obtain a contradiction as in
\eqref{Conndndm}. Hence \eqref{OSOSOSO} holds for $\xX_2$, thus also holds for $\aA_3$ by \eqref{LetNSoOP----1}\,(ii).
This proves \eqref{OSOSOSO}.
%
Now up to $O(\eE_1)^2$, we have [noting that $\nn_1=\eE_1^{-1}\gg\lL$ by \eqref{MSmde33333},\,\eqref{ll-dd}$\ssc\,$],
\begin{eqnarray}
\label{WeHam000}
\!\!\!\!\!\!\!\!\!\!\!\!\!\!\!\!\!\!\!\!\!\!\!\!\!\!\!&&
1\ \
\stackrel{{}^{\sc\rm\eqref{LetNSoOP----1}\,(iii)}}{\le}
\ \
\frac{\tilde \xX_1^{8\lL^6(2-\dD)}}{\aA_1\aA_3^{\lL(2+\dD)}\Big(
 \lL+1 - \lL \aA_2^{\lL^4}\aA_3^2\tilde \xX_1^{-16\lL^5}\Big)}
\nonumber\\
\!\!\!\!\!\!\!\!\!\!\!\!\!\!\!\!\!\!\!\!\!\!\!\!\!\!\!&&
\stackrel{{}^{\sc\rm\eqref{cont-B2},\,\eqref{MSMSMS000000}}}{=}
\gamma_1:=
\frac{\eta^{8\lL^6(2-\dD)}
\aA_2^{-\frac{8\lL^6(2-\dD)}{5}
-a_1a_2^{-1}}
}{\aA_3^{\lL(2+\dD)}\Big(
 \lL+1 - \lL \aA_3^2\eta^{-16\lL^5}
 \aA_2^{\lL^4+\frac{16\lL^5}{5}}
 \Big)}
\nonumber\\
\!\!\!\!\!\!\!\!\!\!\!\!\!\!\!\!\!\!\!\!\!\!\!\!\!\!\!&&
\phantom{1}
\stackrel{{}^{\sc\rm\eqref{ImMpP}\,(1)}}{<}
\frac{\eta^{8\lL^6(2-\dD)}
\aA_2^{-\lL(2+\dD)a_3a_2^{-1}-\frac{8\lL^6(2-\dD)}{5}
-a_1a_2^{-1}}
}{
 \lL+1 - \lL \eta^{-16\lL^5}
 \aA_2^{2a_3a_2^{-1}+
 \lL^4+\frac{16\lL^5}{5}}
 }
\nonumber\\
\!\!\!\!\!\!\!\!\!\!\!\!\!\!\!\!\!\!\!\!\!\!\!\!\!\!\!&&
\phantom{1\ \ \ \ }
=\gamma_2:=
\frac{\eta^{8\lL^6(2-\dD)}
\aA_2^{-\frac{\lL^6(1+7\dD+O(\dD)^2)}{5}
}}{
 \lL+1 - \lL \eta^{-16\lL^5}
 \aA_2^{\frac{\lL^5}{5}(16+5\dD+O(\dD)^2)}
 }
 \nonumber\\
\!\!\!\!\!\!\!\!\!\!\!\!\!\!\!\!\!\!\!\!\!\!\!\!\!\!\!&&
\phantom{1\ \ \ \ }
\le\frac{\eta^{8\lL^6(2-\dD)}
}{
 \lL+1 - \lL \eta^{-16\lL^5}
 }=1-8\lL^5\Big(1+O(\dD)^1\Big)\eE_1<1,
\end{eqnarray}
a contradiction, where the second inequality is obtained by noting that $\gamma_1$ is locally a decreasing function on
$\aA_3$ [as $\frac{\ptl\gamma_1}{\ptl\aA_3}|_{(\eta,\aA_2,\aA_3)=(1,1,1)}=-1\ssc\,$] and we have  \eqref{OSOSOSO}, thus we can replace $\aA_3$ by $\aA_2^{a_3a_2^{-1}}$ by \eqref{ImMpP}\,(1); while the last inequality is obtained by noting that
$\gamma_2$ is locally a decreasing function on $\aA_2$
[as $\frac{\ptl\gamma_1}{\ptl\aA_2}|_{(\eta,\aA_2)=(1,1)}=-\frac{2\lL^5}{5}(1+O(\dD)^1)<0\ssc\,$]
and $\aA_2\ge1$. This proves \eqref{ImMpP}\,(2).
Assume $\tilde\xX_1\le\eE_1$ [then $\xX_1\le\eE_1^{\eE_1}$ by \eqref{x1-small}\,(i)$\ssc\,$].
Then \eqref{LetNSoOP----1}\,(ii) with \eqref{A1-A2-cond}\,(b),\,\eqref{ImMpP}\,(1)
implies that $\xX_2\gg\nn_1$ [noting from \eqref{MSmde33333} that $\nn_1\gg\nn$], and then
\eqref{LetNSoOP----1}\,(i) shows that $\zZ\gg\xX_2$. Thus also $\zZ\gg1$ and $\zZ\gg\xX_1$. We obtain
\equa{Conndndm}{\dis
\g_{|x_1|,|x_2|}
=\g_{\xX_1\kk,\xX_2\kk}\stackrel{{}^{\sc\rm\eqref{wePPPP1+}}}{<}
\g_{\zZ\kk,\xX_2\kk}
\stackrel{{}^{\sc\rm\eqref{wePPPP1}}}{\le}\g_{\zZ\kk,\zZ\kk}
\stackrel{{}^{\sc\rm\eqref{AssG}}}{\le}\zZ\g_{\kk,\kk}
\stackrel{{}^{\sc\rm\eqref{TaKa},\,\eqref{SimMMSMS}}}{=}|x_2+y_2|\stackrel{{}^{\sc\rm\eqref{Ak=1}}}{\le}\g_{|x_1|,|x_2|},
}
a contradiction. This proves the first inequality of \eqref{ImMpP}\,(2)
[similarly, we can prove the last inequality of \eqref{ImMpP}\,(3)$\ssc\,$].
Assume $\xX_2\le\eE_1$. Then \eqref{ImMpP}\,(2) with \eqref{ImMpP}\,(1),\,\eqref{LetNSoOP----1}\,(ii) shows that $\aA_2<1$, a contradiction with \eqref{A1-A2-cond}\,(b). This completes the proof of \eqref{ImMpP}\,(3).
By \eqref{ImMpP}\,\,(1),\,(2),\,\eqref{SimMMSMS}\,(a),\,(b), we see that $x_1,x_2,A_3\ne0$, then by \eqref{LetNSoOP----1}\,(i),\,\eqref{SimMMSMS}\,(c), we see that $x_2+y_2\ne0$. This proves \eqref{-EiathA0}.
}\NOUSE{
Note from \eqref{LetNSoOP----1}\,(ii),\,(iii) that we have the following formula,
\equa{A2===A111}{\dis
A_2^{-\b_1\b_2^{-1}}
=\frac{
X_1^{\b_1(1-\dD-\l_2\dD^{20})}X_2^{\b_1(\dD+\l_3\dD^{20})}
\Big(
\frac{A_1}{A_3^{\b_1\b_3^{-1}}X_1^{\b_1}}
+\lL\frac{A_2^{\b_1(1-\frac{\dD^6}{2}-\l_1\dD^{20})}A_3^{2\b_1\b_3^{-1}}}{X_1^{\b_1}}
\Big)
}{(\lL+1)A_3^{2\b_1\b_3^{-1}}}.
}
This with \eqref{cont-B2},\,\eqref{ImMpP}\,(1),\,\eqref{MSMSMS000000} gives, up to $O(\nn_1^{-1})^2$,
\begin{eqnarray}
\lL+1\le
=\eta^{-\b_1(\dD+\l_2\dD^{20})}\xX_2^{\b_1(\dD+\l_3\dD^{20})}
\Big(\aA_2^{a_1a_2^{-1}+3\b_1\b_3^{-1}a_3a_2^{-1}+\b_1\b_2^{-1}+\frac{\lL^3(1-\dD^5)}{3\b_2}(\dD+\l_2\dD^{20})}
+\lL\aA_2^{\b_1(1-\frac{\dD^6}{2}-\l_1\dD^{20})+\frac{\lL^3(1-\dD^5)}{3\b_2}(\dD+\l_2\dD^{20})}
\Big)
\end{eqnarray}
By \eqref{A1-A2-cond} and the last equality of
\eqref{LetNSoOP----2}\,(ii), we immediately obtain \eqref{ImMpP}\,(2)
[if $|X_1|\le\d_0$ then $|A_2|\ll1$; if $|X_1|\ge\ell_0^{2}$ then $|A_2|\gg\ell_1^{\ell_0^3+\d_0}$]. From this and the second equality of
\eqref{LetNSoOP----2}\,(ii), we obtain  \eqref{ImMpP}\,(3).
Then \eqref{ImMpP}\,(4) is obtained from \eqref{LetNSoOP----2}\,(i),\,\eqref{A1-A2-cond} and  \eqref{ImMpP}\,(2).
}\NOUSE{
For convenience, we simply denote, for all possible $i$,
\equa{tildexz}{\dis\tilde\xX_1=|\widetilde X_1|,\ \ \ \
\xX_i=|X_i|,
\ \ \ \ \zZ=|Z|,\ \ \ \ \aA_i=|A_i|.}
 To prove  \eqref{ImMpP}\,(2),  assume  equality occurs in
any  inequality of  \eqref{C+ToSayas+1}\,(c), i.e., $\xX_1=
\l\aA_1^{\frac12-\d_0^5}$ with $\l=1\pm\d$.
We need to use the following formula,  up to $O(\d)^2$,
\begin{eqnarray}
\label{X-2==aa}
&&\!\!\!\!\!\!\!\!\!\!\!\!\!\!\!\!\!\!\!\!\!\!\!\!\!\!\!\!\!\!\!\!\!\!\!\!\!\!\!\!\!\!
\xX_1^{-1}\xX_2^{-\ell_0^2-1}
\stackrel{{}^{\sc\rm\eqref{LetNSoOP----2}\,(iii)}}{=}
\aA_3^{-\ell_0^2-1}\xX_1^{2\ell_0^3+2\ell_0-1}
\stackrel{{}^{\sc\rm\eqref{ImMpP}\,(1)}}{\le}
\l^{2\ell_0^3+2\ell_0-1}
\aA_1^{-(\ell_0^2+1)a_3a_1^{-1}+(2\ell_0^3+2\ell_0-1)(\frac12-\d_0^5)}\!\!\!\!\!\!\!\!\!\!\!\!\!\!\!
\nonumber\\
&&\!\!\!\!\!\!\!\!\!\!\!\!\!\!\!\!\!\!\!\!\!\!\!\!\!\!\!\!\!\!\!\!\!\!\!\!\!\!\!\!\!\!
\ \ \ \ \ \ \ \ \ \ \ \ \ \ \ \ \ \ \
=\ \ \ \
\l^{2\ell_0^3+2\ell_0-1}\aA_1^{-\frac12+
O(\d_0)^1}.\!\!\!\!\!
\end{eqnarray}
Dividing the formula \eqref{LetNSoOP----2}\,(ii) by $\aA_2$ and using $\aA_2=\aA_1^{a_2a_1^{-1}+O(\d)^2}$, up to $O(\d)^2$, 
%
we can compute, 
\begin{eqnarray}
&\!\!\!\!\!\!\!\!\!\!\!\!\!\!\!\!\!\!\!\!\!\!\!\!\!\!\!\!\!\!\!\!&
\label{A2-derf-}
1
\stackrel{{}^{\sc\rm\eqref{LetNSoOP----2}\,(ii)}}{\le}
\l^{-1}\aA_1^{-a_2a_1-\frac12+\d_0^5}\frac{\frac34+\frac14\l^{-2}\aA_1^{2\d_0^5}}
{2-\l^{-2}\aA_1^{2\d_0^5}}
=\frac{\l^{-1}\aA_1^{\frac{7\d_0^4}{16}+O(\d_0)^5}\Big(
\frac34+\frac14\l^{-2}\aA_1^{2\d_0^5}\Big)}
{2-\l^{-2}\aA_1^{2\d_0^5}}
.
\end{eqnarray}
When $\l=1+\d$, the right-hand side of \eqref{A2-derf-} is $<(1+\d)^{-1}$ (as $\aA_1\le1$), which is a contradiction; in particular, this
proves the last inequality of \eqref{ImMpP}\,(2).
Now assume $\l=1-\d$. Then since $\aA_1\le1$, \eqref{A2-derf-} gives us the following important fact,
\equa{FactaA1}{\dis\aA_1=1+O(\d)^1.}
Now dividing the formula \eqref{LetNSoOP----2}\,(i) by $\aA_1$, using
\eqref{X-2==aa}, we obtain
\begin{eqnarray}
&\!\!\!\!\!\!\!\!\!\!\!\!\!\!\!\!\!\!\!\!\!\!\!\!\!\!\!&
\label{aa1+A2-derf-}
1\le\gamma_1:=\l^2\aA_1^{-2\d_0^5}\Big(1-\d_0^4+\d_0^4
\l^{2\ell_0^3+2\ell_0-1}\aA_1^{-\frac12+
O(\d_0)^1}\Big)
\NOUSE{\!\!\!\!\!\!\!\!\!\!\!\!\!\!\!\!\!\!\!\!\!
\nonumber\\
&\!\!\!\!\!\!\!\!\!\!\!\!\!\!\!\!\!\!\!\!\!\!\!\!\!\!\!&
\phantom{1}\le
\l^{-1}\aA_1^{\frac{\d_0^2}{4}+O(\d_0)^3}
\frac{1-\frac{\d_0^4}{2}+\frac{\d_0^4}2\l^{2\ell_0^3-1}\aA_1^{-\ell_0(1+O(\d_0)^1)}}{1+\frac{\d_0^3}2-\d_0^4-
\Big(\frac{\d_0^3}{2}-\d_0^4\Big)
\l^{2\ell_0^3-1}\aA_1^{-\ell_0(1+O(\d_0)^1)}}
}
.
\end{eqnarray}
One can easily observe that $\gamma_1(\l,\aA_1)$ is a local function on $\l,\aA_1$ (which are $1+O(\d)^1$ elements) such that
$\frac{\ptl\gamma_1}{\ptl\aA_1}=-\frac1{}$
\equa{Such-lll}{\dis
{\rm(i)\ }\frac{\ptl\gamma_1}{\ptl\aA_1}|_{(\l,\aA_1)=(1,1)}=-\frac{\d_0^2}{4}<0,\ \ \ \ {\rm(ii)\ }\frac{\ptl\gamma_1}{\ptl\l}|_{(\l,\aA_1)=(1,1)}=-\d_0+O(\d_0)^2<0.}
 Since $0<\d\ll\d_0$ and $\aA_1\ge1$, $\l=1+\d$,
we obtain from \eqref{aa1+A2-derf-} the following, up to $O(\d)^2$,
 \equa{Gamamamama}{\dis
1\le\gamma_1(\l,\aA_1)\le\gamma_1(\l,1)=1+\Big(-\d_0+O(\d_0)^2\Big)\d,}
which is again a contradiction. This proves \eqref{ImMpP}\,(2).
From \eqref{ImMpP}\,(2) and \eqref{A1-A2-cond}, we obtain \eqref{ImMpP}\,(3).
The first inequality of \eqref{ImMpP}\,(4) follows from \eqref{ImMpP}\,(1),\,(2) and \eqref{LetNSoOP----2}\,(iii).
Assume $|X_2|\ge\ell_2$. Then we obtain from \eqref{LetNSoOP----2}\,(i),\,(ii),\,\eqref{ImMpP}\,(2) that $A_1=(1-\d_0^4)X_1^2+O(\d_2)^1$, $A_2=\frac{1-\frac{\d_0^4}{2}}{X_1(1+\frac{\d_0^3}{2}-\d_0^4)}+O(\d_2)^1$, and
[noting from  \eqref{B1-ep} that $a_1+2a_2=\d_0^2+O(\d_0)^3>0$],
\equa{fmfmfm}{\dis
1\stackrel{{}^{\sc\rm\eqref{cont-B2}}}{\le}|A_1A_2^2|
=(1-\d_0^4)\Big(\frac{1-\frac{\d_0^4}{2}}{1+\frac{\d_0^3}{2}-\d_0^4}\Big)^2+O(\d_2)^1<1,
}
a contradiction. This proves \eqref{ImMpP}\,(4).
Finally assume $|1+\frac{\d_0^3}2-\d_0^4-\big(\frac{\d_0^3}{2}-\d_0^4\big)\frac{1}{X_1X_2^{\ell_0^2}}|=\d_2$. Then
\equa{Msdnene}{\dis
\frac{1}{X_1X_2^{\ell_0^2}}=\Big(1+\frac{\d_0^3}2-\d_0^4\Big)
\Big(\frac{\d_0^3}{2}-\d_0^4\Big)^{-1}+O(\d_2)^1=2\ell_0^3\Big(1+O(\d_0)^1\Big),
} and so
$|1-\frac{\d_0^4}{2}+\frac{\d_0^4}{2X_1X_2^{\ell_0^2}}|=
1+O(\d_0)^1$. From this and \eqref{LetNSoOP----2}\,(iii),\,\eqref{ImMpP}\,(2), we obtain that $|A_2|\gg1$, a contradiction with
\eqref{A1-A2-cond}\,(ii).
}\NOUSE{
Since $\aA_1\ge1$, the above proves
\equa{aa1-is-}{\dis \aA_1=1+O(\d)^1.}
Note that we can regard \eqref{LetNSoOP----2}\,(iii) as a linear equation on $A_1$ and solve it to obtain
\equa{New-A1}{\dis
A_1=\frac{X_1^{3\ell_0^2}\Big(7+(13+28\d_0^5)\frac{X_2^{5}}{A_2X_1^{3\ell_0^2}}\Big)}
{A_3^4X_2^2\Big(35-(15-28\d_0^5)\frac{X_2^5}{A_2X_1^{3\ell_0^2}}\Big)}+O(\d)^2.
}
For convenience, we denote $\l=(1+\d)^{\ell_0^2}$.
Using \eqref{New-A1},\,\eqref{LetNSoOP----2}, we obtain, up to $O(\d)^2$
[note that we have successively made the following replacements in the first inequality:
$\xX_2\to\aA_3^{-1}$ by \eqref{LetNSoOP----2}\,(ii); $\aA_3^{-1}\to\aA_1^{\frac{1}{9}}$ by \eqref{ImMpP}\,(1) (observing that only negative powers of $\aA_3$ appear); $\aA_2\to\aA_1^{-\frac13}$ by \eqref{cont-B2}; $\xX_1^{\ell_0^2}\to\l\aA_1^{\frac{7}{27}}\ssc\,$],
\begin{eqnarray}
&\!\!\!\!\!\!\!\!\!\!\!\!\!\!\!\!\!\!\!\!\!\!\!\!\!\!\!&
\label{A1+A2-derf-ANo}
1=\aA_1^{-\frac{71}{35}-\d_0^5}\aA_2^{-1}\left|
\frac{X_1^{3\ell_0^2}\Big(7+(13+28\d_0^5)\frac{X_2^{5}}{A_2X_1^{3\ell_0^2}}\Big)}
{A_3^4X_2^2\Big(35-(15-28\d_0^5)\frac{X_2^5}{A_2X_1^{3\ell_0^2}}\Big)}
\left(\frac{X_1^{3\ell_0^2}}{A_3^4X_2^2\Big(2-\frac1{X_1^{\ell_0^2}}\Big)}
\right)^{\frac{36}{35}+\d_0^5}\right.
\nonumber\\
&\!\!\!\!\!\!\!\!\!\!\!\!\!\!\!\!\!\!\!\!\!\!\!\!\!\!\!&
\phantom{1=\aA_1^{-\frac{71}{35}-\d_0^5}\aA_2^{-1}\,\,|}
\left.\times
\frac{X_2^5\Big(\frac{15}{28}-\d_0^5+\Big(\frac{13}{28}+\d_0^5\Big)\frac{X_1^{3\ell_0^2}}{A_1A_3^4X_2^2}\Big)}
{X_1^{3\ell_0^2}\Big(\frac54-\frac{X_1^{3\ell_0^2}}{4A_1A_3^4X_2^2}\Big)}\right|
\nonumber\\
&\!\!\!\!\!\!\!\!\!\!\!\!\!\!\!\!\!\!\!\!\!\!\!\!\!\!\!&
\phantom{1}
\le\frac{\l^3\aA_1\Big(7{\sc\!}+
{\sc\!}(13{\sc\!}+{\sc\!}28\d_0^5)\frac{\aA_1^{\frac19}}{\l^3}\Big)}
{\aA_1^{\frac{178}{105}+\d_0^5}\Big(35-(15-28\d_0^5)\frac{\aA_1^{\frac19}}{\l^3}\Big)}\Big(\frac{\l^3\aA_1}{2{\sc\!}-{\sc\!}\frac1{\l
\aA_1^{\frac7{27}}}}\Big)^{\frac{36}{35}+\d_0^5}
\frac{\Big(\frac{15}{28}
{\sc\!}-{\sc\!}\d_0^5{\sc\!}+{\sc\!}\Big(\frac{13}{28}{\sc\!}+{\sc\!}\d_0^5\Big)\l^3\Big)}{
\l^{3}\aA_1^{\frac29}\Big(\frac54-
\frac{\l^3}{4}\Big)}
\nonumber\\
&\!\!\!\!\!\!\!\!\!\!\!\!\!\!\!\!\!\!\!\!\!\!\!\!\!\!\!&
\phantom{1}
=\frac{\l^{\frac{108}{35}+3\d_0^5}\aA_1^{\frac19}\Big(7+(13+28\d_0^5)
\l^{-3}\aA_1^{\frac19}\Big)\Big(\frac{15}{28}-\d_0^5+\Big(\frac{13}{28}+\d_0^5\Big)\l^3\Big)}
{\Big(35-(15-28\d_0^5)\l^{-3\aA_1^{\frac19}}\Big)\Big(2-\l^{-1}\aA_1^{-\frac7{27}}\Big)
\Big(\frac54-\frac{\l^3}{4}\Big)}
\end{eqnarray}
On one hand, by \eqref{A2-derf-}, we have [noting that $|A_1|^{-\frac{7\d_0^2}{27}}\le1$],
\begin{eqnarray}
&\!\!\!\!\!\!\!\!\!\!\!\!\!\!\!\!\!\!\!\!\!\!\!\!\!\!\!&
\label{A1+A2-derf-}
c_1:=
\aA_1\aA_3^2\stackrel{{}^{\sc\rm\eqref{A2-derf-}}}{\le}
\frac{\xX_1^{3\ell_0^2}}{2-\xX_1^{-\ell_0^2}}
=\frac{\Big((1+\d)\aA_1^{\frac{7\d_0^2}{27}}\Big)^{3\ell_0^2}}{2-\Big((1+\d)\aA_1^{\frac{7\d_0^2}{27}}\Big)^{-\ell_0^2}}
\le(1+2\ell_0^2\d)\aA_1^{\frac79},
\end{eqnarray}
We have $|A_1A_3A_1^{-3}|\ge|A_1|^{\frac25+O(\d_0)^1}\ge1$ by
\eqref{cont-B2},\,\eqref{A1-A2-cond}\,(a),\,\eqref{ImMpP}\,(a). Thus
\begin{eqnarray}
&\!\!\!\!\!\!\!\!\!\!\!\!\!\!\!\!\!\!\!\!\!\!\!\!\!\!\!&
\label{A2-derf-+1}
1\le
|A_1A_3A_2^{-3}|
\stackrel{{}^{\sc\rm\eqref{LetNSoOP----2}\,(ii),\,(iii)}}{\le}
\frac{|A_1|\xX_2^2\Big(\frac12+\frac{\xX_1^3}{2|A_1|\xX_2^2}\Big)}
{\xX_1^3\Big(\frac54-\frac{\xX_1^3}{4|A_1|\xX_2^2}\Big)}
=
\frac{\Big(\frac12+\frac{(1+\d)^3}{2}\Big)}
{(1+\d)^3\Big(\frac54-\frac{(1+\d)^3}{4}\Big)}
=1-\frac{3\d}{4},
\end{eqnarray}
again a contradiction.
This proves \eqref{ImMpP}\,(2).
Then \eqref{ImMpP}\,(3) follows from \eqref{A1-A2-cond},\,\eqref{LetNSoOP----2}\,(ii) and \eqref{ImMpP}\,(2).
By \eqref{ImMpP}\,(2), we obtain \eqref{ImMpP}\,(3)--(6) by using formula
$\pm(|a|-|b|)\le|a+b|\le|a|+|b|$ for $a,b\in\C$ [noting that $\xX_1>1-\d$ implies that $\xX_1^{-1}\le1+O(\d)^1\ssc\,$],
where the first inequality of \eqref{ImMpP}\,(7) follows from \eqref{LetNSoOP----2}\,(iii),\,\eqref{A1-A2-cond} and \eqref{ImMpP}\,(2)--(6).
}\NOUSE{
The proof of \eqref{ImMpP}\,(4),\,(5) is a little trick, first we give a remark.
\begin{rema}\rm\label{ABCD12}
We will use the fact in the first inequality of \eqref{C+ToSayas+1}\,(d) that
$|C_1|\le1+\d$
to prove that
 equality cannot occur in the last inequality of \eqref{C+ToSayas+1}\,(d), i.e.,
$|B_1|<1+\d$.
Then we will use this fact to
to prove that equality cannot occur in the first inequality of \eqref{C+ToSayas+1}\,(d), i.e.,
$|C_1|<1+\d$.
Note that this is no problem as we only use the defining conditions on $\ol V_0$.
To illustrate our arguments, we add a parameter $\l$ in the coefficients of $\d^4$ in \eqref{C1C2},
 then we set $\l=1$.
\end{rema}
To prove \eqref{ImMpP}\,(4),
 recall from \eqref{LetNSoOP----1}\,(i) that
 $A_1=\frac{\widetilde X_1^8}{\widetilde X_2^{10+\d_0}}\big(\frac34+\frac{\widetilde X_1^4}{4}\big)$, which can be easily rewritten as%
,
\begin{eqnarray}
\label{Rewritt-}
&\!\!\!\!\!\!\!\!\!\!\!\!\!\!\!\!\!\!\!\!\!\!\!\!\!\!\!\!\!\!&
{\rm(i)\ }a+b+c=0\ \mbox{ with}
\nonumber\\[0pt]&\!\!\!\!\!\!\!\!\!\!\!\!\!\!\!\!\!\!\!\!\!\!\!\!\!\!\!\!\!\!&
{\rm(ii)\ }a=\frac{\widetilde X_1^{12}}{4A_1\widetilde X_2^{10+\d_0}}
,\ \ \ \ \
{\rm(iii)\ }b=-1,\ \ \ \ \
{\rm(iv)\ }c=\frac{3\widetilde X_1^8}{4A_1\widetilde X_2^{10+\d_0}}.
\end{eqnarray}
Therefore [the following is obtained by first
regarding \eqref{Rewritt-}\,(i) as the equation $a\eta^2+b\eta+c=0$ and then solving $\eta$
and the setting $\eta=1$],
\begin{eqnarray}
\label{A-product}
\!\!\!\!\!\!\!\!\!\!\!\!\!\!\!\!\!\!\!\!&&
\Big(1+\frac{b}{2a}\Big(1+(1-4acb^{-2})^{\frac12}\Big)\Big)
\Big(1+\frac{b}{2a}\Big(1-(1-4acb^{-2})^{\frac12}\Big)\Big)\stackrel{{}^{\sc\rm\eqref{Rewritt-}\,(i)}}{=}0,
\end{eqnarray}
where, 
\begin{eqnarray}
\label{abc====}
&\!\!\!\!\!\!\!\!\!\!\!\!\!\!\!\!\!\!\!\!\!\!\!\!\!\!&
{\rm(i)\ \ }\frac{b}{2a}\
\stackrel{{}^{\sc\rm\eqref{Rewritt-}\,(ii),\,(iii),\,\eqref{LetNSoOP----1}\,(i)}}{=}\ -\frac{C_1}{2(1-\d_0)},
\nonumber\\
&\!\!\!\!\!\!\!\!\!\!\!\!\!\!\!\!\!\!\!\!\!\!\!\!\!\!&
{\rm(ii)\ \ }4acb^{-2}\ \stackrel{{}^{\sc\rm\eqref{Rewritt-}\,(ii)\mbox{--}(iv),\,\eqref{C1-2and}\,(ii)}}{=}\
4(1-\d_0)\d_0D_1,
\end{eqnarray}
and where we always choose $(1-4acb^{-2})^{\frac12}$ to be the unique element defined by \eqref{bimeformo},
\equa{half-power--}{\dis
(1-4acb^{-2})^{\frac12}=\Big(1-4(1-\d_0)\d_0D_1\Big)^{\frac12}=1+
\mbox{$\sum\limits_{i=1}^\infty$}\binom{\frac12}{i}\Big(-4(1-\d_0)\d_0D_1\Big)^i,}
which converges absolutely  by the fact from \eqref{ImMpP}\,(6) that
$4(1-\d_0)\d_0|D_1|<1$.
Therefore, we obtain from \eqref{A-product},\,\eqref{abc====} that 
either
\begin{eqnarray}
\label{1=bbb}
&\!\!\!\!\!\!\!\!\!\!\!\!\!\!\!\!\!\!\!\!\!\!\!\!\!\!\!\!\!\!&
{\rm(i)\ }1=C_1E_1,
\\\nonumber
&\!\!\!\!\!\!\!\!\!\!\!\!\!\!\!\!\!\!\!\!\!\!\!\!\!\!\!\!\!\!&
{\rm(ii)\ }E_1:=\frac1{2(1-\d_0)}\mbox{\Large$\Big($}1+\Big(1-4(1-\d_0)\d_0D_1\Big)^{\frac12}
\mbox{\Large$\Big)$}\\\nonumber
&\!\!\!\!\!\!\!\!\!\!\!\!\!\!\!\!\!\!\!\!\!\!\!\!\!\!\!\!\!\!&
\ \ \ \ \ \ \ \ \
\stackrel{{}^{\sc\rm\eqref{half-power--}}}{=}\frac1{2(1-\d_0)}
\mbox{\Large$\Big($}2+
\mbox{$\sum\limits_{i=1}^\infty$}\binom{\frac12}{i}\Big(-4(1-\d_0)\d_0D_1\Big)^i\mbox{\Large$\Big)$},\!\!\!\!\!\!\!\!\!\!\!\!\!\!\!\!\!
\end{eqnarray}or
else
\equa{or-else-A}{\dis
1=
\frac{C_1}{2(1-\d_0)}\mbox{\Large$\Big($}1-\Big(1-4(1-\d_0)\d_0D_1\Big)^{\frac12}\mbox{\Large$\Big)$}.}
Assume we have the later case. Then,
where the first inequality follows from the fact that $(-1)^{i+1}\binom{\frac12}{i}$ is positive for all $i\ge1$,
\begin{eqnarray}
\label{A1-latercase}
\!\!\!\!\!\!\!\!\!\!\!\!\!\!\!\!\!\!\!\!&&
|C_1|^{-1}
\stackrel{{}^{\sc\rm\eqref{or-else-A}}}{=}
\frac{1}{2(1{\ssc\!}-{\ssc\!}\d_0)}\mbox{\Large$\Big|$}1{\ssc\!}-
{\ssc\!}\Big(1{\ssc\!}-{\ssc\!}4(1{\ssc\!}-{\ssc\!}\d_0)\d_0D_1\Big)^{\frac12}\mbox{\Large$\Big|$}
\stackrel{{}^{\sc\rm\eqref{half-power--}}}{=}
\frac1{2(1{\ssc\!}-{\ssc\!}\d_0)}
\mbox{\Large$\Big|$}
\mbox{$-\sum\limits_{i=1}^\infty$}\binom{\frac12}{i}\Big(-4(1{\ssc\!}-{\ssc\!}\d_0)\d_0D_1\Big)^i\mbox{\Large$\Big|$}
\!\!\!\!\!\!\!\!\!\!\nonumber\\
\!\!\!\!\!\!\!\!\!\!\!\!\!\!\!\!\!\!\!\!&&\phantom{\Big|\frac{A_1 X_2}{X_1 Z^{1+ \d_0^2}}}
\!\!\!\!\!\le
\frac1{2(1-\d_0)}
\mbox{$\sum\limits_{i=1}^\infty$}(-1)^{i+1}\binom{\frac12}{i}\Big|4(1-\d_0)\d_0D_1\Big|^i\!\!\!\!
\nonumber\\
\!\!\!\!\!\!\!\!\!\!\!\!\!\!\!\!\!\!\!\!&&
\phantom{\Big|AA\Big|^{-1}}
\!\!
\stackrel{{}^{\sc\rm\eqref{ImMpP}\,(6)}}{\le}
\frac1{2(1-\d_0)}
\mbox{$\sum\limits_{i=1}^\infty$}(-1)^{i+1}\binom{\frac12}{i}\Big(4(1-\d_0)\d_0(1+\d^2)\Big)^i
\!\!\!\!\!\!\!\nonumber\\
\!\!\!\!\!\!\!\!\!\!\!\!\!\!\!\!\!\!\!\!&&
\phantom{\Big|AA\Big|}\ \ \
\stackrel{{}^{\sc\rm\eqref{bimeformo}}}{=}\ \
\frac1{2(1-\d_0)}
\mbox{\Large$\Big($}1{\ssc}-{\ssc}\Big(1{\ssc}-{\ssc}
4(1-\d_0)\d_0(1+\d^2)\Big)^{\frac12}\mbox{\Large$\Big)$}
=\d_0+O(\d_0)^2<(1+\d)^{-1},
\end{eqnarray}
 which contradicts   \eqref{C1-2and}\,(ii).
This proves that we can
 only have \eqref{1=bbb}.
Exactly similar to the evaluation in \eqref{A1-latercase}, we can deduce from
\eqref{C1-2and}\,(ii) and the right-hand side of \eqref{1=bbb}\,(ii) the following,
\begin{eqnarray}\label{D1-is}
&\!\!\!\!\!\!\!\!\!\!\!\!\!\!&
|E_1|\stackrel{{}^{\sc\rm\eqref{1=bbb}\,(ii)}}{\ge}
\frac1{2(1-\d_0)}\mbox{\Large$\Big($}2-
\mbox{$\sum\limits_{i=1}^\infty$}(-1)^{i+1}\binom{\frac12}{i}\Big|4(1-\d_0)\d_0D_1\Big|^i\mbox{\Large$\Big)$}%
\nonumber\\
&\!\!\!\!\!\!\!\!\!\!\!\!\!\!&\phantom{|D_1|}\!
\stackrel{{}^{\sc\rm\eqref{C1-2and}\,(ii)}}{\ge}\frac1{2(1-\d_0)}\mbox{\Large$\Big($}2-
\mbox{$\sum\limits_{i=1}^\infty$}(-1)^{i+1}\binom{\frac12}{i}\Big|4(1-\d_0)\d_0(1+\d^2)\Big|^i\mbox{\Large$\Big)$}%
\nonumber\\
&\!\!\!\!\!\!\!\!\!\!\!\!\!\!&\phantom{|D_1|}\ \,\ \stackrel{{}^{\sc\rm\eqref{bimeformo}}}{=}\ \
\frac1{2(1-\d_0)}\mbox{\Large$\Big($}1+\Big(1-4(1-\d_0)\d_0(1+\d^2)\Big)^{\frac12}\mbox{\Large$\Big)$}
 =1+O(\d)^2.
\end{eqnarray}
 Thus
we obtain from \eqref{1=bbb}\,(i) that $|C_1|=|E_1|^{-1}{\ssc}\stackrel{{}^{\sc\rm\eqref{D1-is}}}{\le}{\ssc}1{\ssc}+{\ssc}O(\d)^2
{\ssc}<{\ssc}1{\ssc}+{\ssc}\d$
, which gives \eqref{ImMpP}\,(7).
}


\NOUSE{
Now assume   equality occurs the first inequality of \eqref{C+ToSayas+1}\,(c), i.e., $|\widetilde X_1|=
1-\d$. Similar to \eqref{A2-derf-}, up to $O(\d)^2$, we have
[note that in this case we do not need \eqref{A2-derf-}\,(ii), instead we need \eqref{A2-derf-+1}\,(ii)$\ssc\,$],
\begin{eqnarray}
&\!\!\!\!\!\!\!\!\!\!\!\!\!\!\!\!\!\!\!\!\!\!\!\!\!\!\!&
\label{A2-derf-+1}
{\rm(i)\ }|C_2|\le
1-\frac{(1+\d_0)\d}{1+2\d_0},\ \ \ \ \ \
{\rm(ii)\ }
|1+\d_0-\d_0\widetilde X_1^{\ell_0+1}|\ge
1+(1+\d_0)\d,
\nonumber\\
&\!\!\!\!\!\!\!\!\!\!\!\!\!\!\!\!\!\!\!\!\!\!\!\!\!\!\!&
{\rm(iii)\ }
\Big|1+\d_0^3-\d_0^3\widetilde X_1C_2\Big|\ge1+\frac{\d_0^3(2+3\d_0)\d}{1+2\d_0}.
\end{eqnarray}
Then as in \eqref{A2-derf}, up to $O(\d)^2$, we have
\begin{eqnarray}
&\!\!\!\!\!\!\!\!\!\!\!\!\!\!\!\!\!\!\!\!\!\!\!\!\!\!\!&
\label{A2-derf+1}\dis
1\stackrel{{}^{\sc\rm\eqref{A1-A2-cond}\,(i)}}{\le}
|A_1|^{\d_0^5+O(\d_0)^6}\stackrel{{}^{\sc\rm\eqref{B1-ep}}}{=}|A_1|^
{a_3a_1^{-1}(\d_0+\d_0^2-\d_0^4)-\d_0^4-a_2a_1^{-1}\d_0^2}
\stackrel{{}^{\sc\rm\eqref{cont-B2},\,\eqref{ImMpP}\,(1)}}{\le}
\Big|\frac{A_3^{\d_0+\d_0^2-\d_0^4}}{A_1^{\d_0^4}A_2^{\d_0^2}}\Big|
\nonumber\\
&\!\!\!\!\!\!\!\!\!\!\!\!\!\!\!\!\!\!\!\!\!\!\!\!\!\!\!&
\ \ \ \ \,\stackrel{{}^{\sc\rm\eqref{LetNSoOP----2}}}{=} \ \ \mbox{\Large$\Big|$}\Big(\frac{C_2}{X_2}
\Big)^{\d_0+\d_0^2-\d_0^4}
\frac{X_2^{\d_0(1-\d_0^3)}}{\widetilde X_1(1+\d_0-\d_0\widetilde X_1^{\ell_0+1})}
\cdot\frac{X_2^{\d_0^2}}{(1+\d_0^3-\d_0^3\widetilde X_1C_2)^{\d_0}\widetilde X_1^{\d_0}C_2^{\d_0}}\mbox{\Large$\Big|$}
\nonumber\\
&\!\!\!\!\!\!\!\!\!\!\!\!\!\!\!\!\!\!\!\!\!\!\!\!\!\!\!&
\ \ \ \ \ \ \le\ \ \ \
\frac{\Big(1-\frac{(1+\d_0)\d}{1+2\d_0}\Big)^{\d_0^2-\d_0^4}}
{(1-\d)^{1+\d_0}\Big(1+(1+\d_0)\d\Big)\Big(1+\frac{\d_0^3(2+3\d_0)\d}{1+2\d_0}\Big)}
=
1+\Big(-\d_0^2+O(\d_0)^3\Big)\d,
\end{eqnarray}
which is again a contradiction.
}%
\NOUSE{
Now assume equality occurs in the last inequality of
\eqref{ImMpP}\,(3), i.e., $\tilde\xX_1=(1+\d)|A_1|^{-\d_0^2}$.
Noting from \eqref{B1-ep} that $(a_3+\frac{(2+\d_0^3)a_2}{3+4\d_0})a_1^{-1}=\frac{\d_0}{3}+O(\d_0)^2$, we have
[noting from \eqref{A1-A2-cond}\,(a) that $|A_1|\ge1$],
\begin{eqnarray}
&\!\!\!\!\!\!\!\!\!\!\!\!\!\!\!\!\!\!\!\!\!\!\!\!\!\!\!&
\label{A2-derf+1}\dis
|A_1|^{\frac{\d_0}{3}+O(\d_0)^2}\stackrel{{}^{\sc\rm\eqref{cont-B2},\,\eqref{ImMpP}\,(1)}}{\le}
|A_3A_2^{\frac{2+\d_0^3}{3+4\d_0}}|
\stackrel{{}^{\sc\rm\eqref{LetNSoOP----2}\,(ii),\,(iii)}}{\le}
\frac{1-\d_0^2+\d_0^2\tilde\xX_1}{\tilde \xX_1^2}
=\frac{1-\d_0^2+\d_0^2(1+\d)|A_1|^{-\d_0^2}}{(1+\d)^2|A_1|^{-2\d_0^2}}\!\!\!\!\!\!\!\!
\nonumber\\
&\!\!\!\!\!\!\!\!\!\!\!\!\!\!\!\!\!\!\!\!\!\!\!\!\!\!\!&
\phantom{|A_1|^{\frac{\d_0}{3}+O(\d_0)^2}\ \ \ \ \ \ \ \ }\le\ \ \ \ \
|A_1|^{2\d_0^2}\frac{1-\d_0^2+\d_0^2(1+\d)}{(1+\d)^2}
=|A_1|^{2\d_0^2}\Big(1+\Big(-2+O(\d_0)^2\Big)\d\Big),
\end{eqnarray}
which is again a contradiction. This prove
\eqref{ImMpP}\,(3).
By \eqref{ImMpP}\,(2),\,(3) and
\eqref{tX11111}, we immediately obtain \eqref{ImMpP}\,(4).
}\NOUSE{using the formula $\pm(|a|-|b|)\le|a+b|\le|a|+|b|$ for $a,b\in\C$.
By \eqref{ImMpP}\,(3),\,(4) we see from \eqref{LetNSoOP----2}\,(ii) that
\eqref{ImMpP}\,(5) holds [using convention \eqref{MSmde33333+}$\ssc\,$].
Then by \eqref{ImMpP}\,(2),\,(4),\,(5) and \eqref{LetNSoOP----2}\,(i), we see that $|\widetilde X_1|>\d_2^{\d_1}$.
From this and  \eqref{ImMpP}\,(3),\,\eqref{tX11111}, we obtain \eqref{ImMpP}\,(6).
}\NOUSE{
\equa{FirAssX1}{\dis
|\widetilde X_1|\le\d_1^{\frac1{10}}=O(\d_1)^{\frac1{10}}=O(\d_0)^1.}
 Then
the left-hand side of \eqref{LetNSoOP----2}\,(iii) is a $1+O(\d_0)^1$ element by \eqref{In-A1A2}\,(b), the second term in the right-hand side is an $O(\d_0)^1$ element by \eqref{FirAssX1}, we obtain
\equa{Dwn3n3}{\dis
\frac{(1-\d_0)X_2}{\widetilde X_1^3}=1+O(\d_0)^1.}
 Then \eqref{LetNSoOP----2}\,(ii) gives
\equa{1+memesss}{\dis\!\!\!\!\!
O(\d_0)^2\stackrel{{}^{\sc\rm\eqref{In-A1A2}\,(b),\,\eqref{FirAssX1}}}{=}\widetilde X_1^2A_2
\stackrel{{}^{\sc\rm\eqref{LetNSoOP----2}\,(ii),\,\eqref{Dwn3n3}}}{=}
\frac{1{\ssc\!}+{\ssc\!}\a_1}{1{\ssc\!}-{\ssc\!}\d_0}\Big(1{\ssc\!}+{\ssc\!}O(\d_0)^1\Big)
{\ssc\!}-{\ssc\!}\a_1\stackrel{{}^{\sc\rm\eqref{a0BiggerThen}\,(iii)}}{=}1{\ssc\!}+{\ssc\!}O(\d_0)^1,
\!\!\!\!}
a contradiction. Next assume $|\widetilde X_1|\ge\ell_1^{\frac1{10}}$.
Similar to the proof of \eqref{Dwn3n3}, this time we can use
\eqref{LetNSoOP----2}\,(ii),\,\eqref{In-A1A2}\,(b) to obtain that $\frac{(1+\a_1)X_2}{\widetilde X_1^5}=1+O(\d_0)^1$. Then \eqref{LetNSoOP----2}\,(iii) shows
\equa{NowSSS}{\dis
O(\d_0)^2=\widetilde X_1^{-2}A_2=\frac{1-\d_0}{1+\a_1}\Big(1+O(\d_0)^1\Big)+O(\d_0)^1=\frac1{1+\a_1}+O(\d_0)^1,
}
which again is a contradiction by \eqref{a0BiggerThen}\,(iii). This proves \eqref{ImMpP}\,(2).
Now assume $|X_2|\le\d_1^{\frac35}$. Then $\frac{X_2}{\widetilde X_1^5},\frac{X_2}{\widetilde X_1^3}$ are
$O(\d_0)^1$ elements by \eqref{ImMpP}\,(2),
we obtain from \eqref{LetNSoOP----2}\,(ii),\,(iii) and \eqref{In-A1A2} that $\frac{\a_1}{\widetilde X_1^2}=1+O(\d_0)^1$,
$\d_0\widetilde X_1=1+O(\d_0)^1$, and
\equa{MaoDuanpp}{\dis
O(\d_0)^2\stackrel{{}^{\sc\rm\eqref{a0BiggerThen}\,(iii)}}{=}\d_0^2\a_1\Big(1+O(\d_0)^1\Big)=(\d_0\widetilde X_1)^2\Big(\frac{\a_1}{\widetilde X_1^2}\Big)=1+O(\d_0)^1,
}
a contradiction. This proves the first inequality of \eqref{ImMpP}\,(3), the second inequality
follows from \eqref{ImMpP}\,(1.b),\,\eqref{In-A1A2}\,(a) and \eqref{B1-ep}.
}\NOUSE{
Assume  equality occurs in the second inequality of \eqref{C+ToSayas+1}\,(d), i.e., $|X_1|=1+\d$.
Then
\begin{eqnarray}
\label{A1---}
\!\!\!\!\!\!\!\!\!\!\!\!\!\!\!\!\!\!\!\!\!\!\!\!\!\!\!\!\!\!\!\!\!\!\!\!\!\!\!\!\!\!
\!\!\!\!\!\!\!\!\!\!\!\!\!\!\!\!\!\!\!\!\!\!\!\!\!\!\!\!\!\!\!\!\!\!&\!\!\!\!\!\!\!\!\!\!\!\!\!\!\!\!\!\!\!\!\!\!\!\!\!\!\!\!\!\!\!\!\!\!\!&
1{\ssc\!}+{\ssc\!}O(\d)^2
\stackrel{{}^{\sc\rm\eqref{MSMSH1H20}\,(i),\,\eqref{C+ToSayas+1}\,(i)}}{\le}
|A_1|
\stackrel{{}^{\sc\rm\eqref{Poly-x1-x2}\,(i),\,\eqref{ImMpP}\,(1.b)}}{\le}
(1{\ssc\!}+{\ssc\!}\d)^{-\ell_0^2}\Big(\d_0{\ssc\!}+{\ssc\!}(1{\ssc\!}
-{\ssc\!}\delta_0)(1{\ssc\!}+{\ssc\!}\d)^{\ell_0^2(1+\d_0)}\Big){\ssc\!}+{\ssc\!}O(\d)^2
\!\!\!\!\!\!\!\!\!\!\!\!\!\!\!\!\!\!\!\!\!\!\!\!\!\!\!\!\!\!\!\!\!\!\!\!\!\!\!\!\!\!\!\!\!\!\!\!\!\!\!\!\!\!\!\!\!\!\!\!\!\!\!\!\!\!\!\!\!\!\!\!\!\!\!\!\!\!\!\!\!\!\!\!\!\!\!\!\!\!\!\!\!\!\!\!
\nonumber\\
\!\!\!\!\!\!\!\!\!\!\!\!\!\!\!\!\!\!\!\!\!\!\!\!\!\!\!\!\!\!\!\!\!\!\!\!\!\!\!\!\!\!\!\!\!\!\!\!\!
\!\!\!\!\!\!\!\!\!\!\!\!\!\!\!\!\!\!\!\!\!\!\!\!\!\!\!&\!\!\!\!\!\!\!\!\!\!\!\!\!\!\!\!\!\!\!\!\!\!\!\!\!\!\!\!\!\!\!\!\!\!\!&
\phantom{1+O(\d)^2{\ssc\!}{\ssc\!}\ \ \ \ \ \ \ \, \  }=\ \ \ \ \ \ \ 1-\d+O(\d)^2,
\!\!\!\!\!\!\!\!\!\!\!\!\!\!\!\!\!\!\!\!\!
\end{eqnarray}
a contradiction, which proves \eqref{ImMpP}\,(2).
Assume  equality occurs in the second inequality of \eqref{C+ToSayas+1}\,(e), i.e., $|H_0|=1+\d$.
Then
\begin{eqnarray}
\label{H0---}
\!\!\!\!\!\!\!\!\!\!\!\!\!\!\!\!\!\!\!\!\!\!\!\!\!\!\!&&
1+O(\d)^2
\stackrel{{}^{\sc\rm\eqref{MSMSH1H20}\,(iii),\,\eqref{C+ToSayas+1}\,(i)}}{\le}
|A_2|
\nonumber\\
\!\!\!\!\!\!\!\!\!\!\!\!\!\!\!\!\!\!\!\!\!\!\!\!\!\!\!&&\ \ \ \ \ \ \ \ \ \,
\stackrel{{}^{\sc\rm\eqref{LetNSoOP----1}\,(iii),\,\eqref{ImMpP}\,(1.b)}}{\le}
(1+\d)^{-\ell_0^2}\Big(\d_0+(1-\d_0)(1+\d)^{\ell_0^2(1+\d_0)}\Big)+O(\d)^2\!\!\!\!\!\!\!\!\!\!\!\!
\nonumber\\
\!\!\!\!\!\!\!\!\!\!\!\!\!\!\!\!\!\!\!\!\!\!\!\!\!\!\!&&\ \ \ \ \ \ \ \ \ \
\ \ \ \ \ \ \ \ \ \ \ \
=\ \ \ \ \ \ 1-\d+O(\d)^2,
\end{eqnarray}
a contradiction, which proves \eqref{ImMpP}\,(3).
We have
\begin{eqnarray}
\label{H0---+1}
\!\!\!\!\!\!\!\!\!\!\!\!\!\!\!\!\!\!\!\!\!\!\!\!\!\!\!&&
|H_0|
\stackrel{{}^{\sc\rm\eqref{LetNSoOP----1}\,(iii),\,\eqref{ImMpP}\,(1.b),\,(3)}}{\le}
\Big|\frac{A_1^{\ell_0}X_1^{\ell_0^3}}{A_2}\Big|^{\d_0^2}
\Big(\d_0+(1-\d_0)(1+\d^3)^{\ell_0^2(1+\d_0)}\Big)^{\d_0^2}
\\
\nonumber\!\!\!\!\!\!\!\!\!\!\!\!\!\!\!\!\!\!\!\!\!\!\!\!\!\!\!&&\ \ \ \ \ \ \
\stackrel{{}^{\sc\rm\eqref{MSMSH1H20},\,\eqref{cont-B2},\,\eqref{ImMpP}\,(2)}}{\le}
|A_1|^{\ell_0-b_2b_1^{-1}}+O(\d)^3
\stackrel{{}^{\sc\rm\eqref{B1-ep}}}{=}|A_1|^{-\ell_0^2(1+O(\d_0)^1)}+O(\d)^3\le1+O(\d)^3.\!\!\!\!\!\!\!\!\!\!\!\!
\end{eqnarray}
Then
\begin{eqnarray}
\label{H0---+2}
\!\!\!\!\!\!\!\!\!\!\!\!\!\!\!\!\!\!\!\!\!\!\!\!\!\!\!&&
1+\d_0^2-\d_0^2|X_1|^{\ell_0}
\stackrel{{}^{\sc\rm\eqref{ImMpP}\,(1.b)}}{\le}
1+\d_0^2-\d_0^2|X_1^{\ell_0}X_2^{\ell_0(1+\d_0)}|+O(\d)^3
\\
\nonumber\!\!\!\!\!\!\!\!\!\!\!\!\!\!\!\!\!\!\!\!\!\!\!\!\!\!\!&&
\phantom{1+\d_0^2-\d_0^2|X_1|^{\ell_0}\ \ \,}\stackrel{{}^{\sc\rm\eqref{Poly-x1-x2}\,(ii)}}{\le}|H_0|+O(\d)^3
\stackrel{{}^{\sc\rm\eqref{H0---+1}}}{\le}1+O(\d)^3.
\!\!\!\!\!\!\!\!\!\!\!\!
\end{eqnarray}
The above implies \eqref{ImMpP}\,(4).
Equ.~\eqref{ImMpP}\,(2) with \eqref{ImMpP}\,(4) means that $|X_1|=1+O(\d)^3$. Then \eqref{H0---+2} shows that $|H_0|\ge1+O(\d)^3$, which with
\eqref{H0---+1} implies \eqref{ImMpP}\,(5).
%
%
%
By \eqref{ImMpP}\,(1.b) and the facts that $|X_1|=1+O(\d)^2$
and $|A_1|=|A_1|+O(\d)^5\ge1+O(\d)^5$, we  obtain from \eqref{C1-2and}\,(iii) that $|D_1|\le1+O(\d)^3$, which implies
 \eqref{ImMpP}\,(6).
\NOUSE{
By \eqref{C1-2and}, \eqref{cont-B2} [which obviously also holds for elements in $\ol V_0$
by \eqref{C+LetNSoOP}\,(i)$\ssc\,$], \eqref{ImMpP}\,(1), we obtain%
,
\begin{eqnarray}
\label{1=bbb+1+1}
&\!\!\!\!\!\!\!\!\!\!\!\!\!\!\!\!\!\!\!\!\!\!\!\!\!\!\!\!\!\!&
|A_2''|\stackrel{{}^{\sc\rm\eqref{C1-2and}}}{=}\Big|\frac{A_1^{4\ell_0^4}Z^4}{A_2^2}\Big|
\stackrel{{}^{\sc\rm\eqref{cont-B2},\,\eqref{ImMpP}\,(1)}}{\le}|A_1|^{-4\ell_0^4-8a_1^{-1}-2a_2a_1^{-1}}
+O(\d)^3
\nonumber\\&\!\!\!\!\!\!\!\!\!\!\!\!\!\!\!\!\!\!\!\!\!\!\!\!\!\!\!\!\!\!&
\phantom{{\rm(i)\ }|C_1|}\,\ \stackrel{{}^{\sc\rm\eqref{B1-ep}}}{=}|A_1|^{-2\ell_0^7(1+O(\d_0)^1)}+O(\d)^3
\stackrel{{}^{\sc\rm\eqref{A1andA2}\,(i)}}{<}1+\d^2
.
\end{eqnarray}
This proves \eqref{ImMpP}\,(4).
}
\NOUSE
{
Note that we may often need to use the fact  from \eqref{LetNSoOP}\,(i) that $|A_{\rOnE}|\le1$ and to apply \eqref{A0===+1} frequently. We also often need to use the fact from
\eqref{ImMpP}
\,(3),\,(6)
, 
\equa{X202020}{\dis \!\!\!\!\!\!\!\!\!
{\rm(i)\ } X_2{\ssc}={\ssc}\tildeX_1^{-1}A_3^{-1}{\ssc}+{\ssc}O(\d)^3, \ \ \ \ \ {\rm(ii)\ }  |A_3|^{-1}{\ssc}<{\ssc}|A_1|^{-c_3c_1^{-1}}+{\ssc}O(\d)^3.\!\!\!\!\!\!\!\!\!\!}
%
%
%
%
Now assume $(1-\d^3)| A_{\rOnE}|^{a_1}\ge|\tildeX_{\ZeRo}|$. Then
by \eqref{LetNSoOP}\,(ii),
\equa{x0==wo3499}{\dis
|\tildeX_{\ZeRo}|=(1-\d^3)| A_{\rOnE}|^{a_1}.}
}\NOUSE
{%
We want to prove \eqref{X0sm,m-1}
.
Noting 
from the second and third inequalities of \eqref{LetNSoOP}\,(i) that, as mentioned in Remark \ref{Rema-V2-def}\,(ii),
$|A_2|=|A_1|^{-\frac{587}{20}+O(\d)^3}=|A_1|^{-\frac{587}{20}}+O(\d)^3$
. From this and \eqref{ImMpP}\,(4),\,(5),\,\eqref{x0==wo3499}, 
we have%
,}\NOUSE
{%
By \eqref{ImMpP}\,(4),\,\eqref{X202020}, we have [the following is the reason of defining $a_1$ in \eqref{defi-a1}$\ssc\,$],
\begin{eqnarray}
\label{X0sm,m-1}
&\!\!\!\!\!\!\!\!\!\!\!\!\!\!\!\!\!\!\!\!\!\!\!\!\!\!\!\!\!\!\!&
1+O(\d)^3=|A_1^{-1}(\tildeX_1A_3)^{-1}\tildeX_1^{\ell_0}|
\nonumber\\&\!\!\!\!\!\!\!\!\!\!\!\!\!\!\!\!\!\!\!\!\!\!\!\!\!\!\!\!\!\!\!&
\phantom{1+O(\d)^3}
\le(1-\d^3)^{\ell_0-1}|A_1^{-1-(a_1+c_3c_1^{-1})+\ell_0a_1}|+O(\d)^3=(1-\d^3)^{\ell_0-1}+O(\d)^3,
\!\!\!\!\!\!\end{eqnarray}
which is a contradiction. This proves \eqref{ImMpP}\,(9).
}}{
\NOUSE
{%
The proof of \eqref{ImMpP}\,(2) is a little tricky, as we already mentioned in Remark \ref{B1-C1}, $B_1,C_1$ are related to $A_1$ in such a way that
\eqref{LetNSoOP----1}\,(iv) is equivalent to \eqref{Re-Writtt}, from which, we can obtain
\eqref{1=bbb-remark}, which is equivalent to \eqref{1=bbb} below.
We give another remark below.
\begin{rema}\rm\label{ABCD12}For any element $(p_1,p_2)$ in $\ol V_0$, by Lemma \ref{lemm-condition-XZ},
we will use the condition in the last inequality of  \eqref{C+ToSayas+1}\,(c) that
$1\le(1+\eE_1^2)|C_1|$
to prove that
 equality cannot occur in the first  inequality of \eqref{C+ToSayas+1}\,(c), i.e.,
$|B_1|<1+\eE_1$.
Then we will use the condition in the first inequality of \eqref{C+ToSayas+1}\,(c) that
$|B_1|\le1+\eE_1$
to prove that equality cannot occur in the last inequality of \eqref{C+ToSayas+1}\,(c), i.e.,
$1<(1+\eE_1^2)|C_1|$.
Note that this is no problem as we only use the defining conditions \eqref{ToSayas+1} on $V_0$ [which implies  \eqref{C+ToSayas+1} by Lemma \ref{lemm-condition-XZ}].
\end{rema}
}
\NOUSE
{%
\begin{lemm}\label{procesi-lemm3}For $(p_1,p_2)\in\ol V_0$, we have
\begin{eqnarray}
\label{Formmeme}\dis
1-\d<\bB_1
\stackrel{{}^{\sc\rm\eqref{LetNSoOP----1-redefine++}\,(iv)}}{\equiv}
\aA_1\aA_3^{-(\ell_0^2+2\ell_0+1)}
\tilde\xX_1^{-\ell_0^4(2+4\d_0-3\d_0^2-6\d_0^3)}
.
\end{eqnarray}
\end{lemm}
\noindent{\it Proof.~}
\NOUSE{
Write $B_1^{-1}$ as  its polar decomposition $B_1^{-1}=r e^{i\theta}$ for
some $r\in\R_{>0},\,\th\in\R$ with $0\le\th<2\pi$.
\begin{clai}\label{procesi-clai1}
The function  of $\theta$ $($when $r$ is fixed$)$ given by
\equa{funccccc}{\dis
|5 -4B_1^{-1}|=   |5 -4 r  e^{ i\theta}|,}  has a minimum
$  |5 -4 r   |$ for $\theta=0$.
\end{clai}
To prove the claim, note that the function  of $\theta$ given by $  |5 -4 r  e^{ i\theta}|$   has a minimum when its
square $ |5 -4 r  e^{ i\theta}|^2$  does. Then
\begin{eqnarray}
\label{Claim--11111}
&\!\!\!\!\!\!\!\!\!\!\!\!\!\!\!\!\!\!\!\!\!\!\!\!&
{\rm(i)\ }
 |5 -4 r  e^{ i\theta}|^2=\big(5-4r \cos(\theta)\big)^2+\big(4r  \sin(\theta)\big)^2
 =25-40r \cos(\theta)+r^2\big(4\cos(\theta)\big)^2+r^2\big(4 \sin(\theta)\big)^2
 \!\!\!\!\!\!\!\!\!\!
\nonumber\\
&\!\!\!\!\!\!\!\!\!\!\!\!\!\!\!\!\!\!\!\!\!\!\!\!&
\phantom{{\rm(i)\ }  |5 -4 r  e^{ i\theta}|^2}
 =
  25+16r^2-  40 r \cos(\theta),
\nonumber\\
&\!\!\!\!\!\!\!\!\!\!\!\!\!\!\!\!\!\!\!\!\!\!\!\!&
{\rm(ii)\ } \min\big(25+16r^2-  40 r \cos(\theta)\big)=\big(25+16r^2-  40 r \cos(\theta)\big)|_{\cos(\th)=1}
\nonumber\\
&\!\!\!\!\!\!\!\!\!\!\!\!\!\!\!\!\!\!\!\!\!\!\!\!&
\phantom{{\rm(ii)\ } \min\big(25+16r^2-  40 r \cos(\theta)\big)}
=\big(25+16r^2-  40 r \cos(\theta)\big)|_{\th=0} =25+16r^2 -4 0 r,\ \ \  \implies
\nonumber\\
&\!\!\!\!\!\!\!\!\!\!\!\!\!\!\!\!\!\!\!\!\!\!\!\!&
{\rm(iii)\ } \min |5 -4 r  e^{ i\theta}|= \sqrt{25+16r^2 -40  r }= |5-4r|.
\end{eqnarray}
This proves the claim.
Now by \eqref{CCCC===},\,\eqref{--Re-Writtt},
}%
We have
\begin{eqnarray}\label{hB1}
&\!\!\!\!\!\!\!\!\!\!\!\!\!\!\!\!\!\!\!\!\!\!\!\!\!\!\!\!\!\!\!&
h_1(B_1):=\bB_1^{-\frac{2(1-2\d_0^2)}{2-3\d_0^2}}\cC_1
\stackrel{{}^{\sc\rm\eqref{Re-Writtt}}}{=}
 \bB_1^{-\frac{2-4\d_0^2}{2-3 \d_0^2}+1}
 |\ell_0+1{\sc}-
 {\sc}\ell_0B_1
 |
 \nonumber\\
&\!\!\!\!\!\!\!\!\!\!\!\!\!\!\!\!\!\!\!\!\!\!\!\!\!\!\!\!\!\!\!&
\phantom{h_1(B_1)}
 ={\sc}\bB_1^{\frac{\d_0^2}{2-3 \d_0^2}}
 |\ell_0+1-\ell_0B_1
|\stackrel{{}^{\sc\rm \eqref{LetNSoOP----1-redefine++}\,(viii)}}{\le}1+O(\d)^2
.
\end{eqnarray}
Using the well-known formula $\pm(|a|-|b|)\le|a-b|$ for any $a,b\in\C$, i.e., $\big||a|-|b|\big|\le|a-b|$, we
obtain that
$|\ell_0+1-\ell_0\bB_1
|\le|\ell_0+1-\ell_0B_1
|$. Thus
%
we have, 
for any $(p_1,p_2)\in\ol V_0$,
\equa{g(r)}{\dis
 g_1(\bB_1):= \bB_1^{\frac{ \d_0^2}{2-3\d_0^2}}|\ell_0+1 -\ell_0 \bB_1
|\leq h_1(B_1)\le  1+O(\d) ^2
 .}
 If $ \bB_1\ge
1+\d_0>
 1$ (recalling that $\d_0\ell_0=1$), then 
in particular
we have \eqref{Formmeme}, so we may assume   $ \bB_1<1+\d_0
$. 
 We have then \equa{g(r)1}{\dis g_1(\bB_1)=\bB_1^{\frac{ \d_0^2}{2-3 \d_0^2}}(\ell_0+1-\ell_0\bB_1
 )  \leq  1+O(\d) ^2
,}
{    and 
\begin{eqnarray}
\label{d-g(r)=}
&&\!\!\!\!\!\!\!\!\!\!\!\!\!\!\!\!\!\!\!\!\!\!
\frac{d g_1}{d \bB_1}=\frac{ \d_0^2}{2-3 \d_0^2}\bB_1^{\frac{ \d_0^2}{2-3 \d_0^2}-1}
(\ell_0+1 -\ell_0 \bB_1
) -\ell_0\bB_1^{\frac{ \d_0^2}{2-3 \d_0^2}}
\nonumber\\
&&\!\!\!\!\!\!\!\!\!\!\!\!\!\!\!\!\!\!\!\!\!\!
\phantom{\frac{d g_1}{d \bB_1}}
=
\frac{1}{2-3\d_0^2}\bB_1^{-\frac{2-4\d_0^2}{2-3\d_0^2}}
\Big(\d_0^2(\ell_0+1-\ell_0\bB_1)-\ell_0(2-3\d_0^2)\bB_1\Big)
\nonumber\\
&&\!\!\!\!\!\!\!\!\!\!\!\!\!\!\!\!\!\!\!\!\!\!
\phantom{\frac{d g}{d \bB_1}}
=\frac1{2-3\d_0^2}\bB_1^{-\frac{2-4\d_0^2}{2-3\d_0^2}}\Big(
\d_0+\d_0^2-2(\ell_0-\d_0)\bB_1\Big),\end{eqnarray}
so $g_1(\bB_1)$
  is 
a strictly decreasing function on $\bB_1$ when $\frac{\d_0^2}{2(1-\d_0)}\le\bB_1\le1$%
.
  By \eqref{C+LetNSoOP}\,(iii), we have $\bB_1\ge\frac12$.
Since $g_1(\bB_1)$ is
  strictly decreasing when $\frac12\le\bB_1\le1$ and 
\equa{g(r)11111}{\dis\!\!\!\!\!\!\!\!\!
  g_1\big(1{\sc\!}-{\sc\!}\d\big)
\stackrel{{}^{\sc\rm\eqref{g(r)1}}}{=}
\Big(1{\sc\!}-{\sc\!}\d\Big)^{\frac{\d_0^2}{2-3\d_0^2}}
\Big(\ell_0{\sc\!}+{\sc\!}1{\sc\!}-{\sc\!}\ell_0(1{\sc\!}-{\sc\!}\d)\Big)
  {\sc\!}=
  {\sc\!}1{\sc\!}+{\sc\!}\Big(\ell_0{\sc\!}+{\sc\!}O(\d_0)^0\Big)\d
  {\sc\!}+{\sc\!}O(\d)^2{\sc\!}>{\sc\!}1{\sc\!}+{\sc\!}O(\d)^2,\!\!\!\!\!\!}
  in order for \eqref{g(r)1} to hold, we must have $\bB_1>1-\d
  $.
This proves \eqref{Formmeme}.}\hfill$\Box$\vskip7pt
}\NOUSE
{
To prove \eqref{ImMpP}\,(2), we first deduce \eqref{1=bbb-remark} in another way.
Recall from definition \eqref{LetNSoOP----1}\,(iv)
 that
\equa{A1-rewww}{\dis  A_1=
\frac{\a_1^{101}A_3^{101\lL (1-20\dD)}\widetilde X_1^{2828}}{Z^{101}}
\Big(\frac15 +\frac{4A_2^{\lL-1}}{5
 A_3^{202\lL (1-10\dD)}\widetilde X_1^{3636}}\Big),} which can be 
rewritten as the following [noting from \eqref{Amsmene} and notation \eqref{denote-t2} that $A_1$ is invertible in $\ol V_0$]%
,
\begin{eqnarray}
\label{Rewritt-}
&\!\!\!\!\!\!\!\!\!\!\!\!\!\!\!\!\!\!\!\!\!\!\!\!\!\!\!\!\!\!&
{\rm(i)\ }a+b+c=0\ \mbox{ with}
\nonumber\\[0pt]&\!\!\!\!\!\!\!\!\!\!\!\!\!\!\!\!\!\!\!\!\!\!\!\!\!\!\!\!\!\!&
{\rm(ii)\ }a=\frac{4\a_1^{101}A_2^{\lL(1 - \dD)}}{5A_1 A_3^{101\lL}\widetilde X_1^{808} Z^{101}}
,\ \ \ \ \
{\rm(iii)\ }b=-1,\ \ \ \ \
{\rm(iv)\ }c=\frac{\a_1^{101}A_3^{101\lL (1-20\dD)}\widetilde X_1^{2828}}{5A_1Z^{101}}
.
\end{eqnarray}
Note that \eqref{Rewritt-}\,(i) is equivalent to the following (noting that $a,b,c$ are all invertible),
%
\begin{eqnarray}
\label{A-product}
\!\!\!\!\!\!\!\!\!\!\!\!\!\!\!\!\!\!\!\!&&
{\rm(i)\ } 1-4acb^{-2}=1+4a(a+b)b^{-2}=(1+2ab^{-1})^2,\ \ \ \implies\ \ \
\nonumber\\\!\!\!\!\!\!\!\!\!\!\!\!\!\!\!\!\!\!\!\!&&
{\rm(ii)\ }1=-\frac{b}{2a}\Big(1+(1-4acb^{-2})^{\frac12}\Big),\mbox{ \ or else \ }
{\rm(iii)\ }1=-\frac{b}{2a}\Big(1-(1-4acb^{-2})^{\frac12}\Big),
\end{eqnarray}
where [noting from \eqref{LetNSoOP----1}\,(vi) that $C_1$ is invertible in $\ol V_0$ as $A_1$ is invertible], 
\begin{eqnarray}
\label{abc====}
&\!\!\!\!\!\!\!\!\!\!\!\!\!\!\!\!\!\!\!\!\!\!\!\!\!\!&
{\rm(i)\ \ }\frac{b}{2a}\
\stackrel{{}^{\sc\rm\eqref{Rewritt-}\,(ii),\,(iii),\,\eqref{LetNSoOP----1}\,(v)}}{=}
\ -\frac{5B_1}{8},
\nonumber\\
&\!\!\!\!\!\!\!\!\!\!\!\!\!\!\!\!\!\!\!\!\!\!\!\!\!\!&
{\rm(ii)\ \ }4acb^{-2}\ \stackrel{{}^{\sc\rm\eqref{Rewritt-}\,(ii)\mbox{--}(iv),\,\eqref{LetNSoOP----1}\,(vi)}}{=}\
\frac{16}{25C_1},
\end{eqnarray}
and where we always choose $(1-4acb^{-2})^{\frac12}$ to be the unique element defined by \eqref{bimeformo},
\equa{half-power--}{\dis
(1-4acb^{-2})^{\frac12}=\Big(1-\frac{16}{25C_1}\Big)^{\frac12}=1+
\mbox{$\sum\limits_{i=1}^\infty$}\binom{\frac12}{i}\Big(-\frac{16}{25C_1}\Big)^i,}
which converges absolutely  by the fact from \eqref{C+ToSayas+1}\,(c) that
$\frac{16}{25|C_1|}\le\frac{16(1+\eE_1^2)}{25}<1$.
Therefore, we obtain from \eqref{A-product}\,(ii) or \eqref{A-product}\,(iii), and \eqref{abc====} the following:
either [noting that the following is equivalent to \eqref{1=bbb-remark}$\ssc\,$]
\begin{eqnarray}
\label{1=bbb}
&\!\!\!\!\!\!\!\!\!\!\!\!\!\!\!\!\!\!\!\!\!\!\!\!\!\!\!\!\!\!\!\!\!\!\!\!\!&
{\rm(i)\ }1=B_1D_1,
\ \ \
{\rm(ii)\ }D_1:=\frac5{8}\mbox{\Large$\Big($}1+\Big(1-\frac{16}{25C_1}\Big)^{\frac12}
\mbox{\Large$\Big)$}
\stackrel{{}^{\sc\rm\eqref{half-power--}}}{=}\frac58
\mbox{\Large$\Big($}2+
\mbox{$\sum\limits_{i=1}^\infty$}\binom{\frac12}{i}\Big(-\frac{16}{25C_1}
\Big)^i\mbox{\Large$\Big)$},\!\!\!\!\!
\end{eqnarray}or
else
\equa{or-else-A}{\dis
1=
\frac{5B_1}{8}\mbox{\Large$\Big($}1-\Big(1-\frac{16}{25C_1}\Big)^{\frac12}\mbox{\Large$\Big)$}.}
Let us assume we have the later case \eqref{or-else-A}. Then we have the following,
where the first inequality follows from the fact that $(-1)^{i+1}\binom{\frac12}{i}$ is positive for all $i\ge1$,
\begin{eqnarray}
\label{A1-latercase}
\!\!\!\!\!\!\!\!\!\!\!\!\!\!\!\!\!\!\!\!&&
|B_1|^{-1}
\stackrel{{}^{\sc\rm\eqref{or-else-A}}}{=}
\frac{5}{8}\mbox{\Large$\Big|$}1-
\Big(1-\frac{16}{25C_1}\Big)^{\frac12}\mbox{\Large$\Big|$}
\stackrel{{}^{\sc\rm\eqref{half-power--}}}{=}
\frac5{8}
\mbox{\Large$\Big|$}
\mbox{$-\sum\limits_{i=1}^\infty$}\binom{\frac12}{i}\Big(\frac{16}{25C_1}\Big)^i\mbox{\Large$\Big|$}
\!\!\!\!\!\!\!\!\!\!\nonumber\\
\!\!\!\!\!\!\!\!\!\!\!\!\!\!\!\!\!\!\!\!&&\phantom{|B_1|^{-1}\ \  \ \ \ }\!\!\!\!\!\le
\ \ \ \frac58
\mbox{$\sum\limits_{i=1}^\infty$}(-1)^{i+1}\binom{\frac12}{i}\Big(\frac{16}{25|C_1|}\Big)^i
\stackrel{{}^{\sc\rm\eqref{C+ToSayas+1}\,(c)}}{\le}
\frac58
\mbox{$\sum\limits_{i=1}^\infty$}(-1)^{i+1}\binom{\frac12}{i}\Big(\frac{16(1+\eE_1^2)}{25}\Big)^i
\!\!\!\!\!\!\!\nonumber\\
\!\!\!\!\!\!\!\!\!\!\!\!\!\!\!\!\!\!\!\!&&
\phantom{\Big|AA\Big|}\ \ \
\stackrel{{}^{\sc\rm\eqref{bimeformo}}}{=}\ \
\frac58
\mbox{\Large$\Big($}1{\ssc}-{\ssc}\Big(1{\ssc}-{\ssc}
\frac{16(1+\eE_1^2)}{25}\Big)^{\frac12}\mbox{\Large$\Big)$}
=\frac14+\frac{\eE_1^2}{3}+O(\eE_1)^4<(1+\eE_1)^{-1},
\end{eqnarray}
 which is a contradiction with the first inequality of
    \eqref{C+ToSayas+1}\,(c).
This proves that we can
 only have \eqref{1=bbb}, which is equivalent to \eqref{1=bbb-remark}.
Exactly similar to the evaluation in \eqref{A1-latercase}, we can deduce from
the right-hand side of \eqref{1=bbb}\,(ii) the following,
\begin{eqnarray}\label{D1-is}
&\!\!\!\!\!\!\!\!\!\!\!\!\!\!&
|D_1|\stackrel{{}^{\sc\rm\eqref{1=bbb}\,(ii)}}{\ge}
\frac58
\mbox{\Large$\Big($}2-
\mbox{$\sum\limits_{i=1}^\infty$}(-1)^{i+1}\binom{\frac12}{i}\Big(\frac{16}{25|C_1|}\Big)^i
\mbox{\Large$\Big)$}%
\nonumber\\
&\!\!\!\!\!\!\!\!\!\!\!\!\!\!&\phantom{|D_1|}\ \,\ \stackrel{{}^{\sc\rm\eqref{bimeformo}}}{=}\ \
\frac58
\mbox{\Large$\Big($}1+\Big(1-\frac{16(1+\eE_1^2)}{25}\Big)^{\frac12}\mbox{\Large$\Big)$}
 =1-\frac{\eE_1^2}{3}+O(\eE_1)^4.
\end{eqnarray}
 Thus
we obtain from \eqref{1=bbb}\,(i) the following
\equa{B1b1b,s,s}{\dis|B_1|=|D_1|^{-1}
{\ssc}\stackrel{{}^{\sc\rm\eqref{D1-is}}}{\le}{\ssc}1{\ssc}+{\ssc}\frac{\eE_1^{2}}{3}+O(\eE_1)^4,
}
which in particular gives the first inequality of \eqref{ImMpP}\,(2).
}%
\NOUSE{
For convenience, we denote, 
\equa{DeD12}{\dis
{\rm(i)\ }D_2=\frac{X_2^{\ell_0-2}}{\widetilde X_1^2Z^{\ell_0}},
\ \ \ \ \ \ \ \
{\rm(ii)\ }D_3=
\frac{\widetilde X_1^{2} Z^{\ell_0+1}}{X_2^{\ell_0}}.
}
Then
\begin{eqnarray}
\label{aaasss--LetNSoOP----1}
&\!\!\!\!\!\!\!\!\!\!\!\!\!\!\!\!\!\!\!\!\!\!\!\!\!\!\!\!\!\!\!\!\!\!\!\!\!\!\!\!\!\!\!\!\!\!\!\!\!\!\!
&
{\rm(i)\ }A_3\stackrel{{}^{\sc\rm\eqref{LetNSoOP----1}\,(i)}}{=}
\frac{\widetilde X_1^2 Z^{\ell_0+1}}{X_2^{\ell_0}\Big(\frac15 +
\frac{4 X_2^{\ell_0-2}}{5\widetilde X_1^2 Z^{\ell_0}}\Big)}
\stackrel{{}^{\sc\rm\eqref{DeD12}}}{=}\frac{D_3}{\frac15+\frac{4D_2}{5}\Big)},
\!\!\!\!\!\!\!\!\!
\!\!\!\!\!\!\!\!\!
\!\!\!\!\!\!\!\!\!
\!\!\!\!\!\!\!\!\!
\nonumber\\
&\!\!\!\!\!\!\!\!\!\!\!\!\!\!\!\!\!\!\!\!\!\!\!\!\!\!\!\!\!\!\!\!\!\!\!\!\!\!\!\!\!\!\!\!\!\!\!\!\!\!\!
&
{\rm(ii)\ }
A_2\stackrel{{}^{\sc\rm\eqref{LetNSoOP----1}\,(iii)}}{=}
\frac{\widetilde X_1^2Z^{\ell_0}}{A_3^{10}X_2^{\ell_0-2}
\Big(2 - \d_0 - \frac{(1 - \d_0) X_2^{\ell_0-2}}{
\widetilde X_1^2 Z^{\ell_0}}\Big)}
\stackrel{{}^{\sc\rm\eqref{DeD12}}}{=}\frac1{A_3^{10}}D_2\Big(2-\d_0 -(1-\d_0)D_2\Big),
\!\!\!\!\!\!\!\!\!
\!\!\!\!\!\!\!\!\!
\!\!\!\!\!\!\!\!\!
\!\!\!\!\!\!\!\!\!
\nonumber\\
&\!\!\!\!\!\!\!\!\!\!\!\!\!\!\!\!\!\!\!\!\!\!\!\!\!\!\!\!\!\!\!\!\!\!\!\!\!\!\!\!\!\!\!\!\!\!\!\!\!\!\!
&
{\rm(iv)\ }
A_1^{-1}\stackrel{{}^{\sc\rm\eqref{aaasss--LetNSoOP----1}\,(i),\,(ii)}}{=}\frac{D_{1}^{-1}}{3-\frac2{A_3D_2}},\ \ \ \ \ \ \ {\rm(v)\ }
D_2^4\stackrel{{}^{\sc\rm\eqref{LetNSoOP----1}\,(iv),\,(v),\,\eqref{DeD12}\,(ii)}}{=}C_1B_1^{-2}.
\!\!\!\!\!\!\!\!\!
\end{eqnarray}
}
\NOUSE{Similar to \eqref{LetNSoOP----1-redefine++}, we have
\begin{eqnarray}
\label{{LetNSoOP----1-redefine++000}}
&\!\!\!\!\!\!\!\!\!\!\!\!\!\!\!\!\!\!\!\!\!\!\!\!\!\!\!\!\!\!\!\!\!\!
&
{\rm(i)\ }1\stackrel{{}^{\sc\rm\eqref{LetNSoOP----1}\,(iii),\,\eqref{A2-rewrite}\,(ii)}}{=}
\bB_2^{-1}
 \Big|\frac{2 +2 \d_0}{1 + 2 \d_0} -
\frac{B_1}{
1+2\d_0}\Big|,
\!\!\!\!\!\!\!\!\!\!\!\!\!\!\nonumber\\
&\!\!\!\!\!\!\!\!\!\!\!\!\!\!\!\!\!\!\!\!\!\!\!\!\!\!\!\!\!\!\!\!\!\!
&
{\rm(ii)\ }
\bB_2=
\frac
{\aA_2\aA_3^{\ell_0^2}\tilde \xX_1^{\ell_0^4(1+2\d_0-\d_0^2-3\d_0^3-2\d_0^4)}}
{\xX_2^{\ell_0^2(1+2\d_0+\d_0^2)}}
\equiv\frac{\aA_2}{\aA_3^{2\ell_0+1}\tilde\xX_1^{\ell_0^4(2+4\d_0-2\d_0^2-9\d_0^3-4\d_0^4)}}.
\!\!\!\!\!\!\!\!\!\!\!\!\!\!\nonumber\\
&\!\!\!\!\!\!\!\!\!\!\!\!\!\!\!\!\!\!\!\!\!\!\!\!\!\!\!\!\!\!\!\!\!\!
&{\rm(iii)\ }
\cC_2\stackrel{{}^{\sc\rm\eqref{A2-rewrite}\,(iii),\,\eqref{equa-Case6-lemm}\,(iv)}}{\equiv}
\frac{\aA_1\aA_2\aA_3^{\ell_0^2-1}\tilde \xX_1^{\ell_0^4(2+4\d_0-4\d_0^2-9\d_0^3-2\d_0^4)}}
{\xX_2^{\ell_0^2(2+4\d_0+1)}}\equiv
\frac{\aA_1\aA_2}{\aA_3^{\ell_0^2+4\ell_0+2}
\tilde X_1^{\ell_0^4(4+8\d_0-5\d_0^2-15\d_0^3-4\d_0^4)}}
,
\!\!\!\!\!\!\!\!\!\!\!\!\!\!\nonumber\\
&\!\!\!\!\!\!\!\!\!\!\!\!\!\!\!\!\!\!\!\!\!\!\!\!\!\!\!\!\!\!\!\!\!\!
&
{\rm(v)\ }
\aA_2\stackrel{{}^{\sc\rm\eqref{Related-A1-A2}}}{\le}
\frac{X_2^{\ell_0^2(1+2\d_0+\d_0^2)}}{A_3^{\ell_0^2}\widetilde X_1^{\ell_0^4(1+2\d_0-\d_0^2-3\d_0^3-2\d_0^4)}}
 \Big(\frac{2 +2 \d_0}{1 + 2 \d_0} -
\frac{A_1\widetilde X_1^{\ell_0^4(1+2\d_0-3\d_0^2-6\d_0^3)}}{
(1+2\d_0)A_3X_2^{\ell_0^2(1+2\d_0-2\d_0^2)}Z^2}\Big),
=
\tilde \xX_1^6 \xX_2^{-10}\Big(\frac15+
\frac{4}{5}\aA_1\aA_3^{-2}\tilde  \xX_1^{20\ell_0+5}\xX_2\Big)+O(\d)^2
\!\!\!\!\!\!\!\!\!\!\!\!\!\!\nonumber\\
&\!\!\!\!\!\!\!\!\!\!\!\!\!\!\!\!\!\!\!\!\!\!\!\!\!\!\!\!\!\!\!\!\!\!
&
\phantom{{\rm(iv)\ }}
\stackrel{{}^{\sc\rm\eqref{LetNSoOP----1-redefine++}\,(ii)}}{=}
\frac15
\aA_3^{-10}\tilde\xX_1^{-800\ell_0-494}
+\frac{4}{5}\aA_1\aA_3^{-11}\tilde \xX_1^{-700\ell_0-439}+O(\d)^2,
\!\!\!\!\!\!\!\!\!\!\!\!\!\!\nonumber\\
&\!\!\!\!\!\!\!\!\!\!\!\!\!\!\!\!\!\!\!\!\!\!\!\!\!\!\!\!\!\!\!\!\!\!
&
{\rm(vi)\ }
\bB_1^{-1}\stackrel{{}^{\sc\rm\eqref{LetNSoOP----1}\,(iv)}}{=}
\aA_1^{-1}\aA_3^2\tilde \xX_1^{-20\ell_0-5}\xX_2^{-70\ell_0-5}
\zZ^{70\ell_0+4}
=\aA_1^{-1}\aA_3^2\tilde \xX_1^{-20\ell_0-5}\xX_2^{-1}+O(\d)^2
\!\!\!\!\!\!\!\!\!\!\!\!\!\!\nonumber\\
&\!\!\!\!\!\!\!\!\!\!\!\!\!\!\!\!\!\!\!\!\!\!\!\!\!\!\!\!\!\!\!\!\!\!
&
\phantom{{\rm(vi)\ }\ \ \ \ }\ \ \ \ \
=\ \ \ \aA_1^{-1}\aA_3\tilde \xX_1^{-100\ell_0-55}+O(\d)^2
,
\!\!\!\!\!\!\!\!\!\!\!\!\!\!
\\[4pt]
&\!\!\!\!\!\!\!\!\!\!\!\!\!\!\!\!\!\!\!\!\!\!\!\!\!\!\!\!\!\!\!\!\!\!
&
{\rm(vii)\ }
\cC_1
=\aA_1\aA_3^{-3}\tilde \xX_1^{20\ell_0-5}\xX_2^{70\ell_0+5}\zZ^{-70\ell_0-4}
=\aA_1\aA_3^{-3}\tilde \xX_1^{20\ell_0-5}\xX_2+O(\d)^2
=\aA_1\aA_3^{-2}\tilde \xX_1^{100\ell_0+45}+O(\d)^2,
\!\!\!\!\!\!\!\!\!\!\!\!\!\!
\!\!\!\!\!\!\!\!\!\!\!\!\!\!
\nonumber\\
\nonumber&\!\!\!\!\!\!\!\!\!\!\!\!\!\!\!\!\!\!\!\!\!\!\!\!\!\!\!\!\!\!\!\!\!\!\!\!\!\!\!
&
{\rm(viii)\ }
\bB_1^{-100\ell_0-45}\cC_1^{100\ell_0+55}
=\aA_1^{10}\aA_3^{-100\ell_0-65}+O(\d)^2
\stackrel{{}^{\sc\rm\eqref{ImMpP}}}{\le}
\aA_1^{-\frac{200\ell_0}{3}(1+O(\d_0)^{\frac12})}+O(\d)^2<1+O(\d)^2.
\!\!\!\!\!\!\!\!\!\!\!\!\!\!\!\!\!\!\!\!\!\!\!\!\!\!\!\!
\end{eqnarray}
}\NOUSE
{%
\begin{lemm}\label{Ano-lemm}For $(p_1,p_2)\in\ol V_0$, we have
\begin{eqnarray}\label{ImMpP+1}\label{x1-x2==}
&\!\!\!\!\!\!\!\!\!\!\!\!\!\!\!\!\!\!
\!\!\!\!\!\!\!\!\!\!
&
{\rm(i)\ }
\bB_1<1+2\d_0+\d_0^2,\ \ \
\ \
{\rm(ii)\ }\dD_0
\stackrel{{}^{\sc\rm\eqref{More-de-1}\,(ii)}}{=}
\tilde \xX_1^{\ell_0^2(1+2\d_0)}\aA_3^{-\ell_0^{10}}\le1+O(\d)^1,
\ \ \ \ \ {\rm(iii)\ }
1+O(\d)^1\le\Big|\frac32 -\frac{ 1}{2 A_3\widetilde X_1^{10}}\Big|
\le2+O(\d)^1,\!\!\!\!\!\!\!\!\!\!\!
\nonumber\\
&\!\!\!\!\!\!\!\!\!\!\!\!\!\!\!\!\!\!
\!\!\!\!\!\!\!\!\!\!
&
{\rm(iv)\ }
1+O(\d)^1\le\Big|\frac75 - \frac{2}{5 A_3\widetilde X_1^{10}}\Big|
\le\frac95+O(\d)^1,
\ \ \ \ \
{\rm(v)\ }\d_2^{\d_0^2}<\tilde\xX_1<(1+\d)\aA_1\stackrel{{}^{\sc\rm\eqref{A1-A2-cond}\,(a)}}{\le}
\ell_1+O(\d)^1
,
\!\!\!\!\!\!\!\!\!\!\!\!\!\!\!\!\!\!\!
\nonumber\\
&\!\!\!\!\!\!\!\!\!\!\!\!\!\!\!\!\!\!
\!\!\!\!\!\!\!\!\!\!
&
 {\rm(vi)\ }\d_2<\xX_2<\ell_2
\NOUSE{,
\!\!\!\!
\!\!\!\!\!\!\!\!\!\!\!\!\!\!
\\[-0pt]\nonumber
\!\!\!\!\!\!\!\!\!\!\!\!
&\!\!\!\!\!\!\!\!\!\!\!\!\!\!\!\!\!\!\!\!\!\!\!\!\!\!\!\!\!\!\!\!\!\!\!\!\!\!\!\!\!\!\!\!\!
&
{\rm(v)\ }
\d_2^2<\tilde\xX_1<\ell_2^2,\ \ \ \
{\rm(vi)\ }\d_2^{100}<\Big|\frac43 - \frac{X_2^{84} Z^8}{3\widetilde X_1^{88}}\Big|<\ell_2^{100},
\ \ \ \
{\rm(vii)\ }\d_2^{100}<\Big|\frac13 + \frac{2 X_2^{84} Z^8}{3\widetilde X_1^{88}}\Big|
<\ell_2^{100}
}
.
\end{eqnarray}
\end{lemm}
\noindent{\it Proof.~}
If $\bB_1\ge3$, then 
by \eqref{g(r)}, we have
$g(\bB_1)\ge\bB_1^{\frac{2 \d_0}{20 + 11 \d_0}}(2\bB_1-3)>3$,
a contradiction with \eqref{g(r)}. This proves \eqref{x1-x2==}\,(i).
We have
\equa{Msms838383}{\dis
\aA_3^{-1}\tilde\xX_1^{-10}
\stackrel{{}^{\sc\rm\eqref{LetNSoOP----1}\,(iv),\,(v)}}{=}
\cC_1\bB_1^{-1}=
\bB_1^{-\frac{20+9\d_0}{20+11\d_0}}\cC_1\bB_1^{-\frac{2\d_0}{20+11\d_0}}
\tilde\xX_1\stackrel{{}^{\sc\rm\eqref{Formmeme},\,\eqref{hB1}}}{\le}1+O(\d)^1.
}
This proves  \eqref{x1-x2==}\,(ii). From this and the well-known formula $\pm(|a|-|b|)\le|a-b|\le|a|+|b|$ for $a,b\in\C$,
we obtain \eqref{x1-x2==}\,(iii),\,(iv).
If $\tilde\xX_1\le\d_2^{\d_0^2}$ (recalling that $\d_2\ll\d_0$), then by \eqref{LetNSoOP----1-redefine++}\,(iii),\,\eqref{A1-A2-cond}\,(a),\,\eqref{x1-x2==}\,(iii),
 we see that
$\aA_3\ll1$, a contradiction with \eqref{ImMpP}. Thus we have the first inequality of
\eqref{x1-x2==}\,(v).
Assume the second inequality of  \eqref{x1-x2==}\,(v) does not hold.
Then $\tilde\xX_1=(1+\d)\aA_1$ by \eqref{C+LetNSoOP}\,(iv).
Up to $O(\d)^2$, we have
\begin{eqnarray}
&&\!\!\!\!\!\!\!\!\!\!\!\!\!\!\!\!\!\!\!\!\!\!\!\!\!\!\!\!\!\!
\label{msm3474575}
1\ \ \ \stackrel{{}^{\sc\rm\eqref{LetNSoOP----1-redefine++}\,(v)}}{\le}
\frac15
(1+\d)^{-800\ell_0-494}\aA_2^{-1}\aA_3^{-10}\aA_1^{-800\ell_0-494}
+\frac{4}{5}(1+\d)^{-700\ell_0-439}\aA_2^{-1}\aA_1\aA_3^{-11}\aA_1^{-700\ell_0-439}
\!\!\!\!\!\!\!\!\!\!\!\!\!\!\!\!\!\!\!\!\nonumber\\
&&\!\!\!\!\!\!\!\!\!\!\!\!\!\!\!\!\!\!\!\!\!\!\!\!\!\!\!\!\!\!
\ \ \ \stackrel{{}^{\sc\rm\eqref{cont-B2},\,\eqref{ImMpP}}}{\le}
\frac15
(1+\d)^{-800\ell_0-494}\aA_1^{-800\ell_0(1+O(\d_0)^{\frac12})}
+\frac{4}{5}(1+\d)^{-700\ell_0-439}\aA_1^{-700\ell_0(1+O(\d_0)^{\frac12})}
\!\!\!\!\!\!\!\!\!\!\!\!\!\!\!\!\!\!\nonumber\\
&&\!\!\!\!\!\!\!\!\!\!\!\!\!\!\!\!\!\!\!\!\!\!\!\!\!\!\!\!\!\!
\phantom{1}
\ \ \
\stackrel{{}^{\sc\rm\eqref{A1-A2-cond}\,(a)}}{\le}
\frac15
(1+\d)^{-800\ell_0-494}
+\frac{4}{5}(1+\d)^{-700\ell_0-439}
=1-720\ell_0\Big(1+O(\d_0)^1\Big)\d
,
\end{eqnarray}
which is a contradiction. This proves \eqref{x1-x2==}\,(v).
Finally, by \eqref{A1-A2-cond}\,(b),\,\eqref{x1-x2==}\,(iv),\,(v) and the second equality of
\eqref{LetNSoOP----1-redefine++}\,(iv), we obtain \eqref{x1-x2==}\,(vi).\hfill$\Box$\vskip7pt
%
%
%
%
%
\NOUSE{
We already have the first inequity of \eqref{ImMpP+1}\,(ii) by \eqref{ImMpP}. Assume
$\aA_3\ge\ell_2^{\frac1{11}+\d_1}$. Then
by \eqref{C+ToSayas+1}\,(d),\,\eqref{A1-A2-cond}, we have
$\aA_2^{-44}\aA_3^{44}\xX_2^4\ge\d_1^{\frac{66}{25}}
\ell_2^{4+44\d_1}\d_2^4>\ell_2^{43\d_1}\gg\ell_1$ [recalling from \eqref{MSmde33333} that $\ell_2=\d_2^{-1}\gg\ell_1=\d_1^{-1}\gg1$],
and we obtain from 
\eqref{LetNSoOP----1+ThisOne}\,(iv) that $\aA_3\ll1$, a contradiction. This proves \eqref{ImMpP+1}\,(iii).
We have,
\begin{eqnarray}
\label{memxx2}
&&\!\!\!\!\!\!\!\!\!\!\!\!\!\!\!\!\!\!\!\!\!\!\!\!\!\!\!\!\!\!\!\!\!\!\!\!\!\!\!
{\rm(i)\ }\xX_2\ \ \ \ \ \ \ \ \stackrel{{}^{\sc\rm\eqref{LetNSoOP----1+ThisOne}\,(v)}}{ =}\ \
(|B_1|\aA_1^{-1}\aA_2^{-\frac{89}{2}}\aA_3^{\frac{87}{2}})^{\frac14}{\sc}+{\sc}O(\d)^2
\!\!\!\!\!\!\!\!\!\!\!\!\!\!\!\!\!\!\!\!\!\!\!\!\!\!\!\nonumber\\
&&\!\!\!\!\!\!\!\!\!\!\!\!\!\!\!\!\!\!\!\!\!\!\!\!\!\!\!\!\!\!\!\!\!\!\!\!\!\!\!
\ \ \ \ \ \ \ \ \stackrel{{}^{\sc\rm\eqref{cont-B2},\,\eqref{ImMpP},\,\eqref{ImMpP+1}\,(i)}}{\ge}
\Big(\frac13\aA_1^{-1-\frac{267}{100}+\frac{609}{100}}\Big)^{\frac14}{\sc}+{\sc}O(\d)^2{\sc}>{\sc}\d_2,
\!\!\!\!\!\!\!\!\!\!\!\!\!\!\!\!\!\!\!\!\!\!\!\!\!\!\!\nonumber\\
&&\!\!\!\!\!\!\!\!\!\!\!\!\!\!\!\!\!\!\!\!\!\!\!\!\!\!\!\!\!\!\!\!\!\!\!\!\!\!\!
{\rm(ii)\ }\xX_2
\stackrel{{}^{\sc\rm\eqref{cont-B2},\,\eqref{Formmeme},\,\eqref{ImMpP+1}\,(iii)}}{\le}
\Big(\aA_1^{-1-\frac{267}{100}}\ell_2^{\frac{87}{22}+\frac{87\d_1}{2}}\Big)^{\frac14}+O(\d)^1<\ell_2,
\end{eqnarray}
where we have used the equality of (i) in the first inequality of (ii).
This proves \eqref{ImMpP+1}\,(iv). By \eqref{LetNSoOP----1+ThisOne}\,(ii),\,\eqref{ImMpP+1}\,(iii),\,(iv), we obtain
\eqref{ImMpP+1}\,(v).
Finally, by \eqref{LetNSoOP----1}\,(i),\,(iii),\,\eqref{A1-A2-cond}\,(a),\,\eqref{ImMpP+1}\,(iii), we obtain
\eqref{ImMpP+1}\,(vi),\,(vii).
}\NOUSE
{
Now if $\xX_2\le\d_2$, then $\zZ\le\d_2+O(\d)^2$ by \eqref{equa-Case6-lemm}\,(iv) and
$\big|\ell_0+1-\frac{\ell_0}{X_2^{\ell_0}}\big|\ge\ell_0\xX_2^{-\ell_0}-\ell_0-1
\ge\ell_0\ell_2^{\ell_0}-\ell_0-1\gg\ell_2$ by \eqref{A1-A2-cond}\,(b) and the fact that $\ell_2=\d_2^{-1}\gg\ell_0=\d_0^{-1}$,
and so $\aA_2\ll1$ by \eqref{equa-Case6-lemm}\,(iv),\,\eqref{LetNSoOP----1}\,(i). We obtain a contradiction with \eqref{ImMpP}\,(i). This prove
the first inequality of \eqref{x1-x2==}\,(iii). If $\xX_2\ge\ell_2$, then
$\big|\ell_0+1-\frac{\ell_0}{X_2^{\ell_0}}\big|\ge\ell_0+1-\ell_0\d_2^{\ell_0}\ge\ell_0$, and
$\aA_3\ll1$ by \eqref{LetNSoOP----1}\,(i),\,\eqref{equa-Case6-lemm}\,(iv).
We obtain a contradiction again.
This proves \eqref{x1-x2==}\,(iii).
}\NOUSE{
By \eqref{equa-Case6-lemm}\,(iv),\,\eqref{(T0o(eE)1},\,\eqref{LetNSoOP----1}\,(i),
up to $O(\d)^2$, we have $X_2=Z=A_3\widetilde X_1^{-1}$. Thus
by \eqref{LetNSoOP----1}\,(ii),\,(iii), we have
\begin{eqnarray}
\label{anoth-LetNSoOP----1}
&\!\!\!\!\!\!\!\!\!\!\!\!\!\!\!\!\!\!\!\!\!\!\!\!\!\!\!\!\!\!\!\!\!\!\!\!\!\!\!\!\!\!\!\!\!\!\!\!\!\!\!
&
{\rm(i)\ }
A_1=\frac{ X_2^9}{A_3^{20}\widetilde  X_1^{5}}\Big(\frac15 + \frac{4 A_3^{10}\widetilde X_1^5}{5}\Big)
=
\frac1{A_3^{11}\widetilde X_1^{14}}\Big(\frac15+\frac{4A_3^{10}\widetilde X_1^5}{5}\Big)+O(\d)^2
,
\nonumber\\
&\!\!\!\!\!\!\!\!\!\!\!\!\!\!\!\!\!\!\!\!\!\!\!\!\!\!\!\!\!\!\!\!\!\!\!\!\!\!\!\!\!\!\!\!\!\!\!\!\!\!\!
&
{\rm(ii)\ }
A_2=\frac{A_3^{115}\widetilde X_1^{11}}{X_2^{100}}
=A_3^{15}\widetilde X_1^{111}+O(\d)^2
.
\end{eqnarray}
By \eqref{LetNSoOP----1}\,(iii), up to $O(\d)^1$, we have
\begin{eqnarray}
\label{11TiX1==}
\!\!\!\!\!\!\!\!\!\!\!\!\!\!\!\!\!\!\!\!\!\!\!\!\!\!\!&&
{\rm(i)\ }
\aA_2\stackrel{{}^{\sc\rm\eqref{LetNSoOP----1}\,(iii)}}{=}\aA_3^{115}\tilde\xX_1^{11}\xX_2^{100}
\stackrel{{}^{\sc\rm\eqref{C1B1}\,(i)}}{=}
\aA_3^{15}\tilde\xX_1^{111},\ \ \implies\ \ {\rm(ii)\ }
 \tilde\xX_1=\aA_2^{\frac1{111}}\aA_3^{-\frac{15}{111}}\mbox{ \ [up to $O(\d)^1\ssc\,$]},
\nonumber\\
\!\!\!\!\!\!\!\!\!\!\!\!\!\!\!\!\!\!\!\!\!\!\!\!\!\!\!&&
{\rm(iii)\ }
|B_1|\stackrel{{}^{\sc\rm\eqref{C1B1}\,(ii)}}{=}
\aA_1\aA_3\tilde\xX_1^9
\stackrel{{}^{\sc\rm\eqref{11TiX1==}\,(ii)}}{=}
\aA_1\aA_2^{\frac3{37}}\aA_3^{-\frac8{37}}
\end{eqnarray}
If $|B_1|\le\frac12$, then $r=|B_1|^{-1}\ge2$ and by \eqref{g(r)}, $g(r)= r^{\frac8{17}}(4r-5)>3$, a contradiction.
This proves \eqref{x1-x2==}\,(ii).
Assume the first inequality of \eqref{x1-x2==}\,(ii) does not hold, then by \eqref{C+ToSayas+1}\,(d), we have
\equa{equality==}{\dis
\tilde\xX_1=(1-\d)\aA_1^{-\frac15}.
}
Then up to $O(\d)^2$, we have
\begin{eqnarray}
\label{Connnn}
&&\!\!\!\!\!\!\!\!\!\!\!\!\!\!\!\!\!\!\!\!\!\!\!\!\!\!\!
1\ \ \ \ \ \ \stackrel{{}^{\sc\rm\eqref{LetNSoOP----1}\,(iii)}}{\le}\ \ \
\aA_2^{-1}\aA_3^{45}\xX_2^{-45}\tilde\xX_1^6\Big(\frac16+\frac56\aA_1\aA_3^{-114}\tilde\xX_1^{18}\xX_2^{66}\Big)
 \stackrel{{}^{\sc\rm\eqref{X2===}}}{=}
\aA_2^{-1}\tilde\xX_1^{51}\Big(\frac16+\frac56\aA_1\aA_3^{-48}\tilde\xX_1^{-48}\Big)\!\!\!\!\!\!\!\!\!\!\!\!\!\!\!\!
\nonumber\\&&\!\!\!\!\!\!\!\!\!\!\!\!\!\!\!\!\!\!\!\!\!\!\!\!\!\!\!
\ \ \stackrel{{}^{\sc\rm\eqref{cont-B2},\,\eqref{ImMpP},\,\eqref{equality==}}}{\le}
(1-\d)^{51}\aA_1^{-\frac{18}{41}-\frac{51}5}\Big(\frac16+\frac56(1-\d)^{-48}
\aA_1^{1-\frac{48}{246}+\frac{48}{5}}\Big)
\nonumber\\&&\!\!\!\!\!\!\!\!\!\!\!\!\!\!\!\!\!\!\!\!\!\!\!\!\!\!\!
\ \ \ \ \ \ \ \ \ \ \ \ =\ \ \ \ \ \ (1-\d)^{51}\aA_1^{-\frac{2181}{205}}\Big(\frac16+
\frac56(1-\d)^{-48}\aA_1^{\frac{2549}{246}}\Big)
=\frac16(1-\d)^{51}\aA_1^{-\frac{21815}{205}}+\frac56(1-\d)^3\aA_1^{-\frac{48}{205}}
\nonumber\\&&\!\!\!\!\!\!\!\!\!\!\!\!\!\!\!\!\!\!\!\!\!\!\!\!\!\!\!
\ \ \ \ \ \ \ \ \ \ \stackrel{{}^{\sc\rm\eqref{A1-A2-cond}}}{\le}\ \ \ \ \ \
1-11\d\mbox{ \ \ [up to $O(\d)^2$]},
\end{eqnarray}
a contradiction, which proves the first inequality of \eqref{x1-x2==}\,(ii).
Now up to $O(\d)^1$, we have
\begin{eqnarray}
\label{tilde-x1}
&&\!\!\!\!\!\!\!\!\!\!\!\!\!\!\!\!\!\!\!\!\!\!\!\!\!\!\!
{\rm(i)\ }
\tilde\xX_1\stackrel{{}^{\sc\rm\eqref{C1B1}\,(i)}}{=}
(|B_1|\aA_1^{-1}\aA_3^{-102})^{\frac1{102}}+O(\d)^1
\stackrel{{}^{\sc\rm\eqref{A1-A2-cond}\,(a),\,\eqref{ImMpP},\,\eqref{Formmeme}}}{\le}1+O(\d)^1,
\nonumber\\
&&\!\!\!\!\!\!\!\!\!\!\!\!\!\!\!\!\!\!\!\!\!\!\!\!\!\!\!
{\rm(ii)\ }\xX_2\stackrel{{}^{\sc\rm\eqref{X2===}}}{=}\aA_3\tilde\xX_1^{-1}+O(\d)^1\stackrel{{}^{\sc\rm\eqref{ImMpP},\,\eqref{tilde-x1}\,(i)}}{\ge}1+O(\d)^1,
\nonumber\\
&&\!\!\!\!\!\!\!\!\!\!\!\!\!\!\!\!\!\!\!\!\!\!\!\!\!\!\!
{\rm(iii)\ }\xX_2\stackrel{{}^{\sc\rm\eqref{C1B1}\,(i)}}{=}
\big(|B_1|\aA_1^{-1}\aA_3^{-36}\tilde\xX_1^{-168}\big)^{\frac1{66}}+O(\d)^1\stackrel{{}^{\sc\rm\eqref{C+ToSayas+1}\,(c),\,\eqref{ImMpP},\,\eqref{Formmeme}}}{\le}
\big((1-\d)\aA_1^{-\frac18}\big)^{-\frac{28}{11}}+O(\d)^1\!\!\!\!\!\!\!\!\!\!\!\!\!\!\!\!\!\!\!\!\!\!
\nonumber\\
&&\!\!\!\!\!\!\!\!\!\!\!\!\!\!\!\!\!\!\!\!\!\!\!\!\!\!\!
\phantom{{\rm(iii)\ }\ \ \ }\stackrel{{}^{\sc\rm\eqref{cont-B2}\,(a)}}{<}\ell_2
.
\end{eqnarray}
This proves \eqref{x1-x2==}\,(ii),\,(iii).
}
}%
\NOUSE{%
Now note from \eqref{LetNSoOP----1}\,(i),\,(ii) that we have
\equa{A2-inv}{\dis
 A_2^{-1} \tilde X_1^{-2\ell-1} X_2^{10\ell}=
\alpha_1^{-4\ell} A_3^{-200\ell}Z^{-10\ell-\lL}=A_3^{-\frac{1605\ell}{8}-\frac{\lL}{16}}+O(\eE_1)^1,}
  so that we have (i) below,
\begin{eqnarray}
\label{A2-B1aaaaaa}
&\!\!\!\!\!\!\!\!\!\!\!\!\!\!\!\!\!\!\!\!\!\!\!\!\!\!\!\!\!\!&
 {\rm(i)\ }\aA_2^{-1}\aA_3^{\frac{1605\ell}{8}+\frac{\lL}{16}}
 \tilde\xX_1^{-2\ell-1}\xX_2^{10\ell}\stackrel{{}^{\sc\rm\eqref{A2-inv},\,\eqref{(T0o(eE)1}}}{=
}1+O(\eE_1)^1,\ \ \
\\\nonumber
&\!\!\!\!\!\!\!\!\!\!\!\!\!\!\!\!\!\!\!\!\!\!\!\!\!\!\!\!\!\!&
{\rm(ii)\ }\frac12+O(\eE_1)^1\stackrel{{}^{\sc\rm\eqref{x1-x2==}\,(ii)}}{\le}
\aA_1\aA_3^{101\lL+\frac{101}{16}}\tilde \xX_1^{808}\aA_2^{-\lL+1}
\stackrel{{}^{\sc\rm\eqref{C1B1}}}{=}|B_1|+O(\eE_1)^1
\stackrel{{}^{\sc\rm\eqref{Formmeme}}}{\le}1+O(\eE_1)^1.
\!\!\!\!\!\!\!\!\!\!\!\!\!\!\!\!
\end{eqnarray}
\NOUSE{where the last inequality of (ii) is obtained from \eqref{1=bbb} by noting the following
\equa{D1====}{\dis\!\!\!\!\!\!|B_1|^{-1}\stackrel{{}^{\sc\rm\eqref{1=bbb}\,(i)}}{=}
|D_1|\stackrel{{}^{\sc\rm\eqref{1=bbb}\,(ii),\,\eqref{ImMpP}\,(2)}}{\le}\frac58\Big(1+\big(1+\frac{16(1+\eE_1^2)}{25}\big)^{\frac12}\Big)
=\frac18\Big(5{\sc\!}+{\sc\!}\sqrt{41}\Big){\sc\!}+{\sc\!}O(\eE_1)^2.\!\!\!}
}%
To prove the last inequalities
of \eqref{x1-x2==}\,(iii),\,(iv), we need to evaluate upper bounds for $\tilde\xX_1,\xX_2$.  This is  done in the next
\eqref{SMSmene},\,\eqref{X2==mmmmm}, using \eqref{A2-B1aaaaaa}. The upper bounds
  are powers of $\aA_1$ [which is $\le\nn$ by \eqref{A1-A2-cond}$\ssc\,$]  by  some exponents  independent of $
\nn_0$.  Thus
the last inequalities
of \eqref{x1-x2==}\,(iii),\,(iv) hold.
\NOUSE{
To prove the last inequalities
of \eqref{ImMpP}\,(3),\,(4), we need to evaluate upper bounds for $\tilde\xX_1,\xX_2$, which are done in
\eqref{SMSmene},\,\eqref{X2==mmmmm}. We remark that using \eqref{A2-B1aaaaaa}, one can easily obtain that $\tilde\xX_1,\xX_2$ are upper bounded by some
number which is a power of $\aA_1$ such that the power is independent of $\nn_0$, thus
the last inequalities
of \eqref{ImMpP}\,(3),\,(4) hold. In this sense, readers do not need to check \eqref{SMSmene},\,\eqref{X2==mmmmm} carefully. However, to be precise, we give
explicit evaluations in \eqref{SMSmene},\,\eqref{X2==mmmmm}.
}
The last inequality of  \eqref{A2-B1aaaaaa}\,(ii) shows that, up to $O(\eE_1)^1$,
we have [noting again that in order to use \eqref{ImMpP-(1)}, it is required that the power of $\aA_3$ is non-positive],
\begin{eqnarray}
\label{SMSmene}
&\!\!\!\!\!\!\!\!\!\!\!\!\!\!\!\!\!\!\!\!\!\!\!\!\!\!\!\!\!\!\!\!\!\!\!\!&
\tilde\xX_1\stackrel{{}^{\sc\rm\eqref{A2-B1aaaaaa}\,(ii)}}{\le}
\aA_1^{-\frac1{808}}\aA_2^{\frac{\lL-1}{808}}\aA_3^{-\frac{\lL(16+\dD)}{128}}
\stackrel{{}^{\sc\rm\eqref{cont-B2},\,\eqref{a2a1a3},\,\eqref{ImMpP-(1)}}}{\le}
\aA_1^{-\frac1{808}+\frac{\lL-1}{808}(2\dD+\frac{2130\dD^2}{101})-\frac{\lL(16+\dD)}{128}(\frac{\dD}{101}
 + \frac{1931 \dD^2}{10201})+O(\dD)^2}\!\!\!\!\!\!\!\!\!\!\!\!\!\!\!\!\!\!\!\!\!\!\!\!\!\!\!\!
\nonumber\\
&\!\!\!\!\!\!\!\!\!\!\!\!\!\!\!\!\!\!\!\!\!\!\!\!\!\!\!\!\!\!\!\!\!\!\!\!&
\ \ \ \ \ \ \ \ 
=\ \ \ \aA_1^{-\frac{149\dD}{1305728}+O(\dD)^2}
\mbox{ \ [up to $O(\eE_1)^1$]}.
\end{eqnarray}
Using this in \eqref{A2-B1aaaaaa}\,(i), we obtain, up to $O(\eE_1)^1$
(recalling that $0<\dD=\lL^{-1}\ll\d=\ell^{-1}$),
\begin{eqnarray}
\label{X2==mmmmm}
\!\!\!\!\!\!\!\!\!
\!\!\!\!\!\!\!\!\!
\!\!\!\!\!\!\!\!\!
\!\!\!\!\!\!\!\!\!
&&
\xX_2\ \ \ \ \ \stackrel{{}^{\sc\rm\eqref{A2-B1aaaaaa}\,(i)}}{=}
\aA_2^{\frac{\d}{10}}\aA_3^{-(\frac{\d\lL}{160}+\frac{321}{16})}\tilde\xX_1^{\frac15+\frac{\d}{10}}
\nonumber\\
\!\!\!\!\!\!\!\!\!
\!\!\!\!\!\!\!\!\!
\!\!\!\!\!\!\!\!\!
\!\!\!\!\!\!\!\!\!
&&
\ \
\stackrel{{}^{\sc\rm\eqref{cont-B2},\,\eqref{SMSmene},\,\eqref{ImMpP-(1)}}}
{\le}
\aA_1^{\frac1{10}(2\dD+\frac{2130\dD^2}{101})-(\frac{\d\lL}{160}+\frac{321}{16})(\frac{\dD}{101}
 + \frac{1931 \dD^2}{10201})+(\frac15+\frac{\d}{10})(-\frac{149\dD}{1305728})+O(\dD)^2}\!\!\!\!\!\!\!\!\!\!\!\!\!\!\!\!\!\!
\nonumber\\
\!\!\!\!\!\!\!\!\!
\!\!\!\!\!\!\!\!\!
\!\!\!\!\!\!\!\!\!
\!\!\!\!\!\!\!\!\!
&& \ \ \ \ \
\ \  \ \ \ \ 
\ \ = \ \ \ \ \ \ \aA_1^{-\frac{\d}{16160}+O(\dD)^1} \mbox{ \ [up to $O(\eE_1)^1$]}
.\!\!\!\!\!\!\!\!\!
\end{eqnarray}
Then \eqref{A1-A2-cond} with \eqref{SMSmene},\,\eqref{X2==mmmmm} shows that the last inequalities
of \eqref{x1-x2==}\,(iii),\,(iv) hold.
To prove the first inequality of \eqref{x1-x2==}\,(iii), let us
assume conversely, \equa{xx1smallerth}{\tilde\xX_1\le\eE_0.}
Recall from \eqref{MSmde33333}
that $\nn_0=\eE_0^{-1}\gg\nn$ and we can also require that $\nn_0^{\dD}\gg\nn$, therefore
$\tilde\xX_1^{-\frac{8\dD}{16+\dD}}\ge\nn_0^{\frac{8\dD}{16+\dD}}\gg\nn$, which with  \eqref{x1-x2==}\,(ii) gives (i) below, where the second equality of (i)
follows from
the equality in \eqref{A2-B1aaaaaa}\,(ii), and where the first equality of (ii) follows from
the first equality of \eqref{X2==mmmmm}. Thus up to $O(\eE_1)^1$, we have
\begin{eqnarray}
\label{mmmmoooo}
&\!\!\!\!\!\!\!\!\!\!\!\!\!\!\!\!\!\!&
{\rm(i)\ }\zZ\stackrel{{}^{\sc\rm\eqref{LetNSoOP----1}\,(i)}}{=}
\aA_3^{\frac1{16}}\stackrel{{}^{\sc\rm\eqref{A2-B1aaaaaa}\,(ii)}}{=}
\Big(|B_1|
\aA_1^{-1}\tilde \xX_1^{-808}\aA_2^{\lL+1}\Big)^{\frac1{101(16\lL+1)}}
%
%
\stackrel{{}^{\sc\rm\eqref{A1-A2-cond},\,\eqref{xx1smallerth},\,\eqref{x1-x2==}\,(ii)}}{\gg}1,
\nonumber\\
&\!\!\!\!\!\!\!\!\!\!\!\!\!\!\!\!\!\!&
{\rm(ii)\ }\xX_2\stackrel{{}^{\sc\rm\eqref{X2==mmmmm}}}{=}
\aA_2^{\frac{\d}{10}}\aA_3^{-(\frac{\d\lL}{160}+\frac{321}{16})}\tilde\xX_1^{\frac15+\frac{\d}{10}}
\stackrel{{}^{\sc\rm\eqref{A1-A2-cond}\,(b),\,\eqref{xx1smallerth},\,\eqref{mmmmoooo}\,(i)}}{\ll}1,
\nonumber\\
&\!\!\!\!\!\!\!\!\!\!\!\!\!\!\!\!\!\!&
{\rm(iii)\ }\xX_1\stackrel{{}^{\sc\rm\eqref{x1-small}\,(i)}}{<}
\xX_2^{\frac15}
\tilde \xX_1^{\frac1{100}}\stackrel{{}^{\sc\rm\eqref{xx1smallerth},\,\eqref{mmmmoooo}\,(ii)}}{\ll}1.
\end{eqnarray}
%
%
Now assume $\xX_2\le\eE_0$.
We have $\tilde\xX_1\le1+O(\eE_1)^1$ by \eqref{SMSmene}. Thus $\xX_1<1$ exactly as in \eqref{mmmmoooo}\,(iii).
Using the facts that $\xX_2\le\eE_0$ and $\tilde\xX_1>\eE_0\ge\xX_2$ (which is just proven), we obtain from the equality of
\eqref{mmmmoooo}\,(ii) the following [noting from \eqref{MSmde33333} that we can choose $\nn_0=\eE_0^{-1}$ such that $\nn_0^{\dD}\gg1$ and observing that
$(\frac{\d\lL}{160}+\frac{321}{16})^{-1}=\frac{160\ell\dD}{1+3210\ell\dD}$],
\equa{aa3-is-bigg}{\dis\!\!\!\!\!\!\!\!\!\!
\aA_3\stackrel{{}^{\sc\rm\eqref{mmmmoooo}\,(ii)}}{=}
(\aA_2^{\frac{\d}{10}}\xX_2^{-1}\tilde\xX_1^{\frac15+\frac{\d}{10}})^{\frac{160\ell\dD}{1+3210\ell\dD}}>
(\aA_2^{\frac1{10}}\xX_2^{-\frac45+\frac{\d}{10}})^{\frac{160\ell\dD}{1+3210\ell\dD}}
\stackrel{{}^{\sc\rm\eqref{A1-A2-cond}\,(b)}}{\ge}\xX_2^{-(\frac45-\frac{\d}{10})\frac{160\ell\dD}{1+3210\ell\dD}}
\gg1.\!\!\!\!\!\!\!\!}
Thus $\zZ\gg1$ by \eqref{LetNSoOP----1}\,(i). In particular, $\zZ>1$ and $\xX_1<\zZ,
\,\xX_2<\zZ$. Once again, we have \eqref{Conndndm}, which is a contradiction. This proves
\eqref{x1-x2==}\,(iv).
Finally, by \eqref{LetNSoOP----1}\,(i),\,\eqref{ImMpP-(1)},\,\eqref{A1-A2-cond}\,(a) and notation \eqref{MSmde33333}, we see that the first inequality of
\eqref{x1-x2==}\,(v) holds.
If $\zZ\ge\nn_0$, then $\zZ\gg1$ and by  \eqref{x1-x2==}\,(iii),\,(iv),\,\eqref{x1-small}\,(i), we  see that $\zZ>\xX_1$, $\zZ>\xX_2$, and we obtain a contradiction exactly as in \eqref{Conndndm}.
This proves \eqref{x1-x2==}\,(v).
\NOUSE{
Similarly,
\eqref{LetNSoOP----1}\,(ii) says that $A_2=\frac{X_1}{ A_3^{2- 2 \d_0^2}}\big(\frac{1 - \d_0^4}{2}
  +\frac{ (1 + \d_0^4)A_3^{2-\d_0^{11}}}{2 X_1^2}\big)$,
which can be easily \vspace*{-5pt}rewritten as
\equa{X1===-1}
{\dis\tilde a X_1^{-2}+\tilde b X_1^{-1}+\tilde c=0,
\mbox{ \ where }\tilde a=\frac{(1+\d_0^4)A_3^{2\d_0^2-\d_0^{11}}}{2A_2},\ \ \
\tilde b=-1,\ \ \
\tilde c=\frac{1-\d_0^4}{2A_2A_3^{2-2\d_0^2}}.
}
Recalling from \eqref{LetNSoOP----1}\,(v),\,\eqref{C1-2and}\,(ii) that
$A_2=\frac{A_2}{ A_3^{2 \d_0^2-\d_0^{11}}},\,C_2=\frac{1}{A_2^2A_3^{2- 4 \d_0^2+\d_0^{11}}}$, we obtain
\equa{tilde-a-ab}
{\dis
\frac{\tilde b}{2\tilde a}\stackrel{{}^{\sc\rm\eqref{X1===-1}}}{=}-\frac{ A_2}{1+\d_0^4},\ \ \ \ \ \ \ \
4\tilde a\tilde c\tilde b^{-2}\stackrel{{}^{\sc\rm\eqref{X1===-1}}}{=}(1-\d_0^8)C_2.
}
Thus using the same arguments after \eqref{abc====}, we can deduce
\NOUSE{
observe that we can rewrite  \eqref{LetNSoOP----1}\,(i),\,(ii),\,(iv)--(vii) as follows,
\begin{eqnarray}
\label{Rewritt}
&\!\!\!\!\!\!\!\!\!\!\!\!\!\!\!\!\!\!\!\!\!\!\!\!\!\!\!\!\!\!&
{\rm(i)\ }A_1'\stackrel{{}^{\sc\rm\eqref{LetNSoOP----1}\,(i)}}{=}\a_1\Big(\frac{1-\d_0^2}{2}+\frac{(1+\d_0^2)A_1'^2}{2}\b_1\Big), \ \ A_1':=\frac{X_1Z^{1+\d_0^2}}{X_2},\ \
{\rm(ii)\ }\a_1:=\frac{1}{A_1A_2^{\d_0^2}},
\  \ {\rm(iii)\ }\b_1:=\frac1{A_2^2},\!\!\!\!\!\!\!\!\!\!\!\!\!\!\!\!\!\!\!\!\!\!
\nonumber\\[0pt]&\!\!\!\!\!\!\!\!\!\!\!\!\!\!\!\!\!\!\!\!\!\!\!\!\!\!\!\!\!\!&
 {\rm(iv)\ }A_1\stackrel{{}^{\sc\rm\eqref{LetNSoOP----1}\,(iv)}}{=}(\a_1\b_1)^{-1},\  \ \ \ \ {\rm(v)\ }C_1\stackrel{{}^{\sc\rm\eqref{LetNSoOP----1}\,(v)}}{=}\a_1^2\b_1,
\nonumber\\[4pt]&\!\!\!\!\!\!\!\!\!\!\!\!\!\!\!\!\!\!\!\!\!\!\!\!\!\!\!\!\!\!&
{\rm(vi)\ }\tildeX_1^{-1}\stackrel{{}^{\sc\rm\eqref{LetNSoOP----1}\,(ii)}}{=}\a_2\Big(\frac{1-\d_0^4}{2}+\frac{(1+\d_0^4)X_1^{-2}}{2}\b_2\Big),
\ \ \ \ {\rm(vii)\ }\a_2:=\frac1{A_2A_3^{2-2\d_0^2}},
\ \ \ \ {\rm(viii)\ }\b_2:=A_3^2,\!\!\!\!\!\!\!\!\!\!\!\!\!\!\!\!\!\!\!\!\!
\nonumber\\&\!\!\!\!\!\!\!\!\!\!\!\!\!\!\!\!\!\!\!\!\!\!\!\!\!\!\!\!\!\!&
{\rm(ix)\ }A_2\stackrel{{}^{\sc\rm\eqref{LetNSoOP----1}\,(vi)}}{=}(\a_2\b_2)^{-1},\ \ \ \ \ \ {\rm(x)\ }C_2\stackrel{{}^{\sc\rm\eqref{LetNSoOP----1}\,(vii)}}{=}\a_2^2\b_2.
\end{eqnarray}
This enables us to rewrite \eqref{Rewritt}\,(i),\,(vi) as the following
(here is the reason we define $A_1,C_1,A_2,$ $C_2$),
where we always regard the element $\big(1-(1-\d_0^4)C_1\big)^{\frac12}$, respectively
$\big(1-(1-\d_0^8)C_2\big)^{\frac12}$,
 as the unique element, defined by the formula  \eqref{bimeformo}, which is a power series of
 $C_1$, respectively $C_2$, converging absolutely by condition \eqref{ToSayas+1}\,(d),
\begin{eqnarray}
\!\!\!\!\!\!\!\!\!\!\!\!\!\!\!\!\!\!\!\!\!\!\!\!\!\!\!\!\!\!\!&&
\label{re-B2}
{\rm(i)\ }\mbox{\Large$\Big($}A_1'-\frac{A_1}{1+\d_0^4}\Big(1+\big(1-(1-\d_0^4)C_1\big)^{\frac12}\Big)\mbox{\Large$\Big)$}
\mbox{\Large$\Big($}A_1'-\frac{A_1}{1+\d_0^4}\Big(1-\big(1-(1-\d_0^4)C_1\big)^{\frac12}\Big)\mbox{\Large$\Big)$}=0,
\\\nonumber
\!\!\!\!\!\!\!\!\!\!\!\!\!\!\!\!\!\!\!\!\!\!\!\!\!\!\!\!\!\!\!&&
{\rm(ii)\ }
\mbox{\Large$\Big($}\tildeX_1^{-1}-
\frac{ A_2}{1+\d_0^8}\Big(1+\big(1-(1-\d_0^8)C_2\big)^{\frac12}\Big)\mbox{\Large$\Big)$}
\mbox{\Large$\Big($}\tildeX_1^{-1}-\frac{A_2}{1+\d_0^8}\Big(1-\big(1-(1-\d_0^8)C_2\big)^{\frac12}\Big)\mbox{\Large$\Big)$}
=0.\!\!\!\!\!\!\!\!\!\!\!\!\!\!\!
\end{eqnarray}
\eqref{re-B2}\,(ii) we obtain
}%
\begin{eqnarray}
\label{1=bbb+1}
&\!\!\!\!\!\!\!\!\!\!\!\!\!\!\!\!\!\!\!\!\!\!\!\!\!\!\!\!\!\!&
{\rm(i)\, }\tildeX_1^{-1}{\ssc}={\ssc}A_2D_2, \
\\&\!\!\!\!\!\!\!\!\!\!\!\!\!\!\!\!\!\!\!\!\!\!\!\!\!\!\!\!\!\!&
\nonumber{\rm(ii)\, }D_2=\frac1{1+\d_0^4}\mbox{\Large$\Big($}1+\Big(1-(1-\d_0^8)C_2\Big)^{\frac12}\mbox{\Large$\Big)$}\stackrel{{}^{\sc\rm\eqref{bimeformo}}}{=}\frac1{1+\d_0^4}
\mbox{\Large$\Big($}2+
\mbox{$\sum\limits_{i=1}^\infty$}\binom{\frac12}{i}\Big(-(1-\d_0^8)C_2\Big)^i\mbox{\Large$\Big)$},\!\!\!\!\!\!\!\!\!\!\!\!\!
\end{eqnarray}
and the last inequality of \eqref{ImMpP}\,(8).
}\NOUSE{Now we write condition \eqref{ToSayas+1}\,(f) as follows (we add a parameter $\l$, cf.~Remark \ref{ABCD12}),
\equa{C1C2}{\dis(1-\l\d^4)|C_1|\le1\le\frac{1+\l\d^4}{|C_2|}.}
By \eqref{1=bbb}\,(ii),\,\eqref{1=bbb+1}\,(ii), we have
\begin{eqnarray}
\label{111AMsnsns}
\!\!\!\!\!\!\!\!\!\!\!\!\!\!\!\!\!\!\!\!&&
{\rm(i)\ }|D_1|\ge
\frac58\Big(1+\big(1-\frac{16}{25|C_1|}\big)^{\frac12}\Big)
\ge\frac58\Big(1+\big(1-\frac{16}{25}(1-\l\d^4)^{-1}\big)^{\frac12}\Big)=1-\frac{\l\d^4}{3}+O(\d)^8,
\nonumber\\
\!\!\!\!\!\!\!\!\!\!\!\!\!\!\!\!\!\!\!\!&&
{\rm(ii)\ }|D_2|
\ge
\frac1{2(1{\ssc}-{\ssc}\d_0)}\Big(1{\ssc}+{\ssc}\big(1{\ssc}-
{\ssc}4 \d_0 (1{\ssc} -{\ssc} \d_0)(1+\l\d^4)\big)^{\frac12}\Big)
=1-\frac{\d_0\l\d^4}{1-2\d_0}+O(\d)^8.
\end{eqnarray}
Using the above in \eqref{1=bbb}\,(i),\,\eqref{1=bbb+1}\,(i),
 we can deduce
[we also need to use \eqref{cont-B2},\, \eqref{ImMpP}\,(6)--(10)$\ssc\,$],
\begin{eqnarray}
\label{A1A2Must}
\!\!\!\!\!\!\!\!\!\!\!\!\!\!\!\!\!\!\!\!&&
{\rm(i)\ }
1\le E_1:=\frac{|D_1|}{1-\frac{\l\d^4}{3}}+O(\d)^3
=\Big|\frac{(1+\frac{\l\d^4}{3})X_1^{80}}{A_1 A_2^{29} A_3^{74}}\Big|+O(\d)^3
=\frac{(1+\frac{\l\d^4}{3})+O(\d)^3}{|A_1 A_2^{29} A_3^{154} X_2^{80}|},
\nonumber
\\
\!\!\!\!\!\!\!\!\!\!\!\!\!\!\!\!\!\!\!\!&&
 {\rm(ii)\ }
 1\le E_2:=\frac{|D_2|}{1-\frac{\d_0\l\d^4}{1-2\d_0}}+O(\d)^3
 =\Big|\frac{\big(1-\frac{\d_0\l\d^4}{1-2\d_0}\big)X_2^{48}}{A_2 A_3^5}\Big|+O(\d)^3.
\!\!\!\!\!\!
\end{eqnarray}
From this and \eqref{cont-B2},\,\eqref{ImMpP}\,(3),\,(6), we can deduce
\begin{eqnarray}
\label{A1A2Must+}
\!\!\!\!\!\!\!\!\!\!\!\!\!\!\!\!\!\!\!\!\!\!\!\!\!\!\!\!\!&&
|C_1|=\Big|\frac{A_2^{87}A_3^{172}}{A_1^2X_1^{190}}\Big|+O(\d)^3\le
E_1^{\frac{19}{8}}|A_1^{-2+\frac{203(113-\d_0)}{1170}} A_3^{362} X_2^{190}|+O(\d)^3
\le\frac{1+\frac{19\l\d^4}{24}+O(\d)^3}{|A_1|^{\frac{581}{5616}+O(\d_0)^1}}
\nonumber\\
\!\!\!\!\!\!\!\!\!\!\!\!\!\!\!\!\!\!\!\!\!\!\!\!\!\!\!\!\!&&
\phantom{|C_2|}
<(1-\l\d^4)^{-1},
\end{eqnarray}
i.e., we have the last inequality of \eqref{ImMpP}\,(12) when $\l=1$.
From \eqref{A1A2Must+}, one sees that the power of $E_1$ plays a crucial role in obtaining the last inequality
(as long as the power is small than $3$, the last strict inequality holds).
By \eqref{ImMpP}\,(10), we have, 
\begin{eqnarray}
\label{A1A2Must+11}
\!\!\!\!\!\!\!\!\!\!\!\!\!\!\!\!\!\!\!\!\!\!\!\!\!\!\!&&
|C_2|=\frac1{|A_2^2 A_3^9|}=\frac{1}{|A_1|^{\frac{7}{3510}+O(\d_0)^1}}+O(\d)^3<
1+\d^4
,
\end{eqnarray}
i.e., we have the first inequality of \eqref{ImMpP}\,(12).
}%
}\NOUSE{
Now assume  equality occurs in the first inequality of \eqref{ToSayas+1}\,(f), i.e., $1=(1-\d^4)|A_1'|$.
Finally using
\eqref{ImMpP}\,(3),\,(5),\,(6), we obtain
\equan{Finana}{\dis\!\!\!\!\!\!\!\!\!\!\!\!
|\tildeX_1^{22\ell_0+3} X_2^{14}A_3^{-2}|
{\ssc\!}={\ssc\!}|\tildeX_1^{26\ell_0+3}A_3^{-2}{\ssc\!}+{\ssc\!}O(\d)^3
{\ssc\!}={\ssc\!}|A_3|^{-20-\frac{3\d_0}{2}}{\ssc\!}+{\ssc\!}O(\d)^3
{\ssc\!}<{\ssc\!}|A_1|^{\frac5{137}(1-\d_0)(20+\frac{3\d_0}{2})}{\ssc\!}+{\ssc\!}O(\d)^3.
\!\!\!\!\!\!\!\!\!\!\!\!\!\!\
}
From this and \eqref{ImMpP}\,(8),\,\eqref{cont-B2}  with the fact that $|A_1|\ge1$, we obtain \eqref{ImMpP}\,(12)
as follows
,
\begin{eqnarray}
\label{MSMA2}
&\!\!\!\!\!\!\!\!\!\!\!\!\!\!\!\!\!\!\!\!\!\!\!\!\!\!\!\!\!\!&
|A_2|\le|A_1^{\frac5{137}(1-\d_0)(20+\frac{3\d_0}{2})-2+33\frac5{137}}
|
+O(\d)^3
=|A_1|^{-\frac9{137}+O(\d_0)^1}+O(\d)^3<1+\d^3.
\end{eqnarray}
}
\NOUSE{To prove \eqref{ImMpP}\,(8),
we already have that $X_1,H_0$ are $1+O(\d)^1$ elements and so is $H_1$  by \eqref{MSMSH1H20}\,(ii).
Using the facts that $|X_2|\le1+O(\d)^5$ and $|A_1|=|A_1|+O(\d)^5\ge1+O(\d)^5$, we obtain from
\eqref{C1-2and}\,(i),\,(ii) that $A_1,X_2$ must be $1+O(\d)^1$ elements. Then $A_1$ is a $1+O(\d)^5$ element by
\eqref{MSMSH1H20}\,(i). Then $A_2,A_2$ are  $1+O(\d)^1$ elements by \eqref{cont-B2},\,\eqref{MSMSH1H20}\,(iii).
}\NOUSE{
 recall that we have notation \eqref{MSmde33333} and \eqref{A1andA2}.
First assume $|X_1|\le\d_2$. Then $|X_1|=\d_2$ by \eqref{C+LetNSoOP}\,(iv) and $|X_1'|=\ell_0^4+1+O(\d_2)^1$ by
\eqref{ImMpP}\,(8). Then \equa{|X-2====|}
{\d_2^{\d_1}<|X_2|<\ell_2^{\d_1},}
which is obtained from \eqref{A1andA2}\,(i),\,\eqref{ImMpP}\,(2) as follows:
if $|X_2|\le\d_2^{\d_1}$ then we see from \eqref{ImMpP}\,(2)
that the second term inside the brackets is $\gg1$ by noting that
$\ell_2^{\d_1}\gg\ell_1$ and $|A_1|\gg|X_2|^{-\ell_0^2}\gg\ell_1$, a contradiction with \eqref{A1andA2}\,(i);
if $|X_2|>\ell_2^{\d_1}$ then $|A_1|\gg|X_2|^{\ell_0^3}\gg\ell_1$, a contradiction again.
Now by \eqref{|X-2====|}, we see from  \eqref{ImMpP}\,(3) that
the second term inside  the brackets is $\gg1$ and $|A_2|\gg\ell_1^{3\ell_0^4}$, a contradiction with
\eqref{A1andA2}\,(ii).
Now assume $|X_1|\ge\ell_2$.
Then $|X_1|=\d_2$ by \eqref{C+LetNSoOP}\,(iv) and $|X_1'|=\ell_0^4\ell_2+O(\d_2)^1$ by
\eqref{ImMpP}\,(8). Exactly similar to the proof of \eqref{|X-2====|}, we obtain from
\eqref{ImMpP}\,(2) the following,
\equa{|X-2====|+}
{\d_2^{\d_1}<\Big|\frac{X_2^{\ell_0^4}}{X_1'^{\ell_0}}\Big|<\ell_2^{\d_1}.}
Then as above we can obtained a contradiction from  \eqref{ImMpP}\,(3). This proves that $\d_2<|X_1|<\ell_2$.
Similarly, we can prove $\d_2<|X_2|<\ell_2$ (see Remark \ref{mekr-beww} below). This proves \eqref{ImMpP}\,(6).
\begin{rema}\rm\label{mekr-beww}
In fact if we use Notation \ref{nota-sim}, then we always obtain from \eqref{ImMpP}\,(2) the following,
\equa{MSMSMSMSMMM}{\dis
\Big|\frac{X_2^{\ell_0^4}}{X_1'^{\ell_0}}\Big|\sim_{\ell_2\ssc\,}1.}
Then we can obtain a contradiction from \eqref{ImMpP}\,(3) if either of $|X_1|,|X_2|$ is not $\sim_{\ell_2\ssc\,}1$.
\end{rema}
To prove \eqref{ImMpP}\,(7), assume $|X_1'|\le\d_2$. Then \eqref{ImMpP}\,(8) shows that $|X_1|=1+O(\d_0)^4$.
Then we obtain from \eqref{ImMpP}\,(3) that $|X_2|\sim_{\ell_2\ssc\,}1$. Then \eqref{ImMpP}\,(2) shows that
$|A_1|\succ_{\ell_2\ssc\,}1$, a contradiction. Thus $\d_2<|X_1'|$.
Finally assume  $|X_1'|\ge(1+\d)|X_2|^{\ell_0^3}$. Then \eqref{C+LetNSoOP}\,(v) gives
\equa{MEMEMMEmememe}{\dis
|X_1'|=(1+\d)|X_2|^{\ell_0^3}.}
We have
\begin{eqnarray}
\label{MFINN}
&\!\!\!\!\!\!\!\!\!\!\!\!\!\!\!\!\!\!\!\!\!\!\!\!\!\!&
1\stackrel{{}^{\sc\rm\eqref{A1andA2}\,(i)}}{\le}|A_1|
\stackrel{{}^{\sc\rm\eqref{ImMpP}\,(2),\,\eqref{MEMEMMEmememe}}}{\le}
(1+\d)^{-\ell_0}\Big(\d_0+(1-\d_0)(1+\d)^{\ell_0(1+\d_0)}\Big)
\nonumber\\&\!\!\!\!\!\!\!\!\!\!\!\!\!\!\!\!\!\!\!\!\!\!\!\!\!\!&
\ \ \ \ \ \ =\ \ \ 1-(\ell_0-\d_0+\d_0^2)\d+O(\d)^2<1,
\end{eqnarray}
a contradiction.
}\NOUSE{first we have the following,
where the first equality follows from the fact in \eqref{LetNSoOP----1}\,(iv) \vspace*{-0pt}that $A_1=A_1 A_2^2 A_3^{\d_0^2}$ and
$\frac{X_2}{X_1Z^{1+\d_0^2}}=\frac{1}{X_1X_2^{\d_0^2}}+O(\d)^3=\frac{A_3^{\d_0^2}}{X_1}+O(\d)^3$ by \eqref{ImMpP}\,(2),\,(6), while the two inequalities follows from \eqref{C+LetNSoOP}\,(iii), and  the last equality from definition \eqref{LetNSoOP----1}\,(v).
\begin{eqnarray}\label{1-dSmallerThan}
\!\!\!\!\!\!\!\!\!\!\!\!\!\!\!\!\!\!\!\!\!\!\!\!\!\!\!\!\!\!\!\!&&
(1-\d)\Big|\frac{A_1A_3^{2\d_0^2}}{A_2^6X_1}\Big|+O(\d)^3\ \,
\stackrel{{}^{\sc \rm\eqref{LetNSoOP----1}\,(iv),\,\eqref{ImMpP}\,(2),\,(6)}}{
=}\ \, (1-\d)\Big|\frac{A_1X_2}{X_1Z^{1+\d_0^2}}\Big|
\stackrel{{}^{\sc\rm \eqref{C+LetNSoOP}\,(iii)}}
{\le}1
\nonumber\\
\!\!\!\!\!\!\!\!\!\!\!\!\!\!\!\!\!\!\!\!\!\!\!\!\!\!\!\!\!\!\!\!\!\!\!&&
\phantom{(1-\d)\Big|\frac{A_1A_2^2A_3^{2\d_0^2}}{X_1}\Big|+O(\d)^3\ \ \ \ \ }\ \ \ \ \
\stackrel{{}^{\sc\rm\eqref{C+LetNSoOP}\,(iii)}}{\le}
\ \
\frac{1+\d}{|A_2X_1|}
\stackrel{{}^{\sc\rm\eqref{LetNSoOP----1}\,(v)}}{
=}(1+\d)\Big|\frac{A_3^{2\d_0^2-\d_0^{11}}}{A_2X_1}\Big|+O(\d)^3.
\end{eqnarray}
Using \eqref{MSM1111}
and convention \eqref{MSmde33333+}, by
\eqref{+C+LetNSoOP}, 
we see that
$|A_{\rOnE}|\sim_{\ell_2\ssc\,}1\sim_{\ell_2\ssc\,}|A_2|$.
Thus by \eqref{1-dSmallerThan}, we obtain \equa{x1xxxxxxx}{\mbox{ ${\rm(i)\ }
\Big|\frac{A_3^{2\d_0^2}}{X_1}\Big|\stackrel{{}^{\sc\rm\eqref{1-dSmallerThan}}}{\preceq_{\ell_2\ssc\,}}1
\stackrel{{}^{\sc\rm\eqref{1-dSmallerThan}}}{\preceq_{\ell_2\ssc\,}}\Big|\frac{A_3^{2\d_0^2}}{X_1}\Big|,
\mbox{ \ \ \ i.e., \ \ (ii) }
|X_1|
{\sim_{\ell_2\ssc\,}}|A_3|^{2\d_0^2}$.}} Then
\eqref{ImMpP}\,(4) gives the following,
\equa{AAAA22222}{\dis
1\sim_{\ell_2\ssc\,}|A_1|\stackrel{{}^{\sc\rm\eqref{ImMpP}\,(4),\,\eqref{x1xxxxxxx}\,(ii) }}{\preceq_{\ell_2\ssc\,}}\frac1{|A_3|^{2\d_0^2}}\Big(\frac{1-\d_0^2}{2}+
\frac{(1+\d_0^2)|A_3|^{2\d_0^2}}{2}\Big)=
\frac{1-\d_0^2}{2|A_3|^{2\d_0^2}}+
\frac{1+\d_0^2}{2}
,} which implies that $|A_3|\preceq_{\ell_2\ssc\,}1$. However by
\eqref{ImMpP}\,(3), we have $|A_3|\succeq_{\ell_2\ssc\,}1$. This with \eqref{ImMpP}\,(6) gives that
 $|X_2|\sim_{\ell_2\ssc\,}
|A_3|\sim_{\ell_2\ssc\,}1$, and so $|X_1|\sim_{\ell_2\ssc\,}1$ by \eqref{x1xxxxxxx} and \eqref{MSmde33333}.
In particular, we have
}\NOUSE{
By the last inequality of \eqref{LetNSoOP}\,(iii), we have $|X_1|\preceq_{\ell_2\ssc\,}1$,
then by \eqref{ImMpP}\,(5),
$|X_2|\preceq_{\ell_2\ssc\,}1$.
However
by \eqref{ImMpP}\,(3),\,(6), we have $|X_2|\succeq_{\ell_2\ssc\,}1$. Thus $|X_2|\sim_{\ell_2\ssc\,}1$, and then
 $|\tildeX_1|\sim_{\ell_2\ssc\,}1$ by \eqref{ImMpP}\,(5).
In particular
}
}
%
\NOUSE{By \eqref{A1A2Must}, we obtain
\begin{eqnarray}
\label{A1A2Must+1}
\!\!\!\!\!\!\!\!\!\!\!\!\!\!\!\!\!\!\!\!&&
{\rm(i)\ }
|X_1|\ge|A_1A_2^{29}A_3^{74}|^{\frac1{80}}+O(\d)^3
\ge|A_1^{\frac{13499}{3510}+O(\d_0)^1}|^{\frac1{80}}+O(\d)^3>1-\d,
\nonumber
\\
\!\!\!\!\!\!\!\!\!\!\!\!\!\!\!\!\!\!\!\!&&
 {\rm(ii)\ }
|X_2|\ge |A_2A_3^5|^{\frac1{48}}+O(\d)^3\ge|A_1^{-\frac{14}{585}}|^{\frac1{48}}+O(\d)^3>(1-\d)|A_1|^{-1}
.
\!\!\!\!\!\!
\end{eqnarray}
}
\NOUSE{Now assume  equality occurs in the last inequality of
of \eqref{LetNSoOP}\,(iii), i.e., $|X_1|=1+
\d$. By
\eqref{ImMpP}\,(4) with \eqref{ImMpP}\,(3),\,(6), we immediately obtain%
,
\equa{MEmeme}
{\dis\!\!\!\!\!\!\!
1\le\frac{(1+\d)^{-2\ell_0^4}}{|A_1A_3^{3\ell_0^3+\d_0}|}
\Big(\frac{1+\d_0^2}{2(1+\d)^{\ell_0^2}}+\frac{(1-\d_0^2)(1+\d)^{\ell_0^2}}{2}\Big)+O(\d)^3
\le\frac{1-(2\ell_0^4+1)\d}{|A|^{2\ell_0(1+O(\d_0)^1)}}+O(\d)^3
,
}
which is smaller than $1$,
a contradiction.
}%
}%
%

\NOUSE{
Further, by {\rm\eqref{TaKa},\,\eqref{SimMMSMS},\,\eqref{ImMpP}\,(2.b) and \eqref{LetNSoOP}\,(vi)},
  we have \eqref{-EiathA0}\,(2) as follows [recalling from \eqref{MSmde33333} that $0<\d\ll\d_2\ll\d_0\ll1{\ssc\,}$],
\begin{eqnarray}\label{x2+yyyy2}
&\!\!\!\!\!\!\!\!\!\!\!\!\!\!\!\!\!\!\!\!\!\!\!\!\!\!\!\!&
  |x_2+y_2|
\stackrel{{}^{\rm\eqref{TaKa},\,\eqref{SimMMSMS}}}{=}
 |\widetilde Z(\bar x_2+\bar y_2)|\stackrel{{}^{\rm\eqref{TaKa},\,\eqref{SimMMSMS}}}{=}
  \g_{\kk,\kk}|Z|^{\d_0^2}\stackrel{{}^{\rm\eqref{ImMpP}\,(2.b)}}{=}
\g_{\kk,\kk}\big|X_2+O(\d)^3\big|^{\d_0^2}
\nonumber\\&\!\!\!\!\!\!\!\!\!\!\!\!\!\!\!\!\!\!\!\!&  \phantom{|x_2+y_2|}\ \ \, \,=
    \g_{\kk,\kk}\big(|X_2|^{\d_0^2}+O(\d)^3\big)
\stackrel{{}^{\rm\eqref{LetNSoOP}\,(vi)}}{\ge}\g_{\kk,\kk}\big(\d_2^{\d_0^2}+O(\d)^3\big)\stackrel{{}^{\sc\rm\eqref{MSmde33333}}}{>}0.
   \end{eqnarray}
}

\NOUSE{
For later convenience, we denote, for all possible $i$,
\equa{DDeeXX}{\dis
\aA_i=|A_i|,\ \ \ \ \xX_i=|X_i|,\ \ \ \ \yY=|Z|.
}
Now assume  equality occurs in the first inequality of \eqref{LetNSoOP}\,(v), i.e.,
$\xX_1=(1-\d)\aA_1^{\ell_0^2}$.
By \eqref{ImMpP}\,(5) (and replacing $\xX_2$ by $\xX_1^{-1}\bA_3^{-1}$),
and using \eqref{cont-B2},\,\eqref{A1A2Must}\,(ii), we have,
where thoughout the rest proof of this lemma, we conduct computations up to $O(\d)^3$, thus $O(\d)^3$ is omitted;
also note  that
$\aA_1\le1\le\aA_3$ by \eqref{LetNSoOP}\,(i),\,\eqref{ImMpP}\,(3),
\begin{eqnarray}
\label{X1Smallerthan-A2}
&\!\!\!\!\!\!\!\!\!\!\!\!\!\!\!\!\!\!&
5{\ssc\!}\le
{\ssc\!}\aA_1^{-1}\aA_3^{-88}(1{\ssc\!}+{\ssc\!}4\aA_3^{85}\xX_1^{-5})E_2^{\frac{1}{11}+\d_0}
{\ssc\!}={\ssc\!}\xX_1^{22\d_0}\aA_1^{-\frac{699}{814}+O(\d_0)^1}
\Big(\xX_1^{5+O(\d_0)^1}\aA_3^{-85+O(\d_0)^1}{\ssc\!}+{\ssc\!}4\aA_3^{\frac{364\d_0}{11}+O(\d_0)^2}
\Big)\!\!\!\!\!\!\!
\nonumber\\
&\!\!\!\!\!\!\!\!\!\!\!\!\!\!\!\!\!\!&\phantom{
5}
\le
\xX_1^{21\d_0}\aA_1^{-\frac{699}{814}+O(\d_0)^1}
\Big(\xX_1^{5+O(\d_0)^1}\aA_3^{-85+O(\d_0)^1}{\ssc\!}+{\ssc\!}4\xX_1^{\d_0}\aA_3^{\frac{364\d_0}{11}+O(\d_0)^2}
\Big(\aA_1^{-1}\aA_3^{-88}(1{\ssc\!}+{\ssc\!}4\aA_3^{85}\xX_1^{-5})\Big)^{\frac{365\d_0}{11\times 88}}
\Big),
\end{eqnarray}
Note that $f(\aA_3)$ is a decreasing function of $\aA_3$ when $\aA_3\le\aA_0:=(21\ell_0)^{\frac{1}{84+\d_0}}=
(21\ell_0)^{\frac1{84}}\big(1+O(\d_0)^1\big)$.
Therefore if $\aA_3\le\aA_0$, then $f(\aA_3)\le f(\aA_1^{-\frac{1 + \d_0}{37 + 19 \d_0}})$ by
\eqref{ImMpP}\,(3), and then the right-hand side of \eqref{X1Smallerthan-A2}\,(i) is smaller than $5$, a contradiction.
Thus $\aA_3>\aA_0$, and we obtain that $f(\aA_3)=4\aA_3^{4\d_0}+O(\d_0)^1$ by \eqref{X1Smallerthan-A2}\,(ii).
Thus \eqref{X1Smallerthan-A2}\,(i) gives
\equa{bB3isBigger}
{\dis
\aA_3\ge\Big(\frac54+O(\d_0)^1\Big)^{\frac{\ell_0}{4}}\aA_1^{-\frac{21\ell_0^2}{4}}.
}
Thus by the first inequalities of \eqref{A1A2Must+},\,\eqref{A1A2Must+11}, we see that $|C_1|,|C_2|<\d_0$.
Then $|D_1|=\frac54+O(\d_0)^1,$ $|D_2|=\frac{4}{3}+O(\d_0)^1$ by
\eqref{1=bbb}\,(ii),\,\eqref{1=bbb+1}\,(ii). We obtain from
\eqref{1=bbb}\,(i),\,\eqref{1=bbb+1}\,(i) the following
[noting from \eqref{A1A2Must} that $\dD_1:=|D_1|=|E_1|,\,\dD_2:=|D_2|=|E_2|$ up to $O(\d)^4\ssc\,$],
\begin{eqnarray}
\label{A1A2Must+++}
&\!\!\!\!\!\!\!\!\!\!\!\!\!\!\!\!\!\!&
{\rm(i)\ }
\frac54+O(\d_0)^1=\dD_1=\aA_1^{-1}\aA_3^{-1}\xX_2^2,
\ \ \ \
{\rm(ii)\ }
\frac{4}{3}+O(\d_0)^1=\dD_2=\aA_1^{\frac{115}{74}+O(\d_0)^1}\aA_3^{11+\d_0}\xX_2^{-22},\ \ \
\nonumber\\
&\!\!\!\!\!\!\!\!\!\!\!\!\!\!\!\!\!\!&
{\rm(iii)\ }
1=\xX_1\xX_2\aA_3+O(\d)^1=\xX_2\bA_1^{\ell_0^2}\aA_3+O(\d)^1
.\end{eqnarray}
Finally assume  equality occurs in the first inequality of \eqref{LetNSoOP}\,(iii), i.e., $|X_2|=(1-\d)|A_1|^{\frac1{30}}$.
Similarly to \eqref{X1Smallerthan-A2-1},
by \eqref{ImMpP}\,(4) [this time we replace $\tildeX_1$ by
$(X_2A_3)^{-1}\ssc\,$], we have
\begin{eqnarray}
\label{X1Smallerthan-A2-1}
&\!\!\!\!\!\!\!\!\!\!\!\!\!\!\!\!\!\!&
1\le
\frac{|X_2|^{24}}{|A_1A_3^{517}|}
\Big(\frac15 +\frac{4}{5}|A_3|^{\frac{571}{5}}\Big)
+O(\d)^3
\le
(1-\d)^{24}|A_1|^{\frac{25173}{3385}}
\Big(\frac15 + \frac{4}{5
|A_1|^{\frac{5710}{677}}}\Big)\le1-24\d+O(\d)^3,
\end{eqnarray}
a contradiction.
By \eqref{A1A2Must}, we obtain
\begin{eqnarray}
\label{X1Smallerthan}
&\!\!\!\!\!\!\!\!\!\!\!\!\!\!\!\!\!\!&
|\tildeX_1|<(1+\d^3)|A_1^{-\frac{9}{89}}A_3^{-51}|^{\frac1{42}}+O(\d)^3
<(1+\d^3)|A_1^{-\frac{9}{89}+51\times\frac{25}{801}}|^{\frac1{42}}=(1+\d^3)|A_1|^{\frac{199}{5607}},
\nonumber\\&\!\!\!\!\!\!\!\!\!\!\!\!\!\!\!\!\!\!&
|\tildeX_1|>(1-\d^3)|A_1^{\frac{20}{89}-10\times\frac{25}{801}}|^{\frac1{10}}=
=(1-\d^3)|A_1|^{-\frac{7}{801}}.
\end{eqnarray}
This proves \eqref{ImMpP}\,(13).
Finally assume  equality occurs in the inequality of \eqref{LetNSoOP}\,(iii), i.e., $|X_2|=
(1-\d)|A_1|^{-100}$.
By \eqref{ImMpP}\,(5) and using \eqref{ImMpP}\,(3),\,(6),\,\eqref{cont-B2},\,\eqref{A1A2Must}, we obtain
\begin{eqnarray}
\label{X1Smallerthan--}
&\!\!\!\!\!\!\!\!\!\!\!\!\!\!\!\!\!\!&
1\le E_2^{\frac{21}{10}}|\tildeX_1^2X_2^2A_3^2|
\frac{|A_3|^9}{|A_1^{\frac{20}{89}}\tildeX_1^{23}X_2|}\Big(\frac14{\ssc\!}+{\ssc\!}
\Big|\frac{3\tildeX_1^{32}}{4A_3^{20}}\Big|E_1^{\frac{16}{21}}\Big){\ssc\!}+{\ssc\!}O(\d)^3
{\ssc\!}={\ssc\!}\frac{|X_2|}{|A_1^{\frac{62}{89}}A_3^{10}|}\Big(\frac15{\ssc\!}+{\ssc\!}
\frac{1}{|A_1^{\frac{48}{623}}
A_3^{\frac{412}{7}}|}\Big){\ssc\!}+{\ssc\!}O(\d)^3\!\!\!\!\!\!\!\!\!
\nonumber\\&\!\!\!\!\!\!\!\!\!\!\!\!\!\!\!\!\!\!&\phantom{1}
\le\frac{(1-\d)|A_1|^{-100}}
{|A_1^{\frac{62}{89}-10\times\frac{25}{801}}|}
\Big(\frac15+\frac{1}{|A_1^{\frac{48}{623}-\frac{412}{7}\times\frac{25}{801}}|}\Big)+O(\d)^3
\le1-\d+O(\d)^3<1,\end{eqnarray}
which is a contradiction.
This completes the proof of \eqref{ImMpP}\,(14).
}\NOUSE
{%
Assume the last inequality of \eqref{LetNSoOP}\,(ii) occurs, i.e., $|\tildeX_1|=(1+\d)|A_1|$.
By \eqref{ImMpP}\,(6),\,\eqref{A1A2Must}, we deduce from
\eqref{ImMpP}\,(5) the following, where for convenience, we denote $\xX_1=|\tildeX_1|,\,\xX_2=|X_2|,\,
\bA_i=|A_i|$, $i=1,2,3$, and we  also use
\eqref{cont-B2} and the fact from \eqref{ImMpP}\,(3) that $\bA_3^{-1}\le1+O(\d)^3$,
\begin{eqnarray}
\label{CCCaaaaaaa}
&\!\!\!\!\!\!\!\!\!\!\!\!\!\!\!\!\!\!\!\!\!\!\!\!\!\!\!&
4\le \Big(\frac{\xX_2\bA_3}{\xX_1}\Big)^{\frac{62}{69}}
\Big(\frac{\xX_1^{33}\xX_2^{21}}{\bA_1\bA_3^{600}}\Big)^{\frac{1}{207}}
\frac{\bA_3^2}{\bA_2\xX_1^{10}\xX_2}\Big(1+\frac{3\xX_1^{12}}{\bA_3^4}\Big)+O(\d)^3
\le\bA_1^{-\frac1{207}}\bA_2^{-1}\xX_1^{-\frac{247}{23}}(1+3\xX_1^{12})+O(\d)^3
\!\!\!\!\!\!\!\!\!\!\!\!\!\!\!\nonumber\\
&\!\!\!\!\!\!\!\!\!\!\!\!\!\!\!\!\!\!\!\!\!\!\!\!\!\!\!&
\phantom{4}
=(1+\d)^{-\frac{247}{23}}\bA_1^{-1+\frac7{83}-\frac{247}{23}}
(1+3(1+\d)^{12}\bA_1^{12})+O(\d)^3
=f(\bA_1)-\frac{712}{23}\d+O(\d)^3,
\end{eqnarray}
where $f(\bA_1):=\bA_1^{-\frac{183143}{17181}}(1+3\bA_1^{12})$
is a decreasing function of $\bA_1$ at $\bA_1=1$ (i.e., $\frac{df}{d\bA_1}|_{\bA_1=1}<0$).
Since $\bA_1\ge1$ and when $\bA_1=1$, \eqref{CCCaaaaaaa} cannot hold, we must have $\bA_1>1$. We deduce that
$\bA_1>1.200$ [noting that when $\bA_1\le\bA_0:=\big(\frac{183143}{69087}\big)^{\frac1{12}}\simeq 1.084$,
$f(\bA_1)$
is a decreasing function, and when $\bA_1\ge\bA_0$ it is a increasing function, and we have
$f(1.200)\simeq4\times0.993<1$, thus we must have $\bA_1>1.200$].
Then using the fact that $\bA_3=\xX_1\xX_2^{-1}=\bA_1\xX_2^{-1}+O(\d)^3$, we obtain
from \eqref{A1A2Must}\,(i) the following,
\begin{eqnarray}
\label{CCCaaaaaaa+}
&\!\!\!\!\!\!\!\!\!\!\!\!\!\!\!\!\!\!\!\!\!\!\!\!\!\!\!&
{\rm(i)\ }1+O(\d)^1
<|A_1''|=\xX_2^{54+567}\bA_1^{-1-567}+O(\d)^1=
\xX_2^{621}\bA_1^{568}+O(\d)^1, \ \mbox{ i.e., }\ \xX_2\ge\bA_1^{\frac{568}{621}}+O(\d)^1,
\nonumber\\ &\!\!\!\!\!\!\!\!\!\!\!\!\!\!\!\!\!\!\!\!\!\!\!\!\!\!\!&
{\rm(ii)\ }
|A_2'|=\bA_1^{\frac{14}{83}-8}\xX_2^{-1}+O(\d)^1
\le\bA_1^{-\frac{497938}{51543}}+O(\d)^1\simeq0.171.
\end{eqnarray}
}\NOUSE{
First by \eqref{ImMpP}\,(5),\,(6), we have $
|A_3|=|\tildeX_1|^{-2}+O(\d)^3$. Using this and
\eqref{ImMpP}\,(5),\,\eqref{cont-B2}, we deduce
from \eqref{ImMpP}\,(4) the following (noting that $|A_1|\le1$),
\begin{eqnarray}
\!\!\!\!\!\!\!\!\!\!\!\!\!\!\!\!\!\!&&
\label{Bsnwnw}
5|A_1|\le|\tildeX_1|^{134}
+3|A_1^{\frac{660}{533}}
\tildeX_1^{4\ell_0(\ell_0+2)+8\ell_0(4\ell_0+1)-4\ell_0-4\ell_0(4\ell_0+1)+4\ell_0}|
+O(\d)^3\!\!\!
\nonumber\\
\!\!\!\!\!\!\!\!\!\!\!\!\!\!\!\!\!\!&&
\phantom{4|A_1|}
\le(1-\d^3)|A_1|^2+O(\d)^3,
\end{eqnarray}
a contradiction, which proves \eqref{ImMpP}\,(9).
}%
%
%

\NOUSE
{
Next assume
$(1-\d^3)|A_{\rOnE}| \ge|X_{\OnE}|$.
Then by
 \eqref{LetNSoOP}\,(iii), we have
\equa{Msuxxxxx}{\dis
|X_{\OnE}|=(1-\d^3)|A_{\rOnE}|.}
%
From this and \eqref{ImMpP}\,(4),\,\eqref{X202020}, we have,
\begin{eqnarray}
\label{X0sm,m+1-}
&\!\!\!\!\!\!\!\!\!\!\!\!\!\!\!\!\!\!\!\!\!\!\!\!\!&
1+O(\d)^3=|(X_2A_3)^{-1}X_2^4A_1^{-1}|\le(1-\d^3)^3+O(\d)^3=1-3\d^3+O(\d)^3,\end{eqnarray}
a contradiction.
 This proves
\eqref{ImMpP}\,(8).
With $A_2$ defined in \eqref{LetNSoOP----1}\,(v), using \eqref{ImMpP}\,(2),\,(8) and the fact that $|A_2|=|A_1|^{-\frac{587}{20}}+O(\d)^3$,  we obtain%
,
\begin{eqnarray}
\label{A2====}
&\!\!\!\!\!\!\!\!\!\!\!\!\!\!\!\!\!\!\!\!\!\!\!\!&
|A_2|=
\frac{1}{|A_1^{47} A_1^{-\frac{587}{10}}X_2^{18}|}+O(\d)^3
\le(1-\d^3)^{-18}|A_1|^{-\frac{63}{10}}+O(\d)^3\le1+O(\d)^3
.\end{eqnarray}
This proves \eqref{ImMpP}\,(9).
%
}\NOUSE%
{%
\begin{rema}\rm\label{Twoformula-rema}We would like to mention that \eqref{1=bbb} and Lemma \ref{X1Small} are ones of the most difficult parts for us to obtain, which, as stated in Remark \ref{AboutX1}\,(iii)\,(c), causes the definition of $V_2$ to be complicated.
\end{rema}
}\NOUSE{%
\begin{lemm}\label{X1Small}
The equality cannot occur in the last inequality of {\rm\eqref{LetNSoOP}\,(ii)}, i.e.,
\equa{equ-X1Small}{\dis(1-\d^3)|A_1|^{-1}<|\tildeX_1|.}
\end{lemm}
\noindent{\it Proof.~}Assume
\equa{equ-X1Small+1}{\dis(1-\d^3)|A_1|^{-1}=|\tildeX_1|.}
Using the facts from \eqref{ImMpP}\,(5),\,(6) that $|X_2|=|\tildeX_1|^{\ell_0}+O(\d)^3,\,
|A_3|=|\tildeX_1|^{-2\ell_0}+O(\d)^3$,
we obtain from \eqref{ImMpP}\,(4) the following, where for convenience we denote
$\xX:=|\tildeX_1|,\,\bB:=|A_1|\ge1$,
\begin{eqnarray}
\label{X0sm,m-1+1}
&\!\!\!\!\!\!\!\!\!\!\!\!\!\!\!\!\!\!\!\!\!\!\!\!\!\!\!\!\!\!\!&
5\le
\bB^{-1}\xX^{40\ell_0+4}+4\bB^{-1+33\frac{5}{137}}\xX^{-1}+O(\d)^3
\!\!\!\!
\nonumber\\[0pt]&\!\!\!\!\!\!\!\!\!\!\!\!\!\!\!\!\!\!\!\!\!\!\!\!\!\!\!\!\!\!\!&
\phantom{5}
=f(\bB)+O(\d)^3,\ \ f(\bB)=(1-\d^3)^{40\ell_0+4}\bB^{-40\ell_0-5}
+4(1-\d^3)^{-1}\bB^{\frac{165}{137}}.
\end{eqnarray}
Note that $f(\bB)$ is a decreasing function of $\bB$ when $\bB\le\bA_0$, where \equa{b3p}{\dis
\bA_0=\Big(\frac{137(1-\d^3)^{40\ell_0+5}(40\ell_0+5)}{4\times165}\Big)^{\frac1{40\ell_0+5+\frac{165}{137}}}=
\Big(\frac{137(40\ell_0+5)}{660}\Big)^{\frac{\d_0}{40}+O(\d_0)^2}+O(\d)^3.
}
If $\bB\le\bA_0$ then $5\le f(\bB)+O(\d^4)\le f(1)+O(\d^4)<5-\d^3$, a contradiction. Thus $\bB>\bA_0$.
From this, using \eqref{ImMpP}\,(8) and formulas $|A_3|=\xX^{-2\ell_0}+O(\d)^3,\,|X_2|=\xX^{\ell_0}$, we have
\equa{A2===1m2m2}{\dis\!\!\!\!\!\!\!\!\!\!\!\!
|A_2|{\ssc\!}={\ssc\!}\bB^{-2+\frac{165}{137}}\xX^{40\ell_0+3}
{\ssc\!}+{\ssc\!}O(\d^4){\ssc\!}={\ssc\!}\bB^{-40\ell_0(1+O(\d_0)^1)}{\ssc\!}+{\ssc\!}O(\d)^3{\ssc\!}\le
{\ssc\!}\bA_0^{-40\ell_0(1+O(\d_0)^1)}{\ssc\!}+{\ssc\!}O(\d)^3{\ssc\!}={\ssc\!}O(\d_0)^1.
\!\!\!\!\!\!\!\!\!}
Thus by \eqref{1=bbb},\,\eqref{ImMpP}\,(7), we have
\begin{eqnarray}
\label{X1-ell02}
\!\!\!\!\!\!\!\!\!\!\!\!\!\!\!\!\!\!\!\!\!\!\!\!&&
\xX^{4\ell_0^2}=|A_1|\Big(\frac54+O(\d_0)^1\Big)+O(\d)^3=
\bB^{1-\frac{165}{137}}\xX^{4\ell_0^2+1}\Big(\frac54+O(\d_0)^1\Big)+O(\d)^3.
\end{eqnarray}
Using \eqref{equ-X1Small} in \eqref{X1-ell02}, we solve that $\bB=\big(\frac54\big)^{\frac{137}{28}}+O(\d_0)^1$ (recall that $0<\d\ll\d_0$). However similarly to \eqref{X0sm,m-1+1}, we have
\begin{eqnarray}
\label{X0sm,m-1+1+1}
&\!\!\!\!\!\!\!\!\!\!\!\!\!\!\!\!\!\!\!\!\!\!\!\!\!\!\!\!\!\!\!&
5\ge 4(1-\d^3)^{-1}\bB^{\frac{165}{137}}-(1-\d^3)^{40\ell_0+4}\bB^{-40\ell_0-5}+O(\d)^3
\!\!\!\!
\nonumber\\[0pt]&\!\!\!\!\!\!\!\!\!\!\!\!\!\!\!\!\!\!\!\!\!\!\!\!\!\!\!\!\!\!\!&
\phantom{5}
\ge4\big(\frac54\big)^{\frac{165}{28}}-1+O(\d_0)^1>5,
\end{eqnarray}
a contradiction.
First assume $\bA_3\le\bA_3'$. Then \eqref{ImMpP}\,(3) shows that $f(\bA_3)\le f(\bA_1^{c_3c_1^{-1}+\d^4})$.
Using this in the right-hand side of \eqref{X0sm,m-1+1}\,(i), we obtain that it is $\le 2(1+\d^3)^{-1}+O(\d)^3
=2-2\d^3+O(\d)^3$, a contradiction.\hfill$\Box$
}\NOUSE
{%
\begin{lemm}\label{holo}Every branch of any of
$A_1,A_2,A_3$ is 
a locally holomorphic function of $\tilde X_1,X_2,Z$  when $(p_1,p_2)\in V_0$.
\end{lemm}\noindent{\it Proof.~}~We see from definition \eqref{LetNSoOP----1} that it suffices to prove that $A_1$ is locally holomorphic.
Recall Remark \ref{B123-unique}\,(iii). Assume $\frac{\partial A_4}{\partial A_1}=0$. Then by \eqref{LetNSoOP----1}\,(i), we can easily compute,
\equa{pB4}{\dis\!\!\!\!\!\!\!\!\!\! 0{\ssc}=
{\ssc}A_1\frac{\partial A_4}{\partial A_1}{\ssc}={\ssc}\big(1{\ssc}-{\ssc}\d_0-\frac{\d_0^2}{2}\big)
A_1^{1-\d_0-\frac{\d_0^2}{2}}{\ssc}-{\ssc}2
\frac{Z \Big({\ssc}\frac{2-\d_0^2}{4}\tildeX_{\ZeRo}^{2}A_{\rOnE}^{2}\Big)}
{X_{\OnE}^{1+\d_{\rZeRo}^2+\d_{\rZeRo}^3}\tildeX_{\ZeRo}A_3^{a_3}},
\!\!\!\!\!\!\!\!}
 or   $
A_1^{1-\d_0-\frac{\d_0^2}{2}}{\ssc}=\big({\ssc}\frac{2}{1{\ssc}-{\ssc}\d_0-\frac{\d_0^2}{2}}\big) \frac{Z \big({\ssc}\frac14(2
{\ssc}-{\ssc}\d_{\rZeRo}^2)\tildeX_{\ZeRo}^2A_{\rOnE}^{2}\big)}
{X_{\OnE}^{1+\d_{\rZeRo}^2+\d_{\rZeRo}^3}\tildeX_{\ZeRo}A_3^{a_3}}$.
Using this in the first term of $A_4$ in \eqref{LetNSoOP----1}\,(i), we immediately obtain the equality below,
where the inequality follows from \eqref{LetNSoOP}\,(ii),
\begin{eqnarray*}
\label{imm_b}
&\!\!\!\!\!\!\!\!\!\!\!\!\!\!\!\!\!\!\!\!& \frac{2+\d_0^2}{4}=
\Big(\frac{2}{1{\ssc}-{\ssc}\d_0-\frac{\d_0^2}{2}}-1\Big)\frac{2-\d_0^2}{4}|\tildeX_1^{2}A_1^{2}|
\ge
\frac12\Big(1+2\d_0+O(\d_0)^2\Big)(1-\d)=\frac12+\d_0+O(\d_0)^2,
\end{eqnarray*}
a contradiction. This proves that $\frac{\partial A_4}{\partial A_1}$ cannot take the zero value when $(p_1,p_2)\in
V_0$. By the implicit functions theorem, $A_1$ is 
a locally holomorphic function when $(p_1,p_2){\ssc}\in{\ssc} V_0$.\hfill$\Box$
}%

\begin{lemm}\label{Prop(ii)Holds}
We have $\ol V_0=V_0$
.
\end{lemm}
\noindent{\it Proof.~}Let $(p_{\ZeRo},p_{\OnE})\in \ol V_{\rZeRo}$.
Using notations in \eqref{denote-t2}, we already see from
  \eqref{ImMpP}
~that all strict inequalities in
  \eqref{ToSayas+1}\,(c)
[or \eqref{LetNSoOP}\,(iii)$\ssc\,$]
~hold
. In particular, $(p_{\ZeRo},p_{\OnE})\in S_2$ (cf.~proof of Lemma \ref{lemm-condition-XZ}).

By Lemma \ref{lemm-condition-XZ}, first assume  equality occurs in the first inequality of \eqref{C+ToSayas+1}\,(a).
Using notations in \eqref{denote-t2}, we have, where $T_2$ is defined in \eqref{denote-t2},
\NOUSE{ $|A_1|=|A_2|=1$.
we see that 
\eqref{LetNSoOP}\,(iii)--(v) hold. Thus $(p_1,p_2)\in S_2$ by Definition \ref{def-S0-S1}.
 Now we divide the rest 
of the proof 
into two steps.
\vskip5pt\noindent{\it Step 1:
First assume in {\rm\eqref{C+LetNSoOP}\,(i)},
 equality occurs in the first inequality
.}
}\NOUSE{
\begin{rema}\rm\label{Initial-Rema}
We will see below, as mentioned in Remark \ref{AboutX1}\,(iii)\,(a),
that
the powers appearing in
\eqref{LetNSoOP----1}\,(i)--(iii) play key roles
.
\end{rema}
}
%
\begin{eqnarray}
\label{1bab-that}
&\!\!\!\!\!\!\!\!\!\!\!\!\!\!\!\!\!\!\!\!\!\!\!\!\!\!\!&
{\rm(i)\ }\aA_{\rOnE}=\aA_2=1, \ \ \ \ \  {\rm(ii)\ }\xX_2
\stackrel{{}^{\sc\rm\eqref{ImMpP}\,(i)}}{>}1
\NOUSE{
,\ \ \ \ \
{\rm(iii)\ }\xX_2\stackrel{{}^{\sc\rm\eqref{We0o0o0o}\,(i)}}{\le}1+O(\d)^1,
\nonumber\\
&\!\!\!\!\!\!\!\!\!\!\!\!\!\!\!\!\!\!\!\!\!\!\!\!\!\!\!&
{\rm(iv)\ }\tilde\xX_1\stackrel{{}^{\sc\rm\eqref{ImMpP}\,(iv),\,\eqref{1bab-that}\,(ii),\,(iii)}}{=}1+O(\d)^1,\ \
\stackrel{{}^{\sc\rm\eqref{equa-Case6-lemm}\,(iv)}}{ \implies}\ \ \ \ \
{\rm(v)\ }a=1+O(\d)^1\mbox{ for }a\in T_2
}
.
\end{eqnarray}
Then we can use the definitions of $A_1,A_2$ in \eqref{LetNSoOP----1-re-give} to obtain
\begin{eqnarray}
&&\!\!\!\!\!\!\!\!\!\!\!\!\!\!\!\!\!\!\!\!\!\!\!\!\!\!\!
\label{1+1bab-that}
{\rm(i)\ }
2\stackrel{{}^{\sc\rm\eqref{LetNSoOP----1-re-give}\,(i)}}{\le}
\frac1{\tilde\xX_1}+\frac{\aA_1}{\xX_2^{\ell_0}}
\stackrel{{}^{\sc\rm\eqref{1bab-that}}}{<}
\frac1{\tilde\xX_1}+1,\ \ \implies\ \ \
{\rm(ii)\ }\tilde\xX_1<1,
\\\nonumber
&&
\!\!\!\!\!\!\!\!\!\!\!\!\!\!\!\!\!\!\!\!\!\!\!\!\!\!\!
{\rm(iii)\, }
2
\stackrel{{}^{\sc\rm\eqref{LetNSoOP----1-re-give}\,(ii)}}{\le}
\frac{\tilde\xX_1}{\aA_1^2}
+\frac{\aA_2\zZ}{\xX_2^{\ell_0+1}}
\stackrel{{}^{\sc\rm\eqref{1+1bab-that}\,(ii),\,\eqref{1bab-that}}}{<}
1+\frac{\zZ}{\xX_2},
\ \ \implies\ \ \  {\rm(iv)\ }1
\stackrel{{}^{\sc\rm\eqref{1bab-that}\,(ii)}}{<}\xX_2<\zZ.
\!\!\!\!\!\!\!\!\!\!\!\!\!\!\!\!\!
\end{eqnarray}
\NOUSE{%
Regarding $\gamma_1$ as a function on $(\xX_2,\tilde\zZ)$ defined
in the  small neighborhood ${\mathcal O}$ of $(1,1)$, where
${\mathcal O}=\big\{(\xX_2,\tilde\zZ)\,\big|\,1-\d_2<\xX_2,\tilde\zZ<1+\d_2\big\}
$
 [recalling from \eqref{MSmde33333} that $0<\d=\ell^{-1}\ll\d_2=\ell_2^{-1}\ll\d_1=\ell_1^{-1}\ll1$],
we have
$\frac{\partial\gamma_1}{\ptl\xX_2}|_{(\xX_2,\tilde\zZ)=(1,1)}=-\ell_0^{15}+\ell_0-1<0$,
  which imply
\equa{mamam}{\dis
\frac{\partial\gamma_1}{\ptl\xX_2}<0
\ \mbox{ for \ }(\xX_2,\tilde\zZ)\in{\mathcal O}.
}
This implies that $\gamma_1$ is a strictly decreasing function on $\xX_2$
when $(\xX_2,\tilde\zZ)\in{\mathcal O}$, which
with the fact in \eqref{1bab-that}\,(ii) that $\xX_2<1$ implies
\begin{eqnarray}
&&\!\!\!\!\!\!\!\!\!\!\!\!\!\!\!\!\!\!\!\!\!\!\!\!\!\!\!
\label{1+1bab-thatAAA}
1+\d_0^5\stackrel{{}^{\sc\rm\eqref{1+1bab-that}\,(iii)}}{<}
\gamma_1<\gamma_1|_{\xX_2=1}=1+\d_0^5\tilde\zZ
.
\end{eqnarray}
Thus $\tilde\zZ>1$, and so $1<\xX_2<\zZ$ and $\tilde\xX_1<z$ by \eqref{1bab-that}\,(ii),\,\eqref{1+1bab-that}\,(ii),\,(iv).
}%
By \eqref{tX1==}, we also have $\xX_1=\a_0^{-1}\tilde\xX_1<\tilde\xX_1<1<\zZ$.
\NOUSE{%
One can easily observe that $\tilde\gamma_1$ is a strictly increasing function on $\tilde\xX_1$ when $|\tilde\xX_1-1|<\d_2$
(recalling that $0<\d\ll\d_2\ll\d_0=\ell_0^{-1}\ll1$). From this and \eqref{1bab-that}\,(v),\,\eqref{1+1bab-thatAAA}\,(i),
we obtain that $\tilde\xX_1>1$.
This with \eqref{1+1bab-thatAAA}\,(ii) and the fact that $\tilde\gamma_2$ is a strictly decreasing on $\tilde\xX_1$ when
$(\tilde\xX_1,\xX_2,\tilde\zZ)\in{\mathcal O}$, we obtain that $2<\tilde\gamma_2|_{\tilde\xX_1=1}=2\tilde\zZ$., i.e., $\tilde\zZ>1$.
Therefore,
\equa{amasmsdmsdm}{\dis
{\rm(i)\ }1\stackrel{{}^{\sc\rm\eqref{1bab-that}\,(ii)}}{<}\xX_2<\xX_2\tilde\zZ\stackrel{{}^{\sc\rm\eqref{1+1bab-that}\,(iii)}}{=}\frac{\zZ}{\tilde\xX_1}<\zZ,\ \ \
{\rm(ii)\ }\xX_1\stackrel{{}^{\sc\rm\eqref{tX1==}}}{<}\tilde\xX_1<\tilde\xX_1\tilde\zZ\stackrel{{}^{\sc\rm\eqref{1+1bab-that}\,(iii)}}{=}\frac{\zZ}{\xX_2}<\zZ.}
}%
We obtain
,
where the first equality follows by noting from \eqref{TaKa},\,\eqref{SimMMSMS},\,\eqref{denote-t2} that
 $|x_1|=\xX_1\kk,\,|x_2|=\xX_2\kk$%
,
\equa{Conndndm}{\dis
\g_{|x_1|,|x_2|}
=\g_{\xX_1\kk,\xX_2\kk}\stackrel{{}^{\sc\rm\eqref{wePPPP1+}}}{<}
\g_{\zZ\kk,\xX_2\kk}
\stackrel{{}^{\sc\rm\eqref{wePPPP1}}}{\le}\g_{\zZ\kk,\zZ\kk}
\stackrel{{}^{\sc\rm\eqref{AssG}}}{\le}\zZ\g_{\kk,\kk}
\stackrel{{}^{\sc\rm\eqref{TaKa},\,\eqref{SimMMSMS}}}{=}|x_2+y_2|\stackrel{{}^{\sc\rm\eqref{Ak=1}}}{\le}\g_{|x_1|,|x_2|},
}which is
a contradiction.
\NOUSE{
First assume $\zZ<1$.
We want to verify whether or not \eqref{TheFasss} holds. To do this, we use notation there
by denoting $\hat k=\zZ$ and so $|x_2+y_2|=
\hat k\gamma_{\kk,\kk}$ by
\eqref{TaKa},\,\eqref{SimMMSMS}, and $k_1:=\xX_1$ with $|x_1|=k_1\kk$, and $k_2:=\xX_2$
with $|x_2|=k_2\kk$. Then we have
\begin{eqnarray}
&&\!\!\!\!\!\!\!\!\!\!\!\!\!\!\!\!\!\!\!\!\!\!\!\!\!\!\!\!\!\!
\label{mamama0009ws9w9}
{\rm(i)\ }\l\stackrel{{}^{\sc\rm\eqref{ImMpP}\,(iv)}}{<}
k_2=\xX_2\stackrel{{}^{\sc\rm\eqref{1+1bab-that}\,(ii)}}{<}1\mbox{ with }\l=1-\d, \ \ \ \ \ {\rm(ii)\ }\hat k=\zZ<1,
\nonumber\\
&&\!\!\!\!\!\!\!\!\!\!\!\!\!\!\!\!\!\!\!\!\!\!\!\!\!\!\!\!\!\!
{\rm(iii)\ }k_2=\xX_2\stackrel{{}^{\sc\rm\eqref{1+1bab-that}\,(ii)}}{\le}
\hat k^{1+\l_2}\mbox{ with }\l_2=1,\ \ \ {\rm(iv)\ }k_1=\xX_1\stackrel{{}^{\sc\rm\eqref{1+1bab-that}\,(i)}}{<}1
.
\end{eqnarray}
This shows that \eqref{TheFasss} holds,
which contradicts the assumption that Case 5 does not occur (at this point it may be worth mentioning that when we say Case 5 does not occur, it means that if it occurs we have already proved that Proposition \ref{real00-inj} holds).
Next assume $\zZ\ge1$ and $\xX_2\le\zZ$.
Then $\xX_1<1\le\zZ$, and
we obtain
,
where the first equality follows by noting from \eqref{TaKa},\,\eqref{SimMMSMS},\,\eqref{denote-t2} that
 $|x_1|=\xX_1\kk,\,|x_2|=\xX_2\kk$
[where in the third inequality,
 in order to apply \eqref{AssG}, we require $\zZ>1$; however in our case here when $\zZ=1$ the third inequality is an equality],
\equa{Conndndm}{\dis
\g_{|x_1|,|x_2|}
=\g_{\xX_1\kk,\xX_2\kk}\stackrel{{}^{\sc\rm\eqref{wePPPP1+}}}{<}
\g_{\zZ\kk,\xX_2\kk}
\stackrel{{}^{\sc\rm\eqref{wePPPP1}}}{\le}\g_{\zZ\kk,\zZ\kk}
\stackrel{{}^{\sc\rm\eqref{AssG}}}{\le}\zZ\g_{\kk,\kk}
\stackrel{{}^{\sc\rm\eqref{TaKa},\,\eqref{SimMMSMS}}}{=}|x_2+y_2|\stackrel{{}^{\sc\rm\eqref{Ak=1}}}{\le}\g_{|x_1|,|x_2|},
}which is
a contradiction.
Now assume $\zZ\ge1$ and $\xX_2>\zZ$. Then by \eqref{1+1bab-that}, we have
\equa{++1+1bab-that}{\dis
{\rm(i)\ }1-2\d\stackrel{{}^{\sc\rm\eqref{ImMpP}\,(iii)}}{<}\xX_1<1,\ \ \ \ \ {\rm(ii)\ }1<\xX_2<\zZ^{1+\eE_1},\ \ \ \ \
{\rm(iii)\ }1=\aA_1^{\eE_1}\aA_2^{2\eE_1}=
\xX_2^{2-\eE_1}\zZ^{-2}\stackrel{{}^{\sc\rm\eqref{ImMpP}\,(i)}}{=}\Big(1+O(\d)^2\Big)\xX_2^{-\eE_1}.
}
\begin{clai}\label{last-claim}
Equ.~\eqref{++1+1bab-that} implies Proposition \ref{real00-inj}.
\end{clai}
To prove the claim, assume \eqref{++1+1bab-that} holds for some element denoted as
$(\hat p_1,\hat p_2)=\big((\hat x_1,\hat y_1),(\hat x_2,\hat y_2)\big)\in V$.
Then \eqref{++1+1bab-that} gives that
$k_1:=\kk^{-1}|\hat x_1|<1,\,1<k_2:=\kk^{-1}|\hat x_2|<\hat k^{1+\eE_1}$ with $\hat k=\gamma_{\kk,\kk}^{-1}|x_2+y_2|$.
Observe from \eqref{LetNSoOP----1}\,(ii),\,(iii) that we have the following relation between $A_1$ and $A_2$,
\begin{eqnarray}
\!\!\!\!\!\!\!\!\!\!\!\!\!\!\!\!\!\!\!\!\!\!\!\!\!\!\!&&
\label{NewA1==}
{\rm(i)\ }A_1=
\frac{\widetilde X_1^2 X_2}{Z}\Big (\frac{2 + 2 \d_0}{1 + 2 \d_0} -\frac{A_2 X_2^2}{
   (1 + 2 \d_0)A_3\widetilde X_1^2 Z^2}\Big),\ \ \ \implies
\nonumber\\
\!\!\!\!\!\!\!\!\!\!\!\!\!\!\!\!\!\!\!\!\!\!\!\!\!\!\!&&
{\rm(ii)\ }
\frac{2+2\d_0}{1+2\d_0}\le
\frac{\aA_2\xX_2^2}{(1+2\d_0)\aA_3\tilde\xX_1^{2}\zZ^2}+\frac{\aA_1\zZ}{\tilde\xX_1^2\xX_2}.
\stackrel{{}^{\sc\rm\eqref{1bab-that}}}{<}\gamma_1:=
\frac{1}{(1+2\d_0)\tilde\xX_1^{2}\tilde\zZ^2}+\frac{\tilde\zZ}{\tilde\xX_1^2},\ \ \ \tilde\zZ=\frac{\zZ}{\xX_2}.
\end{eqnarray}
We claim, where $T_2$ is defined in \eqref{denote-t2},
\equa{1+O(d)-element}{\dis
a=1+O(\d)^1\mbox{ \ for \ }a\in T_2.
}
We see from \eqref{Formmeme},\,\eqref{x1-x2==}\,(v),\,\eqref{1bab-that} that \eqref{1+O(d)-element} holds for $a=\aA_1,\aA_2,\tilde\xX_1$, and
$\aA_3\ge1+O(\d)^1$. Then \eqref{LetNSoOP----1-redefine++}\,(iv) shows that $\aA_3\le1+O(\d)^1$. Thus
\eqref{1+O(d)-element} holds for $\aA_3$, and hence also for $a=\xX_2,\zZ,\xX_1$ by
\eqref{LetNSoOP----1-redefine++}\,(i),\,\eqref{equa-Case6-lemm}\,(iv),\,\eqref{tX1==}.
Finally by \eqref{LetNSoOP----1}\,(iv),\,(v), we see
\eqref{1+O(d)-element} holds for all $a\in T_2$.
%
%
%
Now we have
,
\begin{eqnarray}
\label{LetNSoOP----1-redefine-----}
&\!\!\!\!\!\!\!\!\!\!\!\!\!\!\!\!\!\!\!\!\!\!\!\!\!\!\!\!\!\!\!\!\!\!
&
{\rm(i)\ }\tilde\xX_1\stackrel{{}^{\sc\rm\eqref{LetNSoOP----1}\,(i)}}{=}
\Big(\aa_3^{-1}
\zZ^{ 280\ell_0+376} \xX_2^{-280\ell_0-375}\Big)^{\frac1{80\ell_0+50}}
<
\Big(
\zZ^{ 280\ell_0+376} \xX_2^{-280\ell_0-375}\Big)^{\frac1{80\ell_0+50}}
,\nonumber\\
&\!\!\!\!\!\!\!\!\!\!\!\!\!\!\!\!\!\!\!\!\!\!\!\!\!\!\!\!\!\!\!\!\!\!
&
{\rm(ii)\ }
3\stackrel{{}^{\sc\rm\eqref{More-de-1}}}{\le}
\aA_3^{-1}\tilde\xX_1^{-10}+
2\aA_1\tilde  \xX_1^{20\ell_0+5}\xX_2^{70\ell_0+5}\aA_3^{-2}\zZ^{-70\ell_0-4}
<
\gamma_1:=
\tilde\xX_1^{-10}+
2\tilde  \xX_1^{20\ell_0+5}\xX_2^{70\ell_0+5}\zZ^{-70\ell_0-4}
,\!\!\!\!\!\!\!\!\!\!\!\!\!\!\!\!\!\!\!\!\!\!\!\!\!\!\!\!\!\!
\nonumber\\
&\!\!\!\!\!\!\!\!\!\!\!\!\!\!\!\!\!\!\!\!\!\!\!\!\!\!\!\!\!\!\!\!\!\!
&
{\rm(iii)\ }
7\stackrel{{}^{\sc\rm\eqref{LetNSoOP----1}\,(iii)}}{\le}
2\aA_3^{-1}\tilde\xX_1^{-10}+5\aA_2\zZ^{90}\tilde \xX_1^{-6}\xX_2^{-80}
<\gamma_2:=
2\tilde\xX_1^{-10}+5\zZ^{90}\tilde \xX_1^{-6}\xX_2^{-80},
\nonumber\\
&\!\!\!\!\!\!\!\!\!\!\!\!\!\!\!\!\!\!\!\!\!\!\!\!\!\!\!\!\!\!\!\!\!\!
&
{\rm(iv)\ }
5=5\aA_2\stackrel{{}^{\sc\rm\eqref{Related-A1-A2}}}{\le}
\tilde \xX_1^6 \xX_2^{80}\zZ^{-90}(1+
4\aA_1\aA_3^{-2}\tilde  \xX_1^{20\ell_0+5}\xX_2^{70\ell_0+5}\zZ^{-70\ell_0-4}\Big)
\nonumber\\
&\!\!\!\!\!\!\!\!\!\!\!\!\!\!\!\!\!\!\!\!\!\!\!\!\!\!\!\!\!\!\!\!\!\!
&
\phantom{{\rm(iv)\ }
1}
<
\Big(
\zZ^{ 280\ell_0+376} \xX_2^{-280\ell_0-375}\Big)^{\frac{6}{80\ell_0+50}}
\xX_2^{80}\zZ^{-90}
%
\nonumber\\
&\!\!\!\!\!\!\!\!\!\!\!\!\!\!\!\!\!\!\!\!\!\!\!\!\!\!\!\!\!\!\!\!\!\!
&
\phantom{{\rm(iv)\ }
1<\gamma_{20}:=}\mbox{\Large$\times$}
%
\mbox{\Large$\Big($}1+
4\Big(
\zZ^{ 280\ell_0+376} \xX_2^{-280\ell_0-375}\Big)^{\frac{20\ell_0+5}{80\ell_0+50}}
\xX_2^{70\ell_0+5}\zZ^{-70\ell_0-4}\mbox{\Large$\Big)$}
\nonumber\\
&\!\!\!\!\!\!\!\!\!\!\!\!\!\!\!\!\!\!\!\!\!\!\!\!\!\!\!\!\!\!\!\!\!\!
&
\phantom{{\rm(iv)\ }1}
=\xX_2^{59+O(\d_0)^1}\zZ^{-69+O(\d_0)^1}\Big(1+4
\xX_2^{-\frac{125}{2}+O(\d_0)^1}\zZ^{\frac{255}{4}+O(\d_0)^1}
\Big)
\nonumber\\
&\!\!\!\!\!\!\!\!\!\!\!\!\!\!\!\!\!\!\!\!\!\!\!\!\!\!\!\!\!\!\!\!\!\!
&
\phantom{{\rm(iv)\ }1}
=\gamma_{20}:=
\xX_2^{59+O(\d_0)^1}\zZ^{-69+O(\d_0)^1}+4
\xX_2^{-\frac{7}{2}+O(\d_0)^1}\zZ^{-\frac{21}{4}+O(\d_0)^1}
.
\end{eqnarray}
%
Regarding $\gamma_3,\gamma_1,\gamma_2,\gamma_{20}$ as  functions on $\tilde\xX_1,\xX_2,\zZ$, one can easily compute
that at point $(\tilde\xX_1,\xX_2,\zZ)=(1,1,1)$, we have
\begin{eqnarray}
\label{FiGamma1}
&&\!\!\!\!\!\!\!\!\!\!\!\!\!\!\!\!\!\!\!\!\!\!\!\!
{\rm(i)\ \ }\frac{\ptl\gamma_1}{\ptl\tilde\xX_1}
=40\ell_0,
\ \ \ \ \ \
\frac{\ptl\gamma_1}{\ptl\xX_2}
=2(70\ell_0+5),
\  \ \ \ \ \
\frac{\ptl\gamma_1}{\ptl\zZ}
=-2(70\ell_0+4),
\nonumber\\
&&\!\!\!\!\!\!\!\!\!\!\!\!\!\!\!\!\!\!\!\!\!\!\!\!
{\rm(ii)\ \ }\frac{\ptl\gamma_2}{\ptl\tilde\xX_1}
=-50,
\ \ \ \ \ \
\frac{\ptl\gamma_2}{\ptl\xX_2}
=-400,
\  \ \ \ \ \
\frac{\ptl\gamma_2}{\ptl\zZ}
=450.
\end{eqnarray}
We claim,
\equa{Claimse-111}{\dis
\zZ^{70+4\d_0}\xX_2^{-(70+5\d_0)}<
\tilde\xX_1^{20+s\d_1}, \mbox{ \ where $s=1$ if $\tilde\xX_1\ge1$ or $s=-1$ if $\tilde\xX_1<1$}.
}
Assume this does not hold, i.e., $\xX_2^{70+5\d_0}\zZ^{-(70+4\d_0)}\le
\tilde\xX_1^{-(20+s\d_1)}$. Then
\equa{22==we=}{\dis
3\stackrel{{}^{\sc\rm\eqref{LetNSoOP----1-redefine-----}\,(ii)}}{<}
\gamma_{10}:=\tilde\xX_1^{-10}+2\tilde\xX_1^{5-\ell_0s\d_1}.
}
Regarding $\gamma_{1}$ as a function on $\tilde\xX_1$
in the small neighborhood $\mathcal O_1=\{\tilde\xX_1\,|\,1-\d_2<\tilde\xX_1<1+\d_2\}$ of $1$
[recalling from \eqref{MSmde33333} that $0<\d=\ell^{-1}\ll\d_2=\ell_2^{-1}\ll\d_1=\ell_1^{-1}\ll1$]%
,
one can easily see from \eqref{LetNSoOP----1-redefine-----}\,(ii),
\equa{SeeeotoSee}{\dis
\frac{d\gamma_{10}}{d\tilde\xX_1}\mbox{\Large$\Big|$}_{\tilde\xX_1=1}=
-10+2(5-\ell_0s\d_1)=-2\ell_0s\d_1>0\mbox{  if $s=-1$, \ \ or \ \ $<0$ if $s=1$.}
}
Thus in ${\mathcal O}_1$, we have%
,
\equa{SeeeotoSee+}{\dis
\frac{d\gamma_{10}}{d\tilde\xX_1}>0\mbox{ if $s=-1$, \ \ or \ \ $<0$ if $s=1$.}
}
Note from \eqref{1+O(d)-element} that we have $\tilde\xX_1\in{\mathcal O}_1$. Thus
\eqref{22==we=} with \eqref{SeeeotoSee+} shows that
$\tilde\xX_1>1$ if $s=-1$, or $\tilde\xX_1<1$ if $s=1$, which is a contradiction with our definition of $s$ in
\eqref{Claimse-111}. This proves \eqref{Claimse-111}, which implies in general (i) below
[noting that we have this formula simply because we have \eqref{FiGamma1}\,(i)$\ssc\,$],
where (ii) is obtained from
\eqref{LetNSoOP----1-redefine-----}\,(ii) by using exactly the same arguments [which can be observed from
\eqref{FiGamma1}\,(ii)$\ssc\,$],
\begin{eqnarray}
\label{Nmsmgama}
&&\!\!\!\!\!\!\!\!\!\!\!\!\!\!\!\!\!\!\!\!\!\!\!\!\!\!\!\!\!\!
{\rm(i)\ }1\le\gamma_{11}:=\tilde\xX_1\Big(\xX_2^{70+5\d_0}\zZ^{-(70+4\d_0)}\Big)^{\frac1{20}+O(\d_1)^1},\ \ \
{\rm(ii)\ }
1<\gamma_{21}:=\tilde\xX_1^{-1}(\xX_2^{-8}\zZ^{9})^{1+O(\d_1)^1},\ \ \ \implies
\!\!\!\!\!\!\!\!\!\!\!\!\nonumber\\
&&\!\!\!\!\!\!\!\!\!\!\!\!\!\!\!\!\!\!\!\!\!\!\!\!\!\!\!\!\!\!
{\rm(iii)\ }1<\gamma_{11}\gamma_{21}
=\xX_2^{-\frac92+\frac{\d_0}{4}+O(\d_1)^1}\zZ^{\frac{11}{2}-\frac{\d_0}{5}+O(\d_1)^1},\ \ \ \implies
\ \ \
{\rm(iv)\ }\xX_2<\zZ^{\frac{11}{9}+O(\d_0)^1}.
\end{eqnarray}
Observe from \eqref{LetNSoOP----1-redefine-----}\,(iv) that $\frac{\ptl\gamma_{20}}{\ptl\xX_2}|_{(\xX_2,\zZ)=(1,1)}=
45+O(\d_0)^1$, $\frac{\ptl\gamma_{20}}{\ptl\zZ}|_{(\xX_2,\zZ)=(1,1)}=
-90+O(\d_0)^1$. From this, as in \eqref{Nmsmgama}\,(i), we can obtain that $\zZ<\xX_2^{\frac12+O(\d_0)^1}$.
This with \eqref{Nmsmgama}\,(iv) shows that $\xX_2,\zZ<1$.
We want to verify whether or not \eqref{TheFasss} holds. To do this, we use notation there
by denoting $\hat k=\zZ$ and so $|x_2+y_2|=
\hat k\gamma_{\kk,\kk}$ by
\eqref{TaKa},\,\eqref{SimMMSMS}, and $k_1:=\xX_1$ with $|x_1|=k_1\kk$, and $k_2:=\xX_2$
with $|x_2|=k_2\kk$. Then we have
\begin{eqnarray}
&&\!\!\!\!\!\!\!\!\!\!\!\!\!\!\!\!\!\!\!\!\!\!\!\!\!\!\!\!\!\!
\label{mamama0009ws9w9}
{\rm(i)\ }\l\stackrel{{}^{\sc\rm\eqref{meme}\,(ii)}}{<}
k_2=\xX_2<1\mbox{ with }\l=\d_2^4, \ \ \ \ \ {\rm(ii)\ }\hat k=\zZ<1,
\nonumber\\
&&\!\!\!\!\!\!\!\!\!\!\!\!\!\!\!\!\!\!\!\!\!\!\!\!\!\!\!\!\!\!
{\rm(iii)\ }k_2=\xX_2\stackrel{{}^{\sc\rm\eqref{Nmsmgama}\,(iv)}}{\le}
\hat k^{1+\l_2}\mbox{ with }\l_2=\frac19,
\nonumber\\
&&\!\!\!\!\!\!\!\!\!\!\!\!\!\!\!\!\!\!\!\!\!\!\!\!\!\!\!\!\!\!
{\rm(iv)\ }
k_1=\xX_1\stackrel{{}^{\sc\rm\eqref{tX1==}}}{<}
\tilde\xX^{\frac1{10}}\xX_2
\stackrel{{}^{\sc\rm\eqref{LetNSoOP----1-redefine-----}\,(i)}}{<}
\Big(\xX_2^{-\frac72+O(\d_0)^1}\zZ^{\frac72+O(\d_0)^1}\Big)^{\frac1{10}}\xX_2
\nonumber\\
&&\!\!\!\!\!\!\!\!\!\!\!\!\!\!\!\!\!\!\!\!\!\!\!\!\!\!\!\!\!\!
\phantom{{\rm(iv)\ }
k_1}=\xX_2^{\frac{13}{20}+O(\d_0)^1}\zZ^{\frac{7}{20}+O(\d_0)^1}
\stackrel{{}^{\sc\rm\eqref{Nmsmgama}\,(iv)}}{<}\zZ^{\frac{103}{90}+O(\d_0)^1}<1.
\end{eqnarray}
This shows that \eqref{TheFasss} holds,
which contradicts the assumption that Case 5 does not occur (at this point it may be worth mentioning that when we say Case 5 does not occur, it means that if it occurs we have already proved that Proposition \ref{real00-inj} holds).
\NOUSE{
This proves that \eqref{T2-isO(ee)} holds for $\tilde\xX_1$.
Now regarding $\gamma_{20}$ in \eqref{LetNSoOP----1-redefine-----}\,(iv) as  a function on
$(\tilde\xX_1,\xX_2,\zZ)$
in the small neighborhood $\mathcal O_2$ of $(1,1,1)$, where
\equa{MSMSMO-2}{\dis
{\mathcal O}_2=\Big\{(\tilde\xX_1,\xX_2,\zZ)\,\Big|\,|\tilde\xX_1-1|<\d_2,\,|\xX_2-1|<\d_2,\,|\zZ-1|<\d_2\Big\},}
one can easily see from \eqref{LetNSoOP----1-redefine-----}\,(iv),
\equa{SeeeotoSee}{\dis
\frac{\ptl\gamma_{20}}{\ptl\zZ}\mbox{\Large$\Big|$}_{(\tilde\xX_1,\xX_2,\zZ)=(1,1,1)}=
48+23\d_1-12=36+23\d_1>0.
}
Thus  we have%
,
\equa{g-210-SeeeotoSee+}{\dis
\frac{\ptl\gamma_{20}}{\ptl\zZ}>0\mbox{ in }{\mathcal O}_2,
}
which means that $\gamma_{21}$ is a strictly
increasing function on $\zZ$
 when $(\tilde\xX_1,\xX_2,\zZ)\in{\mathcal O}_2$.
 By \eqref{ImMpP}\,(iii), we indeed have $(\tilde\xX_1,\xX_2,\zZ)\in{\mathcal O}_2$.
Thus \eqref{LetNSoOP----1-redefine-----}\,(iv) with \eqref{Nmsmgama}\,(iv) gives
\equa{mdmem-0000}{\dis
1\le\gamma_{20}<\gamma_{20}\mbox{\Large$\Big|$}_{\zZ=\xX_2^{\frac{14}{13}+O(\d_1)^3}}
=\gamma_4:=
\tilde\xX_1^{3-13\d_1}
\xX_2^{-\frac{69}{13}-\frac{159\d_1}{13}+O(\d_1)^3}\Big(\frac25+\frac35\tilde\xX_1^{-5}\xX_2^{\frac{45}{13}+O(\d_1)^3}\Big).
}
Now $\gamma_4$ is a strictly decreasing function on $\xX_2$ when $(\tilde\xX_1,\xX_2,1)\in{\mathcal O}_2$
as \equa{ptl-gamma4}{\dis
\frac{\ptl\gamma_4}{\ptl\xX_2}\mbox{\Large$\Big|$}_{(\tilde\xX_1,\xX_2)=(1,1)}
=-\frac{69}{13}-\frac{159\d_1}{13}+O(\d_1)^3+\frac{27}{13}=-\frac{42}{13}-\frac{159\d_1}{13}+O(\d_1)^3<0.
}
From this and \eqref{mdmem-0000}, we can obtain
\equa{X2mdemdm}{\dis\!\!\!\!\!
{\rm(i)\ }\xX_2^{\frac{42}{13}+\frac{159\d_1}{13}+O(\d_1)^3}
{\ssc\!}<{\ssc\!}
\tilde\xX_1^{-13\d_1+\d_1^2}\mbox{ if }\tilde\xX_1{\ssc\!}\ge{\ssc\!}1,
\mbox{ or }
{\rm(ii)\ }\xX_2^{\frac{42}{13}+\frac{159\d_1}{13}+O(\d_1)^3}
{\ssc\!}<{\ssc\!}\tilde\xX_1^{-13\d_1-\d_1^2}\mbox{ if }\tilde\xX_1{\ssc\!}<{\ssc\!}1.\!\!\!\!\!
}
To see this, we consider case (i) by assuming conversely that
$\tilde\xX_1\ge1$ but
$\xX_2\ge\tilde\xX_1^{\frac{-13\d_1+\d_1^2}{{\frac{42}{13}+\frac{159\d_1}{13}+O(\d_1)^3}}}$.
Then by replacing $\xX_2$ in the right-hand side of
\eqref{mdmem-0000} by
$\tilde\xX_1^{\frac{-13\d_1+\d_1^2}{{\frac{42}{13}+\frac{159\d_1}{13}}}+O(\d_1)^3}$, we obtain
a function, denoted as $\gamma_5$, on $\tilde\xX_1$. One can easily observe that
$\frac{d\gamma_5}{\d\tilde\xX_1}|_{\tilde\xX_1=1}=-\d_1^2<0$,
 which means that $\gamma_5$ is a strictly decreasing function
on $\tilde\xX_1$ when $|\tilde\xX_1-1|\le\d_2$ (recalling that $0<\d_2\ll\d_1$), but $\gamma_5$ takes value $>1$ at $\tilde\xX_1=1$ by \eqref{mdmem-0000}.
This implies that
 $\tilde\xX_1<1$, a contradiction
with the assumption. This proves \eqref{X2mdemdm}\,(i). Similarly, we have \eqref{X2mdemdm}\,(ii).
By \eqref{X2mdemdm}, we can in particular obtain (i) below,
\begin{eqnarray}
&&\!\!\!\!\!\!\!\!\!\!\!\!\!\!\!\!\!\!\!\!\!\!\!\!\!\!\!
\label{msmememew}
{\rm(i)\ }
1<\gamma_6:=\tilde\xX_1^{-13\d_1}\xX_2^{-\frac{42}{13}+O(\d_1)^1},\ \ \
{\rm(ii)\ }1<\gamma_6\gamma_{11}^{13\d_1}\stackrel{{}^{\sc\rm\eqref{Nmsmgama}\,(i)}}{=}\xX_2^{-\frac{42}{13}+O(\d_1)^1}\zZ^{-52\d_1+O(\d_1)^2},\ \ \ \implies
\!\!\!\!\!\!\!\!\!\!\!\!\!\nonumber\\
&&\!\!\!\!\!\!\!\!\!\!\!\!\!\!\!\!\!\!\!\!\!\!\!\!\!\!\!
{\rm(iii)\ }\xX_2<\zZ^{-\frac{338\d_1}{21}+O(\d_1)^2.}
\end{eqnarray}
%
%
Then by \eqref{LetNSoOP----1-redefine-----}\,(ii) and \eqref{MDNDNN9999}
[noting from that we
indeed have $(\tilde\xX_1,\xX_2,\zZ)\in{\mathcal O}_1$ as $0<\d\ll\d_2$], we have
\equa{smemem334}{\dis\!\!\!\!\!
2\le\gamma_1<\gamma_1|_{\tilde\xX_1=(\xX_2^{5}\zZ^{-4})^{-1}}=
\gamma_{10}:=\tilde\xX_1^{10}+\tilde\
}
One can easily compute,
\begin{eqnarray}
\label{dm,ememe}
&\!\!\!\!\!\!\!\!\!\!\!\!\!\!\!\!\!\!\!\!\!&
\frac{d\gamma_{10}}{d\zZ}\mbox{\Large$\Big|$}_{\zZ=1}=\frac1{10}\Big(-\frac{223}{3}-20s\d_1\Big)
+\frac9{10}\Big(\frac{47}{3}-20s\d_1\Big)+\frac4{45}(-75+90s\d_1)
\nonumber\\
&\!\!\!\!\!\!\!\!\!\!\!\!\!\!\!\!\!\!\!\!\!&
\phantom{\frac{d\gamma_{10}}{d\zZ}\mbox{\Large$\Big|$}_{\zZ=1}}
=-12s\d_1,
\mbox{ which is $<0$ if $\zZ\ge1$, or
$>0$ if $\zZ<1.$}
\end{eqnarray}
Note that $0<\d_2\ll\d_1$, we see from
\eqref{dm,ememe} that when $|\zZ-1|<\d_2$, $\gamma_{10}$ is a strictly decreasing function on $\zZ$ if $\zZ\ge1$, or a
strictly increasing function on $\zZ$ if $\zZ<1$.
Thus in any case we obtain from \eqref{smemem334} 
that
$1<\gamma_{10}\le\gamma_{10}|_{\zZ=1}=1$, a contradiction.
This proves \eqref{Claimse-111}, which in particular implies that $\xX_2<\zZ^{\frac16+O(\d_1)^1}$.
Similarly, we can use \eqref{1dew2}\,(i) to prove that $\zZ<\xX_2^{10+O(\d_1)^1}<\zZ^{\frac53+O(\d_1)^1}$.
Thus we obtain that $1<\xX_2<\zZ$. We also have
$\xX_1=\a_0^{-1}\tilde\xX_1<\aA_3^{-1}<1<\zZ$ by \eqref{tX1==},\,\eqref{LetNSoOP----1}\,(i),\,\eqref{1bab-that}\,(ii).
By \eqref{LetNSoOP----1}\,(iv),\,\eqref{Formmeme}, we obtain that $\xX_1
=\aA_2^{\frac{\scep}{2}}|B_2|^{-\frac{\scep}{2}}\ge1$, and
by \eqref{LetNSoOP----1}\,(i),\,(iii), $\xX_2=|A_1|=1$, $\zZ=
\aA_3\xX_1^3\xX_2^2>\xX_1^3\ge\xX_1\ge1$.
\NOUSE{
Then we have, where $\dD_1=|D_1|,\,\dD_2=|D_2|$,
\begin{eqnarray}
\label{aaasss--LetNSoOP----1-}
&\!\!\!\!\!\!\!\!\!\!\!\!\!\!\!\!\!\!\!\!\!\!\!\!\!\!\!\!\!\!\!\!\!\!\!\!\!\!\!\!\!\!\!\!\!\!\!\!\!\!\!
&
{\rm(i)\ }\frac43
\stackrel{{}^{\sc\rm\eqref{aaasss--LetNSoOP----1}\,(i)}}{\le}
\aA_3^{-1}\dD_2^{-1}+\frac{\dD_2^4}{3}\stackrel{{}^{\sc\rm\eqref{1bab-that}\,(ii)}}{<}
\gamma_{30}(\dD_2):=\dD_2^{-1}+\frac{\dD_2^4}{3},
\!\!\!\!\!\!\!\!\!
\!\!\!\!\!\!\!\!\!
\!\!\!\!\!\!\!\!\!
\!\!\!\!\!\!\!\!\!
\nonumber\\
&\!\!\!\!\!\!\!\!\!\!\!\!\!\!\!\!\!\!\!\!\!\!\!\!\!\!\!\!\!\!\!\!\!\!\!\!\!\!\!\!\!\!\!\!\!\!\!\!\!\!\!
&
{\rm(ii)\ }1=\aA_2^{-1}\stackrel{{}^{\sc\rm\eqref{LetNSoOP----1}\,(ii)}}{=}
\gamma_2:=\aA_3^{-1}\tilde\xX_1^{-2} \zZ^2\xX_2^4,
\!\!\!\!\!\!\!\!\!
\!\!\!\!\!\!\!\!\!
\!\!\!\!\!\!\!\!\!
\!\!\!\!\!\!\!\!\!
\nonumber\\
&\!\!\!\!\!\!\!\!\!\!\!\!\!\!\!\!\!\!\!\!\!\!\!\!\!\!\!\!\!\!\!\!\!\!\!\!\!\!\!\!\!\!\!\!\!\!\!\!\!\!\!
&
{\rm(iv)\ }3\stackrel{{}^{\sc\rm\eqref{aaasss--LetNSoOP----1}\,(iv)}}{\le}
\aA_1\dD_{1}^{-1}+2\aA_3^{-1}\dD_2^{-1}\stackrel{{}^{\sc\rm\eqref{1bab-that}}}{<}
\dD_{1}^{-1}+2\dD_2^{-1}.
\!\!\!\!\!\!\!\!\!
\end{eqnarray}
Observe that $\gamma_{30}(\dD_2)$ is a strictly increasing function on $\dD_2$ when $\dD_2>\big(\frac34\big)^{\frac15}
$, and $\gamma_{30}(1)=\frac43$. Since $\dD_2\ge\aA_1^{\frac{11}{40}}+O(\d)^1=1+O(\d)^1$ by
\eqref{ImMpP+1}\,(ii),\,\eqref{1bab-that}\,(i), we obtain from \eqref{aaasss--LetNSoOP----1-}\,(i),\,(iv) the following,
\begin{eqnarray}
\label{D333aaa}
&&\!\!\!\!\!\!\!\!\!\!\!\!\!\!\!\!\!\!\!\!\!\!\!\!\!\!
{\rm(i)\ }1\stackrel{{}^{\sc\rm \eqref{aaasss--LetNSoOP----1-}\,(i)}}{<}
\dD_2\stackrel{{}^{\sc\rm\eqref{DeD12}\,(ii)}}{=}\gamma_3:=
\xX_2^{21} \zZ^2\tilde \xX_1^{-22},
\nonumber\\
&&\!\!\!\!\!\!\!\!\!\!\!\!\!\!\!\!\!\!\!\!\!\!\!\!\!\!
{\rm(ii)\ }1\stackrel{{}^{\sc\rm \eqref{aaasss--LetNSoOP----1-}\,(iv),\,\eqref{D333aaa}\,(i)}}{<}
\dD_1^{-1}\stackrel{{}^{\sc\rm\eqref{DeD12}\,(i)}}{=}\gamma_1:=\aA_3\tilde\xX_1\xX_2,
\nonumber\\
&&\!\!\!\!\!\!\!\!\!\!\!\!\!\!\!\!\!\!\!\!\!\!\!\!\!\!
{\rm(iii)\ }1<\gamma_2\gamma_1^2\stackrel{{}^{\sc\rm\eqref{aaasss--LetNSoOP----1-}\,(ii),\,\eqref{D333aaa}\,(ii)}}
{=}\gamma_4:=\aA_3\xX_2^2\zZ^{-2},
\nonumber\\
&&\!\!\!\!\!\!\!\!\!\!\!\!\!\!\!\!\!\!\!\!\!\!\!\!\!\!
{\rm(iv)\ }1<\gamma_3\gamma_1^{22}
\stackrel{{}^{\sc\rm\eqref{D333aaa}\,(i),\,(ii)}}
=\gamma_5:=\aA_3^{22}\xX_2^{-1}\zZ^2,
\nonumber\\
&&\!\!\!\!\!\!\!\!\!\!\!\!\!\!\!\!\!\!\!\!\!\!\!\!\!\!
{\rm(v)\ }1<\gamma_4\gamma_5=\aA
\end{eqnarray}
We claim
\equa{Clamsmsm}{\dis
a=1+O(\d)^1\mbox{ \ for \ }a\in T_0:=\{\tilde\xX_1,\xX_2,\zZ,\aA_3\}.
}
This holds for $a=\xX_2$ by \eqref{ImMpP}\,(iii),\,{1bab-that}\,(i),
thus also holds for $a=\zZ$ by \eqref{equa-Case6-lemm}\,(iv).
By \eqref{LetNSoOP----1++++++}\,(ii), we have $1=\aA_2\le\frac1{\aA_3(10-9\aA_2^{-20})}+O(\d)^1$, which implies that $\aA_3\le1+O(\d)^1$.
From this and \eqref{1bab-that}\,(ii), we see that  \eqref{Clamsmsm} also holds for $a=\aA_3$, thus also holds for $a=\tilde\xX_1$ by
\eqref{LetNSoOP----1}\,(i). Namely, we have \eqref{Clamsmsm}.
Note that we have in fact frequently used \eqref{equa-Case6-lemm}\,(iv) to
obtain \eqref{Clamsmsm}
. To obtain more precisely results, we cannot use
\eqref{equa-Case6-lemm}\,(iv), instead we need the precise definitions in \eqref{LetNSoOP----1} and the first equality of \eqref{MEMEME-A-inv}.
We have, 
\begin{eqnarray}
\label{1dew2}
&\!\!\!\!\!\!\!\!\!\!\!\!\!\!\!\!\!\!\!\!\!\!\!\!\!\!&
{\rm(i)\ }
1=\aA_2\stackrel{{}^{\sc\rm\eqref{LetNSoOP----1}\,(ii)}}{\le}
\frac{\aA_3^{-10}\xX_2^{90}\zZ^{-90}}{10-9\xX_2^{90}\aA_3^{-20}}
\stackrel{{}^{\sc\rm\eqref{1bab-that}\,(ii)}}{<}
\gamma_2:=\frac{\xX_2^{90}\zZ^{-90}}{10-9\xX_2^{90}}
,\nonumber\\ &\!\!\!\!\!\!\!\!\!\!\!\!\!\!\!\!\!\!\!\!\!\!\!\!\!\!&
{\rm(ii)\ }1=\aA_1^{-1}\stackrel{{}^{\sc\rm\eqref{MEMEME-A-inv}}}{\le}
\frac{\frac{\aA_3^{-10}\xX_2^{-20}\zZ^{-71}}{10}+\frac{9\xX_2^{-20}\zZ^{19}}{10}}{\frac{49}{45}-\frac{4\aA_3^{-30}\xX_2^{90}\zZ^{-90}}{45}}
\stackrel{{}^{\rm\sc\eqref{1bab-that}\,(ii)}}{<}
\gamma_1:=
\frac{\frac{\xX_2^{-20}\zZ^{-71}}{10}+\frac{9\xX_2^{-20}\zZ^{19}}{10}}{\frac{49}{45}-\frac{4\xX_2^{90}\zZ^{-90}}{45}}
.\end{eqnarray}
}\NOUSE{
Then
\equa{B1====}{\dis
\xX_2\stackrel{{}^{\sc\rm\eqref{C1B1}\,(i)}}{=}
(|B_1|^{-1}\aA_1\aA_2^{12}\tilde\xX_1^{-18\ell_0^2-11})^{\frac14}
\stackrel{{}^{\sc\rm\eqref{Formmeme},\,\eqref{1bab-that}}}>1-\d.
}
First assume $\xX_2\ge1$.
Then
\equa{MSMSmAA3}{\dis
\ell_0+1\stackrel{{}^{\sc\rm\eqref{LetNSoOP----1}\,(i)}}{\le}\aA_2^{-1}\zZ\xX_2^{-2}+\ell_0\xX_2^{-\ell_0}
\stackrel{{}^{\sc\rm\eqref{1bab-that}}}{<}\zZ\xX_2^{-2}+\ell_0.
}
This implies that $\xX_2^2<\zZ$ and so $\zZ>1$ and $\xX_2<\zZ$.
We have $\xX_1=\a_0^{-1}\tilde\xX_1=\a_0^{-1}<1<\zZ$ by \eqref{tX1==}.
%
%
We obtain
,
where the first equality follows by noting from \eqref{TaKa},\,\eqref{SimMMSMS},\,\eqref{denote-t2} that $|x_1|=\xX_1\kk,\,|x_2|=\xX_2\kk$,
\equa{Conndndm}{\dis
\g_{|x_1|,|x_2|}
=\g_{\xX_1\kk,\xX_2\kk}\stackrel{{}^{\sc\rm\eqref{wePPPP1+}}}{<}
\g_{\zZ\kk,\xX_2\kk}
\stackrel{{}^{\sc\rm\eqref{wePPPP1}}}{\le}\g_{\zZ\kk,\zZ\kk}
\stackrel{{}^{\sc\rm\eqref{AssG}}}{\le}\zZ\g_{\kk,\kk}
\stackrel{{}^{\sc\rm\eqref{TaKa},\,\eqref{SimMMSMS}}}{=}|x_2+y_2|\stackrel{{}^{\sc\rm\eqref{Ak=1}}}{\le}\g_{|x_1|,|x_2|},
}which is
a contradiction.
Thus $\xX_2<1$. Then $\xX_2$ is a $1+O(\d)^1$ element by \eqref{B1====}. We want to prove
\equa{x2za2}{\dis
a=1+O(\d)^1\mbox{ \ for \ }a\in T_0:=\{\xX_2,\zZ,\aA_2\}.
}
This holds for $a=\zZ$ by \eqref{equa-Case6-lemm}\,(iv). We also have
\equa{aa2===}{\dis
1\stackrel{{}^{\sc\rm\eqref{1bab-that}\,(ii)}}{<}\aA_2\stackrel{{}^{\sc\rm\eqref{LetNSoOP----1}\,(i)}}{\le}
\frac{\zZ}{\xX^{2}(\ell_0+1-\ell_0\xX_2^{-\ell_0})}=1+O(\d)^1.
}
This proves \eqref{x2za2}.
By \eqref{LetNSoOP----1}\,(ii), we have
\equa{AA111}{\dis
1\stackrel{{}^{\sc\rm\eqref{1bab-that}\,(i)}}{=}\aA_1\stackrel{{}^{\sc\rm\eqref{LetNSoOP----1}\,(ii),\,\eqref{1bab-that}\,(i)}}{\le}
\gamma_1:=\aA_2^6\xX_2^{-6}\Big(\frac15+\frac{4\xX_2^{10}\aA_2^{-18}}{5}\Big).
}
}\NOUSE{
Then
\begin{eqnarray}
\!\!\!\!\!\!\!\!\!\!\!\!\!\!
&&
\label{msmA2p}{\rm(i)\ }
|X_1|^{-1}\stackrel{{}^{\sc\rm\eqref{a0BiggerThen}}}{=}|\a_0X_1|^{-1}+O(\d)^3
\stackrel{{}^{\sc\rm\eqref{1bab-that}}}{<}\Big|\frac{A_2}{A_1^{\ell_0^4}(\a_0X_1)Z^3}\Big|+O(\d)^3
\nonumber\\\!\!\!\!\!\!\!\!\!\!\!\!\!\!
&&
\phantom{{\rm(i)\ }|X_1|^{-1}}\stackrel{{}^{\sc\rm\eqref{LetNSoOP----1}\,(iv)}}{=}
|A_2'|+O(\d)^3\stackrel{{}^{\sc\rm\eqref{ImMpP}\,(5)}}{\le}1+O(\d)^1,
\nonumber\\[4pt] \!\!\!\!\!\!\!\!\!\!\!\!\!\!\!\!\!\!\!\!\!&&
{\rm(ii)\ }
\ell_0^4+1-\ell_0^4|X_1|^{-1}\stackrel{{}^{\sc\rm\eqref{ImMpP}\,(8)}}{\le}
|X_1'|+O(\d)^3
\nonumber\\ \!\!\!\!\!\!\!\!\!\!\!\!\!\!\!\!\!\!\!\!\!&&
\phantom{{\rm(ii)\ }ell_0^4+1-\ell}
\stackrel{{}^{\sc\rm\eqref{ImMpP}\,(8),\,\eqref{equa-Case6-lemm}\,(iv)}}{\le}
(1+\d)|Z|^{\ell_0^3}+O(\d)^3\stackrel{{}^{\sc\rm\eqref{1bab-that}}}{\le}1+O(\d)^1.
\end{eqnarray}
This proves \eqref{mdmdthree}\,(i).
The left-hand side of \eqref{msmA2p}\,(ii) is, by \eqref{msmA2p}\,(i), bigger than $1+O(\d)^1$,
which proves \eqref{mdmdthree}\,(ii).
To prove \eqref{mdmdthree}\,(iii), we have
\begin{eqnarray}
\!\!\!\!\!\!\!\!\!\!\!\!\!\!
&&
\label{1+msmA2p}
1=|A_2|\stackrel{{}^{\sc\rm\eqref{ImMpP}\,(3),\,\eqref{mdmdthree}\,(i)}}{\le}
|X_2|\Big(\frac14+\frac{3|X_2|^2}{3}\Big)+O(\d)^3,
\end{eqnarray}
thus $|X_2|\ge1+O(\d)^1$, which with \eqref{equa-Case6-lemm}\,(iv),\,\eqref{1bab-that}\,(ii) gives \eqref{mdmdthree}\,(iii). Thus
\equa{mdmdthree}{\dis\!\!\!\!
{\rm(i)\ }|X_1|=1+O(\d)^1, \ \ {\rm(ii)\ }|X_1'|=1+O(\d)^1,\ \ {\rm(iii)\ }
|X_2|\stackrel{{}^{\sc\rm\eqref{equa-Case6-lemm}\,(iv)}}{=} |Z|+O(\d)^1=1+O(\d)^1.\!\!\!
}
}\NOUSE{
 by ,\,\eqref{1-dSmallerThan}, we obtain
[cf.~\eqref{x1xxxxxxx}$\ssc\,$],
\equa{that-XXX1}{\mbox{ $|X_1|\stackrel{{}^{\sc\rm\eqref{1bab-that},\,\eqref{1-dSmallerThan}}}{=}\ \,|A_3|^{2\d_0^2}+O(\d)^1$.}} Then
\eqref{ImMpP}\,(4) gives [using \eqref{A0===+1}$\ssc\,$],
\equa{AAAA22222+1}{\dis
1=|A_1|\ \,\stackrel{{}^{\sc\rm\eqref{ImMpP}\,(4),\,\eqref{that-XXX1}}}{\le}\ \,\frac1{|A_3|^{2\d_0^2}}\Big(\frac{1-\d_0^2}{2}+
\frac{(1+\d_0^2)|A_3|^{2\d_0^2}}{2}\Big)+O(\d)^1,}
which shows that $|A_3|\le1+O(\d)^1$. This with \eqref{1bab-that}\,(ii)
 gives that $|A_3|=1+O(\d)^1$. This
with \eqref{ImMpP}\,(2),\,(6),\,\eqref{that-XXX1} proves that $|X_1|,|X_2|,|Z|$ are $1+O(\d)^1$ elements.
}%
\NOUSE{
 In particular, by \eqref{LetNSoOP----1}\,(iii), we have
\equa{A3isBigger}{\dis
|X_2|^{1-\d_0+\d_0^2}=|A_3^{\d_0(1-\d_0)}Z|<|Z|.}
By \eqref{ImMpP}\,(3),\,(6), we have $|X_2|\ge1+O(\d)^1$, but by \eqref{ImMpP}\,(7), we have
$|X_2|\le1+O(\d)^1$. Thus $|X_1|=1+O(\d)^1$ and so $|X_2|=1+O(\d)^1$ by
\eqref{ImMpP}\,(5).
This with \eqref{ImMpP}\,(2) gives
}}\NOUSE{
By \eqref{LetNSoOP----1}\,(ii),\,(iii), we obtain
\begin{eqnarray}
\label{x0s,dm}
&\!\!\!\!\!\!\!\!\!\!\!\!\!\!\!\!\!\!\!\!\!\!\!\!\!&
1<|A_3|=\frac{1}{|\tildeX_1X_2|}.
\end{eqnarray}
Further, by the first inequalities of
\eqref{ImMpP}\,(13),\,(14), we obtain that $|\tildeX_1|,\,|X_2|\ge1+O(\d)^1$. This together with
\eqref{x0s,dm},$ $\eqref{ImMpP}\,(2) shows that we have the following,
}%
\NOUSE{From this and by the first inequality of \eqref{x1-x2==}\,(ii) and \eqref{tilde-x1}\,(i), we see that $\tilde\xX_1$ is a $1+O(\d)^1$ element.
This with \eqref{tilde-x1}\,(ii) and the equality in \eqref{tilde-x1}\,(iii) shows that $\xX_2$ is also a $1+O(\d)^1$. So is $\zZ$ by \eqref{equa-Case6-lemm}\,(iv). Then by
\eqref{equa-Case6-lemm}\,(i) and \eqref{1bab-that}\,(i), we obtain
%
%
\equa{T2-isO(ee)}{\dis a=1+O(\d)^1\mbox{ \ for \ }a\in T_3:=\{\tilde\xX_1,\,
\xX_2,\,\zZ,\,\aA_3,\a_1\}.}
}\NOUSE{which obviously holds for $\a_1$ as it is $1+O(\ep)^1$ element  by \eqref{LetNSoOP----1}\,(i).
To see it holds for $\tilde\xX_1$, first 
by
the last inequality of \eqref{A2-B1aaaaaa}\,(ii)
and \eqref{1bab-that}, we obtain that $\tilde\xX_1\le1+O(\eE_1)^1$,
more precisely,
\equa{mm4m5m5}{\dis
\tilde\xX_1\stackrel{{}^{\sc\rm\eqref{A2-B1aaaaaa}\,(ii)}}{<}(1+\eE_1)^{\frac1{808}}<1+\eE_1
\mbox{ \ [up to $O(\eE_1)^2$]}.}
To prove that $\tilde\xX_1$ is a $1+O(\eE_1)^1$ element, it remains to prove that $\tilde\xX_1\ge1+O(\eE_1)^1$.
It suffices to  prove that $\tilde\xX_1\ge1-6\eE_1.$ Thus  assume conversely,  \equa{Assumemmmmm}{\dis
\tilde\xX_1<1-6\eE_1<1.} Then the first equality in \eqref{X2==mmmmm} with \eqref{1bab-that},\,\eqref{Assumemmmmm} shows
\equa{Smhhwx222is}{\dis \xX_2\stackrel{{}^{\sc\rm\eqref{X2==mmmmm},\,\eqref{(T0o(eE)1}}}{=}
\aA_2^{\frac{\d}{10}}\aA_3^{-(\frac{\d\lL}{160}+\frac{321}{16})}\tilde\xX_1^{\frac15+\frac{\d}{10}}
\stackrel{{}{\sc\rm\eqref{1bab-that}}}{\le}\tilde \xX_1^{\frac15+\frac{\d}{10}}+O(\eE_1)^2
\stackrel{{}^{\sc\rm\eqref{Assumemmmmm}}}{<}1-\eE_1<1.}
If $\zZ\ge1$, then this together with \eqref{Assumemmmmm},\,\eqref{Smhhwx222is} implies that $\xX_1<\zZ,\,\xX_2<\zZ$ and
we can obtain a contradiction as in \eqref{Conndndm} (cf.~Remark \ref{Need-to-Rem} when
$\zZ=1$).
Thus $\zZ<1$.
We want to verify whether or not \eqref{TheFasss} holds. To do this, we use notation there
by denoting $\hat k=\zZ$ and so $|x_2+y_2|=
\hat k\gamma_{\kk,\kk}$ by
\eqref{TaKa},\,\eqref{SimMMSMS}, and $k_1:=\xX_1$ with $|x_1|=k_1\kk$, and $k_2:=\xX_2$
with $|x_2|=k_2\kk$, and further
 we can obtain, where the last equality of (i) follows from the arguments after \eqref{1bab-that},
\begin{eqnarray}
\label{Ckemmmdmdmd}
&\!\!\!\!\!\!\!\!\!\!\!\!\!\!\!\!\!\!\!\!&
{\rm(i)\ }1>\hat k=\zZ\stackrel{{}^{\sc\rm\eqref{LetNSoOP----1}\,(i),\,\eqref{1bab-that}\,(ii)}}{=}
1+O(\ep)^1,
\nonumber\\
&\!\!\!\!\!\!\!\!\!\!\!\!\!\!\!\!\!\!\!\!&
{\rm(ii)\ }k_2=\xX_2\stackrel{{}^{\sc\rm\eqref{Smhhwx222is}}}{<}
1-\eE_1\stackrel{{}^{\sc\rm\eqref{Ckemmmdmdmd}\,(i)}}{<}
\hat k^{1+\l_2}\mbox{ \ \ with, say,  \ }\l_2=1,\ \ \ \implies\ \ \ {\rm(iii)\ }|x_2|\stackrel{{}^{\sc\rm\eqref{TaKa},\,\eqref{SimMMSMS}}}{=}k_2\kk<\kk,\!\!\!\!\!\!\!\!\!
\nonumber\\
&\!\!\!\!\!\!\!\!\!\!\!\!\!\!\!\!\!\!\!\!&
{\rm(iv)\ }k_1=\xX_1\stackrel{{}^{\sc\rm\eqref{x1-small}\,(i)}}{<}
\xX_2^{\frac15}
\tilde \xX_1^{\frac1{100}}\stackrel{{}^{\sc\rm\eqref{Assumemmmmm},\,\eqref{Smhhwx222is}}}{<}
1,\ \ \ \implies\ \ \ {\rm(v)\ }|x_1|\stackrel{{}^{\sc\rm\eqref{TaKa},\,\eqref{SimMMSMS}}}{=}k_1\kk<\kk,
\nonumber\\
&\!\!\!\!\!\!\!\!\!\!\!\!\!\!\!\!\!\!\!\!&
{\rm(vi)\ }k_2\ge\l:=\frac12,
\end{eqnarray}
where \eqref{Ckemmmdmdmd}\,(vi) is proven as follows (cf.~Remark \ref{Remmma}): Assume $k_2<\l$, then by \eqref{TaKa},\,\eqref{SimMMSMS} and notation \eqref{denote-t2}, we have $|x_2|=k_2\kk<\frac{\kk}{2}$, and we obtain
[noting that we can choose $\kk$ in \eqref{MSmde33333} to be sufficiently larger than $\SS_7$ with $\SS_7$ being the number in Proposition \ref{Sect3-Lemm5} so that
by \eqref{msms-a-a-a}\,(i) we can apply
\eqref{mqp1234-2}; where in (ii), $\d=\ell^{-1}$ is the number defined in \eqref{MSmde33333}$\ssc\,$],
\begin{eqnarray}
\label{msms-a-a-a}
&\!\!\!\!\!\!\!\!\!\!\!\!\!\!\!\!\!\!\!\!\!\!\!\!\!\!\!&
{\rm(i)\ }\kk+O(\ep)^1\stackrel{{}^{\sc\rm\eqref{Ckemmmdmdmd}\,(i)}}{=}\hat k\kk\stackrel{{}^{\sc\rm\eqref{ggggg}}}{<}\hat k\g_{\kk,\kk}\stackrel{{}^{\sc\rm\eqref{TaKa},\,\eqref{SimMMSMS}}}{=}
|x_2+y_2|\le|x_2|+|y_2|\stackrel{{}^{\sc\rm\eqref{dp0p1}}}{\le} h_{p_1,p_2},
\nonumber\\
&\!\!\!\!\!\!\!\!\!\!\!\!\!\!\!\!\!\!\!\!\!\!\!\!\!\!\!&
{\rm(ii)\ }|y_i|\stackrel{{}^{\sc\rm\eqref{mqp1234-2}}}{<}h_{_{\sc p_1,p_2}}^{^{\sc\frac m{m+1}}}\stackrel{{}^{\sc\rm\eqref{msms-a-a-a}\,(i),\,\eqref{MSmde33333}}}{<}\d h_{p_1,p_2}\mbox{ \ for \ }i=1,2,
\nonumber\\
&\!\!\!\!\!\!\!\!\!\!\!\!\!\!\!\!\!\!\!\!\!\!\!\!\!\!\!&
{\rm(iii)\ }h_{p_1,p_2}\stackrel{{}^{\sc\rm\eqref{dp0p1}}}{=}|x_1|+|x_2|+|y_1|+|y_2|\stackrel{{}^{\sc\rm\eqref{Ckemmmdmdmd}\,(iii),\,(v),\,\eqref{msms-a-a-a}\,(ii)}}{<}
2\kk+2\d h_{p_1,p_2},
\nonumber\\
&\!\!\!\!\!\!\!\!\!\!\!\!\!\!\!\!\!\!\!\!\!\!\!\!\!\!\!&
{\rm(iv)\ }\frac{\kk}{2}\stackrel{{}^{\sc\rm\eqref{msms-a-a-a}\,(i)}}{<}h_{p_1,p_2}\stackrel{{}^{\sc\rm\eqref{msms-a-a-a}\,(iii)}}{<}3\kk,
\nonumber\\
&\!\!\!\!\!\!\!\!\!\!\!\!\!\!\!\!\!\!\!\!\!\!\!\!\!\!\!&
{\rm(v)\ }|y_2|\ge|x_2+y_2|-|x_2|\stackrel{{}^{\sc\rm\eqref{msms-a-a-a}\,(i),\,}}{\ge}
\kk-\l_1\kk+O(\ep)^1\ge\frac{\kk}{2}+O(\ep)^1\stackrel{{}^{\sc\rm\eqref{msms-a-a-a}\,(iv)}}{>}h_{_{\sc p_1,p_2}}^{^{\sc\frac m{m+1}}},
\end{eqnarray}
where the last inequality of \eqref{msms-a-a-a}\,(ii) is obtained by noting from \eqref{MSmde33333}
that we may assume $0<\kk^{-\frac1{m+1}}\ll\d=\ell^{-1}$.
In particular, \eqref{msms-a-a-a}\,(v) is a contradiction with \eqref{mqp1234-2}. This proves \eqref{Ckemmmdmdmd}\,(vi).
Now \eqref{Ckemmmdmdmd} shows that \eqref{TheFasss} holds, which contradicts the assumption that Case 5 does not occur (at this point it may be worth mentioning that when we say Case 5 does not occur, it means that if it occurs we have already proved that Proposition \ref{real00-inj} holds).
This proves that \eqref{T2-isO(ee)} holds for $\tilde\xX_1$.
Using \eqref{1bab-that},\,\eqref{A2-B1aaaaaa}\,(i)
and the fact that $\tilde\xX_1=1+O(\eE_1)^1$, we see that $\xX_2\le1+O(\eE_1)^1$.
If, say, $\xX_2<
1-10000\eE_1
$, then this with
\eqref{mm4m5m5},\,\eqref{A2-B1aaaaaa}\,(i),\,\eqref{LetNSoOP----1}\,(i) shows
that $\zZ>1>\xX_2$ and then by \eqref{x1-small}\,(i), we also have
$\zZ>\xX_1$, and we obtain a contradiction again
as in \eqref{Conndndm}. Thus $\xX_2$ is also a
$1+O(\eE_1)^1$ element. Then by \eqref{A2-B1aaaaaa}\,(i),\,\eqref{LetNSoOP----1}\,(i), $\aa_3,\zZ$ are
$1+O(\eE_1)^1$ elements too, i.e., we have \eqref{T2-isO(ee)}.
}{\NOUSE{
\eqref{ImMpP}\,(2), $\tilde\xX_1\le1+O(\eE_1)^1$.
If $\tilde\xX_1<1-\eE_1$, then $\xX_2>1$ and $\xX_2>\a_1^{\lL^5(2 + \dD)}\tilde \xX_1^{2\lL^4  (2 - \dD + \dD^2)}$ by
\eqref{LetNSoOP----1}\,(ii),\,\eqref{1bab-that}\,(ii),
 and we obtain from \eqref{LetNSoOP----1}\,(iii),\,\eqref{1bab-that}\,(ii)
that $\zZ>\xX_2$. Then we obtain a contradiction as in \eqref{Conndndm}. Thus
\eqref{OSOSOSO} holds for $\tilde\xX_1$. Then by
\eqref{1bab-that}\,(ii), we see that $\xX_2\ge1+O(\eE_1)^1$.
If, say, $\xX_2>(1+\eE_1)^{\lL^2}$, then as above, we again have that $\zZ>\xX_2>\xX_1$ and $\zZ>1$
and we obtain a contradiction as in
 \eqref{Conndndm}. Thus \eqref{OSOSOSO} holds for $\xX_2$. Then
\eqref{LetNSoOP----1}\,(ii) with \eqref{1bab-that}\,(ii)  shows that
\eqref{OSOSOSO} holds for $\aA_3$. Thus \eqref{OSOSOSO} holds.
Of course $\a_1>1$ is a $1+O(\eE_1)^1$ element by \eqref{LetNSoOP----1}\,(v).
}%
\NOUSE{Now we need the precise definitions in \eqref{LetNSoOP----1}. We have
\begin{eqnarray}
\label{LetNSoOP----1**}
&\!\!\!\!\!\!\!\!\!\!\!\!\!\!\!\!\!\!\!\!\!\!\!\!\!\!\!\!\!\!\!\!\!\!\!\!\!\!
&
{\rm(i)\ }1<\aA_3\stackrel{{}^{\sc\rm\eqref{LetNSoOP----1}\,(i)}}{=}
\gamma_3:=\zZ\tilde\xX_1
,
\!\!\!\!\!\!\!\!\!\!\!\!\!\!\!\!
\nonumber\\
&\!\!\!\!\!\!\!\!\!\!\!\!\!\!\!\!\!\!\!\!\!\!\!\!\!
&{\rm(ii)\ }
1\ \ =\ \ \aA_2\stackrel{{}^{\sc\rm\eqref{LetNSoOP----1}\,(iii)}}{\le}
\tilde\gamma_{20}:=\a_1^4\aA_3^{45}\tilde\xX_1^6\xX_2^{-45}\Big(\frac16+\frac56\a_1^{-6}\aA_3^{-114}\tilde\xX_1^{18}\xX_2^{66}\Big),
\nonumber\\
&\!\!\!\!\!\!\!\!\!\!\!\!\!\!\!\!\!\!\!\!\!\!\!\!\!
&\phantom{{\rm(ii)\ }}
\stackrel{{}^{\sc\rm\eqref{1bab-that}\,(ii)}}{<}\gamma_{20}:=
\a_1^4\tilde\xX_1^6\xX_2^{-45}\Big(\frac16+\frac56\a_1^{-6}\tilde\xX_1^{18}\xX_2^{66}\Big),
\!\!\!\!\!\!\!\!\!\!\!\!\!\!\!\!
\nonumber\\
&\!\!\!\!\!\!\!\!\!\!\!\!\!\!\!\!\!\!\!\!\!\!\!\!\!
&
 {\rm(iii)\ }1=\aA_1\stackrel{{}^{\sc\rm\eqref{LetNSoOP----1}\,(ii)}}{\le}\gamma_{10}:=
\a_1^6\aA_3^{114}\tilde\xX_1^{-18}\xX_2^{-66}\Big(\frac15+\frac45\aA_3^{-150}\tilde\xX_1^{-150}\Big)
\!\!\!\!\!\!\!\!\!\!\!\!\!\!\!\!
\nonumber\\
&\!\!\!\!\!\!\!\!\!\!\!\!\!\!\!\!\!\!\!\!\!\!\!\!\!
&
\phantom{{\rm(iii)\ }1}
<\gamma_0:=
\a_1^6\tilde\xX_1^{-18}\xX_2^{-66}\Big(\frac15+\frac45\tilde\xX_1^{-150}\Big)
 ,\end{eqnarray}
where the strict inequality in (iii) is obtained as follows:
}%
\NOUSE{
Regarding $\gamma_{1}$ is a function on $\xX_2,\zZ$
in the small neighborhood $\mathcal O_1$ of $(1,1)$, where ${\mathcal O}_1$ is defined below%
,
\equa{define-O1}{\dis
{\mathcal O}_1=\Big\{(\xX_2,\zZ)\in\C^2\ \Big|\ |\xX_2-1|<\d_2,\,|\zZ-1|
<\d_2\Big\},
}
one can easily see from \eqref{1dew2}\,(ii),
\equa{SeeeotoSee}{\dis
\frac{\ptl\gamma_{1}}{\ptl\xX_2}\mbox{\Large$\Big|$}_{(\xX_2,\zZ)=(1,1)}=
\frac1{10}\times(-20)+\frac9{10}\times(-20)+\frac4{45}\times(90)=-2-18+8=-12<0.
}
Thus
[noting from \eqref{MSmde33333} that $0<\d=\ell^{-1}\ll\d_2=\ell_2^{-1}\ll1$],
\equa{MDNDNN9999}{\dis \frac{\ptl\gamma_{1}}{\ptl\xX_2}<0 \mbox{ \ \ in \ \ }{\mathcal O}_1,}
which means that $\gamma_{1}$ is a strictly
decreasing function on $\xX_2$
 when $(\xX_2,\zZ)\in{\mathcal O}_1$. 
We claim,
\equa{Claimse-111}{\dis
\xX_2<\zZ^{\frac16+s\d_1}, \mbox{ \ where $s=1$ if $\zZ\ge1$ or $s=-1$ if $\zZ<1$}.
}
Assume this does not hold. Then by \eqref{1dew2}\,(ii) and \eqref{MDNDNN9999} [noting from \eqref{Clamsmsm} that we
indeed have $(\xX_2,\zZ)\in{\mathcal O}$], we have
\equa{smemem334}{\dis
1<\gamma_1<\gamma_{10}:=\gamma_1|_{\xX_2=\zZ^{\frac16+s\d_1}}=
\frac{\frac{\zZ^{-\frac{223}{3}-20s\d_1}}{10}+\frac{9\zZ^{\frac{47}{3}-20s\d_1}}{10}}{\frac{49}{45}-\frac{4\zZ^{-75+90s\d_1}}{45}}.
}
One can easily compute,
\begin{eqnarray}
\label{dm,ememe}
&\!\!\!\!\!\!\!\!\!\!\!\!\!\!\!\!\!\!\!\!\!&
\frac{d\gamma_{10}}{d\zZ}\mbox{\Large$\Big|$}_{\zZ=1}=\frac1{10}\Big(-\frac{223}{3}-20s\d_1\Big)
+\frac9{10}\Big(\frac{47}{3}-20s\d_1\Big)+\frac4{45}(-75+90s\d_1)
\nonumber\\
&\!\!\!\!\!\!\!\!\!\!\!\!\!\!\!\!\!\!\!\!\!&
\phantom{\frac{d\gamma_{10}}{d\zZ}\mbox{\Large$\Big|$}_{\zZ=1}}
=-12s\d_1,
\mbox{ which is $<0$ if $\zZ\ge1$, or
$>0$ if $\zZ<1.$}
\end{eqnarray}
Note that $0<\d_2\ll\d_1$, we see from
\eqref{dm,ememe} that when $|\zZ-1|<\d_2$, $\gamma_{10}$ is a strictly decreasing function on $\zZ$ if $\zZ\ge1$, or a
strictly increasing function on $\zZ$ if $\zZ<1$.
Thus in any case we obtain from \eqref{smemem334} 
that
$1<\gamma_{10}\le\gamma_{10}|_{\zZ=1}=1$, a contradiction.
This proves \eqref{Claimse-111}, which in particular implies that $\xX_2<\zZ^{\frac16+O(\d_1)^1}$.
Similarly, we can use \eqref{1dew2}\,(i) to prove that $\zZ<\xX_2^{10+O(\d_1)^1}<\zZ^{\frac53+O(\d_1)^1}$.
Thus we obtain that $1<\xX_2<\zZ$. We also have
$\xX_1=\a_0^{-1}\tilde\xX_1<\aA_3^{-1}<1<\zZ$ by \eqref{tX1==},\,\eqref{LetNSoOP----1}\,(i),\,\eqref{1bab-that}\,(ii).
}%
}%
\NOUSE{We obtain
,
where the first equality follows by noting from \eqref{TaKa},\,\eqref{SimMMSMS},\,\eqref{denote-t2} that $|x_1|=\xX_1\kk,\,|x_2|=\xX_2\kk$,
\equa{Conndndm}{\dis
\g_{|x_1|,|x_2|}
=\g_{\xX_1\kk,\xX_2\kk}\stackrel{{}^{\sc\rm\eqref{wePPPP1+}}}{<}
\g_{\zZ\kk,\xX_2\kk}
\stackrel{{}^{\sc\rm\eqref{wePPPP1}}}{\le}\g_{\zZ\kk,\zZ\kk}
\stackrel{{}^{\sc\rm\eqref{AssG}}}{\le}\zZ\g_{\kk,\kk}
\stackrel{{}^{\sc\rm\eqref{TaKa},\,\eqref{SimMMSMS}}}{=}|x_2+y_2|\stackrel{{}^{\sc\rm\eqref{Ak=1}}}{\le}\g_{|x_1|,|x_2|},
}which is
a contradiction.
}\NOUSE{
By \eqref{x2za2}, we indeed have
$(\aA_2,\xX_2)\in{\mathcal O}_1$. Thus by
\eqref{MDNDNN9999} and the fact that
$\aA_3>1$, we obtain
\equa{DMDMememe}{\dis
1\le\gamma_{1}<\gamma_{1}|_{\aA_3=1}=\gamma_0:=
\xX_2^{-6}\Big(\frac15+\frac{4\xX_2^{10}}{5}\Big).}
Note that
$\gamma_0$ is a strictly increasing function on $\xX_2$ when $1-\d_2<\xX_2<1+\d_2$.
Thus \eqref{DMDMememe} implies that $\xX_2>1$, a contradiction.
}\NOUSE{
, i.e., we have the strict inequality in
\eqref{LetNSoOP----1**}\,(iii).
The proof of the last equality of \eqref{LetNSoOP----1**}\,(ii) is exactly similar by noting that
$\frac{\ptl\tilde\gamma_{20}}{\ptl\aA_3}|_{(\aA_1,\aA_3,\tilde\xX_1,\xX_2)=(1,1,1,1)}=-50$.
Now we regard $\gamma_0$ in \eqref{LetNSoOP----1**}\,(iii) as a  function on $\a_1,\tilde\xX_1,\xX_2$
in the small neighborhood ${\mathcal O}_2$ of $(1,1,1)$,
 where ${\mathcal O}_2$ is defined below%
,
\equa{define-O}{\dis
{\mathcal O}_2=\Big\{(\a_1,\tilde\xX_1,\xX_2)\in\C^3\ \Big|\ |\a_1-1|<\d_2,\,|\tilde\xX_1-1|<\d_2,\,|\xX_2-1|<\d_2\}.
}
Then  at point
$(1,1,1)$, we have
\equa{ap1-aa3-xx1}{\dis
\frac{\ptl\gamma_0}{\ptl\a_1}
=6,\ \ \ \ \
\frac{\ptl\gamma_0}{\ptl\tilde\xX_1}
=-138,\ \ \ \ \
\frac{\ptl\gamma_0}{\ptl\xX_2}
=-66.
}
This suggests us to denote
(where $s_{1}=1$ if $\tilde\xX_1\ge1$ and $s_{1}=-1$ otherwise, and $s_{2}=1$ if $\xX_2\ge1$ and $s_{2}=-1$ otherwise),
\equa{Dmemennna}{\dis{\rm(i)\ }
\gamma_1:=\a_1^{6+\d_1}
\tilde\xX_1^{-138+s_{1}\d_1}\xX_2^{-66+s_{2}\d_1},\ \ \ \ \
{\rm(ii)\ }\tilde\gamma_1=\gamma_1^{-1}\gamma_0.
}
Then $\tilde\gamma_1$ as a function on $\a_1,\tilde\xX_1,\xX_2$ in  ${\mathcal O}_2$,
we see from \eqref{ap1-aa3-xx1},\,\eqref{Dmemennna} that when $(\a_1,\tilde\xX_1,\zZ)\in{\mathcal O}$, %
  $\tilde\gamma_1$ is 
a strictly decreasing function on $\a_1$, and it is a strictly decreasing function on $\tilde\xX_1$ if $\tilde\xX_1\ge1$ or a
 strictly increasing function on $\tilde\xX_1$ if $\tilde\xX_1<1$, by our choice of $s$
 [noting again from \eqref{MSmde33333} that $0\ll\d\ll\d_2\ll\d_1\ll1$], and further
  it is a strictly decreasing function on $\xX_2$ if $\xX_2\ge1$ or a
 strictly increasing function on $\xX_2$ if $\xX_2<1$
   Thus in any case, we obtain that
$\tilde\gamma_1<\tilde\gamma_1|_{(\a_1,\tilde\xX_1,\xX_2)=(1,1,1)}=1$ [noting from \eqref{LetNSoOP----1}\,(ii) that $\a_1>1$
and that we indeed have $(\a_1,\tilde\xX_1,\xX_2)\in{\mathcal O}$ by \eqref{T2-isO(ee)}$\ssc\,$].
This with \eqref{LetNSoOP----1**}\,(iii),\,\eqref{Dmemennna}\,(ii) proves
\equa{SMSMSMnen3n3333}{\dis
1<\gamma_1=
\a_1^{6+\d_1}
\tilde\xX_1^{-138+s_{1}\d_1}\xX_2^{-66+s_{2}\d_1}.
}
Now by \eqref{LetNSoOP----1**}\,(ii), at point $(1,1,1)$, we have
\equa{gamm2-0-}{\dis\frac{\ptl\gamma_{20}}{\ptl\a_1}=-1,\ \ \ \
\frac{\ptl\gamma_{20}}{\ptl\tilde\xX_1}=21,\ \ \ \ \frac{\ptl\gamma_{20}}{\ptl\xX_2}=10.
}
Thus as above, we obtain
\equa{+SMSMSMnen3n3333}{\dis
1<\gamma_2:=
\a_1^{-1+\d_1}
\tilde\xX_1^{21+s_{1}\d_1}\xX_2^{10+s_{2}\d_1}\stackrel{{}^{\sc\rm\eqref{LetNSoOP----1}\,(ii)}}{<}
\tilde\xX_1^{21+s_{1}\d_1}\xX_2^{10+s_{2}\d_1}.
}
Now we can evaluate the following [where (iii) is obtained by using (ii) in the right-hand side of \eqref{+SMSMSMnen3n3333}$\ssc\,$],
\begin{eqnarray}
\label{LetNSoOP----1**+1}
&\!\!\!\!\!\!\!\!\!\!\!\!\!\!\!\!\!\!\!\!\!\!\!\!\!\!\!\!\!\!\!\!\!\!\!\!\!\!
&
{\rm(i)\ }1<
\gamma_1^{1-\d_1}\gamma_2^{6+\d_1}
\stackrel{{}^{\sc\rm\eqref{SMSMSMnen3n3333},\,\eqref{+SMSMSMnen3n3333}}}{=}
\gamma_{12}:=\tilde\xX_1^{-12+O(\d_1)^{1}}\xX_2^{-6+O(\d_1)^1}
,\ \  \implies\  \ {\rm(ii)\ }\xX_2<\tilde\xX_1^{-2+O(\d_1)^1},
\!\!\!\!\!\!\!\!\!\!\!\!\!\!\!\!
\nonumber\\
&\!\!\!\!\!\!\!\!\!\!\!\!\!\!\!\!\!\!\!\!\!\!\!\!\!
&
 {\rm(iii)\ }1\stackrel{{}^{\sc\rm\eqref{+SMSMSMnen3n3333},\,\eqref{LetNSoOP----1**+1}\,(ii)}}{<}\widetilde\xX_1,
\ \ \ \ \ \ \ \ {\rm(iv)\ }\xX_2\stackrel{{}^{\sc\rm\eqref{LetNSoOP----1**+1}\,(ii),\,(iii)}}{<}1,
\!\!\!\!\!\!\!\!\!\!\!\!\!\!\!\!
\nonumber\\
&\!\!\!\!\!\!\!\!\!\!\!\!\!\!\!\!\!\!\!\!\!\!\!\!\!
&
{\rm(v)\ }
\xX_1\stackrel{{}^{\sc\rm\eqref{tX1==}}}{=}\a_0^{-1}\tilde\xX_1\xX_2^3<\tilde\xX_1\xX_2^3
\stackrel{{}^{\sc\rm\eqref{LetNSoOP----1**+1}\,(ii)}}{<}\xX_2^{-\frac12+O(\d_1)^1+3}\stackrel{{}^{\sc\rm\eqref{LetNSoOP----1**+1}\,(iv)}}{<}1.
\end{eqnarray}
First assume  $\zZ\ge1$. Then $\xX_1<\zZ,\,\xX_2<\zZ$ by \eqref{LetNSoOP----1**+1}\,(iv),\,(v).
We obtain (see Remark \ref{Need-to-Rem}),
where the first equality follows by noting from \eqref{TaKa},\,\eqref{SimMMSMS},\,\eqref{denote-t2} that $|x_1|=\xX_1\kk,\,|x_2|=\xX_2\kk$,
\equa{Conndndm}{\dis
\g_{|x_1|,|x_2|}
=\g_{\xX_1\kk,\xX_2\kk}\stackrel{{}^{\sc\rm\eqref{wePPPP1+}}}{<}
\g_{\zZ\kk,\xX_2\kk}
\stackrel{{}^{\sc\rm\eqref{wePPPP1}}}{\le}\g_{\zZ\kk,\zZ\kk}
\stackrel{{}^{\sc\rm\eqref{AssG}}}{\le}\zZ\g_{\kk,\kk}
\stackrel{{}^{\sc\rm\eqref{TaKa},\,\eqref{SimMMSMS}}}{=}|x_2+y_2|\stackrel{{}^{\sc\rm\eqref{Ak=1}}}{\le}\g_{|x_1|,|x_2|},
}which is
a contradiction.
%
%
\begin{rema}\rm\label{Need-to-Rem}
We remark that as long as $\xX_1,\xX_2<\zZ$ and $\zZ\ge1$, we can obtain
\eqref{Conndndm}. Note that we require $\zZ>1$ to apply \eqref{AssG}, however in our case here, the third inequality of \eqref{Conndndm} is simply an
equality when  $\zZ=1$.\end{rema}
Thus $\zZ<1$.
We want to verify whether or not \eqref{TheFasss} holds. To do this, we use notation there
by denoting $\hat k=\zZ$ and so $|x_2+y_2|=
\hat k\gamma_{\kk,\kk}$ by
\eqref{TaKa},\,\eqref{SimMMSMS}, and $k_1:=\xX_1$ with $|x_1|=k_1\kk$, and $k_2:=\xX_2$
with $|x_2|=k_2\kk$, and further
 we can obtain, where the last equality of (i) follows from the arguments after \eqref{1bab-that},
\begin{eqnarray}
\label{Ckemmmdmdmd}
&\!\!\!\!\!\!\!\!\!\!\!\!\!\!\!\!\!\!\!\!&
{\rm(i)\ }1>\hat k=\zZ
\stackrel{{}^{\sc\rm\eqref{T2-isO(ee)}}}{=}
1+O(\d)^1,
\nonumber\\
&\!\!\!\!\!\!\!\!\!\!\!\!\!\!\!\!\!\!\!\!&
{\rm(ii)\ }k_2=\xX_2\stackrel{{}^{\sc\rm\eqref{LetNSoOP----1**+1}\,(ii)}}{<}
\tilde\xX_1^{-2+O(\d_1)^1}
\stackrel{{}^{\sc\rm\eqref{LetNSoOP----1}\,(i)}}
{=}(\zZ\aA_3^{-1})^{2+O(\d_1)^1}
\stackrel{{}^{\sc\rm\eqref{1bab-that}\,(iii)}}{<}
\hat k^{1+\l_2}\mbox{ \ \ with, say,  \ }\l_2=\frac12
,\  \implies\!\!\!\!\!\!\!\!\!\!\!\!\!\!
\nonumber\\
&\!\!\!\!\!\!\!\!\!\!\!\!\!\!\!\!\!\!\!\!&
{\rm(iii)\ }|x_2|\stackrel{{}^{\sc\rm\eqref{TaKa},\,\eqref{SimMMSMS}}}{=}k_2\kk<\kk,\!\!\!\!\!\!\!\!\!
\nonumber\\
&\!\!\!\!\!\!\!\!\!\!\!\!\!\!\!\!\!\!\!\!&
{\rm(iv)\ }k_1=\xX_1\stackrel{{}^{\sc\rm\eqref{LetNSoOP----1**+1}\,(v)}}{<}
1,\ \ \ \implies\ \ \ {\rm(v)\ }|x_1|\stackrel{{}^{\sc\rm\eqref{TaKa},\,\eqref{SimMMSMS}}}{=}k_1\kk<\kk,
\nonumber\\
&\!\!\!\!\!\!\!\!\!\!\!\!\!\!\!\!\!\!\!\!&
{\rm(vi)\ }k_2=\xX_2\stackrel{{}^{\sc\rm\eqref{equa-Case6-lemm}\,(iv)}}{=}\zZ+O(\d)^2\ge\l:=\frac12.
\end{eqnarray}
}\NOUSE{where \eqref{Ckemmmdmdmd}\,(vi) is proven as follows (cf.~Remark \ref{Remmma}): Assume $k_2<\l$, then by \eqref{TaKa},\,\eqref{SimMMSMS} and notation \eqref{denote-t2}, we have $|x_2|=k_2\kk<\frac{\kk}{2}$, and we obtain
[noting that we can choose $\kk$ in \eqref{MSmde33333} to be sufficiently larger than $\SS_7$ with $\SS_7$ being the number in Proposition \ref{Sect3-Lemm5} so that
by \eqref{msms-a-a-a}\,(i) we can apply
\eqref{mqp1234-2}; where in (ii), $\d=\ell^{-1}$ is the number defined in \eqref{MSmde33333}$\ssc\,$],
\begin{eqnarray}
\label{msms-a-a-a}
&\!\!\!\!\!\!\!\!\!\!\!\!\!\!\!\!\!\!\!\!\!\!\!\!\!\!\!&
{\rm(i)\ }\kk+O(\ep)^1\stackrel{{}^{\sc\rm\eqref{Ckemmmdmdmd}\,(i)}}{=}\hat k\kk\stackrel{{}^{\sc\rm\eqref{ggggg}}}{<}\hat k\g_{\kk,\kk}\stackrel{{}^{\sc\rm\eqref{TaKa},\,\eqref{SimMMSMS}}}{=}
|x_2+y_2|\le|x_2|+|y_2|\stackrel{{}^{\sc\rm\eqref{dp0p1}}}{\le} h_{p_1,p_2},
\nonumber\\
&\!\!\!\!\!\!\!\!\!\!\!\!\!\!\!\!\!\!\!\!\!\!\!\!\!\!\!&
{\rm(ii)\ }|y_i|\stackrel{{}^{\sc\rm\eqref{mqp1234-2}}}{<}h_{_{\sc p_1,p_2}}^{^{\sc\frac m{m+1}}}\stackrel{{}^{\sc\rm\eqref{msms-a-a-a}\,(i),\,\eqref{MSmde33333}}}{<}\d h_{p_1,p_2}\mbox{ \ for \ }i=1,2,
\nonumber\\
&\!\!\!\!\!\!\!\!\!\!\!\!\!\!\!\!\!\!\!\!\!\!\!\!\!\!\!&
{\rm(iii)\ }h_{p_1,p_2}\stackrel{{}^{\sc\rm\eqref{dp0p1}}}{=}|x_1|+|x_2|+|y_1|+|y_2|\stackrel{{}^{\sc\rm\eqref{Ckemmmdmdmd}\,(iii),\,(v),\,\eqref{msms-a-a-a}\,(ii)}}{<}
2\kk+2\d h_{p_1,p_2},
\nonumber\\
&\!\!\!\!\!\!\!\!\!\!\!\!\!\!\!\!\!\!\!\!\!\!\!\!\!\!\!&
{\rm(iv)\ }\frac{\kk}{2}\stackrel{{}^{\sc\rm\eqref{msms-a-a-a}\,(i)}}{<}h_{p_1,p_2}\stackrel{{}^{\sc\rm\eqref{msms-a-a-a}\,(iii)}}{<}3\kk,
\nonumber\\
&\!\!\!\!\!\!\!\!\!\!\!\!\!\!\!\!\!\!\!\!\!\!\!\!\!\!\!&
{\rm(v)\ }|y_2|\ge|x_2+y_2|-|x_2|\stackrel{{}^{\sc\rm\eqref{msms-a-a-a}\,(i),\,}}{\ge}
\kk-\l_1\kk+O(\ep)^1\ge\frac{\kk}{2}+O(\ep)^1\stackrel{{}^{\sc\rm\eqref{msms-a-a-a}\,(iv)}}{>}h_{_{\sc p_1,p_2}}^{^{\sc\frac m{m+1}}},
\end{eqnarray}
where the last inequality of \eqref{msms-a-a-a}\,(ii) is obtained by noting from \eqref{MSmde33333} that we may assume $0<\kk^{-\frac1{m+1}}\ll\d=\ell^{-1}$.
In particular, \eqref{msms-a-a-a}\,(v) is a contradiction with \eqref{mqp1234-2}. This proves \eqref{Ckemmmdmdmd}\,(vi).
}%
\NOUSE{Now \eqref{Ckemmmdmdmd} shows that \eqref{TheFasss} holds, which contradicts the assumption that Case 5 does not occur (at this point it may be worth mentioning that when we say Case 5 does not occur, it means that if it occurs we have already proved that Proposition \ref{real00-inj} holds).
This proves that \eqref{T2-isO(ee)} holds for $\tilde\xX_1$.
Now if $\zZ\ge1$, then
 \eqref{LetNSoOP----1**+1}\,(vi),\,(vii) imply that $\xX_1<\zZ,\,\xX_2<\zZ$ and we obtain a
contradiction as in \eqref{Conndndm}.
Hence $\zZ<1$.
Then
$\zZ$ is a $1+O(\ep)^1$ element as in \eqref{Ckemmmdmdmd}\,(i), and we have
\eqref{Ckemmmdmdmd}\,(iii) by \eqref{LetNSoOP----1**+1}\,(vii).
By \eqref{LetNSoOP----1**+1}\,(iv),
$\xX_2<\zZ^{1+\frac{3}{20}+O(\dD)^{20}}<\zZ^{1+\l_2}$ with
$\l_2=\frac1{10}$. The same arguments after
\eqref{Ckemmmdmdmd} show that
we still have
$\xX_2\ge\l:=\frac12$,
i.e.,
 \eqref{TheFasss} holds, which again contradicts the assumption that Case 5 does not occur.
}\NOUSE{
From this, \eqref{1bab-that}\,(ii) and \eqref{LetNSoOP----2}\,(iii), we obtain that $\xX_2\ge1+O(\d)^1$, but by
\eqref{LetNSoOP----2}\,(i), we have $\xX_2\le1+O(\d)^1$. This with \eqref{equa-Case6-lemm}\,(iv) proves
%
%
%
\begin{eqnarray}
\label{euqansnde}
&\!\!\!\!\!\!\!\!\!\!\!\!\!\!\!\!\!\!\!\!\!\!\!\!\!\!\!\!\!\!&
a=1+O(\d)^1\mbox{  for  }a\in T_2
:=
\Big\{
\xX_1=|X_1|,\,
\xX_2=|X_2|,\,\zZ=|Z|,\,\tilde\xX_1=\a_0\xX_1
\Big\}.\!\!\!\!\!\!\!\!
\end{eqnarray}
\NOUSE{
Equ.~\eqref{euqansnde}\,(i) holds for $\tilde\xX_1$ by \eqref{ImMpP}\,(3), thus also holds for $\cC_2,\bB_2$ by
\eqref{LetNSoOP----2}\,(iii) (using the formula $|a+b|\le|a|+|b|$ for $a,b\in\C$). Then it also holds for $\xX_2$ by
\eqref{ImMpP}\,(2), thus also for $\zZ$ by \eqref{equa-Case6-lemm}\,(iv). This proves
\eqref{euqansnde}\,(i).
}\NOUSE{
By \eqref{ImMpP}\,(2),\,(3), we have
[noting from \eqref{LetNSoOP----2}\,(iii) that when \eqref{mememrmmr}\,(i) holds, we have $|C_2|\le1+O(\d)^1\ssc\,$],
\equa{mememrmmr}{\dis {\rm(i)\ }\tilde\xX_1\stackrel{{}^{\sc\rm\eqref{ImMpP}\,(3)}}{=}1+O(\d)^1,\ \ \ \
{\rm(ii)\ }\xX_2\stackrel{{}^{\sc\rm\eqref{ImMpP}\,(2),\,\eqref{LetNSoOP----2}\,(iii)}}{\le}1+O(\d)^1.}
We have
\begin{eqnarray}
\label{AAAaaws}
&\!\!\!\!\!\!\!\!\!\!\!\!\!\!\!\!\!\!\!\!\!\!\!\!&
1+\d_0+O(\d)^1\stackrel{{}^{\sc\rm\eqref{LetNSoOP----2}\,(i)}}{\le}
\d_0\tilde\xX_1^{\ell_0+1}
+|A_1|^{\d_0^4}\xX_2^{\d_0(1-\d_0^3)}\tilde\xX_1^{-1}+O(\d)^1
\nonumber\\
&\!\!\!\!\!\!\!\!\!\!\!\!\!\!\!\!\!\!\!\!\!\!\!\!&
\ \ \ \ \ \ \ \ \ \ \ \ \ \ \ \ \ \ \,\ \stackrel{{}^{\sc\rm\eqref{mememrmmr}\,(i)}}{=}
\ \ \ \d_0
+\xX_2^{\d_0(1-\d_0^3)}+O(\d)^1.
\end{eqnarray}
This with \eqref{mememrmmr}\,(ii) shows that $\xX_2=1+O(\d)^1$. Thus we have \eqref{euqansnde} by \eqref{equa-Case6-lemm}\,(iv),\,\eqref{tX11111}\,(ii),\,(iii).
By \eqref{LetNSoOP----2}\,(iii),\,\eqref{mememrmmr}\,(i), we also obtain
\equa{BandC}{\dis{\rm(i)\ }\cC_2:=|C_2|\le1+O(\d)^1,\ \ \ \
{\rm(ii)\ }\bB_2:=|B_2|\le1+O(\d)^1.}
}\NOUSE{
Note that $\gamma_1$ is a strictly increasing (respectively, decreasing)
function on $\tilde\xX_1$ when $a\le\tilde\xX_1<1$ (respectively $0<\tilde\xX_1<a$) with $a=
(1+\d_0)^{-\d_0}$. Since $\gamma_1|_{\tilde \xX_1=1-\d}=1+\d_0-\d_0\d+O(\d)^2$, in order for \eqref{AAAaaws}\,(i) to
hold, we must have $\tilde\xX_1<a$. Since $\gamma_1|_{\tilde\xX_1=1-2\d_0^2+5\d_0^3}=1+\d_0-\frac{\d_0^4}{3}<1+\d_0,$ we obtain
\equa{mememe}{\tilde\xX_1<1-2\d_0^2+5\d_0^3.}
which with \eqref{mememrmmr} implies that $\tilde \xX_1,\xX_2\ge1+O(\d)^1$. Thus $\tilde\xX_1,\xX_2$ [and also $\zZ$ by
\eqref{equa-Case6-lemm}\,(iv)$\ssc\,$] are $1+O(\d)^1$ elements. Then \eqref{LetNSoOP----2}\,(i),\,(ii) show
that $\hH,\hH_1$ are $1+O(\d)^1$ elements.
This proves \eqref{euqansnde}.
}
Note that we have in fact frequently used \eqref{equa-Case6-lemm}\,(iv) to
obtain \eqref{euqansnde}
. To obtain more precisely results, we cannot use
\eqref{equa-Case6-lemm}\,(iv), instead we need the precise definitions in \eqref{LetNSoOP----1}.
By \eqref{LetNSoOP----1},\,\eqref{LetNSoOP----2}, we obtain
\begin{eqnarray}
\label{1dew2}
&\!\!\!\!\!\!\!\!\!\!\!\!\!\!\!\!\!\!\!\!\!\!\!\!\!\!&
{\rm(i)\ }
1=\aA_2\stackrel{{}^{\sc\rm\eqref{LetNSoOP----1}\,(ii)}}{\le}\tilde\xX_1^{-1}
\frac1{\Big(\frac{\ell_0}{2}-\Big(\frac{\ell_0}{2}-1\Big)\xX_2^{-1}\Big)},\ \ \
\implies\ \ \tilde\xX_1\le1,
\nonumber\\ &\!\!\!\!\!\!\!\!\!\!\!\!\!\!\!\!\!\!\!\!\!\!\!\!\!\!&
{\rm(ii)\ }1<\aA_3^{\d_0}\stackrel{{}^{\sc\rm\eqref{LetNSoOP----1}\,(iii)}}{=}\tilde\xX_1^2\zZ^{1+\d_0}\xX_2
\stackrel{{}^{\rm\sc\eqref{1dew2}\,(i)}}{\le}\zZ^{1+\d_0}\xX_2^{-1},\ \ \ \implies\ \ \ \xX_2<\zZ^{1+\d_0},
\nonumber\\ &\!\!\!\!\!\!\!\!\!\!\!\!\!\!\!\!\!\!\!\!\!\!\!\!\!\!&
{\rm(iii)\ }
1\ \ \ \ =\ \ \ \ |A_1|\stackrel{{}^{\sc\rm\eqref{LetNSoOP----2}\,(iii)}}{\le}
\gamma_2:=\tilde\xX_1^2\Big(1-\d_0^4+\d_0^4\tilde\xX_1^{-1}\xX_2^{\ell_0^3(1+\d_0)}\zZ^{-\ell_0^3(1+2\d_0)}\Big)
\nonumber\\ &\!\!\!\!\!\!\!\!\!\!\!\!\!\!\!\!\!\!\!\!\!\!\!\!\!\!&\ \ \ \ \ \ \ \ \
\stackrel{{}^{\sc\rm\eqref{1dew2}\,(i)}}{\le}
1-\d_0^4+\d_0^4\xX_2^{\ell_0^3(1+\d_0)}\zZ^{-\ell_0^3(1+2\d_0)},\ \ \ \implies\ \ \
\zZ^{1+2\d_0}\le\xX_2^{1+\d_0}
,
\end{eqnarray}
where the second inequality of \eqref{1dew2}\,(iii) is obtained by noting that $\gamma_2$ is (locally) a strictly increasing
function on $\tilde\xX_1$.
Note that $\gamma_1$ is a function of $\tilde\xX_1,\xX_2$ such that
\equa{msmsmsmsm222}{\dis
{\rm(i)\ }\frac{\ptl\gamma_1}{\ptl\tilde\xX_1}\Big|_{(\tilde\xX_1,\xX_2)=(1,1)}
=2,\ \ \ \ \ {\rm(ii)\ }\frac{\ptl\gamma_1}{\ptl\xX_2}\Big|_{(\tilde\xX_1,\xX_2)=(1,1)}=-3.
}
We claim
\equa{Cmm2222}{\dis
{\rm(i)\ }\xX_2\le\tilde\xX_1^{\frac23+\d_1}\mbox{ if }\tilde\xX_1\ge1,
\ \ \ \ \ \
{\rm(ii)\ }\xX_2^3\le\tilde\xX_1^{\frac23-\d_1}\mbox{ if }\tilde\xX_1<1.
}
To prove the claim, we start the case with $\tilde\xX_1\ge1$ and assume
$\tilde\xX_1^{\frac23+\d_1}<\xX_2$.
Then \eqref{1dew2} gives
\equa{1dew2+1}{\dis
1<\gamma_3:=\gamma_1|_{\tilde\xX_2=\tilde\xX_1^{\frac23+\d_1}}=
\frac{\tilde\xX_1^{\frac35-2\d_1}}{2-\tilde\xX_1^{-\frac53-\d_1}}.
}Note that
$\gamma_3$ is a strictly decreasing function of $\tilde\xX_1$ when \eqref{euqansnde} holds
(as $\frac{d\gamma_2}{d \tilde\xX_1}|_{x=1}=-3\d_1<0$ and $0<\d\ll\d_1$).
Thus \eqref{1dew2+1} shows that $\tilde\xX_1<1$, which is a contradiction
 with the case we started. This proves \eqref{Cmm2222}\,(i).
Similarly, we have \eqref{Cmm2222}\,(ii).
Thus in general,
\equa{InGenem}{\dis
\xX_2\le\tilde\xX_1^{\frac23+O(\d_1)^1}.
}
Similarly, using the fact in \eqref{1dew2}\,(iii) that
\equa{msmsmsmsm222+}{\dis
{\rm(i)\ }\frac{\ptl\gamma_2}{\ptl\tilde\xX_1}\Big|_{(\tilde\xX_1,\xX_2)=(1,1)}
=-\frac34,\ \ \ \ \ {\rm(ii)\ }\frac{\ptl\gamma_1}{\ptl\xX_2}\Big|_{(\tilde\xX_1,\xX_2)=(1,1)}=-\frac32,
}
we obtain
\equa{InGenem+1}{\dis
\xX_2\le\tilde\xX_1^{\frac23+O(\d_1)^1}.
}
Next, by \eqref{LetNSoOP----1}\,(iii), we have [using \eqref{1dew2}$\ssc\,$],
\equa{ssssss0000}{\dis
1<|A_3|
\le\gamma_3:=(1 - \d_0^2 +
   \d_0^2 \tilde \xX_1\gamma_1^{2\ell_0^2(2-\d_0^2)})\tilde \zZ^{2 + \d_0^3}\tilde \xX_1^{-2},
}
Using \eqref{msmsmsmsm222}, at the point $(\tilde\xX_1,\tilde\xX_2,\tilde\zZ)=(1,1,1)$, we have
\equa{gamma-3prime}{\dis
{\rm(i)\ }\frac{\ptl\gamma_3}{\ptl\tilde\xX_1}=\d_0^2
\Big(1+2\ell_0^2(2-\d_0^2)\frac{\ptl\gamma_1}{\ptl\tilde\xX_1}\Big)-2=
2 + 3 \d_0^2 - 2 \d_0^4,\ \ \
{\rm(ii)\ }\frac{\ptl\gamma_3}{\ptl\tilde\xX_2}=2(2-\d_0^2)\frac{\ptl\gamma_1}{\ptl\tilde\xX_2}
=-2(2-\d_0^2)(1-\d_0^3),
}
Now we consider the case with $\bB_2<1$ and assume
$\tilde\xX_1^{12\ell_0^3-1}\zZ^{\ell_0^3}<\bB_2^{1+\d_0^3}$. Then
 \eqref{1dew2} gives
\equa{1dew2+1+1}{\dis
1+\d_0^3<
\d_0^3\bB_2+\bB_2^{\d_0^3},
}
which immediately implies that $1<\bB_2$, a contradiction with the assumption
 $\bB_2<1$.
This proves \eqref{Cmm2222}\,(ii).
Thus in general we have
(noting the important fact that $\a_1>1$, from  which we obtain the first inequality below),
\equa{In-ggg}{ \dis\bB_1\stackrel{{}^{\sc\rm\eqref{LetNSoOP----1}\,(iii),\,\eqref{a0BiggerThen}\,(iv)}}{<}
\bB_2\le(\tilde\xX_1^{12\ell_0^3-1}\zZ^{\ell_0^3})^{1+O(\d_0)^3}.
}
We also have
\begin{eqnarray}
&\!\!\!\!\!\!\!\!\!\!\!\!\!\!\!\!\!\!\!\!\!\!\!\!\!\!\!&
\label{AlsoWj}
1+\d_0^3\stackrel{{}^{\sc\rm\eqref{LetNSoOP----1}\,(i)}}{\le}
\d_0^3\tilde\xX_1^{3\ell_0^3}+|A_1|^{\d_0^3}\tilde\xX_1^{-15}\zZ^{-1}\Big(\frac{1+\d_0^3+\d_0^3
\tilde\xX_1^{3\ell_0^3}}{1+2\d_0^3}\Big)^{\d_0^3}
\nonumber\\
&\!\!\!\!\!\!\!\!\!\!\!\!\!\!\!\!\!\!\!\!\!\!\!\!\!\!\!&
\phantom{1+\d_0^3\ \ \ \ }
=\ \ \ \gamma_2:=\d_0^3\tilde\xX_1^{3\ell_0^3}+\tilde\xX_1^{-15}\zZ^{-1}\Big(\frac{1+\d_0^3+\d_0^3
\tilde\xX_1^{3\ell_0^3}}{1+2\d_0^3}\Big)^{\d_0^3}.
\end{eqnarray}
\begin{clai}\label{clai-1111}
We have
\equa{Now=1+O}{\dis
{\rm(i)\ }\bB_2>1-\ell_1\d,\ \ \ \ \ {\rm(ii)\ }
\tilde \xX_1>1-\d,\ \ \ \ \ {\rm(iii)\ }a=1+O(\d)^1\mbox{ \ for \ }a\in T_2.
}
\end{clai}
To prove the claim, assume conversely \eqref{Now=1+O}\,(i) does not hold, i.e.,
\equa{converse}{\dis
\bB_2\le1-\ell_1\d.}
Then
\begin{eqnarray}
\!\!\!\!\!\!\!\!\!\!\!\!\!\!\!\!\!\!\!\!\!\!&&
\label{664646}
{\rm(i)\ }\xX_2\stackrel{{}^{\sc\rm\eqref{LetNSoOP----1}\,(iv)}}{=}|A_3|^{-1}\bB_2
\stackrel{{}^{\sc\rm\eqref{1bab-that}\,(ii)}}{<}\bB_2
\stackrel{{}^{\sc\rm\eqref{converse}}}{\le}1-\ell_1\d,
\nonumber\\
\!\!\!\!\!\!\!\!\!\!\!\!\!\!\!\!\!\!\!\!\!\!&&
{\rm(ii)\ }\xX_1:=|X_1|
\stackrel{{}^{\sc\rm\eqref{tX11111}}}{=}
\a_0^{-\frac12}\tilde\xX_1^{\frac{\d_1}{2}}\xX_2^{\frac12}
\stackrel{{}^{\sc\rm\eqref{a0BiggerThen}\,(ii)}}{<}
\tilde\xX_1^{\frac{\d_1}{2}}\xX_2^{\frac12}
\stackrel{{}^{\sc\rm\eqref{ImMpP}\,(3)}}{<}(1+\d)^{\frac{\d_1}{2}}\xX_2^{\frac12}
\stackrel{{}^{\sc\rm\eqref{664646}\,(i)}}{<}1.
\end{eqnarray}
First assume $\zZ\ge1$. Then $\xX_1,\xX_2<\zZ$.
Noting from \eqref{TaKa},\,\eqref{SimMMSMS} that $|x_1|=\xX_1\kk,\,|x_2|=\xX_2\kk,$ $|x_2+y_2|=\zZ\g_{\kk,\kk}$,
we have
[where the third inequality holds trivially if $\zZ=1$ and it follows from \eqref{AssG} in case $\zZ>1$],
\equa{Conndndm}{\dis
\g_{|x_1|,|x_2|}
=\g_{\xX_1\kk,\xX_2\kk}\stackrel{{}^{\sc\rm\eqref{wePPPP1+}}}{<}
\g_{\zZ\kk,\xX_2\kk}
\stackrel{{}^{\sc\rm\eqref{wePPPP1}}}{\le}\g_{\zZ\kk,\zZ\kk}
\stackrel{{}^{\sc\rm\eqref{AssG}}}{\le}\zZ\g_{\kk,\kk}
=|x_2+y_2|\stackrel{{}^{\sc\rm\eqref{Ak=1}}}{\le}\g_{|x_1|,|x_2|},
}
a contradiction.
Next assume $\zZ<1$. In this case we want to prove that $\xX_2\le\zZ^{\frac{\ell_0^3}{2}}$. Assume conversely
\equa{Convvvnn}{\dis
\zZ^{\frac{\ell_0^3}{2}}<\xX_2.
}
Then
\equa{xxx2222}{\dis
\xX_2\stackrel{{}^{\sc\rm\eqref{664646}\,(i)}}{<}\bB_2
\stackrel{{}^{\sc\rm\eqref{In-ggg},\,\eqref{Convvvnn},\,\eqref{ImMpP}\,(3)}}{<}
\Big((1+\d)^{12\ell_0^3-1}\xX^2\Big)^{1+O(\d_0)^3},
}
which implies that $\xX_2>(1+\d)^{-13\ell_0^2}$, a contradiction with
\eqref{664646}\,(i).
This proves that $\xX_2\le\zZ^{\frac{\ell_0^3}{2}}$.
Now we want to verify whether or not \eqref{TheFasss} holds. To do this, we use  notations in
\eqref{TheFasss}, and obtain
\equa{1-k2k2k3mmm}{\dis{\rm(i)\ }
\hat k:=|Z|=\zZ<1,\ \ \
{\rm(ii)\ }k_2:=|X_2|=\xX_2
\le
\hat k^{1+\l_2},
\mbox{ where $\l_2=\frac{\ell_0^3}2-1$.}
}
We also have $k_1:=|X_1|=\xX_1<1$ by \eqref{664646}\,(ii).
Further, $k_2\ge\l$ with $\l=\d_2^{\d_1}$ by \eqref{ImMpP}\,(5).
Note that $\l,\l_2$ are independent of $\kk$ by
Remark \ref{rema3.1}\,(i), by
 \eqref{TaKa},\,\eqref{SimMMSMS}, we  see that
 \eqref{TheFasss} holds, which contradicts the assumption that Case 5 does not occur.
This proves \eqref{Now=1+O}\,(i), and hence $\bB_2$ is a $1+O(\d)^1$ element [and so is $\bB_1$ by \eqref{LetNSoOP----2}\,(iii)$\ssc\,$].
Next assume $\tilde\xX_1\le1-\d$.
 Then $\bB_2<1$ by \eqref{LetNSoOP----1}\,(iii) (using formula
 $|a+b|\le|a|+|b|$ for $a,b\in\C$ and noting that $0<\ep\ll\d$).
 Then $\xX_2<1$ by the first inequality of \eqref{664646}\,(i), and also $\xX_1<1$ by the first inequality
 of \eqref{664646}\,(ii).
 Now if $\zZ\ge1$, we obtain a contradiction as in \eqref{Conndndm}; if $\zZ<1$ then $\xX_2\le\zZ^{\frac{\ell_0^3}{2}}$
 by the first inequality of \eqref{664646}\,(i) and the last inequality of \eqref{In-ggg}, and we obtain
 a contradiciton as in \eqref{1-k2k2k3mmm}.
This proves \eqref{Now=1+O}\,(i) and so $\xX_1$ is also a $1+O(\d)^1$ element.
Then \eqref{1dew2} shows that $\zZ\ge1+O(\d)^1$, which implies that it is a $1+O(\d)^1$ element by
\eqref{euqansnde}\,(i). So is $\xX_2$ by \eqref{equa-Case6-lemm}\,(iv). This proves Claim \ref{clai-1111}.
\vskip5pt
Now
by \eqref{Now=1+O}\,(iii), as in the proof of \eqref{Cmm2222}\,(i), we can now refine
the result in \eqref{In-ggg}
by using \eqref{1dew2},\,\eqref{msmsmsmsm222} to obtain the following more precisely result,
\equa{1+In-ggg}{ \dis\bB_1<
\bB_2\le(\tilde\xX_1^{12\ell_0^3-1}\zZ^{\ell_0^3})^{\frac1{1-\d_0^3}+O(\d_0)^9}.
}
Now by \eqref{AlsoWj}, we have
\equa{msmsm000000}{\dis
{\rm(i)\ }
\frac{\ptl\gamma_2}{\ptl\tilde\xX_1}\Big|_{(\tilde\xX_1,\zZ)=(1,1)}=-12+3\d_0^3+O(\d_0)^6,\ \ \ \ \ \
{\rm(ii)\ }\frac{\ptl\gamma_2}{\ptl\zZ}\Big|_{(\tilde\xX_1,\zZ)=(1,1)}=-1.
}
From this, as in the proof of \eqref{Cmm2222}, we obtain from
\eqref{AlsoWj} the following,
\equa{Anothermmm}{\dis
\zZ\tilde\xX_1^{12-3\d_0^3+O(\d_0)^6}\le1.
}
\NOUSE{
Then \eqref{AlsoWj} shows that $\zZ\ge1+O(\d)^1$, thus it
is a $1+O(\d)^1$ element and so is $\xX_2$ by \eqref{equa-Case6-lemm}\,(iv).
This proves \eqref{Now=1+O}.
Now by \eqref{Now=1+O}\,(ii), as in the proof of \eqref{Cmm2222}, we obtain
from \eqref{AlsoWj} the following,
\equa{ANn-0000}{\dis
\tilde\xX_1^{3+O(\d_0)^3}<\zZ.
}
}%
Further,
\equa{1==AlsoWj}{\dis\!\!\!\!\!\!\!\!\!\!
(\tilde\xX_1^{12\ell_0^3-1}\zZ^{\ell_0^3})^{\frac1{1-\d_0^3}+O(\d_0)^9}
 \stackrel{{}^{\sc\rm\eqref{1+In-ggg}}}{>}
\bB_1=\bB_1|A_1|^{\d_0^3}
\stackrel{{}^{\sc\rm\eqref{LetNSoOP----1}\,(i),\,(ii)}}{\ge}
\frac{\tilde\xX_1^{15}\zZ\Big(\frac{(1+\d_0^3)^2-\d_0^6\xX_1^{6\ell_0^3}}{1+2\d_0^3}\Big)}
{\Big(\frac{1+\d_0^3+\d_0^3\tilde\xX_1^{3\ell_0^3}}{1+2\d_0^3}\Big)^{\d_0^3}}
,
\!\!\!\!\!\!\!\!}
which gives
\equa{AlsoWj+1}{\dis
1\le\gamma_3:=\frac{(\tilde\xX_1^{12\ell_0^3-1}\zZ^{\ell_0^3})^{\frac1{1-\d_0^3}+O(\d_0)^9}\Big(\frac{1+\d_0^3+\d_0^3\tilde\xX_1^{3\ell_0^3}}{1+2\d_0^3}\Big)^{\d_0^3}}{\tilde\xX_1^{15}\zZ\Big(\frac{(1+\d_0^3)^2-\d_0^6\xX_1^{6\ell_0^3}}{1+2\d_0^3}\Big)}
.
}
We have
\equa{gamm3}{\dis
{\rm(i)\ }\frac{\ptl\gamma_3}{\ptl\tilde\xX_1}\Big|_{(\tilde\xX_1,\zZ)=(1,1)}=12\ell_0^3-4+20\d_0^3+O(\d_0)^6,\ \ \ \
{\rm(ii)\ }\frac{\ptl\gamma_3}{\ptl\zZ}\Big|_{(\tilde\xX_1,\zZ)=(1,1)}=\ell_0^3+\d_0^3+O(\d_0)^6.
}
From this, one can easily see,  as in \eqref{In-ggg},
\equa{xxxab}{\dis1<\zZ\tilde\xX_1^{12-4\d_0^3+O(\d_0)^6}.}
This with \eqref{Anothermmm} implies (i) below, then (ii) is obtained from (i) and \eqref{xxxab},
\equa{i-iiiii}{\dis
{\rm(i)\ }\tilde\xX_1<1,\ \ \ \ {\rm(ii)\ }\zZ>1.
}
We have
\begin{eqnarray}
&\!\!\!\!\!\!\!\!\!\!\!\!\!\!\!\!\!\!\!\!\!\!\!\!\!\!\!\!\!\!&
\label{mmxxX2}
\xX_2\stackrel{{}^{\sc\rm\eqref{664646}\,(i)}}{<}\bB_2
\stackrel{{}^{\sc\rm\eqref{1+In-ggg}}}{\le}
\Big(\tilde\xX_1^{12\ell_0^3-1-(12-3\d_0^3+O(\d_0)^6)\ell_0^3}(\zZ\tilde\xX_1^{12-3\d_0^3+O(\d_0)^6})^{\ell_0^3}\Big)^{\frac1{1-\d_0^3}+O(\d_0)^9}
\nonumber\\
&\!\!\!\!\!\!\!\!\!\!\!\!\!\!\!\!\!\!\!\!\!\!\!\!\!\!\!\!\!\!&
\ \ \ \ \ \stackrel{{}^{\sc\rm\eqref{Anothermmm}}}{\le}\ \ \Big(\tilde\xX_1^{2+O(\d_0)^3}\Big)^{\frac1{1-\d_0^3}+O(\d_0)^9}
\stackrel{{}^{\sc\rm\eqref{i-iiiii}\,(i)}}{<}1.
\end{eqnarray}
By \eqref{i-iiiii}\,(i),\,\eqref{mmxxX2}, as in the proof of
\eqref{664646}\,(ii), we obtain that $\xX_1<1$. This with
\eqref{i-iiiii}\,(ii),\,\eqref{mmxxX2} gives \eqref{Conndndm}, which is a contradiction.
\NOUSE{
where (ii) is obtained as follows:
\equa{Didiek}{\dis
1
\stackrel{{}^{\sc\rm\eqref{1bab-that}\,(ii),\,\eqref{In-ggg}}}{<}
\frac{|A_3|\tilde\zZ^{\d_0+O(\d_0)^2}}{\tilde\xX_1\bB_2}
\stackrel{{}^{\sc\rm\eqref{LetNSoOP----1}\,(iv)}}{=}\frac{\tilde\zZ^{\d_0+O(\d_0)^2}}{\tilde\xX_2}
.}
By notation \eqref{euqansnde}\,(ii),
we obtain from \eqref{i-iiiii},\,\eqref{tX11111}\,(ii),\,(iii) the following,
\equa{1-2-3---}{\dis
{\rm(i)\ }\xX_2=\tilde\zZ^{-\frac15}\zZ^{\frac65}<\zZ^{\frac65},\ \ \ \ \ {\rm(ii)\ }\zZ^2
=\tilde\xX_2\xX_2<\xX_2^{1+O(\d_0)^1},\ \ \ \ \ {\rm(iii)\ }\xX_2,\zZ<1,
}
where (iii) is obtained from (i) and (ii).
We want to verify whether or not \eqref{TheFasss} holds. To do this, we use  notations in
\eqref{TheFasss}, and obtain
\equa{k2k2k3mmm}{\dis{\rm(i)\ }\hat k:=|Z|=\zZ\stackrel{{}^{\sc\rm\eqref{1-2-3---}\,(iii)}}{<}1,\ \ \
{\rm(ii)\ }k_2:=|X_2|=\xX_2\stackrel{{}^{\sc\rm\eqref{1-2-3---}\,(i)}}{<}\hat k^{1+\l_2}\mbox{ \ with \ }\l_2=\frac15.
}
We also have (noting the important fact that $\a_0>1$, from  which we obtain the first inequality below)
\equa{MMMMMSmsms}{\dis
k_1:=|X_1|\stackrel{{}^{\sc\rm\eqref{tX11111}\,(i)}}{=}
(\a_0^{-1}\xX_1^{\d_1}\xX_2^3)^{\frac15}
\stackrel{{}^{\sc\rm\eqref{a0BiggerThen}\,(ii)}}{<}( \xX_1^{\d_1}\xX_2^3)^{\frac15}
\stackrel{{}^{\sc\rm\eqref{xxxab}\,}}{<}( \tilde\zZ^{(\d_0+O(\d_0)^2)\d_1}\xX_2^3)^{\frac15}
<1,
}
where the last inequality is obtained by \eqref{1-2-3---}\,(iii) and noting that $\tilde\zZ=\zZ^5\xX_2^{-4}$.
We also have $k_2\ge\l$ with $\l=\d_2^{\d_1}$ by \eqref{ImMpP}\,(5). Note that $\l,\l_2$ are independent of $\kk$ by
Remark \ref{rema3.1}\,(i), by
 \eqref{TaKa},\,\eqref{SimMMSMS}, we  see that
 \eqref{TheFasss} holds, which contradicts the assumption that Case 5 does not occur.
}%
}%
}This proves that the first strict inequality of \eqref{LetNSoOP}\,(i) holds for $(p_1,p_2)$.

\NOUSE{
Similarly from \eqref{LetNSoOP----1}\,(i), we obtain
\equa{1dew2+2}{\dis
2
\le \xX_1^9+\tilde\xX_1\xX_2^3\zZ^{-1}
.
}
We claim \equa{clq111}{\dis
\xX_2^3\zZ^{-1}\ge1,}
 otherwise \eqref{1dew2+2} shows that $\xX_1>1$ but \eqref{1dew2} [with condition \eqref{euqansnde}$\ssc\,$] shows that $\xX_1<1$, a contradiction.
Equ.~\eqref{clq111} with \eqref{1bab-that}\,(ii),\,\eqref{LetNSoOP----1} gives that $\xX_2^2<\zZ\le\xX_2^3$. Thus
\equa{mdnrnrn}{\dis
1<\xX_2<\zZ.
}
Further \eqref{1dew2} with \eqref{euqansnde} gives that $\tilde \xX_1\le(\xX_2^3\zZ^{-1})^{\frac{\d_0}{11}+\d_0^2}$.
This with \eqref{tX11111} proves
\equa{x11,,}{\dis
\xX_1:=|X_1|=\Big(\a_0^{-1}\tilde \xX_1^{\frac1{20}}\xX_2^3\Big)^{\frac15}
\stackrel{{}^{\sc\rm\eqref{a0BiggerThen}}}{<}\Big(\tilde \xX_1^{\frac1{20}}\xX_2^3\Big)^{\frac15}<
\zZ.
}
Using the notations in \eqref{SimMMSMS},\,\eqref{euqansnde}, we have $|x_1|=\xX_1\kk,\,|x_2|=\xX_2\kk,\,|x_2+y_2|=
\zZ\g_{\kk,\kk}$. Then we obtain the following contradiction,
\equa{Commene}{\dis
\g_{|x_1|,|x_2|}\stackrel{{}^{\sc\rm\eqref{Ak=1}}}{\ge}|x_2+y_2|=\zZ\g_{\kk,\kk}
\stackrel{{}^{\sc\rm\eqref{AssG}}}{>}\g_{\zZ\kk,\zZ\kk}
\stackrel{{}^{\sc\rm\eqref{wePPPP1+}}}{>}\g_{\xX_1\kk,\zZ\kk}
\stackrel{{}^{\sc\rm\eqref{wePPPP1}}}{\ge}\g_{\xX_1\kk,\xX_2\kk}
=\g_{|x_1|,|x_2|}.
}
Next assume in {\rm\eqref{C+LetNSoOP}\,(i),}
 equality occurs in the last inequality, i.e., $\aA_{\rOnE}=\nn$.
Then $\tilde\xX_1<1$ by
Denote $\aA_i=|A_i|$ for $i=1,2,3$ and using notation in \eqref{euqansnde}\,(ii).
Set $\b_1=|2-\widetilde X_1^{9}|$, $\b_2=|2-\widetilde X_1^{11}A_1|$ [which satisfies \eqref{ImMpP}\,(4)$\ssc\,$].
Then \eqref{LetNSoOP----2} gives
\equa{Fibebe}{\dis
\aA_1^{-22}\aA_3^{-44}\b_1^{22}\stackrel{{}^{\sc\rm\eqref{LetNSoOP----2}\,(i)}}{=}
\tilde\xX_1^{22}\stackrel{{}^{\sc\rm\eqref{LetNSoOP----2}\,(ii)}}{=}(\aA_2\aA_3^{-2})^{\d_0}\aA_1^{-2}\b_1^{-1}.
}
%
%
%
%
%
}\NOUSE{Using \eqref{LetNSoOP----1}\,(i) and \eqref{1bab-that}
,
we obtain,
\begin{eqnarray}
\label{B1-0-0}
&\!\!\!\!\!\!\!\!\!\!\!\!\!\!\!\!\!\!\!\!\!\!\!\!\!\!\!\!\!&
1\stackrel{{}^{\sc\rm\eqref{1bab-that}}}{=}|A_1|\stackrel{{}^{\sc\rm\eqref{LetNSoOP----1}\,(i),\,\eqref{euqansnde}\,(ii)}}{\le}
\gamma_1:=\zZ^{\ell_0(1+\d_0)}\xX_1'^{-\ell_0}\Big(\d_0+(1-\d_0)\xX_1'^{\ell_0(1 + \d_0)}\zZ^{-\ell_0^4(1 + \d_0)}\Big)
\nonumber\\
&\!\!\!\!\!\!\!\!\!\!\!\!\!\!\!\!\!\!\!\!\!\!\!\!\!\!\!\!\!&
\phantom{1}\ \
<\xX_1'^{-\ell_0}\Big(\d_0+(1-\d_0)\xX_1'^{\ell_0(1 + \d_0)}\Big),
\end{eqnarray}
where the  strict inequality follows from the facts that $\zZ<1$,\,$\zZ=1+O(\d)^1$ and
$\gamma_1$ is a strictly increasing function on $\zZ$ when $\xX_2=1+O(\d)^1$ is fixed as $\frac{\ptl\gamma_1}{\ptl\zZ}|_{(\xX_1',\zZ)=(1,1)}=\ell_0^2>0$. Now the right-hand side of \eqref{B1-0-0} is a strictly decreasing function on $\xX_1'$, which implies (i) below; then we obtain
\equa{MEmememn}{\dis
{\rm(i)\ }\xX_1'<1,\ \ \ \ \ {\rm(ii)\ }|X_1|=\xX_1\stackrel{{}^{\sc\rm\eqref{a0BiggerThen}}}{<}
\a_0\xX_1\stackrel{{}^{\sc\rm\eqref{LetNSoOP----1}\,(ii),\,
\eqref{MEmememn}\,(i)}}{<}1.
}
Then
\begin{eqnarray}
\label{B1-0-0+1}
&\!\!\!\!\!\!\!\!\!\!\!\!\!\!\!\!\!\!\!\!\!\!\!\!\!\!\!\!\!&
1=|A_2|\stackrel{{}^{\sc\rm\eqref{LetNSoOP----1}\,(iii),\,\eqref{euqansnde}\,(ii)}}{\le}\gamma_2:=(\a_0\xX_1)^{-1}
\zZ^2\yY^{-1}\Big(\frac14+
\frac{3(\a_0\xX_1)^{2}\yY^{-2}}{4}\Big)<
\yY^{-1}\Big(\frac14+
\frac{3\yY^{2}}{4}\Big),
\end{eqnarray}
where the strict inequality follows from \eqref{msmA2p}\,(i),\,\eqref{MEmememn}\,(ii) and the fact that $\gamma_2$ is a strictly increasing on $\a_0\xX_1$ when \eqref{euqansnde}\,(i) holds. The right-hand side of \eqref{B1-0-0+1} is a strictly increasing function on $\yY$ when $\yY=1+O(\d)^1$, which implies that $\yY>1$, i.e.,
\equa{msmsnyYYYY}{\dis |X_2|=\xX_2<\zZ^2=|Z|^2.}
}\NOUSE{
\begin{eqnarray}
&\!\!\!\!\!\!\!\!\!\!\!\!\!\!\!\!\!\!\!\!\!\!\!\!\!\!\!\!\!&
\aA_3^{-\d_0^2}\Big((1{\ssc\!}-{\ssc\!}\d_0^2)\yY{\ssc\!}+{\ssc\!}(1{\ssc\!}+{\ssc\!}\d_0^2)\yY^{-1}\Big)
\stackrel{{}^{\sc\rm\eqref{1bab-that}\,(ii)}}{<}
(1{\ssc\!}-{\ssc\!}\d_0^2)\yY{\ssc\!}+{\ssc\!}(1{\ssc\!}+{\ssc\!}\d_0^2)\yY^{-1}
,\!\!\!\!\!\!\!\!\!\!\!\!\!\!
\nonumber\\\nonumber
&\!\!\!\!\!\!\!\!\!\!\!\!\!\!\!\!\!\!\!\!\!\!\!\!&
{\rm(ii)\ }2\  \ \,\stackrel{{}^{\sc\rm\eqref{1bab-that}}}{=}2|A_2|
\stackrel{{}^{\sc\rm\eqref{LetNSoOP----1}\,(ii),\,\eqref{euqansnde}\,(ii)}}{\le}\gamma_0:=
\aA_3^{-2+2\d_0^2}\Big(
(1-\d_0^4)\xX_1+(1+\d_0^4)\aA_3^2\xX_1^{-1}\Big)
\\
&\!\!\!\!\!\!\!\!\!\!\!\!\!\!\!\!\!\!\!\!\!\!\!\!&
\phantom{{\rm(ii)\ }2}
\!\!\!\!\!\!\stackrel{{}^{\sc\rm\eqref{1bab-that},\,\eqref{euqansnde}\,(i)}}{<}
(1-\d_0^4)\xX_1+(1+\d_0^4)\xX_1^{-1}
,\!\!\!\!\!\!\!\!
\end{eqnarray}
where the strict inequality in (i) is obvious as $\aA_3>1$, while that in (ii) is
obtained by the fact that $\gamma_0$ is
a strictly decreasing function of $\aA_3$ when $\xX_1$ is fixed and when \eqref{euqansnde}\,(i) holds.
Now the right-hand sides of \eqref{B1-0-0}\,(i),\,(ii) are decreasing functions of $\yY$ and $\xX_1$ respectively
when \eqref{euqansnde}\,(i) holds (noting that $0<\d\ll\d_0$),
which proves that $\yY<1,\,\xX_1<1$. Thus, by notations in \eqref{euqansnde}\,(ii) and \eqref{LetNSoOP----1}\,(iii), we have \equa{SMEMEMEm}{\dis
|X_1|<1,\ \ \ \ \ |X_2|<|Z|^{1+\d_0^2}\stackrel{{}^{\sc\rm\eqref{LetNSoOP----1}\,(iii)}}{=}\aA_3^{-1-\d_0^2}<1.}
}\NOUSE{We claim
, while the strict inequality in (i) is obtained by noting that $\gamma_0$ is
a strictly decreasing function of $\aA$
(since $\frac{\ptl\gamma_0}{\ptl\aA}|_{(\aA,\yY)=(1,1)}=-10<0$,
$\frac{\ptl\tilde\gamma_0}{\ptl\aA}|_{(\aA,\xX_1,\xX_2)=(1,1,1)}=-15<0$)
when
other variables are fixed and when all variables satisfy
\eqref{euqansnde} (thus $\gamma_0<\gamma_0|_{\aA=1}=\gamma_1$).
Note that $\gamma_1(\xX),\,\gamma_2(\yY)$ are increasing functions when $\xX,\yY\in T_2$. Thus \eqref{B1-0-0} shows
that $\xX,\yY>1$.
This with \eqref{x0s,dm} gives,
\equa{This-is-ok}{\dis
|X_2|<|Z|^{\frac{7}{5}}<|X_2|^{\frac{7}{10}}.
}
By the next lemma, we may assume this does not occur.
}\NOUSE
{
First we use \eqref{B1-0-0}\,(ii) to prove 
\begin{eqnarray}
\label{Bsdbebe00}
\!\!\!\!\!\!\!\!\!\!\!\!\!\!\!\!\!\!\!\!&&
{\rm(i)\ }\bA_1<\bA_2^{\frac{\ell_0^7}{2}+\d_1}
\mbox{ \ if \ }\bA_2\ge1,\mbox{ \ \ or \ \ }
{\rm(ii)\ }\bA_1<\bA_2^{\frac{\ell_0^7}{2}-\d_1}\mbox{ \ if \ }\bA_2<1
.
\end{eqnarray}
We start with the case (i) that $\bA_2\ge1$. Assume conversely $\bA_1\ge\bA_2^{\frac{\ell_0^7}{2}+\d_1}$.
Noting that $\gamma_2$ 
is a strictly decreasing function of $\bA_1$
when $\bA_2$ is  fixed and when all variables satisfy \eqref{euqansnde} (since
$\frac{\ptl{\ssc}\gamma_2}{\ptl\bA_1}|_{(\bA_1,\bA_2)=(1,1)}=-2\d_0^{11}<0$ and $0<\d\ll\d_0$),
we obtain from  \eqref{B1-0-0}\,(i),
\equan{beta2}{\dis
1<\gamma_2(\bA_1,\bA_2)\le\gamma_2(\bA_2^{\frac{\ell_0^7}{2}+\d_1},\yY)=
\tilde\gamma_2(\bA_2):=\bA_2^{(\frac{\ell_0^7}{2}+\d_1)(\d_0^2-\d_0^6-2\d_0^{11})}
\Big(\frac{(1+\d_0^4)\bA_2}{2}
+\frac{1-\d_0^4}{2\bA_2^{2(\frac{\ell_0^7}{2}+\d_1)\d_0^2+1}}\Big)
.}
Noting that $\tilde\gamma_2$  is a strictly decreasing function on $\bA_2$ when $\bA_2$
satisfies \eqref{euqansnde} [as $\frac{d{\ssc}\tilde\gamma_2}{d\bA_2}|_{\bA_2=1}=-2\d_0^{11}\d_1<0$
and $0<\d\ll\d_1\ll\d_0$ by \eqref{MSmde33333}$\ssc\,$],
we obtain from the above that $\bA_2<1$, which a contradiction with the case we started. This proves
\eqref{Bsdbebe00}\,(i). Similarly, we have \eqref{Bsdbebe00}\,(ii).
%
Using the fact from \eqref{LetNSoOP----1}\,(iv),\,(v) that $\bA_1=\xX_1^{-\ell_0^2}$, $\bA_2=\xX_1\yY$, and noting that $0<\d_1\ll\d_0$, we obtain from \eqref{Bsdbebe00} the following,
\equa{giveIngeneral}{\dis
{\rm(i)\ }1<\bA_1^{-2}\bA_2^{\ell_0^7+O(\d_1)^1}=\xX_1^{\ell_0^7(1+O(\d_0)^5)}\yY^{\ell_0^7+O(\d_1)^1},\mbox{ \ \ or \ \ }
{\rm(ii)\ }1<\xX_1\yY^{1+O(\d_0)^5}.
}
If we denote $\bA_1'=\bA_1\yY^{-\frac{\d_0^4}{2}}$
then as in the proof of \eqref{Bsdbebe00}, using \eqref{B1-0-0}\,(i), we can prove that $\bA_2'^{2\ell_0^4}<\bA_1^{\d_0^2+O(\d_1)^1}$ [noting that $\frac{\partial\gamma_1}{\partial\bA_1}|_{(\bA_1,\bA_2,\yY)=(1,1,1)}
=\d_0^2$]. Thus
\equa{msmsme888}{\dis
1<\bA_2'^{-2}\bA_1^{\d_0^6+O(\d_1)^1}=\xX^{-2-\d_0^4+O(\d_1)^1}\yY^{-2+\d_0^4}.
}
This with \eqref{giveIngeneral}\,(ii) shows that $\yY>1$, i.e., $|Z|<|X_2|^{1-\d_0}$, and  $\xX<1$ by \eqref{msmsme888}.
\NOUSE{
Since $\gamma_1$ defined in
\eqref{B1-0-0}\,(i) is a strictly increasing function of $\bA_1$ when \eqref{euqansnde} holds,
we obtain \equan{memen9329}{\dis
1{\sc\!}<{\sc\!}
\gamma_1(\bA_1,\bA_2,\yY){\sc\!}<{\sc\!}
\tilde \gamma_1(\bA_2,\yY){\sc\!}:={\sc\!}\gamma_1(\bA_2^{\frac{\ell_0^7}{2}+O(\d_1)^2},\bA_2,\yY)
{\sc\!}={\sc\!}\yY\bA_2^{-2\ell_0^4}\Big(\frac{(1{\sc\!}+{\sc\!}\d_0^2)\bA_2^{\frac{\ell_0^7}{2}+O(\d_1)^2}}{2}
 {\sc\!}+{\sc\!}\frac{1{\sc\!}-{\sc\!}\d_0^2}{2 \bA_2^{\frac{\ell_0^7}{2}+O(\d_1)^2}}\Big)
.}
Note that $\tilde \gamma_1$ is locally a strictly increasing function on both variables, as in the proof of
\eqref{Bsdbebe00}, we can prove [noting that $\frac{\partial\tilde\gamma_1}{\partial\bA_2}|_{(\bA_2,\yY)
=(1,1)}=\frac{\ell_0^5}{2}-2\ell_0^4+O(\d_1)^2\ssc\,$],
\equa{e,e,e}{{\rm(i)\ }
\yY^{-1}<\bA_2^{\frac{\ell_0^5}{2}-2\ell_0^4+O(\d_1)^1},
}
(noting that $0<\d_1\ll\d_0$). Thus using the fact from \eqref{LetNSoOP----1}\,(iv) that
 $\bA_2=\xX_1\yY$, we obtain,
\equa{Mememe}{\dis
1<\bA_2\yY^{\ell_0^2(1+O(\d_0)^1)}=\xX_1\yY^{\ell_0^2(1+O(\d_0)^1)}.}
Noting  the fact from \eqref{LetNSoOP----1}\,(iii) that
$\xX_1=\bA_1^{-\d_0}\aA_3^{-\ell_0}<\bA_1^{-\d_0^2}$, we obtain
from \eqref{Bsdbebe00}\,(i),\,(ii) that
$\xX_1<\yY^{1+O(\d_1)^1}$, which with \eqref{Mememe} gives that $\yY>1$, i.e., $|Z|<|X_2|^{1-\d_0}$.
}%
This with \eqref{A3isBigger} gives, where
\equa{Masmwmw00000}{\dis
|X_2|^{1-\d_0+\d_0}<|Z|<|X_2|^{1-\d_0}.}
 Thus $|X_2|,|Z|<1$.
}{%
\NOUSE{Next we want to verify whether or not \eqref{TheFasss} holds. To do this, we use  notations in
\eqref{TheFasss}, 
\equa{denotttt}{k_1:=|X_1|\stackrel{{}^{\sc\rm\eqref{MEmememn}\,(ii)}}{
<}1,\ \ \ \hat k:=|Z|<1,\ \ \ k_2:=|X_2|\stackrel{{}^{\sc\rm\eqref{msmsnyYYYY}}}{<}\hat k^{1+\l_2}<1,\ \ \ \l_2=1.}
We have $k_2\ge\l_2:=\d_2$.
\begin{eqnarray}
&\!\!\!\!\!\!\!\!\!\!\!\!\!\!\!\!\!\!\!\!&
k_1\stackrel{{}^{\sc\rm\eqref{denotttt}}}{=}
|\widetilde X_1|\stackrel{{}^{\sc\rm\eqref{ImMpP}\,(1)}}{<}|X_1|^{\d_0^2\d}
\stackrel{{}^{\sc\rm\eqref{SMEMEMEm}}}{<}1,
\nonumber\\
&\!\!\!\!\!\!\!\!\!\!\!\!\!\!\!\!\!\!\!\!&
k_2\stackrel{{}^{\sc\rm\eqref{denotttt},\,\eqref{SimMMSMS}\,(v)}}{=}|X_2|^{\d_0^2}
\stackrel{{}^{\sc\rm\eqref{SMEMEMEm}}}{<}|Z|^{\d_0^2(1+\d_0^2)}
\stackrel{{}^{\sc\rm\eqref{SimMMSMS}\,(vi),\,\eqref{denotttt}}}{=}
\hat k^{1+\l_2},\ \ \l_2=\d_0^2,
\nonumber\\
&\!\!\!\!\!\!\!\!\!\!\!\!\!\!\!\!\!\!\!\!&
k_2=|X_2|^{\d_0^2}\stackrel{{}^{\sc\rm\eqref{XXXX2222}}}{>}\l,\ \ \ \l=\d_2^{\d_0^2}.
\end{eqnarray}
}%
\NOUSE{
Note that $\l,\l_2$ are independent of $\kk$.
Thus \eqref{TheFasss} holds, which contradicts the assumption that Case 5 does not occur.
This proves that the first strict inequality of \eqref{LetNSoOP}\,(i) holds for $(p_1,p_2)$.
in general we have $1<\xX_2^{-1}\xX_1^{1+O(\d_1)^1}$ and
$1<\xX_1^{-1}\yY^{3+O(\d_1)^1}$, i.e.,
we have by notations in \eqref{euqansnde},
\equa{i.e.,d,d}{\dis
{\rm(i)\ }1<\Big|\frac{\tildeX_1}{X_2^{1+O(\d_1)^1}}\Big|,
\ \ \ \ {\rm(ii)\ }1<\Big|\frac{X_2^{3+O(\d_1)^1}Z^{4+O(\d_1)^1}}{\tildeX_1}\Big|
.}
Multiplying \eqref{i.e.,d,d}\,(i) with \eqref{x0s,dm} and \eqref{i.e.,d,d}\,(ii) respectively, we obtain
\equa{X2-Z}{\dis
|X_2^{1+O(\d_1)^1}Z|<1<
|X_2^{2+O(\d_1)^1}Z^{4+O(\d_1)^1}|.}
This proves that
$k_2:=|X_2|<1$ and $k:=|Z|>1$,
 which
together with \eqref{ImMpP}\,(1),\,\eqref{x0s,dm}  implies that
  $k_{\ZeRo}:=|X_{\ZeRo}|<|\tildeX_{\ZeRo}|^{\d}<|Z ^{-1}|^\d<
 1$.
By 
\eqref{TaKa}, \eqref{SimMMSMS}\,(c),
we obtain,
\equa{mwemnexy-1}{\mbox{$|x_{\OnE}|=k_{\OnE}\kk<\kk<k\kk$, \ \
$|x_{\ZeRo}|=k_{\ZeRo}\kk<\kk<k\kk$, \ \ $|x_{\OnE}+y_{\OnE}|=|Z(\bar x_2+\bar y_2)|=k\g_{\kk,\kk}$.}}
Thus, we have [where the first two inequalities follow from Lemma \ref{g-1-lemm} and \eqref{wePPPP1},
while the last from definition \eqref{Ak=1}$\ssc\,$],
\equa{111--ByDsllo-0??}{\g_{k\kk,k\kk}>\g_{|x_{\ZeRo}|,k\kk}\ge\g_{|x_{\ZeRo}|,|x_{\OnE}|}\ge|x_{\OnE}+y_{\OnE}|
= k{\ssc\,}\g_{\kk,\kk},}
which is a contradiction with \eqref{AssG}.
}%
}%
}\vskip5pt{
Next assume in {\rm\eqref{C+LetNSoOP}\,(i),}
 equality occurs in the last inequality, i.e., $\aA_1=|A_{\rOnE}|=\d_1$, and $\aA_2=\d_1+O(\d)^1$
 by \eqref{cont-B2}.
}%

First assume $\tilde\xX_1\le\d_1^{\d_0}$. Then by recalling that $\d\ll\d_1\ll\d_0$, up to $O(\d)^1$, we have
\begin{eqnarray}
\label{Mndndnpp}
&\!\!\!\!\!\!\!\!\!\!\!\!\!\!\!\!\!\!\!\!\!\!\!\!&
{\rm(i)\ }
\Big|2-\frac1{\tilde X_1}\Big|\ge\tilde \xX_1^{-1}-2\ge\ell_1^{\d_0}-2>1,
\nonumber\\
&\!\!\!\!\!\!\!\!\!\!\!\!\!\!\!\!\!\!\!\!\!\!\!\!&
{\rm(ii)\ }\d_1=\aA_1\stackrel{{}^{\sc\rm\eqref{LetNSoOP----1-re-give}\,(i)}}{=}
\xX_2^{\ell_0}
\Big|2-\frac1{\tilde X_1}\Big|
\stackrel{{}^{\sc\rm\eqref{ImMpP}\,(i)}}{\ge}
\Big|2-\frac1{\tilde X_1}\Big|
\stackrel{{}^{\sc\rm\eqref{Mndndnpp}\,(i)}}{>}
1,
\end{eqnarray}
which is a contradiction.
Thus $\tilde\xX_1>\d_1^{\d_0}$. Then $\frac{\tilde\xX_1}{\aA_1^2}\ge\ell_1^{2-\d_0}$, and
up to $O(\d)^1$,
\begin{eqnarray}
\label{Mndndnpp++}
&\!\!\!\!\!\!\!\!\!\!\!\!\!\!\!\!\!\!\!\!\!\!\!\!\!\!\!\!\!\!\!\!\!\!\!\!\!&
{\rm(i)\ }
\Big|2-\frac{\tilde X_1}{A_1^2}\Big|>\ell_1^{2-\d_0}-2>1,
\ \ \ \ \ {\rm(ii)\ }\d_1=\aA_2
\stackrel{{}^{\sc\rm\eqref{We0o0o0o++1}}}{>}
\xX_2^{\ell_0}\Big |2- \frac{\widetilde  X_1}{A_1^2}\Big|
\stackrel{{}^{\sc\rm\eqref{ImMpP}\,(i),\,\eqref{Mndndnpp++}\,(i)}}{\ge}
1,\!\!\!\!
\end{eqnarray}
which is again a contradiction.
%
\NOUSE{%
\begin{eqnarray}
\label{Finanan}
&\!\!\!\!\!\!\!\!\!\!\!\!\!\!\!\!\!\!\!\!\!\!\!\!\!\!\!\!\!\!\!\!\!&
{\rm(i)\ }\tilde\xX_1
\stackrel{{}^{\sc\rm\eqref{ImMpP}\,(iv)}}{\le}
\aA_1^{-(1-\d_0^2)}=\d_1^{1-\d_0^2}<O(\d_1)^{\frac12},
\!\!\!\!\!\!\!\!\!\!\!\!\!\!\!\!\!\!\!\!\!\!\!\!\!\!\!\!\!\!\!\!\!
\\
&\!\!\!\!\!\!\!\!\!\!\!\!\!\!\!\!\!\!\!\!\!\!\!\!\!\!\!\!\!\!\!\!\!&
{\rm(ii)\ }
\aA_1\stackrel{{}^{\sc\rm\eqref{LetNSoOP----1}\,(i)}}{=}
\tilde \xX_1^3 \xX_2\Big |5 - \frac{4\widetilde X_1 }{X_2^{\ell_0}}\Big|
\stackrel{{}^{\sc\rm\eqref{Finanan}\,(i)}}{=}\tilde \xX_1^3 \xX_2\Big(5+O(\d_1)^{\frac12}\Big)
\stackrel{{}^{\sc\rm\eqref{Finanan}\,(i)}}{\ll}
\xX_2
.\NOUSE{,
\!\!\!\!\!\!\!\!\!\!\!\!\!\!\!\!\!\!\!\!\!\!\!\!\!\!\!\!\!\!\!\!\!
\nonumber\\
&\!\!\!\!\!\!\!\!\!\!\!\!\!\!\!\!\!\!\!\!\!\!\!\!\!\!\!\!\!\!\!\!\!&
{\rm(iii)\ }
\frac{\aA_1 \xX_2^{\ell_0-1}}{\tilde \xX_1^4}
\stackrel{{}^{\sc\rm\eqref{ImMpP}\,(iv)}}{\ge}
\aA_1^{5-4\d_0^2}\xX_2^{4\ell_0^2+\ell_0-1}
\stackrel{{}^{\sc\rm\eqref{ImMpP}\,(i)}}{\ge}
\aA_1^{5-4\d_0^2}=\ell_1^{5-4\d_0^2},
\!\!\!\!\!\!\!\!\!\!\!\!\!\!\!\!\!\!\!\!\!\!\!\!\!\!\!\!\!\!\!\!\!
\nonumber\\
&\!\!\!\!\!\!\!\!\!\!\!\!\!\!\!\!\!\!\!\!\!\!\!\!\!\!\!\!\!\!\!\!\!&
{\rm(iv)\ }
\Big|\ell_0^5+1-\frac{\ell_0^5A_1 X_2^{\ell_0-1}}{\widetilde X_1^4}\Big|
\stackrel{{}^{\sc\rm\eqref{Finanan}\,(iii)}}{=}
\frac{\aA_1 \xX_2^{\ell_0-1}}{\tilde \xX_1^4}\Big(\ell_0^5+O(\d_1)^1\Big)
=\frac{\eta_2\aA_1 \xX_2^{\ell_0-1}}{\tilde \xX_1^4},
\!\!\!\!\!\!\!\!\!\!\!\!\!\!\!\!\!\!\!\!\!\!\!\!\!\!\!\!\!\!\!\!\!
\nonumber\\
&\!\!\!\!\!\!\!\!\!\!\!\!\!\!\!\!\!\!\!\!\!\!\!\!\!\!\!\!\!\!\!\!\!&
{\rm(iv)\ }
\d_1^{\ell_0^8}=\aA_2
\stackrel{{}^{\sc\rm\eqref{We0o0o0o++1}}}{=}
\tilde \xX_1^{\ell_0^8 - 5\ell_0^5} \xX_2^{\ell_0^{20}}
\Big |\ell_0^5+1 - \frac{\ell_0^5A_1 X_2^{\ell_0-1}}{\widetilde X_1^4}\Big|
\!\!\!\!\!\!\!\!\!\!\!\!\!\!\!\!\!\!\!\!\!\!\!\!\!\!\!\!\!\!\!\!\!
\nonumber\\
&\!\!\!\!\!\!\!\!\!\!\!\!\!\!\!\!\!\!\!\!\!\!\!\!\!\!\!\!\!\!\!\!\!&
{\ \ \ \ \ \ }
\stackrel{{}^{\sc\rm\eqref{ImMpP}\,(iii)}}{=}
\eta_2\tilde \xX_1^{\ell_0^8 - 5\ell_0^5} \xX_2^{\ell_0^{20}}
\frac{\aA_1 \xX_2^{\ell_0-1}}{\tilde \xX_1^4}
=\eta_2\eta_1^{-\frac13(\ell_0^8 - 5\ell_0^5-4)}\aA_1(\eta_1\tilde\xX_1^3\xX_2)^{\frac13(\ell_0^8 - 5\ell_0^5-4)} \xX_2^{\ell_0^{20}+O(\d_0)^{-19}}
\!\!\!\!\!\!\!\!\!\!\!\!\!\!\!
\!\!\!\!\!\!\!\!\!\!\!\!\!\!\!\!\!\!\!\!\!\!\!\!\!\!\!\!\!\!\!\!\!
\nonumber\\
&\!\!\!\!\!\!\!\!\!\!\!\!\!\!\!\!\!\!\!\!\!\!\!\!\!\!\!\!\!\!\!\!\!&
\ge\ell_0\tilde \xX_1^{2\ell_0}\xX_2^{2\ell_0^2}-(\ell_0+1)
\stackrel{{}^{\sc\rm\eqref{ImMpP}\,(iv)}}{\ge}\ell_0\aA_1^{-\d_0^3}-(\ell_0+1)>1,
\!\!\!\!\!\!\!\!\!\!\!\!\!\!\!\!\!\!\!\!\!\!\!\!\!\!\!\!\!\!\!\!\!\nonumber\\
&\!\!\!\!\!\!\!\!\!\!\!\!\!\!\!\!\!\!\!\!\!\!\!\!\!\!\!\!\!\!\!\!\!&
{\rm(ii)\ }
\d_1\ \ =\ \aA_1\stackrel{{}^{\sc\rm\eqref{LetNSoOP----1}\,(i)}}{=}
\tilde\xX_1^{2\ell_0^2-1}\xX_2^{\ell_0^6}\Big|\ell_0+1 -\ell_0\widetilde X_1^{2\ell_0}X_2^{2\ell_0^2}\Big|
\stackrel{{}^{\sc\rm\eqref{ImMpP}\,(iv),\,\eqref{Finanan}\,(i)}}{\ge}\aA_1^{-\d_0^3(2\ell_0^2-1)}\xX_2^{\ell_0^6+O(\d_0)^{-5}}
\!\!\!\!\!\!\!\!\!\!\!\!\!\!\!\!\!\!\!\!\!\!\!\!\!\!\!\!\!\!\!\!\!\nonumber\\
&\!\!\!\!\!\!\!\!\!\!\!\!\!\!\!\!\!\!\!\!\!\!\!\!\!\!\!\!\!\!\!\!\!&
\phantom{{\rm(ii)\ }\ }
\stackrel{{}^{\sc\rm\eqref{ImMpP}\,(i)}}{\ge}
\aA_1^{-2\d_0+O(\d_0)^2}>\ell_1^{2\d_0+O(\d_0)^2}\gg\d_1.
}
\nonumber
\end{eqnarray}
Equ.~\eqref{Finanan}\,(ii) is a contradiction with \eqref{We0o0o0o}\,(i). 
}\NOUSE{
Then by \eqref{cont-B2}, $\aA_2=\d_1^{\frac{32}{5}+O(\d_0)^{\frac12}}$, and
up to $O(\d)^1$,
we have, where
$\b_1=|\frac32-\frac1{2A_3\widetilde X_1^{10}}|,\,\b_2=|\frac75-\frac{2}{5A_3\widetilde X_1^{10}}|$
[cf.~\eqref{x1-x2==}\,(iii),\,(iv)$\ssc\,$],
\begin{eqnarray}
\label{XX1==mamwmw=}
&&\!\!\!\!\!\!\!\!\!\!\!\!\!\!\!\!\!\!\!\!\!\!\!\!\!\!\!\!\!\!\!\!\!
{\rm(i)\ }\aA_3\stackrel{{}^{\sc\rm\eqref{LetNSoOP----1-redefine++}\,(iii)}}{=}
\b_1^{-1}\ell_1\tilde\xX_1^{100\ell_0+55}
\stackrel{{}^{\sc\rm\eqref{LetNSoOP----1-redefine++}\,(iv)}}{=}
\b_1^{-1}\ell_1\Big(\b_2\ell_1^{\frac{32}{5}+O(\d_0)^{\frac12}}\aA_3^{-10}\Big)^{\frac{100\ell_0+55}{800\ell_0+494}}
\nonumber\\
&&\!\!\!\!\!\!\!\!\!\!\!\!\!\!\!\!\!\!\!\!\!\!\!\!\!\!\!\!\!\!\!\!\!
\phantom{{\rm(i)\ }\aA_3\ \ \ \ }=\ \
\b_1^{-1}\b_2^{\frac18+O(\d_0)^1}\ell_1^{\frac95+O(\d_0)^{\frac12}}\aA_3^{-\frac{5}{4}+O(\d_0)^1},\ \ \ \implies
\nonumber\\
&&\!\!\!\!\!\!\!\!\!\!\!\!\!\!\!\!\!\!\!\!\!\!\!\!\!\!\!\!\!\!\!\!\!
{\rm(ii)\ }\aA_3=\b_3
\end{eqnarray}
Thus as above, we can obtain \eqref{Conndndm}, which is again a contradiction.
}%
\NOUSE{
$\aA_2\aA_3^{30}>\ell_1$ by \eqref{cont-B2},\,\eqref{ImMpP}\,(i).
Thus from \eqref{LetNSoOP----1++++++}\,(ii), we obtain, as in \eqref{msdmseme444}, the following,
\equa{----msdmseme444}{\dis{\rm(i)\ }
\frac{X_2^{90}}{A_3^{20}}=\frac{10}{9}+O(\d_1)^1,\ \ \ \
{\rm(ii)\ }\frac15+\frac{4X_2^{90}}{5A_3^{20}}\stackrel{{}^{\sc\rm\eqref{msdmseme444}\,(i)}}{=}\frac{49}{45}+O(\d_1)^1,
\ \ \ \
{\rm(iii)\ }\xX_2\stackrel{{}^{\sc\rm\eqref{msdmseme444}\,(i)}}{>}\aA_3^{\frac29}.
}
Then by \eqref{LetNSoOP----1++++++}\,(i),
\equa{smmsmsmene040430}{\dis
\ell_1=\aA_1=\Big(\frac{10}{9}\Big)^{\frac{91}{90}}\aA_3^{\frac29}
}
where $\b_1=\frac45,\,\b_2=\frac23$ (recalling that $0<\d=\ell^{-1}\ll\d_1=\ell_1^{-1}\ll1$),
\begin{eqnarray}
\label{Ne-compaaa}
\!\!\!\!\!\!\!\!\!\!\!\!\!\!\!\!\!\!\!\!\!\!\!\!&&
{\rm(i)\ }
\aA_3^{150}\tilde\xX_1^{150}
\stackrel{{}^{\sc\rm\eqref{C1B1}\,(i)}}{=}
(|B_1|\aA_1^{-1})^{\frac{25}{17}}+O(\d)^1
\stackrel{{}^{\sc\rm\eqref{Formmeme}}}{\le}
\d_1,
\nonumber\\
\!\!\!\!\!\!\!\!\!\!\!\!\!\!\!\!\!\!\!\!\!\!\!\!&&
{\rm(ii)\ }A_1\stackrel{{}^{\sc\rm \eqref{LetNSoOP----1}\,(ii)}}{=}
\frac{\a_1^6 A_3^{114}}{\widetilde X_1^{18} X_2^{66}}\Big(\frac15 + \frac4{5 A_3^{150}\widetilde X_1^{150}}\Big)
=\frac{4}{5A_3^{36}\widetilde X_1^{168} X_2^{66}}\Big(\frac{A_3^{150}\widetilde X_1^{150}}4+1\Big)+O(\d)^1
\nonumber\\
\!\!\!\!\!\!\!\!\!\!\!\!\!\!\!\!\!\!\!\!\!\!\!\!&&
\ \ \ \ \ \ \ \ \ \stackrel{{}^{\sc\rm\eqref{Ne-compaaa}\,(i)}}{=}
\frac{\b_1}{A_3^{36}\widetilde X_1^{168} X_2^{66}}\Big(1+O(\d_1)^1\Big)+O(\d)^1
\nonumber\\
\!\!\!\!\!\!\!\!\!\!\!\!\!\!\!\!\!\!\!\!\!\!\!\!&&
\ \ \ \ \
\stackrel{{}^{\sc\rm\eqref{equa-Case6-lemm}\,(iv),\,\eqref{LetNSoOP----1}\,(i)}}{=}
\frac{\b_1}{A_3^{102}\widetilde X_1^{102}}\Big(1+O(\d_1)^1\Big)+O(\d)^1
,
\nonumber\\
\!\!\!\!\!\!\!\!\!\!\!\!\!\!\!\!\!\!\!\!\!\!\!\!&&
{\rm(iii)\ }
A_2\ \ \ \ \ \,  \ \stackrel{{}^{\sc\rm \eqref{LetNSoOP----1}\,(iii)}}{=}\ \ \ \
\frac{\a_1^4A_3^{45}\widetilde X_1^6}{ X_2^{45}}\Big (\frac16 + \frac{5 A_1\widetilde X_1^{18} X_2^{66}}{
   6 \a_1^6A_3^{114}}\Big)
\stackrel{{}^{\sc\rm\eqref{LetNSoOP----1}\,(ii)}}{=}
\frac{\a_1^4 A_3^{45} \widetilde X_1^6}{X_2^{45}}\Big (\frac13+\frac{2}{3A_3^{150}\widetilde X_1^{150}}\Big)
\nonumber\\
\!\!\!\!\!\!\!\!\!\!\!\!\!\!\!\!\!\!\!\!\!\!\!\!&&
\phantom{{\rm(iii)\ }\ \ \ \ \ \ \ \ \, \ \ \ \ \  }
=\frac{2}{3A_3^{105}\widetilde X_1^{144} X_2^{45}}\Big(\frac{A_3^{150}\widetilde X_1^{150}}2+1\Big)+O(\d)^1
\nonumber\\
\!\!\!\!\!\!\!\!\!\!\!\!\!\!\!\!\!\!\!\!\!\!\!\!&&
\phantom{{\rm(iii)\ } }
\stackrel{{}^{\sc\rm\eqref{equa-Case6-lemm}\,(iv),\,\eqref{LetNSoOP----1}\,(i),\,\eqref{Ne-compaaa}\,(i)}}{=}
\frac{\b_2}{A_3^{150}\widetilde X_1^{99}}\Big(1+O(\d_1)^1\Big)+O(\d)^1,
\nonumber\\
\!\!\!\!\!\!\!\!\!\!\!\!\!\!\!\!\!\!\!\!\!\!\!\!&&
{\rm(iv)\ }\aA_1^{\frac{741}{41}+O(\d)^1}=\aA_1^{33-34\times\frac{18}{41}+O(\d)^1}
\stackrel{{}^{\sc\rm\eqref{cont-B2}}}{=}
\aA_1^{33}\aA_2^{-34}
\nonumber\\
\!\!\!\!\!\!\!\!\!\!\!\!\!\!\!\!\!\!\!\!\!\!\!\!&&
\phantom{{\rm(iv)\ }\ \ \ \ \ \ \ }
\stackrel{{}^{\sc\rm\eqref{Ne-compaaa}\,(ii),\,(iii)}}{=}
\b_1^{33}\b_2^{-44}
\aA_3^{34\times150-33\times102}\Big(1+O(\d_1)^1\Big)+O(\d)^1
\nonumber\\
\!\!\!\!\!\!\!\!\!\!\!\!\!\!\!\!\!\!\!\!\!\!\!\!&&
\phantom{{\rm(iv)\ }\ \ \ \ \ \ \ }
\ \ \ \ \ \
=\ \ \ \ \b_1^{33}\b_2^{-44}\aA_3^{1734}\Big(1+O(\d_1)^1\Big)+O(\d)^1
\stackrel{{}^{\sc\rm\eqref{ImMpP}}}{>}\b_1^{33}\b_2^{-44}\aA_1^{\frac{1734}{246}}\Big(1+O(\d_1)^1\Big)+O(\d)^1.
\end{eqnarray}
By \eqref{Ne-compaaa}\,(i),\,(ii),\,(iv), we obtain a contradiction as in
\eqref{Conndndm}.
}%
\NOUSE{%
Using the fact in \eqref{cont-B2} that $|A_1^{\d_0^3}A_2^{\d_0^3}|=1+O(\d)^2$, we immediately obtain from
\eqref{LetNSoOP----2}\,(i),\,(ii) and \eqref{ImMpP}\,(4) that $\tilde\xX_1=1+O(\d_0)^3$.
Then \eqref{LetNSoOP----2}\,(i) shows that $\xX_2=|A_1|^{\d_0^3}\big(1+O(\d_0)^3\big)$
Using \eqref{ImMpP}\,(4), we obtain from \eqref{LetNSoOP----2}\,(i),\,(ii) the following
\equa{mdmeme====}{\dis
|\widetilde X_1|=|A_2X_2|^{\d_0}\Big(+O(\d_0)^1\Big)=|A_1^{\d_0^4}X_2^{\d_0(1-\d_0^3)}|\Big(1+O(\d_0)^1\Big).
}
Thus we have [noting that $\ell_1\gg\ell_0$ (for instance we can assume $\ell_1>\ell_0^{\ell_0^5}$) and  $1-\ell_0^3a_2a_1^{-1}=3\d_0+O(\d_0)^2$ by \eqref{B1-ep}$\ssc\,$],
\equa{mem4m4m4m}{\dis
|X_2|=\Big|A_1A_2^{-\ell_0^3}\Big|\Big(1+O(\d_0)^1\Big)^{\ell_0^4}
\stackrel{{}^{\sc\rm\eqref{cont-B2},\,\eqref{B1-ep}}}{=}
\ell_1^{3\d_0+O(\d_0)^2}\Big(1+O(\d_0)^1\Big)^{\ell_0^4}\gg1,
}
a contradiction with \eqref{ImMpP}\,(2).
\NOUSE{\begin{rema}\rm\label{FirstIneq-rema}
We will see, as mentioned in Remark \ref{AboutX1}\,(iv), that
there will exist an inconsistence in  \eqref{ToSayas+1}
. As we already explained in Remark \ref{Antotm}\,(ii), the inconsistence arises from the fact that we choose different  $u,\tildev $ in \eqref{LetF3rst-final}.
\end{rema}
We have \equa{X1pSmall}{\dis
|X_1'|\stackrel{{}^{\sc\rm\eqref{C+LetNSoOP}\,(7),\,\eqref{equa-Case6-lemm}\,(iv)}}{\le}|Z|^{\ell_0^3}+O(\d)^1
\stackrel{{}^{\sc\rm\eqref{C+LetNSoOP}\,(1)}}{\le}\ell_1^{-2\ell_0^3 a_1^{-1}}+O(\d)^1
\stackrel{{}^{\sc\rm\eqref{B1-ep}}}{<}\d_1.
}}\NOUSE{
Thus $|X_1|=1+O(\d_0)^1$ by \eqref{ImMpP}\,(8).
By \eqref{1=bbb+1+1},\,$|A_2''|<\d_1$, thus $|A_2'''|=\frac43+O(\d_1)^1$ by \eqref{1=bbb}\,(ii), and $|A_2'|=\frac43+O(\d_1)^1$.
We obtain from \eqref{LetNSoOP----1}\,(iv)
that
\equa{Z33333}{\dis
|Z|^3=\frac{3|A_2A_1^{-\ell_0^4}|}{4}+O(\d_0)^1
=\frac{3|A_1|^{a_2a_1^{-1}-\ell_0^4}}{4}+O(\d_0)^1\gg1.}
By \eqref{1=bbb+1+1}, we have $|C_1|,|C_2|\prec_{\ell_1\ssc\,}1$.
Thus $|D_1|,|D_2|\sim_{\ell_1\ssc\,}1$ by \eqref{1=bbb}\,(ii),\,\eqref{1=bbb+1}\,(ii).
Then \eqref{1=bbb}\,(i),\,\eqref{1=bbb+1}\,(i) give the following
[cf.~\eqref{1-dSmallerThan}, where  the inequalities there now become ``\,$\sim_{\ell_1\ssc\,}$\,''],
\equa{1-dSmallerThan+1}{\dis
\Big|\frac{A_1A_2^2A_3^{2\d_0^2}}{X_1}\Big|\sim_{\ell_1\ssc\,}
\Big|\frac{A_1}{A_1'}\Big|\stackrel{{}^{\sc\rm\eqref{1=bbb}\,(i)}}{\sim_{\ell_1\ssc\,}}1\stackrel{{}^{\sc\rm\eqref{1=bbb+1}\,(i) }}{\sim_{\ell_1\ssc\,}}
\frac1{|A_2X_1|}\sim_{\ell_1\ssc\,}
\Big|\frac{A_3^{2\d_0^2}}{A_2X_1}\Big|.
}
This implies that $1\sim_{\ell_1\ssc\,}|A_1A_2^3|$, but by \eqref{B1-ep}, \eqref{cont-B2}
 [which obviously also holds for elements in $\ol V_0$
by \eqref{C+LetNSoOP}\,(i)$\ssc\,$],
\equa{wemamamamsn}{\mbox{ $|A_1A_2^3|\stackrel{{}^{\sc\rm\eqref{cont-B2}}}{\sim_{\ell_1\ssc\,}}
|A_1|^{1+3a_2a_1^{-1}}\stackrel{{}^{\sc\rm\eqref{B1-ep}}}{\sim_{\ell_1\ssc\,}}|A_1|^{-2+O(\d_0)^1}\prec_{\ell_1\ssc\,}1$,}} a contradiction.
}%
}%
This proves that the last strict inequality of \eqref{LetNSoOP}\,(i) holds for $(p_1,p_2)$.

\NOUSE{
By \eqref{ImMpP}\,(7), we have
$\frac{1+ \d_0^2}{2|X_1^{\ell_0^2}X_2^{\ell_0^3}|}\ge\frac{1+\d_0^2}{2}+O(\d)^1$ and
$\frac{(1- \d_0^2)|X_1^{\ell_0^2}X_2^{\ell_0^3}|}{2}\le\frac{1-\d_0^2}{2}+O(\d)^1$. Thus by
\eqref{ImMpP}\,(4), we have (i) below; by \eqref{ImMpP}\,(5) and \eqref{cont-B2}, we have (ii), then using
\eqref{ImMpP}\,(3),\,(6), we have,
noting from
\eqref{B1-ep} that $a_2a_1^{-1}=-10\ell_0^2(1+O(\d_0)^1)$,
\begin{eqnarray}
\label{nenemmmmm}
&\!\!\!\!\!\!\!\!\!\!\!\!\!\!\!\!\!\!\!\!\!\!\!\!&
{\rm(i)\ }
1\sim_{\ell_1\ssc\,}\a_1:=
\frac1{|A_1X_1^{2\ell_0^4+\ell_0^2}X_2^{3\ell_0^3+\d_0}|},\ \ \ \
{\rm(ii)\ }
1\sim_{\ell_1\ssc\,}\a_2=\Big|
\frac{X_1^{\d_0^4+\d_0^9}}{A_1^{a_2a_1^{-1}} X_2^{\d_0}}
\Big|,
\\\nonumber
&\!\!\!\!\!\!\!\!\!\!\!\!\!\!\!\!\!\!\!\!\!\!\!\!&
{\rm(iii)\ }
1\sim_{\ell_1\ssc\,}\a_1\a_2^{\ell_0^4}\sim_{\ell_1\ssc\,}
|A_1^{10\ell_0^4(1+O(\d_0)^1)}X_2^{-\ell_0^{12}(1+O(\d_0)^1)}|\sim_{\ell_1\ssc\,}
|A_1^{10\ell_0^4(1+O(\d_0)^1)}A_3^{-\ell_0^{12}(1+O(\d_0)^1)}|
.\end{eqnarray}
Noting that $a_3a_1^{-1}=1+O(\d_0)^1,$ we obtain from \eqref{nenemmmmm}\,(iii) that $|A_3|\sim_{\ell_1\ssc\,}|A_1|^{10\d_0+O(\d_0)^2}\prec_{\ell_1\ssc\,}|A_1|^{a_3a_1^{-1}}$, a contradiction with
\eqref{ImMpP}\,(3).
}\NOUSE{
By \eqref{A1A2Must+},\,\eqref{A1A2Must+11}, we see that $|C_1|,|C_2|\prec_{\ell_1\ssc\,}1$. Thus $|D_1|,|D_2|\sim_{\ell_1\ssc\,}1$ by \eqref{1=bbb}\,(ii),\,\eqref{1=bbb+1}\,(ii).
By \eqref{ImMpP}\,(7),\,(9), we have
\begin{eqnarray}
\label{Finaalalal}
&\!\!\!\!\!\!\!\!\!\!\!\!\!\!\!\!\!\!\!\!\!&
{\rm(i)\ }1\sim_{\ell_1\ssc\,}|D_1|\sim_{\ell_1\ssc\,}|\tildeX_1^{-60}A_1^{-1}|\sim_{\ell_1\ssc\,}
\sim_{\ell_1\ssc\,}
\Big|\frac{A_2^4 X_2^{21}}{A_1A_3^{30}\tildeX_1^{21}}\Big|\sim_{\ell_1\ssc\,}\frac1{A_1^{\frac{9}{89}}A_3^{51}
\tildeX_1^{42}},
\nonumber\\&\!\!\!\!\!\!\!\!\!\!\!\!\!\!\!\!\!\!\!\!\!&
{\rm(ii)\ }1\sim_{\ell_1\ssc\,}|D_2|\sim_{\ell_1\ssc\,}|\tildeX_1^{3}A_2^{-1}|\sim_{\ell_1\ssc\,}
\sim_{\ell_1\ssc\,}
\Big|\frac{\tildeX_1^9}{A_2A_3^{11} X_2}\Big|
\sim_{\ell_1\ssc\,}\Big|\frac{\tildeX_1^{10}}{A_1^{\frac{20}{89}}A_3^{10}
}\Big|,
\nonumber\\&\!\!\!\!\!\!\!\!\!\!\!\!\!\!\!\!\!\!\!\!\!&
{\rm(iii)\ }1\sim_{\ell_1\ssc\,}|D_1^{5}D_2{21}|\sim_{\ell_1\ssc\,}|A_1^{-\frac{465}{89}}A_3^{-465}|.
\end{eqnarray}
We obtain that $|A_3|\sim_{\ell_1\ssc\,}|A_1|^{-\frac1{89}}\prec_{\ell_1\ssc\,}|A_1|^{-\frac{25}{801}}$, which is a contradiction with \eqref{ImMpP}\,(3).
}\NOUSE{%
By \eqref{MSMA2}, we have $|A_2|\sim_{\ell_1\ssc\,}|A_1|^{-1+O(\d_0)^1}\preceq_{\ell_1\ssc\,}1$.
Thus by \eqref{1=bbb}, we deduce that $|A_3|\sim_{\ell_1\ssc\,}|A_2|$ and so $|\tildeX_1|^{4\ell_0}\sim_{\ell_1\ssc\,}|A_1A_2|$.
Using
\eqref{ImMpP}\,(5),\,(6), we have $|A_3|\sim_{\ell_1\ssc\,}|\tildeX_1|^{-4\ell_0-1}$, then using
\eqref{ImMpP}\,(7),\,(8), we  obtain the following,
\equa{X1====}{\!\!\!\!\!\!\!\!\!|\tildeX_1|^{4\ell_0}\sim_{\ell_1\ssc\,}|A_1A_2|\sim_{\ell_1\ssc\,}
|A_1^{-1}X_2^{\ell_0+2}A_3^{\ell_0+1}|
\sim_{\ell_1\ssc\,}|A_1^{-1}\tildeX_1^{4\ell_0(\ell_0+2)-(\ell_0+1)(4\ell_0+1)}|
=|A_1^{-1}\tildeX_1^{3\ell_0-1}|.\!\!\!\!\!\!\!\!\!
}
We obtain that $|\tildeX_1|\sim_{\ell_1\ssc\,}|A_1|^{-\frac{1}{\ell_0+1}}$, and so
 \equa{InSOOO}{\dis
|A_3|\sim_{\ell_1\ssc\,}|A_1|^{\frac{4\ell_0+1}{\ell_0+1}}\prec_{\ell_1\ssc\,}
|A_1|^{\frac{3\d_0}{5}+\frac{4\d_0^2}{25}},}
 which is a contradiction with \eqref{ImMpP}\,(3).
}\NOUSE{
Note that the first equalities of \eqref{X0sm,m-1}\,(i),\,(ii) always hold [without condition \eqref{x0==wo3499}$\ssc\,$].
Thus $|\a_1|\sim_{\ell_1\ssc\,}$
First observe that the  second term, denoted as $\a_2$, inside the brackets of the right-hand side of \eqref{ImMpP}\,(5.a) is equal to $\frac{9\tildeX_1^{10}A_1^5}{10}+O(\d)^3$ by using \eqref{ImMpP}\,(4.a). Then note from \eqref{MSM1111} and convention \eqref{MSmde33333+} that $|A_1|\prec_{\ell_1\ssc\,}1$ by the assumption in Step 2.
By \eqref{ImMpP}\,(7), we have that
$|\tildeX_1A_1^{\frac12}|\succeq_{\ell_1\ssc\,}1$.
If $|\tildeX_1A_1^{\frac12}|\sim_{\ell_1\ssc\,}1$, then $|\a_2|\sim_{\ell_1\ssc\,}1$ and as in
\eqref{X0sm,m-1}, we see that $\a_1\sim_{\ell_1\ssc\,}|A_1|^{\frac{11}4}\prec_{\ell_1\ssc\,}1$, which implies that the right-hand side of
\eqref{ImMpP}\,(5.a) is  $\prec_{\ell_1\ssc\,}1$, a contradiction. Thus
$|\tildeX_1A_1^{\frac12}|\succ_{\ell_1\ssc\,}1$, and so $|\a_2|\succ_{\ell_1\ssc\,}1$.
Then  \eqref{ImMpP}\,(5) shows that we must have $|\a_1\a_2|\sim_{\ell_1\ssc\,}1$, i.e.,
\begin{eqnarray}
\label{B2-IMMM+1}
&\!\!\!\!\!\!\!\!\!\!\!\!\!\!\!\!\!\!\!\!\!\!\!\!\!\!\!\!\!\!\!\!&
1\sim_{\ell_1\ssc}|\a_1\a_2|\sim_{\ell_1\ssc\,}|A_2^{-1}A_3^{-11}A_1^5X_2^{-5}|\sim_{\ell_1\ssc\,}|\tildeX_1^{11}X_2^6A_1^{\frac{21}{4}}|
,
\end{eqnarray}
where the last $\sim_{\ell_1}$ is obtained by using the facts that
$A_3\sim_{\ell_1\ssc\,}(\tildeX_1X_2)^{-1},$ $|A_2|\sim_{\ell_1\ssc\,}|A_1|^{-\frac14}$.
Further, we see from \eqref{A2====} that $|A_2|\preceq_{\ell_1\ssc\,}|A_1|^{\frac52}\prec_{\ell_1\ssc\,}1$. Thus \eqref{ToSayas+1}\,(h) with \eqref{ToSayas+2}\,(i) shows that
$A_3\sim_{\ell_1\ssc\,}A_2+O(A_2)^2\sim_{\ell_1\ssc\,}A_2$, and we obtain from \eqref{ToSayas+1}\,(e) that $|\tildeX_1|^{10}\sim_{\ell_1\ssc\,}|A_1A_2|$,
i.e. [using \eqref{ImMpP}\,(2),\,(6) and the fact that $|A_2|=|A_1|^{-\frac14}+O(\d)^3\ssc\,$],
\begin{eqnarray}
\label{Anenex1}
&\!\!\!\!\!\!\!\!\!\!\!\!\!\!\!\!\!\!\!\!\!\!\!\!\!\!\!\!\!\!\!\!\!\!\!\!\!\!\!\!\!\!\!\!\!\!\!\!\!\!\!\!\!&
1{\ssc}\sim_{\ell_1}|\tildeX_1^{-10}A_1A_2|\sim_{\ell_1\ssc\,}|\tildeX_1^{-10}A_2^{-1}A_3^{-11}X_2^{-5}|
\sim_{\ell_1\ssc\,}|\tildeX_1X_2^6A_2^{-1}|\sim_{\ell_1\ssc\,}
|\tildeX_1X_2^6A_1^{\frac14}|
.
\!\!\!\!\!\!\end{eqnarray}
Dividing \eqref{B2-IMMM+1} by \eqref{Anenex1} we obtain that $|\tildeX_1^{10}A_1^5|\sim_{\ell_1\ssc\,}1$, a contradiction.
\vskip5pt
}%
The above with \eqref{C+LetNSoOP} shows that \eqref{LetNSoOP} is satisfied by $(p_1,p_2)$. Thus \eqref{ToSayas+1} is satisfied by $(p_1,p_2)$,
which implies that $(p_1,p_2)\in V_0$. This proves that  $\ol V_0=V_0$.
\hfill$\Box$\vskip5pt

Now Proposition \ref{real00-inj}
follows from Lemmas \ref{V0NOT0}--\ref{Prop(ii)Holds} together with  the following facts.
When $(p_1,p_2)\in V_0$, by \eqref{ImMpP}, $x_1,x_2$ are bounded,
thus $V_0$ is bounded by Proposition \ref{pro-also}. Further by \eqref{SimMMSMS},\,\eqref{equa-Case6-lemm}\,(iii),\,
\eqref{A1-A2-cond},\,\eqref{ImMpP}
, we see that
 $x_1,x_2,x_2+y_2\ne0$, i.e., \eqref{-EiathA0} holds.
This completes the proof
.\hfill$\Box$

\appendix
\section{Another proof of surjectivity of $\pi_1$ 
}\label{sect7}
To see our original idea on how to obtain the surjectivity of $\pi_1$,  for your reference, below we would like to present our original  proof of surjectivity of $\pi_1$  by simply using  fundamental mathematical methods.

Assume conversely $\pi_1$ is not surjective. Fix any $(\xi_1,\xi_2)\in\C^2\bs\pi_1(V)$, and define
\equa{Equ-Theo-3}{L_{p_1,p_2}=|x_1-\xi_1|^2+|x_2-\xi_2|^2\mbox{ for }(p_1,p_2)=\big((x_1,y_1),(x_2,y_2)\big)\in V.}
Then $L_{p_1,p_2}>0$ for $(p_1,p_2)\in V$. Fix any $(\bar p_1,\bar p_2)\in V$, and denote
\equa{V-2}{V_3=\{(p_1,p_2)\in V\,|\,L_{p_1,p_2}\le L_{\bar p_1,\bar p_2}\}.}
Then $V_3\ne\emptyset.$ By \eqref{Equ-Theo-3} and Theorem \ref{Theo-2}, $V_3$ is compact in $\C^4$. Then the following proposition (which says that the continuous function $L_{p_1,p_2}$ does not have the minimal value on the compact set $V_3$) immediately gives a contradiction, which proves the  surjectivity of $\pi_1$.
\begin{prop}\label{Main-prop}
For any $(p_1,p_2)\in V_3$, in every neighborhood of $(p_{\ZeRo},p_{\OnE})$ there exists $(q_1,q_2)=\big((\dot x_1,\dot y_1),(\dot x_2,\dot y_2)\big)\in V_3$ such that
\equa{Equ-Theo-3+1}{L_{q_1,q_2}<L_{p_1,p_2}.}
\end{prop}
\noindent{\it Proof.~}Let $(p_{\ZeRo},p_{\OnE})=\big((x_{\ZeRo},y_{\ZeRo}),(x_{\OnE},y_{\OnE})\big)\in V$ and set
(and define $G_{\ZeRo},G_{\OnE}$ similarly)
\begin{eqnarray}
\label{Newf0F1++00}&\!\!\!\!\!\!\!\!\!\!\!\!\!\!\!\!\!\!\!\!\!\!\!\!\!\!&
F_{\ZeRo}= F(x_{\ZeRo}+\a_{\ZeRo}x,y_{\ZeRo}+y)-F(x_1,y_1),  \ \,
\ \ F_{\OnE}= F(x_{\OnE}+\a_{\OnE}x,y_{\OnE}+y)-F(x_2,y_2),\mbox{ \  where}\\[0pt]
\label{0+Newf0F1++00}&\!\!\!\!\!\!\!\!\!\!\!\!\!\!\!\!\!\!\!\!\!\!\!\!\!\!&
\a_{\ZeRo}=\left\{\begin{array}{cl}1&\mbox{if }x_{\ZeRo}=\xi_{\ZeRo},\\[0pt]
x_{\ZeRo}-\xi_{\ZeRo}&\mbox{else},\end{array}\right.\ \ \ \ \ \ \ \ \ \ \ \ \ \ \, \
\a_{\OnE}=\left\{\begin{array}{cl}1&\mbox{if }x_{\OnE}=\xi_{\OnE},\\[0pt]
x_{\OnE}-\xi_{\OnE}&\mbox{else.}\end{array}\right.
\end{eqnarray}
Define $q_{\ZeRo},q_{\OnE}$ accordingly [cf.~\eqref{q0q1} and \eqref{1++??+q0q1}${\ssc\,}$],
\equa{1+++q0q1+Add}{
q_{\ZeRo}:=(\dot x_{\ZeRo},\dot y_{\ZeRo})=(x_{\ZeRo}+\a_{\ZeRo} s\ep, y_{\ZeRo}+t\ep),\ \ \ \ \ q_{\OnE}:=(\dot x_{\OnE}, \dot y_{\OnE})=
( x_{\OnE}
+\a_{\OnE}u
\ep,y_{\OnE}+v\ep).}
We need to choose $(q_1,q_2)$ to be in $V$ [then it is automatically in $V_3$ when \eqref{Equ-Theo-3+1} holds].

Let $A_1$ be the unique invertible $2\times 2$ complex matrix such that for Jacobian pair
 \equa{NAnene}{\mbox{$(\bar F_i,\bar G_i):=(F_i,G_i)A_1$ for $i=1,2$, }}
 the linear parts of $\bar F_1,\bar G_1,\bar F_2,\bar G_2$, denoted as $\bar F_1\lin,\bar G_1\lin,\bar F_2\lin,\bar G_1\lin$, are of the following forms,
\equa{formmmamsm}{\bar F_1\lin=x,\ \ \bar G_1\lin=y,\ \ \ \ \ \bar F_2\lin=ax+by,\ \ \bar G_2\lin=cx+dy,}
for some $a,b,c,d\in\C$ with $ad-bc=1$.
Note from \eqref{Newf0F1++00},\,\eqref{1+++q0q1+Add} and \eqref{NAnene} that the equation \eqref{detmm} is equivalent to the equation
\equa{1+2=detmm}{\big(\bar F_1(s\ep,t\ep),\bar G_1(s\ep,t\ep)\big)=\big(\bar F_2(u\ep,v\ep),\bar G_2(u\ep,v\ep)\big).}
Then we can solve from \eqref{1+2=detmm} to obtain
%
\equa{AsIn2.2}{s=au+bv+O(\ep)^1.}
Note 
that $(x_{\ZeRo},x_{\OnE})\ne(\xi_{\ZeRo},\xi_{\OnE})$ by the assumption.
First suppose $x_{\ZeRo}\ne\xi_{\ZeRo},\,x_{\OnE}\ne\xi_{\OnE}$ (then $\a_{\ZeRo}=x_{\ZeRo}-\xi_{\ZeRo},\,\a_{\OnE}=x_{\OnE}-\xi_{\OnE}$).
In this case, by \eqref{Equ-Theo-3},\,\eqref{Equ-Theo-3+1}, we need to choose $u,v$ such   that,
\equa{2DDDbar1}{C_{\ZeRo}:=\b_{\ZeRo}|1+s\ep|^2+\b_{\OnE}|1+u\ep|^2-(\b_{\ZeRo}+\b_{\OnE})<0,} where
$\b_{\ZeRo}=|x_{\ZeRo}-\xi_{\ZeRo}|^2>0,\,  \b_{\OnE}=|x_{\OnE}-\xi_{\OnE}|^2>0$.
Using \eqref{AsIn2.2} in \eqref{2DDDbar1}, we immediately see (by comparing the coefficients of $\ep^1$) that if $b\ne0$ or $a\ne-{\b_{\ZeRo}}{\b_{\OnE}^{-1}}$, then we have a solution for \eqref{2DDDbar1}. Thus assume \equa{Assb==0}{b=0,\ \ \ \ \ \ a=-{\b_{\OnE}}{\b_{\ZeRo}^{-1}}\in\R_{\ne0}.}
Then we can write
\begin{eqnarray}
\label{4-3}&\!\!\!\!\!\!\!\!\!\!\!\!\!\!\!\!\!\!\!\!\!\!\!\!\!\!\!\!\!\!&
{\rm(i)\ }\bar F_{\ZeRo}=
\mbox{$\sum\limits_{i\ge2}a_iy^i+x\Big(1+\sum\limits_{i\ge1}\hat a_iy^i\Big)+\cdots, \ \ \bar F_{\OnE}= \sum\limits_{i\ge2}b_iz^i+ax\Big(1+\sum\limits_{i\ge1}\hat b_iz^i\Big)+\cdots$},
\\[-2pt]
\nonumber
&\!\!\!\!\!\!\!\!\!\!\!\!\!\!\!\!\!\!\!\!\!\!\!\!\!\!\!\!\!\!&
{\rm(ii)\, } \bar G_{\ZeRo}=y+\mbox{$\sum\limits_{i\ge2}c_iy^i+x\sum\limits_{i\ge1}c'_iy^i+\cdots,\ \ \ \ \ \bar G_{\OnE}=z+\sum\limits_{i\ge2}d_iz^i+x\sum\limits_{i\ge1}d'_iz^i+\cdots
,\ z=cx+d y,$\!\!\!\!\!}
\end{eqnarray}
for some $a_i,\hat a_i,b_i,\hat b_i,c_i,c'_i,d_i,d'_i\in\C$,  and where we regard $\bar F_{\OnE},\bar G_{\OnE}$ as polynomials on $x,z$ and we omit  terms with $x$-degree $\ge2$.
\begin{lemm}\label{6-q10}
There exists some $i\ge2$ such that 
{\mbox{$(a_i,c_i)\ne(b_i,d_i)$.}}
\end{lemm}\noindent{\it Proof.~}~Denote $(\bar F,\bar G)=(F,G)A_1$, and we use the same symbols with a bar to denote the associated elements corresponding to the pair $(\bar F,\bar G)$. Then by definition of Keller maps, we have
\equa{Keller-map}{\dis \bar \sigma(p)=\big(\bar F(p),\bar G(p)\big)=\big(F(p),G(p)\big)A_1=\sigma(p)A_1\mbox{ for }p\in\C^2.}
 Assume $(a_i,c_i)=(b_i,d_i)$ for $i\ge2$. Then by the bar version of \eqref{Newf0F1++00}, we  obtain  (and the like for $\bar G$),
\begin{eqnarray}
\label{ObataDDD}&&\!\!\!\!\!\!\!\!\!\!\!\!\!\!\!\!\!\!\!\!\!\!\!\!
\bar F(x_{\ZeRo}, y_{\ZeRo}+\kk)=\bar F_{\ZeRo}\big|_{(x,y)=(0,\kk)}=
 \bar F_{\OnE}\big|_{(x,z)=(0,\kk)}=\bar F( x_{\OnE},y_{\OnE}+d^{-1}\kk),
 \end{eqnarray}
i.e., $\bar \si(\hat p_{\ZeRo})=\bar \si(\hat p_{\OnE})$ and so $\si(\hat p_{\ZeRo})=\si(\hat p_{\OnE})$ by
\eqref{Keller-map} for
all $\kk\gg1$, where \equa{hat-ppppp}{\dis
\hat p_{\ZeRo}:=(\hat x_{\ZeRo},\hat y_{\ZeRo})=(x_{\ZeRo},y_{\ZeRo}+\kk),\ \ \ \ \hat p_{\OnE}:=(\hat x_{\OnE},\hat y_{\OnE})=(x_{\OnE},y_{\OnE}+d^{-1}\kk).}
Since $p_{\ZeRo}\ne p_{\OnE}$, we have $\hat p_{\ZeRo}\ne\hat p_{\OnE}$ when $\kk\gg1$, i.e., $(\hat p_{\ZeRo},\hat p_{\OnE})\in V$.
Then \eqref{hat-ppppp} gives that $\dH_{\hat p_{\ZeRo},\hat p_{\OnE}}\sim_{\ssc\kk\,}\kk
$ and  $|\hat y_{\OnE}|\sim_{\ssc\kk\,}\kk
\succ_{\ssc\kk\,} \dH_{_{\sc \hat p_{\ZeRo},\hat p_{\OnE}}}^{^{\sc \frac{\ssTH m}{\ssTH m+1}}}$, a contradiction with  \eqref{mqp1234-2}.
\hfill$\Box$
\begin{lemm}
\label{5-lemm3}Let $i_{\rZeRo}\ge2$ be the minimal number satisfying Lemma $\ref{6-q10}$.
We can assume
\equa{ai=bi=0}{\dis
{\rm(i)\ }a_i=0,\, i=2,...,2i_{\rZeRo},\  {\rm(ii)\ }b_i=0,\, i=2,...,2i_{\rZeRo},\
{\rm(iii)\ }c_i=d_i,\, i=2,...,i_{\rZeRo}-1,\ c_{i_{\rZeRo}}\ne d_{i_{\rZeRo}}.
}
\end{lemm}
%
\noindent{\it Proof.~}%
To solve $s$ from \eqref{1+2=detmm}, we can  replace $ \bar F_j$ by $\bar  F_j+\sum_{i=2}^{2i_{\rZeRo}}\b_i \bar G^{{\ssc\,}i}_j$ for some $\b_i\in\C$ and $j=\ZeRo,\OnE$
(
observe that $i_{\rZeRo}$ is still the minimal number satisfying  Lemma \ref{6-q10} after the replacement
by considering either $c_{i_{\rZeRo}}\ne d_{i_{\rZeRo}}$ or $c_{i_{\rZeRo}}=d_{i_{\rZeRo}},\,a_{i_{\rZeRo}}\ne b_{i_{\rZeRo}}$), thanks to the term $y$ in $ \bar G_{\ZeRo}$, we can then suppose \eqref{ai=bi=0}\,(i) holds.
Assume $b_k\ne0$ for some $k\le 2i_{\rZeRo}$. Take minimal such $k\ge2$.
Setting [noting from \eqref{4-3} that this amounts to setting $x=u\ep=\check u\ep^k$, $z=w\ep$ in $\bar F_{\OnE},\,\bar G_{\OnE}$, and setting
$x= s\ep,\,y=t\ep$ in $\bar F_{\ZeRo},\,\bar G_{\ZeRo}$, and letting $\bar F_{\ZeRo}=\bar F_{\OnE},\,\bar G_{\ZeRo}=\bar G_{\OnE}$ to solve $s,\,t$],
%
\equa{[MAMMS]}{u=\check u\ep^{k-1},\ \ \ \ \ \ \ \
\dis v=d^{-1}(w-cu),}  and regarding $\check u,w$ as new variables,
by \eqref{4-3}\,(i), we have
\equa{F0=F1==}{\dis\!\!\!\!\!\!\!\!\!
\bar F_{\OnE}|_{(x,z)=(\check u\scep^k,w\scep)}{\sc\!}={\sc\!}(b_kw^k{\sc\!}+{\sc\!}a\check u)\ep^k{\sc\!}+{\sc\!}O(\ep)^{k+1}{\sc\!}={\sc\!}
\bar F_{\ZeRo}(s\scep,t\scep){\sc\!}={\sc\!}O(\ep)^{k+1}{\sc\!}+{\sc\!}s\ep\Big(1{\sc\!}+{\sc\!}O(t\ep)^1\Big){\sc\!}+{\sc\!}O(s\ep)^2.
\!\!\!\!\!\!}
This shows that we have $s\ep=O(\ep)^k$ and the right-hand side of \eqref{F0=F1==} becomes $s\ep+O(\ep)^{k+1}$. Hence,
\equa{s=a}{s=(b_kw^k+ a\check u)\ep^{k-1}+O(\ep)^k.}
Using this, \eqref{Assb==0} and the first equation of \eqref{[MAMMS]} in
\eqref{2DDDbar1}, we obtain by choosing $\check u=0$ and $w\in\C$ with $(b_kw^k)\re<0$,
\begin{eqnarray*}
&&
C_{\rZeRo}{ }:=
{ }\b_{\ZeRo}\Big|1{ }+{ }(-\b_{\ZeRo}^{-1}\b_{\OnE}\check u
{ }+{ }b_kw^k)\ep^{k}{ }+{ }O(\ep)^{k+1}\Big|^2{ }+{ }\b_{\OnE}|1{ }+{ }\check u\ep^k|^2{ }
-{ }(\b_{\ZeRo}{ }+{ }\b_{\OnE})
\nonumber\\[0pt]&&\phantom{C_{\ZeRo}}
{ }={ }2\b_{\ZeRo}(b_kw^k)\re\ep^k{ }+{ }O(\ep)^{k+1}{ }<{ }0.
\end{eqnarray*}
This proves  Proposition \ref{Main-prop} in this case. Therefore, we can assume
\eqref{ai=bi=0}\,(ii) holds. Then we have \eqref{ai=bi=0}\,(iii) by Lemma \ref{6-q10}.\hfill$\Box$
\begin{lemm}\label{5-ficici}We have
 $\hat a_i=\hat b_i$ for $1\le i\le i_{\rZeRo}-2$ and $\kappa_{\rZeRo}:=\hat b_{i_{\rZeRo}-1}-\hat a_{i_{\rZeRo}-1}=i_{\rZeRo}(c_{i_{\rZeRo}}-d_{i_{\rZeRo}})\ne0$.
\end{lemm}
\noindent{\it Proof.~}%
For $1\le i\le i_{\rZeRo}-1$, by \eqref{4-3} we have,
%
\begin{eqnarray}
\label{CsssII}
&\!\!\!\!\!\!\!\!\!\!\!\!\!\!\!\!\!\!\!\!\!\!&
\mbox{$\hat a_i+\sum\limits_{{\rOnE}\le j<i}(j+1)\hat a_{j-i}c_{j+1}+(i+1)c_{i+1}=\Coeff\big(J(\bar F_{\ZeRo},\bar G_{\ZeRo}),x^0y^i\big)=0$, }
\nonumber\\[-4pt]
&\!\!\!\!\!\!\!\!\!\!\!\!\!\!\!\!\!\!\!\!\!\!&
\mbox{$\hat b_i+\sum\limits_{{\rOnE}\le j<i}(j+1)\hat b_{j-i}d_{j+1}+(i+1)d_{i+1}=\dis \frac{\Coeff\big(J(\bar F_{\OnE},\bar G_{\OnE}),x^0y^i\big)}{ad}=0$. }
\end{eqnarray}
Induction on
$i$ for $1\le i\le i_{\rZeRo}$, we obtain the lemma by \eqref{ai=bi=0}\,(iii).\hfill$\Box$
\vskip7pt
Now we set,
\equa{AMSMSM23333}{\dis u=u_{\rOnE}\ii\ep^{i_{\rZeRo}-1},\ \ \
\ \ \ \ \  v=d^{-1}(w-cu),}
for $u_{\rOnE}\in\R_{\ne0}$. As in \eqref{F0=F1==}, one can see that $s\ep=O(\ep)^{i_{\rZeRo}}$.
Thus by \eqref{4-3}\,(ii), we have
\begin{eqnarray}
\label{CsssII+1}
&\!\!\!\!\!\!\!\!\!\!\!\!\!\!\!\!\!\!\!\!\!\!\!\!\!\!\!\!\!\!&
\bar G_{\ZeRo}(s\ep,t\ep){\sc\!}={\sc\!}t\ep{\sc\!}+{\sc\!}\mbox{$\sum\limits_{i=2}^{i_{\rZeRo}}$}c_i(t\ep)^i{\sc\!}+{\sc\!}O(\ep)^{i_{\rZeRo}+1}
{\sc\!}={\sc\!}
\bar G_{\OnE}|_{(x,z)=(u_{\rOnE}\ii\scep^{i_{\rZeRo}},w\scep)}{\sc\!}={\sc\!}w\ep{\sc\!}+{\sc\!}
\mbox{$\sum\limits_{i=2}^{i_{\rZeRo}}$}d_i(w\ep)^i+O(\ep)^{i_{\rZeRo}+1}
.\end{eqnarray}
Hence $t\ep=w\ep+O(\ep)^{i_{\rZeRo}}$. Using this and \eqref{4-3}\,(i), we obtain,
\begin{eqnarray}
\label{CsssII+2}
&\!\!\!\!\!\!\!\!\!\!\!\!\!\!\!\!\!\!\!\!\!\!\!\!\!\!\!\!\!\!&
\bar F_{\ZeRo}(s\ep,t\ep){\sc}={\sc}s\ep\Big(
1{\sc}+{\sc}\mbox{$\sum\limits_{i=1}^{i_{\rZeRo}-1}$}\hat a_i(t\ep)^i\Big){\sc}+{\sc}O(\ep)^{2i_{\rZeRo}}
\nonumber\\[-8pt]
&\!\!\!\!\!\!\!\!\!\!\!\!\!\!\!\!\!\!\!\!\!\!\!\!\!\!\!\!\!\!&
\phantom{F_{\ZeRo}(s\ep,t\ep)}
{\sc}={\sc}
\bar F_{\OnE}|_{(x,z)=(u_{\rOnE}\ii\scep^{i_{\rZeRo}},w\scep)}{\sc}={\sc}au_{\rOnE}\ii\ep^{i_{\rZeRo}}\Big(1{\sc}+{\sc}
\mbox{$\sum\limits_{i=1}^{i_{\rZeRo}-1}$}\hat b_i(w\ep)^i\Big)+O(\ep)^{2i_{\rZeRo}}
,\end{eqnarray}
which with Lemma \ref{5-ficici} gives that $s\ep=au_{\rOnE}\ii\ep^{i_{\rZeRo}}(1+\kappa_{\ZeRo}w^{i_{\rZeRo}-1}\ep^{i_{\rZeRo}-1})+O(\ep)^{2i_{\rZeRo}}$.
Thus $C_{\rZeRo}$ defined in \eqref{2DDDbar1} becomes,
\begin{eqnarray}
\label{C0-Filll}
\!\!\!\!\!\!\!\!\!\!\!\!\!\!\!\!\!\!\!\!\!\!\!\!\!\!\!\!&&
C_{\rZeRo}=\b_{\ZeRo}\Big|1-\b_{\OnE}\b_{\ZeRo}^{-1}u_{\rOnE}\ii\ep^{i_{\rZeRo}}(1+\kappa_{\ZeRo}w^{i_{\rZeRo}-1}\ep^{i_{\rZeRo}-1})+O(\ep)^{2i_{\rZeRo}}\Big|^2+\b_{\OnE}|1+u_{\rOnE}\ii
\ep^{i_{\rZeRo}}|^2-\b_{\ZeRo}-\b_{\OnE}
\nonumber\\
\!\!\!\!\!\!\!\!\!\!\!\!\!\!\!\!\!\!\!\!\!\!\!\!\!\!\!\!&&
\phantom{C_{\rZeRo}}=2\b_{\OnE}u_{\rOnE}(\kappa_{\ZeRo}w^{i_{\rZeRo}-1})\im\ep^{2i_{\rZeRo}-1}+O(\ep)^{2i_{\rZeRo}},
\end{eqnarray}
which is negative if we choose $u_{\rOnE}=1$ and $w\in\C$ with $(\kappa_{\ZeRo}w^{i_{\rZeRo}-1})\im<0.$
This proves Proposition \ref{Main-prop} in case $\b_{\ZeRo}>0,\,\b_{\OnE}>0$.

Now if $x_{\ZeRo}=\xi_{\ZeRo}$ (then $x_{\OnE}\ne\xi_{\OnE}$), then the first term of $C_{\rZeRo}$ in \eqref{2DDDbar1} becomes $|s\ep|^2=O(\ep)^2$ and we can easily choose any $u$ with $u\re<0$ to satisfy that $C_{\rZeRo}<0$. Similarly, if $x_{\OnE}=\xi_{\OnE}$ (then $x_{\ZeRo}\ne\xi_{\ZeRo}$), then the second term of $C_{\rZeRo}$ in \eqref{2DDDbar1} becomes $|u\ep|^2=O(\ep)^2$ and using \eqref{AsIn2.2}, we can easily choose $u,v$ with $(au)\re<0$ and $v=0$ (in case $a\ne0$) or  with $u=0$ and $(bv)\re<0$ (in case $b\ne0$) to satisfy that $C_{\rZeRo}<0$.
This proves Proposition \ref{Main-prop}.\hfill$\Box$
%

\section{Another proof of Proposition \ref{real00-inj+1}}
\label{our-prop33}

To see our original idea on how to obtain Proposition \ref{real00-inj+1},  for your reference, below
 we wish to present our original  proof of Proposition \ref{real00-inj+1}.
\vskip5pt\noindent{\it Another proof of Proposition \ref{real00-inj+1}.~}~We can more precisely rewrite \eqref{i.e.1@such111that=2} as follows, for some $\tilde\a_i,\tilde\b_i\in\C$,
\begin{eqnarray}\label{i.e.1@such111that=2+A}
\!\!\!\!\!\!\!\!\!\!\!\!\!\!\!\!\!\!\!\!\!\!\!\!\!
&
\!\!\!\!\!\!\!\!\!\!&
{\rm(i)\ }C_{\rOnE}{\ssc}:={\ssc}
\big|1+(\tilde\a_{\rOnE}u+\tilde \a_2\tildev )\ep_1+(\tilde \b_{\rOnE}u^2+\tilde \b_2u\tildev +\tilde \b_3\tildev ^2)\ep_1^2
\big|
-1+O(\ep_1)^3\ge0,
\nonumber\\[0pt] 
\nonumber\!\!\!\!\!\!\!\!\!\!\!\!\!\!\!\!\!\!\!\!\!\!\!\!\!\!\!\!&
\!\!\!\!\!\!\!\!\!\!&
{\rm(ii)\, }C_2{\sc}:={\sc}
\kappa_{\ZeRo}'\big|1{\sc}+{\sc}(\tilde\a_3u{\sc}+{\sc}\tilde \a_4\tildev )\ep_1
{\sc}+{\sc}(\tilde \b_4u^2{\sc}+{\sc}\tilde \b_5u\tildev {\sc}+{\sc}\tilde \b_6\tildev ^2)\ep_1^2
\big|
{\sc}+{\sc}|1{\sc}+{\sc}u\ep_1|{\sc}
\\[-2pt] 
\!\!\!\!\!\!\!\!\!\!\!\!\!\!\!\!\!\!\!\!\!\!\!\!\!\!\!\!&
\!\!\!\!\!\!\!\!\!\!&
\phantom{{\rm(ii)\, }C_2:=}
+{\sc}\kappa'_{\OnE}|1{\sc}+{\sc}\tildev \ep_1|{\sc}-{\sc}
(\kappa'_{{\ZeRo}}{\sc}+{\sc}1+\kappa'_{\OnE}){\sc}+{\sc}O(\ep_1)^3{\sc}>{\sc}0.\!\!\!\!\!\!\!\!\!\!\!\!\!\!\!\!\!
\end{eqnarray}
First assume $\tilde\a_{\rOnE}\ne0$. 
Looking at \eqref{i.e.1@such111that=2+A}\,(i), we can take
\equa{take-u+A}{\mbox{ $u=-\tilde\a_{\rOnE}^{-1}\tilde\a_2\tildev +(\b_{\rOnE}\tildev ^2+\b_2w)\ep_1$
for some $\b_i,w\in\C$ with $w\re>0$,}}
so that $C_{\rOnE}$ has the form [the fact is that we can use \eqref{i.e.1@such111that=2+A}\,(i),\,\eqref{take-u+A} to solve $\b_1,\b_2$ precisely from the first equality below],
\equa{C1=====+A}{\mbox{$\dis \phantom{a}\ \ \ \ \ C_{\rOnE}=|1+w\ep_1^2|-1+O(\ep_1)^3=w\re\ep_1^2+O(\ep_1)^3>0$,}} i.e.,
\eqref{i.e.1@such111that=2+A}\,(i) holds.
%
Using \eqref{take-u+A} in \eqref{i.e.1@such111that=2+A}\,(ii), we see  that $C_2$ becomes  the following form, for some $\tilde\a_i\in\C$,
\begin{eqnarray}
\label{i.e.1@such111that=2++A}
\!\!\!\!\!\!\!\!\!\!\!\!\!\!\!\!\!\!\!\!\!\!\!\!\!
&
\!\!\!\!\!\!\!\!\!\!&
C_2=
\kappa_{\ZeRo}'\big|1+\tilde\a_5\tildev \ep_1+(\tilde \a_6\tildev ^2+\tilde\a_7 w)\ep_1^2\big|
+\big|1+\tilde\a_8\tildev \ep_1+(\tilde \a_9\tildev ^2+\tilde\a_{{\rOnE}0} w)\ep_1^2\big|
\nonumber\\[0pt]
\!\!\!\!\!\!\!\!\!\!\!\!\!\!\!\!\!\!\!\!\!\!\!\!\!
&
\!\!\!\!\!\!\!\!\!\!&
\phantom{C_2=}
+\kappa'_{\OnE}|1+\tildev \ep_1|-(\kappa'_{{\ZeRo}}+1+\kappa'_{\OnE})+O(\ep_1)^3>0.
\end{eqnarray}
Observe that for any $\a=\a\re+\a\im\ii\in\C$ [recall Convention \ref{conv1}\,(1)$\ssc\,$], we have, 
\begin{eqnarray}
\label{AOnnn+A}
\nonumber
&\!\!\!\!\!\!\!\!\!\!\!\!\!\!\!\!\!\!\!\!\!\!\!\!&
{\rm(i)\ }
|1+\a\ep_1|=\sqrt{(1+\a\re\ep_1)^2+(a\im\ep_1)^2}=1+\a\re\ep_1+\frac{(\a\im)^2}2\ep_1^2+O(\ep_1)^3,
\\&\!\!\!\!\!\!\!\!\!\!\!\!\!\!\!\!\!\!\!\!\!\!\!\!&
{\rm(ii)\ }A_0:=|1+\a\tildev ^2\ep_1^2|=1+\Big(\a\re\big((\tildev \re)^2-(\tildev \im)^2\big)-2\a\im\tildev \re\tildev \im\Big)\ep_1^2+O(\ep_1)^3.
\end{eqnarray}
We see that $C_2$ is an $O(\ep_1)^1$ element and we can easily compute
\begin{eqnarray}
\label{i.e.1@such111that=2+1+A}
\!\!\!\!\!\!\!\!\!\!\!\!\!\!\!\!\!\!\!\!\!\!\!\!\!
&
\!\!\!\!\!\!\!\!\!\!&
\Coeff(C_2,\ep_1)=
\kappa_{\ZeRo}'(\tilde\a_5\tildev )\re
+(\tilde\a_8\tildev )\re
+\kappa'_{\OnE}\tildev \re=(c_0\tildev )\re,\mbox{ \ where \ }c_0=\kappa'_{\ZeRo}\tilde\a_5+\tilde\a_8+\kappa'_{\OnE}.
\end{eqnarray}
Thus if  $c_{\rZeRo}\ne0$, we can always choose
$\tildev \in\C$ with $(c_{\rZeRo}\tildev )\re>0$ to satisfy \eqref{i.e.1@such111that=2++A}.

Hence assume $c_{\rZeRo}=0$. Then
$C_2$ in \eqref{i.e.1@such111that=2++A} becomes an $O(\ep_1)^2$ element.
Our purpose is to compute $\tilde \b$ defined in \eqref{t-beta+A} below. First using notation \eqref{nene34bb4b44b}, one can observe from \eqref{AOnnn+A}\,(ii) the facts that $\Coeff\big(A_0,(\tildev \re)^2\ep_1^2\big)=\a\re$ and $\Coeff\big(A_0,(\tildev \im)^2\ep_1^2\big)=-\a\re$, which imply the following important fact,
\equa{a0---fff+A}{\Coeff\big(A_0,(\tildev \re)^2\ep_1^2\big)+\Coeff\big(A_0,(\tildev \im)^2\ep_1^2\big)=\a\re-\a\re=0.} From this and \eqref{i.e.1@such111that=2++A},\,\eqref{AOnnn+A}\,(ii), one can observe
that $\tilde\a_6,\,\tilde\a_7,\,\tilde\a_9,\,\tilde\a_{{\rOnE}0}$ do not contribute
to $\tilde \b$  defined in \eqref{t-beta+A}. Then by \eqref{AOnnn+A}\,(i) and noting that $\big((\tilde\a_5\tildev )\im\big)^2=(\tilde\a_{5\rm\,re}v\im+\tilde\a_{5\rm\,im}v\re)^2$, we obtain
\begin{eqnarray}
\label{2ccoeff+A}
&\!\!\!\!\!\!\!\!\!\!\!\!\!\!\!\!\!\!\!\!\!\!\!\!\!&
\Coeff(2C_2,\ep_1^2)=\big(\kappa'_1(\tilde \a_{5\rm\,re})^2+(\tilde\a_{8\rm\,re})^2+\kappa'_2\big)(\tildev \im)^2
\nonumber\\
&\!\!\!\!\!\!\!\!\!\!\!\!\!\!\!\!\!\!\!\!\!\!\!\!\!&
\phantom{\Coeff(2C_2,\ep_1^2)=}
+
\big(\kappa'_1(\tilde\a_{5\rm\,im})^2+(\tilde\a_{8\rm\,im})^2\big)(\tildev \re)^2+2(\kappa'_1\tilde\a_{5\rm\,re}\tilde\a_{5\rm\,im}
+\tilde\a_{8\rm\,re}\tilde\a_{8\rm\,im})\tildev \re\tildev \im+\cdots,
\end{eqnarray}
where $\cdots$ means terms contributed by $\tilde\a_6,\,\tilde\a_7,\,\tilde\a_9,\,\tilde\a_{{\rOnE}0}$, which do not contribute the following,
therefore, we obtain the following crucial fact,
\begin{eqnarray}
\label{t-beta+A}
\!\!\!\!\!\!\!\!\!\!\!\!\!\!\!\!\!\!\!\!\!\!\!\!\!\!\!\!
&&
\tilde\b{\ssc}:=
{\ssc}\Coeff\big(2C_2,(\tildev \re)^2\ep_1^2\big){\ssc}+{\ssc}\Coeff\big(2C_2,(\tildev \im)^2\ep_1^2\big)
\nonumber\\[-0pt]\!\!\!\!\!\!\!\!\!\!\!\!\!\!\!\!\!\!\!\!\!\!\!\!\!&&\phantom{\tilde\b}{\ssc}
{\ssc}={\ssc}
\kappa'_1\big((\tilde \a_{5\rm\,re})^2+
(\tilde\a_{5\rm\,im})^2\big)+(\tilde\a_{8\rm\,re})^2+(\tilde\a_{8\rm\,im})^2+\kappa'_2\ge\kappa'_2>0.
\end{eqnarray}
This will ensure us to achieve our goal:
 we can  choose $\tildev $
with  $(\tildev \re)^2$ being sufficiently larger than $(\tildev \im)^2$ if
$\Coeff\big(C_2,(\tildev \re)^2\ep_1^2\big)>0$ or
 with  $(\tildev \im)^2$ being sufficiently larger than $(\tildev \re)^2$ if
$\Coeff\big(C_2,(v\im)^2\ep_1^2\big)>0$,
to guarantee that \eqref{i.e.1@such111that=2++A} holds (when $w$ is fixed).
 This proves
 Proposition {\rm\ref{real00-inj+1}}  for the case that $\tilde\a_{\rOnE}\ne0$.

Now assume $\tilde\a_{\rOnE}=0$. By the symmetry of $u$ and $\tildev $ in \eqref{i.e.1@such111that=2+A},
we may also assume $\tilde\a_2=0$.
Then we have one of the following,
\equa{mwm2m2m2+A}{\dis{\rm(i)\ } C_{\rOnE}=0,\ \ \ \ \mbox{ or } \ \ \ \
{\rm(ii)\ }C_{\rOnE}
=|1+g(u,\tildev )\ep_1^k|-1+O(\ep_1)^{k+1},}
for some nonzero homogeneous polynomial $g(u,\tildev )$ of $u,\tildev $ with degree $k\in\Z_{\ge2}$,
and we assume we have \eqref{mwm2m2m2+A}\,(ii) as
\eqref{i.e.1@such111that=2+A}\,(i) holds trivially for case \eqref{mwm2m2m2+A}\,(i).

In case $c_{\rOnE}:=\kappa'_{\ZeRo}\tilde\a_3+1\ne0$ [see \eqref{i.e.1@such111that=2+A}\,(ii)$\ssc\,$], we can solve the problem as follows:
Let $\a\in\C$ be determined later, and take $\tildev =\a u$. Then by \eqref{AOnnn+A}\,(i), as in \eqref{i.e.1@such111that=2+1+A}, we can compute from \eqref{i.e.1@such111that=2+A}\,(ii), \eqref{mwm2m2m2+A}\,(ii),
\begin{eqnarray}
\label{C-oeddd+A}
&\!\!\!\!\!\!\!\!\!\!\!\!\!\!\!\!\!\!\!\!&
{\rm(i)\ }\Coeff(C_2,\ep_1)=\Big(\big(\kappa'_1(\tilde\a_3+\tilde\a_4\a)+1+\kappa'_2\a\big)u\Big)\re=(\tilde c_1u)\re\mbox{ \ with \ }\tilde c_1=c_1+(\kappa'_1\tilde\a_4+\kappa'_2)\a,
\nonumber\\&\!\!\!\!\!\!\!\!\!\!\!\!\!\!\!\!\!\!\!\!&
{\rm(ii)\ }\Coeff(C_1,\ep_1^k)=\big(g_1(\a)u^k\big)\re\mbox{ \ with \ }g_0(\a)=g(1,\a).
\end{eqnarray}
Note that $g_0(\a)$ is a nonzero polynomial of $\a$ with degree $\le k$. We can always $\a\in\C$ such that $\tilde c_1\ne0$ and $g_0(\a)\ne0$.
Then we choose $u\in\C$ with $(\tilde c_{\rOnE}u)\re>0$ so that \eqref{i.e.1@such111that=2+A}\,(ii) holds by \eqref{C-oeddd+A}\,(i),
and further
$\big(g_0(\a)u^k\big)\re>0$ (this can be always done since $k\ge2$), i.e., $C_1>0$ by \eqref{mwm2m2m2+A}\,(ii),\,\eqref{C-oeddd+A}\,(ii).

If $c_2:=\kappa'_{\ZeRo}\tilde\a_4+\kappa'_{\OnE}\ne0$, we can
solve the problem symmetrically.

Now assume $c_{\rOnE}=c_2=0$, i.e., \equa{c1c3==+A}{\dis
\tilde\a_3=-\kappa_1'^{-1}\in\R,\ \ \ \ \ \tilde \a_4=-\kappa'_2\kappa'^{-1}_1\in\R.}
Then $C_2$ is an $O(\ep_1)^2$ element. One can  easily compute as in \eqref{t-beta+A},
\begin{eqnarray}
\label{<ASW,s+A}
\dis
\!\!\!\!\!\!\!\!\!\!\!\!\!\!\!\!\!\!\!\!\!\!\!\!&&
\Coeff\big(2C_2,(\tildev \re)^2\ep_1^2\big)+\Coeff\big(2C_2,(\tildev \im)^2\ep_1^2\big)
=\kappa'_1\a_4^2+\kappa_2'=\kappa'_2+\kappa'^2_2\kappa'^{-1}_1>0.
\!\!\!\!\!\!\!\!\!\!
\end{eqnarray}
If $g(u,\tildev )$ does not depend on $\tildev $ [i.e., $g(u,\tildev )=b''u^k$ for some $b''\in\C_{\ne0}$], then
we can first choose $u\in\C $ to satisfy that
$g(u,\tildev )\re=(b''u^k)\re>0$ then
 choose $\tildev \in\C$ with  $(\tildev \re)^2$ being sufficiently larger than $(\tildev \im)^2$ if
$\Coeff\big(C_2,(\tildev \re)^2\ep_1^2\big)>0$ or
 with  $(\tildev \im)^2$ being sufficiently larger than $(\tildev \re)^2$ if
$\Coeff\big(C_2,(\tildev \im)^2\ep_1^2\big)>0$, to guarantee that $\Coeff(C_2,\ep_1^2)>0$, i.e.,
\eqref{i.e.1@such111that=2+A}\,(ii) holds.

Thus assume $g(u,\tildev )$ depend on $\tildev $.
We set $\tildev =\a u$ with $\a,u\in\C$ being determined later
.
Then \eqref{mwm2m2m2+A}\,(ii) and \eqref{i.e.1@such111that=2+A}\,(ii) become the following forms,
for some $\tilde\a_{{\rOnE}1}\in\C$, and the non-constant  polynomial $g_{\rZeRo}(\a):=g(1,\a)$ of $\a$ with degree $\le k$
[where the number $\tilde\a:=-\kappa_{\ZeRo}'^{-1}(1+\a)$ in $C_2$ is obtained from
 $\tilde\a_3u+\tilde \a_4\tildev $ using \eqref{c1c3==+A}, and $\tilde\a_{11}=\tilde\beta_4+\a\tilde\b_5+\a^2\tilde\b_6\ssc\,$],
\begin{eqnarray}
\label{i.e.1@such111that=2+4+A}
\!\!\!\!\!\!\!\!\!\!\!\!\!\!\!\!\!\!\!
&
\!\!\!\!\!\!\!\!\!\!&
{\rm(i)\ }C_{\rOnE}=|1+g_{\rZeRo}(\a)u^k\ep_1^k|-1+O(\ep_1)^{k+1},
\nonumber\\[0pt]
\!\!\!\!\!\!\!\!\!\!\!\!\!\!\!\!\!\!\!\!\!\!\!\!\!
&
\!\!\!\!\!\!\!\!\!\!&
{\rm(ii)\, }
C_2{\sc\!}=
{\sc\!}\kappa'_{\ZeRo}\big|1{\sc\!}+{\sc\!}\tilde\a u\ep_1
{\sc\!}+{\sc\!}\tilde\a_{{\rOnE}1}u^2\ep_1^2\big|{\sc\!}+{\sc\!}|1{\sc\!}+{\sc\!}u\ep_1|
{\sc\!}+{\sc\!}\kappa'_{\OnE}
|1{\sc\!}+{\sc\!}\a u\ep_1|{\sc\!}-{\sc\!}(\kappa'_{\ZeRo}{\sc\!}+{\sc\!}1{\sc\!}+{\sc\!}\kappa'_{\OnE})
{\sc\!}+{\sc\!}O(\ep_1)^3{\sc\!}>{\sc\!}0.
\end{eqnarray}
Since $g_0(\a)$ is a non-constant polynomial of $\a$, we can choose $\a\in\C$ satisfying (i) below, then we can choose $\theta,u\in\C$ satisfying, where $\d>0$ is sufficiently small,
\begin{eqnarray}
\label{alpha-theta-u+A}
&\!\!\!\!\!\!\!\!\!\!\!\!\!\!\!\!\!\!\!\!\!\!\!\!&
{\rm(i)\ }g_0(\a)=1\mbox{  if  $k\equiv 0,1,2,7$ ({\rm mod\,}$8$), and \ $g_0(\a)=-1$ else,}
\nonumber\\
&\!\!\!\!\!\!\!\!\!\!\!\!\!\!\!\!\!\!\!\!\!\!\!\!&
{\rm(ii)\ } \mbox{$\theta=1$ if $(\tilde\a_{11})\re\le0$, and $\theta=-1$  else,}
\nonumber\\
&\!\!\!\!\!\!\!\!\!\!\!\!\!\!\!\!\!\!\!\!\!\!\!\!&
{\rm(iii)\ } u=\sqrt{2}e^{\frac{\theta(1-\d)\pi\ii}{4}},
\nonumber\\
&\!\!\!\!\!\!\!\!\!\!\!\!\!\!\!\!\!\!\!\!\!\!\!\!&
{\rm(iv)\ }u\re=1+O(\d)^1,\ \  u\im=\theta+O(\d)^1, \ \ \ (u\re)^2=(u\im)^2+O(\d)^1=1+O(\d)^1,
 \end{eqnarray}
where (iv) simply follows from (iii).
We can compute from \eqref{i.e.1@such111that=2+4+A} to obtain the following,
where  (i) follows from \eqref{alpha-theta-u+A}\,(i) and the facts that $\theta=\pm1$ and $\d>0$ is sufficiently small, while (ii) is obtained, as
 in \eqref{2ccoeff+A},  by  using  \eqref{AOnnn+A},\, \eqref{alpha-theta-u+A},\,(iv) and the fact from \eqref{alpha-theta-u+A}\,(ii) that $\theta(\tilde\a_{11})\re\ge0$\vspace*{-3pt},
\begin{eqnarray}\label{C2-ep2==+A}
&\!\!\!\!\!\!\!\!\!\!\!\!\!\!\!\!\!\!\!\!&
\dis{\rm(i)\ }
\Coeff(C_1,\ep_1^k)=\sqrt{2}\Big(g_0(\a)e^{\frac{k\theta(1-\d)\pi\ii}{4}}\Big)\re=\sqrt{2}g_0(\a)\cos\Big(\frac{k(1-\d)\pi\ii}{4}\Big)>0,
\nonumber\\&\!\!\!\!\!\!\!\!\!\!\!\!\!\!\!\!\!\!\!\!&
{\rm(ii)\ }
\Coeff(C_2,\ep_1^2)
=\frac12\Big(\kappa'_1\big((\tilde\a\re)^2+(\tilde\a\im)^2\big)+1+\kappa'_2\big((\a\re)^2+(\a\im)^2\big)\Big)-2\theta(\tilde\a_{11})\re+O(\d)^1
\nonumber\\&\!\!\!\!\!\!\!\!\!\!\!\!\!\!\!\!\!\!\!\!&
\phantom{{\rm(i)\ }\Coeff(C_2,\ep_1^2)}
\ge\frac12+O(\d)^1>0.
\end{eqnarray}
Thus $C_1>0$ and $C_2>0$ by \eqref{C2-ep2==+A}.
This proves Proposition \ref{real00-inj+1}.\hfill$\Box$

\section{A proposition provided by 
Bin Xu}
\label{sect8}
We wish to thank professor  Bin Xu from University of Science and Technology of China, who provides us the following
material.

\begin{prop}\label{prop1}
Assume that $V$ satisfies {\rm $($C1$)$--(C3),\,$($C4$)'$} stated in the beginning of section $3$. Then the projection $\pi_1:V\to\C^2$, $(p_1,p_2)\mapsto(x_1,x_2)$ for $(p_1,p_2)\in V$, is not surjective.
\end{prop}

For a proof, we need the following lemma in page 7 of \cite{Na92}.

\begin{lemm}
\label{lem:cov}
Let $X, Y$ be two Hausdorff, locally compact topological spaces, and let $\pi:X\to Y$ be both surjective and a local homeomorphism
$($i.e. for any $\a\in X$, there exists an open neighbourhood $X_\a$ of $\a$ such that $Y_\a=\pi(X_\a)$ is open in $Y$ and $\pi|_{X_\a}$ is a homeomorphism onto $Y_\a)$. Then $\pi$ is a finite covering $($i.e. for any $y_0\in Y$, there is an open neighbourhood $Y_{y_0}$ of $y_0$ such that $\pi^{-1}(Y_{y_0})$ is a disjoint union $\cup_{j=1}^n\, X_j$ of open sets $X_j$ of $X$ with the property that $\pi|_{X_j}$ is a homeomorphism onto $Y_{y_0}$ for all $1\leq j\leq n)$ if and only if it is proper $($i.e. for any compact set $K\subset Y$, the inverse image $\pi^{-1}(K)$ is compact in $X)$.
\end{lemm}

\noindent{\it Proof of Proposition \ref{prop1}.~}~Assume that $\pi_1$ is surjective.~%
By (C3),\,(C4)$'$,~$
\pi_1$ is a local homeomorphism.
By (C2)
, it is proper (see also Proposition \ref{pro-also}). Hence, by Lemma \ref{lem:cov}, $
\pi_1$ is a covering map. Since ${\Bbb C}^2$ is simply connected and $
V$ is connected, $
\pi_1$ is a holomorphic homeomorphism. Then there exist two holomorphic functions
$\phi_1, \,\phi_2$ on ${\Bbb C}^2$ such \vspace*{-3pt}that
\[
V=\Big\{\Big(\big(x_1,\phi_1(x_1,x_2)\big),\big(x_2,\phi_2(x_1,x_2)\big)\Big)\,\Big|\,(x_1,x_2)\in {\Bbb C}^2\Big\}.\]
By using both (C2)
\ and the Taylor developments of $\phi_1$ and $\phi_2$ on ${\Bbb C}^2$, we see that
$\phi_1,\phi_2$ are constant,
which contradicts that in $
V$, $x_1,y_1$ are locally expressed by  holomorphic functions on $x_2,y_2$\vspace*{-3pt}.\hfill$\Box$

{ \small\footnotesize \lineskip=2pt

}
\end{document}